\documentclass[11pt]{article} 
 
\usepackage{latexsym} 
\usepackage{amsfonts} 
\usepackage{amsmath} 
\usepackage{amsthm,amssymb,stackrel,esint} 
\usepackage{mathrsfs} 
\usepackage{dsfont,yfonts} 
\usepackage{bbold} 
\usepackage{caption} 
\usepackage{epsfig} 
\usepackage{float} 
\usepackage{subfigure,mathtools} 
\usepackage{psfrag} 
\usepackage{graphicx} 
\usepackage{epsfig} 
\usepackage{amsbsy} 
\usepackage{hyperref,color}
\usepackage{xcolor}
\usepackage{enumitem}
\usepackage{quoting}
\usepackage{scalerel,stackengine}
\stackMath
\newcommand\reallywidehat[1]{%
	\savestack{\tmpbox}{\stretchto{%
			\scaleto{%
				\scalerel*[\widthof{\ensuremath{#1}}]{\kern-.6pt\bigwedge\kern-.6pt}%
				{\rule[-\textheight/2]{1ex}{\textheight}}%WIDTH-LIMITED BIG WEDGE
			}{\textheight}% 
		}{0.5ex}}%
	\stackon[1pt]{#1}{\tmpbox}%
}
\DeclareMathAlphabet{\mathpgoth}{OT1}{pgoth}{m}{n}
\DeclareMathAlphabet{\mathpzc}{OT1}{pzc}{m}{it}

\DeclareMathOperator{\Val}{\matV} 
\DeclareMathOperator{\meas}{meas}

\newtheorem{theorem}{Theorem} 
\newtheorem*{prop*}{Theorem} 
\newtheorem{theo}[theorem]{Theorem} 
\newtheorem{coro}[theorem]{Corollary} 
\newtheorem{defi}[theorem]{Definition} 
\newtheorem{lemma}[theorem]{Lemma} 
\newtheorem{prop}[theorem]{Proposition} 
\newtheorem{rmk}[theorem]{Remark}

\newcommand{\zerarcounters}{\setcounter{equation}{0}\setcounter{theorem}{0}} 
 
\newcommand{\ZZZ}{\mathds{Z}} 
\newcommand{\CCC}{\mathds{C}} 
\newcommand{\NNN}{\mathds{N}} 
\newcommand{\QQQ}{\mathds{Q}} 
\newcommand{\RRR}{\mathds{R}} 
\newcommand{\TTT}{\mathds{T}} 
\newcommand{\uno}{\mathds{1}} 

\newcommand{\LLL}{\mathds{L}}

\newcommand{\calA}{{\mathcal A}} 
\newcommand{\BB}{{\mathcal B}} 
\newcommand{\CCCC}{{\mathcal C}} 
\newcommand{\DD}{{\mathcal D}} 
\newcommand{\calE}{{\mathcal E}} 
\newcommand{\calF}{{\mathcal F}} 
\newcommand{\matG}{{\mathcal G}}  %\newcommand{\calG}{{\mathcal G}} 

\newcommand{\calJ}{{\mathcal J}} 
 
\newcommand{\LL}{{\mathcal L}} 
\newcommand{\MM}{{\mathcal M}}

\newcommand{\calP}{{\mathcal P}} 
 
\newcommand{\RR}{{\mathcal R}} 
\newcommand{\ttR}{{\mathtt R}}

\newcommand{\SSSS}{{\mathbcalboondox R}\!\Theta} %\newcommand{\SSSS}{{\mathcal S}} 
\newcommand{\TT}{{\mathcal T}} 
 
\newcommand{\VV}{{\mathcal V}}

\newcommand{\calmC}{{\mathscr C}}

\newcommand{\calmK}{{\mathscr K}}

\newcommand{\calmN}{{\mathscr N}} 
 
\newcommand{\calmP}{{\mathscr P}}

\newcommand{\calmT}{{\mathscr T}}

\newcommand{\calmX}{{\mathscr X}}

\DeclareFontFamily{U}{BOONDOX-calo}{\skewchar\font=45 }
\DeclareFontShape{U}{BOONDOX-calo}{m}{n}{
  <-> s*[1.05] BOONDOX-r-calo}{}
\DeclareMathAlphabet{\mathcalboondox}{U}{BOONDOX-calo}{m}{n}
\DeclareMathAlphabet{\mathbcalboondox}{U}{BOONDOX-calo}{b}{n}

\newcommand{\gota}{{\mathfrak a}} 
\newcommand{\gotb}{{\mathfrak b}}

\newcommand{\gote}{{\mathfrak e}} 
\newcommand{\gotf}{{\mathfrak f}} 
\newcommand{\gotg}{{\mathfrak g}}

\newcommand{\gotj}{{\mathfrak j}}

\newcommand{\gotn}{{\mathfrak n}}

\newcommand{\gotr}{{\mathfrak r}}

\newcommand{\gotw}{{\mathfrak w}} 
 
\newcommand{\gotA}{{\mathfrak A}} 
\newcommand{\BBB}{{\mathfrak B}} 
\newcommand{\gotB}{{\mathbcalboondox B}} 
\newcommand{\gotC}{{\mathfrak C}} 
\newcommand{\gotD}{{\mathfrak D}}
\newcommand{\gD}{{\mathbcalboondox D}} 
\newcommand{\gotE}{{\mathfrak E}} 
\newcommand{\gotF}{{\mathfrak F}} 
 
\newcommand{\gotG}{{\mathfrak G}} 
\newcommand{\gotH}{{\mathfrak H}} 
\newcommand{\gotI}{{\mathfrak I}}

\newcommand{\gotN}{{\mathfrak N}}

\newcommand{\gotP}{{\mathbcalboondox A}} 
 
\newcommand{\gQ}{{\mathbcalboondox Q}} 
\newcommand{\gotR}{{\RR}} %\newcommand{\gotR}{{\mathbcalboondox R}} 
\newcommand{\fR}{{\mathfrak R}}
\newcommand{\gotS}{{\mathbcalboondox R}\!{\mathbcalboondox T}} %\newcommand{\gotS}{{\mathbcalboondox S}} 
\newcommand{\fS}{{\mathbcalboondox R}\fR} %\newcommand{\fS}{{\mathfrak S}}
\newcommand{\gotT}{{\mathbcalboondox T}} 
\newcommand{\fT}{{\mathfrak T}}

\newcommand{\gotW}{{\mathfrak W}} 
 
\newcommand{\fC}{{\mathbcalboondox C}}
 
\newcommand{\matC}{{\mathscr C}} 
\newcommand{\matD}{{\mathscr D}} 
 
\newcommand{\matF}{{\mathscr F}} 
\newcommand{\calG}{{\mathscr G}} %\newcommand{\matG}{{\mathscr G}} 
\newcommand{\matH}{{\mathscr H}}

\newcommand{\matK}{{\mathscr K}} 
\newcommand{\matL}{{\mathscr L}} 
\newcommand{\matM}{{\mathscr M}} 
 
\newcommand{\matO}{{\mathscr O}}

\newcommand{\matR}{{\mathscr R}}

\newcommand{\matU}{{\mathscr U}} 
\newcommand{\matV}{{\mathscr V}}

\newcommand{\matW}{{\mathscr W}}

\newcommand{\und}{\underline}

\newcommand{\ol}{\overline} 
 
\newcommand{\Fullbox}{{\rule{2.0mm}{2.0mm}}} 
 
\newcommand{\EP}{\hfill\Fullbox\vspace{0.2cm}} 
\newcommand{\prova}{\noindent\textit{Proof. }} 
\newcommand{\io}{\infty} 
\newcommand{\e}{\varepsilon} 
\newcommand{\al}{\alpha} 
\newcommand{\de}{\delta} 
\newcommand{\be}{\beta} 
\newcommand{\n}{\nu} 
\newcommand{\m}{\mu} 
\newcommand{\eps}{\epsilon}
\newcommand{\x}{\xi} 
 
\newcommand{\ka}{\kappa} 
\newcommand{\g}{\gamma} 
\newcommand{\om}{\omega} 
\newcommand{\h}{\eta} 
\newcommand{\ze}{\zeta} 
\newcommand{\la}{\lambda} 
\newcommand{\f}{\varphi} 
\newcommand{\s}{\sigma} 
 
\newcommand{\del}{\partial}

\newcommand{\av}[1]{\langle #1 \rangle}

\newcommand{\starF}{\calF^*}

\newcommand{\der}{{\rm d}} 
\newcommand{\ii}{{\rm i}} 
\newcommand{\To}{{\mathring{T}}} 
\newcommand{\TTo}{{\mathring{\TT}}} 
\newcommand{\PPTTo}{{\mathring{\calP}_{\TT}}} 
\newcommand{\PPTo}{{\mathring{\calP}_{T}}} 

\newcommand{\thetao}{{\mathring{\vartheta}}} 
\newcommand{\thetabo}{{\mathring{\breve\vartheta}}}

\newcommand{\jap}[1]{\langle #1 \rangle}
\newcommand{\nmax}{\ol{n}}
\newcommand{\nmin}{\underline{n}}
\newcommand{\pow}{q}
\newcommand{\Uperp}{U_{\!\perp}}

\makeatletter
\newcommand{\oset}[3][0ex]{%
  \mathrel{\mathop{#3}\limits^{
    \vbox to#1{\kern-1\ex@
    \hbox{$\scriptstyle#2$}\vss}}}}
\makeatother

\def\st#1{\oset{*}{#1}}

\def\tilde#1{\widetilde{#1}}
 
\def\ins#1#2#3{\vbox to0pt{\kern-#2 \hbox{\kern#1 #3}\vss}\nointerlineskip}

\begin{document}

%%%%%%%%%%%%%%%%%%%%%%%%%%%%%%%%%%%%%%%%%%%%%%%%%%%%%%%%%%%%%%%%%%%%%%%%%% 
%%%%%%%%%%%%%%%%%%%%%%%%%%%%%%%%%%%%%%%%%%%%%%%%%%%%%%%%%%%%%%%%%%%%%%%%%% 
\title{\textbf{Almost-periodic solutions to the NLS equation \\ with smooth convolution potentials}}
%%%%%%%%%%%%%%%%%%%%%%%%%%%%%%%%%%%%%%%%%%%%%%%%%%%%%%%%%%%%%%%%%%%%%%%%%% 
 %%%%%%%%%%%%%%%%%%%%%%%%%%%%%%%%%%%%%%%%%%%%%%%%%%%%%%%%%%%%%%%%%%%%%%%%%% 

\author{\textbf{Livia Corsi, Guido Gentile, Michela Procesi}\\
\\
\small Dipartimento di Matematica e Fisica, Universit\`a Roma Tre, Roma, 
00146, Italy\\
\small e-mail: livia.corsi@uniroma3.it,
guido.gentile@uniroma3.it,  michela.procesi@uniroma3.it}

\date{\today} 

\maketitle

%%%%%%%%%%%%%%%%%%%%%%%%%%%%%%%%%%%%%%%%%%%%%%%%%%%%%%%%%%%%%%%%%%%%%%%%%% 
%%%%%%%%%%%%%%%%%%%%%%%%%%%%%%%%%%%%%%%%%%%%%%%%%%%%%%%%%%%%%%%%%%%%%%%%%% 
\begin{abstract} 
We consider the one-dimensional NLS equation with a convolution potential and a quintic nonlinearity. We prove that,
for most choices of potentials with polynomially decreasing Fourier coefficients, there exist  almost-periodic solutions 
in the Gevrey class with frequency %$\om_j$ 
satisfying a Bryuno non-resonance condition.
This allows convolution potentials of class $C^p$, for any integer $p$: 
as far as we know this is the first result where the regularity of the potential is arbitrarily large and not compensated
by a corresponding smoothing of the nonlinearity.
\smallskip


\smallskip

\end{abstract} 
%%%%%%%%%%%%%%%%%%%%%%%%%%%%%%%%%%%%%%%%%%%%%%%%%%%%%%%%%%%%%%%%%%%%%%%%%% 
%%%%%%%%%%%%%%%%%%%%%%%%%%%%%%%%%%%%%%%%%%%%%%%%%%%%%%%%%%%%%%%%%%%%%%%%%% 

 \tableofcontents
 
%%%%%%%%%%%%%%%%%%%%%%%%%%%%%%%%%%%%%%%%%%%%%%%%%%%%%%%%%%%%%%%%%%%%%%%%%% 
%%%%%%%%%%%%%%%%%%%%%%%%%%%%%%%%%%%%%%%%%%%%%%%%%%%%%%%%%%%%%%%%%%%%%%%%%% 
\zerarcounters 
\section{Introduction} 
\label{intro} 
%%%%%%%%%%%%%%%%%%%%%%%%%%%%%%%%%%%%%%%%%%%%%%%%%%%%%%%%%%%%%%%%%%%%%%%%%% 
%%%%%%%%%%%%%%%%%%%%%%%%%%%%%%%%%%%%%%%%%%%%%%%%%%%%%%%%%%%%%%%%%%%%%%%%%% 
  
In the analysis of Hamiltonian PDEs  modelling wave propagation in spatially confined systems,
a useful paradigm is to investigate the existence of invariant manifolds, particularly finite- and infinite-dimensional invariant tori.
Typically one studies nonlinear PDEs on a compact manifold with an elliptic fixed point at zero, where the linear part $L$
has a purely imaginary pure point spectrum. Thus,  in the basis of eigenfunctions of $L$, 
the PDE is represented by an infinite chain of harmonic oscillators coupled by the nonlinearity:
all solutions of the linear system lie on invariant tori and one may wonder
which ones -- if any -- survive the onset of the nonlinearity.
   
If one considers a finite truncation of the problem, standard KAM theory ensures that,
under some non-degeneracy assumption, most tori survive and most solutions are quasi-periodic in time,
with frequency vectors satisfying some sufficiently strong irrationality condition.
Then a natural question is whether one can give similar statements in the infinite-dimensional case.
The first results in this direction date back to the last decade of the last century, and deal with the
existence of periodic and quasi-periodic solutions to some fundamental nonlinear PDEs,
first in space dimension $D=1$~\cite{K,K2,K3,W,CW,P,KP,B1,BBM},
then in higher space dimension~\cite{B2,EK,BB,GY,GXY,BCP,CM}.
However, such solutions are in no way typical for infinite-dimensional systems:
even in integrable models, typical solutions are supported on infinite-dimensional tori, and hence are \emph{almost-periodic} in time --
here and henceforth, following Bohr \cite{bohr,bohr1},  the set of almost-periodic functions on a Banach space $(X,\|\cdot\|_X)$ 
is meant as the closure w.r.t.~the uniform topology of the set of trigonometric polynomials with values in $X$.

With this in mind, to extend KAM theory to the infinite-dimensional case, one is interested in inquiring
whether ``typical" solutions of nearly integrable PDEs  are almost-periodic. In this context,
a crucial point is the choice of the phase space and of the measure defined on it.
Unfortunately, up to now, there are no results on this question: the only available results
are on the simpler problem of existence of ``some'' almost-periodic solutions,
without even attempting at studying which fraction of phase space they fill.
In fact, there are relatively few results on the subject,
all on somewhat cooked-up models.

A first main difficulty is that the small divisors one has to deal with are expected to be ``extremely small'':
it is well known that in constructing finite-dimensional invariant tori the non-resonance conditions
become weaker and weaker as the dimension increases.
To simplify the small divisor problem, imposing -- say -- some strong Diophantine condition
on the eigenvalues of the linear part, it is customary to consider parameter families of PDEs
with as many free parameters as the dimension of the invariant tori that one wishes to construct. 
In the case of the nonlinear Schr\"odinger equation or, more generally, PDEs whose linear part is constant coefficients,
the simplest way of introducing such parameters is through a Fourier multiplier operator~\cite{CP,Bgafa,Po,Bjfa,GX,BMP2}.

Before stating our result, let us describe more precisely the problem we aim to study, and 
briefly review the literature on the subject.

We  consider a simplified model for the nonlinear Schr\"odinger (NLS) equation on the circle with \emph{external parameters}, i.e.
\begin{equation}\label{nls}
\ii u_t - u_{xx} + T_V u + \e |u|^4u =0 ,\qquad x\in\TTT := \RRR / 2\pi \ZZZ ,
\end{equation}
where $\e$ is a small real parameter and $T_V$ is the \emph{Fourier multiplier operator} with Fourier multiplier (or symbol) $V\!:\ZZZ\to\RRR$.
The latter means that, if we write
\begin{equation}\label{giustoper}
u(x) = \sum_{j\in\ZZZ} u_j  e^{\ii jx} ,  
\end{equation}
then we have
\begin{equation}\label{giustoper1}
(T_V u) (x) :=  \sum_{j\in\ZZZ} V_ju_j \,  e^{\ii jx} .
\end{equation}
We consider a quintic nonlinearity, as in ref.~\cite{Bjfa}, to deal with an explicit non-integrable case,
but we expect our result to extend to any
equation with $|u|^4$ replaced by an analytic function of $|u|^2$;
however this would introduce further technical intricacies
which would only make the computations more cumbersome without really shedding more light on the relevant features.

Note that $T_V$ in \eqref{nls} is a bounded operator if $\{V_j\}_{j\in\ZZZ}\in \ell^\io(\ZZZ,\RRR)$.
If suitable decay properties are assumed on $\{V_j\}_{j\in\ZZZ}$, then one can write $T_V u= \calF^{-1}(V) \ast u$,
where $\calF$ is the Fourier transform operator and $\ast$ denotes convolution:
in that case, one says that \eqref{nls} describes the NLS equation with a \emph{convolution potential} $\calF^{-1}(V)$.
For our purposes, it is convenient to identify the Fourier multiplier $V$ with the set of external parameters $V=\{V_j\}_{j\in\ZZZ}$;
for simplicity's sake, we refer to $V$ as the \emph{potential} \footnote{Strictly speaking,
$V$ is rather either the Fourier transform of the potential -- when the latter exists -- or the Fourier multiplier.} of the NLS equation \eqref{nls}.  
	
In dealing with infinite-dimensional tori, where one expects to need modulation of infinitely many parameters,
the regularity of the potential plays a preeminent r\^ole: less regularity implies stronger modulation and
hence stronger non-resonance conditions, and this in turn entails results which are easier to obtain, but weaker.
A complementary point of view, again with the aim of simplifying the problem, is to look for infinite-dimensional tori
which are ``very close'' to a finite-dimensional one, for instance by considering a NLS equation
with a forcing term which is almost-periodic in time and whose Fourier coefficients decay very fast.
This corresponds to proving the existence of almost-periodic solutions with very high regularity.

In the light of the considerations above, it is not surprising that the first results in the field provided
very regular solutions for model PDEs, with potentials with very low regularity in the case of the NLS equation.
In this regard, we recall the results obtained by Bourgain and by P\"oschel, more precisely, the existence of:
\begin{itemize}[topsep=1ex]
\itemsep0em
\item at least super-exponentially decaying almost-periodic solutions to the nonlinear
wave equation with analytic multiplicative potentials~\cite{Bgafa};
\item at least super-exponentially decaying solutions to the NLS equation with smoothing nonlinearities and 
Fourier multipliers in $\ell^2(\ZZZ,\RRR)$~\cite{Po};
\item Gevrey solutions to the NLS equation with Fourier multipliers in $\ell^\infty(\ZZZ,\RRR)$~\cite{Bjfa}.
\end{itemize}

Informally, all the results of the papers mentioned above can be rephrased as follows.

\vspace{-0.1cm}
\begin{quoting}
\textit{Consider a Hamiltonian PDE with external parameters and with an elliptic fixed point at zero. 
Fix an invariant torus of the linearized equation which satisfies the very special condition that the corresponding actions 
admit some upper and lower bounds.\footnote{
For instance, Bourgain in ref.~\cite{Bjfa} requires the actions of the considered torus to be such that, with the notation \eqref{giustoper},
$r e^{-s\sqrt{|j|}} \le |u_j|^2 \le 2r e^{-s\sqrt{|j|}}$ for all $j\in\ZZZ$ and some $r,s>0$.}
Then for a positive measure set of external parameters, 
there exists  an almost-periodic solution to the non-linear PDE  supported 
on an infinite-dimensional invariant torus which is close to that of the linear equation.}
\end{quoting}
\vspace{-0.1cm}

We note here that, so far, in the study of almost-periodic solutions,
multiplicative potentials have been handled only in the case of the nonlinear wave equations, whereas
all results in the case of the NLS equation, unless one allows a smoothing nonlinearity~\cite{Po},
require at best a convolution potential.
Since the aforementioned seminal results there have been some advances both with the aim of strengthening the regularity assumptions
on the Fourier multipliers and in the direction of finding less regular solutions -- a phenomenon distinctive of infinite-dimensional systems,
which may occur already in the integrable cases~\cite{KM,GK,GKT}.
For the first kind of results see for instance refs.~\cite{GX} and \cite{BMP2}, which extend the results of ref.~\cite{Bjfa}
to a positive measure set of analytic and Gevrey, respectively, almost-periodic solutions,
while we refer to ref.~\cite{BMP3} for the second topics. Finally we mention ref.~\cite{CY2} for the wave equation and
ref.~\cite{CLSY},
where nonlinear stability is studied for the infinite-dimensional tori of ref.~\cite{Bjfa}, and
refs.~\cite{RLZ,RL1,Liu,CMP}, where almost-periodically forced PDEs are considered.

In the present paper we prove the existence of solutions with Gevrey regularity for the  NLS equation \eqref{nls}
with a convolution potential of arbitrarily high regularity. 
To describe more in detail our result, let us introduce the scale of weighted sequence spaces
\[
\ell^{N,\io} = \ell^{N,\io}(\ZZZ,\RRR) := \Big\{  \{V_j\}_{j\in\ZZZ} : V_j\in\RRR,\  
\|V\|_{N,\io}:=\sup_{j\in\ZZZ}|V_j|\jap{j}^N <\io \Big\},
\]
where $N\in\ZZZ_+$ and
\begin{equation} \label{jap}
\jap{j}:=\max\{1,|j|\} ,
\end{equation}
and consider \eqref{nls} with $V\in \ell^{N,\io}$;
note that, for any $N$ fixed, the space $\ell^{N,\io}$ is trivially isomorphic to $\ell^\io(\ZZZ,\RRR)$.
Then an informal statement of our result reads as follows.

\vspace{-0.1cm}
\begin{quoting}
\textit{Consider equation \eqref{nls} and fix an invariant torus of the linearized equation.
Then, for a positive probability measure set of external parameters in $\ell^{N,\io}$,
there exists  an almost-periodic solution to the non-linear PDE  supported 
on an infinite-dimensional invariant torus which is close to that of the linear equation.}
\end{quoting}
\vspace{-0.1cm}

The probability measure on $\ell^{N,\io}$ is defined as follows. Given any normed space $(X,\|\cdot\|_X)$, for any $\rho>0$ let $\matU_\rho(X)$ denote
the open ball  of radius $\rho$ centered at zero in $X$, and let $\ol{\matU}_\rho(X)$ denote the closure  of $\matU_\rho(X)$.
If $X=\ell^{N,\io}$ we endow $\ol{\matU}_\rho(\ell^{N,\io})$ with
the probability measure induced by the product measure on the ball $\ol{\matU}_\rho(\ell^\io)$. This means that,
for any
\[
 A = \prod_{j\in\ZZZ} A_j\subseteq  \prod_{j\in\ZZZ} [-\jap{j}^{-N}\rho,\jap{j}^{-N}\rho] ,
\]
we define the probability measure of $A$ as
\begin{equation} \label{eq:measure}
\meas(A) := \lim_{h\to\io} \prod_{j=-h}^{h} \meas (A_j) ,\qquad \meas(A_j)=\frac{m(A_j)}{2\jap{j}^{-N}\rho},
\end{equation}
with $m(\cdot)$ denoting the Lebesgue measure.

Let us now specify the type of the almost-periodic functions we consider. 
We regard an almost-periodic function $t\mapsto f(t)\in X$, with $(X,\|\cdot\|_X)$ a Banach space,
as the trace of a function of infinitely many angles $\f\in\TTT^\ZZZ$
computed at $\f=\om t$, for some $\om\in\RRR^\ZZZ$; in this case we say that $\om$ is the \emph{frequency vector} of $f$.
In other words, setting\footnote{Here and in the following $\|\cdot\|_p$ denotes as usual the $\ell^p$ norm;
in particular $\|\cdot\|_\io$ is the sup-norm. The subscript $f$ stands for ``finite support''.}
$\ZZZ^{\ZZZ}_{f}:=  \{ \nu \in \ZZZ^\ZZZ : \|\nu\|_1 < \infty\}$ and writing, for any $\nu\in\ZZZ^{\ZZZ}_{f}$ and any $\f\in\TTT^\ZZZ$,
\begin{equation} \label{scalarproduct}
\nu \cdot \f := \sum_{i\in\ZZZ} \nu_i \f_ i ,
\end{equation}
we have
\[
f(t)=F(\om t) \,,\qquad\quad
F(\f) = \sum_{\nu\in\ZZZ^{\ZZZ}_{f}}
f_{\nu} \,  e^{\ii\nu\cdot\f} , \qquad\quad f_\nu\in X, \qquad\quad \sum_{\nu\in\ZZZ^{\ZZZ}_{f}}
\|f_{\nu} \|_X<\io .
\]
In our context, $X$ is a space of periodic functions in $x$ with Gevrey regularity, and $\om$ is restricted to the set
\begin{equation}\label{rettangolo}
{ { \gQ} := \Big\{ \om\in \RRR^\ZZZ \; :\; |\om_j-j^2| \le 1/2 \Big\} .} 
\end{equation}

We look for an almost-periodic solution $u(x,t)$ to \eqref{nls} with frequency vector $\om\in\gQ$ which is
Gevrey both in space and in time. To be more precise,
defining,  for any $\nu\in\ZZZ^{\ZZZ}_{f}$ and any $\mu\ge 0$,\footnote{Note that $|\nu|_0=\|\nu\|_1$.}
\begin{equation}\label{ric}
|\nu|_\mu := \sum_{i\in \ZZZ}  \langle i\rangle ^\mu |\nu_i| ,
\end{equation}
we look for a solution $u(x,t)$ which can be written as
\begin{equation}\label{fury}
u(x,t)= U(x,\omega t) \,,\qquad\qquad
U (x,\f) = \sum_{j\in\ZZZ} 
\sum_{\nu\in\ZZZ^{\ZZZ}_{f}}
u_{j,\nu} \, e^{\ii jx} e^{\ii\nu\cdot\f}\,,
\end{equation}
with $U$ defined on\footnote{{Using the notations and results of \cite{MP}, $\TTT^\ZZZ$
is a Banach manifold modelled on $\ell^\io$}; see Appendix \ref{gevgev} for further details. } 
$\TTT\times\TTT^\ZZZ$ and such that there exist $\al\in(0,1)$ and $s_1,s_2>0$ such that
\begin{equation}\label{espo}
\|U\|_{s_1,s_2,\al}:=\sup_{\substack{j\in\ZZZ \\ \nu\in\ZZZ^{\ZZZ}_{f}}}|u_{j,\nu}|e^{s_1|\nu|_\al}e^{s_2\jap{j}^\al} < \infty .
\end{equation}
Thus
the problem of existence of an almost-periodic solution  is reduced to the problem of
existence of an infinite-dimensional invariant torus parameterized by $\f\in\TTT^\ZZZ$.

%%%%%%%%%%%%%%%%%%%%%%%%%%%%%%%%%%%%%%%%%%%%%%%%%%%%%%%%%%%%%%%%%%%%%%%%%% 
\begin{rmk} \label{gevrey}
\emph{
A function $U\!:\TTT\times\TTT^\ZZZ\to\CCC$
such that \eqref{espo} holds is Gevrey of index $1/\al$ in $x$ and analytic in each argument $\f_j$, and hence, for $\om\in\gQ$,
the corresponding function $u$ defined as in \eqref{fury} is Gevrey with index $2/\al$ in $t$ (see Appendix \ref{gevgev}).
}
\end{rmk}
%%%%%%%%%%%%%%%%%%%%%%%%%%%%%%%%%%%%%%%%%%%%%%%%%%%%%%%%%%%%%%%%%%%%%%%%%% 

For $\e=0$, the NLS equation \eqref{nls} reduces to the linear equation $\ii u_t - u_{xx} + T_V u=0$ on the circle, whose solutions
\begin{enumerate}[topsep=0ex]
\itemsep0em
\item
are of the form $u(x,t)= U_{0}(x,\omega t;c)$, with 
\begin{equation} \label{sollinear}
U_{0}(x,\varphi;c):=\ \sum_{j\in \ZZZ} c_j e^{\ii (j x +\varphi_j)} , \qquad \omega =\{j^2+V_j\}_{j\in \ZZZ},\qquad c=\{c_j\}_{j\in\ZZZ},
\end{equation}
\item are almost-periodic for all choices of the sequences $\{c_j\}_{j\in\ZZZ}$ and $\{V_j\}_{j\in\ZZZ}$, 
\item are Gevrey in $x,t$ if $c \in \mathtt g({s,\alpha})$, with
\begin{equation}
\label{gsa}
\mathtt g({s,\alpha}):=\{c\,:\, \|c\|_{s,\alpha}:= \sup_{j\in\ZZZ} |c_j| e^{s\jap{j}^\alpha} <\io\} .
\vspace{-.5cm}
\end{equation}
\end{enumerate}

For future convenience, given any $c\in \mathtt g({s,\alpha})$ we set $|c|^2:=\{|c_j|^2\}_{j\in\ZZZ}$.

%%%%%%%%%%%%%%%%%%%%%%%%%%%%%%%%%%%%%%%%%%%%%%%%%%%%%%%%%%%%%%%%%%%%%%%%%% 
\begin{rmk} \label{gevrey0}
\emph{
Note that $\|U_0\|_{s_1,s_2,\alpha}<\io$ for all $s_1,s_2> 0$ such that $s_1+s_2\le s$.
}
\end{rmk}
%%%%%%%%%%%%%%%%%%%%%%%%%%%%%%%%%%%%%%%%%%%%%%%%%%%%%%%%%%%%%%%%%%%%%%%%%% 

Our main result ensures that, for all Gevrey $U_0$ and for many $V\in\ell^{N,\io}$, there exists a
solution to \eqref{nls} of the form \eqref{fury}, branching from $U_0$, provided that $\e\neq 0$ is small enough.

%%%%%%%%%%%%%%%%%%%%%%%%%%%%%%%%%%%%%%%%%%%%%%%%%%%%%%%%%%%%%%%%%%%%%%%%%% 
{\begin{theo}[\textbf{Main Theorem}]\label{main}
	For any $s>0$, any $\alpha\in(0,1)$ and any $N\in\ZZZ_+$, there exists  $\ol{\e}=\ol{\e}(s,\al,N)>0$
	such that for all $\e\in(-\ol{\e},\ol{\e})$ and for all $c\in \matU_1(\mathtt g({s,\al}))$ 
	there exists  a set ${\matC_N(\e,c,s,\al)}\subseteq\ol{\matU}_{1/4}(\ell^{N,\io})$, such that
\begin{enumerate}[topsep=0ex]
\itemsep0em
\item one has
	\begin{equation}\label{sipuo}
	\lim_{\e\to0}
	\meas(\matC_N(\e,c,s,\al))=1,
	\end{equation}
\item ${\matC_N(\e,c,s,\al)}$ depends on $c$ only through $|c|^2$,
\item for all $V\in {\matC_N(\e,c,s,\al)}$ the NLS equation \eqref{nls} 
	admits an almost-periodic solution $u(x,t)$ of the form \eqref{fury}, which is Gevrey in both $t$ and $x$ in the following sense:
	there is $s_0=s_0(\e)\in(0,s)$ and, for all $s_2\in(s_0,s)$, a constant $K= K(s_2,\al,N)>0$
	such that
		\begin{equation} \label{UU0}
		\|U-U_{0}\|_{s_1,s_2,\alpha} \le K |\e| ,
		\end{equation}
		with $s_1=s-s_2$ and $U_0$ as in \eqref{sollinear}. 
\end{enumerate}
\end{theo}}
%%%%%%%%%%%%%%%%%%%%%%%%%%%%%%%%%%%%%%%%%%%%%%%%%%%%%%%%%%%%%

%%%%%%%%%%%%%%%%%%%%%%%%%%%%%%%%%%%%%%%%%%%%%%%%%%%%%%%%%%%%%
\begin{rmk} \label{reg+meas}
\emph{
The bound \eqref{UU0} implies that the more regular the solution in time, the less regular is it in the space variable. Note
that, as in ref.~\cite{BMP2}, the values for the actions allowed fill a ball about the origin and hence are not confined to a zero-measure set.
 }
\end{rmk}
%%%%%%%%%%%%%%%%%%%%%%%%%%%%%%%%%%%%%%%%%%%%%%%%%%%%%%%%%%%%%%%%%%%%%%%%%% 

The strategy we follow in order to prove Theorem \ref{main}
is to reformulate the problem of persistence of solutions with infinitely many frequencies 
as a counterterm problem, in the spirit of Moser's approach in the finite-dimensional case \cite{M}.
There are many advantages in this approach, the most evident being that 
it allows to assume since the beginning appropriate non-resonance conditions on the frequency vector $\om$ of the solution.\footnote{In the 
infinite-dimensional case, this is deeply exploited already in ref.~\cite{BMP1}.} 

Precisely we proceed as follows. First, we fix any infinite-dimensional frequency vector $\om$ which satisfies some non-resonance condition:
the condition we require is weaker than the Diophantine conditions used in the literature (see Definition \ref{wbr}).
Next, we study a  modified equation (see \eqref{mtoro}), obtained by replacing the linear operator $-\del_{xx}+T_V$ with $T_\om + T_\h$,
where $\eta=\{\eta_j\}_{j\in \ZZZ}$ is a sequence of free parameters: we prove in Theorem \ref{moser} that it is
possible to fix the sequence $\eta$ at a suitable value, that we call the ``counterterm'', so that %\eqref{mtoro}
the modified equation admits  a solution $U(x,\om t)$ satisfying \eqref{UU0}.
%to a modified equation (see \eqref{mtoro}) where the operator $-\del_{xx}+T_V$ is replaced with $T_\om + T_\h$,
%where the parameters $\eta=\{\eta_j\}_{j\in \ZZZ}$ are called the ``counterterm''; see also Theorem \ref{moser}.
%
%where an additional term -- the ``counterterm''  $\eta=\{\eta_j\}_{j\in \ZZZ}$ -- is added to the linear part 
%(see \eqref{mtoro} and Theorem \ref{moser}).
For $U(x,\om t)$ to solve the original equation as well, the two linear operators must coincide: this induces a relation
between $\om$ and $V$ (see \eqref{perv}). Then we prove that for many choices of $V\in\ell^{N,\io}$ -- in the sense of \eqref{sipuo} --
%and hence we have to fix the potential $V$ appropriately (see \eqref{perv}),
%controlling at the same time that 
the corresponding $\om$ is among those satisfying the non-resonance condition;
this is achieved by studying an implicit function problem and showing that the latter can be solved for $V$
in a subset of $\ell^{N,\io}$ with large probability measure (see Proposition \ref{bohboh}).
This is a rather standard strategy, and it appears for instance also in ref.~\cite{BMP2},
where obtaining a regular counterterm would require a smoothing nonlinearity as in ref.~\cite{Po}:
indeed, the results of ref.~\cite{BMP2}, in our notation, show that $\eta$ belongs to $\ell^\io$ and hence,
since the regularity of $V$ is related to the regularity of the counterterm $\eta$,
they only imply that $V$ is in $\ell^\io$ as well.\footnote{In the direction of considering more regular potentials,
we mention ref.~\cite{GX}, where a T\"oplitz-Lipschitz property is showed to hold, so as to allow $V\in \ell^2$.}

In our paper the higher regularity of the potential $V$ is obtained by  
showing that, if the frequency vector $\om=\{\om_j\}_{j\in \ZZZ}$ has a suitable asymptotic expansion as $j\to \io,$ then
the counterterm admits a similar  asymptotic expansion (see Proposition \ref{decadeta} for a formal statement).
%and the heuristic discussion in Subsection \ref{rrmk} for further comments).
This is the crucial point that allows us to use the aforementioned implicit function theorem argument 
which leads to the regularity of $V$. 
%Note that the asymptotic expansion provided by Proposition \ref{decadeta}
%We stress that this is the crucial argument allowing
%us to deal with regular potential. 
An intuitive demonstration of why the asymptotic expansion is expected to hold
is given in Subsection \ref{rrmk}; however the argument described therein provides an expansion whose 
coefficients are formal power series, so one is left with the quite non-trivial problem of proving the convergence of such series.
By following a Renormalization Group approach to tackle such a task, we are able to exploit the symmetries of
the problem in order to obtain bounds which ensure the convergence.

We stress that, with respect to the existing literature, the main novelty of our result resides in the higher regularity of the Fourier multiplier operators,
which allows convolution potentials: it is the first result in which the regularity of the potential is not correlated to a
correspondingly smoothing nonlinearity. This is an important step in the direction of dealing with physically relevant models, where
the potential is very regular, and often the nonlinearity is even unbounded. 
Moreover it gives an insight on how to deal with the case of zero potential. 
We emphasize that in Theorem \ref{main} one has $\ol{\e}(s,\al,N)\to 0$
as $N\to+\io$; this is a typical feature in asymptotic expansions of pseudo-differential 
Sch\"odinger operators (see for instance refs.~\cite{Mar,Ric}), thus our result may well be optimal.

%Furthermore, as in refs.~\cite{GX,BMP2}, the values for the actions allowed fill a ball
%about the origin and hence are not confined to a zero-measure set. Finally
% the set of frequencies we consider is larger,
%since we introduce a Bryuno-like non-resonance condition which is weaker than the Diophantine condition considered in refs.~\cite{Bjfa,BMP2,BMP3}.

Another novelty is that we introduce an %infinite-dimensional
Bryuno-like non-resonance condition, which is weaker than the Diophantine
condition considered in refs.~\cite{Bjfa,BMP2,BMP3}, so our result holds for a larger set of frequencies; we refer to Definitions \ref{br} and
\ref{wbr}, and Appendix \ref{confronto} for further details.

\smallskip

\noindent
{\it Acknowledgements}. We  thank L.~Biasco, R.~Feola, E.~Haus and J.E.~Massetti for many useful comments.

%%%%%%%%%%%%%%%%%%%%%%%%%%%%%%%%%%%%%%%%%%%%%%%%%%%%%%%%%%%%%%%%%%%%%%%%%%
%%%%%%%%%%%%%%%%%%%%%%%%%%%%%%%%%%%%%%%%%%%%%%%%%%%%%%%%%%%%%%%%%%%%%%%%%% 
\section{Functional setting and scheme of the proof} %{An infinite-dimensional counterterm theorem} 
\label{setup} 
%%%%%%%%%%%%%%%%%%%%%%%%%%%%%%%%%%%%%%%%%%%%%%%%%%%%%%%%%%%%%%%%%%%%%%%%%% 
%%%%%%%%%%%%%%%%%%%%%%%%%%%%%%%%%%%%%%%%%%%%%%%%%%%%%%%%%%%%%%%%%%%%%%%%%% 

In this Section we introduce the technical setup and we state the key results from which we deduce the main Theorem \ref{main},
by following the strategy briefly outlined after Remark \ref{reg+meas}.
%Moreover we give an heuristic presentation about the validity of our technical results, and a brief presentation of the strategy of the proof.

If compared with refs.~\cite{Bjfa,BMP2}, our approach is more direct: we  look for almost-periodic solutions to \eqref{nls}
which are written explicitly in their Fourier series representation, namely which, if $\om\in\RRR^\ZZZ$
denotes the frequency vector, can be written as in \eqref{fury},
with the sum being totally convergent in a rather strong topology.
More precisely, the search of Gevrey almost-periodic solutions with frequency vector $\om\in\RRR^\ZZZ$ is reduced
to finding solutions $U(x,\f)$ to the \emph{torus equation}
\begin{equation}\label{toro}
\ii\om\cdot\del_\f U - \del_{xx} U + T_V U + \e |U|^4 U =0
\end{equation}
in the scale of Gevrey spaces
\begin{equation}\label{gev}
\mathtt{G}(s_1,s_2,\al):= \{ U:\TTT\times\TTT^\ZZZ\to\CCC\ :\ \|U\|_{s_1,s_2,\al}<\io\} ,
\end{equation}
that is in the scale of Banach spaces of functions $\TTT\times\TTT^\ZZZ\mapsto\CCC$ equipped with the Gevrey norm \eqref{espo},
with $\al\in(0,1)$, for suitable $s_1,s_2>0$.

%%%%%%%%%%%%%%%%%%%%%%%%%%%%%%%%%%%%%%%%%%%%%%%%%%%%%%%%%%%%%%%%%%%%%%%%%%
%%%%%%%%%%%%%%%%%%%%%%%%%%%%%%%%%%%%%%%%%%%%%%%%%%%%%%%%%%%%%%%%%%%%%%%%%% 
\subsection{First step: an infinite-dimensional counterterm theorem} 
\label{setup.p1} 
%%%%%%%%%%%%%%%%%%%%%%%%%%%%%%%%%%%%%%%%%%%%%%%%%%%%%%%%%%%%%%%%%%%%%%%%%% 
%%%%%%%%%%%%%%%%%%%%%%%%%%%%%%%%%%%%%%%%%%%%%%%%%%%%%%%%%%%%%%%%%%%%%%%%%% 

As already mentioned,
we  prove the existence of an almost-periodic solution by using a ``counterterm'' approach.
This means that we
 start by considering, instead of \eqref{toro}, the \emph{modified torus equation}
\begin{equation}\label{mtoro}
\ii\om\cdot\del_\f U + (T_\om + T_\h) U + \e |U|^4 U =0,
\end{equation}
where, not to overwhelm the notation and consistently with the definition of $V$, we are identifying the Fourier multipliers
$\om$ and $\eta$ with the corresponding sequences
$\{\om_j\}_{j\in\ZZZ}\in \gQ$ %\RRR^\ZZZ$ 
and $\{\h_j\}_{j\in\ZZZ}\in \ell^\io(\ZZZ,\RRR)$, respectively,
and look for a suitable sequence 
$\{\h_j\}_{j\in\ZZZ}\in \ell^\io(\ZZZ,\RRR)$,
that we call \emph{counterterm},
such that there exists a solution to \eqref{mtoro} in the space \eqref{gev}, for suitable $s_1$ and $s_2$.

%%%%%%%%%%%%%%%%%%%%%%%%%%%%%%%%%%%%%%%%%%%%%%%%%%%%%%%%%%%%%%%%%%%%%%%%%% 
\begin{rmk}\label{chestoaffa}
\emph{
The modified torus equation \eqref{mtoro} reduces to the torus equation \eqref{toro}
if the sequences $\om,\h,V$ satisfy
\begin{equation}\label{perv}
\om_j + \h_j = j^2 + V_j , \qquad j \in \ZZZ .
\end{equation}
}
\end{rmk}
%%%%%%%%%%%%%%%%%%%%%%%%%%%%%%%%%%%%%%%%%%%%%%%%%%%%%%%%%%%%%%%%%%%%%%%%%% 

In order to find a solution to \eqref{mtoro} we need to impose some non-resonance condition on the frequency vector $\om$.

%%%%%%%%%%%%%%%%%%%%%%%%%%%%%%%%%%%%%%%%%%%%%%%%%%%%%%%%%%%%%%%%%%%%%%%%%% 
\begin{defi}[\textbf{Bryuno condition}]
\label{br}
 Fix $\al\in(0,1)$. For $\om\in{\gQ}$, consider the non-increasing function $\be_\om \! : [1,+\io) \to\RRR$ given by
\begin{equation}\label{beta}
\beta_\om(x):= \inf_{\substack{\nu\in\ZZZ^{\ZZZ}_{f} \vspace{-.1cm}\\ 0 <|\nu|_{\al/2}\le x}}  |\om\cdot\nu| , 
\qquad\qquad
\om\cdot\nu := \sum_{i\in\ZZZ} \om_i \nu_i .
\end{equation}
Let $\gotR$ be the set of diverging increasing sequences $r=\{r_m\}_{m\ge0}$ on $[1,+\io)$. 
For any given $r\in\gotR$, introduce the Bryuno function
\begin{equation}\label{juno}
\BB_\om(r) := \sum_{m\ge1}\frac{1}{r_{m-1}}\log \left(\frac{1}{\be_\om(r_m)}\right)
\end{equation}
and set
\begin{equation}\label{set}
\gotB:=\bigcup_{r\in\gotR}\{\om\in{\gQ}\;:\; \BB_\om(r)<\io\}.
\end{equation}
We say that $\om\in{\gQ}$ satisfies the infinite-dimensional \emph{Bryuno condition}  if  $\om\in\gotB$.
Equivalently, we say that $\om\in\gotB$ is a (infinite-dimensional) \emph{Bryuno vector}.
\end{defi}
%%%%%%%%%%%%%%%%%%%%%%%%%%%%%%%%%%%%%%%%%%%%%%%%%%%%%%%%%%%%%%%%%%%%%%%%%% 

%%%%%%%%%%%%%%%%%%%%%%%%%%%%%%%%%%%%%%%%%%%%%%%%%%%%%%%%%%%%%%%%%%%%%%%%%%   
\begin{rmk} \label{r2m}
\emph{
As a special case of sequence in $\gotR$ one can take $r_m=2^m$ $\forall m\ge 0$;
this is a common choice in the finite-dimensional case. In fact, the Bryuno condition as stated 
in Definition \ref{br} is equivalent to requiring that $\BB_\om(\bar r) < +\io$ for $\bar r:=\{2^m\}_{m\ge 0}$
(see Appendix \ref{confronto}). The advantage in using Definition \ref{br} is that it allows more general sequences,
possibly adapted to the frequency vector $\om$; for instance, in the finite-dimensional case with $d=2$,
if $\om=(1,\la)$, a natural choice is to take $r_m=q_m$, where $\{q_m\}_{m\ge 0}$ is the sequence of the denominators
of the best approximants of $\la$ (see \eqref{originale} in Appendix \ref{confronto}).
}
\end{rmk}
%%%%%%%%%%%%%%%%%%%%%%%%%%%%%%%%%%%%%%%%%%%%%%%%%%%%%%%%%%%%%%%%%%%%%%%%%%   

%%%%%%%%%%%%%%%%%%%%%%%%%%%%%%%%%%%%%%%%%%%%%%%%%%%%%%%%%%%%%%%%%%%%%%%%%%   
\begin{rmk} \label{serveunnome}
\emph{
The definition \eqref{rettangolo} of $\gQ$ implies $\be_\om(x) \le 1/2$ for all $x\ge 1$.
}
\end{rmk}
%%%%%%%%%%%%%%%%%%%%%%%%%%%%%%%%%%%%%%%%%%%%%%%%%%%%%%%%%%%%%%%%%%%%%%%%%%   

%%%%%%%%%%%%%%%%%%%%%%%%%%%%%%%%%%%%%%%%%%%%%%%%%%%%%%%%%%%%%%%%%%%%%%%%%%   
\begin{rmk}
\emph{
Compared to the finite-dimensional analogue,
the set $\gotB$ depends on the value of $\alpha$ appearing in the definition of the function $\be_\om$ in \eqref{beta};
this is not surprising, since when dealing with infinite dimensional vector spaces all results depend strongly on the norms.
In particular, as $\alpha$ increases the value $\BB_\om(r)$ becomes smaller and hence,
in principle, the set $\gotB$ may become larger.
On the other hand, it turns out that the closer $\al$ to $0$, the closer to 0 is the value of $\ol{\e}$ in Theorem \ref{main}
(see also Remarks \ref{alpha0}, \ref{loss} and \ref{toner}).
In fact, the value of $\al$ is related to the regularity of the solution and the value $\al=0$ is not allowed (see \eqref{espo}).
}
\end{rmk}
%%%%%%%%%%%%%%%%%%%%%%%%%%%%%%%%%%%%%%%%%%%%%%%%%%%%%%%%%%%%%%%%%%%%%%%%%%   

The following definition is adapted from ref.~\cite{Bgafa}; see also refs.~\cite{Bjfa,BMP1,BMP2,MP,CMP}.

%%%%%%%%%%%%%%%%%%%%%%%%%%%%%%%%%%%%%%%%%%%%%%%%%%%%%%%%%%%%%%%%%%%%%%%%%% 
\begin{defi}[\textbf{Diophantine condition}]\label{dio}
Set, for $\g>0$ and $\tau>1/2$,
\begin{equation} \label{diofantinoBIS}
\gD(\g,\tau)  :=\Bigl\{ \omega\in {\gQ} : |\omega\cdot \nu|> \g
\prod_{i\in \ZZZ}\frac{1}{(1+ \jap{i}^{2}|\nu_i|^{2})^\tau} \; \forall \nu\in \ZZZ^{\ZZZ}_{f} \setminus\!\{0\} \Bigr\} ,
\end{equation}
and define
\begin{equation} \label{diofantinoBIS2}
\gD := \!\!\!\! \bigcup_{\substack{\g\in(0,1) \\ \tau >1/2}} \!\!\! \gD(\g,\tau) .
\end{equation}
We call $\gD$ the set of \emph{Diophantine vectors} and $\gD(\g,\tau)$
the set of Diophantine vectors with parameters $\g,\tau$.
\end{defi}
%%%%%%%%%%%%%%%%%%%%%%%%%%%%%%%%%%%%%%%%%%%%%%%%%%%%%%%%%%%%%%%%%%%%%%%%%% 

%%%%%%%%%%%%%%%%%%%%%%%%%%%%%%%%%%%%%%%%%%%%%%%%%%%%%%%%%%%%%%%%%%%%%%%%%% 
\begin{rmk}
\emph{
In Appendix \ref{confronto} we compare Definition \ref{br} with more common finite dimensional Bryuno conditions in the literature.
Moreover, since  $\gD \subseteq \gotB$ (see Lemma \ref{diovsbrj}), then a result proved by Bourgain in ref.~\cite{Bjfa}
(see also Lemma \ref{misurazza}) implies that $\gD$ -- and hence $\gotB$ -- has full probability measure.
}
\end{rmk}
%%%%%%%%%%%%%%%%%%%%%%%%%%%%%%%%%%%%%%%%%%%%%%%%%%%%%%%%%%%%%%%%%%%%%%%%%% 

%%%%%%%%%%%%%%%%%%%%%%%%%%%%%%%%%%%%%%%%%%%%%%%%%%%%%%%%%%%%%%%%%%%%%%%%%%%% 
%%\begin{rmk}
%%\emph{
%%As it turns out, we actually need the convergence of the sum for a suitable subsequence $\{r_{m_n}\}_{n\ge0}$ constructed recursively:
%% $r_{m_{n}}$ is the smallest element of the sequence $\{r_m\}_{m\ge1}$ such that
%% $\be()$
%% }
%%\end{rmk}
%%%%%%%%%%%%%%%%%%%%%%%%%%%%%%%%%%%%%%%%%%%%%%%%%%%%%%%%%%%%%%%%%%%%%%%%%%%% 

Let us come back to the modified torus equation \eqref{mtoro}. Set
\begin{equation}\label{diffop}
\gotD(\om):= (-\ii\om\cdot\del_\f - T_\om) ,
\end{equation}
so that \eqref{mtoro} reads
\begin{equation}\label{gene}
\gotD(\om) \, U = T_\h U + \e |U|^4U.
\end{equation}
For $\al\in(0,1)$, $\om\in\gotB$ and $s_1,s_2>0$, the kernel of the operator $\gotD(\om)$ in $\mathtt{G}(s_1,s_2,\al)$ is 
\begin{equation}\label{nucleo}
\begin{aligned}
\mathtt{K}(s_1+s_2,\al)&:=
\mbox{Ker}(\gotD(\om))\cap \mathtt{G}(s_1,s_2,\al) \phantom{\Big\{} \\
&= \Big\{ U_0(x,\f;c) = \sum_{j\in\ZZZ}c_je^{\ii(jx+\f_j)}\ : \ c\in\CCC^\ZZZ,\ \sup_{j\in\ZZZ}|c_j|e^{(s_1+s_2)\jap{j}^\al}<\io \Big\} ,
\end{aligned}
\end{equation}
so that, if $\gote_j \in \ZZZ^{\ZZZ}_{f}$ is such that $(\gote_j)_i=\de_{ji}$ and
\begin{equation}\label{range}
\mathtt{W}(s_1,s_2,\al):= \Big\{ U\in \mathtt{G}(s_1,s_2,\al) \ : \  U(x,\f)= \sum_{j\in\ZZZ} \sum_{{\nu\in\ZZZ^{\ZZZ}_{f}\setminus{ \{\gote_j\}}}} 
u_{j,\nu} \, e^{\ii (j x+\nu\cdot\f) } \Big\} 
\end{equation}
 denotes the
range of $\gotD(\om)$ in $\mathtt{G}(s_1,s_2,\al)$,
one has 
\begin{equation}\label{nulker}
\mathtt{G}(s_1,s_2,\al)=: \mathtt{K}(s_1+s_2,\al)\oplus \mathtt{W}(s_1,s_2,\al) .
\end{equation}

%%%%%%%%%%%%%%%%%%%%%%%%%%%%%%%%%%%%%%%%%%%%%%%%%%%%%%%%%%%%%%%%%%%%%%%%%% 
\begin{rmk}\label{iso}
\emph{
For all $s>0$ and all $\al>0$, the space $\mathtt{K}(s,\al)$ is isometrically identified with the sequence space
$\mathtt g(s,\al)$ in \eqref{gsa}.
}
\end{rmk}
%%%%%%%%%%%%%%%%%%%%%%%%%%%%%%%%%%%%%%%%%%%%%%%%%%%%%%%%%%%%%%%%%%%%%%%%%% 

%%%%%%%%%%%%%%%%%%%%%%%%%%%%%%%%%%%%%%%%%%%%%%%%%%%%%%%%%%%%%%%%%%%%%%%%%% 
\begin{rmk} \label{extrabryuno}
\emph{
By the gauge-covariance property of the NLS equation
(that is the covariance under the transformation $u(x,t) \mapsto e^{\ii \la} \, u(x,t) $ for any $\la\in\RRR$),
if $U$ solves \eqref{gene}, then in \eqref{range} there appear only Fourier labels satisfying the constraint
that the sum of the components equals one. Such a property turns out to be crucial
in order to prove the forthcoming Proposition \ref{bohboh}.
}
\end{rmk}
%%%%%%%%%%%%%%%%%%%%%%%%%%%%%%%%%%%%%%%%%%%%%%%%%%%%%%%%%%%%%%%%%%%%%%%%%% 

In view of Remark \ref{extrabryuno}, it is natural to introduce the following set of frequency vectors.

%%%%%%%%%%%%%%%%%%%%%%%%%%%%%%%%%%%%%%%%%%%%%%%%%%%%%%%%%%%%%%%%%%%%%%%%%% 
\begin{defi}[\textbf{Weak Bryuno condition}]
\label{wbr}
 Fix $\al\in(0,1)$. For $\om\in{\gQ}$, consider the non-increasing function $\be^{(0)}_\om \! : [1,+\io) \to\RRR$ given by
\begin{equation}\label{beta*}
\beta^{(0)}_\om(x):= \inf_{\substack{\nu\in\ZZZ^{\ZZZ}_{f} \vspace{-.1cm}\\ 0 <|\nu|_{\al/2}\le x \vspace{-.1cm}\\ {{\sum_{i\in\ZZZ}} \nu_i = 0}}}  |\om\cdot\nu| .
\end{equation}
For any given $r\in\gotR$, introduce the \emph{weak Bryuno function}

\begin{equation}\label{wjuno}
\BB_\om^{(0)}(r) := \sum_{m\ge1}\frac{1}{r_{m-1}}\log \left(\frac{1}{\be_\om^{(0)}(r_m)}\right)
\end{equation}
and set
\begin{equation}\label{wset}
\gotB^{(0)}:=\bigcup_{r\in\gotR}\{\om\in{\gQ}\;:\; \BB_\om^{(0)}(r)<\io\}.
\end{equation}
We say that $\om\in{\gQ}$ satisfies the infinite-dimensional \emph{weak Bryuno condition}  if  $\om\in\gotB^{(0)}$.
\end{defi}
%%%%%%%%%%%%%%%%%%%%%%%%%%%%%%%%%%%%%%%%%%%%%%%%%%%%%%%%%%%%%%%%%%%%%%%%%% 

%%%%%%%%%%%%%%%%%%%%%%%%%%%%%%%%%%%%%%%%%%%%%%%%%%%%%%%%%%%%%%%%%%%%%%%%%% 
\begin{defi}[\textbf{Separate analyticity}]
\label{ana}
{Let $(Z,\|\cdot\|_Z)$ and $(Y,\|\cdot\|_Y)$ be two complex Banach spaces. We say that a function $F:Z\to Y$ is 
\emph{separately analytic} if there is an analytic function ${F}_{\rm ext}\!:Z\times Z\to Y$ such that $F(z)={F}_{\rm ext}(z,\ol{z})$, where $\ol{z}$
is the complex conjugate of $z$.}
\end{defi}
%%%%%%%%%%%%%%%%%%%%%%%%%%%%%%%%%%%%%%%%%%%%%%%%%%%%%%%%%%%%%%%%%%%%%%%%%% 

%%%%%%%%%%%%%%%%%%%%%%%%%%%%%%%%%%%%%%%%%%%%%%%%%%%%%%%%%%%%%%%%%%%%%%%%%% 
\begin{rmk} \label{sepa-ana}
\emph{
Note that the definition of separate analyticity given above differs from the usual definition according to which a function $F\!:Z \times Z \to Y$
is separately analytic if it is analytic in each of the two variables while the other one is kept fixed \cite{JP}.
}
\end{rmk}
%%%%%%%%%%%%%%%%%%%%%%%%%%%%%%%%%%%%%%%%%%%%%%%%%%%%%%%%%%%%%%%%%%%%%%%%%% 

The following result, which, as already mentioned, is the analogue of Moser's counterterm theorem for the finite-dimensional case,
is the first step towards the proof of Theorem \ref{main}.

%%%%%%%%%%%%%%%%%%%%%%%%%%%%%%%%%%%%%%%%%%%%%%%%%%%%%%%%%%%%%%%%%%%%%%%%%% 
\begin{theo}[\textbf{Counterterm Theorem}]\label{moser}
Fix $s>0$, $\al\in(0,1)$ and $\om\in\gotB^{(0)}$. Then there exists $\e_0=\e_0(s,\al,\om)>0$ such that for all $\e\in(-\e_0,\e_0)$
and all $c\in \matU_1(\mathtt{g}(s,\al))$ there exist
a sequence $\h(c,\om,\e)$  and a function $U(x,\f;c,\om,\e)$
of the form \eqref{fury} which solve \eqref{mtoro} and satisfy the following properties.

\begin{enumerate}[topsep=0ex]
\itemsep0.0em
\item One has $ \h(c,\om,\e)=\{\h_j(c,\om,\e)\}_{j\in\ZZZ}\in\ell^\io(\ZZZ,\RRR)$ and there is a positive constant $\Phi_1=\Phi_1(s,\al,\om)$, independent of both $\e$ and $c$, such that
\begin{equation}\label{oh1}
\null\hspace{-.5cm}
\|\h(c,\om,\e)\|_{\io} \le  {\Phi_1}|\e|  . 
\end{equation}
\item Writing
\begin{equation} \label{Udecomposto}
U(x,\f;c,\om,\e) =U_0(x,\f;c)+ \Uperp(x,\f;c,\om,\e), 
\end{equation}
with $U_0(x,\f;c)$ as in \eqref{sollinear}, there exists a value $s'_0=s'_0(\e)\in(0,s)$ such that, for all $s_2\in(s'_0,s)$ and $s_1=s-s_2$,  one has
 $\Uperp(x,\f;c,\om,\e)\in \mathtt{W}(s_1,s_2,\al)$ and
there is a positive constant $\Phi_2=\Phi_2(s_2,\al,\om)$, independent of both $\e$ and $c$, such that
\begin{equation}\label{oh2}
\null\hspace{-.5cm}
\|\Uperp(c,\om,\e) \|_{s_1,s_2,\al} \le  {\Phi_2}|\e| . 
\end{equation}
\item The sequence $\eta(c,\om,\e)$ depends on $c$ only through $|c|^2$.
\item Both $\eta(c,\om,\e)$ and $U(x,\f;c,\om,\e)$ are separately analytic in $c,\ol{c}$ for $c\in \matU_1(\mathtt g(s,\alpha))$.
In particular $\eta(c,\om,\e)$ is analytic in $|c|^2$.
\end{enumerate}
\end{theo}
%%%%%%%%%%%%%%%%%%%%%%%%%%%%%%%%%%%%%%%%%%%%%%%%%%%%%%%%%%%%%%%%%%%%%%%%%% 

%
%%%%%%%%%%%%%%%%%%%%%%%%%%%%%%%%%%%%%%%%%%%%%%%%%%%%%%%%%%%%%%%%%%%%%%%%%%% 
%\begin{rmk} \label{perdita}
%\emph{
%Since $\gotD(\om)^{-1}$ is an unbounded operator, a loss of regularity in the angles is expected:
%this explains why one cannot take $s_1=s$, i.e.~$s_2=0$, in Theorem \ref{moser}.
%}
%\end{rmk}
%%%%%%%%%%%%%%%%%%%%%%%%%%%%%%%%%%%%%%%%%%%%%%%%%%%%%%%%%%%%%%%%%%%%%%%%%%% 

%%%%%%%%%%%%%%%%%%%%%%%%%%%%%%%%%%%%%%%%%%%%%%%%%%%%%%%%%%%%%%%%%%%%%%%%%% 
\begin{rmk} \label{attilio}
\emph{
The bound \eqref{oh2} ensures that $U(x,\f;c,\om,\e)$ is analytic in $\f\in\TTT^\ZZZ$,
with an explicit control on the analyticity strip (see Appendix \ref{gev}).
}
\end{rmk}
%%%%%%%%%%%%%%%%%%%%%%%%%%%%%%%%%%%%%%%%%%%%%%%%%%%%%%%%%%%%%%%%%%%%%%%%%% 

%%%%%%%%%%%%%%%%%%%%%%%%%%%%%%%%%%%%%%%%%%%%%%%%%%%%%%%%%%%%%%%%%%%%%%%%%% 
\begin{rmk} \label{cosmetics}
\emph{
The existence of a solution to the modified equation \eqref{mtoro} can be deduced also from~\cite{BMP2}. 
The improvement here, besides the stronger result that the frequency vector can vary in a larger set and the 
functions $\h$ and $U_\perp$ are proved to be separately analytic in $c,\ol{c}$ and analytic in $\f$, 
is that the proof we give allows us to
 obtain an explicit control on the Taylor expansions of $\h$ and $U_\perp$.
This turns  out to be fundamental in the proof of the forthcoming Proposition \ref{decadeta}
(see Subsection \ref{rrmk} for further comments), which is in turn is the key to proving Theorem \ref{main}.
}
\end{rmk}
%%%%%%%%%%%%%%%%%%%%%%%%%%%%%%%%%%%%%%%%%%%%%%%%%%%%%%%%%%%%%%%%%%%%%%%%%% 

%%%%%%%%%%%%%%%%%%%%%%%%%%%%%%%%%%%%%%%%%%%%%%%%%%%%%%%%%%%%%%%%%%%%%%%%%% 
\begin{rmk} \label{normally}
\emph{
In fact, the functions $\eta$ and $U$ in Theorem \ref{moser} are found to be
\emph{normally separately analytic} in $c,\ol{c}$ for $c\in \matU_1(\mathtt g(s,\alpha))$,
according to the standard terminology \cite{GKT0} adapted to our case (see Appendix \ref{gevgevgevgev}).
}\end{rmk}
%%%%%%%%%%%%%%%%%%%%%%%%%%%%%%%%%%%%%%%%%%%%%%%%%%%%%%%%%%%%%%%%%%%%%%%%%% 

%%%%%%%%%%%%%%%%%%%%%%%%%%%%%%%%%%%%%%%%%%%%%%%%%%%%%%%%%%%%%%%%%%%%%%%%%% 
\begin{rmk} \label{mettiamoloqui}
\emph{
Besides the gauge-covariance property mentioned in Remark \ref{extrabryuno},
another remarkable property of the NLS equation \eqref{nls} 
is the translation-covariance, that is the covariance under the transformation $u(x,t) \mapsto u(x+\la,t)$ for any $\la\in\RRR$.
While the gauge-covariance property, in our approach, is used only at the end for the measure estimates (see Section \ref{provafinale}),
on contrast, we use the translation-covariance property at length all along the paper (see also the comments in  Remark \ref{importante}).
In particular the traslation-covariance implies item 3 in Theorem \ref{moser} and that 
\[
U(x,\f,c,\om,\e)= U(0,0, c(x,\f),\om,\e)\,,\qquad c(x,\f):= \{c_j e^{\ii j x +\ii \f_j}\}_{j\in\ZZZ}\,.
\] 
}
\end{rmk}
%%%%%%%%%%%%%%%%%%%%%%%%%%%%%%%%%%%%%%%%%%%%%%%%%%%%%%%%%%%%%%%%%%%%%%%%%% 

%%%%%%%%%%%%%%%%%%%%%%%%%%%%%%%%%%%%%%%%%%%%%%%%%%%%%%%%%%%%%%%%%%%%%%%%%% 
\begin{rmk} \label{dovesimette}
\emph{
By inserting \eqref{Udecomposto} and \eqref{nucleo} into \eqref{nls}, one realizes immediately that the solution $U(x,\f;c,\om,\e)$
satisfies the symmetry 
$$
U(x,\f; \la c,\om,\la^{-4}\e)= \la \, U(x,\f;c,\om,\e), \qquad \h(\la c,\om,\la^{-4}\e) = \h(c,\om,\e). 
$$
}
\end{rmk}
%%%%%%%%%%%%%%%%%%%%%%%%%%%%%%%%%%%%%%%%%%%%%%%%%%%%%%%%%%%%%%%%%%%%%%%%%% 

%%%%%%%%%%%%%%%%%%%%%%%%%%%%%%%%%%%%%%%%%%%%%%%%%%%%%%%%%%%%%%%%%%%%%%%%%% 
\begin{rmk} \label{lecostantibis}
\emph{
By looking at the proof of Theorem \ref{moser}, one relizes that $\Phi_1=\tilde\Phi(s,\al,\om)$ and $\Phi_2=\tilde\Phi(s_2,\al,\om)$
for a suitable function $\tilde\Phi$.
}
\end{rmk}
%%%%%%%%%%%%%%%%%%%%%%%%%%%%%%%%%%%%%%%%%%%%%%%%%%%%%%%%%%%%%%%%%%%%%%%%%% 

Our scheme to prove Theorem \ref{moser} is the following. 
We start with  a Lyapunov-Schmidt decomposition, namely we split
%\[
%\mathtt{G}(s,s',\al)=: \mathtt{K}(s+s',\al)\oplus \mathtt{W}(s,s',\al) ,
%\]
%which implicitly defines the spaces $\mathtt{W}(s,s',\al)$.
%We write
%$U=U_0+\Uperp$, with $U_0\in \mathtt{K}=\mathtt{K}(s,\al)$ and 
%$\Uperp\in\mathtt{W}=\mathtt{W}(s_1,s_2,\al)$, with $s_1+s_2=s$; thus, we split
\eqref{gene} into the so-called \textit{range} and \textit{kernel} 
equations, i.e. %~we write
\begin{subequations}\label{sub}
\begin{align}
& \gotD(\om)\Uperp =   T_\h \Uperp+ \e\Pi_{\mathtt{W}}|U|^4U\,, \phantom{\big)}
\label{ran} \\
& T_\h U_0 + \e\Pi_{\mathtt{K}}|U|^4U =0 , \phantom{\Big)}
\label{bif}
\end{align}
\end{subequations}
with $\Pi_{\mathtt{W}}$ and $\Pi_{\mathtt{K}}$ denoting the projectors on $\mathtt{W}=\mathtt{W}(s_1,s_2,\al)$ and 
$\mathtt{K}=\mathtt{K}(s_1+s_2,\al)$, respectively.

%%%%%%%%%%%%%%%%%%%%%%%%%%%%%%%%%%%%%%%%%%%%%%%%%%%%%%%%%%%%%%%%%%%%%%%%%% 
\begin{rmk}\label{fastidio}
\emph{
In principle, a solution $\h$ to \eqref{bif} -- if any exists -- is in $\CCC^\ZZZ$. On the other hand, if $\ol{\h}$ denotes the complex conjugate
of $\h$ and $\ol{\h}\ne \h$, the counter\-term
cannot be interpreted as a correction to the frequency vector $\om$ so as to satisfy \eqref{perv}.
Thus, we have to check that the counterterm $\h$ in \eqref{bif} is real.
}
\end{rmk}
%%%%%%%%%%%%%%%%%%%%%%%%%%%%%%%%%%%%%%%%%%%%%%%%%%%%%%%%%%%%%%%%%%%%%%%%%%

Let us write
\begin{equation}\label{decay}
c_j =e^{-s\jap{j}^\al}\g_j,\qquad |\g_{j}| < 1\,,\qquad j\in\ZZZ.
\end{equation}

We look for  the solution $U(x,\f;c,\om,\e)$ and the counterterm $\h(c,\om,\e)$
as formal series in $\e$ (\emph{Lindstedt series}), namely we write
\begin{subequations} \label{formale+fout}
\begin{align}
U(x,\f;c,\om,\e) & = \sum_{j\in\ZZZ} c_j e^{\ii (j x +\f_j)}+  
\sum_{k\ge1}\e^k
\sum_{j\in\ZZZ} \sum_{\substack{\nu\in\ZZZ^{\ZZZ}_{f}\\ \nu\ne \gote_j}} 
u_{j,\nu}^{(k)}(c,\om) \, e^{\ii (j x+\nu\cdot\f)} ,
\label{fout} \\
\h_j (c,\om,\e) & = \sum_{k\ge1} \e^k \h_j^{(k)} (c,\om) .
\label{formale}
\end{align}
\end{subequations}

Inserting \eqref{formale+fout} into \eqref{sub} we compute
recursively the coefficients $u^{(k)}_{j,\nu}(c,\om)$ and $\h_j^{(k)}(c,\om)$ to all orders $k\ge 1$. 
Then, to conclude the proof of Theorem \ref{moser},
 we need to prove that the series are absolutely convergent, uniformly in $x,\f,c$, in the appropriate Gevrey class,
and that the counterterm is real (see Remark \ref{fastidio}) and depends on $c$ only through $|c|^2$.
Sections \ref{esplicito} to \ref{convergenza} are devoted to this task: 
Theorem \ref{moser} follows combining Proposition \ref{bellaprop} and Corollary \ref{contenti}.

%%%%%%%%%%%%%%%%%%%%%%%%%%%%%%%%%%%%%%%%%%%%%%%%%%%%%%%%%%%%%%%%%%%%%%%%%% 
%%%%%%%%%%%%%%%%%%%%%%%%%%%%%%%%%%%%%%%%%%%%%%%%%%%%%%%%%%%%%%%%%%%%%%%%%% 
\subsection{Second step: asymptotic expansion of the counterterm} 
\label{saline} 
%%%%%%%%%%%%%%%%%%%%%%%%%%%%%%%%%%%%%%%%%%%%%%%%%%%%%%%%%%%%%%%%%%%%%%%%%% 
%%%%%%%%%%%%%%%%%%%%%%%%%%%%%%%%%%%%%%%%%%%%%%%%%%%%%%%%%%%%%%%%%%%%%%%%%% 

In proving Theorem \ref{moser} we obtain a rather detailed expression for both the solution \eqref{fout} and the counterterm \eqref{formale}
as a ``sum over trees''. We use such an  expression of the counterterm $\h(c,\om,\e)$
to prove that, if we assume the frequencies to have a special dependence on $j$,
then the counterterms turn out to have the same dependence as well. 
More precisely the subset of frequencies in $\gQ$ we consider is as follows. 
Define
\begin{equation} \label{bienne}
\matW_N:= \begin{cases}
\ol{\matU}_{1/2}(\ell^{\io}),& N=0 \\
[-1/4,1/4] \times \ol{\matU}_{{1/2}}(\ell^{1,\io}), & N=1 \\
[-1/4,1/4]^{N-1}\times \ol{\matU}_{{1/2}}(\ell^{N,\io}),& N\ge2,
\end{cases}
\end{equation}
and write $\ze:=\xi$ for $\ze\in \matW_0$ and $\ze:=(\ka,\xi)$ for $\ze \in \matW_N$, with $N\ge 1$,
where $\xi\in \ol{\matU}_{{1/2}}(\ell^{N,\io})$, while $\ka=\ka_0\in[-1/4,1/4]$ if $N=1,2$ and
$\ka=(\ka_0,\ka_2,\ldots,\ka_{N-1})\in[-1/4,1/4]^{N-1}$ if $N \ge 3$.
For all $j\in\ZZZ$, define also $\om_j(\ze):=\om_j(\xi)=j^2+\x_j$ for $\ze\in \matW_0$,
while, for $N\ge1$ and $\ze\in\matW_N$, define\footnote{Here and henceforth we are using the convention that the sum over the empty set is $0$.}
\begin{equation}\label{espome}
        \om_j(\ze) := \om_j(\ka,\xi)=
	\begin{cases}
	\displaystyle{ j^2 +\kappa_0+ \sum_{q=2}^{N-1} \frac{\kappa_q}{j^q} + 
	\xi_j} \, ,& \qquad j \neq 0 , \\
	\ka_0+\xi_0,&\qquad j=0, 
	\end{cases}
\end{equation}
and set
\begin{equation} \label{KN}
\matK_N:=\{ \ze\in \matW_N:\,\omega(\ze)\in \gotB^{(0)}\}
\end{equation}
for $N\ge 0$. Then the following result is proved in Section \ref{labifsec}.

%%%%%%%%%%%%%%%%%%%%%%%%%%%%%%%%%%%%%%%%%%%%%%%%%%%%%%%%%%%%%%%%%%%%%%%%%% 
\begin{prop}\label{decadeta}
For any $s>0$, any $\al\in(0,1)$, any $N\ge 1$ and any $\ze\in \matK_N$ there exists ${\e_1=\e_1(s,\al,\ze,N)\in(0,\e_0)}$,
with $\e_0$ as in Theorem \ref{moser}, such that for all $\e\in(- \e_1,\e_1)$ 
and all $c\in\matU_1(\mathtt{g}(s,\al))$
%$c\in\mathtt{g}(s,\al)$ 
%such that $\|c\|_{s,\al}\le 1$
%
 there exist {real} coefficients $\gota_q(c,\ze,\e)$, with $q=0,2,\ldots,N-1$ if $N\ge 3$ and $q=0$ if $N=1,2$, and a sequence
$\gotr(c,\ze,\e)= \{\gotr_j(c,\ze,\e)\}_{j\in\ZZZ} \in \ell^{N,\io}(\ZZZ,\RRR)$, {depending on $c$ only through $|c|^2$} and such that 
\begin{equation} \label{nemmenounnome!}
\h_j(c,\om(\ze),\e) = 
\begin{cases}
\gota_0(c,\ze,\e) + \gotr_0(c,\ze,\e),&\qquad j=0 \\
\phantom{aaa}\\
\displaystyle{\gota_0(c,\ze,\e)+ \sum_{q=2}^{N-1} \frac{\gota_q(c,\ze,\e)}{j^q} + {\gotr_j(c,\ze,\e)},}& \qquad j \neq 0 .
\end{cases}
\end{equation}
Moreover there is a positive constant $C_0=C_0(s,\al,\ze)$, independent of both $\e$ and $c$, such that
\begin{equation}\label{tosse}
\sup_{q=0,2,\ldots,N-1}|\gota_q(c,\ze,\e)| \le C_0|\e|,\qquad \|\gotr(c,\ze,\e)\|_{N,\io}\le C_0| \e| ,
\end{equation}
where only $q=0$ has to be considered in the first bound if $N=1,2$. 
\end{prop}
%%%%%%%%%%%%%%%%%%%%%%%%%%%%%%%%%%%%%%%%%%%%%%%%%%%%%%%%%%%%%%%%%%%%%%%%%% 

Proposition \ref{decadeta} above is the key result allowing us to prove the existence of almost-periodic solutions to the NLS equation \eqref{nls}
with smooth convolution potentials. For a formal argument that shows why the result can be expected to hold, 
we refer to Subsection \ref{rrmk}, where we also highlight the difficulties that one encounters in making the formal argument rigorous,
and provide a motivation for choosing a Renormalization Group approach over KAM or Nash-Moser algorithms.

In particular, the Renormalization Group method provides a systematic tool to compute, at least in principle, the coefficients $\gota_q(c,\ze,\e)$,
$q=0,2,\ldots,N$.
For instance, if, relying on ref.~\cite{MP} and following the notations therein, for $F\in\mathtt G(s_1,s_2,\alpha)$ we write 
\begin{subequations} \label{mediaaaa}
\begin{align}
\langle F(\cdot,\f)\rangle_\TTT&=
\fint_{\TTT} F:= \frac{1}{2\pi}   \int_{0}^{2\pi} F(x,\varphi) \der x \,,
\label{media1} \\
 \langle F(x,\cdot)\rangle_{\TTT^\ZZZ}
 &= \fint_{\TTT^\ZZZ} F :=
\lim_{l\to\infty}\fint_{\TTT^{2l+1}}F(x,\varphi) \, \der\varphi_{-l}\ldots \der \varphi_l\,, 
\label{media2} \\
\langle  F \rangle_{\TTT\times \TTT^\ZZZ}&=
\langle  F(\cdot,\cdot) \rangle_{\TTT\times \TTT^\ZZZ}:= \fint_{\TTT\times \TTT^\ZZZ}F ,
\label{media}
\end{align}
\end{subequations}
we have the following result.

%%%%%%%%%%%%%%%%%%%%%%%%%%%%%%%%%%%%%%%%%%%%%%%%%%%%%%%%%%%%%%%%%%%%%%%%%% 
\begin{lemma}\label{carinello}
Under the assumptions of Proposition \ref{decadeta}, 
%let  $U$ be the solution of \eqref{mtoro} as in Theorem \ref{moser}. Then 
one has
$$
\gota_0(c,\ze,\e) = - 3\,\e\av{|U(\cdot,\cdot;c,\om(\ze),\e)|^4}_{\TTT\times\TTT^\ZZZ}.
$$
\end{lemma}
%%%%%%%%%%%%%%%%%%%%%%%%%%%%%%%%%%%%%%%%%%%%%%%%%%%%%%%%%%%%%%%%%%%%%%%%%% 

%However, before doing that,
%%in Subsection \ref{commenti}
%in the rest of this section 
%we describe how to conclude the proof of Theorem \ref{main}.

%%%%%%%%%%%%%%%%%%%%%%%%%%%%%%%%%%%%%%%%%%%%%%%%%%%%%%%%%%%%%%%%%%%%%%%%%% 
%%%%%%%%%%%%%%%%%%%%%%%%%%%%%%%%%%%%%%%%%%%%%%%%%%%%%%%%%%%%%%%%%%%%%%%%%% 
\subsection{Third step: the implicit function problem} 
\label{commenti} 
%%%%%%%%%%%%%%%%%%%%%%%%%%%%%%%%%%%%%%%%%%%%%%%%%%%%%%%%%%%%%%%%%%%%%%%%%% 
%%%%%%%%%%%%%%%%%%%%%%%%%%%%%%%%%%%%%%%%%%%%%%%%%%%%%%%%%%%%%%%%%%%%%%%%%% 

Theorem \ref{moser} is an abstract Moser-like counterterm theorem which ensures the existence of the solution $U$ and 
of the counterterm $\eta$ for all $c\in \matU_1(\mathtt{g}(s,\alpha))$ and all $\om\in\gotB^{(0)}$:
in particular the radius of convergence $\e_0$ is uniform in $c$, but it depends on the frequency vector $\om$.
The same happens for $\e_1$ as a function of $\zeta$ in Proposition \ref{decadeta}.
However, in order to apply Theorem \ref{moser} to the NLS equation \eqref{nls}, we need
to  solve the compatibility equation  \eqref{perv} and  express the frequency vector
$\om$ and the counterterm $\h$ in terms of the potential $V$; this is an implicit
function problem and, in order to solve it, we need Lipschitz regularity on some open set
where one has uniform bounds on $\e_1$.

To achieve this, first we restrict $\ze$ to a  subset $\matK_N(\g)$ of $\matK_N$
in which, the radius of convergence $\e_1$ is bounded uniformly in $\ze$.
Then, we prove that the coefficients $\gota_q,\gotr_j$ are Lipschitz-continuous w.r.t.~the $\ell^\io$-norm 
in $\matK_N(\g)$, so that, by relying on standard continuation arguments, they can be extended to functions
defined in $\matW_N$ and admitting there the same bounds. 
%Moreover we show that
%the probability measure of $\matK_N \setminus \matK_N(\g)$ is of order $\g$.
Then, once all this has been done,
an implicit function argument allows us to obtain $\ze=\ze(V)$.
We finally need to show that the set of  potentials $V$ such that
$\ze(V)\in \matK_N(\g)$ has the following two properties:
it is measurable and it has positive measure  in $\ol{\matU}_{1/4}(\ell^{N,\io})$.
Once the first property has been proved, the second one can be easily obtained from the Lipschitz regularity.
Thus, before anything else, we have to prove that $\matK_N(\g)$ is a measurable set:
since this requires $\ze(V)$ to be continuous w.r.t.~the product topology, in the argument above we need to ensure that 
the coefficients $\gota_q,\gotr_j$ and their extentions satisfy the same property.

More precisely, to implement the scheme described above, we proceed as follows.
First of all, in order to define the set $\matK_N(\g)$ we need to introduce some notation.
Fix $\g>0$ and $\tau>1/2$, and let $\gotR^*\subseteq\gotR$ be the set of sequences $r^*\in\gotR$ such that, setting
\begin{equation} \label{notte}
	\be^*(x) = \be^*(x,\g,\tau) := \g \!\!\!\! 
	\inf_{\substack{\nu\in\ZZZ^{\ZZZ}_{f} \vspace{-.1cm}\\ 0 <|\nu|_{\al/2}\le x \vspace{-.1cm}\\ {\sum_{i\in\ZZZ} \nu_i = 0}}} 
	%\inf_{\substack{\nu\in\ZZZ^{\ZZZ}_{f} \vspace{-.1cm}\\ 0 <|\nu|_{\al/2}\le x}}
	\prod_{i\in\ZZZ} \frac{1}{(1+\jap{i}^2|\nu_i|^2)^\tau} ,
	%\qquad \BB(r,\g) := \sum_{m\ge1}\frac{1}{r_{m-1}}\log\frac{1}{\be^*(r_m,\g)} , 
\end{equation}
one has
\begin{equation} \label{adessoce}
	\BB(r^*) = \BB(r^*,\g,\tau) := \sum_{m\ge1}\frac{1}{r^*_{m-1}}\log\frac{1}{\be^*(r_m^*,\g,\tau)}<\io.
\end{equation}
%

%%%%%%%%%%%%%%%%%%%%%%%%%%%%%%%%%%%%%%%%%%%%%%%%%%%%%%%%%%%%%%%%%%%%%%%%%% 
\begin{rmk}\label{mopro}
\emph{
In Appendix \ref{confronto} we show that $\{2^m\}_{m\ge0}\in\gotR^*$. In particular $\gotR^*$ is non-empty
and it does not depend on $\g,\tau$ (see Lemma \ref{diovsbrj} and Remark \ref{primanoncera}).
}
\end{rmk}
%%%%%%%%%%%%%%%%%%%%%%%%%%%%%%%%%%%%%%%%%%%%%%%%%%%%%%%%%%%%%%%%%%%%%%%%%% 

Then the appropriate subset $\matK_N(\g)$ is defined as follows.

%%%%%%%%%%%%%%%%%%%%%%%%%%%%%%%%%%%%%%%%%%%%%%%%%%%%%%%%%%%%%%%%%%%%%%%%%% 
\begin{defi}[\textbf{Good parameters}]
\label{unibr}
Let $\gotB(\gamma,\tau)\subset\gotB^{(0)}$ denote the set of weak Bryuno vectors $\om$ with the property that
there are $r\in\gotR$ and $r^*\in\gotR^*$ such that $\BB_\om^{(0)}(r) \le \BB(r^*,\g,\tau)$.
%$r_m\ge r_m^*$ and $\be_\om(r_m)\ge\be^*(r_m^*,\g)$ for all $m\ge0$.
For $N\ge 0$ the set of \emph{good parameters} is, 
\begin{equation} \label{bgamma}
\matK_N(\g):=\{ \ze \in \matW_N : \omega(\ze) \in \gotB(\gamma,N+1)\}.
\end{equation}
\end{defi}
%%%%%%%%%%%%%%%%%%%%%%%%%%%%%%%%%%%%%%%%%%%%%%%%%%%%%%%%%%%%%%%%%%%%%%%%%% 

%%%%%%%%%%%%%%%%%%%%%%%%%%%%%%%%%%%%%%%%%%%%%%%%%%%%%%%%%%%%%%%%%%%%%%%%%% 
\begin{rmk}\label{alpha0}
\emph{
In Remark \ref{primanoncera} we show that $\BB(r^*,\g,\tau)$ diverges as either $\al\to0$ or $\g\to 0$ or $\tau\to 1/2$.
Therefore the smaller $\g$ or $\al$, the larger is the bound on $\BB^{(0)}_\om(r)$ and hence the smaller is the estimated value of $\ol{\e}$
in Theorem \ref{main} (see Remark \ref{loss} for more details).
}
\end{rmk}
%%%%%%%%%%%%%%%%%%%%%%%%%%%%%%%%%%%%%%%%%%%%%%%%%%%%%%%%%%%%%%%%%%%%%%%%%% 

%%%%%%%%%%%%%%%%%%%%%%%%%%%%%%%%%%%%%%%%%%%%%%%%%%%%%%%%%%%%%%%%%%%%%%%%%% 
\begin{rmk}\label{misuyrabgt}
\emph{
In Appendix \ref{confronto} we also prove that, for all $\g>0$ and all $\tau>1/2$, one has
$\gD(\g,\tau)\subseteq\gotB(\g,\tau)$ and hence 
there is a positive constant $C$ such that $\meas(\gQ\setminus\gotB(\g,\tau))\le C \gamma$. 
}
\end{rmk}
%%%%%%%%%%%%%%%%%%%%%%%%%%%%%%%%%%%%%%%%%%%%%%%%%%%%%%%%%%%%%%%%%%%%%%%%%% 

Next, we extend the functions of interest outside the set $\matK_N(\g)$.
In order to do that, we have to prove first that such functions satisfy suitable Lipschitz bounds inside $\matK_N(\g)$. 
Given any subset $\matU\subseteq\matW_N$, with $N\ge 1$, and any map $f:\matU\to E$, 
with $(E,|\cdot|_E)$ some Banach space, we define the Lipschitz norm of $f$ as
\begin{equation}\label{lipnorm}
	|f|_{\matU,E}^{{\rm{Lip}}}:= \sup_{\ze\in \matU} | f(\ze) |_E + 
	\sup_{\substack{ \ze,\ze'\in \matU \\ \ze\neq\ze'}} \frac{|f(\ze) -f(\ze') |_E}{\|\ze-\ze'\|_{\io}} ,
\end{equation}
where $\|\ze\|_{\io} := \max\{ \|\kappa\|_{\io} , \|\x'\|_{\io} \}$ for $N\ge 1$ and $\|\ze\|_{\io} = \|\x'\|_{\io}$ for $N=0$,
and write
\begin{equation}\label{torta}
	|f|_{E}^{{\rm{Lip}}} := |f|_{\matW_N,E}^{{\rm{Lip}}}.
\end{equation}

The following result is proved in Subsection \ref{lip}.

%%%%%%%%%%%%%%%%%%%%%%%%%%%%%%%%%%%%%%%%%%%%%%%%%%%%%%%%%%%%%%%%%%%%%%%%%% 
\begin{prop}\label{fabrizio}
Fix $N\ge 0$ and $\g>0$.
For any $s>0$ and any $\al\in(0,1)$ %and all $s_1 \ge 0$ and $s_2>0$ such that $s_1+s_2=s$
there exist $\e_2=\e_2(s,\al,N,\g)>0$ and a positive constant $C_1=C_1(s,\al,N,\g)$ such that,
for all  $\e\in(-\e_2, \e_2)$ and all $c\in \matU_1(\mathtt{g}(s,\al))$, one has
\begin{subequations} \label{dolore!}
\begin{align}
|\h(c,\om(\cdot),\e)|^{{\rm{Lip}}}_{\matK_N(\g), \ell^\io} & \le C_1 |\e| , 
\label{dolore!a} \\
%|\Uperp(c,\om(\cdot),\e) |^{{\rm{Lip}}}_{\matK_N(\g),\mathtt{W}(s_1,s_2,\al)}  & \le C |\e| , 
%\label{dolore!b} \\
|\gota_q(c,\cdot,\e)|^{{\rm{Lip}}}_{\matK_N(\g),\RRR} & \le C_1 |\e|  , \qquad q=0,2,\ldots,N-1, 
\label{muffa1} \\
|\gotr(c,\cdot,\e)|^{{\rm{Lip}}}_{\matK_N(\g),\ell^{N,\io}} & \le C_1 |\e|  ,
\label{muffa2}
\end{align}
\end{subequations}
if $N\ge3$.
If $N=0,1,2$, bounds analogous to \eqref{dolore!} 
hold, %for $|\e|\le \e_2$, for some positive constants $\ol{\e}$ and $C$,
where \eqref{muffa1} and \eqref{muffa2} are missing if $N=0$ and \eqref{muffa1} holds only for $q=0$ if $N=1,2$.
Moreover, for all $\e\in(-\e_2,\e_2)$  there exists a value $s''_0=s''_0(\e)\in(0,s)$ such that, for all $s_2\in(s''_0,s)$ and $s_1=s-s_2$,
there is a positive constant $C_2=C_2(s_2,\al,N,\g)$ such that
\begin{equation}\label{dolore!b}
|\Uperp(c,\om(\cdot),\e) |^{{\rm{Lip}}}_{\matK_N(\g),\mathtt{W}(s_1,s_2,\al)}   \le C_2 |\e|.
\end{equation}
\end{prop}
%%%%%%%%%%%%%%%%%%%%%%%%%%%%%%%%%%%%%%%%%%%%%%%%%%%%%%%%%%%%%%%%%%%%%%%%%% 

In Subsection \ref{lipext}, we use McShane's theorem \cite{mcshane} and extend the functions in Proposition \ref{fabrizio} to functions
defined in the whole $\matW_N$, continuous w.r.t.~the product topology and satisfying therein the  bounds  \eqref{dolore!}.

Finally, we study the implicit function equation \eqref{perv}.
Let us consider first the case $N=0$. Recall that  $\matW_0= \ol{\matU}_{1/2}(\ell^\io)$ and $\om_j(\xi)=j^2+\xi_j$.  
Following the discussion above, for $|\e|< \e_2$, 
we construct an infinite-dimensional invariant torus and a counterterm for all $\xi\in \matK_0(\g)$. 
We extend the counterterm $\h_j(c,\om(\xi),\e)$ to a function $\h_j^{\rm Ext}(c,\xi,\e)$ defined for all $\x\in\matW_0$.
The compatibility equation \eqref{perv} reads
\begin{equation} \label{pervexplicit}
	\xi_j + \h_j^{\rm Ext}(c,\xi,\e) = V_j , \qquad j \in \ZZZ .
\end{equation}
We then use a fixed point argument to obtain $\x=\x(V)$, and hence, for all $V\in \ol{\matU}_{1/4}(\ell^\io)$, we obtain $\om=\om(\xi(V))\in \gQ$, as a
function continuous w.r.t.~the product topology and satisfying Lipschitz estimates similar to \eqref{dolore!};
to conclude, we use the latter properties to prove that the
set of potentials $V$ for which $\om(\x(V))\in\gotB(\gamma,1)$ has positive probability measure. This covers the result of \cite{BMP2}.

Let us now consider the case $N\ge1$.
For $|\e|< \e_2$ and $c \in \matU_1(\mathtt g(s,\alpha))$, consider the coefficients $\gota_q,\gotr_j$
defined in \eqref{nemmenounnome!} for $\ze=(\kappa,\xi)\in \matK_N(\g)$, and call $A_q, R_j$, respectively,
the extended functions obtained by applying McShane's theorem.
Then  \eqref{perv} becomes
\begin{equation} \label{system}
	\begin{cases}
		\kappa_q + A_q(c,\ze,\e)  = 0 , \qquad q=0, 2,\ldots, N-1 , \\
		%\xi_0 + R_0(c,\ze,\e) = V_ 0 , \\
		\x_j + R_j(c,\ze,\e) = V_j , \qquad j \in \ZZZ , %\neq 0 , 
	\end{cases}
\end{equation}
where only $q=0$ must be considered in the first line if $N=1,2$. The implicit function equation \eqref{system}
can be solved once more via a fixed point argument (see Subsection \ref{implicit}).
For fixed $c\in\matU_1( \mathtt{g}(s,\al))$ and $\e$ small enough, we write, for $j \neq 0$,
\begin{equation}
	\label{orrore}
	\om_j(\ze(V)) := j^2 + \ka_0(V) + \sum_{q=2}^{N-1} \frac{\ka_q(V)}{j^q} + {\xi_j(V)} , 
\end{equation}
{and verify that for ${\e}$ small enough $\om(\ze(V))\in \gQ$ for all $V \in \ol{\matU}_{1/4}(\ell^{N,\io})$}.
Then we define
\begin{equation} \label{rettangolino}
	\widehat \matC_N(\e,c,\g)= \left\{ {V}\in\ol{\matU}_{1/4}(\ell^{N,\io}) \, : \, \om(\ze(V)) \in {\gotB(\g,N+1)} \right\} .
\end{equation}
The following result is proved in Subsection \ref{misura}.

%%%%%%%%%%%%%%%%%%%%%%%%%%%%%%%%%%%%%%%%%%%%%%%%%%%%%%%%%%%%%%%%%%%%%%%%%%   
\begin{prop} \label{bohboh}
There is an absolute constant $C$ such that, if $\g< C^{-1}$, the set $\widehat \matC_N(\e,c,\g)$ defined in \eqref{rettangolino} satisfies
	\[
	\meas(\widehat \matC_N(\e,c,\g))\ge1-C\g.
	\]
\end{prop}
%%%%%%%%%%%%%%%%%%%%%%%%%%%%%%%%%%%%%%%%%%%%%%%%%%%%%%%%%%%%%%%%%%%%%%%%%%   

%%%%%%%%%%%%%%%%%%%%%%%%%%%%%%%%%%%%%%%%%%%%%%%%%%%%%%%%%%%%%%%%%%%%%%%%%% 
%%%%%%%%%%%%%%%%%%%%%%%%%%%%%%%%%%%%%%%%%%%%%%%%%%%%%%%%%%%%%%%%%%%%%%%%%   
\subsection{Conclusion of the proof}
\label{provamain} 
%%%%%%%%%%%%%%%%%%%%%%%%%%%%%%%%%%%%%%%%%%%%%%%%%%%%%%%%%%%%%%%%%%%%%%%%%   
%%%%%%%%%%%%%%%%%%%%%%%%%%%%%%%%%%%%%%%%%%%%%%%%%%%%%%%%%%%%%%%%%%%%%%%%%% 

The results stated in Subsections \ref{setup.p1}, \ref{saline} and \ref{commenti} are proved in Sections 
\ref{esplicito} to \ref{convergenza}, \ref{labifsec} and \ref{provafinale}, respectively.
Here we show that such results allow us to conclude immediately the proof of our main result Theorem \ref{main}.

Fix $s>0$, $\al\in(0,1)$ and $N\ge0$. 
Let $\e_2(s,\al,N,\g)$ be as in Proposition \ref{fabrizio} and note that by Remark \ref{alpha0} one has
\[
\lim_{\g\to0}\e_2(s,\al,N,\g)=0.
\]
For $\e$ sufficiently small, define
\[
\g_0(\e,s,\al,N):=2\inf\{ \g \, :\, |\e|<\e_2(s,\al,N,\g)\}.
\]
Precisely, since 
\[
\lim_{\e\to0}\g_0(\e,s,\al,N)=0 ,
\]
we need to fix $\ol{\e}=\ol{\e}(s,\al,N)$ in such a way that for all $\e\in(-\ol{\e},\ol{\e})$ one has $\g_0(\e,s,\al,N)<C^{-1}$, with $C$ the absolute constant
appearing in Proposition \ref{bohboh}. Now, for $\e\in(-\ol{\e},\ol{\e})$, set $\matC_N(\e,c,s,\al):=\widehat \matC_N(\e,c,\g_0(\e,s,\al,N))$,
so that \eqref{sipuo} is satisfied. 

Fix $\e\in(-\ol{\e},\ol{\e})$  and $c\in\matU_1(\mathtt{g}(s,\al))$, and note that by construction one has
 $|\e|<\e_2(s,\al,N,\g_0(\e,s,\al,N))$ and $\om(\ze(V))\in\gotB(\g_0(\e,s,\al,N),N+1)$
for all $V\in\matC_N(\e,c,s,\al)$. Thus, for any value $V\in\matC_N(\e,c,s,\al)$ we can apply 
Proposition \ref{fabrizio}  and conclude that \eqref{mtoro} with $\om=\om(\ze(V))$
admits a solution $U(x,\f;c,\om(\ze(V)),\e)$ which belongs to $\mathtt G(s_1,s_2,\al)$ for all $s_2\in(s_0(\e),s)$, and $s_1=s-s_2$.
This is the desired solution, and it is Gevrey in $x,t$ following Remark \ref{gevrey}.

Then,  Proposition \ref{decadeta} combined with \eqref{system} implies that, for all $V\in\matC_N(\e,c,s,\al)$, one has
$$
 \om_j(\ze(V))+\h_j(c,\om(\ze(V),\e) = j^2 + V_j  , \qquad j \in \ZZZ ,
 $$
and, therefore, the function $u(x,t):=U(x,\om(\ze(V))t;c,\om(\ze(V)),\e)$ solves \eqref{nls} as well.

%%%%%%%%%%%%%%%%%%%%%%%%%%%%%%%%%%%%%%%%%%%%%%%%%%%%%%%%%%%%%%%%%%%%%%%%%% 
%%%%%%%%%%%%%%%%%%%%%%%%%%%%%%%%%%%%%%%%%%%%%%%%%%%%%%%%%%%%%%%%%%%%%%%%%   
\section{A few comments on the strategy} 
\label{finrrmk} 
%%%%%%%%%%%%%%%%%%%%%%%%%%%%%%%%%%%%%%%%%%%%%%%%%%%%%%%%%%%%%%%%%%%%%%%%%   
%%%%%%%%%%%%%%%%%%%%%%%%%%%%%%%%%%%%%%%%%%%%%%%%%%%%%%%%%%%%%%%%%%%%%%%%%% 

As already pointed out in Remark \ref{cosmetics}, Theorem \ref{moser}, in principle, could be deduced from the result in ref.~\cite{BMP2}.
In particular, also by using a KAM approach one could show that the counterterm $\eta$ is well defined and real, and satisfies the identity \eqref{bif}. 
However, this is not enough if one wishes to solve the implicit function equation \eqref{perv} for $V\in \ell^{N,\io}$.
In this section, we explain where difficulties arise when trying to prove Proposition \ref{decadeta}
after proving a counterterm result like Theorem \ref{moser} through a KAM scheme, and why we resort to the Renormalization Group method
in order to overcome such difficulties.

%%%%%%%%%%%%%%%%%%%%%%%%%%%%%%%%%%%%%%%%%%%%%%%%%%%%%%%%%%%%%%%%%%%%%%%%%% 
%%%%%%%%%%%%%%%%%%%%%%%%%%%%%%%%%%%%%%%%%%%%%%%%%%%%%%%%%%%%%%%%%%%%%%%%%   
\subsection{A formal argument for the asymptotic expansion of the counterterm} 
\label{rrmk} 
%%%%%%%%%%%%%%%%%%%%%%%%%%%%%%%%%%%%%%%%%%%%%%%%%%%%%%%%%%%%%%%%%%%%%%%%%   
%%%%%%%%%%%%%%%%%%%%%%%%%%%%%%%%%%%%%%%%%%%%%%%%%%%%%%%%%%%%%%%%%%%%%%%%%% 

Let us fix $\ze\in\calmK_N$, with $\calmK_N$ as in \eqref{KN}, and consider a frequency vector $\omega=\omega(\ze)$ of the form \eqref{espome}. 
If we introduce the functions
\[
\begin{array}{lccc}
{\gotf_{j,\nu}=\gotf_{j,\nu}(x,\f):= e^{\ii j x+\ii\nu\cdot\f}} , & & j\in\ZZZ, &\nu\in\ZZZ^\ZZZ_f, \\
{\gotf_j=\gotf_j(x,\f):=\gotf_{j,\gote_j}(x,\f)= e^{\ii j x+\ii\f_j}} , & & j\in\ZZZ, &
\end{array}
\]
and recall \eqref{media}, then \eqref{bif} becomes
\[
\eta_j= - \frac{\e}{c_j} \langle  \ol{\gotf}_j |U|^4 U \rangle_{\TTT\times \TTT^\ZZZ},
\]
where $U=U_0+U_\perp$ is the function given by Theorem \ref{moser} and satisfying \eqref{ran}.

As a byproduct of the proof of Theorem \ref{moser},
not only the function $\eta$ is separately analytic in $c,\ol{c}$, but also, for all $j\in\ZZZ$, we have
(see Section \ref{convergenza} -- and \eqref{labase} in particular -- for details)
\[
\eta_j(c,\omega,\e)= \widehat\eta_j(c,\omega,\e)+ |c_j|^2 \, \widetilde\eta_j(c,\omega,\e) ,
\]
where $\widehat\eta_j$ does not depend on $c_j$, while $|c_j|^2\widetilde\eta_j(c,\omega,\e)$ 
is of order $\e e^{-\delta_0 \jap{j}^{\alpha_0}}$, for some positive constants $\delta_0,\alpha_0$. 
Thus, in studying the expansion of $\eta_j$ we can essentially ignore any contribution
which contains a factor $|c_j|^2$.
Similarly, again as a consequence of the proof of Theorem \ref{moser}, if one writes $U$ according to \eqref{fury}, i.e.~as
\[
U (x,\f) = \sum_{j'\in\ZZZ}\sum_{{\nu\in\ZZZ^\ZZZ_f }} u_{j',\nu} \, e^{\ii j'x} e^{\ii\nu\cdot\f}\,,
\]
then for any coefficient $u_{j',\nu}$ such that $\nu_j\ne 0$, we can factorize 
\[
u_{j',\nu} = |c_j|^{|\nu_j|} \widetilde  u_{j',\nu} ,
\]
with $\widetilde  u_{j',\nu}$ still separately analytic in $c,\ol{c}$,
and also such a contribution turns out to be sub-exponentially small in $j$.\footnote{Likely, the decay properties
of the corrections to the leading terms can be obtained also by carefully analyzing the proof of Theorem 3 in ref.~\cite{BMP2}.}
As a consequence, if we set
\begin{equation}
 \label{stiz}
 \mathtt Z_j:= \{(j',\nu)\in\ZZZ\times\ZZZ^\ZZZ_f: \nu_j=0\} \,,\qquad \Pi_{ \mathtt Z_j}U := \sum_{(j',\nu)\in \mathtt Z_j} 
 u_{j',\nu} \gotf_{j',\nu} ,
 \end{equation}
 we find that $U-  \Pi_{ \mathtt Z_j}U$ is of order $e^{-\delta_0\jap{j}^{\al_0}}$,
 possibly redefining the constants $\delta_0,\alpha_0$.
 This is a very informal statement, 
 but it can be made precise  in terms of an appropriate norm $\|\cdot\|_{s_1,s_2,\alpha}$ (see \eqref{espo}), with $s_1+s_2 \le s-\delta_0$.
Then we can prove that
\begin{equation} \label{mariposa}
\begin{aligned}
\!\!\!\!
\eta_j(c,\omega,\e) &=  -  \frac{\e}{c_j} 
\left\langle  \ol{\gotf}_j
\big(3|\Pi_{\mathtt Z_j}U|^4 U + 2   (\Pi_{\mathtt Z_j}U)^3 (\Pi_{\mathtt Z_j}\ol{U} )\ol{U}\big) 
\right\rangle_{\TTT\times \TTT^\ZZZ}
+ O(e^{-\delta_0\jap{j}^{\al_0}})  \\
&= - \frac{3\e}{c_j} \left\langle   \ol{\gotf}_j %\fint_{\TTT^\ZZZ\times \TTT}
|\Pi_{\mathtt Z_j}U|^4 U_0  \right\rangle_{\TTT\times \TTT^\ZZZ} \\ 
& \quad  - \frac{\e}{c_j} \left\langle \ol{\gotf}_j %\fint_{\TTT^\ZZZ\times \TTT}
\big(3|\Pi_{\mathtt Z_j}U|^4 U_\perp +
2   (\Pi_{\mathtt Z_j}U)^3 (\Pi_{\mathtt Z_j}\ol{U} )\ol{U}_\perp\big) 
\right\rangle_{\TTT\times \TTT^\ZZZ}
+ O(e^{-\delta_0\jap{j}^{\al_0}}) 
\\&= - 3\e \left\langle  %\fint_{\TTT^\ZZZ\times \TTT}
|U|^4 \right\rangle_{\TTT\times \TTT^\ZZZ} \\ 
& \quad  - \frac{\e}{c_j} \left\langle \ol{\gotf}_j %\fint_{\TTT^\ZZZ\times \TTT}
\big(3|\Pi_{\mathtt Z_j}U|^4 U_\perp +
2   (\Pi_{\mathtt Z_j}U)^3 (\Pi_{\mathtt Z_j}\ol{U} )\ol{U}_\perp\big) 
\right\rangle_{\TTT\times \TTT^\ZZZ}
+ O(e^{-\delta_0\jap{j}^{\al_0}}) ,
\end{aligned}
\end{equation}
and,
recalling that $U_\perp$ is in the range $\mathtt{W}$ (see \eqref{ran}), we obtain the identity
\begin{equation}\label{agito}
\begin{aligned} 
U_\perp & =\gotD^{-1}(\om)\,\Pi_{\mathtt{W}}(T_\eta U_\perp + \e\Pi_{\mathtt{W}}|U|^4 U ) \\
& = \gotD^{-1}(\om)\,\Pi_{\mathtt{W}}(T_\eta U_\perp +3\e|\Pi_{\mathtt Z_j}U|^4 U +
2 \e(\Pi_{\mathtt Z_j}U )^2  (\Pi_{\mathtt Z_j}\ol{U}) \ol{U} ) + c_j O(e^{-\delta_0\jap{j}^{\al_0}})  \\
& = \gotD^{-1}(\om)\,\Pi_{\mathtt{W}}(3\e|\Pi_{\mathtt Z_j}U|^4 U_0 +
2 \e(\Pi_{\mathtt Z_j}U )^2  (\Pi_{\mathtt Z_j}\ol{U}) \ol{U}_0 ) \\ 
&\qquad +
 \gotD^{-1}(\om)\,\Pi_{\mathtt{W}}(T_\eta U_\perp +3\e|\Pi_{\mathtt Z_j}U|^4 U_\perp +
2 \e(\Pi_{\mathtt Z_j}U )^2  (\Pi_{\mathtt Z_j}\ol{U}) \ol{U}_\perp ) 
+ c_j O(e^{-\delta_0\jap{j}^{\al_0}}).
\end{aligned}
\end{equation}

The intuition behind the expansion \eqref{nemmenounnome!} of the counterterm is  that the operator $\gotD^{-1}(\om)\Pi_{\mathtt{W}}$
should ``gain one derivative in space''. This implies that, if this intuition were correct, by inserting \eqref{agito}
into the last line in \eqref{mariposa} one  would obtain a factor $O(j^{-1})$ and hence, by dropping the projection $\Pi_{\mathtt Z_j}$
in the second line of \eqref{mariposa} (which may be done up to sub-exponentially small corrections),
one would find $\gota_0=-3\e\av{|U|^4}_{\TTT\times\TTT^\ZZZ}$.
Note that, with a different approach based on the so-called T\"oplitz-Lipschitz estimates, 
in ref.~\cite{GX} it was proved that $\h_j(c,\om,\e)$ is a constant up to corrections $O(j^{-1})$, but without providing the value of the constant
and under the stronger assumption that the coefficients $c_j$ decay super-exponentially.

Unfortunately, the intuition mentioned above is only partially correct: when inserting \eqref{agito} into the last line of \eqref{mariposa},
other contributions to $\gota_0$ appear, as the following argument shows.
Let us introduce some notation, by recalling \eqref{media1}, decomposing 
$\mathtt{G}(s_1,s_2,\al)  =\mathtt{G}_a + {\mathtt{G}}_n$, with
\[
\begin{aligned}
{\mathtt{G}}_a &={\mathtt{G}}_a(s_1,s_2,\al):=\{ F\in \mathtt{G}(s_1,s_2,\al) : \av{F(\cdot,\f)}_\TTT=0\},\\
\mathtt{G}_n &= \mathtt{G}_n(s_1,s_2,\al):=\{ F\in \mathtt{G}(s_1,s_2,\al) : F(x,\f)=\av{F(\cdot,\f)}_\TTT\}, 
\end{aligned}
\]
and defining the operators
\begin{subequations}
\begin{align}
\mathcal{G}_a[F]&:=
\ol{\gotf}_j\gotD^{-1}(\om)\, \Pi_{\mathtt{W}} [\,\gotf_j \,\Pi_{{\mathtt{G}_a}}F ]
=\sum_{\substack{j'\ne0 \\ \nu\neq \gote_{j+j'}-\gote_j}}
 \frac{F_{j',\nu}}{\omega\cdot\nu + \omega_j-\omega_{j +j'}}  \gotf_{j',\nu}  ,
 \label{sumabovep}
\\ %
\mathcal{G}_n[F]&:= \ol{\gotf}_j
 \gotD^{-1}(\om)\, \Pi_{\mathtt{W}} [\, \gotf_j\, \Pi_{{\mathtt{G}_n}}F ] =
\sum_{\nu\ne 0} \frac{F_{0,\nu}} {\omega\cdot\nu} \gotf_{0,\nu}
=  (-\ii\omega\cdot\partial_{\varphi})^{-1}\av{F - \av{F}_{\TTT^\ZZZ}}_\TTT ,
 \label{sumabove0}
\end{align}
\end{subequations}
where $ \Pi_{{\mathtt{G}_a}}$ and $ \Pi_{{\mathtt{G}_n}}$ are the projectors onto 
$\mathtt{G}_a$ and $\mathtt{G}_n$, respectively.
In \eqref{sumabovep}, 
for $\om=\om(\ze)$ of the form \eqref{espome}, at least in the regime $|j|\gg \max\{ |j'|,|\omega\cdot\nu|\}$, we may Taylor expand the factors
$(\omega\cdot\nu + \omega_j-\omega_{j +j'})^{-1}$ in powers of $j^{-1}$ up to order $N-1$, so as to obtain
\begin{equation}\label{estiqa}
\begin{aligned} 
\frac{1}{\omega\cdot\nu +\omega_j-\omega_{j + j'}} 
%\nonumber \\ & \qquad 
& = \frac{1}{\displaystyle{\omega\cdot\nu -2j j'- (j')^2 + 
\sum_{q=2}^{N-1} \kappa_q \Bigl( \frac{1}{j^q}-\frac{1}{(j + j')^q} \Bigr) +\xi_j- \xi_{j+ j'}} } \\
& = \frac{1}{2j j'} - \frac{\omega\cdot\nu - (j')^2}{2j^2(j')^2} + \ldots+O(j^{-N})\,. 
\end{aligned}
\end{equation}
This means that, provided we show that the terms with $|j|\lesssim \max\{ |j'|,|\omega\cdot\nu|\}$ give a sub-exponentially small contribution,
we have
\begin{equation} \label{pantrosa}
\mathcal{G}_a[F]=   - \frac1{2j}\partial_x^{-1} \Pi_{{\mathtt{G}_a}}F -
  \frac{1}{4j^2}\bigl(-\ii \omega\cdot\partial_\varphi + \partial_{xx} \bigr) \partial_x^{-2}  \Pi_{{\mathtt{G}_a}}F  + \ldots+O(j^{-N}) .
\end{equation}

On the other hand, not only the factors $j^{-1}$ are obtained at the price of applying to $F$ an unbounded operator,
but, in addition, in \eqref{sumabove0} no derivative is gained w.r.t.~$x$. Then, when inserting \eqref{agito}
into the last line of \eqref{mariposa}, we obtain, among others, a contribution of the form 
\begin{equation}\label{elvis}
\begin{aligned}
- 9\e^2 \left\langle \ol{\gotf}_j\,
\big( |\Pi_{\mathtt Z_j}U|^4 
\gotD^{-1}(\om)\,\Pi_{\mathtt{W}}\big)  \big[|\Pi_{\mathtt Z_j}U|^4\, \gotf_j \big] 
\right\rangle_{\TTT\times \TTT^\ZZZ},
\end{aligned}
\end{equation}
coming from the first summand in the third line of \eqref{agito}. Dropping the projector $\Pi_{\mathtt{Z}_j}$,
which produces again a sub-exponentially small correction, 
\eqref{elvis} reads
\begin{equation}\label{cash}
-9\e^2\av{|U|^4(\mathcal{G}_a + \mathcal{G}_n)[|U|^4] }_{\TTT\times\TTT^\ZZZ}=
-9\e^2\av{|U|^4\mathcal{G}_n[|U|^4] }_{\TTT\times\TTT^\ZZZ}-9\e^2\av{|U|^4\mathcal{G}_a[|U|^4] }_{\TTT\times\TTT^\ZZZ} ,
\end{equation}
with only the first summand contributing to $\gota_0$.
In contrast, the second summand is of order $O(j^{-1})$;
precisely, using \eqref{pantrosa} we arrive at
\[
\begin{aligned}
\av{|U|^4\mathcal{G}_a[|U|^4] }_{\TTT\times\TTT^\ZZZ} 
&= - \frac1{2j}\av{|U|^4\partial_x^{-1} \Pi_{{\mathtt{G}_a}}|U|^4 }_{\TTT\times\TTT^\ZZZ} \\
 &\qquad -  \frac{1}{4j^2} \av{ |U|^4
 \bigl(-\ii \omega\cdot\partial_\varphi + \partial_{xx} \bigr) \partial_x^{-2} \Pi_{{\mathtt{G}_a}}|U|^4 }_{\TTT\times\TTT^\ZZZ}
+ \ldots+O(j^{-N}) \\
&= - \frac1{2j}\av{ (\Pi_{{\mathtt{G}_a}} |U|^4)\partial_x^{-1} \Pi_{{\mathtt{G}_a}}|U|^4 }_{\TTT\times\TTT^\ZZZ} \\
&\qquad -  \frac{1}{4j^2} \av{ |U|^4
\bigl(-\ii \omega\cdot\partial_\varphi + \partial_{xx} \bigr)  \partial_x^{-2} \Pi_{{\mathtt{G}_a}}|U|^4 }_{\TTT\times\TTT^\ZZZ}
 + \ldots+O(j^{-N}) ,
\end{aligned}
\]
where the first term is found to vanish integrating by parts, while the others contribute to $\gota_2,\ldots,\gota_{N-1}$.
In fact, it turns out that also the contribution to $\gota_0$ in \eqref{cash} vanishes: indeed we have
\[
\begin{aligned}
\av{|U|^4\mathcal{G}_n[|U|^4] }_{\TTT\times\TTT^\ZZZ}&= \fint_{\TTT\times\TTT^\ZZZ} 
\Big(|U|^4  (-\ii\omega\cdot\partial_\varphi)^{-1} \fint_{\TTT}(|U|^4 - \fint_{\TTT^\ZZZ} |U|^4)\Big) \\
&= \fint_{\TTT\times\TTT^\ZZZ} \Big((|U|^4 - \fint_{\TTT^\ZZZ} |U|^4) (-\ii\omega\cdot\partial_\varphi)^{-1}
 \fint_{\TTT}(|U|^4 - \fint_{\TTT^\ZZZ} |U|^4)\Big) \\
&= \fint_{\TTT^\ZZZ}\big(\fint_{\TTT} (|U|^4 - \fint_{\TTT^\ZZZ} |U|^4)\big) (-\ii\omega\cdot\partial_\varphi)^{-1} 
\fint_{\TTT}\big(|U|^4 - \fint_{\TTT^\ZZZ} |U|^4\big) =0,
\end{aligned}
\]
where, once more, the last identity is obtained integrating by parts. 

The problem is that we can iterate \eqref{agito} an arbitrary number of times before inserting it into
\eqref{mariposa}, so obtaining, among others, terms of the form
\begin{equation} \label{stonza}
\begin{aligned}
& - \frac{(3\e)^{h+1}}{c_j} 
\left\langle  \ol{\gotf}_j %\fint_{\TTT^\ZZZ\times \TTT}
\big( |\Pi_{\mathtt Z_j}U|^4 
\gotD^{-1}(\om)\,\Pi_{\mathtt{W}}\big)^h\big[|\Pi_{\mathtt Z_j}U|^4U_0\big] 
\right\rangle_{\TTT\times \TTT^\ZZZ} \\
& \qquad = - {(3\e)^{h+1}} \left\langle  \ol{\gotf}_j %\fint_{\TTT^\ZZZ\times \TTT}
\big( |U|^4 
\gotD^{-1}(\om)\,\Pi_{\mathtt{W}}\big)^h  \big[|U|^4  \gotf_j\big] 
\right\rangle_{\TTT\times \TTT^\ZZZ} + O(e^{-\de_0\jap{j}^{\al_0}}) \\
&\qquad =  - {(3\e)^{h+1}} \langle  \big(|U|^4 (\mathcal{G}_n + \mathcal{G}_a) \big)^h [|U|^4]\rangle_{\TTT\times\TTT^\ZZZ}  + O(e^{-\de_0\jap{j}^{\al_0}}) ,
\end{aligned}
\end{equation}
for any $h\ge1$.
Thus, for any $h\ge1$, we obtain a contribution to $\gota_0$ of the form
\begin{equation}\label{print}
- {(3\e)^{h+1}} \langle  \big(|U|^4 \mathcal{G}_n  \big)^h [|U|^4]\rangle_{\TTT\times\TTT^\ZZZ},
\end{equation}
as well as contributions to $\gota_q$, for all $q\ge q_a$, where $q_a$ is the number of times the operator $\mathcal{G}_a$
appears when expanding $(|U|^4 (\mathcal{G}_n + \mathcal{G}_a))^h$  in  \eqref{stonza}. 
 
In conclusion, in agreement with the heuristic argument above, each coefficient $\gota_q$, with $q=0,2,\ldots,N-1$,
is found to admit a formal power series expansion in $\e$, so we are left with the problem of showing the convergence of
the series -- as well as that of proving that all terms we claimed to be sub-exponential remainders are indeed so. 
What makes the problem harder is that we need to prove the convergence of formal power series in the presence of both small denominators, 
due to the unbounded operators $\mathcal{G}_n$ in \eqref{sumabove0}, and large numerators, due to the unbounded operators 
 in \eqref{pantrosa}:
this is where the Renormalization Group method shows all its power. 

As it turns out, it is convenient to use the Renormalization Group since the beginning to prove Theorem \ref{moser}, 
since this provides us with $U$ in \eqref{fout} as a convergent power series in $\e$ with explicit expressions for the coefficients.
Then, using such a convergent series, for all $j\ne0$ we are able to write the counterterm in the form
\[
\h_j (c,\om(\ze),\e) = \sum_{q=0}^{N-1} \frac{\gota_q(c,\ze,\e,j)}{j^q} + {\gotA(c,\ze,\e,j)}, \qquad {\gotA(c,\ze,\e,j)} = O(j^{-N}) ,
\]
where the sequences $\{\gota_q(c,\ze,\e,j)\}_{j\in\ZZZ}$ and $\{\gotA(c,\ze,\e,j)\}_{j\in\ZZZ}$ are formal power series in $\e$
constructed as follows. First, we iterate the procedure above infinitely many times, so as to obtain 
contributions to $\h_j(c,\om(\ze),\e)$ where there appear the operators $\mathcal{G}_a$
and $\mathcal{G}_n$. Next, we take the Taylor expansion of $\mathcal{G}_a$ obtained by inserting \eqref{estiqa} into \eqref{sumabovep}.
The difference, with respect to the approach outlined above, is that the formal power series can be studied and proved to converge
with the same technique we used to prove the convergence of the series for $U$; see Subsection \ref{spiegotree} below for 
further comments.
Finally we prove that the sequences $\{\gota_q(c,\ze,\e,j)\}_{j\in\ZZZ}$ are
Cauchy sequences, % with $|\gota_q(j)-\gota_q(j+1)|\lesssim O(e^{-\de_0\jap{j}^{\al_0}})$, 
which allows us to define
\[
\gota_q(c,\ze,\e):=\lim_{j\to\io}\gota_q(c,\ze,\e,j),
\]
and we show that the limit is reached at sub-exponential rate.

Another advantage of the approach based on the Renormalization Group is that it provides
very explicit expressions also for the coefficients $\gota_q(c,\ze,\e,j)$. This allows us to prove, for instance,
that $\gota_1(c,\ze,\e)=0$ (see Lemma \ref{auno}) and that $\gota_0(c,\ze,\e)=-3\e\av{|U|^4}_{\TTT\times\TTT^\ZZZ}$ 
(see Subsection \ref{provocarinello})).

%A drawback of making the argument rigorous is that one loses track of the nice interpretation
%of the coefficient $\gota_q(c,\ze,\e)$ in terms of integrals of the solution.

%%%%%%%%%%%%%%%%%%%%%%%%%%%%%%%%%%%%%%%%%%%%%%%%%%%%%%%%%%%%%%%%%%%%%%%%%%
%%%%%%%%%%%%%%%%%%%%%%%%%%%%%%%%%%%%%%%%%%%%%%%%%%%%%%%%%%%%%%%%%%%%%%%%%%
\subsection{The Renormalization Group approach}
\label{spiegotree}
%%%%%%%%%%%%%%%%%%%%%%%%%%%%%%%%%%%%%%%%%%%%%%%%%%%%%%%%%%%%%%%%%%%%%%%%%%
%%%%%%%%%%%%%%%%%%%%%%%%%%%%%%%%%%%%%%%%%%%%%%%%%%%%%%%%%%%%%%%%%%%%%%%%%%

The use of a graphical representation for the coefficients of perturbative series in terms of  \emph{diagrams} is quite
common in Quantum Field Theory (QFT) and dates back to Feynman \cite{F1,F2}. In the context of classical
KAM theory it was originally introduced by Gallavotti  \cite{Galla}, inspired by the pioneering
work by Eliasson \cite{Elia1} (see ref.~\cite{GM0} for a comparison of the approaches), 
and thereafter it has been used in many other related papers (see ref.~\cite{Gen} for a review).\footnote{Eliasson's paper
inspired also other authors (see for instance refs.~\cite{FT,CF}), who, however, did not push forward 
with the application of the diagrammatic techniques to the KAM theory.
In particular the analogy with the QFT was first stressed in ref.~\cite{FT}, but only in refs.~\cite{Galla,GGM,GM0}
the Renormalization Group ideas were fully implemented to prove the convergence of the perturbative series.}
In the context of Hamiltonian PDEs the tree formalism has been fruitfully used in refs.~\cite{GMP,GP,MaPro} to prove the existence
of periodic solutions to the NLS equation and the beam equation.\footnote{A different Renormalization Group approach,
consisting in solving iteratively a sequence of fixed point equations, is used in refs.~\cite{BGK,BKS} to study the existence
of quasi-periodic solutions for both finite- and infinite-dimensional Hamiltonian systems.} 
More recently, Deng and Hani obtained a rigorous derivation of the wave kinetic equation 
from the cubic NLS equation by using trees to characterize the expressions appearing in the formal expansion of solutions
to NLS equation (see ref.~\cite{HD} and references therein).

With respect to most of the aforementioned literature, the present paper is closer in spirit to the finite-dimensional
classical KAM, where one looks for full-dimensional invariant tori. 
Indeed, when passing to the Fourier side w.r.t.~the space variable,
the NLS equation \eqref{nls} can be seen as a system of infinitely many weakly coupled harmonic oscillators written in complex variables.
In principle, it would be natural to pass to (infinitely many) action-angle variables before implementing any KAM scheme,
however it is well known that the action-angle  variables are singular near the origin -- this
is essentially the reason why Bourgain in ref.~\cite{Bjfa} requires the conditions $|u_j|^2\ge re^{-s\sqrt{|j|}}$.
To overcome such a difficulty, we work with the original coordinates instead of using action-angle variables,
by following the same strategy as in ref.~\cite{CGP} for the finite-dimensional case.
In fact, our proof of Theorem \ref{moser} relies on the same key ingredients as in ref.~\cite{CGP},
combined with two technical results on integer numbers, due to Bourgain, 
%which we report in Appendix \ref{dibou} for completeness, and
which are crucial to deal with the infinite-dimensional case.

The idea of the proof of Theorem \ref{moser} is as follows.
%%%The operator $ \gotD(\om)\circ\Pi_{\mathtt{W}}$ acting on the
%%% pair $(U,\ol{U})$  is block-diagonal in the Fourier basis, and it is of the form
%%%%
%%%\[
%%%\gotD(\om)\circ\Pi_{\mathtt{W}}={\diag_{\substack{(j,\nu)\in\ZZZ\times\ZZZ^\ZZZ_f \\ \nu\ne\gote_j}}}
%%%\begin{pmatrix}
%%%x_{j,\nu} & 0 \cr
%%%0 & x_{j,\nu}
%%%\end{pmatrix}
%%%\qquad\qquad
%%%x_{j,\nu}:=\om\cdot\nu - \om_j.
%%%\]
%%%%
For all $k\ge1$, all $j\in\ZZZ$ and all $\nu\in\ZZZ^\ZZZ_f$, 
the coefficients $u^{(k)}_{j,\nu}(c,\om)$ and $\h^{(k)}_{j}(c,\om)$ of the series \eqref{formale+fout} 
are computed recursively in $k$, and found to be expressed as infinite sums of terms,
each of which is a monomial in $c$ and $\ol{c}$ times the product of the so-called \emph{small divisors} arising from the
symbol $(\om\cdot\nu - \om_j)^{-1}$ of the operator $(\gotD(\om))^{-1} \Pi_{\mathtt W}$ in \eqref{ran}.
Thanks to Bourgain's results mentioned above, that we report in Appendix \ref{dibou} for completeness
(see Lemmata \ref{constance} and \ref{costanza}), one can show that
the coefficients $u^{(k)}_{j,\nu}(c,\om)$ and $\h^{(k)}_{j}(c,\om)$
are defined as absolutely convergent series. 
Nevertheless, due to the presence of terms where the small divisors accumulate,
a crude bound gives
\begin{equation}\label{jar}
|u^{(k)}_{j,\nu}(c,\om)|\le (C(k))^k e^{-s_1|\nu|_\al}e^{-s_2\jap{j}^\al} 
\end{equation}
and an analogous bound for $\h^k_j(c,\om)$,
with $C(k)$ diverging as $k$ goes to infinity (see Lemma \ref{welldefined}).
According to such a bound, if on the one hand the sums over $j\in\ZZZ$ and $\nu\in\ZZZ^\ZZZ_f$ turn out to be absolutely convergent,
on the other hand one can only deduce that both the solution and the counterterm are defined as formal power series in $\e$,
and the main difficulty is indeed to prove that the series in $\e$ is absolutely convergent as well.

%%%%%%%%%%%%%%%%%%%%%%%%%%%%%%%%%%%%%%%%%%%%%%%%%%%%%%%%%%%%%%%%%%%%%%%%%%
% FIGURA 19-2
%%%%%%%%%%%%%%%%%%%%%%%%%%%%%%%%%%%%%%%%%%%%%%%%%%%%%%%%%%%%%%%%%%%%%%%%%%
\begin{figure}[ht]
%\vspace{-.2cm}
\centering
%\null
%\hspace{-.6cm}
%\ins{070pt}{-128pt}{$\ell_1$}
%\ins{170pt}{-100pt}{$\ell_2$}
\ins{270pt}{-120pt}{$T$}
\ins{116pt}{-052pt}{$\ell_{2}$}
\ins{238pt}{-068pt}{$\ell_{3}$}
\ins{034pt}{-078pt}{$\ell_{1}$}
\ins{382pt}{-068pt}{$\ell_{4}$}
%\ins{346pt}{-058pt}{$\ell^{in}_\TT$}
\subfigure{\includegraphics*[width=6in]{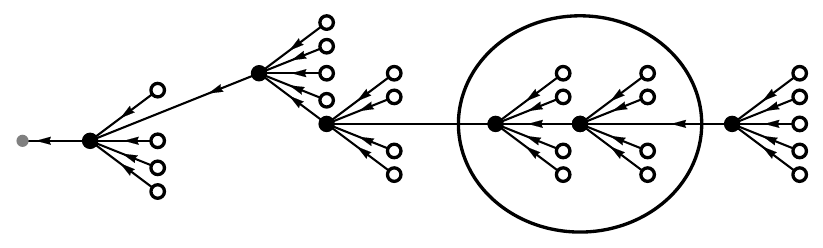}}
\vspace{-.2cm}
\caption{\small An example of tree containing a resonant cluster.}
\label{intro1}
\end{figure}
%%%%%%%%%%%%%%%%%%%%%%%%%%%%%%%%%%%%%%%%%%%%%%%%%%%%%%%%%%%%%%%%%%%%%%%%%%

The phenomenon outlined above is easier to visualize by means of a graphical representation of each term via a labelled rooted tree, 
which helps us to recognize the dangerous ones.
Shortly, to each line $\ell$ of a rooted tree
we assign labels $(j_\ell,\nu_\ell)\in\ZZZ\times\ZZZ^\ZZZ_f$ and we associate with $\ell$ a small divisor
$x_{\ell}:= \om\cdot\nu_\ell - \om_{j_\ell}$. It may happen that the tree contains subgraphs
with the following structure (see Figure \ref{intro1}): the subgraph is between lines whose small divisors are equal in absolute value, while
the absolute values of the small divisors associated with all the lines of the subgraph are much larger.
A subgraph of such a kind is called a resonant cluster -- or self-energy cluster, since it plays the same r\^ole as
the self-energy diagram in the context of QFT. For instance in Figure \ref{intro1} the resonant cluster $T$ is
between the two lines $\ell_3$ and $\ell_4$: the absolute values of the small divisors
associated with the lines $\ell_3$ and $\ell_4$ are equal to each other and are much smaller than the absolute values
of the small divisors associated with all the lines internal to $T$.

%%%%%%%%%%%%%%%%%%%%%%%%%%%%%%%%%%%%%%%%%%%%%%%%%%%%%%%%%%%%%%%%%%%%%%%%%%
% FIGURA 19-2
%%%%%%%%%%%%%%%%%%%%%%%%%%%%%%%%%%%%%%%%%%%%%%%%%%%%%%%%%%%%%%%%%%%%%%%%%%
\begin{figure}[H]
%\vspace{-.2cm}
\centering
%\null
%\hspace{-.6cm}
\ins{070pt}{-142pt}{$T_1$}
\ins{140pt}{-142pt}{$T_2$}
\ins{270pt}{-142pt}{$T$}
\ins{116pt}{-074pt}{$\ell_{2}$}
\ins{238pt}{-091pt}{$\ell_{3}$}
\ins{034pt}{-100pt}{$\ell_{1}$}
\ins{382pt}{-091pt}{$\ell_{4}$}
%\ins{346pt}{-058pt}{$\ell^{in}_\TT$}
\subfigure{\includegraphics*[width=6in]{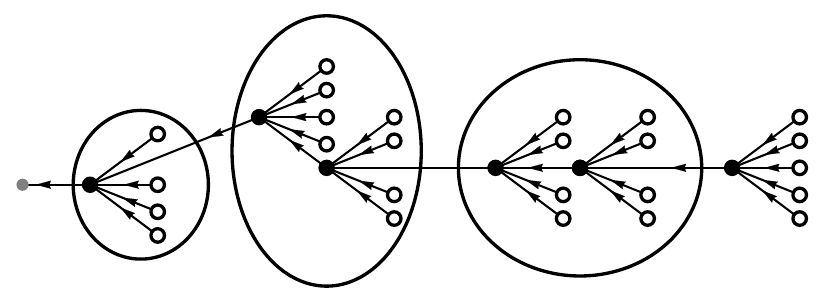}}
\vspace{-.2cm}
\caption{\small An example of tree containing a chain of resonant clusters: the small divisors of the the lines $\ell_1$, $\ell_2$, $\ell_3$ and $\ell_4$ are equal in absolute value,
and the absolute values of the small divisors associated with all the lines internal to the subgraphs $T$, $T_1$ and $T_2$
are much larger.}
\label{intro2}
\end{figure}
%%%%%%%%%%%%%%%%%%%%%%%%%%%%%%%%%%%%%%%%%%%%%%%%%%%%%%%%%%%%%%%%%%%%%%%%%%

In a rooted tree there is accumulation of small divisors whenever
it contains a path of lines oriented towards the root such that $p+1$ of them, for some $p\ge 1$,
have the same small divisors in absolute value
and the $p$ subgraphs between them are resonant clusters: we say in that case that the tree contains a chain
of $p$ resonant clusters. See Figure \ref{intro2} for an example of a chain of $p$ resonant clusters, with $p=3$.
Indeed, when a chain of resonant clusters arises in a tree,
the term of the sum corresponding to such a tree contains a factor $O(|x_{\ell}|^{-p})$ and, since $p$ can be $O(k)$,
the factor can not be bounded better than by $(C(k))^k$, with $C(k)$ as in \eqref{jar} (see the discussion in Subsection \ref{nuovo}, and
in particular the bound \eqref{stimeorrendec}).
The crucial point is that, if we carefully analyze the trees containing resonant clusters,
we find that, when grouping them together appropriately, there are suitable cancellations which allows us to improve drastically the bond \eqref{jar}.
In the case of classical KAM theorem in finite dimension, as it has been shown in \cite{BGK},
such cancellations are ultimately due to the 
invariance of the action functional w.r.t.~the translation $\om t\mapsto\om t + \be$ for any $\be\in\RRR^d$, where $d$ is the dimension of the 
invariant torus and $\om$ is the frequency vector of the solution supported on such a torus.
Here the cancellation mechanism %which we show essentially by hand,
is less transparent, both due to the infinite dimensional context
and to the fact that we are not working in action-angle variables.

%%%%%%%%%%%%%%%%%%%%%%%%%%%%%%%%%%%%%%%%%%%%%%%%%%%%%%%%%%%%%%%%%%%%%%%%%%
%%%%%%%%%%%%%%%%%%%%%%%%%%%%%%%%%%%%%%%%%%%%%%%%%%%%%%%%%%%%%%%%%%%%%%%%%%
\subsection{Content of the paper}
\label{paragrafone}
%%%%%%%%%%%%%%%%%%%%%%%%%%%%%%%%%%%%%%%%%%%%%%%%%%%%%%%%%%%%%%%%%%%%%%%%%%
%%%%%%%%%%%%%%%%%%%%%%%%%%%%%%%%%%%%%%%%%%%%%%%%%%%%%%%%%%%%%%%%%%%%%%%%%%

The rest of the paper is organized as follows. 

We introduce in Section \ref{esplicito} the tree formalism which provides the \emph{tree expansion}, that is 
 a graphical representation, of the solution to \eqref{mtoro}.
In Subsection \ref{sec:range} we describe
the tree expansion of the solution of the range equation \eqref{ran} alone, 
with the counterterm left as a free parameter: this is used later (see Subsection \ref{simmetrie}) %, in particular Proposition \ref{emmeeladerivatadigi}) 
to prove the symmetry property at the basis of the cancellation mechanism mentioned at the end of
Subsection \ref{spiegotree}.
In Subsection \ref{labif}
we describe the tree expansion of  the solution to the full equation \eqref{sub}, that is of both the function $U$ and the sequence $\h$ 
in Theorem \ref{moser}.

In Section \ref{sec} we define rigorously the resonant clusters introduced and roughly described in Subsection \ref{spiegotree}.

In Section
\ref{stimazze} we prove first of all that the coefficients  $u^{(k)}_{j,\nu}(c,\om)$ 
and $\h_j^{(k)}(c,\om)$ are well defined to all orders $k\ge 1$ and satisfy the bounds \eqref{jar}, so that the series \eqref{formale+fout} 
are well defined as formal power series in $\e$.
Next we prove that the resonant
clusters are indeed the only possible obstruction to the convergence, 
by showing that if we could ignore them, then
the series \eqref{formale+fout} would converge as power series in $\e$.

In Section \ref{convergenza} we conclude the proof of Theorem \ref{moser}.
Specifically in Subsection \ref{simmetrie} we prove a crucial symmetry property,
which is used thereafter, in Subsection \ref{LFR}, to
prove that the terms of the expansion corresponding to the trees which contain
resonant clusters, when summed together, admit bounds which ensure the convergence of the series in $\e$.

In Section \ref{labifsec}, by using a slightly modified version of the previous tree expansions, we prove Proposition \ref{decadeta} 
about the asymptotic expansion of the couterterm.

Finally in Section \ref{provafinale} we conclude the proof of Theorem \ref{main} 
by providing the Lipschitz estimates needed to use MacShane's theorem, so as to solve the implicit function problem \eqref{perv},
and finally deriving the measure estimates on the set $\matC_N(\e)$.

%%%%%%%%%%%%%%%%%%%%%%%%%%%%%%%%%%%%%%%%%%%%%%%%%%%%%%%%%%%%%%%%%%%%%%%%%%
%%%%%%%%%%%%%%%%%%%%%%%%%%%%%%%%%%%%%%%%%%%%%%%%%%%%%%%%%%%%%%%%%%%%%%%%%%
\section{Tree expansion}
\label{esplicito}
\zerarcounters
%%%%%%%%%%%%%%%%%%%%%%%%%%%%%%%%%%%%%%%%%%%%%%%%%%%%%%%%%%%%%%%%%%%%%%%%%%
%%%%%%%%%%%%%%%%%%%%%%%%%%%%%%%%%%%%%%%%%%%%%%%%%%%%%%%%%%%%%%%%%%%%%%%%%%

In this section we describe how to represent graphically in terms of trees both the coefficients $u^{(k)}_{j,\nu}(c,\om)$ of the formal
solution \eqref{fout} and the elements $\h_j^{(k)}(c,\om)$ of the counterterm \eqref{formale}.

%%%%%%%%%%%%%%%%%%%%%%%%%%%%%%%%%%%%%%%%%%%%%%%%%%%%%%%%%%%%%%%%%%%%%%%%%%
\subsection{Graphs and trees}
\label{alberi}
%%%%%%%%%%%%%%%%%%%%%%%%%%%%%%%%%%%%%%%%%%%%%%%%%%%%%%%%%%%%%%%%%%%%%%%%%%

%%%%%%%%%%%%%%%%%%%%%%%%%%%%%%%%%%%%%%%%%%%%%%%%%%%%%%%%%%%%%%%%%%%%%%%%%%
\begin{defi}[\textbf{Graph}]
A \emph{graph} $G$ is a pair $G=(V,L)$, where $V=V(G)$ is a
non-empty set whose elements are called \emph{vertices} 
and $L=L(G)$ a family of couples of unordered elements of $V(G)$, whose
elements are called \emph{lines} (or \emph{edges}).
\end{defi}
%%%%%%%%%%%%%%%%%%%%%%%%%%%%%%%%%%%%%%%%%%%%%%%%%%%%%%%%%%%%%%%%%%%%%%%%%%

Given two vertices $v,w\in V(G)$, the line $\ell=(v,w)$ is said to be \emph{incident} on both $v$ and $w$.
We say that a graph $G$ is \emph{simple} if two lines can not be incident on the same couples of vertices
A graph $G$ is \emph{finite} if $|V(G)|,|L(G)|<\io$;\footnote{Here and henceforth, $|A|$ denotes the cardinality of the set $A$.} 
in the following we tacitly assume all the graphs to be finite. We can represent a graph $G$ by drawing a set of points (the vertices)
and lines connecting them; couples appearing more than once are drawn as distinct lines.
A \emph{planar} graph is a graph $G$ which can be drawn on a plane
without lines crossing.

Given a graph $G$ and two different vertices $v,w\in V(G)$ 
we shall say that $v,w$ are \emph{connected} if
either $(v,w)\in L(G)$, or there exist
$v_{0}=v,v_{1},\ldots,v_{n-1},v_{n}=w\in V(G)$ such that
$\calP=\calP(v,w):=\{(v_{0},v_{1}),\ldots,(v_{n-1},v_{n})\}\subseteq L(G)$; we shall say
that $\calP$ is a \emph{path} connecting $v$ to $w$. A graph $G$ is \emph{connected} if for any $v,w\in V(G)$
there exists a path of lines connecting them.
A graph $G$ has a \emph{loop} if there exists $v\in V(G)$
such that either $(v,v)\in L(G)$ or there exists a path $\calP$ connecting $v$ to itself.

%%%%%%%%%%%%%%%%%%%%%%%%%%%%%%%%%%%%%%%%%%%%%%%%%%%%%%%%%%%%%%%%%%%%%%%%%%
\begin{defi}[\textbf{Subgraph}]
A \emph{subgraph} $S$ of a graph $G$ is a pair $S=(V,L)$
such that $V=V(S)\subseteq V(G)$, $L=L(S)\subseteq L(G)$ and $S$ is itself
a graph; if $S$ is a subgraph of $G$, we say $S$ is \emph{contained} in $G$.
\end{defi}
%%%%%%%%%%%%%%%%%%%%%%%%%%%%%%%%%%%%%%%%%%%%%%%%%%%%%%%%%%%%%%%%%%%%%%%%%%

%%%%%%%%%%%%%%%%%%%%%%%%%%%%%%%%%%%%%%%%%%%%%%%%%%%%%%%%%%%%%%%%%%%%%%%%%%
\begin{defi}[\textbf{Oriented graph}]
An \emph{oriented} graph is a graph $G$, where each couple in $L(G)$ is ordered.
If $G$ is an oriented graph, we say that the line $\ell=(v,w)
\in L(G)$ \emph{exits} the vertex $v$ and \emph{enters} $w$. We can represent
an oriented line $\ell=(v,w)$ by drawing the line with an arrow from $v$ to $w$ superimposed.
\end{defi}
%%%%%%%%%%%%%%%%%%%%%%%%%%%%%%%%%%%%%%%%%%%%%%%%%%%%%%%%%%%%%%%%%%%%%%%%%%

Given an oriented graph $G$, any subgraph $S$ inherits
the orientation of $G$. We say that a line $\ell=(v,w)\in L(G)$
\emph{exits} a subgraph $S$ of $G$ if $v\in V(S)$, while $w\in V(G)\setminus V(S)$; analogously we say that $\ell'=(v',w')$ \emph{enters} $S$
if $v'\in V(G)\setminus V(S)$, while $w'\in V(S)$. Note that if $\ell$ enters or exits $S$ then $\ell\notin L(S)$.

%%%%%%%%%%%%%%%%%%%%%%%%%%%%%%%%%%%%%%%%%%%%%%%%%%%%%%%%%%%%%%%%%%%%%%%%%%
\begin{defi}[\textbf{Rooted tree graph}]
A \emph{tree graph}  $\vartheta$ is a planar connected graph with no loops.
A \emph{rooted tree graph} is a tree graph with a (unique) special vertex $r=r_\vartheta$, called the \emph{root}, 
which has only one incident line, called the \emph{root line} and denoted by $\ell_{\vartheta}$.
\end{defi}
%%%%%%%%%%%%%%%%%%%%%%%%%%%%%%%%%%%%%%%%%%%%%%%%%%%%%%%%%%%%%%%%%%%%%%%%%%

The root induces a natural orientation on the lines of the graph {toward} the root:
for any $v\in V(\vartheta)\setminus\{r\}$ either $(v,r)\in L(\vartheta)$ or there are $v_1,\ldots,v_n\in V(\vartheta)$ such that
$\calP(v,r)=\{(v,v_1),(v_1,v_2),\ldots ,(v_n,r)\}\subseteq L(\vartheta)$. 
Thus a rooted tree is an oriented graph.

%%%%%%%%%%%%%%%%%%%%%%%%%%%%%%%%%%%%%%%%%%%%%%%%%%%%%%%%%%%%%%%%%%%%%%%%%%
\begin{rmk}\label{ellev}
\emph{
Given a rooted tree graph $\vartheta$, a line $\ell$ is uniquely determined by the
vertex $v$ which it exits, so we may write $\ell=\ell_{v}$.
}
\end{rmk}
%%%%%%%%%%%%%%%%%%%%%%%%%%%%%%%%%%%%%%%%%%%%%%%%%%%%%%%%%%%%%%%%%%%%%%%%%%

The orientation on a rooted tree graph $\vartheta$ provides a partial ordering relation on the vertices:
given $v,w\in V(\vartheta)$, we write $v\preceq w$ if $\calP(v,w)\subseteq L(\vartheta)$.
If $\ell=\ell_{v}$ and $\ell'=\ell_{v'}$ for some vertices $v,v'\in V(\vartheta)$,
we write $v\prec\ell'$, $\ell\prec v'$  and $\ell\prec \ell'$ if $v\prec v'$.
Given two distinct lines $\ell,\ell'\in L(\vartheta)$ such that $\ell\prec\ell'$, let
$\calP(\ell,\ell')$ denote the unique path connecting the vertex $\ell$ enters to the vertex $\ell'$ exits.
If $\ell'$ exits the same vertex $\ell$ enters then $\calP(\ell,\ell')=\emptyset$.
Given a path $\calP$ we say that a node $v$
is in $\calP$ if either $\ell_v\in\calP$ or there is a line $\ell\in\calP$ entering $v$.

Given a rooted tree graph $\vartheta$, for each vertex $v\in V(\vartheta)$ let $s_v$ be the number of lines entering $v$
(\emph{branching label}); we call \emph{nodes} the vertices $v\in V(\vartheta)$ with $s_v\ge1$ except the root $r$ and set
$$
N(\vartheta):=\{v\in V(\vartheta)\ :\ s_v\ge1\}\setminus\{r\}.
$$ 
We call \emph{leaves} the vertices $v$ with $s_v=0$ and set
$$
\Lambda(\vartheta):=\{v\in V(\vartheta)\ :\ s_v=0\}.
$$
Finally, If $T$ is a subgraph of $\vartheta$, we set $N(T):=N(\vartheta)\cap V(T)$ and $\Lambda(T):=\Lambda(\vartheta)\cap V(T)$.

%%%%%%%%%%%%%%%%%%%%%%%%%%%%%%%%%%%%%%%%%%%%%%%%%%%%%%%%%%%%%%%%%%%%%%%%%%
\begin{rmk} \label{siamopignoli}
\emph{%remark
%In fact
Most of the time 
we use the letter $v$ only for nodes, and the letter $\la$ for leaves.
}%remark
\end{rmk}
%%%%%%%%%%%%%%%%%%%%%%%%%%%%%%%%%%%%%%%%%%%%%%%%%%%%%%%%%%%%%%%%%%%%%%%%%%

%%%%%%%%%%%%%%%%%%%%%%%%%%%%%%%%%%%%%%%%%%%%%%%%%%%%%%%%%%%%%%%%%%%%%%%%%%
\begin{rmk}\label{drawing}
\emph{
It turns out to be convenient for us to draw the trees with the root on the far left and the arrows all pointing from right to left.
In particular, for any node $v$ with branching label $s_v\ge 2$, the lines $\ell_1,\ldots,\ell_{s_v}$ entering $v$ are all drawn
to the right of $v$ on the top of each other, and a path $\calP(v,w)$ connects the vertex $v$ to the right with the vertex $w$ to the left
(see Figure \ref{superalcolici} below for an explicit example).
}
\end{rmk}
%%%%%%%%%%%%%%%%%%%%%%%%%%%%%%%%%%%%%%%%%%%%%%%%%%%%%%%%%%%%%%%%%%%%%%%%%%

%%%%%%%%%%%%%%%%%%%%%%%%%%%%%%%%%%%%%%%%%%%%%%%%%%%%%%%%%%%%%%%%%%%%%%%%%%
\begin{defi}[\textbf{Labelled rooted tree graph}] \label{indicetti}
A \emph{labelled rooted tree graph} %(or simply \emph{tree} in the following)
is a rooted tree graph together with a label function
defined on $V(\vartheta)$ and $L(\vartheta)$. 
The labels are assigned according to the following rules:
\begin{enumerate}[topsep=0ex]
\itemsep0.0em
\item with each leaf $\la\in \Lambda(\vartheta)$ we associate a \emph{component} label $j_\la\in\ZZZ$ 
and a \emph{sign} label $\s_\la\in\{\pm\}$, and, for all $j\in\ZZZ$ and $\s\in\{\pm\}$, we call $\Lambda_{j,\s}(\vartheta)$
the set of leaves $\la\in\Lambda(\vartheta)$ such that $j_\la=j$ and $\s_\la=\s$,
{while, for all $j\in\ZZZ$, we set $\Lambda_{j}(\vartheta)=\Lambda_{j,+}(\vartheta)\cup\Lambda_{j,-}(\vartheta)$};
\item with each node $v\in N(\vartheta)$ we associate an \textit{order} label $k_v\ge0$,
a \textit{component} label $j_v\in\ZZZ$ and a sign label $\s_v\in\{\pm\}$;
%
%%We represent with black bullets the nodes 
%%with $s_v=5$ and with white squares the nodes with $s_v=1$; finally we represent with white bullets
%%the leaves and with a grey bullet the root (see Figure \ref{superalcolici}).
%
\item with the root $r$ we associate a \emph{component} label $j_r=j_{\ell_\vartheta}$;
\item with each line $\ell\in L(\vartheta)$ we associate a \emph{component} label $j_{\ell}\in\ZZZ$, a \emph{sign} label $\s_\ell\in\{\pm\}$,
a \textit{scale} label $n_\ell \in \ZZZ$ and a \emph{momentum} label $\nu_\ell\in\ZZZ^{\ZZZ}_{f}$.
\end{enumerate}
\end{defi}
%%%%%%%%%%%%%%%%%%%%%%%%%%%%%%%%%%%%%%%%%%%%%%%%%%%%%%%%%%%%%%%%%%%%%%%%%%

%%%%%%%%%%%%%%%%%%%%%%%%%%%%%%%%%%%%%%%%%%%%%%%%%%%%%%%%%%%%%%%%%%%%%%%%%%
\begin{rmk}\label{tree}
\emph{
In the rest of the paper we will consider only labelled rooted tree graphs. We will refer to them simply as trees, but one has to keep in mind
that henceforth ``trees'' are meant as ``labelled rooted tree graphs".
}
\end{rmk}
%%%%%%%%%%%%%%%%%%%%%%%%%%%%%%%%%%%%%%%%%%%%%%%%%%%%%%%%%%%%%%%%%%%%%%%%%%

%%%%%%%%%%%%%%%%%%%%%%%%%%%%%%%%%%%%%%%%%%%%%%%%%%%%%%%%%%%%%%%%%%%%%%%%%%
\begin{defi}[\textbf{Order of a tree}]
We define the \emph{order} of a tree $\vartheta$ as 
\begin{equation}\label{ordine}
k(\vartheta) := \sum_{v\in N(\vartheta)} k_v .
\end{equation}
\end{defi}
%%%%%%%%%%%%%%%%%%%%%%%%%%%%%%%%%%%%%%%%%%%%%%%%%%%%%%%%%%%%%%%%%%%%%%%%%%

%%%%%%%%%%%%%%%%%%%%%%%%%%%%%%%%%%%%%%%%%%%%%%%%%%%%%%%%%%%%%%%%%%%%%%%%%%
\begin{rmk}\label{stoper}
\emph{
The  number of vertices in a tree $\vartheta$ is bounded as
\[
|V(\vartheta)|\le 1 + \sum_{v\in N(\vartheta)}s_v,
\]
as it can be easily seen by induction. In particular, if there exists $p\ge1$ such that $s_v\le p$ for all $v\in N(\vartheta)$
one has $|V(\vartheta)|\le 1+p |N(\vartheta)|$ and hence $|\Lambda(\vartheta)|\le 1+(p-1) |N(\vartheta)|$.
Moreover, the number of trees with $V$ vertices is bounded by $2^{2V}$ \cite[Lemma A1]{GM}, so that, 
if $s_v\le p$ for all $v\in N(\vartheta)$, the number of trees with $N$ nodes is bounded by $4^{(p+1)N}$.
}
\end{rmk}
%%%%%%%%%%%%%%%%%%%%%%%%%%%%%%%%%%%%%%%%%%%%%%%%%%%%%%%%%%%%%%%%%%%%%%%%%%

%%%%%%%%%%%%%%%%%%%%%%%%%%%%%%%%%%%%%%%%%%%%%%%%%%%%%%%%%%%%%%%%%%%%%%%%%%
\begin{defi}[\textbf{Subtree}] \label{subtree}
Given a tree $\vartheta$ and a node $v \in N(\vartheta)$, let $\vartheta'$ be the subgraph  
of $\vartheta$ such that
%\vspace{-.2cm}
\begin{enumerate}[topsep=0ex]
\itemsep0em
\item $V(\vartheta')$ contains all the vertices $v' \preceq v$ and the vertex ${v}_0$ which $\ell_v$ enters,
\item $L(\vartheta')$ contains  all the lines $\ell' \preceq \ell_v$.
\end{enumerate}
%\vspace{-.2cm}
Then $\vartheta'$ is a tree and  ${v}_0$ is its root. We say that $\vartheta'$ is a \emph{subtree} of the tree $\vartheta$,
and write $\vartheta'\subseteq\vartheta$, with $\vartheta'\subset\vartheta$ except when $v_0$ is the root of $\vartheta$.
\end{defi}
 %%%%%%%%%%%%%%%%%%%%%%%%%%%%%%%%%%%%%%%%%%%%%%%%%%%%%%%%%%%%%%%%%%%%%%%%%%

 %%%%%%%%%%%%%%%%%%%%%%%%%%%%%%%%%%%%%%%%%%%%%%%%%%%%%%%%%%%%%%%%%%%%%%%%%%
\begin{defi}[\textbf{Grafting}]
Given two trees $\vartheta$ and $\vartheta'$ we say that the tree $\vartheta''$ is obtained
by grafting $\vartheta'$ to a node $v\in N(\vartheta)$ if we identify
$r_{\vartheta'}$ with $v$ in such a way that $\ell_{\vartheta'}$ enters $v$
{and the value of $s_v$ increases by one}. % and $\vartheta'\subset\vartheta''$.
\end{defi}
 %%%%%%%%%%%%%%%%%%%%%%%%%%%%%%%%%%%%%%%%%%%%%%%%%%%%%%%%%%%%%%%%%%%%%%%%%%
 
 %%%%%%%%%%%%%%%%%%%%%%%%%%%%%%%%%%%%%%%%%%%%%%%%%%%%%%%%%%%%%%%%%%%%%%%%%%
\begin{defi}[\textbf{Cutting}]
If $\vartheta'$ is a subtree of $\vartheta$ and $v\in N(\vartheta)$ is the node $\ell_{\vartheta'}$ enters,
 we say that the tree $\vartheta''$ is obtained from $\vartheta$ by cutting
 $\vartheta'$ if
$N(\vartheta'')=N(\vartheta)\setminus N(\vartheta')$ and $L(\vartheta'') = L(\vartheta)\setminus L(\vartheta')$,
in such a way that $s_v$ decreases by one.
\end{defi}
 %%%%%%%%%%%%%%%%%%%%%%%%%%%%%%%%%%%%%%%%%%%%%%%%%%%%%%%%%%%%%%%%%%%%%%%%%%

%%%%%%%%%%%%%%%%%%%%%%%%%%%%%%%%%%%%%%%%%%%%%%%%%%%%%%%%%%%%%%%%%%%%%%%%%%
\begin{rmk} \label{nodiefoglie}
\emph{
We draw as black bullets the nodes $v$ with $s_v > 2$, as white squares the nodes $v$ with $s_v\le 2$,
as white bullets the leaves and as a grey bullet the root (see Figure \ref{superalcolici} for an explicit example).
The lines entering a node are drawn so that the lines with sign ``$-$''
are above %to the left of
the lines with sign ``$+$'
(see Figure \ref{caffe1} for an explicit example).
}
\end{rmk}
%%%%%%%%%%%%%%%%%%%%%%%%%%%%%%%%%%%%%%%%%%%%%%%%%%%%%%%%%%%%%%%%%%%%%%%%%%

Our aim is to represent both the solution and the counterterm as ``sum over trees'' (in a sense that will be
clear later on). We shall introduce several sets of trees, which differ because of
the numbers of lines allowed to enter the nodes and the constraints to be imposed on the labels.

%%%%%%%%%%%%%%%%%%%%%%%%%%%%%%%%%%%%%%%%%%%%%%%%%%%%%%%%%%%%%%%%%%%%%%%%%%
\subsection{The range equation: unexpanded trees}
\label{sec:range}
%%%%%%%%%%%%%%%%%%%%%%%%%%%%%%%%%%%%%%%%%%%%%%%%%%%%%%%%%%%%%%%%%%%%%%%%%%

Fix once and for all $s\ge0$ and $\al\in(0,1)$ and consider any $c\in \mathtt{g}(s,\al)$ with $\|c\|_{s,\al}\le1$.
Let us consider first the range equation \eqref{ran} alone, and look for a solution for any fixed $\h$: the exact value of $\h$ is found
by requiring the kernel equation \eqref{bif} to be solved as well order by order. For this reason
it is more convenient to think of $\h$  as a sequence in both indices $j$ and $k$, i.e.
\[
\h :=\{\h^{(k)}_j\}_{j\in\ZZZ, k\ge1}\in\CCC^{\ZZZ\times\NNN},
\]
and look for a formal solution of \eqref{ran} with $\h\in\CCC^{\ZZZ\times\NNN}$ as a free parameter.

To avoid confusion
we call $w=w(x,\f;c,\om,\e,\h)$ the solution to \eqref{ran} alone, that we write as a formal Lindstedt series 
\begin{equation}\label{fw}
w (x,\f;c,\om,\h,\e)= \sum_{j\in\ZZZ}c_je^{\ii j x} + 
\sum_{j\in\ZZZ} \sum_{\substack{\nu\in\ZZZ^{\ZZZ}_{f}\\ \nu\ne \gote_j}} 
\sum_{k\ge1}\e^k
w_{j,\nu}^{(k)}(c,\om,\h) \, e^{\ii (j x+\nu\cdot\f)}\,,
\end{equation}
Note that $w$ depends, among others, on the parameter $\h$, and, in general, if we set ${U}=w$
then $U$ does not solve \eqref{gene}. 
For this to happen we shall need to fix $\h$ in a proper way.

Consider any $\om\in\gotB^{(0)}$ and fix any $r\in\gotR$ such that $\BB_\om^{(0)}(r)<\io$.
Inserting \eqref{fw} into \eqref{ran}, for $k=1$ we obtain, for $\nu\ne\gote_j$,
\begin{equation}\label{ran1}
w_{j,\nu}^{(1)} (c,\om,\h) =   \frac{1}{(\om \cdot \nu  - \om_j) }
\sum_{\substack{ \gote_{j_1}-\gote_{j_2} + \gote_{j_3} - \gote_{j_4} + \gote_{j_5}  =\nu \\ j_1-j_2+j_3 - j_4 +j_5 = j }}
c_{j_1}\ol{c}_{j_2}c_{j_3}\ol{c}_{j_4}c_{j_5}\,.
\end{equation}

%%%%%%%%%%%%%%%%%%%%%%%%%%%%%%%%%%%%%%%%%%%%%%%%%%%%%%%%%%%%%%%%%%%%%%%%%%
\begin{rmk} \label{+-} 
\emph{%remark
Throughout the paper, we shall use equivalently the notation $z=z^+$ and 
$\ol{z}=z^-$ for any $z\in\CCC$, where $\ol{z}$ is the complex conjugate of $z$.
}%remark
\end{rmk}
%%%%%%%%%%%%%%%%%%%%%%%%%%%%%%%%%%%%%%%%%%%%%%%%%%%%%%%%%%%%%%%%%%%%%%%%%%

%
It may be convenient to write explicitly also the expression for $\ol{w}^{(1)}_{j,\nu}$ as
\begin{equation}\label{u1bar}
\ol{w}_{j,\nu}^{(1)} (c,\om,\h) =   \frac{1}{(\om \cdot \nu  - \om_j) }
\sum_{\substack{ - \gote_{j_1} +\gote_{j_2} - \gote_{j_3} + \gote_{j_4} - \gote_{j_5}  = - \nu \\ 
-j_1+ j_2 - j_3 + j_4 - j_5 = - j }}
\ol{c}_{j_1}{c}_{j_2}\ol{c}_{j_3}{c}_{j_4}\ol{c}_{j_5}\,.
\end{equation}

%%%%%%%%%%%%%%%%%%%%%%%%%%%%%%%%%%%%%%%%%%%%%%%%%%%%%%%%%%%%%%%%%%%%%%%%
\begin{rmk} \label{serve?}
\emph{
In fact, the coefficients $w^{(k)}_{j,\nu}(c,\om,\h)$ with $k=1$ depend only on $\om$ and $c$;
on the contrary, for $k \ge 2$, they do depend explicitly also on $\h$, as \eqref{ran2} below shows already for $k=2$.
}
\end{rmk}
%%%%%%%%%%%%%%%%%%%%%%%%%%%%%%%%%%%%%%%%%%%%%%%%%%%%%%%%%%%%%%%%%%%%%%%%

To simplify the forthcoming discussion, we introduce also the coefficients %$w^{(0)}_{j,\nu}(c,\om,\h)$ by setting
\begin{equation} \label{questoserve}
w^{(0)}_{j,\nu}(c,\om,\h) := \begin{cases}
c_j , & \quad \nu = \gote_j , \\
0 , & \quad \nu\neq\gote_j ,
\end{cases}
\end{equation}
so that, if we further require $w^{(k)}_{j,\gote_j}(c,\om,\h)$ to vanish for all $j\in\ZZZ$ and all $k\ge 1$, we can rewrite \eqref{fw} as
\begin{equation} \nonumber 
w (x,\f;c,\om,\h,\e)= \sum_{j\in\ZZZ} \sum_{\nu\in\ZZZ^{\ZZZ}_{f}} 
\sum_{k\ge0}\e^k
w_{j,\nu}^{(k)}(c,\om,\h) \, e^{\ii (j x+\nu\cdot\f)} \, .
\end{equation}

%%%%%%%%%%%%%%%%%%%%%%%%%%%%%%%%%%%%%%%%%%%%%%%%%%%%%%%%%%%%%%%%%%%%%%%%
\begin{rmk}\label{formalpoly}
\emph{%remark
Both coefficients \eqref{ran1} and \eqref{u1bar}, as well as the forthcoming coefficients \eqref{ran2} and \eqref{rank},
are not finite sums, and indeed they have to be considered as sums of infinitely many polynomials in $c$ of finite degree.
Of course, at least in \eqref{ran1} and \eqref{u1bar}, the condition \eqref{decay} is enough to ensure convergence.
}%remark
\end{rmk}
%%%%%%%%%%%%%%%%%%%%%%%%%%%%%%%%%%%%%%%%%%%%%%%%%%%%%%%%%%%%%%%%%%%%%%%%

We introduce a graphical representation of the coefficients $w^{(k)}_{j,\nu}(c,\om,\h)$ in terms of trees.
The aim is to express the coefficients as sums of several contributions, each of which is given by the
``value'' of a tree, in such a way to make easier to get bounds and single out the trees
which produce bad bounds and hence hamper the convergence of the series.

\begin{figure}[H]
\vspace{.2cm}
\centering
\ins{37pt}{-010pt}{$w^{(k)}_{j,\nu}(c,\om,\h)=$}
\ins{145pt}{-022pt}{$j \,\gote_j+$}
\ins{193pt}{-000pt}{$k$}
\ins{248pt}{-014pt}{$c_j^\s=$ }
\ins{301pt}{-022pt}{$j \,\gote_j\s$}
\ins{355pt}{-002pt}{$j\s$}
%\ins{180pt}{-064pt}{$=$}
\null
\hspace{.5cm}
\subfigure{\includegraphics*[width=1.4in]{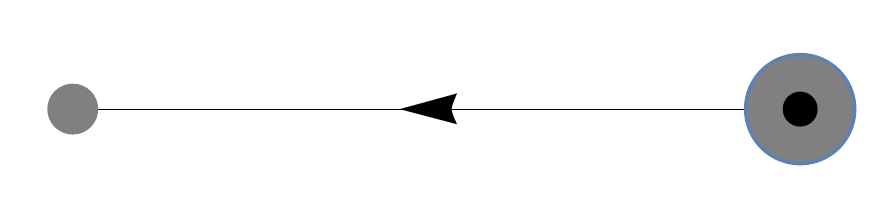}}
\hspace{2.0cm}
\subfigure{\includegraphics*[width=1.4in]{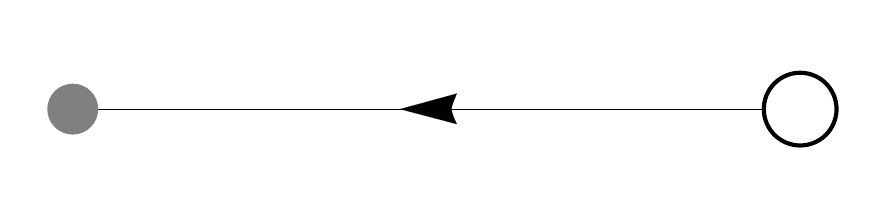}}
\caption{\small Graphical representation of  $w^{(k)}_{j,\nu}(c,\om,\h)$ and $c^\s_j$.}
\label{venticello}
\end{figure}

We start by representing each coefficient $w^{(k)}_{j,\nu}(c,\om,\h)$ as a tree with a line exiting a grey bullet 
with a smaller black bullet inside, and each coefficient $c_j^\s$ as a tree with a line exiting a leaf,
drawn as a white bullet (see Figure \ref{venticello}).

%%%%%%%%%%%%%%%%%%%%%%%%%%%%%%%%%%%%%%%%%%%%%%%%%%%%%%%%%%%%%%%%%%%%%%%%
% FIGURA 19-1
%%%%%%%%%%%%%%%%%%%%%%%%%%%%%%%%%%%%%%%%%%%%%%%%%%%%%%%%%%%%%%%%%%%%%%%%
%%%%%%%%%%%%%%%%%%%%%%%%%%%%%%%%%%%%%%%%%%%%%%%%%%%%%%%%%%%%%%%%%%%%%%%%

\begin{figure}[H]
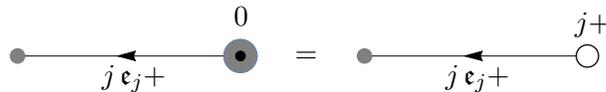

\vspace{.2cm}
\centering
\ins{145pt}{-022pt}{$j \,\gote_j+$}
\ins{195pt}{-000pt}{$0$}
\ins{218pt}{-017pt}{$=$}
\ins{275pt}{-022pt}{$j \,\gote_j+$}
\ins{324pt}{-002pt}{$j+$}
%\ins{180pt}{-064pt}{$=$}
\null
\hspace{-.6cm}
\subfigure{\includegraphics*[width=1.4in]{figura0a}}
\hspace{.8cm}
\subfigure{\includegraphics*[width=1.4in]{figura0}}
\caption{\small Graphical representation of the idenity $w^{(0)}_{j,\nu}(c,\om,\h)=c^+_j$.}
\label{mancava}
\end{figure}

Then we may represent the identity $w^{(0)}_{j,\gote_j}(c,\om,\h)=c^+_j$ in \eqref{questoserve} as in Figure \ref{mancava}.

%%%%%%%%%%%%%%%%%%%%%%%%%%%%%%%%%%%%%%%%%%%%%%%%%%%%%%%%%%%%%%%%%%%%%%%%
% FIGURA 19-1
%%%%%%%%%%%%%%%%%%%%%%%%%%%%%%%%%%%%%%%%%%%%%%%%%%%%%%%%%%%%%%%%%%%%%%%%
%%%%%%%%%%%%%%%%%%%%%%%%%%%%%%%%%%%%%%%%%%%%%%%%%%%%%%%%%%%%%%%%%%%%%%%%

\begin{figure}[ht]
\vspace{.2cm}
\centering
\ins{120pt}{-068pt}{$j \,\nu \,+$}
\ins{169pt}{-046pt}{$1$}
\ins{195pt}{-064pt}{$=$}
\ins{240pt}{-068pt}{$j\,\nu\,+$}
\ins{327pt}{-035pt}{$j_2$}
\ins{335pt}{-031pt}{$\gote_{j_2}$}
\ins{341pt}{-024pt}{$-$}
\ins{327pt}{-050pt}{$j_4$}
\ins{336pt}{-048pt}{$\gote_{j_4}$}
\ins{344pt}{-042pt}{$-$}
\ins{327pt}{-066.5pt}{$j_1 \gote_{j_1} \! +$}
\ins{327pt}{-084pt}{$j_3$}
\ins{336pt}{-091pt}{$\gote_{j_3}$}
\ins{346pt}{-097pt}{$+$}
\ins{322pt}{-098pt}{$j_5$}
\ins{331pt}{-106pt}{$\gote_{j_5}$}
\ins{341pt}{-112pt}{$+$}
\ins{352pt}{+002pt}{$j_2 -$}
\ins{352pt}{-023pt}{$j_4 -$}
\ins{352pt}{-050pt}{$j_1 +$}
\ins{352pt}{-078pt}{$j_3 +$}  
\ins{352pt}{-102pt}{$j_5 +$}
\null
\hspace{-.6cm}
\subfigure{\includegraphics*[width=1.4in]{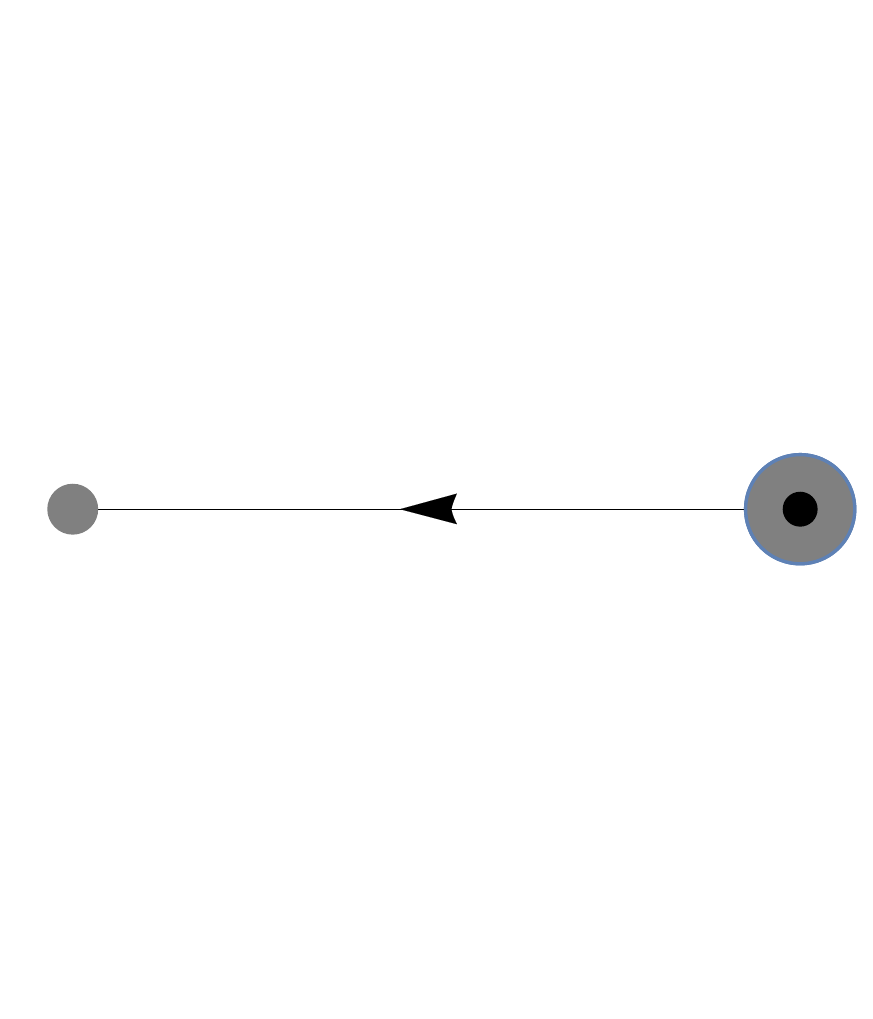}}
\hspace{.8cm}
\subfigure{\includegraphics*[width=2.1in]{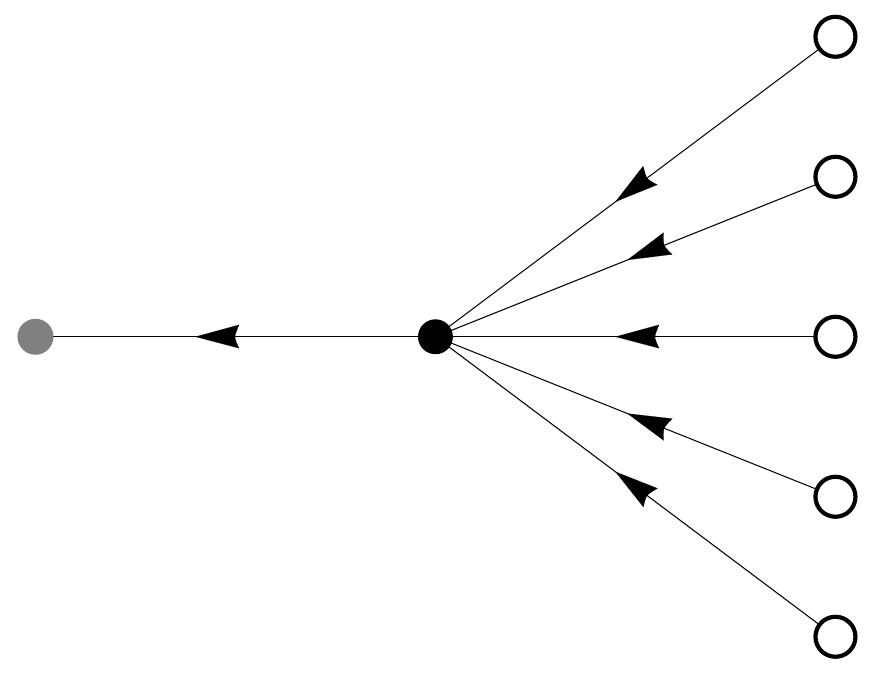}}
\caption{\small Graphical representation of \eqref{ran1}.}
\label{caffe1}
\end{figure}

We represent graphically the r.h.s of \eqref{ran1} as a tree with 1 node $v$ with $k_v=1$ and 5 leaves  (see Figure \ref{caffe1}). 
We represent the node as a black bullet and again the leaves as white bullets, according to Remark \ref{nodiefoglie}.
By construction, such a tree is of order $1$.
The sum appearing in \eqref{ran1} is not explicitly written in Figure \ref{caffe1}.

%%%%%%%%%%%%%%%%%%%%%%%%%%%%%%%%%%%%%%%%%%%%%%%%%%%%%%%%%%%%%%%%%%%%%%%%
% FIGURA 19-2
%%%%%%%%%%%%%%%%%%%%%%%%%%%%%%%%%%%%%%%%%%%%%%%%%%%%%%%%%%%%%%%%%%%%%%%%
%%%%%%%%%%%%%%%%%%%%%%%%%%%%%%%%%%%%%%%%%%%%%%%%%%%%%%%%%%%%%%%%%%%%%%%%

\begin{figure}[ht]
\vspace{-.2cm}
\centering
\ins{046pt}{-068pt}{$j \,\nu \,+$}
\ins{097pt}{-046pt}{$k$}
\ins{116pt}{-064.5pt}{$=$}
\ins{154pt}{-068pt}{$j\,\nu\,+$}
\ins{219pt}{-068pt}{$j\,\nu\,+$}
\ins{314pt}{-068pt}{$j\,\nu\,+$}
\ins{184pt}{-046pt}{$k_1\,j\,+$}
\ins{258pt}{-046pt}{$k_2$}
\ins{276.5pt}{-063pt}{$+$}
\ins{397pt}{-038pt}{$j_2$}
\ins{405pt}{-035pt}{$\nu_2$}
\ins{411pt}{-027pt}{$-$}
\ins{395pt}{-054pt}{$j_4$}
\ins{403pt}{-052pt}{$\nu_4$}
\ins{412pt}{-045pt}{$-$}
\ins{390pt}{-067.5pt}{$j_1 \nu_1 +$}
\ins{397pt}{-083.5pt}{$j_3$}
\ins{406pt}{-092pt}{$\nu_3$}
\ins{414pt}{-097pt}{$+$}
\ins{392pt}{-095pt}{$j_5$}
\ins{401pt}{-104pt}{$\nu_5$}
\ins{408pt}{-109pt}{$+$}
\ins{427pt}{+000pt}{$k_2$}
\ins{427pt}{-025pt}{$k_4$}
\ins{427pt}{-049pt}{$k_1$}
\ins{427pt}{-075pt}{$k_3$}  
\ins{427pt}{-099pt}{$k_5$}
\null
\hspace{-.6cm}
\subfigure{\includegraphics*[width=1.4in]{figura1a}}
\hspace{.4cm}
\subfigure{\includegraphics*[width=2.0in]{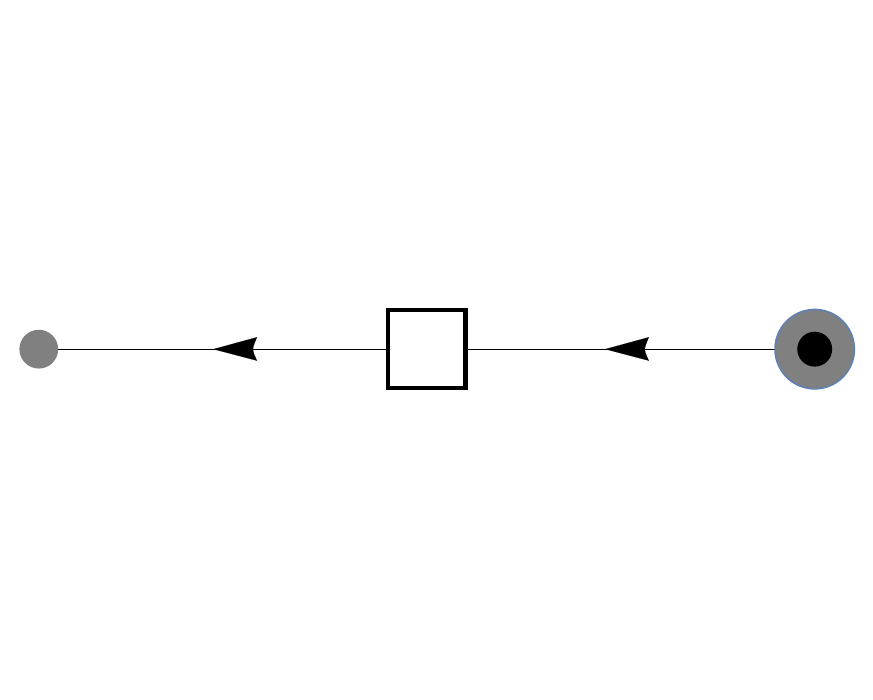}}
\hspace{.4cm}
\subfigure{\includegraphics*[width=2.0in]{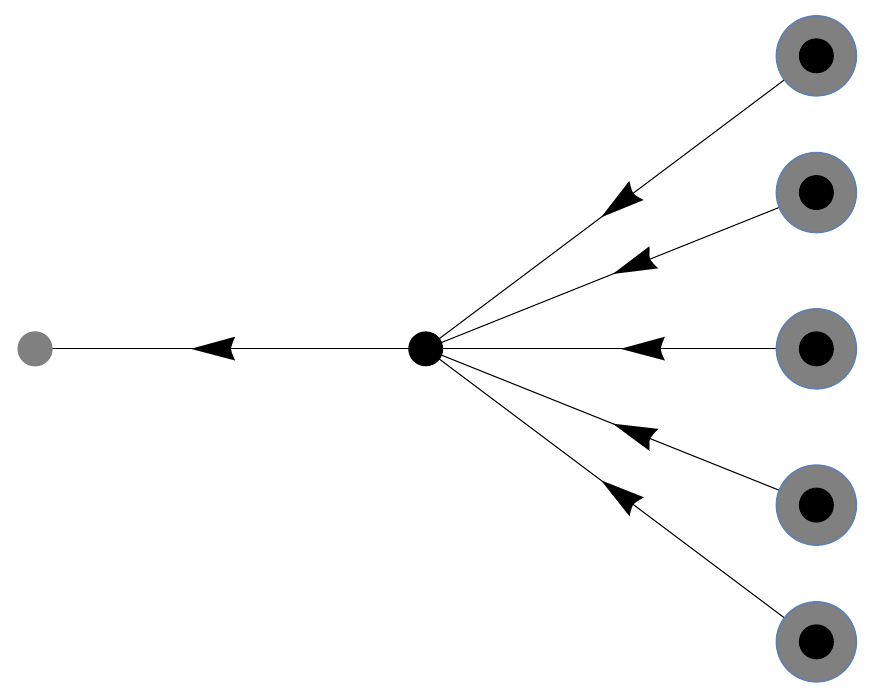}}
\caption{\small Graphical representation of \eqref{rank} (and of \eqref{ran2} for $k=2$).}
\label{no!}
\end{figure}

Things start getting more interesting for $k=2$; indeed we find
\begin{equation}\label{ran2}
w_{j,\nu}^{(2)}  = \frac{1}{(\om\cdot\nu-\om_j)}\Bigg[\h^{(1)}_j w_{j,\nu}^{(1)} +
\sum_{\substack{ \nu_1 -\nu_2 + \nu_3-\nu_4 +\nu_5 =\nu \\ j_1-j_2+j_3 - j_4 +j_5 = j  \\
k_1+k_2+k_3+k_4+k_5=1}}
w_{j_1,\nu_1}^{(k_1)}\ol{w}_{j_2,\nu_2}^{(k_2)}w_{j_3,\nu_3}^{(k_3)}\ol{w}_{j_4,\nu_4}^{(k_4)}w_{j_5,\nu_5}^{(k_5)}
\Bigg]\, ,
\end{equation}
where we have not made explicit the dependence of the coefficients on $c,\om,\h$ in order not to overwhelm further the notation.
Note that the constraint on the labels $k_i$ in the sums appearing in \eqref{ran2} implies that only one among
them equals $1$ while the others are zero. We represent
\[
\frac{1}{(\om\cdot\nu-\om_j)}\h^{(1)}_j w_{j,\nu}^{(1)}
\] 
as a tree with the root line exiting a node $v$ with $k_v=1$, drawn as a white square (recall Remark \ref{nodiefoglie}),
and a subtree of order 1 as in Figure \ref{caffe1} entering the node $v$;
such a tree is a special case of the tree represented in the middle of Figure \ref{no!}, with $k_1=k_2=1$.
Instead, we represent 
\[
\frac{1}{(\om\cdot\nu-\om_j)}w_{j_1,\nu_1}^{(k_1)}\ol{w}_{j_2,\nu_2}^{(k_2)}w_{j_3,\nu_3}^{(k_3)}\ol{w}_{j_4,\nu_4}^{(k_4)}w_{j_5,\nu_5}^{(k_5)}
\]
as a tree with 1 node $v$ with $k_v=1$ and 5 subtrees entering $v$; such a tree corresponds to the tree represented
to the right in Figure \ref{no!}, in the case in which one has $k_i=1$ for some $i\in\{1,\ldots,5\}$ and $k_j=0$ for $j\neq i$.
The vertices with $k_i=0$ are leaves and hence are drawn as white bullets, according to Remark \ref{nodiefoglie} and Figure \ref{mancava}.

%%%%%%%%%%%%%%%%%%%%%%%%%%%%%%%%%%%%%%%%%%%%%%%%%%%%%%%%%%%%%%%%%%%%%%%%
% FIGURA 19-1
%%%%%%%%%%%%%%%%%%%%%%%%%%%%%%%%%%%%%%%%%%%%%%%%%%%%%%%%%%%%%%%%%%%%%%%%
%%%%%%%%%%%%%%%%%%%%%%%%%%%%%%%%%%%%%%%%%%%%%%%%%%%%%%%%%%%%%%%%%%%%%%%%

More generally, for $k\ge2$ we obtain
\begin{equation}\label{rank}
w_{j,\nu}^{(k)} = \frac{1}{(\om\cdot\nu-\om_j)}\Bigg[
\sum_{k_1+k_2=k} \h^{(k_1)}_j w^{(k_2)}_{j,\nu} + 
\!\!\!\!\!\!\!\!\!\!
\sum_{\substack{ \nu_1 -\nu_2 + \nu_3-\nu_4 +\nu_5 =\nu \\ j_1-j_2+j_3 - j_4 +j_5 = j  \\
k_1+k_2+k_3+k_4+k_5=k}}
\!\!\!\!\!\!\!\!\!\!
w_{j_1,\nu_1}^{(k_1)}\ol{w}_{j_2,\nu_2}^{(k_2)}w_{j_3,\nu_3}^{(k_3)}\ol{w}_{j_4,\nu_4}^{(k_4)}w_{j_5,\nu_5}^{(k_5)} 
\Bigg] ,
\end{equation}
and we represent graphically \eqref{rank} as in Figure \ref{no!}; again, 
the sum over the labels is understood, and,
when $i\in\{1,\ldots,5\}$ is such that $k_i=0$, then the corresponding vertex is a leaf.
We also associate with the node drawn as a  white square the labels $k_1,j,+$.

Analogous graphical representations can be easily obtained for $\bar w^{(k)}_{j,\nu}$,
by changing the labels ``$+$'' with ``$-$'' and vice versa,
and reverse the order of all the lines entering the same node so that the lines
with sign ``$-$'' are still above %to the left of
the lines with sign ``$+$''
(see the comments at the end of Subsection \ref{alberi}).

%%%%%%%%%%%%%%%%%%%%%%%%%%%%%%%%%%%%%%%%%%%%%%%%%%%%%%%%%%%%%%%%%%%%%%%%%%
\begin{rmk}\label{disegno}
\emph{
As Figure \ref{no!} shows, trees are drawn horizontally with the root to the left and the arrows pointing from the right to the left.
Hence, if a node $v$ precedes a node $w$ (i.e.~$v \preceq w$), then $v$ is to the right of $w$.
}
\end{rmk}
%%%%%%%%%%%%%%%%%%%%%%%%%%%%%%%%%%%%%%%%%%%%%%%%%%%%%%%%%%%%%%%%%%%%%%%%%%

%%%%%%%%%%%%%%%%%%%%%%%%%%%%%%%%%%%%%%%%%%%%%%%%%%%%%%%%%%%%%%%%%%%%%%%%%%
\begin{rmk}\label{momento}
\emph{
One may easily check by induction that, due to the constraints in the sums, $w_{j,\nu}^{(k)}\neq0$ requires
\begin{equation}\label{pigreco}
j = 
\pi(\nu):=\sum_{i\in\ZZZ} i \nu_i\,.
\end{equation}
In particular, this means that $\nu$ fixes $j$. We do not provide here the details, since
the assertion will be proved later as a consequence of Lemmata \ref{nonsocosa} and \ref{boh} (see Remark \ref{noneradetto}).
}
\end{rmk}
%%%%%%%%%%%%%%%%%%%%%%%%%%%%%%%%%%%%%%%%%%%%%%%%%%%%%%%%%%%%%%%%%%%%%%%%%%

\begin{figure}[ht]
\vspace{-.2cm}
\centering
%\null
%\hspace{-.6cm}
\ins{056pt}{-108pt}{$r$}
%\ins{068pt}{-085pt}{$\ell^{out}_\TT$}
%\ins{378pt}{-044pt}{$v_{in}$}
%\ins{346pt}{-058pt}{$\ell^{in}_\TT$}
\subfigure{\includegraphics*[width=5in]{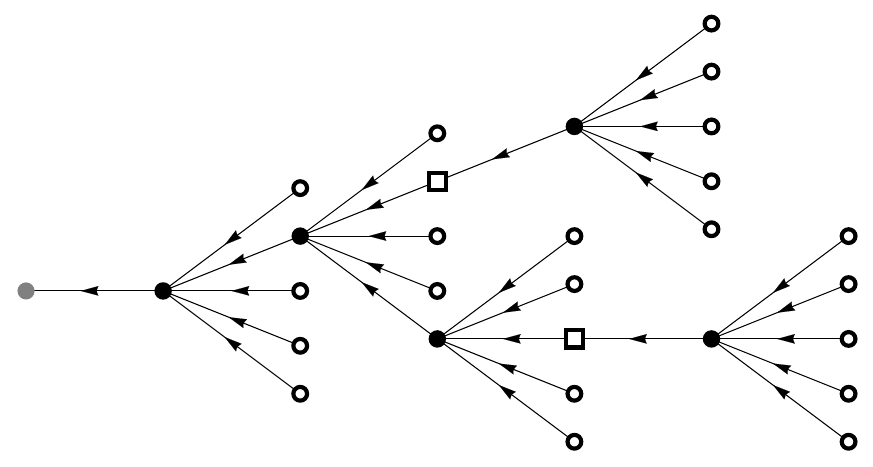}}
\caption{\small A tree $\vartheta\in\Theta$ with 7 nodes and 21 leaves (the labels are not shown).}
\label{superalcolici}
\end{figure}

By iterating the graphical representation in Figure \ref{no!} until only leaves
and subtrees of order 1 are left, which are represented according to Figures \ref{mancava} and \ref{caffe1}, respectively,
we obtain a graphical representation of the coefficients $w^{(k)}_{j,\nu}$ in \eqref{fw} in terms of trees,
as defined in Subsection \ref{alberi}, by requiring further constraints on the labels.
To this aim, we introduce the first set of trees as follows
(refer to Figure \ref{superalcolici} for an example).

%%%%%%%%%%%%%%%%%%%%%%%%%%%%%%%%%%%%%%%%%%%%%%%%%%%%%%%%%%%%%%%%%%%%%%%%%%
% FIGURA 19-2
%%%%%%%%%%%%%%%%%%%%%%%%%%%%%%%%%%%%%%%%%%%%%%%%%%%%%%%%%%%%%%%%%%%%%%%%%%
%%%%%%%%%%%%%%%%%%%%%%%%%%%%%%%%%%%%%%%%%%%%%%%%%%%%%%%%%%%%%%%%%%%%%%%%%%

%%%%%%%%%%%%%%%%%%%%%%%%%%%%%%%%%%%%%%%%%%%%%%%%%%%%%%%%%%%%%%%%%%%%%%%%%%
\begin{defi}[\textbf{Set $\boldsymbol{\Theta}$ of the unexpanded trees}]
\label{thetone}
Let $\Theta$ be the set of trees $\vartheta$ satisfying the following constraints:

\begin{enumerate}[topsep=0ex]
\itemsep0.0em
\item for any node $v\in N(\vartheta)$ one has $s_v=1$ or $s_v=5$;
\item for any node $v\in N(\vartheta)$ one has $k_v\ge1$ if $s_v=1$ and $k_v=1$ if $s_v=5$;
\item\label{thetone0}  if $\ell \in L(\vartheta)$ exits a leaf $\la\in\Lambda(\vartheta)$, then $\s_\ell =\s_\la $, $j_\ell= j_\la$ and $\nu_\ell=\gote_{j_\la}$;
\item\label{thetone1} if $\ell \in L(\vartheta)$ exits a node $v\in N(\vartheta)$ and $\ell$ is not the root line, then $\nu_\ell\ne\gote_{j_\ell}$;
\item\label{thetone2} if $\ell \in L(\vartheta)$ exits a node $v \in N(\vartheta)$, then
\begin{equation}\label{conserva}
\s_{\ell}1 = \sum_{i=1}^{s_v}\s_{\ell_i}1,\qquad
\s_\ell \nu_\ell = \sum_{i=1}^{s_v}\s_{\ell_i}\nu_{\ell_i},\qquad 
\s_\ell j_\ell =
\sum_{i=1}^{s_v} \s_{\ell_i}j_{\ell_i} ,
\end{equation}
where $\ell_1,\ldots,\ell_{s_v}$ are the lines entering $v$;
\item if $\ell \in L(\vartheta)$ exits a node $v\in N(\vartheta)$, then one has $\s_v=\s_\ell$ and $j_v = j_\ell$;
\item for any line $\ell\in L(\vartheta)$ one has $n_\ell=-1$ if $\nu_\ell=\gote_{j_\ell}$ and $n_\ell\ge0$ otherwise.
\end{enumerate}
We call any $\vartheta\in\Theta$ an \emph{unexpanded tree},
and define $\Theta^{(k)}_{j,\nu,\s}$ as the set of the unexpanded trees $\vartheta\in\Theta$ with order $k(\vartheta)=k$
such that the root line has component $j$, momentum $\nu$ and sign $\s$.
\end{defi}
%%%%%%%%%%%%%%%%%%%%%%%%%%%%%%%%%%%%%%%%%%%%%%%%%%%%%%%%%%%%%%%%%%%%%%%%%%

%%%%%%%%%%%%%%%%%%%%%%%%%%%%%%%%%%%%%%%%%%%%%%%%%%%%%%%%%%%%%%%%%%%%%%%%%%
\begin{rmk} \label{noetaleaf}
\emph{
{
The constraints \ref{thetone0}, \ref{thetone1} and \ref{thetone2} in Definition \ref{thetone} imply that
if $\ell \in L(\vartheta)$ exits a leaf $\la\in\Lambda(\vartheta)$ and enters a node $v\in N(\vartheta)$ then $s_v=5$,
the only exception being the tree $\vartheta$ with $|N(\vartheta)|=|\Lambda(\vartheta)|=1$.
}
}
\end{rmk}
%%%%%%%%%%%%%%%%%%%%%%%%%%%%%%%%%%%%%%%%%%%%%%%%%%%%%%%%%%%%%%%%%%%%%%%%%%

%%%%%%%%%%%%%%%%%%%%%%%%%%%%%%%%%%%%%%%%%%%%%%%%%%%%%%%%%%%%%%%%%%%%%%%%%%
\begin{rmk}\label{importante}
\emph{
Iterating the equations in \eqref{conserva} one finds, for all $\ell\in L(\vartheta)$,
\begin{equation}\label{lefoglie}
        \s_\ell 1 = \sum_{{\substack{\la\in\Lambda(\vartheta) \\ \la\preceq \ell}}} \s_{\la} 1 ,\qquad
	\s_\ell \nu_\ell = \sum_{{\substack{\la\in\Lambda(\vartheta) \\ \la\preceq \ell}}}  \s_{\la} \gote_{j_\la},\qquad
	\s_\ell j_\ell=\sum_{{\substack{\la\in\Lambda(\vartheta) \\ \la\preceq \ell}}}  \s_{\la} {j_\la}\,.
\end{equation}
In \eqref{lefoglie}, the first and second relations imply that, for any $\ell\in\Lambda(\vartheta)$,
\begin{equation}\label{tocchera}
\sum_{i\in\ZZZ} (\nu_\ell)_i =1,
\end{equation}
which reflects the gauge-covariance mentioned in Remark \ref{extrabryuno},
while the last relation is due to the translation-covariance and,
as noted in Remark \ref{mettiamoloqui}, plays a crucial role in the following.
}
\end{rmk}
%%%%%%%%%%%%%%%%%%%%%%%%%%%%%%%%%%%%%%%%%%%%%%%%%%%%%%%%%%%%%%%%%%%%%%%%%%

%%%%%%%%%%%%%%%%%%%%%%%%%%%%%%%%%%%%%%%%%%%%%%%%%%%%%%%%%%%%%%%%%%%%%%%%%%
\begin{rmk}\label{cc}
\emph{
Any tree in $\Theta^{(k)}_{j,\nu,-}$ can be obtained
by a tree in $\Theta^{(k)}_{j,\nu,+}$ by changing the sign labels of all nodes, leaves and lines
and reversing the order all the lines entering the same node (so that the lines
with sign ``$-$'' are still above %to the left of 
the lines with sign ``$+$''). 
}
\end{rmk}
%%%%%%%%%%%%%%%%%%%%%%%%%%%%%%%%%%%%%%%%%%%%%%%%%%%%%%%%%%%%%%%%%%%%%%%%%%

Next, we prove a relation which, combined with the forthcoming Lemma \ref{boh}, gives the property mentioned in Remark \ref{momento}.

%%%%%%%%%%%%%%%%%%%%%%%%%%%%%%%%%%%%%%%%%%%%%%%%%%%%%%%%%%%%%%%%%%%%%%%%%%
\begin{lemma}\label{nonsocosa}
Let $\vartheta$ be a tree in $\Theta$.
For all $\ell \in L(\vartheta)$ one has
\begin{equation} \nonumber %\label{pi}
j_{\ell} = \pi(\nu_{\ell}) = \sum_{i\in\ZZZ} i (\nu_\ell)_i\,.
\end{equation}
\end{lemma}
%%%%%%%%%%%%%%%%%%%%%%%%%%%%%%%%%%%%%%%%%%%%%%%%%%%%%%%%%%%%%%%%%%%%%%%%%%

%%%%%%%%%%%%%%%%%%%%%%%%%%%%%%%%%%%%%%%%%%%%%%%%%%%%%%%%%%%%%%%%%%%%%%%%%%
\prova If $\ell$ exists a leaf $\la$, one has $(\nu_\ell)_i=(\nu_\la)_i=(\gote_{j_\la})_i=\de_{j_\la i}$,
so that $\pi(\nu_\ell)=j_\la=j_\ell$. If $\ell$ exits a node,
one proceeds by induction on the order $k$ of the subtree $\vartheta$ which has $\ell$ as root line.

If $k=1$, $\vartheta$ has only 1 node with $5$ entering lines 
which exit as many leaves, so that
\begin{equation}\nonumber
\s_\ell (\nu_\ell)_i = \sum_{\la \in \Lambda(\vartheta)} \s_\la \de_{j_\la i} ,
\qquad
\s_\ell j_\ell = \sum_{\la \in \Lambda(\vartheta)} \s_\la j_\la ,
\end{equation}
and hence
\begin{equation}\nonumber
\s_\ell \pi(\nu_\ell) := \sum_{i \in \ZZZ}  i \s_\ell (\nu_\ell)_i = \sum_{i \in \ZZZ} \sum_{\la \in \Lambda(\vartheta)} i \s_\la \de_{j_\la i} = 
\sum_{\la \in \Lambda(\vartheta)} \s_\la j_\la = \s_\ell j .
\end{equation}

If $k>1$, let $v$ be the node $\ell$ exits:
\begin{itemize}[topsep=0ex]
\itemsep0.0em
\item if $s_v=1$, then the line entering $v$ has the same labels as $\ell$
and is the root line of a tree of order $\le k-1$, so that the result follows from the inductive hypothesis;
\item if $s_v=5$, and $\ell_1,\ldots,\ell_5$ denote the lines entering $v$, by \eqref{conserva} one has
\begin{equation}\nonumber
\s_\ell (\nu_\ell)_i = \sum_{p=1}^5 \s_{\ell_p} (\nu_{\ell_p})_i ,
\qquad
\s_\ell j _\ell= \sum_{p=1}^5 \s_{\ell_p} j_{\ell_p}
%=\sum_{p=1}^5 \s_{\ell_p} \pi(\nu_{\ell_p}) ,
%= \sum_{p=1}^5 \sum_{i\in\ZZZ} i \s_{\ell_p} (\nu_{\ell_p)_i} ,
\end{equation}
which imply, again by using the inductive hypothesis,
\begin{equation}\nonumber
\s_\ell \pi(\nu_\ell) := \sum_{i \in \ZZZ}  i \s_\ell (\nu_\ell)_i = 
\sum_{i \in \ZZZ}  \sum_{p=1}^5 i \s_{\ell_p} (\nu_{\ell_p})_i 
= \sum_{p=1}^5 \s_{\ell_p} \pi(\nu_{\ell_p}) = \sum_{p=1}^5 \s_{\ell_p} j_{\ell_p} = \s_\ell j_\ell .
\end{equation}
\end{itemize}
This proves the assertion.
\EP
%%%%%%%%%%%%%%%%%%%%%%%%%%%%%%%%%%%%%%%%%%%%%%%%%%%%%%%%%%%%%%%%%%%%%%%%%%

Let $\vartheta$ be a tree in $\Theta$. With each leaf $\la\in\Lambda(\vartheta)$ we associate a \textit{leaf factor}
\begin{equation}\label{foglie}
\LL_\la =\LL_\la(c):= c_{j_\la}^{\s_\la}\,,
\end{equation}
with each node $v\in N(\vartheta)$ we associate a \textit{node factor} 
\begin{equation}\label{nodi}
\calF_v = \calF_v(\h) := %\left\{
\begin{cases}
1, & \quad 
s_v =5\\
(\h_{j_v}^{(k_v)})^{\s_v} , & \quad 
s_v=1 ,
\end{cases}
%\right.
\end{equation}
and
with each line $\ell\in L(\vartheta)$ we associate a \textit{small divisor}\footnote{The small divisor \eqref{simbolo}
is the symbol of the differential operator $\gotD(\om)$ defined in \eqref{diffop}.}
%\textit{symbol} 
\begin{equation}\label{simbolo}
x_\ell:=\om\cdot \nu_\ell - \om_{j_\ell}
\end{equation}  
and a \textit{propagator}
\begin{equation}\label{prop}
{\matG_\ell(\om)=\calG_{n_\ell}(x_\ell) , }%\qquad \hbox{stesso simbolo?} }
\end{equation}
where $\calG_{n}:\RRR\to\RRR$ is the function
\begin{equation}\label{procesi}
\calG_{n}(x):=\left\{
\begin{aligned}
	&\frac{\Psi_{n}(x) %- \om_{j_\ell})
	}{x %- \om_{j_\ell}
	} ,\quad &x \ne 0 , \\ %\ne\om_{j_{\ell}}, \\
	&1,  &x= 0 , %\om_{j_{\ell}},
\end{aligned}
\right.	
\end{equation}
and $\{\Psi_n(x) \}_{n\in \ZZZ_+}$ is a partition of unity defined as follows, in terms of the frequency vector $\om$ and the
sequence $r$ which we fixed at the beginning of Subsection \ref{sec:range}.

\begin{figure}[ht] 
\vskip.2truecm 
\centering 
\ins{370pt}{-132pt}{$x$} 
\ins{312pt}{-135pt}{$\displaystyle{\rho_{0}}$} 
\ins{210pt}{-130pt}{$\displaystyle{\frac{\rho_{0}}{2}}$} 
\ins{188pt}{-135pt}{$\displaystyle{\rho_{1}}$} 
\ins{140pt}{-130pt}{$\displaystyle{\frac{\rho_{1}}{2}}$} 
%\ins{120pt}{-130pt}{$\displaystyle{\frac{\al_{m_{2}}}{8}}$} 
\ins{104pt}{-130pt}{$\displaystyle{\frac{\rho_{2}}{2}}$} 
\ins{130pt}{-10pt}{$\Psi_{2}(x)$} 
\ins{210pt}{-10pt}{$\Psi_{1}(x)$} 
\ins{330pt}{-10pt}{$\Psi_{0}(x)$} 
%\ins{098pt}{-10pt}{$\chi_{0}(x)$} 
\includegraphics[width=4in]{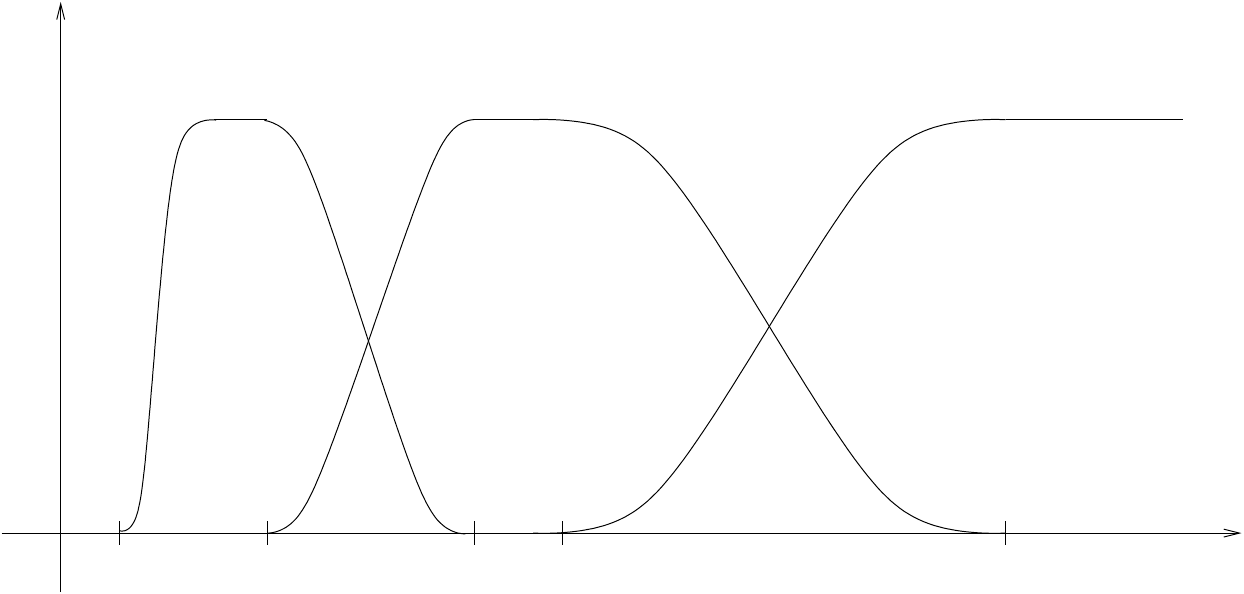} 
\vskip.2truecm 
\caption{Some of the $C^{\io}$ functions $\Psi_{n}(x)$ for $x \ge0$; here 
$\rho_{n}=\frac{1}{8}\be({m_n})$.} 
\label{fig:1} 
\end{figure}

Roughly we would like to say that a line $\ell$ has scale $m$ if $|x_\ell|
\approx \be_{\om}^{(0)}(r_m)$; however, since the sequence $\{\be_{\om}^{(0)}(r_m)\}_{m\ge0}$
is only non-increasing, for the scales to be clearly identified
one must take a decreasing subsequence $\{\be_{\om}^{(0)}(r_{m_n})\}_{n\ge0}$ whose elements are well separated,
and say that $\ell$ has scale $n$ if $|x_\ell|
\approx \be_{\om}^{(0)}(r_{m_n})$. It would be natural to consider a sharp
partition through step functions with support in 
$[\be_{\om}^{(0)}(r_{m_n}),\be_{\om}^{(0)}(r_{m_{n-1}}))$, but later on we shall 
need to take derivatives of the propagators, so it is more convenient to
use a smooth partition.
Thus, following ref.~\cite{CG2}, we introduce the sequence $\{m_n\}_{n\ge0}$ such that
%$m_0=0$ and $m_{n+1}=m_n+p_n+1$
%$m_0=0$ and $m_{n+1}=m_n+p_n$
$m_0=0$ and
%$\be_\om(r_{m_0})<1$
and $m_{n+1}=m_n+p_n$
where
\[
p_n := \min \Bigl\{q\in\NNN\;:\; 
%\be_{\om}(r_{m_n}) > 2 \be_{\om}(r_{m_n+q})\}.
\be_{\om}^{(0)}(r_{m_n+q}) < \frac{1}{2} \be_{\om}^{(0)}(r_{m_n}) \Bigr\}.
\]
In other words $\be_{\om}^{(0)}(r_{m_{n}})$ is the largest element in $\{\be_{\om}^{(0)}(r_{m})\}_{m\ge0}$
smaller than  $\frac{1}{2}\be_{\om}^{(0)}(r_{m_{n-1}})$, i.e.~such that
\begin{equation}\label{ie}
\be_{\om}^{(0)}(r_{m_n}) < \frac{1}{2}\be_{\om}^{(0)}(r_{m_{n-1}}) \le \be_{\om}^{(0)}(r_{m_n -1}) .
\end{equation}
Next, let $\chi:\RRR\to\RRR$ be a $C^{\io}$ even function, non-increasing for $x\ge0$,
%and non-decreasing for $x<0$, 
such that 
\begin{equation}\label{eq:3.10} 
\chi(x)=\left\{ 
\begin{aligned} 
&1,\qquad |x| \le 1/2, \\ 
&0,\qquad |x| \ge 1 , 
\end{aligned}\right.
\end{equation} 
and set 
\begin{subequations} \label{lachi}
\begin{align}
\chi_{-1}(x) & :=1 ,
\label{lachia} \\
\chi_{n}(x) & :=\chi \Bigl( \frac{8}{\be_\om^{(0)}({r_{m_n}})}x \Bigr) , \qquad n\ge 0.
\label{lachib}
\end{align}
\end{subequations}
Finally set $\Psi_{n}(x) :=\chi_{n-1}(x)-\chi_{n}(x)$ for $n\ge 0$ (see Figure \ref{fig:1}).

%%%%%%%%%%%%%%%%%%%%%%%%%%%%%%%%%%%%%%%%%%%%%%%%%%%%%%%%%%%%%%%%%%%%%%%%%%
% Figure 1 
%%%%%%%%%%%%%%%%%%%%%%%%%%%%%%%%%%%%%%%%%%%%%%%%%%%%%%%%%%%%%%%%%%%%%%%%%%
%%%%%%%%%%%%%%%%%%%%%%%%%%%%%%%%%%%%%%%%%%%%%%%%%%%%%%%%%%%%%%%%%%%%%%%%%%

%%%%%%%%%%%%%%%%%%%%%%%%%%%%%%%%%%%%%%%%%%%%%%%%%%%%%%%%%%%%%%%%%%%%%%%%%%
\begin{rmk} \label{infinito}
\emph{
To simplify the notation, in the following we shorten $\be(m):=\be_\om^{(0)}(r_m)$,
so that \eqref{ie} reads $2\be(m_n) < \be(m_{n-1}) \le 2\be(m_{n}-1)$.
It is also convenient to define $\be(m_{-1})=+\io$.
}
\end{rmk}
%%%%%%%%%%%%%%%%%%%%%%%%%%%%%%%%%%%%%%%%%%%%%%%%%%%%%%%%%%%%%%%%%%%%%%%%%%

%%%%%%%%%%%%%%%%%%%%%%%%%%%%%%%%%%%%%%%%%%%%%%%%%%%%%%%%%%%%%%%%%%%%%%%%%%
\begin{rmk} \label{otto}
\emph{
Of course the factor $1/2$ in both \eqref{ie} and \eqref{eq:3.10} is arbitrary, and indeed we could instead use any 
number $a<1$. In turn the factor $8$ in \eqref{lachib} could be replaced with any number $b$ such that $ab>2$, a property
that will be used in Lemma \ref{ovvio2}.
However this would only make the notations unnecessarily more involved, so we decided to fix the values of $a$ and $b$ once and for all.
}
\end{rmk}
%%%%%%%%%%%%%%%%%%%%%%%%%%%%%%%%%%%%%%%%%%%%%%%%%%%%%%%%%%%%%%%%%%%%%%%%%%

%%%%%%%%%%%%%%%%%%%%%%%%%%%%%%%%%%%%%%%%%%%%%%%%%%%%%%%%%%%%%%%%%%%%%%%%%%
\begin{rmk} %[\textbf{Small divisor}] 
\label{smalldivisor}
\emph{%remark
The reason why we call \emph{small divisor} the symbol $x_\ell=\om\cdot\nu_\ell-\om_{j_\ell}$, appearing in the denominator
of the propagator \eqref{prop}, is that all convergence issues are due to the fact that it can be arbitrarily small.
}%remark
\end{rmk}
%%%%%%%%%%%%%%%%%%%%%%%%%%%%%%%%%%%%%%%%%%%%%%%%%%%%%%%%%%%%%%%%%%%%%%%%%%

%%%%%%%%%%%%%%%%%%%%%%%%%%%%%%%%%%%%%%%%%%%%%%%%%%%%%%%%%%%%%%%%%%%%%%%%%%
\begin{defi}[\textbf{Value of an unexpanded tree}]\label{defvalue1}
{For any unexpanded tree  $\vartheta\in \Theta$
we define the \emph{value} of $\vartheta$ as
\begin{equation}\label{val}
\VV(\vartheta;c,\om,\h) := \Big(\prod_{\la\in \Lambda(\vartheta)}\LL_\la(c) \Big) \Big(\prod_{v\in N(\vartheta)} \calF_v(\h)\Big)
\Big(\prod_{\ell\in L(\vartheta)} \matG_\ell(\om) \Big) ,
\end{equation}
with $\LL_\la(c)$, $\calF_v(\h)$ and $\matG_\ell(\om)$ defined according to \eqref{foglie}, \eqref{nodi} and \eqref{prop}, respectively.}
\end{defi}
%%%%%%%%%%%%%%%%%%%%%%%%%%%%%%%%%%%%%%%%%%%%%%%%%%%%%%%%%%%%%%%%%%%%%%%%%%

%%%%%%%%%%%%%%%%%%%%%%%%%%%%%%%%%%%%%%%%%%%%%%%%%%%%%%%%%%%%%%%%%%%%%%%%%%
\begin{rmk}\label{freddo-rmk}
\emph{%remark
For any tree $\vartheta\in\Theta$, once the component and momentum labels have been assigned,
the number of scales that can be associated with a line $\ell\in L(\vartheta)$ in such a way that $\matG_\ell(\om)\ne0$ and hence
$\VV(\vartheta;c,\om,\h)\ne0$ is at most two.
}%remark
\end{rmk}
%%%%%%%%%%%%%%%%%%%%%%%%%%%%%%%%%%%%%%%%%%%%%%%%%%%%%%%%%%%%%%%%%%%%%%%%%%

%%%%%%%%%%%%%%%%%%%%%%%%%%%%%%%%%%%%%%%%%%%%%%%%%%%%%%%%%%%%%%%%%%%%%%%%%%
\begin{rmk}\label{lescale}
\emph{%remark
If a tree $\vartheta\in\Theta$ is such that $\VV(\vartheta;c,\om,\h)\ne0$, then
 for all $\ell\in L(\vartheta)$ with $n_\ell\ge 0$ one has
\begin{equation}\label{taglia}
%\frac{1}{16}\be_{\om}(r_{m_{n_\ell}}) \le |\om\cdot\nu_\ell - \om_{j_\ell}| < \frac{1}{8}\be_{\om}(r_{m_{n_\ell-1}}),
\frac{1}{16}\be({m_{n_\ell}}) < |x_\ell| < \frac{1}{8}\be({m_{n_\ell-1}}) .
\end{equation}
%
%where $\be(m_{-1})$ has to be interpreted as $+\io$.
This means that, whenever $\VV(\vartheta;c,\om,\h)\ne0$, the size of the small divisors are clearly identifiable through the corresponding scales.
}%remark
\end{rmk}
%%%%%%%%%%%%%%%%%%%%%%%%%%%%%%%%%%%%%%%%%%%%%%%%%%%%%%%%%%%%%%%%%%%%%%%%%%

%%%%%%%%%%%%%%%%%%%%%%%%%%%%%%%%%%%%%%%%%%%%%%%%%%%%%%%%%%%%%%%%%%%%%%%%%%
\begin{lemma}\label{boh}
For any $k\ge1$, $j\in\ZZZ$, $\nu\in\ZZZ^{\ZZZ}_{f}$ such that $\nu\ne\gote_j$, one has formally
\begin{equation}\label{uk}
w^{(k)}_{j,\nu}(c,\om,\h) = \sum_{\vartheta\in \Theta^{(k)}_{j,\nu,+}}\VV(\vartheta;c,\om,\h),\qquad
\ol{w}^{(k)}_{j,\nu}(c,\om,\h) = \sum_{\vartheta\in \Theta^{(k)}_{j,\nu,-}}\VV(\vartheta;c,\om,\h).
\end{equation}
\end{lemma}
%%%%%%%%%%%%%%%%%%%%%%%%%%%%%%%%%%%%%%%%%%%%%%%%%%%%%%%%%%%%%%%%%%%%%%%%%%

%%%%%%%%%%%%%%%%%%%%%%%%%%%%%%%%%%%%%%%%%%%%%%%%%%%%%%%%%%%%%%%%%%%%%%%%%%
\prova
By construction.
\EP
%%%%%%%%%%%%%%%%%%%%%%%%%%%%%%%%%%%%%%%%%%%%%%%%%%%%%%%%%%%%%%%%%%%%%%%%%%

%%%%%%%%%%%%%%%%%%%%%%%%%%%%%%%%%%%%%%%%%%%%%%%%%%%%%%%%%%%%%%%%%%%%%%%%%%
\begin{rmk} \label{noneradetto}
\emph{
Thanks to Lemma \ref{boh}, Lemma \ref{nonsocosa} yields the property \eqref{pigreco} in Remark \ref{momento}.
 }
\end{rmk}
%%%%%%%%%%%%%%%%%%%%%%%%%%%%%%%%%%%%%%%%%%%%%%%%%%%%%%%%%%%%%%%%%%%%%%%%%%

%%%%%%%%%%%%%%%%%%%%%%%%%%%%%%%%%%%%%%%%%%%%%%%%%%%%%%%%%%%%%%%%%%%%%%%%%%
\begin{rmk} \label{ultimo}
\emph{
We can also represent formally the r.h.s of \eqref{bif} as sum over trees, in a similar way.
Indeed, setting
\begin{equation}\label{gi}
G_{j,\s}^{(k)}(c,\om,\h) := \sum_{\vartheta\in \Theta^{(k)}_{j,\gote_j,\s}}\VV(\vartheta;c,\om,\h),
\end{equation}
we see that \eqref{bif} can be written order by order as $G^{(k)}_{j,\s}(c,\om,\h)=0$.
}
\end{rmk}
%%%%%%%%%%%%%%%%%%%%%%%%%%%%%%%%%%%%%%%%%%%%%%%%%%%%%%%%%%%%%%%%%%%%%%%%%%

%%%%%%%%%%%%%%%%%%%%%%%%%%%%%%%%%%%%%%%%%%%%%%%%%%%%%%%%%%%%%%%%%%%%%%%%%%
\begin{rmk}\label{formalseries}
\emph{
Both the identity \eqref{uk} and the definition \eqref{gi} are only formal because 
it is not \emph{a priori} obvious that the sums in \eqref{uk} and in \eqref{gi} are well-defined.
However, we shall verify later, as a consequence of the bounds we shall prove,
%see Lemma \ref{welldefined} below)
that, in fact, both $w^{(k)}_{j,\nu}(c,\om,\h)$ and
$G_{j,\s}^{(k)}(c,\om,\h)$ are well-defined for all $k\ge 1$ and all $(j,\nu) \in\ZZZ \times \ZZZ^{\ZZZ}_{f}$, and, furthermore,
the coefficients $w^{(k)}_{j,\nu}(c,\om,\h)$ are summable over $j$ and $\nu$. This yields that the almost-periodic solution \eqref{fw} 
to \eqref{ran} is well-defined to all orders as a formal power series. What really requires some work is to prove that the series converges.
}
\end{rmk}
%%%%%%%%%%%%%%%%%%%%%%%%%%%%%%%%%%%%%%%%%%%%%%%%%%%%%%%%%%%%%%%%%%%%%%%%%%

%%%%%%%%%%%%%%%%%%%%%%%%%%%%%%%%%%%%%%%%%%%%%%%%%%%%%%%%%%%%%%%%%%%%%%%%%%
\subsection{The full equation: expanded trees}
\label{labif}
%%%%%%%%%%%%%%%%%%%%%%%%%%%%%%%%%%%%%%%%%%%%%%%%%%%%%%%%%%%%%%%%%%%%%%%%%%

Now we aim to take into account also \eqref{bif}. If we explicitly solve it for $k=1$ we find
\begin{equation}\label{bif1}
\h_j^{(1)}(c,\om) = - \frac{1}{c_j}
\sum_{\substack{\gote_{j_1}-\gote_{j_2} + \gote_{j_3} - \gote_{j_4} + \gote_{j_5} =\gote_j \\ j_1-j_2+j_3 - j_4 +j_5 = j }}
c_{j_1}\ol{c}_{j_2}c_{j_3}\ol{c}_{j_4}c_{j_5} .
\end{equation}

%%%%%%%%%%%%%%%%%%%%%%%%%%%%%%%%%%%%%%%%%%%%%%%%%%%%%%%%%%%%%%%%%%%%%%%%%%
\begin{rmk} \label{nosingo}
\emph{
The singularity at $c_j=0$ in \eqref{bif1} is only apparent:  from the constraints in the sum
we see that at least one among $j_i$, $i=1,\ldots,5$, equals $j$ and hence
there is at least a factor $c_j$ which compensates for the denominator $c_j$ (see also Lemma \ref{cijei} and Remark \ref{singolare}). 
}
\end{rmk}
%%%%%%%%%%%%%%%%%%%%%%%%%%%%%%%%%%%%%%%%%%%%%%%%%%%%%%%%%%%%%%%%%%%%%%%%%%

%%%%%%%%%%%%%%%%%%%%%%%%%%%%%%%%%%%%%%%%%%%%%%%%%%%%%%%%%%%%%%%%%%%%%%%%%%
\begin{rmk} \label{riserve?}
\emph{%remark
In fact $\h^{(k)}(c,\om)$ depends on $\om$ only from $k=2$ on.
}%remark
\end{rmk}
%%%%%%%%%%%%%%%%%%%%%%%%%%%%%%%%%%%%%%%%%%%%%%%%%%%%%%%%%%%%%%%%%%%%%%%%%%

%%%%%%%%%%%%%%%%%%%%%%%%%%%%%%%%%%%%%%%%%%%%%%%%%%%%%%%%%%%%%%%%%%%%%%%%%%
\begin{rmk}
\emph{%remark
A trivial computation gives
\begin{equation}\label{eta1bis}
	\h_j^{(1)} (c,\om)=  - |c_j|^4 - 6 \sum_{\substack{ i \in \ZZZ \\ i\neq j}} |c_i|^2 |c_j|^2   - 3 \sum_{\substack{ i \in \ZZZ \\ i\neq j}} |c_i|^4  
	- 12 \sum_{\substack{i,k \in \ZZZ \\ i,k \ne j , i\ne k}} |c_i|^2|c_k|^2 ,
\vspace{-.2cm}
\end{equation}
which shows that $\h_j^{(1)}(c,\om)$ is real for all $j\in\ZZZ$. 
}%remark
\end{rmk}
%%%%%%%%%%%%%%%%%%%%%%%%%%%%%%%%%%%%%%%%%%%%%%%%%%%%%%%%%%%%%%%%%%%%%%%%%%

For $k\ge2$, we solve recursively both equations \eqref{sub} simultaneously. This means that we look
for a solution $u=u(x,\f;c,\om,\e)$ and a counterterm $\eta=\eta(c,\om,\e)$ which, inserted into \eqref{sub}, make the equations
to be satisfied order by order. Thus, we obtain 
\begin{subequations}\label{bif2}
\begin{align}
& u_{j,\nu}^{(k)} = \frac{1}{(\om\cdot\nu-\om_j)}\Bigg[
\sum_{k_1+k_2=k} \h^{(k_1)}_j  u^{(k_2)}_{j,\nu} + 
\!\!\!\!\!\!\!\!\!\!\!\!
\sum_{\substack{ \nu_1 -\nu_2 + \nu_3-\nu_4 +\nu_5 =\nu \\ j_1-j_2+j_3 - j_4 +j_5 = j  \\
k_1+k_2+k_3+k_4+k_5=k}}
\!\!\!\!\!\!\!\!\!\!\!\!
u_{j_1,\nu_1}^{(k_1)}\ol{u}_{j_2,\nu_2}^{(k_2)}u_{j_3,\nu_3}^{(k_3)}\ol{u}_{j_4,\nu_4}^{(k_4)}u_{j_5,\nu_5}^{(k_5)}  \Bigg] , 
\label{bif2a} \\
& \h_j^{(k)} = -\frac{1}{c_j}
\sum_{\substack{ \nu_1 -\nu_2 + \nu_3-\nu_4 +\nu_5 =\gote_j \\ j_1-j_2+j_3 - j_4 +j_5 = j  \\
k_1+k_2+k_3+k_4+k_5=k-1}}
u_{j_1,\nu_1}^{(k_1)}\ol{u}_{j_2,\nu_2}^{(k_2)}u_{j_3,\nu_3}^{(k_3)}\ol{u}_{j_4,\nu_4}^{(k_4)}u_{j_5,\nu_5}^{(k_5)}.
\phantom{\sum^{\sum^{\sum^{\sum}}}}
\label{bif2b}
\end{align}
\end{subequations}

%%%%%%%%%%%%%%%%%%%%%%%%%%%%%%%%%%%%%%%%%%%%%%%%%%%%%%%%%%%%%%%%%%%%%%%%%%
\begin{rmk} \label{notationforeta}
\emph{
In \eqref{bif2} as well as in \eqref{rank}, we are not writing the dependence of the coefficients on their arguments $c,\om$.
However, it is important to keep in mind that $u$ does not depend on $\eta$ and
now $\h$ itself is part of the solution and depends on $c,\om$, since it is no longer
the free parameter $\h$ appearing in the series \eqref{fw}. Of course, if we are able to solve \eqref{sub}, then
we deduce that $u(x,\f;c,\om,\e)=w(x,\f;c,\om,\eta(c,\om,\e),\e)$.
}
\end{rmk}
%%%%%%%%%%%%%%%%%%%%%%%%%%%%%%%%%%%%%%%%%%%%%%%%%%%%%%%%%%%%%%%%%%%%%%%%%%

\begin{figure}[H]
\vspace{-.2cm}
\centering
%\null
%\hspace{-.6cm}
\ins{012pt}{-069.5pt}{$\vartheta=$}
\ins{261pt}{-070.5pt}{$\vartheta' =$}
\ins{109pt}{-082.5pt}{$v_1$}
\ins{178pt}{-051pt}{$v_3$}
\ins{360pt}{-084pt}{$v_1$}
%\ins{378pt}{-044pt}{$v_{in}$}
%\ins{346pt}{-058pt}{$\ell^{in}_\TT$}
\subfigure{\includegraphics*[width=3.0in]{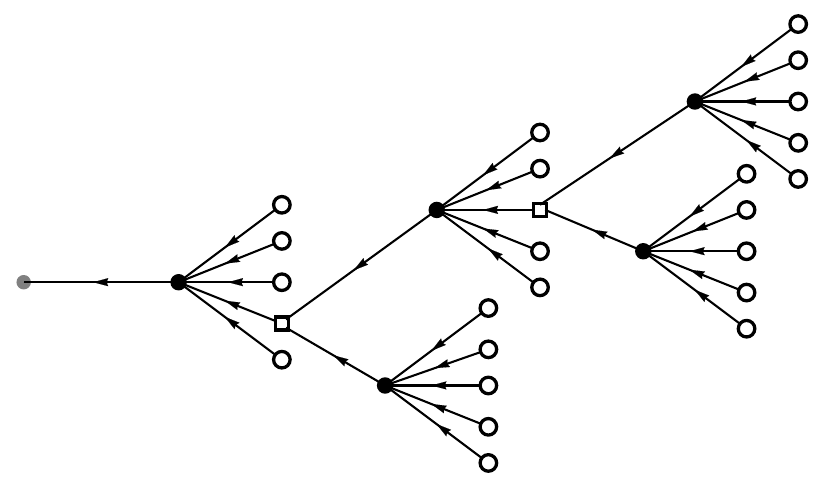}}
\hspace{1cm}
\subfigure{\includegraphics*[width=2.0in]{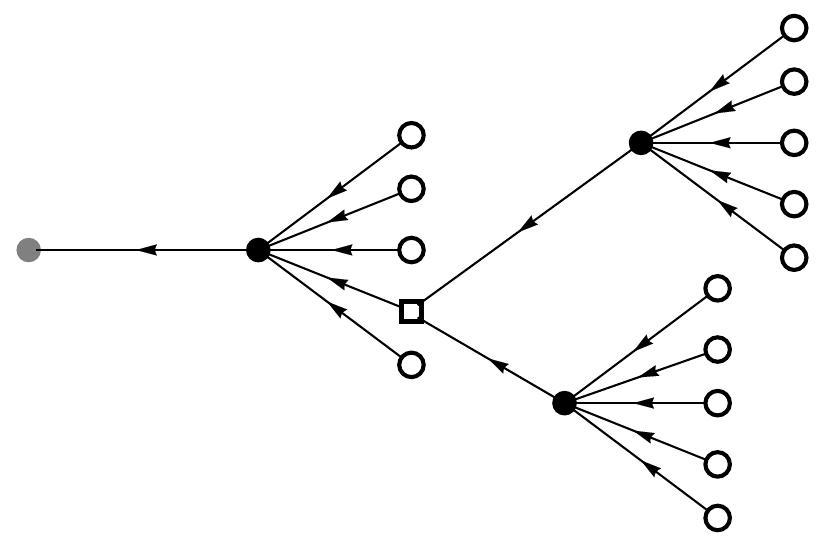}}
\caption{\small Two trees $\vartheta,\vartheta''\in\gotT$ (the labels are not shown).}
\label{fig-expanded}
\end{figure}

\begin{figure}[H]
\vspace{-.2cm}
\centering
%\null
%\hspace{-.6cm}
\ins{002pt}{-020pt}{$\vartheta_1\!=$}
\ins{098pt}{-031pt}{$v_1$}
\ins{376pt}{-034pt}{$v_3$}
\ins{166pt}{-068.5pt}{$\vartheta_2\!=$}
\ins{275pt}{-065.5pt}{$\vartheta_3\!=$}
%\ins{378pt}{-044pt}{$v_{in}$}
%\ins{346pt}{-058pt}{$\ell^{in}_\TT$}
\subfigure{\includegraphics*[width=1.8in]{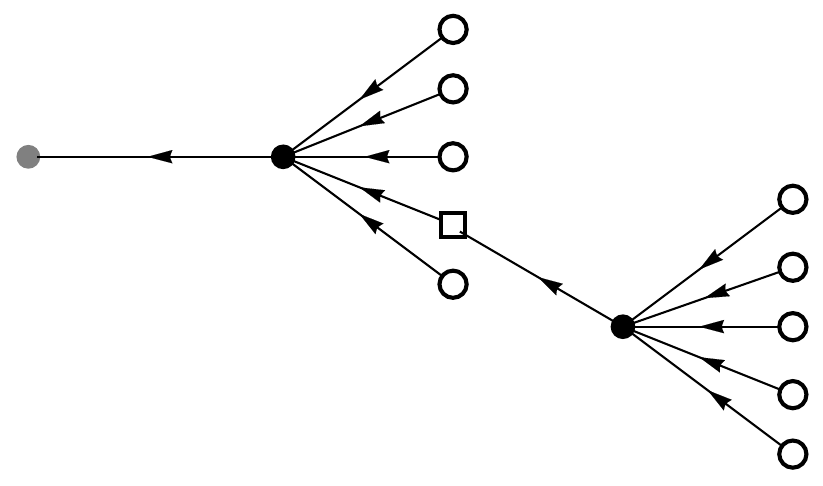}}
\hspace{1cm}
\subfigure{\includegraphics*[width=1.0in]{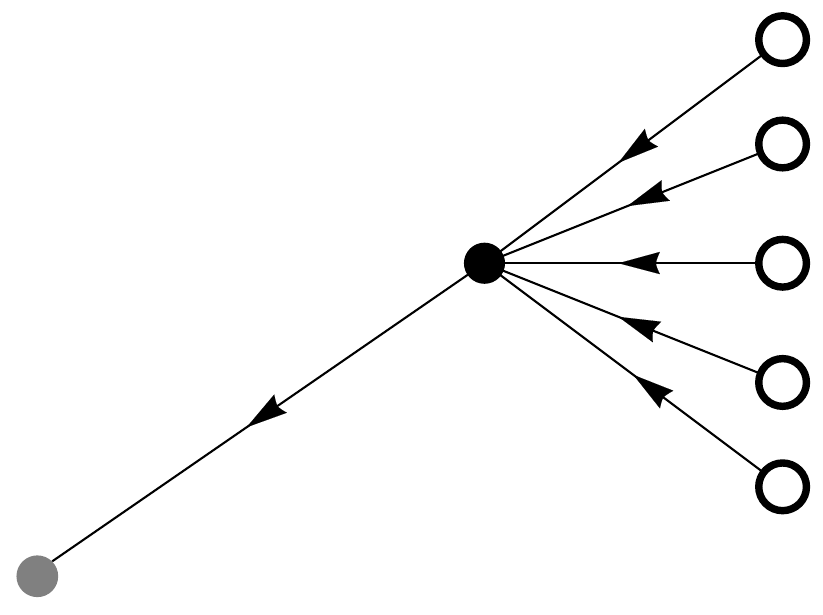}}
\hspace{1cm}
\subfigure{\includegraphics*[width=1.9in]{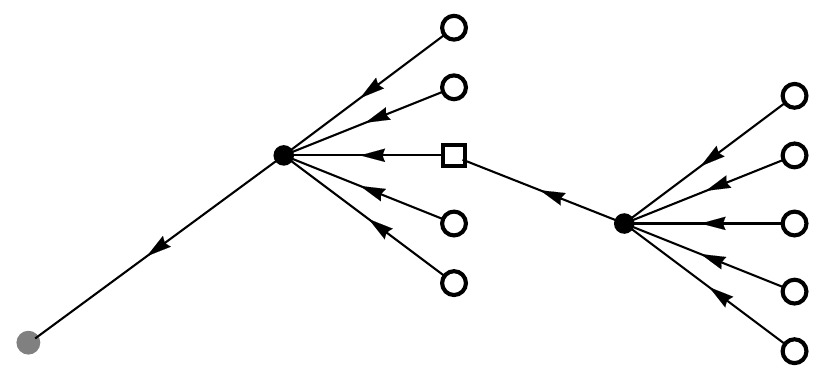}}
\caption{\small The trees $\vartheta_1,\vartheta_2,\vartheta_3\in\Theta$ used to construct the unexpanded trees in Figure \ref{fig-expanded}.}
\label{fig-unexpanded}
\end{figure}

\begin{figure}[H]
\vspace{-.2cm}
\centering
%\null
%\hspace{-.6cm}
\ins{120pt}{-094pt}{$\vartheta''=$}
%\ins{378pt}{-044pt}{$v_{in}$}
%\ins{346pt}{-058pt}{$\ell^{in}_\TT$}
\subfigure{\includegraphics*[width=2.4in]{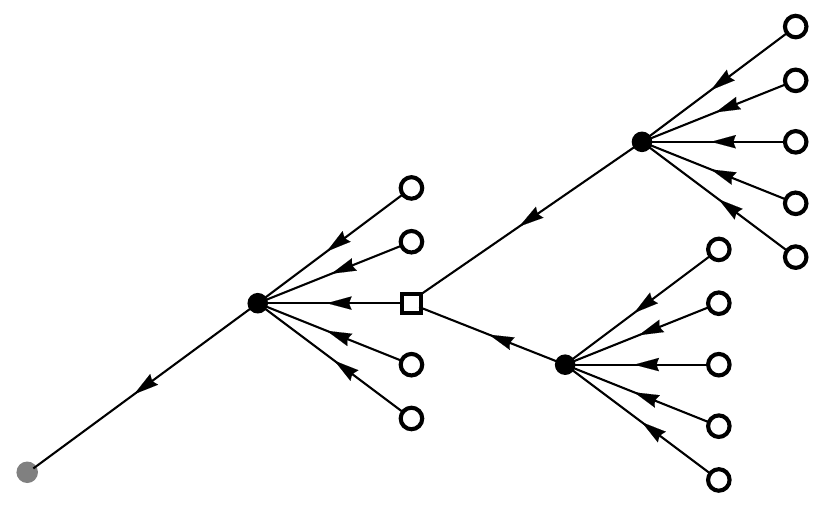}}
\caption{\small The expanded tree $\vartheta_3$ obtained by grafting the tree $\vartheta'$ to the tree $\vartheta''$ in Figure \ref{fig-unexpanded}.}
\label{fig-exp-unexp}
\end{figure}

Again we want to represent graphically the summands in \eqref{bif2}, so we introduce a second set of trees,
which are constructed starting from the set $\Theta$ of Definition \ref{thetone}
through the following recursive procedure (see Figures \ref{fig-expanded} to \ref{fig-exp-unexp} for explicit examples):
\begin{enumerate}[topsep=0ex]
\itemsep0em
\item choose any tree $\vartheta\in\Theta$ {which has not one node and one leaf only};
%{except those constituted by only one node $v$ with $s_v=1$ and a leaf};
%
\item to any node $v\in N(\vartheta)$ with $s_v=1$, graft %represented by a white square we attach
any tree in $\vartheta'\in \Theta^{(k_v)}_{j_v,\gote_{j_v},\s_v}$ %(and hence increase $s_v$ by one), 
in such a way that $\ell_{\vartheta'}$
 remains to the left of the line in $L(\vartheta)$ entering $v$,
and change the value of the order label $k_v$ into $k_v=0$ and the value of the branching label $s_v$ into $s_v=2$;
%
%\item we call $\ell_{\vartheta'}$  a $\eta$\emph{-line};
%
\item for any tree obtained as above, apply once more the procedure {to any node of the grafted trees},
and so on, until no nodes $v$ with $s_v=1$ are left, and eventually each node $v$ has either $s_v=5$ or $s_v=2$.
%
%\item finally, associate with the nodes $v$ with $s_v=2$ a new order label $k_v=0$.
%
\end{enumerate}
This leads naturally to the following definition.

\begin{defi}[\textbf{Set $\boldsymbol{\gotT}\!\!\!\!\!\boldsymbol{\gotT}$ of the expanded trees}]
\label{gotT}
Let $\gotT$ be the set of trees $\vartheta$ satisfying the following constraints:
\begin{enumerate}[topsep=0ex]
\itemsep0.0em
\item for any node $v\in N(\vartheta)$ one has $s_v=2$ or $s_v=5$;
\item for any node $v\in N(\vartheta)$ one has $k_v=0$ if $s_v=2$ and $k_v=1$ if $s_v=5$;
\item if $\ell \in L(\vartheta)$ exits a leaf $\la\in\Lambda(\vartheta)$, then $\s_\ell =\s_\la $, $j_\ell= j_\la$ and $\nu_\ell=\gote_{j_\la}$;
\item if $\ell \in L(\vartheta)$ exits a node $v\in N(\vartheta)$, then one has $\s_v=\s_\ell$ and $j_v = j_\ell$;
\item for any node $v\in N(\vartheta)$ with $s_v=2$, the upper line $\ell'_v$ entering $v$ is called a $\h$-\emph{line} and is such that
$\s_{\ell'_v}=\s_v$, $j_{\ell'_v}=j_v$ and $\nu_{\ell'_v} = \gote_{j_{\ell'_v}}$;
\item {both lines entering a node $v$ with $s_v=2$ do not exit a leaf};
\item if $\ell\in L(\vartheta)$ exits a node $v\in N(\vartheta)$ and $\ell$ is neither the root line nor a $\h$-line,
then $\nu_\ell \neq \gote_{j_\ell}$;
%
%\item for any $\ell,\ell' \in L(\vartheta)$ one writes
%\begin{enumerate}[topsep=0ex]
%\item[7.1] $\ell \st{\preceq} \ell'$ if $\ell\preceq\ell'$ and there is no $\h$-line on $\calP(\ell,\ell')$,
%\item[7.2] $v\st{\preceq} \ell$ if $\ell_v\st{\preceq}\ell$, and $v\st\preceq w$ if $\ell_v\st\preceq\ell_w$;
%\item if the line $\ell \in L(\vartheta)$ exits a node $v\in N(\vartheta)$ with $s_v=5$,
%then $\nu_\ell\ne \s_{\ell}\gote_{j_\ell}$, unless $\ell$ is the root line;
%\end{enumerate}
%
\item if $\ell \in L(\vartheta)$ exits a node $v \in N(\vartheta)$ and $\ell_1,\ldots,\ell_{s_v}$ are the lines entering $v$, then
\begin{equation}\label{lefoglie*}
	{\s_\ell 1 = \sum_{\substack{ \la \in \Lambda(\vartheta) \\\la\st{\preceq} \ell} } \s_{\la} 1,\qquad
	\s_\ell\nu_\ell =  \sum_{\substack{ \la \in \Lambda(\vartheta) \\\la\st{\preceq} \ell} }  \s_{\la} \gote_{j_\la},\qquad
		\s_\ell j_\ell= \sum_{\substack{ \la \in \Lambda(\vartheta) \\\la\st{\preceq} \ell} }  \s_{\la} {j_\la}\,,}
\end{equation}
where 
the ordering relation $\st{\preceq}$ is defined by writing
$\ell \st{\preceq} \ell'$ for any lines $\ell,\ell'\in L(\vartheta)$ such that $\ell\preceq\ell'$ and there is no $\h$-line on $\calP(\ell,\ell')$,
and writing
$v\st{\preceq} \ell$ if $\ell_v\st{\preceq}\ell$, and $v\st\preceq w$ if $\ell_v\st\preceq\ell_w$, for any $v,w\in N(\vartheta)\cup \Lambda(\vartheta)$;
%\item if the line $\ell' \in L(\vartheta)$ enters a node $v \in N(\vartheta)$ with $s_v=1$, then one has
%with $s_v=1$ and $\ell'$ is the line entering $v$, then
%
%\begin{equation}\label{impongo}
%$\s_{\ell'}=\s_v$ and $j_{\ell'}= j_v$;
%\end{equation}
%
\item for any line $\ell\in L(\vartheta)$ one has $n_\ell=-1$ if $\nu_\ell=\gote_{j_\ell}$ and $n_\ell\ge0$ otherwise.
\end{enumerate}
We call any $\vartheta\in\gotT$ an \emph{expanded tree},
and define $\gotT^{(k)}_{j,\nu,\s}$ as the set of expanded trees $\vartheta\in\gotT$ with order $k(\vartheta)=k$
such that the root line has component $j$, momentum $\nu$ and sign $\s$.
\end{defi}
%%%%%%%%%%%%%%%%%%%%%%%%%%%%%%%%%%%%%%%%%%%%%%%%%%%%%%%%%%%%%%%%%%%%%%%%%%

%%%%%%%%%%%%%%%%%%%%%%%%%%%%%%%%%%%%%%%%%%%%%%%%%%%%%%%%%%%%%%%%%%%%%%%%%%
\begin{rmk} \label{mezzormk}
\emph{
%In particular, \eqref{laseconda} shows that 
By comparing Definitions \ref{thetone} and \ref{gotT}, one realizes that any tree $\vartheta\in\gotT$ 
can be constructed recursively as follows.
Take a tree $\vartheta'\in\Theta^{(k)}_{j,\nu,\s}$
and perform the following operations for each node $v\in N(\vartheta')$ with $s_v=1$:
\begin{itemize}[topsep=0ex]
\itemsep0.0em
\item 
first {graft} a subtree $\vartheta_v \in  {\gotT^{(k_v)}_{ j_v ,\gote_{j_v},\s_v}}$ to $v$ in such a way that it remains above %on the left of
the subtree whose root line already enters $v$, and change $s_v$ into $s_v=2$;
\item
then redefine the order label of $v$ as $k_v=0$.
\end{itemize}
Eventually, a tree $\vartheta\in\gotT^{(k)}_{ j ,\nu,\s}$ is obtained.
}
\end{rmk}
%%%%%%%%%%%%%%%%%%%%%%%%%%%%%%%%%%%%%%%%%%%%%%%%%%%%%%%%%%%%%%%%%%%%%%%%%%

%%%%%%%%%%%%%%%%%%%%%%%%%%%%%%%%%%%%%%%%%%%%%%%%%%%%%%%%%%%%%%%%%%%%%%%%%%
\begin{rmk} \label{nonloso?}
\emph{%remark
The conservation law \eqref{lefoglie*} extends to trees $\vartheta\in\Theta$ as well -- and reduces to \eqref{lefoglie} in that case --
since the trees $\vartheta\in\Theta$ contain no $\h$-lines.
}%remark
\end{rmk}
%%%%%%%%%%%%%%%%%%%%%%%%%%%%%%%%%%%%%%%%%%%%%%%%%%%%%%%%%%%%%%%%%%%%%%%%%%

%%%%%%%%%%%%%%%%%%%%%%%%%%%%%%%%%%%%%%%%%%%%%%%%%%%%%%%%%%%%%%%%%%%%%%%%%%
\begin{rmk}\label{kn}
\emph{
For any tree $\vartheta\in\gotT$ one has
\[
|N(\vartheta)|\le 2 k(\vartheta)-1,
\]
as can be proved by induction on the number of nodes.
}
\end{rmk}
%%%%%%%%%%%%%%%%%%%%%%%%%%%%%%%%%%%%%%%%%%%%%%%%%%%%%%%%%%%%%%%%%%%%%%%%%%

Two examples of expanded trees are represented in Figure \ref{fig-expanded}.
To see how they can be constructed starting from trees in $\Theta$,
consider the three unexpanded trees in Figure \ref{fig-unexpanded}.

%%%%%%%%%%%%%%%%%%%%%%%%%%%%%%%%%%%%%%%%%%%%%%%%%%%%%%%%%%%%%%%%%%%%%%%%%%
% FIGURA 
%%%%%%%%%%%%%%%%%%%%%%%%%%%%%%%%%%%%%%%%%%%%%%%%%%%%%%%%%%%%%%%%%%%%%%%%%%
%%%%%%%%%%%%%%%%%%%%%%%%%%%%%%%%%%%%%%%%%%%%%%%%%%%%%%%%%%%%%%%%%%%%%%%%%%

If we graft the tree $\vartheta_2$ to the node $v_1\in N(\vartheta_1)$ with $s_{v_1}=1$ we obtain the tree $\vartheta'$; if, instead
we first graft $\vartheta_3$ to the node $v_1\in N(\vartheta_1)$ and then graft $\vartheta_2$ to the node
$v_3\in N(\vartheta_3)$ with $s_{v_3}=1$ we obtain the expanded tree $\vartheta$.
Here we are ignoring all the labels: the expanded tree inherits the labels of the unexpanded trees starting from
which it is constructed, with the only difference that, for any node $v$ with $s_v=1$, the values of $s_v$ and $k_v$
are changed into $s_v=2$ and $k_v=0$.

%%%%%%%%%%%%%%%%%%%%%%%%%%%%%%%%%%%%%%%%%%%%%%%%%%%%%%%%%%%%%%%%%%%%%%%%%%
% FIGURA 
%%%%%%%%%%%%%%%%%%%%%%%%%%%%%%%%%%%%%%%%%%%%%%%%%%%%%%%%%%%%%%%%%%%%%%%%%%
%%%%%%%%%%%%%%%%%%%%%%%%%%%%%%%%%%%%%%%%%%%%%%%%%%%%%%%%%%%%%%%%%%%%%%%%%%

Note that $\vartheta$ may also be obtained as a result of the following procedure:
we first graft $\vartheta_2$ to the node $v_3\in N(\vartheta_3)$ with $s_{v_3}=1$, so as to construct  the expanded tree $\vartheta''$ in 
Figure \ref{fig-exp-unexp}, and then graft $\vartheta''$ to the node $v_1\in N(\vartheta_1)$ wih $s_v=1$. Note also that,
for instance, just as $\vartheta''$ is obtained by grafting $\vartheta_2$ to the node $v_3$ of $\vartheta_3$, so also
$\vartheta_3$ is obtained by cutting the subtree $\vartheta_2$ of $\vartheta''$.

%%%%%%%%%%%%%%%%%%%%%%%%%%%%%%%%%%%%%%%%%%%%%%%%%%%%%%%%%%%%%%%%%%%%%%%%%%
% FIGURA 
%%%%%%%%%%%%%%%%%%%%%%%%%%%%%%%%%%%%%%%%%%%%%%%%%%%%%%%%%%%%%%%%%%%%%%%%%%
%%%%%%%%%%%%%%%%%%%%%%%%%%%%%%%%%%%%%%%%%%%%%%%%%%%%%%%%%%%%%%%%%%%%%%%%%%

The leaf factors $\LL_\la(c)$ and the propagators $\matG_\ell(\om)$ for trees in $\gotT$
are defined as in \eqref{foglie} and \eqref{prop} respectively,  while the node factors are defined as 
\begin{equation}\label{nodiker}
\starF_v = \starF_v(c) := 
\begin{cases}
\displaystyle{1, \phantom{\raisebox{-.2cm}{\Big(}}} & \quad 
s_v =5\\
-\displaystyle{\frac{1}{c^{\s_v}_{j_v}} }, & \quad 
s_v=2 .
\end{cases}
\end{equation}

%%%%%%%%%%%%%%%%%%%%%%%%%%%%%%%%%%%%%%%%%%%%%%%%%%%%%%%%%%%%%%%%%%%%%%%%%%
\begin{defi}[\textbf{Value of an expanded tree}]\label{defvalue2}
{For any expanded tree $\vartheta\in\gotT$, we define the \emph{value} of $\vartheta$ as
\begin{equation}\label{valeta}
\Val(\vartheta;c,\om):=\Big(\prod_{\la\in \Lambda(\vartheta)}\LL_\la (c)\Big) 
\Big(\prod_{v\in N(\vartheta)} \starF_v(c)\Big)
\Big(\prod_{\ell\in L(\vartheta)} \matG_\ell(\om) \Big) ,
\end{equation}
with $\LL_\la(c)$, $\calF_v^*(c)$ and $\matG_\ell(\om)$ defined according to \eqref{foglie}, \eqref{nodiker} and \eqref{prop}, respectively.}
\end{defi}
%%%%%%%%%%%%%%%%%%%%%%%%%%%%%%%%%%%%%%%%%%%%%%%%%%%%%%%%%%%%%%%%%%%%%%%%%%

Finally set
\begin{subequations}\label{veri}
\begin{align}
u^{(k)}_{j,\nu} (c,\om) 
:= \sum_{\vartheta  \in \gotT^{(k)}_{ j ,\nu,+}} &\Val (\vartheta;c,\om),
\qquad
\ol{u}^{(k)}_{j,\nu} (c,\om):= \sum_{\vartheta  \in \gotT^{(k)}_{ j ,\nu,-}} \Val (\vartheta;c,\om),\label{ukvero}\\
&(\h_j^{(k)}(c,\om))^{\s}:= -\frac{1}{c_j^\s}\sum_{\vartheta  \in \gotT^{(k)}_{ j ,\gote_j,\s}} \Val  (\vartheta;c,\om) \label{etak} \,.
\end{align}
\end{subequations}

We conclude this section by showing that the series \eqref{formale+fout}, with the coefficients given by \eqref{veri},
formally solve \eqref{bif2} and hence \eqref{sub}.
For any $\vartheta\in\Theta\cup\gotT$, set
\begin{equation} \label{Lambdastar}
\Lambda^*(\vartheta) := \{ \la\in \Lambda(\vartheta) : \la\st{\preceq}r_\vartheta\} ,
\qquad
\Lambda^*_{j,\s}(\vartheta) := \{ \la\in \Lambda_{j,\s}(\vartheta) : \la\st{\preceq}r_\vartheta\} ,
\end{equation}
with the notation of Definition \ref{gotT}.
Note that $\Lambda^*(\vartheta)=\Lambda(\vartheta)$ if $\vartheta\in\Theta$.

%%%%%%%%%%%%%%%%%%%%%%%%%%%%%%%%%%%%%%%%%%%%%%%%%%%%%%%%%%%%%%%%%%%%%%%%%%
\begin{lemma}\label{cijei}
For all $k\ge1$, $j\in\ZZZ$ and $\s\in\{\pm\}$ and for any 
$\vartheta\in \Theta^{(k)}_{j,\gote_j,\s}\cup %\breve\gotT^{(k)}_{j,\gote_j,\s}\cup 
\gotT^{(k)}_{j,\gote_j,\s}$ one has
\[
{|\Lambda^*_{j',\s}(\vartheta)| = |\Lambda^*_{j',-\s}(\vartheta)| + \de_{j' j} . }
\]
\end{lemma}
%%%%%%%%%%%%%%%%%%%%%%%%%%%%%%%%%%%%%%%%%%%%%%%%%%%%%%%%%%%%%%%%%%%%%%%%%%

%%%%%%%%%%%%%%%%%%%%%%%%%%%%%%%%%%%%%%%%%%%%%%%%%%%%%%%%%%%%%%%%%%%%%%%%%%
\prova
Let $\ell$ be the root line of $\vartheta$.
The conservation law \eqref{lefoglie*}, with $\nu_{\ell}=\gote_j$, yields
\begin{equation} \nonumber
\begin{aligned}
\sigma \de_{j j'} & =  (\sigma \gote_j)_{j'} = \sum_{\la \in \Lambda^*(\vartheta)} \sigma_\la \, \de_{j' j_\la }  \\
& =
\sigma \sum_{\substack{ \la \in \Lambda^*(\vartheta) \\ \sigma_\la = \sigma}}  \de_{j' j_\la } -
\sigma \sum_{\substack{ \la \in \Lambda^*(\vartheta) \\ \sigma_\la = - \sigma}}  \de_{j' j_\la } 
=
\sigma \Bigl( |\Lambda^*_{j',\sigma}(\vartheta)| - |\Lambda^*_{j',-\sigma}(\vartheta)| \Bigr) ,
\end{aligned}
\end{equation}
which, multiplied by $\sigma$, provides the result.
\EP
%%%%%%%%%%%%%%%%%%%%%%%%%%%%%%%%%%%%%%%%%%%%%%%%%%%%%%%%%%%%%%%%%%%%%%%%%%

%%%%%%%%%%%%%%%%%%%%%%%%%%%%%%%%%%%%%%%%%%%%%%%%%%%%%%%%%%%%%%%%%%%%%%%%%%
\begin{rmk}\label{singolare}
\emph{%remark
%By Lemma \ref{cijei} and Remark \ref{stracon} it follows that $(\h_j^{(k)}(c))^{+}=(\h_j^{(k)}(c))^{-}$, i.e. $\h_j^{(k)}(c)\in\RRR$,
%\[
%\sum_{\vartheta  \in \gotT^{(k)}_{ j ,\gote_j,+}} \Val (\vartheta;c) = \sum_{\vartheta  \in \gotT^{(k)}_{ j ,\gote_j,-}} \Val (\vartheta;c) ,
%\]
Thanks to Lemma \ref{cijei},
$(\h_j^{(k)}(c,\om))^{\s}$, for both $\s=\pm$,
is a formal polynomial in $|c_{j'}|^2$ for all $j'\in\ZZZ$, so that it is
a function of the sequence $|c|^2:=\{|c_j|^2\}_{j\in\ZZZ}$ rather than $c$.
%Analogously  $(\h_j^{(k)}(c,\om,n))^{+}=(\h_j^{(k)}(c,\om,n))^{-}$, i.e. $\h_j^{(k)}(c,\om,n)\in\RRR$.
In particular, one has $(\h_j^{(k)}(c,\om))^{-}=(\h_j^{(k)}(c,\om))^{+}$, and hence $\h_j^{(k)}(c,\om):=(\h_j^{(k)}(c,\om))^{+}$ is real.
}%remark
\end{rmk}
%%%%%%%%%%%%%%%%%%%%%%%%%%%%%%%%%%%%%%%%%%%%%%%%%%%%%%%%%%%%%%%%%%%%%%%%%%

%%%%%%%%%%%%%%%%%%%%%%%%%%%%%%%%%%%%%%%%%%%%%%%%%%%%%%%%%%%%%%%%%%%%%%%%%%
\begin{rmk}\label{calcolo}
\emph{%remark
By Remark \ref{mezzormk}, the two sums \eqref{ukvero} equal the sums  \eqref{uk}
as soon as, for any $\vartheta\in\Theta^{(k)}_{j,\nu,\s}$ and any node $v\in N(\vartheta)$ with $s_v=1$,
one replaces the free parameter $\h^{(k_v)}_{j_v}$
with $\h^{(k_v)}_{j_v}(c,\om)$ in \eqref{etak}. In other words one has
$u^{(k)}_{j,\nu}(c,\om)=w^{(k)}_{j,\nu}(c,\om,\h(c,\om))$.
}%remark
\end{rmk}
%%%%%%%%%%%%%%%%%%%%%%%%%%%%%%%%%%%%%%%%%%%%%%%%%%%%%%%%%%%%%%%%%%%%%%%%%%

%%%%%%%%%%%%%%%%%%%%%%%%%%%%%%%%%%%%%%%%%%%%%%%%%%%%%%%%%%%%%%%%%%%%%%%%%%
\begin{lemma}\label{secondome}
The function \eqref{fout}, with the coefficients $u^{(k)}_{j,\nu}(c,\om)$ given in \eqref{ukvero},
and the counterterm $\h(c,\om,\e)$, with $\h_j(c,\om,\e)$ as in \eqref{formale} and $\h_j^{(k)}(c,\om)$ given in \eqref{etak},
solve \eqref{sub} order by order.
In particular, one has
\[
G^{(k)}_{j,\s}(c,\om,\h(c,\om))\equiv0 , 
\]
where $ G^{(k)}_{j,\s}(c,\om,\h)$ is defined in \eqref{gi}.
\end{lemma}
%%%%%%%%%%%%%%%%%%%%%%%%%%%%%%%%%%%%%%%%%%%%%%%%%%%%%%%%%%%%%%%%%%%%%%%%%%

%%%%%%%%%%%%%%%%%%%%%%%%%%%%%%%%%%%%%%%%%%%%%%%%%%%%%%%%%%%%%%%%%%%%%%%%%%
\proof
By construction.
\EP
%%%%%%%%%%%%%%%%%%%%%%%%%%%%%%%%%%%%%%%%%%%%%%%%%%%%%%%%%%%%%%%%%%%%%%%%%%

%%%%%%%%%%%%%%%%%%%%%%%%%%%%%%%%%%%%%%%%%%%%%%%%%%%%%%%%%%%%%%%%%%%%%%%%%%
%%%%%%%%%%%%%%%%%%%%%%%%%%%%%%%%%%%%%%%%%%%%%%%%%%%%%%%%%%%%%%%%%%%%%%%%%%
\section{Resonances and small divisors}
\label{sec}
\zerarcounters
%%%%%%%%%%%%%%%%%%%%%%%%%%%%%%%%%%%%%%%%%%%%%%%%%%%%%%%%%%%%%%%%%%%%%%%%%%
%%%%%%%%%%%%%%%%%%%%%%%%%%%%%%%%%%%%%%%%%%%%%%%%%%%%%%%%%%%%%%%%%%%%%%%%%%

%Given a tree $\vartheta\in\gotT$, %$\vartheta\in\Theta\cup\gotT$, 
We want to identify the contributions that are a ``real obstruction'' to the convergence of  the series \eqref{veri}.
As usual in KAM theory, all difficulties come from the small divisors, which in our representation appear in the propagators.
In fact, the problem arises from the possible ``accumulation'' of small divisors -- essentially, as we shall see,
from the presence of subgraphs $T$ of the trees with one entering line and one exiting line such that
their propagators are equal in absolute value, while the propagators of the lines in  $T$ are much smaller.
%Roughly the small divisors inside $T$ are larger, thus less dangerous.
We shall discuss in Section \ref{stimazze} why the values of the trees
where the small divisors accumulate are difficult to handle and require a more careful analysis.
Before doing that, we introduce some notation with the aim of identifying the ``dangerous'' trees.

%%%%%%%%%%%%%%%%%%%%%%%%%%%%%%%%%%%%%%%%%%%%%%%%%%%%%%%%%%%%%%%%%%%%%%%%%%
\subsection{{Trimmed trees}}
%%%%%%%%%%%%%%%%%%%%%%%%%%%%%%%%%%%%%%%%%%%%%%%%%%%%%%%%%%%%%%%%%%%%%%%%%%

%%%%%%%%%%%%%%%%%%%%%%%%%%%%%%%%%%%%%%%%%%%%%%%%%%%%%%%%%%%%%%%%%%%%%%%%%%
\begin{defi}[\textbf{Set $\breve{\gotT}$ of the trimmed trees}]\label{potati}
Let $\breve\gotT$ the set of trees defined as $\Theta$ with the only difference that $k_v=0$ for any node $v$ such that $s_v=1$.
The trees in $\breve\gotT$ are called \emph{trimmed trees},
and $\breve\gotT^{(k)}_{j,\nu,\s}$ denotes the set of trimmed trees $\vartheta\in\breve\gotT$ with $|N(\vartheta)|=k$
% order $k(\vartheta)=k$
such that the root line has component $j$, momentum $\nu$ and sign $\s$.
%Given a tree $\vartheta\in\Theta$ we define the \emph{trimmed order} of $\vartheta$ as 
%
%\begin{equation}\label{ordine}
%\breve{k}(\vartheta) := \sum_{\substack{v\in N(\vartheta) \\ s_v=5}} k_v = |\{ v \in N(\vartheta) : s_v=5 \}|.
%\end{equation}
%
%Consider the equivalence class of the trees $\vartheta\in\Theta$ with trimmed order $\breve{k}(\vartheta)=k$ such that the root line
%has component $j$, momentum $\nu$ and sign $\s$, and define \emph{trimmed tree} any the tree obtained from any
%representative of the equivalence class by replacing the values of the labels $k_v$ of the nodes $v$ with 
%$\breve\gotT^{(\breve{k})}_{j,\nu,\s}$.$\breve\gotT^{(k)}_{j,\nu,\s}$
\end{defi}
%%%%%%%%%%%%%%%%%%%%%%%%%%%%%%%%%%%%%%%%%%%%%%%%%%%%%%%%%%%%%%%%%%%%%%%%%%

%%%%%%%%%%%%%%%%%%%%%%%%%%%%%%%%%%%%%%%%%%%%%%%%%%%%%%%%%%%%%%%%%%%%%%%%%%
\begin{rmk}\label{rmkpro}
\emph{
Mind that in Definition \ref{potati} the index $k$ fixes the number of nodes, differently from Definition \ref{gotT} where it
fixes the order.
}
\end{rmk}
%%%%%%%%%%%%%%%%%%%%%%%%%%%%%%%%%%%%%%%%%%%%%%%%%%%%%%%%%%%%%%%%%%%%%%%%%

%%%%%%%%%%%%%%%%%%%%%%%%%%%%%%%%%%%%%%%%%%%%%%%%%%%%%%%%%%%%%%%%%%%%%%%%%%
\begin{rmk}\label{rmkcijei}
\emph{
{Lemma \ref{cijei} extends trivially to the trimmed trees $\vartheta\in\breve\gotT$.}
}
\end{rmk}
%%%%%%%%%%%%%%%%%%%%%%%%%%%%%%%%%%%%%%%%%%%%%%%%%%%%%%%%%%%%%%%%%%%%%%%%%%

%%%%%%%%%%%%%%%%%%%%%%%%%%%%%%%%%%%%%%%%%%%%%%%%%%%%%%%%%%%%%%%%%%%%%%%%%%
\begin{defi}[\textbf{Value of a trimmed tree}]\label{defvalue3}
{For any trimmed tree $\vartheta\in \breve\gotT$ we define the \emph{value} of $\vartheta$ as 
\begin{equation}\label{trival}
\Val(\vartheta;c,\om) := \Big(\prod_{\la\in \Lambda(\vartheta)}\LL_\la(c) \Big) 
\Big(\prod_{\ell\in L(\vartheta)} \matG_\ell(\om) \Big) ,
\end{equation}
with $\LL_\la(c)$ and $\matG_\ell(\om)$ defined according to \eqref{foglie} and \eqref{prop}, respectively.}
\end{defi}
%%%%%%%%%%%%%%%%%%%%%%%%%%%%%%%%%%%%%%%%%%%%%%%%%%%%%%%%%%%%%%%%%%%%%%%%%%

%%%%%%%%%%%%%%%%%%%%%%%%%%%%%%%%%%%%%%%%%%%%%%%%%%%%%%%%%%%%%%%%%%%%%%%%%%
\begin{rmk}\label{didi}
\emph{
One has
 \[
 \Theta\cap\gotT\cap\breve\gotT = \{\vartheta \in\Theta\cup\gotT\cup\breve\gotT\, : \, s_v=5,\ \mbox{ for all }v\in N(\vartheta)\}.
 \] 
In particular, if $\vartheta\in  \Theta\cap\gotT\cap\breve\gotT$ one has
\[
\VV(\vartheta;c,\om,\h)=\Val(\vartheta;c,\om).
\]
}
\end{rmk}
%%%%%%%%%%%%%%%%%%%%%%%%%%%%%%%%%%%%%%%%%%%%%%%%%%%%%%%%%%%%%%%%%%%%%%%%%%

%%%%%%%%%%%%%%%%%%%%%%%%%%%%%%%%%%%%%%%%%%%%%%%%%%%%%%%%%%%%%%%%%%%%%%%%%%
\begin{defi}[\textbf{Trimmed tree associated with an expanded tree}]
\label{trimmed}
%For all $k$, $\nu$, $j$, $\s$
For any expanded tree $\vartheta\in\gotT$, set
\begin{equation} \label{N2theta}
N_2(\vartheta) :=\{ v \in N(\vartheta) :s_{v}=2\}.% \hbox{ and } v \st{\preceq} r_{\vartheta}\}  . 
\end{equation}
For any $v\in N_2(\vartheta)$, if $\vartheta_v$ denotes the subtree of $\vartheta$ entering $v$ whose root line is an $\h$-line, set 
\begin{equation} \label{F2theta}
F_2^*(\vartheta) := \{ \vartheta_v : v \in N_2(\vartheta) \hbox{ and } v \st{\preceq} r_{\vartheta}\}.
\end{equation}
Let $\breve\vartheta$ be the trimmed tree obtained from $\vartheta$ by cutting the subtrees in $F_2^*(\vartheta)$;
we call $\breve\vartheta$ the \emph{trimmed tree associated with the tree $\vartheta$.}
\end{defi}
%%%%%%%%%%%%%%%%%%%%%%%%%%%%%%%%%%%%%%%%%%%%%%%%%%%%%%%%%%%%%%%%%%%%%%%%%%

Examples of trimmed trees can be obtained by looking at Figures \ref{fig-expanded} to \ref{fig-exp-unexp}.
Consider the expanded tree $\vartheta$ in Figure \ref{fig-expanded}.
The trimmed tree $\breve\vartheta$ associated with $\vartheta$ is the unexpanded tree $\theta_1$ in Figure \ref{fig-unexpanded},
while, if $\vartheta_{v_1}$ is the expanded tree grafted to the node $v_1$ of $\vartheta$, then $\vartheta_{v_1}=\vartheta''$
and $\breve\vartheta_{v_1}=\vartheta_3$, with $\vartheta''$ and $\vartheta_3$ as in \ref{fig-unexpanded}
and \ref{fig-exp-unexp}, respectively. Finally, if $\vartheta_{v_3}$ is the expanded tree grafted to $v_3$, 
then $\vartheta_{v_3}=\breve\vartheta_{v_3}=\vartheta_2$, since $\vartheta_{v_3}$ has no node $v$ with $s_v=2$.

%%%%%%%%%%%%%%%%%%%%%%%%%%%%%%%%%%%%%%%%%%%%%%%%%%%%%%%%%%%%%%%%%%%%%%%%%%
\begin{rmk}\label{stracon}
\emph{
By construction,
with the notation of Definition \ref{trimmed}, one has
$$ 
k(\breve\vartheta) = \sum_{v \in N(\breve\vartheta)} k_v =
k(\vartheta) - \sum_{\substack{v\in N_2(\vartheta) \\ v \st{\preceq} r_{\vartheta}}} k(\vartheta_v ),
$$
and 
 \begin{equation}\label{laseconda}
 \Val(\vartheta;c,\om) = \Val(\breve{\vartheta};c,\om)\prod_{\substack{v \in N_2(\vartheta) \\ v \st{\preceq} r_{\vartheta}}}
\Big(- \frac{1}{c_{j_{v}}^{\sigma_{v}}}\Val(\vartheta_{v};c,\om)\Big) 
\end{equation}
which can be iterated so as to obtain
\begin{equation}\label{disposta}
\Val(\vartheta;c,\om)=\Val(\breve{\vartheta};c,\om)\prod_{v \in N_2(\vartheta)}
\Big(- \frac{1}{c_{j_{v}}^{\sigma_{v}}} \Val(\breve\vartheta_v;c,\om) \Big).
\end{equation}
}
\end{rmk}
%%%%%%%%%%%%%%%%%%%%%%%%%%%%%%%%%%%%%%%%%%%%%%%%%%%%%%%%%%%%%%%%%%%%%%%%%%

If $T$ is a subgraph of a tree $\vartheta\in\Theta\cup \gotT \cup \breve\gotT $
we set $N(T):=N(\vartheta)\cap V(T)$, $\Lambda(T):=\Lambda(\vartheta)\cap V(T)$
and $\Lambda_{j,\s}(T):=\Lambda_{j,\s}(\vartheta)\cap V(T)$,
and we define the \emph{order} of $T$ as
\begin{equation}\label{ordineT}
k(T) := \sum_{v\in N(T)} k_v ,
\end{equation}
and set
\begin{equation}\label{n2t}
N_2(T):=N_2(\vartheta)\cap N(T).
\end{equation}

%%%%%%%%%%%%%%%%%%%%%%%%%%%%%%%%%%%%%%%%%%%%%%%%%%%%%%%%%%%%%%%%%%%%%%%%%%
\begin{defi}[\textbf{Self-energy graph}]
\label{self-energy}
A \emph{self-energy graph} (SEG) of a tree $\vartheta$
is a connected subgraph $T$ of $\vartheta$ which 
has only one entering line $\ell'_T$ and one exiting line $\ell_T$.
We set $\calP_T:=\calP(\ell_T',\ell_T)$ and call $\calP_T$
the \emph{path of the self-energy graph} $T$.
\end{defi}
%%%%%%%%%%%%%%%%%%%%%%%%%%%%%%%%%%%%%%%%%%%%%%%%%%%%%%%%%%%%%%%%%%%%%%%%%%

%%%%%%%%%%%%%%%%%%%%%%%%%%%%%%%%%%%%%%%%%%%%%%%%%%%%%%%%%%%%%%%%%%%%%%%%%%
\begin{rmk} \label{seg}
\emph{
The term``self-energy graph'' is borrowed from quantum field theory, where it denotes a special Feynman graph -- called also
\emph{mass graph} -- with two external lines \cite{AGD,NO}.
}
\end{rmk}
%%%%%%%%%%%%%%%%%%%%%%%%%%%%%%%%%%%%%%%%%%%%%%%%%%%%%%%%%%%%%%%%%%%%%%%%%%

%%%%%%%%%%%%%%%%%%%%%%%%%%%%%%%%%%%%%%%%%%%%%%%%%%%%%%%%%%%%%%%%%%%%%%%%%%
\begin{defi}[\textbf{Value of a subgraph}]\label{defvalue4}
{We define the \emph{value} of a subgraph $T$ of a tree $\vartheta\in\Theta$ as
\begin{equation}\label{valT}
	\VV(T;c,\om,\h) := \Big(\prod_{\la\in \Lambda(T)}\LL_\la(c) \Big) \Big(\prod_{v\in N(T)} \calF_v(\h)\Big)
	\Big(\prod_{\ell\in L(T)} \matG_\ell(\om) \Big),
\end{equation}
of a subgraph $T$ of a tree $\vartheta\in\gotT$ as
\begin{equation}\label{valetaT}
\Val(T;c,\om):=\Big(\prod_{\la\in \Lambda(T)}\LL_\la (c)\Big) 
\Big(\prod_{v\in N(T)} \starF_v(c)\Big)
\Big(\prod_{\ell\in L(T)} \matG_\ell(\om) \Big),
\end{equation}
and of a subgraph $T$ of a tree $\vartheta\in\breve\gotT$ as
\begin{equation}\label{valTtri}
	\Val(T;c,\om) := \Big(\prod_{\la\in \Lambda(T)}\LL_\la(c) \Big) 
	\Big(\prod_{\ell\in L(T)} \matG_\ell(\om) \Big) ,
\end{equation}
where $\LL_\la(c)$, $\calF_v(\h)$, $\calF_v^*(c)$ and $\matG_\ell(\om)$ are defined 
according to \eqref{foglie}, \eqref{nodi}, \eqref{nodiker} and \eqref{prop}, respectively.}
\end{defi}
%%%%%%%%%%%%%%%%%%%%%%%%%%%%%%%%%%%%%%%%%%%%%%%%%%%%%%%%%%%%%%%%%%%%%%%%%%

%%%%%%%%%%%%%%%%%%%%%%%%%%%%%%%%%%%%%%%%%%%%%%%%%%%%%%%%%%%%%%%%%%%%%%%%%%
\begin{rmk} \label{aereo}
\emph{
If a connected subgraph $T$ of a tree $\vartheta$ does not contain the root of $\vartheta$, then $T$ has one
and only one exiting line, which we still call $\ell_T$ as in the case of self-energy graphs.
}
\end{rmk}
%%%%%%%%%%%%%%%%%%%%%%%%%%%%%%%%%%%%%%%%%%%%%%%%%%%%%%%%%%%%%%%%%%%%%%%%%%

%%%%%%%%%%%%%%%%%%%%%%%%%%%%%%%%%%%%%%%%%%%%%%%%%%%%%%%%%%%%%%%%%%%%%%%%%%
\begin{defi}[\textbf{Trimmed subgraph of an expanded tree}] 
\label{trimmedsubgraph1}
{Given a subgraph $T$ of an expanded tree $\vartheta$, we define $\breve T$ as the subgraph of $\breve\vartheta$
such that $V(\breve T) = V(T)\cap V(\breve\vartheta)$ and $L(\breve T)= L(T )\cap L(\breve\vartheta)$.
We call $\breve T$ the \emph{trimmed subgraph associated with the subgraph} $T$.}
\end{defi}
%%%%%%%%%%%%%%%%%%%%%%%%%%%%%%%%%%%%%%%%%%%%%%%%%%%%%%%%%%%%%%%%%%%%%%%%%%

%%%%%%%%%%%%%%%%%%%%%%%%%%%%%%%%%%%%%%%%%%%%%%%%%%%%%%%%%%%%%%%%%%%%%%%%%%
\subsection{Resonant clusters}
%%%%%%%%%%%%%%%%%%%%%%%%%%%%%%%%%%%%%%%%%%%%%%%%%%%%%%%%%%%%%%%%%%%%%%%%%%

For any $\vartheta\in \Theta\cup\breve\gotT$ define %\Theta^{(k)}_{j,\nu,\s}$ define
\begin{equation}\label{J}
J(\vartheta):= \sum_{\la\in \Lambda(\vartheta)} \jap{j_\la}^\al - \jap{j}^\al,
\end{equation}
where $j$ is the component associated with the root line of $\vartheta$.

%%%%%%%%%%%%%%%%%%%%%%%%%%%%%%%%%%%%%%%%%%%%%%%%%%%%%%%%%%%%%%%%%%%%%%%%%%
\begin{rmk} \label{altrormk}
\emph{
For any $\vartheta\in\Theta\cup\breve\gotT$
by \eqref{foglie} and \eqref{decay}, %and the third of \eqref{lefoglie}, 
we deduce that
\begin{equation}\label{assai}
\left|\prod_{\la\in \Lambda(\vartheta)} \LL_\la(c) \right|=
\left|\prod_{\la\in \Lambda(\vartheta)} e^{-s\jap{j_\la}^\al} \g_{j_{\la}}\right|\le
e^{- s \sum_{\la\in \Lambda(\vartheta)  }\jap{j_{\la}}^\al} =
e^{-s(\jap{j}^\al +J(\vartheta))}\,.
\end{equation}
}
\end{rmk}
%%%%%%%%%%%%%%%%%%%%%%%%%%%%%%%%%%%%%%%%%%%%%%%%%%%%%%%%%%%%%%%%%%%%%%%%%%

If $T$ is a subgraph of $\vartheta \in \Theta \cup\gotT\cup\breve\gotT$ which does not contain the root of $\vartheta$, 
let $\Lambda^*(T)$ denote the set of leaves $\la\in \Lambda(T)$
such that $\la\st{\preceq}\ell_T$,
 and let $\Lambda^*_{j,\s}(T)$ denote the set of leaves $\la\in \Lambda_{j,\s}(T)$
such that $\la\st{\preceq}\ell_T$.
For any  connected subgraph $T$ of a tree $\vartheta \in \Theta\cup\gotT\cup\breve\gotT $ with only one entering line $\ell'$ and one exiting line $\ell$
%(and $T$ is not a RC) 
we set
\begin{equation}\label{JT}
J(T):=\sum_{\la \in \Lambda^*(T)} \jap{j_\la}^\al - \jap{j_\ell}^\al + \jap{j_{\ell'}}^\al,
\end{equation}
where $\Lambda^*(T)=\Lambda(T)$ if $T\subset\vartheta\in\Theta\cup\breve\gotT$
and $\Lambda^*(T)=\Lambda(\breve T)$ if $T\subset\vartheta\in\gotT$.

%%%%%%%%%%%%%%%%%%%%%%%%%%%%%%%%%%%%%%%%%%%%%%%%%%%%%%%%%%%%%%%%%%%%%%%%%%
\begin{rmk} \label{utile}
\emph{
Given a tree $\vartheta\in \Theta\cup\breve\gotT$, for any subtree {$\vartheta'\subseteq\vartheta$} one has $J(\vartheta') \le J(\vartheta)$;
indeed, if $j$ is the component of the root line of $\vartheta$ and $j'$ is the component of the root line of $\vartheta'$, one has
\begin{equation} \nonumber
J(\vartheta) = \sum_{\la \in \Lambda(\vartheta')} \jap{j_\la}^\al - \jap{j'}^\al +
\sum_{\la \in \Lambda(\vartheta)\setminus\Lambda(\vartheta')} \jap{j_\la}^\al +  \jap{j'}^\al - \jap{j}^\al ,
\end{equation}
with
\begin{equation} \nonumber
\sum_{\la \in \Lambda(\vartheta')} \jap{j_\la}^\al - \jap{j'}^\al = J(\vartheta') , \qquad
\jap{j}^\al \le \sum_{\la \in \Lambda(\vartheta)\setminus\Lambda(\vartheta')} \jap{j_\la}^\al +  \jap{j'}^\al ,
\end{equation}
where the last inequality follows from
the third relation in \eqref{lefoglie*} and the subadditivity of the function $x \mapsto x^{\al}$ for $\al\in(0,1)$.
Similarly, if $\vartheta_1',\ldots,\vartheta_p'\subseteq\vartheta$, then
\begin{equation}\label{i}
\sum_{i=1}^p J(\vartheta_i')\le J(\vartheta).
\end{equation}
}
\end{rmk}
%%%%%%%%%%%%%%%%%%%%%%%%%%%%%%%%%%%%%%%%%%%%%%%%%%%%%%%%%%%%%%%%%%%%%%%%%%

%%%%%%%%%%%%%%%%%%%%%%%%%%%%%%%%%%%%%%%%%%%%%%%%%%%%%%%%%%%%%%%%%%%%%%%%%%
\begin{rmk} \label{utilemaseparato}
\emph{
If $T$ is a connected subgraph of $\vartheta\in \Theta\cup\breve\gotT$
with only one entering line $\ell'$ and one exiting line $\ell$, and $\vartheta'$ is a subtree contained in $T$, one has
$J(\vartheta') \le J(T)$, which follows immediately by writing
\begin{equation} \nonumber
J(T) = J(\vartheta') + \sum_{\la \in \Lambda(T)\setminus\Lambda(\vartheta')} \jap{j_\la}^\al +  \jap{j_{\ell'}}^\al + \jap{j'}^\al - \jap{j}^\al ,
\end{equation}
where $j$ is the component of the line $\ell$ and $j'$ is the component of the root line of $\vartheta'$, and noting that
the root line of $\vartheta'$ can not belong to $\calP_T$ and hence
\begin{equation} \nonumber
\jap{j}^\al \le \sum_{\la \in \Lambda(T)\setminus\Lambda(\vartheta')} \jap{j_\la}^\al +  \jap{j_{\ell'}}^\al + \jap{j'}^\al ,
\end{equation}
by the same argument as in Remark \ref{utile}. By reasoning in a similar way, one finds that if $T'$ is a connected subgraph of $T$
with only one entering line $\ell'_{T'}$ and one exiting line $\ell_{T'}$ both belonging to $\calP_T\cup\{\ell_T,\ell'_{T}\}$, then
$J(T')\le J(T)$. Finally if $T'\subset T \subset \vartheta\in\gotT$ and both $T'$ and $T$ have only one entering line,
respectively $\ell'_{T'}$ and $ \ell'_T$, and one exiting line, 
respectively $\ell_{T'}$ and $\ell_T$, and moreover $\ell_{T'}\st\prec\ell_T$, then $J(T')\le J(T)$.
}
\end{rmk}
%%%%%%%%%%%%%%%%%%%%%%%%%%%%%%%%%%%%%%%%%%%%%%%%%%%%%%%%%%%%%%%%%%%%%%%%%%

%%%%%%%%%%%%%%%%%%%%%%%%%%%%%%%%%%%%%%%%%%%%%%%%%%%%%%%%%%%%%%%%%%%%%%%%%%
\begin{defi}[\textbf{Resonant cluster}]\label{coso}
Let  $C_1=C_1(\al)$ an appropriate constant to be specified in Lemma \ref{constance2}.
Given a tree $\vartheta\in\Theta\cup\gotT\cup\breve\gotT$, %with $\VV(\vartheta;c,\om,\h)\ne0$,
a \textit{resonant cluster} (RC) is a connected subgraph $T$ of $\vartheta$ with only
one entering line $\ell'_T$ and one exiting line $\ell_T$  such that
\begin{enumerate}[topsep=0ex]
\itemsep0.0em
\item
either $T$ consists only of one node $v$ %with $s_v=1$
\item
or, setting
%\begin{equation}\label{maxT}
%$\nmin_T :=\min\{n_{\ell_T},n_{\ell_T'}\}$ and 
$\nmax_T:=\max\{n_\ell :\ell\in L(T)\}$ and $\nmin_T :=\min\{n_{\ell_T},n_{\ell_T'}\}$,
%\end{equation}
the following  holds:
\begin{enumerate}[topsep=0ex]
\itemsep0.0em
\item[2.1.] 
$\nmax_T<\nmin_T$,
\item[2.2.]
one has
\begin{equation}\label{secindipT}
\sum_{\la\in \Lambda^*(T)} \s_\la \gote_{j_\la} + \s_{\ell'_T}\gote_{j_{\ell'_T}}- \s_{\ell_T}\gote_{j_{\ell_T}}=0\,,
\end{equation}
%\begin{equation}\label{secindip}
%\sum_{\la\in E(T)} \s_\la \gote_{j_\la} + \s_{\ell'_T}\gote_{j_{\ell'_T}}- \s_{\ell_T}\gote_{j_{\ell_T}}=0\,,
%\end{equation}
%and $\nmax_\ell  < \nmin_{\ell_T}=\nmin_{\ell'_T}$ for all $\ell\in L(T)$.
%
\item[2.3.] 
for all $\ell\in\calP_T$ one has
$\s_{\ell}(\nu_{\ell} -\gote_{j_{\ell}}) \ne \s_{\ell'_T}(\nu_{\ell'_T} - \gote_{j_{\ell'_T}})$,
\item[2.4.]\label{RCpath} 
if $\vartheta\in\gotT$, then the path $\calP_T$ does not contain any $\h$-line,
\item[2.5.]
$J(T) < C_1 r_{m_{\nmin_T}-1}$,
\item[2.6.]
if $\vartheta\in\gotT$, then $J(\breve\vartheta_v)< C_1 r_{m_{\nmin_T}-1}$ for all $v\in N_2(T)$, with $N_2(T)$ defined in \eqref{n2t}.
\end{enumerate}
\end{enumerate}
We call $\ell_T$ and $\ell'_T$ the \emph{external lines} of $T$, and $\calP_T$ the \emph{resonant path} of $T$.
When $L(T)\neq\emptyset$,
we say that $\nmax_T$ is the \emph{scale} of $T$. Finally we call \emph{trivial} the RC consisting only of one node $v$, %with $s_v=1$,
and we say that it has scale $-1$.
\end{defi}
%%%%%%%%%%%%%%%%%%%%%%%%%%%%%%%%%%%%%%%%%%%%%%%%%%%%%%%%%%%%%%%%%%%%%%%%%%

%%%%%%%%%%%%%%%%%%%%%%%%%%%%%%%%%%%%%%%%%%%%%%%%%%%%%%%%%%%%%%%%%%%%%%%%%%
\begin{rmk} \label{utile?}
\emph{
If $T$ is a RC of a tree $\vartheta\in\gotT$, then, by the condition \ref{RCpath}.4 of Definition \ref{coso},
for any $v\in N_2(\vartheta)$ the tree $\vartheta_v$ -- and hence $\breve\vartheta_v$ as well -- is either external or internal to $T$.
Thus, the cutting of $\vartheta$ naturally produces the trimmed RC $\breve T \in \breve\vartheta$.
}
\end{rmk}
%%%%%%%%%%%%%%%%%%%%%%%%%%%%%%%%%%%%%%%%%%%%%%%%%%%%%%%%%%%%%%%%%%%%%%%%%%

%%%%%%%%%%%%%%%%%%%%%%%%%%%%%%%%%%%%%%%%%%%%%%%%%%%%%%%%%%%%%%%%%%%%%%%%%%
 \begin{rmk}\label{fratellini}
 \emph{
From \eqref{secindipT} and \eqref{lefoglie*} one has
$\s_{\ell_T}(\nu_{\ell_T} -\gote_{j_{\ell_T}}) = \s_{\ell'_T}(\nu_{\ell'_T} - \gote_{j_{\ell'_T}})$
and hence
\begin{equation}\label{quellavera}
{\s_{\ell_T} x_{\ell_T} =  \s_{\ell'_T} x_{\ell'_T} \,. }
\end{equation}
Then, if $T$ is a RC  with $\s_{\ell_T}=\s_{\ell'_T}$,
$\matG_{\ell_T}(\om)$ and $\matG_{\ell'_T}(\om)$ have the same denominator,
whereas,
if the RC $T$ is such that $\s_{\ell_T}=-\s_{\ell'_T}$,
the denominators of $\matG_{\ell_T}(\om)$ and $\matG_{\ell'_T}(\om)$ 
are opposite to each other. 
}
\end{rmk}
%%%%%%%%%%%%%%%%%%%%%%%%%%%%%%%%%%%%%%%%%%%%%%%%%%%%%%%%%%%%%%%%%%%%%%%%%%

%%%%%%%%%%%%%%%%%%%%%%%%%%%%%%%%%%%%%%%%%%%%%%%%%%%%%%%%%%%%%%%%%%%%%%%%%%
\begin{rmk}\label{nuzero}
\emph{%remark
By the conservation law \eqref{lefoglie*}, for any RC $T$ and for any $\ell\in\calP_T$
one has $\s_\ell\nu_{\ell}=\nu_{\ell}^0+ \s_{\ell'_T}\nu_{\ell'_T}$, where
\begin{equation}\label{ellezero}
\nu^0_\ell=\nu_\ell^0(T) := \sum_{\substack{\la\in \Lambda^*(T) \\ \la\prec \ell}} \s_{\la}\gote_{j_{\la}}.
\end{equation}
%Moreover, because of the constraint on the scales in Definition \ref{coso}, for any line $\ell\in L(T)$ -- and hence
%in particular for any $\ell\in\calP_T$ -- one has 
%$$\s_\ell\nu_\ell - \s_\ell\gote_{j_\ell}\ne \s_{\ell'_T}\nu_{\ell'_T} - \s_{\ell'_T}\gote_{j_{\ell'_T}}.$$
}%remark
\end{rmk}
%%%%%%%%%%%%%%%%%%%%%%%%%%%%%%%%%%%%%%%%%%%%%%%%%%%%%%%%%%%%%%%%%%%%%%%%%%

%%%%%%%%%%%%%%%%%%%%%%%%%%%%%%%%%%%%%%%%%%%%%%%%%%%%%%%%%%%%%%%%%%%%%%%%%%
\begin{defi}[\textbf{Resonant line}]
\label{resort}
Given any $\vartheta\in\Theta\cup\gotT\cup\breve\gotT$
we say that a line $\ell \in L(\vartheta)$ is \textit{resonant} if it is both the exiting line of a RC in $\vartheta$
and the entering line of a RC in $\vartheta$,
otherwise we say $\ell$ is \textit{non-resonant}.
\end{defi}
%%%%%%%%%%%%%%%%%%%%%%%%%%%%%%%%%%%%%%%%%%%%%%%%%%%%%%%%%%%%%%%%%%%%%%%%%%

%%%%%%%%%%%%%%%%%%%%%%%%%%%%%%%%%%%%%%%%%%%%%%%%%%%%%%%%%%%%%%%%%%%%%%%%%%
\begin{defi}[\textbf{Chain of resonant clusters}]
\label{chain}
Given a tree $\vartheta\in\Theta\cup\gotT\cup\breve\gotT$, a \emph{chain} of RCs in $\vartheta$
is a set $\gotC=\{T_1,\ldots,T_p\}$, where $T_1,\ldots,T_p$, with $p \ge 2$, are RCs  of $\vartheta$
such that $\ell_{T_{i+1}}=\ell_{T_{i}}'$ for 
$i=1,\ldots,p-1$ and both $\ell_{T_1}$ and $\ell_{T_p}'$ are non-resonant. We call $\ell_\gotC:=\ell_{T_1}$ and $\ell'_\gotC:=\ell'_{T_p}$ the 
\emph{exiting line} and the \emph{entering line} of $\gotC$, respectively, and $\ell_{T_2},\ldots,\ell_{T_{p}}$ the \emph{links} of $\gotC$.
\end{defi}
%%%%%%%%%%%%%%%%%%%%%%%%%%%%%%%%%%%%%%%%%%%%%%%%%%%%%%%%%%%%%%%%%%%%%%%%%%

%%%%%%%%%%%%%%%%%%%%%%%%%%%%%%%%%%%%%%%%%%%%%%%%%%%%%%%%%%%%%%%%%%%%%%%%%%
% FIGURA 19-2
%%%%%%%%%%%%%%%%%%%%%%%%%%%%%%%%%%%%%%%%%%%%%%%%%%%%%%%%%%%%%%%%%%%%%%%%%%
\begin{figure}[ht]
\vspace{-.2cm}
\centering
%\null
%\hspace{-.6cm}
\ins{080pt}{-120pt}{$T_1$}
\ins{154pt}{-120pt}{$T_2$}
\ins{266pt}{-120pt}{$T_3$}
\ins{134pt}{-049pt}{$\ell_{T_2}$}
\ins{243pt}{-065pt}{$\ell_{T_3}$}
\ins{58pt}{-092pt}{$\ell_{T_1}$}
\ins{382pt}{-085pt}{$\ell_{T_3}'$}
%\ins{346pt}{-058pt}{$\ell^{in}_\TT$}
\subfigure{\includegraphics*[width=6in]{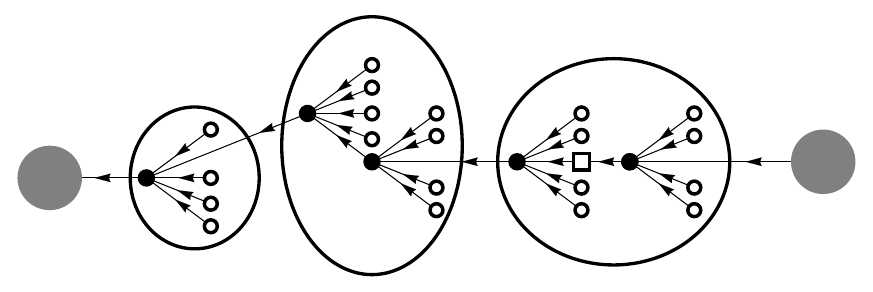}}
\caption{\small A chain of 3 resonant clusters  $\gotC=\{T_1,T_2,T_3\}$; the lines $\ell_{T_2}$ and $\ell_{T_3}$ are resonant.}
\label{catena}
\end{figure}
%%%%%%%%%%%%%%%%%%%%%%%%%%%%%%%%%%%%%%%%%%%%%%%%%%%%%%%%%%%%%%%%%%%%%%%%%%

An example of chain in a tree $\vartheta\in\Theta$ is drawn in Figure \ref{catena}, with the grey circles representing the remaining part
of $\vartheta$; in particular the grey circle to the left contains the root of $\vartheta$.

%%%%%%%%%%%%%%%%%%%%%%%%%%%%%%%%%%%%%%%%%%%%%%%%%%%%%%%%%%%%%%%%%%%%%%%%%%
\begin{defi}[\textbf{Relevant resonant cluster}]
\label{rilevante}
A RC $T$ of a tree $\vartheta$ is \emph{relevant} if there is a chain $\gotC$ in $\vartheta$ such that $T\in\gotC$.
%Given any subgraph $S$ of a tree $\vartheta$, we define $\fT_0(S)$ as the set of relevant RCs contained in $S$. 
\end{defi}
%%%%%%%%%%%%%%%%%%%%%%%%%%%%%%%%%%%%%%%%%%%%%%%%%%%%%%%%%%%%%%%%%%%%%%%%%%

%%%%%%%%%%%%%%%%%%%%%%%%%%%%%%%%%%%%%%%%%%%%%%%%%%%%%%%%%%%%%%%%%%%%%%%%%%
\begin{rmk} \label{percapireacheserveladef}
\emph{
As stressed at the beginning of the section, the product of the small divisors can get out of control when they accumulate.
In turn, for this actually to occur, the same small divisors must appear over and over again. This means that
the resonant clusters are a real source of problems -- and hence are considered ``relevant'' -- not when they are isolated,
but only when they belong to a chain of arbitrarily large length.
}
\end{rmk}
%%%%%%%%%%%%%%%%%%%%%%%%%%%%%%%%%%%%%%%%%%%%%%%%%%%%%%%%%%%%%%%%%%%%%%%%%%

%%%%%%%%%%%%%%%%%%%%%%%%%%%%%%%%%%%%%%%%%%%%%%%%%%%%%%%%%%%%%%%%%%%%%%%%%%
%%%%%%%%%%%%%%%%%%%%%%%%%%%%%%%%%%%%%%%%%%%%%%%%%%%%%%%%%%%%%%%%%%%%%%%%%%
 \section{Well-definedness of the power series}
 \label{stimazze}
\zerarcounters
%%%%%%%%%%%%%%%%%%%%%%%%%%%%%%%%%%%%%%%%%%%%%%%%%%%%%%%%%%%%%%%%%%%%%%%%%%
%%%%%%%%%%%%%%%%%%%%%%%%%%%%%%%%%%%%%%%%%%%%%%%%%%%%%%%%%%%%%%%%%%%%%%%%%%

As mentioned at the beginning of Section \ref{sec}, we now want to show that
the only possible obstruction to the convergence of the series is the 
accumulation of small divisors of the same size which is produced by the presence of chains of RCs in trees, 
(see Remark \ref{percapireacheserveladef}). 

More precisely, in this section we prove that if, in \eqref{trival},
we could neglect the propagators of the resonant lines due to the chains of RCs, we would be able
to bound the value of each trimmed tree -- and hence, thanks to \eqref{disposta}, of each expanded tree --
in such a way to ensure the convergence of the series.
Then, in Section \ref{convergenza}, we shall see how to deal with the presence of the resonant lines.

%%%%%%%%%%%%%%%%%%%%%%%%%%%%%%%%%%%%%%%%%%%%%%%%%%%%%%%%%%%%%%%%%%%%%%%%%%
%%%%%%%%%%%%%%%%%%%%%%%%%%%%%%%%%%%%%%%%%%%%%%%%%%%%%%%%%%%%%%%%%%%%%%%%%%
\subsection{Preliminary bounds}
\label{vedo}
%%%%%%%%%%%%%%%%%%%%%%%%%%%%%%%%%%%%%%%%%%%%%%%%%%%%%%%%%%%%%%%%%%%%%%%%%%
%%%%%%%%%%%%%%%%%%%%%%%%%%%%%%%%%%%%%%%%%%%%%%%%%%%%%%%%%%%%%%%%%%%%%%%%%%

Recall that $s>0$ and $\al\in(0,1)$ are fixed.
Given a tree $\vartheta\in\breve\gotT$, set $\widehat\Lambda(\vartheta):=
\Lambda(\vartheta)\cup\{r_\vartheta\}$ and $p=p(\vartheta):=|\widehat\Lambda(\vartheta)|$.
After numbering the set $\widehat\Lambda(\vartheta)=\{\la_{1},\ldots\la_p\}$,  consider
the sequence $\jap{j_{\la_1}},\ldots,\jap{j_{\la_p}}$ and rearrange it in decreasing order, i.e.~set
\begin{equation}\label{glienne}
\begin{aligned}
&\gotn_1(\vartheta) := \max \{\jap{j_{\la_1}},\ldots,\jap{j_{\la_p}}\} , \\
&\gotn_2(\vartheta) := \max \{\jap{j_{\la_1}},\ldots,\jap{j_{\la_p}}\} \setminus \{\gotn_1(\vartheta)\} , \\
& \gotn_3 (\vartheta):= \max \{\jap{j_{\la_1}},\ldots,\jap{j_{\la_p}}\} \setminus \{\gotn_1(\vartheta),\gotn_2(\vartheta) \} , \\
&\qquad \vdots\\
&\gotn_{p} (\vartheta):= \min \{\jap{j_{\la_1}},\ldots,\jap{j_{\la_p}}\} ,
\end{aligned}
\end{equation}
so that $\gotn_1(\vartheta)\ge \gotn_2(\vartheta)\ge \ldots \ge \gotn_{p}(\vartheta)$.
Recall also the definition of $J(\vartheta)$ in \eqref{J}.

%%%%%%%%%%%%%%%%%%%%%%%%%%%%%%%%%%%%%%%%%%%%%%%%%%%%%%%%%%%%%%%%%%%%%%%%%%
\begin{lemma}\label{meglio}
For any $\vartheta\in \Theta\cup \breve\gotT$, one has
\begin{equation}\label{jean}
J(\vartheta) %=\sum_{\la\in\Lambda(\vartheta)} \jap{j_{\la}}^\al  - \jap{j}^\al
\ge  ({2- 2^\al}) \sum_{i=3}^{p(\vartheta)} \gotn_i^\al (\vartheta) .
\end{equation}
\end{lemma}
%%%%%%%%%%%%%%%%%%%%%%%%%%%%%%%%%%%%%%%%%%%%%%%%%%%%%%%%%%%%%%%%%%%%%%%%%%

%%%%%%%%%%%%%%%%%%%%%%%%%%%%%%%%%%%%%%%%%%%%%%%%%%%%%%%%%%%%%%%%%%%%%%%%%%
\prova
If $\vartheta$ is the trivial tree, then $p(\vartheta)=2$ and \eqref{jean} is trivial. Otherwise $p=p(\vartheta)\ge6$.
By the {third} of \eqref{lefoglie*}, applied to the root line, we conclude that there is a sequence $\hat{\s}_1,\ldots,\hat{\s}_{p}=\pm1,0$
such that 
\[
\sum_{i=1}^{p} \hat{\s}_i \, \gotn_i(\vartheta) =0,
\]
and $\hat{\s}_i=0$ implies $\gotn_i(\vartheta)=1$.
If $\hat{\s}_1=0$ then, since the sequence $\gotn_i(\vartheta)$ is decreasing, we have $\gotn_i(\vartheta)=1$ 
for all $i=1,\ldots,p$ and hence the result is trivial.
Otherwise the result follows from Lemma \ref{constance}.
\EP
%%%%%%%%%%%%%%%%%%%%%%%%%%%%%%%%%%%%%%%%%%%%%%%%%%%%%%%%%%%%%%%%%%%%%%%%%%

For any $k,j,\nu,\s$, given a tree $\vartheta\in %\gotT^{(k)}_{j,\nu,\s}\cup
\breve\gotT^{(k)}_{j,\nu,\s}$, 
let $\calmX(\vartheta)$ denote the set of all trees $\vartheta'\in %\gotT^{(k)}_{j,\nu,\s}\cup
\breve\gotT^{(k)}_{j,\nu,\s}$ which are obtained
from $\vartheta$ by  changing the component labels $j_\la$ but not the sign labels $\s_\la$ of all the leaves $\la\in \Lambda(\vartheta)$, 
and changing all the remaining labels consistently. Similarly, 
for any $k,j,\nu,\s$, given a tree $\vartheta\in\gotT^{(k)}_{j,\nu,\s}$
we call $\calmX(\vartheta)$  the set of all trees $\vartheta'\in \gotT^{(k)}_{j,\nu,\s}$  which are obtained
from $\vartheta$ by  changing the component labels $j_\la$ but not the sign labels $\s_\la$ of all the leaves $\la\in \Lambda(\vartheta)$, 
and changing all the remaining labels consistently, in such a way that the labels $j_v$ and $\s_v$ remain unchanged for all $v\in N_2(\vartheta)$.
By construction, one has $\vartheta'\in \calmX(\vartheta)$ if and only if $\calmX(\vartheta') = \calmX(\vartheta)$; therefore
we may define the equivalence relation
\[
\vartheta' \sim \vartheta \quad : \Longleftrightarrow \quad \calmX(\vartheta') = \calmX(\vartheta) 
\]
and set $\Xi^{(k)}_{j,\nu,\s}:= \gotT^{(k)}_{j,\nu,\s}/\!\sim$ and $\breve\Xi^{(k)}_{j,\nu,\s}:= \breve\gotT^{(k)}_{j,\nu,\s}/\!\sim$.
By taking into account that
\begin{itemize}[topsep=.3em]
\itemsep0em
\item the number of unlabelled trees with $N$ nodes is bounded
by $4^{6N}$ (see Remark \ref{stoper}),
\item if we fix the component labels of the leaves, the number of trees which
are obtained by assigning the sign labels is bounded by $2^L$, if $L$ is the number of leaves,
\end{itemize}
we see that the sets $\Xi^{(k)}_{j,\nu,\s}$ and $\breve\Xi^{(k)}_{j,\nu,\s}$ are finite partitions of $\gotT^{(k)}_{j,\nu,\s}$ and $\breve\gotT^{(k)}_{j,\nu,\s}$,
respectively. We call $X$ any element of the sets $\Xi^{(k)}_{j,\nu,\s}$ and $\breve\Xi^{(k)}_{j,\nu,\s}$,
and $\vartheta^X$ an arbitrary representative of $X$.
Note that, by construction, $X$ is itself a set of trees.

%%%%%%%%%%%%%%%%%%%%%%%%%%%%%%%%%%%%%%%%%%%%%%%%%%%%%%%%%%%%%%%%%%%%%%%%%%
\begin{rmk} \label{rmk:disposta}
\emph{
With the notations of \eqref{disposta}, given a tree $\vartheta \in\gotT$,
any tree in $\calmX(\vartheta)$ fixes a  trimmed tree $\vartheta' \in \calmX(\breve\vartheta)$
and a trimmed tree $\vartheta_v'\in \calmX(\breve\vartheta_v)$ for each $v\in N_2(\vartheta)$.
Moreover $\calmX(\vartheta)$ is obtained by taking all possible trimmed trees $\vartheta'\in \calmX(\breve\vartheta)$ and grafting
any possible trimmed tree $\vartheta_v' \in \calmX(\breve\vartheta_v)$ to any node $v\in N_2(\vartheta)$.
Therefore, any $X\in\Xi^{(k)}_{j,\nu,\s}$ determines uniquely the set $\breve X$ and the $|N_2(\vartheta^X)|$ sets $\breve X_v$
such that $\breve X=\calmX(\breve\vartheta^X)$ and $\breve X_v=\calmX(\breve\vartheta^{X}_{v})$ for each $v\in N_2(\vartheta^X)$.
}
\end{rmk}
%%%%%%%%%%%%%%%%%%%%%%%%%%%%%%%%%%%%%%%%%%%%%%%%%%%%%%%%%%%%%%%%%%%%%%%%%%

%%%%%%%%%%%%%%%%%%%%%%%%%%%%%%%%%%%%%%%%%%%%%%%%%%%%%%%%%%%%%%%%%%%%%%%%%%
\begin{lemma}\label{convergerebbe2}
%Assume there is $C>0$ such  $|\h^{(k)}_j|\le C^k$ for all $k,j$.Then f
For all $s''>0$ there is a positive constant $D_{2}=D_{2}(s'',\al)$ such that,
for any $k\ge 1$, $j\in\ZZZ$, $\nu\in\ZZZ^\ZZZ_f$ and $\s\in\{\pm\}$ 
and any trimmed tree $\vartheta\in \breve\gotT^{(k)}_{j,\nu,\s}$, 
one has
\begin{equation}\label{atutti}
%\sum_{\vartheta\in\breve\gotT^{(k)}_{j,\nu,\s}} e^{- s'' J(\vartheta')} \le D_2^k .
\sum_{\vartheta' \in \calmX(\vartheta)} e^{- s'' J(\vartheta')} \le D_2^k .
\end{equation}
\end{lemma}
%%%%%%%%%%%%%%%%%%%%%%%%%%%%%%%%%%%%%%%%%%%%%%%%%%%%%%%%%%%%%%%%%%%%%%%%%%

%%%%%%%%%%%%%%%%%%%%%%%%%%%%%%%%%%%%%%%%%%%%%%%%%%%%%%%%%%%%%%%%%%%%%%%%%%
\prova
Decompose $\calmX(\vartheta)=\calmX_1(\vartheta)\sqcup \calmX_2(\vartheta)$, such that
there exists $ \la \in \Lambda(\vartheta')$  with $\jap{j_\la} \ge \jap{j}$ if $\vartheta'\in \calmX_1(\vartheta)$,
while one has $\jap{j_\la} < \jap{j}$ $\forall \la \in \Lambda(\vartheta')$ if $\vartheta'\in \calmX_2(\vartheta)$.
Then we write
\[
\sum_{\vartheta'\in \calmX(\vartheta)} e^{-s''J(\vartheta')}  =
\sum_{\vartheta'\in \calmX_1(\vartheta)} e^{-s''J(\vartheta')} + 
\sum_{\vartheta'\in \calmX_2(\vartheta)} e^{-s''J(\vartheta')} ,
\]
and study the two sums separately:
\vspace{-.2cm}
\begin{itemize}
\itemsep0em
\item by using %the definition of $J(\vartheta)$ in \eqref{J},
Remark \ref{stoper}, we bound
\[
\begin{aligned}
& \sum_{\vartheta'\in \calmX_1(\vartheta)} e^{-s''J(\vartheta')} =
\sum_{\vartheta'\in \calmX_1(\vartheta)} e^{s''\jap{j}^\al} \prod_{\la \in \Lambda(\vartheta')} e^{-s'' \jap{j_\la}^{\al}} \\
& \qquad \qquad \le 
\sum_{\la \in \Lambda(\vartheta)} \sum_{\substack{  j_{\la} \in \ZZZ \\ \jap{j_{\la}} \ge \jap{j}}} e^{-s'' (\jap{j_{\la}}^{\al} - \jap{j}^\al)}
\prod_{\la'\in\Lambda(\vartheta)\setminus\{\la\}} \sum_{ j_\la' \in \ZZZ} e^{-s'' \jap{j_{\la'}}^{\al}} \le 5 k \, C^k ,
\end{aligned}
\]
for some positive constant $C$ depending on $\al$ and $s''$;
\item by using Lemma \ref{meglio} and taking into account that $\gotn_1(\vartheta')=|j|$ and that the conservation law \eqref{lefoglie}
fixes $\gotn_2(\vartheta')$ in terms of $\gotn_1(\vartheta'), \gotn_3(\vartheta'),\ldots,\gotn_p(\vartheta')$, with $p=|\Lambda(\vartheta)|+1$, we obtain
\[
\begin{aligned}
\sum_{\vartheta'\in \calmX_2(\vartheta)} e^{-s''J(\vartheta')}  
& \le % \le %\stackrel{\eqref{constance}}{\le}
\sum_{\vartheta'\in \calmX_2(\vartheta)} e^{-s''(2-2^\al) ( \gotn_3^\al(\vartheta') + \ldots + \gotn_p^\al(\vartheta')) } \\
& \le \sum_{n_3 , \dots ,n_{p} \ge 1} e^{-s''(2-2^\al) ( n_3^\al + \ldots + n_p^\al) }\le C^k ,
\end{aligned}
\]
for some positive constant $C$ depending on $\al$ and $s''$.
\end{itemize}
Therefore, in all cases, we obtain \eqref{atutti}
for a suitable constant $D_2=D_2(s'',\al)$. % depending on $\al$ and $s''$.
\EP
%%%%%%%%%%%%%%%%%%%%%%%%%%%%%%%%%%%%%%%%%%%%%%%%%%%%%%%%%%%%%%%%%%%%%%%%%%

The next result is crucial in order to compare $|\om\cdot\nu-\om_j|$, i.e.~the modulus of the small divisor,
with $|\nu-\gote_j|_{\al/2}$; it is a version of \cite[Lemma 7.1]{BMP2} adapted to the present situation.

%%%%%%%%%%%%%%%%%%%%%%%%%%%%%%%%%%%%%%%%%%%%%%%%%%%%%%%%%%%%%%%%%%%%%%%%%%
\begin{lemma}\label{constance2}
There exists a constant $C_1=C_1(\al)$, such that the following holds. 
Consider a tree $\vartheta\in %\Theta^{(k)}_{j,\nu,\s}\cup
\breve\gotT^{(k)}_{j,\nu,\s}$, for some $k\ge 1$, $j\in\ZZZ$, $\nu\in\ZZZ^\ZZZ_f$ and $\s\in\{\pm\}$.
If $0<|\om\cdot \nu -\om_j| <1/2$, then
\begin{equation}\label{zia}
 J(\vartheta) \ge  C_1|\nu -\gote_j |_{\al/2} .
\end{equation}
One can take $C_1(\al):=(2-2^\al)/(2\cdot 3^{\al/2}+1)$.
\end{lemma}
%%%%%%%%%%%%%%%%%%%%%%%%%%%%%%%%%%%%%%%%%%%%%%%%%%%%%%%%%%%%%%%%%%%%%%%%%%

%%%%%%%%%%%%%%%%%%%%%%%%%%%%%%%%%%%%%%%%%%%%%%%%%%%%%%%%%%%%%%%%%%%%%%%%%%
\prova
%Given $\vartheta$, set $\widehat\Lambda(\vartheta):=\Lambda(\vartheta)\cup\{r\}$ and $p:=|\widehat\Lambda(\vartheta)|$. 
Define 
\[
\hat{\s}_\la:=\begin{cases}
\s_\la, &\quad \la\in\Lambda^*(\vartheta),\\
-\s, &\quad \la=r,
\end{cases}
\]
and 
\[
\gotP:=\{(j_\la,\hat{\s}_\la)\}_{\la\in\widehat\Lambda(\vartheta)}.
\]
where $\Lambda^*(\vartheta)$ and $\widehat\Lambda(\vartheta)$ are defined in \eqref{Lambdastar} and at the beginning
of Subsection \ref{vedo}, respectively.

Next, consider the set
$\widehat\gotP\subseteq\gotP$ in which all the elements with $j_\la=0$ and all the pairs $(j_\la,\hat{\s}_{\la}),(j_{\la'},\hat{\s}_{\la'})$
with $j_\la=j_{\la'}$ and $\hat{\s}_\la=-\hat{\s}_{\la'}$ are removed;
in other words one has $\widehat\gotP=\{(j_{1},\hat\s_{1}),\ldots , (j_{\hat{p}},\hat\s_{\hat{p}})\}$, for some $\hat p \le p=p(\vartheta)$,
such that $j_i\ne0$ and if $j_i=j_{i'}$ then $\hat\s_i= \hat\s_{i'}$.
Up to rearranging the indices we may assume $|j_1|\ge |j_2|\ge \ldots \ge |j_{\hat{p}}|\ge1$.
By considering \eqref{lefoglie*}, with $\ell$ the root line of $\vartheta$, and, first subtracting $\s j$ from both sides
the third equation and then $\s\gote_j$ from both sides of the second one, we obtain, respectively,
 \begin{equation}\label{freddo}
 \sum_{i=1}^{\hat{p}} \hat\s_i j_i =0\qquad \mbox{ and }\qquad
 \s(\nu - \gote_j) = \ka_\vartheta \, \gote_0 + \sum_{i=1}^{\hat{p}} \hat\s_i \gote_{j_i}\,,
 \end{equation}
 where
 \begin{equation}\label{kappa}
 \ka_\vartheta := \sum_{\substack{\la\in\widehat\Lambda(\vartheta) \\ \; j_\la=0}} \hat{\s}_{\la}1\,.
 \end{equation}

 Now, by \eqref{rettangolo} we can write $\om_j=j^2+\ze_j$ with $|\ze_j|\le\frac{1}{2}$, so that 
 \begin{equation}\label{nomino}
 \begin{aligned}
|\om \cdot \nu -\om_j| &=|\om\cdot(\nu-\gote_j)| \stackrel{\eqref{freddo}}{=}|\ka_\vartheta \ze_0 + \sum_{i=1}^{\hat{p}}\hat{\s}_i(j^2_i + \ze_{j_i})| \\ &
\ge \left| \sum_{i=1}^{\hat{p}} \hat{\s}_i j_i^2 \right| - \left| \sum_{i\in \ZZZ} \ze_i \left( \nu_i - \de_{ij} \right) \right|
\ge \left| \sum_{i=1}^{\hat{p}} \hat{\s}_i j_i^2 \right| - \frac{1}{2} {\|\nu-\gote_j \|_{1}} %|\nu-\gote_j|_{\ell^1}
\end{aligned}
 \end{equation}
 and hence, since we assumed $|\om\cdot\nu-\om_j|<1/2$  and the right hand side of \eqref{nomino}
 is the difference between an integer and a half-integer, we must have
 \[
 \left|\sum_{i=1}^{\hat{p}}\hat{\s}_i j_i^2 \right|\le {\|\nu-\gote_j \|_{1} \le {p} }. %|\nu-\gote_j|_{\ell^1} \le {p}.
 \]

Note that $\hat{p}\ge3$. Indeed $\hat{p}\ne0$ since we assumed $0<|\om\cdot\nu - \om_j|<1/2$.
If $\hat{p}=1$ then  the first of \eqref{freddo} would imply $(j_1,\s_1)=(0,\s_1)$
but, by definition of $\widehat{\gotP}$, $j_1\ne0$. Finally, if $\hat{p}=2$ again the first of \eqref{freddo}
would imply $(j_2,\s_2)=(j_1,-\s_1)$ and again this cannot happen for pairs in $\hat{\gotP}$. Hence
we can apply Lemma \ref{costanza} to deduce that 
\begin{equation}\label{riunione}
|j_2|\le |j_1| \le 3\sum_{i=3}^{\hat{p}} |j_i|^2 +p-\hat{p} \,.
\end{equation}

Therefore we have
\begin{equation}\label{metro}
\begin{aligned}
|\nu - \gote_j|_{\al/2} & =  \sum_{i\in\ZZZ} \jap{i}^{\al/2}|\nu_i - \de_{ij}| 
 \stackrel{\eqref{freddo}}{=} \ka_\vartheta +  \sum_{i\in\ZZZ\setminus\{0\}} {|i|}^{\al/2}|\nu_i - \de_{ij}|  \\
&\stackrel{\eqref{freddo}}{\le } \ka_\vartheta + \sum_{i=1}^{\hat{p}} |j_i|^{\al/2}
\le \ka_\vartheta + 2|j_1|^{\al/2} + \sum_{i=3}^{\hat{p}} |j_i|^{\al/2} \\
&\stackrel{\eqref{riunione}}{\le} \ka_\vartheta + 2 \left( 3\sum_{i=3}^{\hat{p}} |j_i|^2 + p-\hat{p}\right)^{\al/2} 
 + \sum_{i=3}^{\hat{p}} |j_i|^{\al/2}  \\
 &\le \ka_\vartheta + (2\cdot 3^{\al/2}+1)\sum_{i=3}^{\hat{p}}|j_i|^\al + 2(p-\hat{p}) \\
& \le  (2\cdot 3^{\al/2}+1)\sum_{i=3}^{\hat{p}}|j_i|^\al + 3(p- \hat{p}) \\
 & \le (2\cdot 3^{\al/2}+1) \Big(\sum_{i=3}^{\hat{p}} |j_i|^\al  + \sum_{i=\hat{p}+1}^{p}1\Big)
 \le (2\cdot 3^{\al/2}+1) \sum_{i=3}^{{p}}\gotn_i^\al(\vartheta)
\end{aligned}
\end{equation}
 where we are using the definition \eqref{glienne}.
Hence the assertion follows from Lemma \ref{meglio}, with $C_1(\al):=(2-2^\al)/(2\cdot 3^{\al/2}+1)$.
\EP
%%%%%%%%%%%%%%%%%%%%%%%%%%%%%%%%%%%%%%%%%%%%%%%%%%%%%%%%%%%%%%%%%%%%%%%%%%

%%%%%%%%%%%%%%%%%%%%%%%%%%%%%%%%%%%%%%%%%%%%%%%%%%%%%%%%%%%%%%%%%%%%%%%%%%
\begin{defi}[\textbf{Support property}]
\label{supporto}
Given $\vartheta \in %\Theta\cup
\gotT\cup\breve\gotT$,
we say that a line $\ell\in L(\vartheta)$ verifies the \emph{support property} if either $n_\ell\le0$ or %$n_\ell\ge 1$ and 
\begin{equation}\label{tagliasup}
%\frac{1}{16}\be_{\om}(r_{m_{n_\ell}}) \le |\om\cdot\nu_\ell - \om_{j_\ell}| < \frac{1}{8}\be_{\om}(r_{m_{n_\ell-1}}),
\frac{1}{32}\be({m_{n_\ell}}) < |x_\ell| < \frac{1}{4}\be({m_{n_\ell-1}}) .
\end{equation}
%where $\be(m_{-1})$ has to be interpreted as $+\io$. 
We say that $\vartheta$
verifies the support property if all $\ell\in L(\vartheta)$ verify the support property.
\end{defi}
%%%%%%%%%%%%%%%%%%%%%%%%%%%%%%%%%%%%%%%%%%%%%%%%%%%%%%%%%%%%%%%%%%%%%%%%%%

%%%%%%%%%%%%%%%%%%%%%%%%%%%%%%%%%%%%%%%%%%%%%%%%%%%%%%%%%%%%%%%%%%%%%%%%%%
\begin{rmk}\label{vabbuo}
\emph{
Given $\vartheta\in\gotT\cup\breve\gotT$, if $\Val(\vartheta;c,\om)\ne0$ then $\vartheta$ automatically verifies the support property
(see Remark \ref{lescale}). We introduce the support property
because later we shall need to associate with each line $\ell\in L(\vartheta)$ an \emph{interpolated momentum},
to be denoted $\nu_\ell(\und{y_\ell})$, and we shall see that,
if $n_\ell$ is such that $\calG_{{n_\ell}}(\om\cdot\nu_\ell(\und{y_\ell})-\om_{j_\ell})\ne0$,
then the momentum $\nu_\ell$ is such that $x_\ell=\om\cdot\nu_{\ell} - \om_{j_\ell}$ verifies the bound \eqref{tagliasup}.
}
\end{rmk}
%%%%%%%%%%%%%%%%%%%%%%%%%%%%%%%%%%%%%%%%%%%%%%%%%%%%%%%%%%%%%%%%%%%%%%%%%%

%%%%%%%%%%%%%%%%%%%%%%%%%%%%%%%%%%%%%%%%%%%%%%%%%%%%%%%%%%%%%%%%%%%%%%%%%%
\begin{lemma}\label{quasiB}
Consider $\vartheta\in %\Theta\cup
\gotT\cup\breve\gotT$,
and assume that $\vartheta$ verifies the support property.
If  $\ell \in L(\vartheta)$ has scale $n \ge 1$, then
\[
|\nu_\ell-\gote_{j_\ell}|_{\al/2} \ge r_{m_n-1}
%\ge r_{m_{n-1}}.
\]
\end{lemma}
%%%%%%%%%%%%%%%%%%%%%%%%%%%%%%%%%%%%%%%%%%%%%%%%%%%%%%%%%%%%%%%%%%%%%%%%%%

%%%%%%%%%%%%%%%%%%%%%%%%%%%%%%%%%%%%%%%%%%%%%%%%%%%%%%%%%%%%%%%%%%%%%%%%%%
\prova
If $\ell$ has scale $n\ge1$ we have
\[
%|\om\cdot\nu-\om_j|\le \frac{1}{8}\be_{\om}(r_{m_{n-1}})\le\frac{1}{4}\be_{\om}(r_{m_{n}-1})<\be_{\om}(r_{m_{n}-1})
|x_{\ell}|< \frac{1}{4}\be({m_{n-1}})\le\frac{1}{2}\be({m_{n}-1})<\be({m_{n}-1}),
\]
where we used the definition of the sequence ${m_n}$,
so that the result follows from the definition of $\be_{\om}$ in \eqref{beta}.
\EP
%%%%%%%%%%%%%%%%%%%%%%%%%%%%%%%%%%%%%%%%%%%%%%%%%%%%%%%%%%%%%%%%%%%%%%%%%%

%%%%%%%%%%%%%%%%%%%%%%%%%%%%%%%%%%%%%%%%%%%%%%%%%%%%%%%%%%%%%%%%%%%%%%%%%%
\begin{coro}\label{basta!}
Consider $\vartheta\in %\Theta\cup
\breve\gotT$,
and assume that $\vartheta$ verifies the support property.
If the root line has scale $n\ge1$, then
\[
J(\vartheta) \ge C_1 r_{m_n-1},
\]
where $C_1$ is the constant appearing in Lemma \ref{constance2}.
\end{coro}
%%%%%%%%%%%%%%%%%%%%%%%%%%%%%%%%%%%%%%%%%%%%%%%%%%%%%%%%%%%%%%%%%%%%%%%%%%

%%%%%%%%%%%%%%%%%%%%%%%%%%%%%%%%%%%%%%%%%%%%%%%%%%%%%%%%%%%%%%%%%%%%%%%%%%
\prova 
Combine Lemma \ref{constance2} with Lemma \ref{quasiB}, with $\ell$ being the root line.
\EP

%%%%%%%%%%%%%%%%%%%%%%%%%%%%%%%%%%%%%%%%%%%%%%%%%%%%%%%%%%%%%%%%%%%%%%%%%%
\begin{lemma}\label{consec}
Let $T$ be a SEG of a tree $\vartheta\in %\Theta\cup
\gotT\cup\breve\gotT$, with entering line $\ell'$ and exiting line $\ell$, such that $\calP_T$ does not contain any $\h$-line.
Assume  that $0<|x_{\ell}|<1/2$, $0<|x_{\ell'}|<1/2$ and $\s_\ell x_\ell \ne\s_{\ell'} x_{\ell'}$. Then 
\begin{equation}\label{lezione}
J(T) \ge C_1 |\s_{\ell}(\nu_{\ell} -\gote_{j_\ell}) - \s_{\ell'}(\nu_{\ell'} - \gote_{j_\ell'})|_{\al/2} 
\end{equation}
%
%for some positive constant $C_2=C_2(\al)$.
with the constant $C_1$ as in Lemma \ref{constance2}.
\end{lemma}
%%%%%%%%%%%%%%%%%%%%%%%%%%%%%%%%%%%%%%%%%%%%%%%%%%%%%%%%%%%%%%%%%%%%%%%%%%

%%%%%%%%%%%%%%%%%%%%%%%%%%%%%%%%%%%%%%%%%%%%%%%%%%%%%%%%%%%%%%%%%%%%%%%%%%
\prova
Call $v$ the node $\ell$ exits and $v'$ the node\footnote{The line $\ell'$ cannot exit a leaf since $|x_{\ell'}|>0$.} $\ell'$ exits,
so that  $j_v=j_\ell$ and $j_{v'}=j_{\ell'}$.
%and associate with $v$  a label $j_v=j_\ell$ and with  $v'$ a label $j_{v'}=j_{\ell'}$.
Set $\widehat\Lambda(T):=\Lambda^*(T)\cup\{v,v'\}$, and define 
\[
\hat{\s}_\la:=\begin{cases}
\s_\la, &\quad \la\in\Lambda(T)\cup\{v'\},\\
-\s_\ell , &\quad \la=v,
\end{cases}
\]
and 
\[
\gotP:=\{(j_\la,\hat{\s}_\la)\}_{\la\in \widehat\Lambda(T)}.
\]

%%%
%%% set $\Lambda(\vartheta)=\{\la_{1},\ldots\la_p\}$ and consider the set
%%%$$
%%%%\gotP:=\{(j_{\la_1},\s_{\la_1}),\ldots (j_{\la_p},\s_{\la_p}),(j, - \s)\}
%%%\gotP:=\{(j_{1},\hat{\s}_{1}),\ldots (j_{p},\hat{\s}_{p}),(j_{p+1}, \hat{\s}_{p+1})\}
%%%$$ 
%%%where $j_i=j_{\la_i}$ and $\hat{\s}_i=\s_{\la_i}$ for $i=1,\ldots,p$, while $j_{p+1}=j$ and $\hat{\s}_{p+1}=-\s$.

Let
$\widehat\gotP=\{(j_{1},\hat\s_{1}),\ldots , (j_{\hat{p}},\hat\s_{\hat{p}})\}$ be the largest subset of $\gotP$ such that
$j_i\ne0$ and if $j_i=j_{i'}$ then
$\hat\s_i= \hat\s_{i'}$. We rearrange the indices so that
$|j_1|\ge |j_2|\ge \ldots \ge |j_{\hat{p}}|\ge1$. 
Note that
 \begin{equation}\label{freddoT}
 \sum_{i=1}^{\hat{p}} \hat \s_i j_i =0 , \qquad \qquad
 \s_\ell(\nu_\ell - \gote_{j_\ell}) - \s_{\ell'}(\nu_{\ell'}-\gote_{j_{\ell'}}) = \ka_T\, \gote_0 + 
 \sum_{i=1}^{\hat{p}} \hat \s_i \gote_{j_i}\,,
 \end{equation}
 where
 \begin{equation}\label{kappaT}
 \ka_T := \sum_{\substack{\la\in \widehat\Lambda(T)\\ j_\la=0}} \hat{\s}_{\la}1 .
 \end{equation}
We can repeat essentially word by word the argument of the proof of Lemma
\ref{constance2}, the only difference being that now we use that $|x_\ell \pm x_{\ell'}|<1$.
This gives the same constant $C_1$ as in Lemma \ref{constance2}.
%This gives a constant $C_2=(2-2^\al)/(2\cdot 3^{\al/2}+2)$ instead of $C_1$. 
\EP
%%%%%%%%%%%%%%%%%%%%%%%%%%%%%%%%%%%%%%%%%%%%%%%%%%%%%%%%%%%%%%%%%%%%%%%%%%

%%%%%%%%%%%%%%%%%%%%%%%%%%%%%%%%%%%%%%%%%%%%%%%%%%%%%%%%%%%%%%%%%%%%%%%%%%
\begin{lemma}\label{ovvio2}
Let $\vartheta\in %\Theta\cup
\gotT\cup\breve\gotT$ verify the support property.
Let $T$ be a SEG of $\vartheta$ with entering line $\ell'$ and exiting  line $\ell$,
 such that $\s_{\ell}(\nu_{\ell} -\gote_{j_\ell}) \ne \s_{\ell'}(\nu_{\ell'} - \gote_{j_\ell'})$
and $\calP_T$ does not contain any $\h$-line. If ${n}:=\min\{n_\ell,n_{\ell'}\}\ge1$ then
\begin{equation}\label{sist}
 J(T) \ge C_1  |\s_{\ell}(\nu_{\ell} -\gote_{j_\ell}) - \s_{\ell'}(\nu_{\ell'} - \gote_{j_\ell'})|_{\al/2} \ge C_1 r_{m_n-1} ,
\end{equation}
where  $C_1$ is as in Lemma \ref{constance2}.
%where  $C_2$ is as in Lemma \ref{consec}.
\end{lemma}
%%%%%%%%%%%%%%%%%%%%%%%%%%%%%%%%%%%%%%%%%%%%%%%%%%%%%%%%%%%%%%%%%%%%%%%%%%

%%%%%%%%%%%%%%%%%%%%%%%%%%%%%%%%%%%%%%%%%%%%%%%%%%%%%%%%%%%%%%%%%%%%%%%%%%
\prova
By Remark \ref{lescale} one has
 \[
|x_\ell-x_{\ell'}|\le
|x_\ell | + |x_{\ell'}| 
< \frac{1}{4}\be({m_{n_\ell-1}}) + 
\frac{1}{4}\be({m_{n_{\ell'}-1}}) 
%&\stackrel{n=\min\{n_\ell,n_{\ell'}\}}
\le  \frac{1}{2}\be({m_{n-1}}) \le 
\be({m_{n}-1}),
\]
 so that one has
 \[
|\s_{\ell}(\nu_{\ell} -\gote_{j_\ell}) - \s_{\ell'}(\nu_{\ell'} - \gote_{j_\ell'})|_{\al/2} \ge 
r_{m_n-1}
\]
from the definition of $\be_\om$ in \eqref{beta}, and hence \eqref{sist} follows from Lemma \ref{consec}.
\EP
%%%%%%%%%%%%%%%%%%%%%%%%%%%%%%%%%%%%%%%%%%%%%%%%%%%%%%%%%%%%%%%%%%%%%%%%%%
%
%%%%%%%%%%%%%%%%%%%%%%%%%%%%%%%%%%%%%%%%%%%%%%%%%%%%%%%%%%%%%%%%%%%%%%%%%%%
%\begin{lemma}\label{mortacci}
%Let $\vartheta\in\Theta\cup\breve\gotT$ verify the support property.
%Let $T$ be a subgraph of $\vartheta$ with only one entering line $\ell'$ and one exiting 
%line $\ell$  such that % $\s_{\ell}(\nu_{\ell} -\gote_{j_\ell}) =\s_{\ell'}(\nu_{\ell'} - \gote_{j_\ell'})$ and 
%$\calP_T$ does contain an $\h$-line.
%If $n:=n_{\ell'}\ge1$ then
%\begin{equation} \nonumber %\label{sist}
% J(T) \ge C_1 r_{m_n-1},
%\end{equation}
%where  $C_1$ is as in Lemma \ref{constance2}.
%%where  $C_2$ is as in Lemma \ref{consec}.
%\end{lemma}
%%%%%%%%%%%%%%%%%%%%%%%%%%%%%%%%%%%%%%%%%%%%%%%%%%%%%%%%%%%%%%%%%%%%%%%%%%%
%
%%%%%%%%%%%%%%%%%%%%%%%%%%%%%%%%%%%%%%%%%%%%%%%%%%%%%%%%%%%%%%%%%%%%%%%%%%%
%\prova
%Let $\ell_1$ be the $\h$-line in $\calP_T$ closest to $\ell'$, and let 
%$T'$ be the subgraph of $\vartheta$ having $\ell_1$ as exiting line and $\ell'$ as entering line. Of course $J(T)\ge J(T')$.
%Then $\calP_{T'}$ does not contain any $\h$-line, 
%so we can apply Lemma \ref{ovvio2}.
%\EP
%%%%%%%%%%%%%%%%%%%%%%%%%%%%%%%%%%%%%%%%%%%%%%%%%%%%%%%%%%%%%%%%%%%%%%%%%%%

%%%%%%%%%%%%%%%%%%%%%%%%%%%%%%%%%%%%%%%%%%%%%%%%%%%%%%%%%%%%%%%%%%%%%%%%%%
\begin{defi}[\textbf{Minimum scale}]
\label{scalemin}
For any  $\vartheta\in %\Theta\cup
\gotT\cup\breve\gotT$ and any line $\ell\in L(\vartheta)$ with $\nu_{\ell}\ne \s_\ell\gote_{j_\ell}$,
we define the \emph{minimum scale} as
\begin{equation}\label{minnie}
\nmin_\ell:= \min\{ n\ge0\;:\; \Psi_n(x_\ell )\ne0\} .
\end{equation}
%and the \emph{maximum scale} as
%\begin{equation}\label{topolino}
%\nmax_\ell:= \max\{ n\ge0\;:\; \Psi_n(\om\cdot\nu_\ell-\om_{j_\ell})\ne0\}.
%\end{equation}
\end{defi}
%
%\begin{rmk}\label{stimameglio}
%One easily checks that if $\nmin_\ell=\nmax_\ell$ then
%\begin{equation}\label{tagliameglio}
%\frac{1}{8}\be({m_{\nmin_\ell}})\le  |\om\cdot\nu_\ell - \om_{j_\ell}| \le \frac{1}{16}\be(m_{\nmin_\ell-1}), 
%\end{equation}
%while if $\nmin_\ell<\nmax_\ell$ then
%\begin{equation}\label{tagliameglio2}
%\frac{1}{16}\be({m_{\nmin_\ell}})\le  |\om\cdot\nu_\ell - \om_{j_\ell}| \le \frac{1}{8}\be(m_{\nmin_\ell})
%\end{equation}
%\end{rmk}
%%%%%%%%%%%%%%%%%%%%%%%%%%%%%%%%%%%%%%%%%%%%%%%%%%%%%%%%%%%%%%%%%%%%%%%%%%

Call $L_{\!N\!R}(\vartheta)$ the set of non-resonant lines in $\vartheta$ on scale $\ge0$, and define, for $n\ge 0$,
\begin{equation}\label{Nn}
\gotN_n(\vartheta):=|\{\ell\in L_{\!N\!R}(\vartheta) : \nmin_\ell \ge n \}|\,.
\end{equation}
%where the minimum scale $\nmin_\ell$ is defined in \eqref{minnie}.

Now, we are in the position to prove the following version of the {finite-dimensional Bryuno lemma}, which gives
a bound on the non-resonant lines.

%%%%%%%%%%%%%%%%%%%%%%%%%%%%%%%%%%%%%%%%%%%%%%%%%%%%%%%%%%%%%%%%%%%%%%%%%%
\begin{prop}\label{brjuno}
Let $\vartheta\in %\Theta \cup
\breve\gotT$ %\Theta_{j,\nu,\s}^{(k)}\cup \gotT_{j,\nu,\s}^{(k)}$ for some $k,\nu,j,\s$
verify the support property.
For %any tree $\vartheta$ and 
all $n\ge1$ %\in\Theta_{\nu,j,\s}^{(k)}$ 
one has
% there is a constant $K>0$ such that
 \begin{equation}\label{bru}
 \gotN_n(\vartheta)\le \max\Biggl\{\frac{K}{r_{m_n-1}} J(\vartheta)-2,0\Biggr\},
 \end{equation}
 with $K:= 4/C_1$ and $C_1$ as in Lemma \ref{constance2}. %\max\{4/C_1,2/C_2\}$.
\end{prop}
%%%%%%%%%%%%%%%%%%%%%%%%%%%%%%%%%%%%%%%%%%%%%%%%%%%%%%%%%%%%%%%%%%%%%%%%%%

%%%%%%%%%%%%%%%%%%%%%%%%%%%%%%%%%%%%%%%%%%%%%%%%%%%%%%%%%%%%%%%%%%%%%%%%%%
\prova
The proof is by induction on the number of nodes $k$ of the tree. %$k$.

For $k=1$ only the root line may have scale $n\ge1$, so that $\gotN_n(\vartheta)\le1$,
and the bound follows immediately from Corollary \ref{basta!} since $K > 3/C_1$.

If $k>1$, assume \eqref{bru} to be true for all $k'<k$ and 
consider a tree $\vartheta\in %\Theta_{j,\nu,\s}^{(k)}\cup 
\breve\gotT_{j,\nu,\s}^{(k)}$, for some 
$j,\nu,\s$. Let $\ell$ be the root line of $\vartheta$. If $\nmin_\ell < n$
then the bound is satisfied by induction, using \eqref{i}.  If $\nmin_{\ell}\ge n$, let $\ell_1,\ldots,\ell_p$
be the lines closest to $\ell$ such that $\nmin_{\ell_i}\ge n$ and
call $\vartheta_1,\ldots,\vartheta_p$ the subtrees of $\vartheta$ having $\ell_1,\ldots,\ell_p$ as root lines, respectively.
%We distinguish among three cases.

%\begin{itemize}[topsep=0ex]
%\itemsep0em

%\item[1.] 

If $p=0$ the result follows once more since $\gotN_n(\vartheta)=1$ and $K > 3/C_1$.
%Indeed, by Lemmata \ref{constance2} and \ref{quasiB}, if in $\vartheta$ there is any line
%on minimum scale $\ge n$, then one has $ J(\vartheta) \ge C_1 r_{m_n-1}$.

If $p\ge2$ then
\[
\gotN_n(\vartheta)= 1 + \sum_{i=1}^p \gotN_n(\vartheta_i) \le 1+ \frac{K}{r_{m_n-1}} \sum_{i=1}^p J(\vartheta_i) -2p
\le \frac{K}{ r_{m_n-1}} J(\vartheta)-2 ,
\]
where we have used the inductive hypothesis in the first inequality and in the second one we have exploited \eqref{i} and the fact that $p\ge2$.

%\item[3.] 

If $p=1$, call $T$ the SEG of $\vartheta$ having $\ell$ as exiting line and $\ell_1$ as entering line. We have,
by the inductive hypothesis,
\[
\gotN_n(\vartheta) = 1 + \gotN_n(\vartheta_1) \le \frac{K}{r_{m_n-1}} J(\vartheta_1) -1 \le
 \frac{K}{r_{m_n-1}}(J(\vartheta)-J(T))-1.
\]
where we have used that $J(\vartheta)=J(T)+J(\vartheta_1)$.
Since $K > 1/C_1$, the bound follows if $J(T) \ge C_1 r_{m_n-1}$, which holds necessarily in the following cases:

\begin{itemize}[topsep=0ex]
\itemsep0em
\item if $\s_{\ell}(\nu_{\ell} -\gote_{j_\ell}) \ne \s_{\ell_1}(\nu_{\ell_1} - \gote_{j_{\ell_1}})$,
by Lemma \ref{ovvio2};
\item if $\s_{\ell}(\nu_{\ell} - \gote_{j_\ell}) = \s_{\ell_1}(\nu_{\ell_1} - \gote_{j_{\ell_1}})$
and there is a line $\bar{\ell}\in L(T)$ such that $\nmin_{\bar{\ell}}+1=n_{\bar\ell}=\nmin_{\ell_1}=n$,
because in such a case one has
$\s_{\bar\ell}(\nu_{\bar\ell} -\gote_{j_{\bar\ell}}) \neq \s_{\ell_1}(\nu_{\ell_1} - \gote_{j_{\ell_1}})$, so that
if $\bar\ell\in\calP_T$ then \eqref{bru} follows again by Lemma \ref{ovvio2}, while if $\bar\ell \notin \calP_T$
then \eqref{bru} follows by Corollary \ref{basta!} and Remark \ref{utile}.
\end{itemize}

On the other hand, if $\s_{\ell}(\nu_{\ell} -\gote_{j_\ell}) = \s_{\ell_1}(\nu_{\ell_1} - \gote_{j_{\ell_1}})$,
$n_{\bar\ell} < n$ for all lines $\bar\ell\in L(T)$ and $J(T) <  C_1 r_{m_{n}-1}$,
then one has $J(T) < C_1r_{m_{n}-1}  \le C_1 r_{m_{\nmin_T}-1}$, so that $T$ is a RC.
Let $\ell'_1,\ldots,\ell'_{p'}$ the  lines closest to $\ell_1$ such that $\nmin_{\ell'_i}\ge n$ and
call $\vartheta'_1,\ldots,\vartheta_{p'}'$ the subtrees of $\vartheta$ having $\ell'_1,\ldots,\ell_{p'}'$ as root lines, respectively.

%\begin{enumerate}[topsep=0ex]

%\item[3.2.1] 
If $p'=0$ then $\gotN_n(\vartheta)=2$ and hence the bound follows since $K\ge 4/C_1$.

%\item[3.2.2] 
If $p'\ge2$ we have
\[
\gotN_n(\vartheta)\le 2 + \sum_{i=1}^{p'} \gotN_n(\vartheta'_i) \le 2 + \frac{K}{r_{m_n-1}} \sum_{i=1}^{p'} J(\vartheta'_i) -2p'
\le \frac{K}{ r_{m_n-1}} J(\vartheta)-2.
\]

%\item[3.2.3.] 
If $p'=1$, call $T'$ the SEG of $\vartheta$ having $\ell_1$ as exiting line and $\ell_1'$ as entering line. One has
\[
\gotN_n(\vartheta) \le 2+  \gotN_n(\vartheta_1') \le \frac{K}{r_{m_n-1}} J(\vartheta_1') \le
 \frac{K}{r_{m_n-1}}(J(\vartheta)-J(T')),
\]
because one has $J(\vartheta_1)=J(T')+J(\vartheta_1')$.

Once more the bound \eqref{bru} is obtained whenever one has $J(T') > C_1 r_{m_n-1}$, which, by reasoning as before
and using that $K>2/C_1$, is found to hold true automatically in the following cases:
\begin{itemize}[topsep=0ex]
\itemsep0em
%\item if there are $\h$-lines in $\calP_{T'}$;
\item if %$\calP_{T'}$ does not contain any $\h$-line and 
$\s_{\ell_1}(\nu_{\ell_1} -\gote_{j_{\ell_1}}) \neq \s_{\ell'_1}(\nu_{\ell'_1} - \gote_{j_{\ell'_1}})$;
\item if %$\calP_{T'}$ does not contain any $\h$-line, 
$\s_{\ell_1}(\nu_{\ell_1} -\gote_{j_{\ell_1}}) = \s_{\ell'_1}(\nu_{\ell'_1} - \gote_{j_{\ell'_1}})$
and there is a line $\bar{\ell}\in L(T)$ with $\nmin_{\bar{\ell}}+1=n_{\bar \ell}=\nmin_{\ell'_1}=n$.
% then $\s_{\bar\ell}(\nu_{\bar\ell} -\gote_{j_{\bar\ell}}) - \s_{\ell'_1}(\nu_{\ell'_1} - \gote_{j_{\ell'_1}})\ne0$; if $\bar\ell\notin\calP_T$
%then \eqref{bru} follows by combining Lemmata \ref{constance2} and \ref{quasiB} since $K>1/C_1$, while if
%$\bar\ell\in\calP_T$ then \eqref{bru} follows again from Lemma \ref{ovvio2}, since $K>1/C_2$.
\end{itemize}

Thus, if $J(T') < C_1 r_{m_n-1}$, then %$\calP_T'$ does not contain any $\h$-line, 
$\s_{\ell_1}(\nu_{\ell_1} -\gote_{j_{\ell_1}}) =
\s_{\ell'_1}(\nu_{\ell'_1} - \gote_{j_{\ell'_1}})$  and $n_{\bar{\ell}}<\nmin_{\ell'_1}$ for all $\bar{\ell}\in L(T)$,
so that, using that $C_1 r_{m_n-1} < C_1 r_{m_{\nmin_{T'}}-1}$, one finds that $T'$ is a RC as well as $T$, and hence $\ell_1$ is a resonant line.
Let $\gotC$ be the chain containing $T$ and $T'$, and let
$\ell''$ and $\vartheta''$ be the entering line of $\gotC$ and the subtree of $\vartheta$ which has $\ell''$ as root line, respectively.
Thus $\gotN_n(\vartheta) = 1 +  \gotN_n(\vartheta'')$. 
Let $\ell''_1,\ldots,\ell''_{p''}$
be the lines closest to $\ell''$ such that $\nmin_{\ell''_i}\ge n$ and
call $\vartheta''_1,\ldots,\vartheta_{p''}''$ the subtrees of $\vartheta$ having $\ell''_1,\ldots,\ell_{p''}$ as root lines, respectively.

%\begin{enumerate}[topsep=0ex]

%\item[3.2.3.2.1.]
 
If $p''=0$ then $\gotN_n(\vartheta)=2$ so that the bound follows since $K \ge 4/C_1$.

%\item[3.2.3.2.2.]
If $p''\ge2$ then the bound follows noting that
\[
\gotN_n(\vartheta)\le 2 + \sum_{i=1}^{p''} \gotN_n(\vartheta''_i) \le 2 + \frac{K}{r_{m_n-1}} \sum_{i=1}^{p''} J(\vartheta''_i) -2p''
\le \frac{K}{ r_{m_n-1}} J(\vartheta)-2.
\]

%\item[3.2.3.2.3.]
If $p''=1$, since the SEG $T''$ having $\ell''$ as exiting line and $\ell''_1$ as entering line is not a RC, then,
by reasoning as above, one finds that either of the following conditions is met:
\begin{itemize}[topsep=0ex]
\itemsep0em
\item $\s_{\ell''}(\nu_{\ell''} -\gote_{j_{\ell''}}) \ne \s_{\ell''_1}(\nu_{\ell''_1} - \gote_{j_{\ell''_1}})$, and hence one has $J(T'') > C_1 r_{m_{\nmin_{T''}}-1}$;
\item $\s_{\ell''}(\nu_{\ell''} -\gote_{j_{\ell''}}) = \s_{\ell''_1}(\nu_{\ell''_1} - \gote_{j_{\ell''_1}})$
and there is a line $\bar{\ell}\in\calP_{T''}$ with $\nmin_{\bar{\ell}}+1=n_{\bar\ell}=\nmin_{\ell'_1}=n$, and hence one has
$J(T'') > C_1 r_{m_{\nmin_{T''}}-1}$;
\item the conditions above are not verified but still one has $J(T'') > C_1 r_{m_{\nmin_{T''}}-1}$, otherwise $T''$ would be a RC.
\end{itemize}

Therefore, in all cases one finds 
\[
\gotN_n(\vartheta) \le 2+  \gotN_n(\vartheta_1'') \le \frac{K}{r_{m_n-1}} J(\vartheta_1'') \le
 \frac{K}{r_{m_n-1}}(J(\vartheta)-J(T'')),
\]
so that the bound \eqref{bru} follows once more.
\EP

%%%%%%%%%%%%%%%%%%%%%%%%%%%%%%%%%%%%%%%%%%%%%%%%%%%%%%%%%%%%%%%%%%%%%%%%%%
%%%%%%%%%%%%%%%%%%%%%%%%%%%%%%%%%%%%%%%%%%%%%%%%%%%%%%%%%%%%%%%%%%%%%%%%%%
\subsection{Non-resonant values}
\label{rispezzo}
%%%%%%%%%%%%%%%%%%%%%%%%%%%%%%%%%%%%%%%%%%%%%%%%%%%%%%%%%%%%%%%%%%%%%%%%%%
%%%%%%%%%%%%%%%%%%%%%%%%%%%%%%%%%%%%%%%%%%%%%%%%%%%%%%%%%%%%%%%%%%%%%%%%%%

By relying on Proposition \ref{brjuno} we can now show that the coefficients of the series \eqref{formale+fout} are well-defined to all orders,
and that, if not for the presence of the resonant lines, we could prove the convergence of the series. 

%%%%%%%%%%%%%%%%%%%%%%%%%%%%%%%%%%%%%%%%%%%%%%%%%%%%%%%%%%%%%%%%%%%%%%%%%%
\begin{defi}[\textbf{Non-resonant value of a tree}]
\label{nolabel}
Define the \emph{non-resonant value} of a tree $\vartheta\in\gotT$ as
\begin{equation}\label{valnrnew}
\Val_{\!N\!R}(\vartheta;c,\om) := \Big(\prod_{\la\in \Lambda(\vartheta)}\LL_\la (c)\Big) \Big(\prod_{v\in N(\vartheta)} \calF_v^*(c)\Big)
\Big(\prod_{\ell\in L_{\!N\!R}(\vartheta)} \matG_\ell(\om) \Big),
\end{equation}
and, analogously, the \emph{non-resonant value} of a trimmed tree $\vartheta\in\breve\gotT$ as
\begin{equation}\label{valnr}
\Val_{\!N\!R}(\vartheta;c,\om) := \Big(\prod_{\la\in \Lambda(\vartheta)}\LL_\la (c)\Big) %\Big(\prod_{v\in N(\vartheta)} \calF_v(\h)\Big)
\Big(\prod_{\ell\in L_{\!N\!R}(\vartheta)} \matG_\ell(\om) \Big),
\end{equation}
where we recall that $L_{\!N\!R}(\vartheta)$ is the set of non-resonant lines in $\vartheta$ on scale $\ge0$.
\end{defi}
%%%%%%%%%%%%%%%%%%%%%%%%%%%%%%%%%%%%%%%%%%%%%%%%%%%%%%%%%%%%%%%%%%%%%%%%%%

%%%%%%%%%%%%%%%%%%%%%%%%%%%%%%%%%%%%%%%%%%%%%%%%%%%%%%%%%%%%%%%%%%%%%%%%%%
\begin{lemma}\label{convergerebbe1}
%Assume there is $C>0$ such  $|\h^{(k)}_j|\le C^k$ for all $k,j$.Then f
For all $s_1\ge0$ and $s_2>0$ such that $s_1+s_2 = s$ and all $s'\in(s_1,s)$, 
there is a positive constant {$A_{\!N\!R}=A_{\!N\!R}(s-s',\al,\om)$} such that,
for any $k\ge 1$, $j\in\ZZZ$, $\nu\in\ZZZ^\ZZZ_f$ and $\s\in\{\pm\}$
and for any $\vartheta \in \breve\gotT^{(k)}_{j,\nu,\s}$, one has
\begin{equation}\label{viodio}
|\Val_{\!N\!R}(\vartheta;c,\om)| \le 
A_{\!N\!R}^k  e^{-s_1|\nu|_\al} e^{-s_2\jap{j}^\alpha} e^{-(s'-s_1)J(\vartheta)}.
\end{equation}
\end{lemma}
%%%%%%%%%%%%%%%%%%%%%%%%%%%%%%%%%%%%%%%%%%%%%%%%%%%%%%%%%%%%%%%%%%%%%%%%%%

%%%%%%%%%%%%%%%%%%%%%%%%%%%%%%%%%%%%%%%%%%%%%%%%%%%%%%%%%%%%%%%%%%%%%%%%%%
\prova
%For $\vartheta\in \breve\gotT^{(k)}_{j,\nu,\s}$, f
First of all, by Proposition \ref{brjuno} we can bound, for any $n_0 \ge 1$,
\begin{equation}\label{carbonara}
\begin{aligned}
\null\hspace{-.3cm} \prod_{\ell\in L_{\!N\!R}(\vartheta)}|\matG_{\ell}(\om)|&\le \prod_{n\ge0}\left( 
\frac{16}{\be({m_{n}})}\right)^{\gotN_n(\vartheta)}
\le \left( \frac{16}{\be({m_{n_0}})}\right)^{|L(\vartheta;n_0)|} \!\!
 \prod_{n\ge n_0+1}\left( \frac{16}{\be({m_{n}})}\right)^{\gotN_n(\vartheta)}\\
&\le (A(n_0))^{{|L(\vartheta;n_0)|}} \prod_{n\ge n_0+1}\left( 
\frac{16}{\be({m_{n}})}\right)^{K J(\vartheta) \frac{1}{r_{m_n-1}}} \\
&=  (A(n_0))^{{|L(\vartheta;n_0)|}} \exp\Biggl(K J(\vartheta){\sum_{n\ge n_0+1} \frac{1}{r_{m_{n}-1}}
\log\left( \frac{16}{\be({m_n})} \right)}\Biggr) ,
\end{aligned}
\end{equation}
with $A(n):=16/\be(m_n)$ and $L(\vartheta;n):=\{\ell \in L_{\!N\!R}(\vartheta) : n_{\ell} \le n\}$.
Thus, since $\om\in\gotB$, fixed $s'\in(0,s)$, we can take {$n_0=n_0(s - s')$} large enough so that 
\begin{equation}\label{sn0}
\ka(n_0):= 
K\sum_{n\ge n_0+1}  \frac{1}{r_{m_n-1}}\log
\left(\frac{16}{\be({m_n})} \right) \le s-s'.
\end{equation}
%is smaller than $s$. Set $s':=s-s(n_0)$ and use Lemma \ref{constance2}.

The product  of the leaf factors is bounded in \eqref{assai}.
Moreover, thanks to Remark \ref{stoper}, the overall number of lines -- and hence of non-resonant lines -- on scale $\ge 0$
is bounded by the number of nodes, that is $k$. Then
we obtain 
\[
|\Val_{\!N\!R}(\vartheta;c,\om)| \le (A(n_0))^k  e^{-s'J(\vartheta)}e^{-s\jap{j}^\al}.
\]

Now, take $s_1 \in [0, s')$, and set $s_2:=s-s_1$.
Then, from the definition of $J(\vartheta)$ in \eqref{J} and the fact that
\begin{equation} \label{seserve}
\sum_{\la \in \Lambda(\vartheta)} \jap{j_\la}^\al \ge \sum_{i\in\ZZZ} \jap{i}^\al |\nu_i | = |\nu|_\al ,
\end{equation}
we obtain \eqref{viodio} with $A_{\!N\!R}=A(n_0)$.
\EP
%%%%%%%%%%%%%%%%%%%%%%%%%%%%%%%%%%%%%%%%%%%%%%%%%%%%%%%%%%%%%%%%%%%%%%%%%%

%%%%%%%%%%%%%%%%%%%%%%%%%%%%%%%%%%%%%%%%%%%%%%%%%%%%%%%%%%%%%%%%%%%%%%%%%%
\begin{rmk}\label{loss}
\emph{
The proofs of Lemmata \ref{convergerebbe2} and \ref{convergerebbe1} 
show that both $D_2$ and $n_0$ -- and hence $A(n_0)$ -- diverge as $s_1 \to s$.
This is expected: the larger the analyticity strip in the angles,
the smaller is the radius of convergence of the power series in the perturbation parameter. 
On the other hand, the larger the value of $\BB_\om(r)$, the larger is the value of $n_0$ we have to take in order to satisfy \eqref{sn0}
and, as a consequence, the larger is the value of the constant $A_{\!N\!R}=A(n_0)$.
 }
\end{rmk}
%%%%%%%%%%%%%%%%%%%%%%%%%%%%%%%%%%%%%%%%%%%%%%%%%%%%%%%%%%%%%%%%%%%%%%%%%%

%%%%%%%%%%%%%%%%%%%%%%%%%%%%%%%%%%%%%%%%%%%%%%%%%%%%%%%%%%%%%%%%%%%%%%%%%%
\begin{lemma}\label{11}
Let $s'\in(0,s)$.
For any $k\ge 1$, $j\in\ZZZ$ and $\s\in\{\pm\}$ and for any $\vartheta \in \breve\gotT^{(k)}_{j,\gote_j,\s}$, one has
\begin{equation}\label{viodioeta}
\left|\frac{1}{c_j^\s}\Val_{\!N\!R}(\vartheta;c,\om)\right| \le A_{\!N\!R}^k  e^{-s'J(\vartheta)},
\end{equation}
with $A_{\!N\!R}$ as in Lemma \ref{convergerebbe1}.
\end{lemma}
%%%%%%%%%%%%%%%%%%%%%%%%%%%%%%%%%%%%%%%%%%%%%%%%%%%%%%%%%%%%%%%%%%%%%%%%%%

%%%%%%%%%%%%%%%%%%%%%%%%%%%%%%%%%%%%%%%%%%%%%%%%%%%%%%%%%%%%%%%%%%%%%%%%%% 
\prova{%remark
Recall that, by Lemma \ref{cijei}, for any tree $\vartheta\in\breve\gotT^{(k)}_{j,\gote_j,\s}$ there is at least one leaf, that we call $\la_\vartheta$,
with $j_{\la_\vartheta}=j$ and $\s_{\la_\vartheta}=\s$. Thus, one realizes immediately that
\[
\left| \frac{1}{c_j^\s}\prod_{\la\in \Lambda(\vartheta)} \LL_\la(c)\right| = 
\left| \prod_{\la\in \Lambda(\vartheta) \setminus \la_\vartheta} c_{j_\la}^{\s_\la}\right| \le 
e^{-s J(\vartheta)},
\]
so that, reasoning as in the proof of Lemma \ref{convergerebbe1}, one finds \eqref{viodioeta}.
\EP
%%%%%%%%%%%%%%%%%%%%%%%%%%%%%%%%%%%%%%%%%%%%%%%%%%%%%%%%%%%%%%%%%%%%%%%%%% 

%%%%%%%%%%%%%%%%%%%%%%%%%%%%%%%%%%%%%%%%%%%%%%%%%%%%%%%%%%%%%%%%%%%%%%%%%%
\begin{rmk} \label{convergerebbe}
\emph{
Lemmata \ref{convergerebbe1} and \ref{11}, together with Lemma \ref{convergerebbe2},
%with $s''=s'-s_1$,
imply that, for all $s_1\ge 0,s_2>0$ such that $s_1+s_2 = s$, by fixing $s'=s-s_2/2$, one has %and all $s'\in(s_1,s)$
\begin{equation}\label{sommazza}
\sum_{\vartheta \in \gotT^{(k)}_{j,\nu,\s}}
|\Val_{\!N\!R}(\vartheta;c,\om) |< D_{\!N\!R}^{k} e^{-s_1 |\nu|_\al}e^{-s_2\jap{j}^\al } ,
\end{equation}
with $D_{\!N\!R}=D_{\!N\!R}(s_2,\al,\om):=2^{34}A_{\!N\!R}^2(s_2/2,\al,\om)\,D_2^2(s_2/2,\al)$. 
Recall that $D_2,\,A_{\!N\!R}$ are defined respectively in Lemmata \ref{convergerebbe2} and \ref{convergerebbe1}. Indeed, 
by using \eqref{disposta} in Remark \ref{stracon} and the notations introduced in Remark \ref{rmk:disposta},
we write
\[ 
\begin{aligned}
& \sum_{\vartheta \in \gotT^{(k)}_{j,\nu,\s}} \Val_{\!N\!R}(\vartheta;c,\om) =
\sum_{X \in \Xi^{(k)}_{j,\nu,\s}} 
\sum_{\vartheta \in X} \Val_{\!N\!R}(\vartheta;c,\om) \\
& \quad = \!\!\!\!
\sum_{X \in \Xi^{(k)}_{j,\nu,\s}} 
\Biggl( \sum_{\vartheta'\in \breve X} \Val_{\!N\!R}(\vartheta';c,\om)\! \Biggr)
\prod_{v \in N_2(\vartheta^X)}
\Biggl(  - \sum_{\vartheta_v'\in \breve X_v} 
\frac{1}{c_{j_{v}}^{\sigma_{v}}} \Val_{\!N\!R}(\vartheta_v';c,\om) \! \Biggr) ,
\end{aligned}
\]
so that we bound
\[
\left| \sum_{\vartheta'\in X(\breve\vartheta)} \Val_{\!N\!R}(\vartheta';c,\om)
\prod_{v \in N_2(\vartheta')}
\Biggl(  \sum_{\vartheta_v'\in X(\breve\vartheta_v)} \Val_{\!N\!R}(\vartheta_v';c,\om) \Biggr) \right| \le
A_{\!N\!R}^{2k} D_2^{2k}  e^{-s_1 |\nu|_\al}e^{-s_2\jap{j}^\al } ,
\]
where we have used that $|N(\vartheta)| \le 2k$ for all $\vartheta \in \gotT^{(k)}_{j,\nu,\s}$ by Remark \ref{kn}.
Finally the factor $2^{34}$ arises by counting the number of unlabelled trees ($\le 4^{12k}$) and the number of the sign labels
assigned to the leaves ($\le 2^{10k}$).
Note that $D_{\!N\!R}$ is decreasing in $s_2$ and $D_{\!N\!R}\to \io$ as $s_2\to 0$.
}
\end{rmk}
%%%%%%%%%%%%%%%%%%%%%%%%%%%%%%%%%%%%%%%%%%%%%%%%%%%%%%%%%%%%%%%%%%%%%%%%%%

%%%%%%%%%%%%%%%%%%%%%%%%%%%%%%%%%%%%%%%%%%%%%%%%%%%%%%%%%%%%%%%%%%%%%%%%%%
\begin{lemma} \label{welldefined}
For all $s_1 \ge 0$ and $s_2>0$ such that $s_1+s_2=s$, one has 
\[
\sum_{\vartheta \in \gotT^{(k)}_{j,\nu,\s}}
|\Val(\vartheta;c,\om) | \le (C(k))^{k} e^{-s_1 |\nu|_\al}e^{-s_2\jap{j}^\al } ,
\]
where {$C(k)=C(s_2,\al,\om, k)$} diverges as $k$ goes to $\infty$.  
\end{lemma}
%%%%%%%%%%%%%%%%%%%%%%%%%%%%%%%%%%%%%%%%%%%%%%%%%%%%%%%%%%%%%%%%%%%%%%%%%%

%%%%%%%%%%%%%%%%%%%%%%%%%%%%%%%%%%%%%%%%%%%%%%%%%%%%%%%%%%%%%%%%%%%%%%%%%%
\prova
Fix $s'=s-s_2/2$ and
let $n_0=n_0(s_2/2(2k-1))$ be such that $ \ka(n_0)< s_2/2(2k-1)$, with $\ka(n_0)$ {and $n_0(s)$} defined in \eqref{sn0}.
For any resonant line $\ell$ there are at least two non-resonant lines with the same small divisor:
the exiting line and the entering line of the chain of which $\ell$ is a link. This means that, for any resonant line $\ell$
with $n_\ell\ge n_0+1$ we can bound, using \eqref{carbonara},
\begin{equation} \label{boundresonantline}
| \matG_\ell(\om)| \le \prod_{n \ge n_0+1}
\Biggl( \; \prod_{\substack{\ell' \in L_{\!N\!R}(\vartheta) \\ n_{\ell'} =n}} |\matG_{\ell'}(\om)| \Biggr) \le 
e^{\ka(n_0)J(\vartheta)}.
\end{equation}
Indeed, the propagator of any resonant line on scale $\ge n_0$ can be bounded by the maximum of the propagator
of any non-resonant scale on the same scale, which, in turn, can be bounded by the maximum of the largest propagator of the non-resonant lines,
and the latter is trivially bounded by the product of the maxima of all propagators of the lines on scale $\ge n_0$.

Therefore, an estimate on the product of the propagators of a tree with $k$ nodes
that we can deduce from the estimate \eqref{carbonara} is that
\[
\begin{aligned}
\prod_{\ell\in L(\vartheta)}|\matG_{\ell}(\om)| & \le A^{2k}(n_0) \prod_{n \ge n_0+1}
\Biggl( \;  \prod_{\substack{\ell\in L_{\!N\!R}(\vartheta) \\ n_\ell =n}} |\matG_{\ell}(\om)| \Biggr)
\Biggl( \; \prod_{\substack{\ell\in L(\vartheta) \setminus L_{\!N\!R}(\vartheta) \\ n_\ell =n}} |\matG_{\ell}(\om)| \Biggr) \\
& \le A^{2k}(n_0) \, e^{\ka(n_0)\, J(\vartheta)} e^{(2k-2)\, \ka(n_0)\, J(\vartheta)} ,
\end{aligned}
\]
where we have that
\begin{itemize}[topsep=0ex]
\itemsep-0.2em
\item each resonant line can be bounded as in \eqref{boundresonantline},
\item the number of resonant lines is bounded by $2k-2$, because there are at most $2k$ lines and, in that case,
there is one leaf only,
\item  the root line and the lines exiting the leaves are non-resonant.
\end{itemize}
This implies the assertion, by reasoning as in Remark \ref{convergerebbe}.
\EP
%%%%%%%%%%%%%%%%%%%%%%%%%%%%%%%%%%%%%%%%%%%%%%%%%%%%%%%%%%%%%%%%%%%%%%%%%%

%%%%%%%%%%%%%%%%%%%%%%%%%%%%%%%%%%%%%%%%%%%%%%%%%%%%%%%%%%%%%%%%%%%%%%%%%%
\begin{rmk} \label{welldefinedbis}
\emph{
The bound in Lemma \ref{welldefined} makes the formal power series in $\e$ well-defined to all orders,
but is not quite sufficient to ensures the convergence.
It is not difficult to show that, if one only aimed at proving the coefficients of the power series \eqref{ukvero} and \eqref{etak} 
to be well-defined, one could rely just on Corollary \ref{basta!}, without using Proposition \ref{brjuno}. 
However, Proposition \ref{brjuno} plays a fundamental role when one looks for bounds which ensure the convergence of the series:
indeed, as seen in Remark \ref{convergerebbe}, it implies that if we ignore the resonant lines we obtain summable bounds, 
so that, if we show that the resonant lines can be taken into account without worsening the bounds, then convergence follows.
}
\end{rmk}
%%%%%%%%%%%%%%%%%%%%%%%%%%%%%%%%%%%%%%%%%%%%%%%%%%%%%%%%%%%%%%%%%%%%%%%%%%

Remark \ref{convergerebbe}  and Lemma \ref{welldefined} %, as well as Lemma \ref{osservo},
show that  the resonant lines are the only possible obstruction to convergence of the formal power series \eqref{formale+fout}. 
As we shall see, even though there are trees whose values do not admit
summable bounds (or, at least, for which we are not able to prove that summable bounds exist), nonetheless, thanks
to cancellations related the translation-covariance of the equation, when we sum together the values of all such trees
the overall contribution turns out to be summable as well as those where no resonant lines appear.
Actually, to exploit the cancellations we have to group together in a suitable way the trees with chains of resonant clusters.
How to implement concretely the program sketched above is discussed in Section \ref{convergenza}.

%%%%%%%%%%%%%%%%%%%%%%%%%%%%%%%%%%%%%%%%%%%%%%%%%%%%%%%%%%%%%%%%%%%%%%%%%%
%%%%%%%%%%%%%%%%%%%%%%%%%%%%%%%%%%%%%%%%%%%%%%%%%%%%%%%%%%%%%%%%%%%%%%%%%%
\section{Convergence of the power series}
\label{convergenza}
\zerarcounters
%%%%%%%%%%%%%%%%%%%%%%%%%%%%%%%%%%%%%%%%%%%%%%%%%%%%%%%%%%%%%%%%%%%%%%%%%%
%%%%%%%%%%%%%%%%%%%%%%%%%%%%%%%%%%%%%%%%%%%%%%%%%%%%%%%%%%%%%%%%%%%%%%%%%%

%%%%%%%%%%%%%%%%%%%%%%%%%%%%%%%%%%%%%%%%%%%%%%%%%%%%%%%%%%%%%%%%%%%%%%%%%%
\subsection{Resonant trees}
%%%%%%%%%%%%%%%%%%%%%%%%%%%%%%%%%%%%%%%%%%%%%%%%%%%%%%%%%%%%%%%%%%%%%%%%%%

It is convenient to introduce two further set of trees. The reason is that, in order to improve the bounds on the coefficients of the series \eqref{veri}
with respect to Section  \ref{stimazze}, we have to deal with resonant clusters independently of the trees
in which they are contained; on the other hand the resonant clusters have been defined as subgraphs of given trees.

%%%%%%%%%%%%%%%%%%%%%%%%%%%%%%%%%%%%%%%%%%%%%%%%%%%%%%%%%%%%%%%%%%%%%%%%%%
\begin{defi}[{\textbf{Set $\boldsymbol{\SSSS}$ of the unexpanded resonant  trees}}]
\label{proviamo}
Let $\SSSS$ be the set of trees $\TT$ %rooted
which obey the same rules as the trees in $\Theta$, according to Definition \ref{thetone},
except for the following:
\begin{enumerate}[topsep=0ex]
\itemsep0em
%\item one has $s_v=1,5$ for all $v\in N(\TT)$;
%
\item there are two \emph{special vertices}, the root $r=v_{out}$ and a leaf $v_{in}$, and two \emph{special lines}
$\ell^{out}_\TT$ and  $\ell^{in}_\TT$,  %the line entering $v_{out}$, and $\ell^{in}_\TT$, the line exiting $v_{in}$, such that 
%
%\begin{enumerate}[topsep=0ex]
%\itemsep0.0em
%\item all the lines in $\TT$ are oriented towards $v_{out}$ (so no line exits $v_{out}$);
%
%\item $s_{v_{in}}=0$ (so no line enters $v_{in}$) while $s_{v_{out}}=1$;
such that
\begin{itemize}[topsep=0ex]
\itemsep0em
\item[1.1.] $\ell^{out}_\TT$  enters $v_{out}$ and  $\ell^{in}_\TT$ exits $v_{in}$,
\item[1.2.]  $\ell^{out}_\TT$ exits a node $v$ with $s_v=1$ if and only if $\ell_\TT^{in}$ enters the same node $v$;
\item[1.3.] one has $\nu_{\ell_T^{out}}=\gote_{j_{\ell_T^{out}}}$;
\end{itemize}
%\item
%if $\calP_{\TT}\subseteq L(\TT)$ denotes the path of lines\footnote{Of course $\calP_\TT=\emptyset$ if $\ell_\TT^{in}$ enters
%the same node $\ell_\TT^{out}$ exits.} connecting $\ell_\TT^{in}$ to $\ell_\TT^{out}$ and  $\ell\in\calP_\TT$, then
%$\s_\ell\nu_\ell\ne \s_\ell\gote_{j_\ell}-\s_{\ell^{in}_\TT}\gote_{{\ell^{in}_\TT}}$.
\vspace{-.02cm}
\item the set of lines $L(\TT)$ is redefined by discarding the lines
$\ell_\TT^{in}$ and $\ell_\TT^{out}$, which are called the \emph{external lines} of $\TT$,
so that $\ell_\TT^{in}$ \emph{enters} $\TT$ while $\ell_\TT^{out}$ \emph{exits} $\TT$;
\item the set of leaves $\Lambda(\TT)$ is redefined without counting $v_{in}$;
\item one says that $\TT$ has scale $\nmax_\TT:=\max\{ n_\ell :  \ell\in L(\TT) \}$
%${\nmax}_\TT =\ge-1$ 
if $L(\TT) \ne \emptyset$ and
%\begin{equation}\label{maxTT}
scale $\nmax_\TT=-1$ if it consists only of one node.
%\end{equation}
%\item if $\TT$ is the {trivial RT}, one says that $\TT$ has scale $-1$.
\end{enumerate}
We call the trees $\TT\in\SSSS$  \emph{unexpanded resonant trees}, and define
$\SSSS^{(k)}_{j,\s,j',\s'}$ as the set of all RTs $\TT\in\SSSS$ of order $k$
such that $j_{\ell^{out}_\TT}=j$, $\s_{\ell^{out}_\TT}=\s$,
$j_{\ell^{in}_\TT}=j'$ and $\s_{\ell^{in}_\TT}=\s'$.
\end{defi}
%%%%%%%%%%%%%%%%%%%%%%%%%%%%%%%%%%%%%%%%%%%%%%%%%%%%%%%%%%%%%%%%%%%%%%%%%%

An example of unexpanded resonant tree $\TT\in\SSSS$ is represented in Figure \ref{caffe}.

%%%%%%%%%%%%%%%%%%%%%%%%%%%%%%%%%%%%%%%%%%%%%%%%%%%%%%%%%%%%%%%%%%%%%%%%%%
% FIGURA 19-2
%%%%%%%%%%%%%%%%%%%%%%%%%%%%%%%%%%%%%%%%%%%%%%%%%%%%%%%%%%%%%%%%%%%%%%%%%%
\begin{figure}[ht]
\vspace{-.2cm}
\centering
%\null
%\hspace{-.6cm}
\ins{050pt}{-068pt}{$v_{out}$}
\ins{082pt}{-085pt}{$\ell^{out}_\TT$}
\ins{390pt}{-044pt}{$v_{in}$}
\ins{360pt}{-058pt}{$\ell^{in}_\TT$}
\subfigure{\includegraphics*[width=5in]{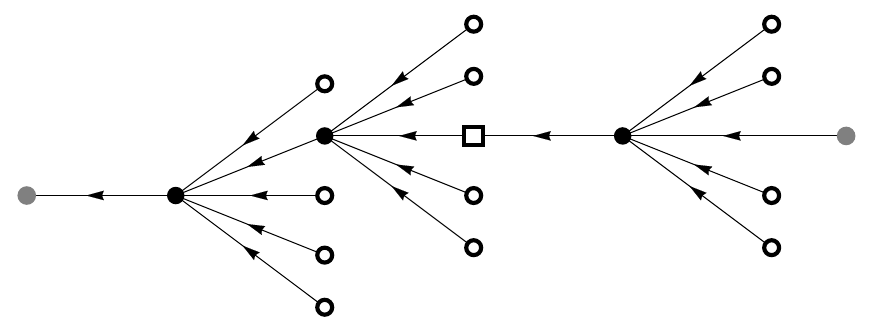}}
\caption{\small A resonant tree $\TT\in\SSSS$ with 4 nodes and 12 leaves (the labels are not shown).}
\label{caffe}
\end{figure}
%%%%%%%%%%%%%%%%%%%%%%%%%%%%%%%%%%%%%%%%%%%%%%%%%%%%%%%%%%%%%%%%%%%%%%%%%%

%Finally we introduce a further set of trees as follows.

%%%%%%%%%%%%%%%%%%%%%%%%%%%%%%%%%%%%%%%%%%%%%%%%%%%%%%%%%%%%%%%%%%%%%%%%%%
\begin{defi}[\textbf{{Set $\boldsymbol{\gotS} %\hspace{-.27cm}\boldsymbol{\gotS}
$ of the expanded resonant trees}}]
\label{speriamo}
Let $\gotS$ be the set of trees $\TT$ which obey the same rules as the trees in $\gotT$, according to Definition \ref{gotT},
except for the following:
\begin{enumerate}[topsep=0ex]
\itemsep0em
\item there are two special vertices, the root $r=v_{out}$ and a leaf $v_{in}$, and two special lines
$\ell^{out}_\TT$ and  $\ell^{in}_\TT$, such that
\begin{itemize}[topsep=0ex]
\itemsep0em
\item[1.1.] $\ell^{out}_\TT$  enters $v_{out}$ and  $\ell^{in}_\TT$ exits $v_{in}$,
\item[1.2.]  $\ell^{out}_\TT$ exits a node $v$ with $s_v=2$ if and only if $\ell_\TT^{in}$ enters the same node $v$;
\item[1.3.] one has $\nu_{\ell_T^{out}}=\gote_{j_{\ell_T^{out}}}$;
\end{itemize}
\vspace{-.02cm}
\item the set of lines $L(\TT)$ is redefined by discarding the lines
$\ell_\TT^{in}$ and $\ell_\TT^{out}$, which are called the \emph{external lines} of $\TT$,
so that $\ell_\TT^{in}$ \emph{enters} $\TT$ while $\ell_\TT^{out}$ \emph{exits} $\TT$;
\item the set of leaves $\Lambda(\TT)$ is redefined without counting $v_{in}$;
\item one says that $\TT$ has scale $\nmax_\TT:=\max\{ n_\ell :  \ell\in L(\TT) \}$ if $L(\TT) \ne \emptyset$ and
scale $\nmax_\TT=-1$ if it consists only of one node.
\end{enumerate}
We call the trees $\TT\in\gotS$ \emph{expanded resonant trees}, and define
$\gotS^{(k)}_{j,\s,j',\s'}$ as the set of RTs $\TT\in\gotS$ of order $k$ with 
$j_{\ell^{out}_\TT}=j$, $\s_{\ell^{out}_\TT}=\s$, %$\nu_{\ell^{in}_T}=\nu'$
$j_{\ell^{in}_\TT}=j'$ and $\s_{\ell^{in}_\TT}=\s'$.
%We call \emph{resonant trees} (RTs) the elements of $\gotS$ as well as the elements of $\SSSS$.
\end{defi}
%%%%%%%%%%%%%%%%%%%%%%%%%%%%%%%%%%%%%%%%%%%%%%%%%%%%%%%%%%%%%%%%%%%%%%%%%%

%%%%%%%%%%%%%%%%%%%%%%%%%%%%%%%%%%%%%%%%%%%%%%%%%%%%%%%%%%%%%%%%%%%%%%%%%%
\begin{rmk} \label{speriamo2}
\emph{
Essentially, Definition \ref{speriamo} repeats Definition \ref{proviamo} up to referring to the set $\gotT$ instead of $\Theta$.
In fact, any expanded resonant tree $\TT$ can be obtained from any RT $\TT'\in\SSSS$ by
\begin{itemize}[topsep=0ex]
\itemsep0em
\item  grafting a subtree $\vartheta\in\gotT^{(k_v)}_{j_v,\gote_{j_v},\s_v}$ to each node $v\in N(\TT')$ with $s_v=1$,
and redefining accordingly the branching label of $v$ as $s_v=2$, and
\item redefining the order label of $v$ as $k_v=0$.
\end{itemize}
Compare %Definition \ref{speriamo} 
with Remark \ref{stracon}, where a similar procedure has been applied to obtain
the set $\gotT$ starting from $\Theta$.
}
\end{rmk}
%%%%%%%%%%%%%%%%%%%%%%%%%%%%%%%%%%%%%%%%%%%%%%%%%%%%%%%%%%%%%%%%%%%%%%%%%%

%%%%%%%%%%%%%%%%%%%%%%%%%%%%%%%%%%%%%%%%%%%%%%%%%%%%%%%%%%%%%%%%%%%%%%%%%%
\begin{defi}[\textbf{Resonant tree}] \label{RT}
We call \emph{resonant tree} (RT) any element of $\SSSS\cup\gotS$, that is any either unexpanded or expanded
resonant tree.
Any $\TT\in \SSSS$ with $N(\TT)=\{v\}$ and $\Lambda(\TT)=\emptyset$ is referred to as a \emph{trivial resonant tree}.
\end{defi}
%%%%%%%%%%%%%%%%%%%%%%%%%%%%%%%%%%%%%%%%%%%%%%%%%%%%%%%%%%%%%%%%%%%%%%%%%%

An example of trivial resonant tree is given in Figure \ref{Tjsk}. 

%%%%%%%%%%%%%%%%%%%%%%%%%%%%%%%%%%%%%%%%%%%%%%%%%%%%%%%%%%%%%%%%%%%%%%%%%%
% FIGURA 12
%%%%%%%%%%%%%%%%%%%%%%%%%%%%%%%%%%%%%%%%%%%%%%%%%%%%%%%%%%%%%%%%%%%%%%%%%%
\begin{figure}[ht]
\vspace{.4cm}
\centering
\ins{124pt}{-000pt}{$\TT_{j,\s}^{(k)} =$}
\ins{212pt}{0012pt}{$k,j,\s$}
\ins{284pt}{0010pt}{$j,\s$}
\ins{180pt}{-012pt}{$j,\gote_j,\s$}
\ins{248pt}{-012pt}{$j,\gote_j,\s$}
\subfigure{\includegraphics*[width=2in]{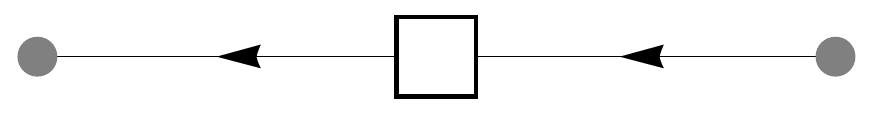}}
\caption{\small The resonant tree $\TT^{(k)}_{j,\s}$ of order $k$ with $j_{\ell_\TT^{out}}=j$ and $\s_{\ell_\TT^{out}}=\s$.}
\label{Tjsk}
\end{figure}
%%%%%%%%%%%%%%%%%%%%%%%%%%%%%%%%%%%%%%%%%%%%%%%%%%%%%%%%%%%%%%%%%%%%%%%%%%

%%%%%%%%%%%%%%%%%%%%%%%%%%%%%%%%%%%%%%%%%%%%%%%%%%%%%%%%%%%%%%%%%%%%%%%%%%
\begin{defi}[\textbf{Resonant path of a resonant tree}] \label{RPT}
Given a resonant tree $\TT$,  let $\calP_{\TT}\subseteq L(\TT)$ denote the path of lines
connecting $\ell_\TT^{in}$ to $\ell_\TT^{out}$ and  $\ell\in\calP_\TT$. We say that $\calP_T$ is the \emph{resonant path} of $\TT$.
\end{defi}
%%%%%%%%%%%%%%%%%%%%%%%%%%%%%%%%%%%%%%%%%%%%%%%%%%%%%%%%%%%%%%%%%%%%%%%%%%

%%%%%%%%%%%%%%%%%%%%%%%%%%%%%%%%%%%%%%%%%%%%%%%%%%%%%%%%%%%%%%%%%%%%%%%%%%
\begin{rmk} \label{speriamo3}
\emph{
Of course $\calP_\TT=\emptyset$ if $\ell_\TT^{in}$ enters the same node $\ell_\TT^{out}$ exits.
Moreover, by construction, for all $\TT\in\gotS$, the path $\calP_\TT$  does not contain any $\h$-line. 
}
\end{rmk}
%%%%%%%%%%%%%%%%%%%%%%%%%%%%%%%%%%%%%%%%%%%%%%%%%%%%%%%%%%%%%%%%%%%%%%%%%%

%%%%%%%%%%%%%%%%%%%%%%%%%%%%%%%%%%%%%%%%%%%%%%%%%%%%%%%%%%%%%%%%%%%%%%%%%%
\begin{rmk}
\emph{
The trivial RT $\TT$ belongs to $\SSSS$ only, and it has  $j_{\ell^{\TT}_{out}}=j_{\ell^{\TT}_{in}}$.}
\end{rmk}
  %%%%%%%%%%%%%%%%%%%%%%%%%%%%%%%%%%%%%%%%%%%%%%%%%%%%%%%%%%%%%%%%%%%%%%%%%%
  
For $\TT\in\SSSS\cup\gotS$, the conservation law \eqref{conserva} leads to
\begin{equation}\label{secindip}
\sum_{\la\in \Lambda^*(\TT)} \s_\la \gote_{j_\la} + \s_{\ell_\TT^{in}}\gote_{j_{\ell_\TT^{in}}}- \s_{\ell_\TT^{out}}\gote_{j_{\ell_\TT^{out}}}=0\,,
\vspace{-.2cm}
\end{equation}
where $\Lambda^*(\TT)$ is the set of leaves $\la\in\Lambda(\TT)$ such that $\la\st{\prec} \ell^{out}_\TT$;
note that $\Lambda^*(\TT)=\Lambda(\TT)$ if $\TT\in\SSSS$.
For all $\ell \in L(\TT)$, set
\begin{equation}\label{ellezeroT}
	\nu^0_\ell=\nu_\ell^0(\TT) := \sum_{\substack{\la\in \Lambda^*(\TT) \\ \la\prec \ell}} \s_{\la}\gote_{j_{\la}} .
\vspace{-.2cm}
\end{equation}

%For all $x\in\RRR$
Introduce the \emph{symbol function} $\x_\ell\!:\RRR\to\RRR$ given by
\begin{equation}\label{lexi}
\x_{\ell}(x) =\begin{cases}
x_\ell , & \ell\notin \calP_\TT, \\
 \s_{\ell}(\om\cdot\nu^0_{\ell} + \s_{\ell^{in}_{\TT}}\om_{j_{\ell^{in}_{\TT}}}+ \s_{\ell^{in}_{\TT}} x ) -\om_{j_{\ell}}    , & \ell\in \calP_\TT ,
\end{cases}
\end{equation}
where $x_\ell$ is the small divisor \eqref{simbolo}.

%%%%%%%%%%%%%%%%%%%%%%%%%%%%%%%%%%%%%%%%%%%%%%%%%%%%%%%%%%%%%%%%%%%%%%%%%%
\begin{defi}[\textbf{Value of a resonant tree }]\label{defvalue5}
For any %unexpanded RT
$\TT\in\SSSS$ and any $x\in\RRR$, we define the \emph{value} of $\TT$ as
\begin{equation}\label{valTvero}
\begin{aligned}
	\VV(\TT;x,c,\om,\h):= & \Big(\prod_{\la\in \Lambda(\TT)}\LL_\la(c) \Big) \Big(\prod_{v\in N(\TT)} \calF_v(\h)\Big)
\Big(\prod_{\ell\in L(\TT)}\calG_{n_\ell}(\x_\ell(x))\Big)	 ,
\end{aligned}
\vspace{-.2cm}
\end{equation}
%$\nu_\ell^0$ is defined as
%
where $\LL_\la(c)$, $\calF_v(\h)$, $\calG_n(x)$ and $\x_\ell(x)$ are defined 
according to \eqref{foglie}, \eqref{nodi}, \eqref{procesi} and \eqref{lexi}, respectively.
Similarly, for any %expanded RT 
$\TT\in\gotS$ and any $x\in\RRR$, we define the \emph{value} of $\TT$ as
\begin{equation}\label{valTeta}
\begin{aligned}
	\Val(\TT;x,c,\om):=&\Big(\prod_{\la\in \Lambda(\TT)}\LL_\la(c) \Big) \Big(\prod_{v\in N(\TT)} \starF_v(c)\Big)
	\Big(\prod_{\ell\in L(\TT)}\calG_{n_\ell}(\x_\ell(x))\Big)	 ,
\end{aligned}
\vspace{-.2cm}
\end{equation}
where $\LL_\la(c)$, $\calF_v^*(c)$, $\calG_n(x)$ and $\x_\ell(x)$ are defined 
according to \eqref{foglie}, \eqref{nodiker}, \eqref{procesi} and \eqref{lexi}, respectively.
\end{defi}
%%%%%%%%%%%%%%%%%%%%%%%%%%%%%%%%%%%%%%%%%%%%%%%%%%%%%%%%%%%%%%%%%%%%%%%%%%

%%%%%%%%%%%%%%%%%%%%%%%%%%%%%%%%%%%%%%%%%%%%%%%%%%%%%%%%%%%%%%%%%%%%%%%%%%
\begin{defi}[\textbf{Resonant tree associated with a resonant cluster}]
\label{TTT}
Consider a tree $\vartheta\in\Theta\cup\gotT$.
With any resonant cluster $T$ in $\vartheta$
we may associate a resonant tree $\TT=\TT_T$ in the following way:
\begin{enumerate}[topsep=0ex]
\itemsep0em
\item replace the subgraph which $\ell_T'$ exits with a leaf $v_{in}$ such that $j_{v_{in}}=j_{\ell_T'}$ and $\s_{v_{in}}=\s_{\ell_T'}$;
\item replace the subgraph which $\ell_T$ enters with a vertex $v_{out}$;
\item modify the momenta of each $\ell\in\calP_T$ so that $\s_\ell \nu_{\ell} = \nu_{\ell}^0+\s_{v_{in}} \gote_{j_{v_{in}}}$.
\end{enumerate}
We call $\TT_T$ the \emph{resonant tree associated with the RC} $T$.
\end{defi}
%%%%%%%%%%%%%%%%%%%%%%%%%%%%%%%%%%%%%%%%%%%%%%%%%%%%%%%%%%%%%%%%%%%%%%%%%

See Figure \ref{RCRT} for an example of resonant tree associated to a resonant cluster of a tree: the resonant tree $\TT_T$
is obtained by replacing the subtree ``on the right" of $T$ with a leaf with the same component label as the root line of the subtree, and
subgraph ``on the left'' of $T$ with a vertex, which represents the root of $\TT_T$.

%%%%%%%%%%%%%%%%%%%%%%%%%%%%%%%%%%%%%%%%%%%%%%%%%%%%%%%%%%%%%%%%%%%%%%%%%%
% FIGURA 
%%%%%%%%%%%%%%%%%%%%%%%%%%%%%%%%%%%%%%%%%%%%%%%%%%%%%%%%%%%%%%%%%%%%%%%%%%
\begin{figure}[ht]
\vspace{-.2cm}
\centering
%\null
%\hspace{-.6cm}
\ins{005pt}{-039pt}{$\vartheta=$}
\ins{140pt}{-020pt}{$T$}
\ins{238pt}{-040pt}{$\TT_T\!=$}
\subfigure{\includegraphics*[width=2.8in]{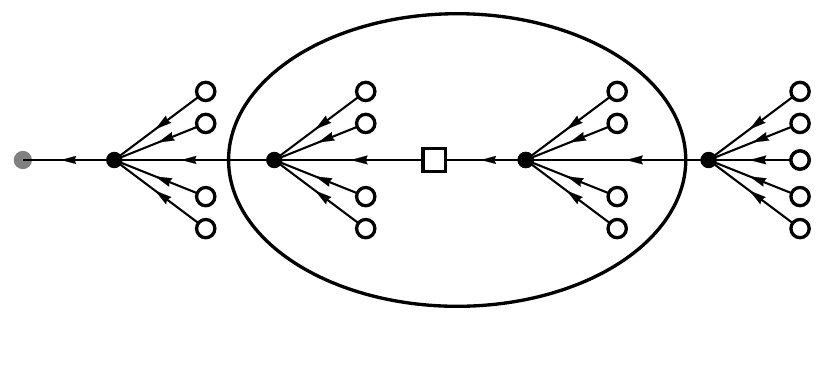}}
\hspace{-.0cm}
\subfigure{\includegraphics*[width=2.8in]{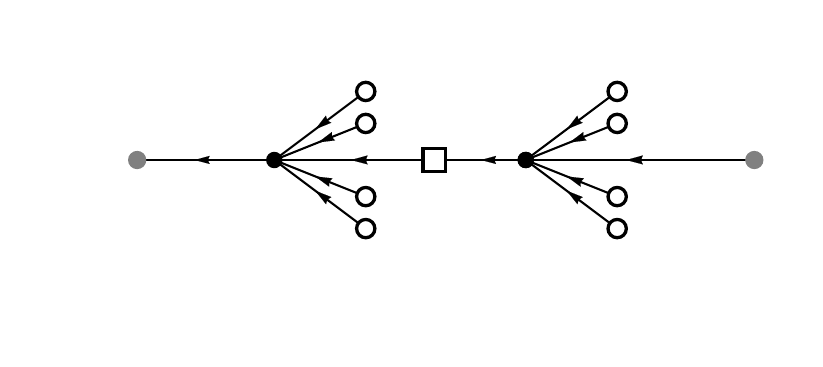}}
\caption{\small A tree $\vartheta\in\Theta$ with a resonant cluster $T$ and the resonant tree $\TT\in\SSSS$ associated with $T$.}
\label{RCRT}
\end{figure}
%%%%%%%%%%%%%%%%%%%%%%%%%%%%%%%%%%%%%%%%%%%%%%%%%%%%%%%%%%%%%%%%%%%%%%%%%%

%%%%%%%%%%%%%%%%%%%%%%%%%%%%%%%%%%%%%%%%%%%%%%%%%%%%%%%%%%%%%%%%%%%%%%%%%%
\begin{defi}[\textbf{Resonance}]
\label{rognoso}
We say that a RC $T$ in a tree $\vartheta\in\Theta\cup\gotT$ and a RC $T'$ in a tree $\vartheta'\in\Theta\cup\gotT$,
possibly coinciding with $\vartheta$,
are \emph{equivalent} if $\TT_T=\TT_{T'}$, and, in that case, we write $T \sim T'$.
We call \emph{resonance} the equivalence class $R:=[T]_\sim$.
\end{defi}
%%%%%%%%%%%%%%%%%%%%%%%%%%%%%%%%%%%%%%%%%%%%%%%%%%%%%%%%%%%%%%%%%%%%%%%%%%

%%%%%%%%%%%%%%%%%%%%%%%%%%%%%%%%%%%%%%%%%%%%%%%%%%%%%%%%%%%%%%%%%%%%%%%%%%
\begin{rmk} \label{pedantissimo}
\emph{%remark
By comparing the two Definitions \ref{proviamo} and \ref{speriamo} with the Definition \ref{coso}, one realizes that for each $\TT$ there exists
a unique resonance $R$ such that  $\TT=\TT_T$, where $T$ is any representative of $R$. We say that such a $R$ is the
\emph{resonance associated} with $\TT$.
}%remark
\end{rmk}
%%%%%%%%%%%%%%%%%%%%%%%%%%%%%%%%%%%%%%%%%%%%%%%%%%%%%%%%%%%%%%%%%%%%%%%%%%

%Let $\SSSS^{(k)}_{j,\s,j',\s'}$ be the set of all RTs $\TT\in\SSSS$ of order $k$, such that $j_{\ell^{out}_\TT}=j$, $\s_{\ell^{out}_\TT}=\s$,
%$j_{\ell^{in}_\TT}=j'$ and $\s_{\ell^{in}_\TT}=\s'$. Similarly 
%let $\gotS^{(k)}_{j,\s,j',\s'}$ denote the set of RTs $\TT\in\gotS$ of order $k$ with  $j_{\ell^{out}_\TT}=j$, $\s_{\ell^{out}_\TT}=\s$, %$\nu_{\ell^{in}_T}=\nu'$
%$j_{\ell^{in}_\TT}=j'$ and $\s_{\ell^{in}_\TT}=\s'$.

%%%%%%%%%%%%%%%%%%%%%%%%%%%%%%%%%%%%%%%%%%%%%%%%%%%%%%%%%%%%%%%%%%%%%%%%%%
\begin{rmk}\label{pedante}
\emph{
If a tree $\vartheta\in\Theta$ has a RC $T$,
by the conservation law \eqref{conserva} we can write
\begin{equation}\label{valTped}
\begin{aligned}
	\VV(T;c,\om,\h)=\Big(&\prod_{\la\in \Lambda(T)}\LL_\la(c) \Big) \Big(\prod_{v\in N(T)} \calF_v(\h)\Big)\\
	&\times\Big(\prod_{\ell\in L(T)\setminus\calP_T}\matG_{\ell}(\om) \Big)
	\Bigg(\prod_{\ell\in \calP_T}
\frac{\Psi_{n_\ell}(\om\cdot (\s_\ell(\nu_\ell^0 +\s_{\ell_T'}\nu_{\ell_T'}))  -\om_{j_\ell})}{\om\cdot (\s_\ell(\nu_\ell^0 +\s_{\ell_T'}\nu_{\ell_T'})) -\om_{j_\ell}}\Bigg).
\end{aligned}
\end{equation}
By comparison with \eqref{valTvero} and recalling \eqref{lexi}, we deduce that, if $\TT_T$ is the RT associated with $T$, then
$\VV(T;c,\om,\h)= \VV(\TT_T;x_{\ell'_T},c,\om,\h)$.
In the same way, if a tree in $\gotT$ has a RC $T$, 
then, if $\TT_T $ is the RT associated with $T$, then
$\Val(T;c,\om)=\Val(\TT_T;x_{\ell'_T},c,\om)$.
}
\end{rmk}
%%%%%%%%%%%%%%%%%%%%%%%%%%%%%%%%%%%%%%%%%%%%%%%%%%%%%%%%%%%%%%%%%%%%%%%%%%

%%%%%%%%%%%%%%%%%%%%%%%%%%%%%%%%%%%%%%%%%%%%%%%%%%%%%%%%%%%%%%%%%%%%%%%%%%
%%%%%%%%%%%%%%%%%%%%%%%%%%%%%%%%%%%%%%%%%%%%%%%%%%%%%%%%%%%%%%%%%%%%%%%%%%
\subsection{Strategy of the proof: truncated trees} %RT and a heuristic argument}
\label{nuovo}
%%%%%%%%%%%%%%%%%%%%%%%%%%%%%%%%%%%%%%%%%%%%%%%%%%%%%%%%%%%%%%%%%%%%%%%%%%
%%%%%%%%%%%%%%%%%%%%%%%%%%%%%%%%%%%%%%%%%%%%%%%%%%%%%%%%%%%%%%%%%%%%%%%%%%

By construction, if $\TT$ is the RT associated with a RC $T$, then both the scale labels of the lines $\ell\in L(\TT)$ and the 
component labels of the leaves $\la\in\Lambda(\TT)$ and of the special vertices of $\TT$ satisfy the same constraints as the
the corresponding labels  of $T$. We have to take into account such constraints if we wish to identify the RTs to be considered
when exploiting the cancellations among the corresponding values. This motivates the definitions introduced in this subsection.

Given a RT $\TT\in\SSSS\cup\gotS$, set
\begin{equation}\label{JTT}
J(\TT):=\sum_{\la \in \Lambda^*(\TT)} \jap{j_\la}^\al - \jap{j_\ell}^\al + \jap{j_{\ell'}}^\al
\end{equation}
and for $\TT\in\gotS$ set
\begin{equation}\label{N2T}
N_2(\TT):=\{v\in N(\TT)\,:\, s_v=2\}.
\end{equation}

%%%%%%%%%%%%%%%%%%%%%%%%%%%%%%%%%%%%%%%%%%%%%%%%%%%%%%%%%%%%%%%%%%%%%%%%%%
\begin{defi}[\textbf{Truncated expanded resonant tree}] \label{frodo}
An expanded resonant tree $\TT\!\in\gotS$ is said to be \emph{truncated} at scale $n$ if 
\begin{enumerate}[topsep=0.5ex]
\itemsep0em
\item $n_\ell\le n$ for all $\ell\in L(\TT)$, 
\item $J(\TT)<C_1r_{m_n-1}$, 
\item $J(\breve{\vartheta}_v)<C_1r_{m_n-1}$ for all $v\in N_2(\TT)$.
\end{enumerate}
We call
$\gotS^{(k)}_{j,\s,j',\s'}(n)$ the set of expanded resonant trees $\TT\in \gotS^{(k)}_{j,\s,j',\s'}$ truncated at scale $n$.
\end{defi}
%%%%%%%%%%%%%%%%%%%%%%%%%%%%%%%%%%%%%%%%%%%%%%%%%%%%%%%%%%%%%%%%%%%%%%%%%%

Define
\begin{equation}\label{emmeta}
%\matM_{j\s j'\s'}^{(k)}(x,c,\om)  := \!\!\!\sum_{\TT\in \gotS^{(k)}_{j,j',\s,\s'}} \!\!\! \Val(\TT;c,x) , \qquad
\matM_{j\s j'\s'}^{(k)}(x,c,\om,n) := \!\!\! \sum_{\TT\in \gotS^{(k)}_{j,j',\s,\s'}(n)} \!\!\!\!\!\! \Val(\TT;x,c,\om).
\end{equation}
For all $j,j'\in\ZZZ$, let us introduce the matrix
\begin{equation}\label{notaz}
\matM_{j j'}^{(k)}(x,c,\om,n) := 
\begin{pmatrix}
\matM_{j+ j'+}^{(k)}(x,c,\om,n) & \matM_{j+ j'-}^{(k)}(-x,c,\om,n) \cr
\matM_{j- j'+}^{(k)}(x,c,\om,n) & \matM_{j- j'-}^{(k)}(-x,c,\om,n)
\end{pmatrix},
\end{equation}
and set\footnote{$P_1$ and $P_3$ are the first and third Pauli matrices.
Usually in the literature the Pauli matrices are denoted $\s_1,\s_2,\s_3$;
we use the notation $P_1,P_2,P_3$ to avoid confusion with the sign labels.}
\begin{equation}\label{matricette}
\CCCC_j :=\begin{pmatrix}
c_j & 0 \cr
0 & \ol{c}_j
\end{pmatrix}, 
\qquad
U:=\begin{pmatrix}
1 & 1 \cr
1& 1
\end{pmatrix},
\qquad
P_1 :=\begin{pmatrix} 0 & 1 \cr 1 & 0 \end{pmatrix}, 
\qquad
P_3:=\begin{pmatrix} 1 & 0 \cr 0 & -1 \end{pmatrix}.
\end{equation}
%

%%%%%%%%%%%%%%%%%%%%%%%%%%%%%%%%%%%%%%%%%%%%%%%%%%%%%%%%%%%%%%%%%%%%%%%%%%
\begin{rmk}\label{fix}
\emph{%remark
Thanks to the constraint {$J(\TT)<C_1r_{m_n-1}$} and $J(\breve{\vartheta}_v)<C_1r_{m_n-1}$ for all $v\in N_2(\TT)$, the set
$\gotS^{(k)}_{j,j',\s,\s'}(n)$ is finite, thus the sum \eqref{emmeta} is finite as well.
}%remark
\end{rmk}
%%%%%%%%%%%%%%%%%%%%%%%%%%%%%%%%%%%%%%%%%%%%%%%%%%%%%%%%%%%%%%%%%%%%%%%%%%

Let us now come back to the problem of bounding the values of the expanded trees.
Consider a tree $\vartheta\in\gotT$ and suppose that there is a chain $\gotC=\{T_1,\ldots,T_{p}\}$ of RCs in $\vartheta$.
Set $k_i=k(T_i)$ and $\ell_i=\ell_{T_i}$ for $i=1,\ldots,p$, and set $\ell_{p=1}=\ell'_{T_p}$.
Set also $j_{i}=j_{\ell_i}$, $\s_i=\s_{\ell_i}$ and $\nu_i=\nu_{\ell_i}$ for $i=1,\ldots,p=1$.
%Note that with each $T_i$ a RT $\TT_i \in \gotS^{(k_i)}_{j_i,\s_i,j_{i+1},\s_{i+1}}$ is associated.
Let $\calF(\vartheta,\gotC)$ be the set of all possible trees  $\vartheta'\in \Theta$ that can be obtained from $\vartheta$
as follows:
\begin{enumerate}[topsep=0ex]
\itemsep0em
\item the RC $T_1$ is replaced with any RC $T_1'$  such that 
$[T_1']_\sim$ is associated with a RT
$\TT_1'\in  %\gotS^{(k_1)}_{j_1,\s_1,j_2,+} \cup \gotS^{(k_1)}_{j_1,\s_1,j_2,-}$ 
%\displaystyle{\bigcup_{{\s}'_2=\pm} $\SSSS^{(k_1)}_{j_1,\s_1,j_2,{\s}'_2}}$
\gotS^{(k_1)}_{j_1,\s_1,j_2,{\s}'_2}$ for some $\s_2'\in\{\pm\}$
and $\nu_2$ is replaced with ${\nu}'_2={\s}'_2\s_1(\nu_1-\gote_{j_1})+\gote_{j_2}$;\footnote{Note that in this way
$\nu_1$ is left unchanged and $T'_1$ is a RC in $\vartheta'$.}
 
\item if $p\ge 3$, for $i=2,\ldots,p-1$  the RC $T_i$ is replaced with any RC $T_i'$ such that
$[T_i']_\sim$ is associated with a RT
$\TT_i'\in
%\displaystyle{\bigcup_{{\s}'_{i+1}=\pm}\SSSS^{(k_i)}_{j_i,{\s}'_i,j_{i+1},{\s}'_{i+1}}}$
\gotS^{(k_i)}_{j_i,{\s}'_i,j_{i+1},{\s}'_{i+1}}$ for some $\s'_{i+1}\in\{\pm\}$
and $\nu_{i+1}$ is replaced with ${\nu}'_{i+1}={\s}'_{i+1}{\s}'_i({\nu}'_i-\gote_{j_i})+\gote_{j_{i+1}}$;\footnote{Note that in this way 
$\nu_{\ell_{T_i'}}={\nu}'_{{i}}$ and $T'_i$ is a RC in $\vartheta'$.}

\item the RC $T_p$ is replaced with any RC $T_p'$ such that 
$[T_p']_\sim$ is associated with a RT
$\TT_p'\in  %\gotS^{(k_1)}_{j_1,\s_1,j_2,+} \cup \gotS^{(k_1)}_{j_1,\s_1,j_2,-}$ \displaystyle{\bigcup_{\s=\pm} 
\gotS^{(k_p)}_{j_p,{\s}'_{p},j_{p+1},\s_{p+1}}$. % and $\nu_{\ell_{T_p'}}=\nu_{{p}}$.
\end{enumerate}

%Note that, in this way, $\nu_1$ is left unchanged and $\nu_{\ell_{T_i'}}={\nu}'_{{i}}$ for $i=2,\ldots,p-1$;
%in particular all the subgraphs $T'_1,\ldots,T'_{p}$ are RCs in $\vartheta'$.

Given a tree in $\calF(\vartheta,\gotC)$,
shortening $x_{i}:=x_{\ell_i}=\om\cdot\nu_{{i}}-\om_{j_i}$ for $i=1,\ldots,p+1$, by Remark \ref{fratellini} we have 
\begin{equation}\label{dormiveglia}
%\begin{aligned}
 \s_i x_{i}=\s_{i+1} x_{i+1} = \s_{p+1}x_{p+1},
\end{equation}
so that, setting $n_i:=\nmin_{T_i}$ and  using Remark \ref{pedante}, we obtain
\begin{equation}\label{sommo}
\begin{aligned}
\sum_{\vartheta'\in\calF(\vartheta,\gotC)} \!\!\! \Val(\vartheta';c,\om) & = A(\vartheta,\gotC) \, \calG_{n_{\ell_{1}}}(x_{p+1}) 
\left[\matM_{j_1 j_{2}}^{(k_1)}(x_{2},c,\om,n_1) \right. \\
%\right.
 & \null \hspace{-.8cm} \left. \times
\prod_{i=2}^p \Big( \calG_{n_{\ell_{i}}}(x_{p+1}) P_3 \matM_{j_{i} j_{i+1}}^{(k_i)}(x_{i+1},c,\om,n_i)\Big) \right]_{\s_1\s_{p+1}} 
\!\!\!\!\!\!  \!\!\!\!\!\! \!\!\!   \calG_{n_{\ell_{p+1}}}(x_{p+1})  \, B(\vartheta,\gotC),
\end{aligned}
\end{equation}
where $A(\vartheta,\gotC)$ and $B(\vartheta,\gotC)$ are common factors and $[\cdot]_{\s\s'}$
denotes the $(\s,\s')$-entry of the matrix in the square bracket, consistently with the notation \eqref{notaz}.

%Set $k:=k_1+\ldots+k_p$. 
The product of the propagators in \eqref{sommo} gives a factor proportional to $|x_{p+1}|^{-(p+1)}$.
This may be a problem, even if all the RT contributing to $ \matM_{j_{i} j_{i+1}}^{(k_i)}(x_{i+1},c,\om,n_i)$ for all $i=1,\ldots,p$
contain no RC, so that a bound in the spirit of \eqref{viodio} can be proved to give
$|\matM_{j_{i} j_{i+1}}^{(k_i)}(x_{i+1},c,\om,n_i)|\le C^{k_i}$ for some $C>0$.
Indeed, in general, %it may happen that  $p =O(k)$ %and $x_{p+1}=O( k^{-1})$,
the factor $|x_{p+1}|^{-p}$ produces a non-summable bound, %of size $k^{O(k)}$, 
thus preventing the convergence of the series. For instance, this happens in the case of the tree $\vartheta$
represented in Figure \ref{pcatena}, where $\gotC=\{T_1,\ldots,T_p\}$ is a chain or RCs 
and $\vartheta'$ is a subtree of order $p$ which does not contain any resonant lines. Then one has
\[
\Val(\vartheta;c,\om) = \prod_{i=1}^p \Biggl( \matG_{\ell_{T_i}} \prod_{\la \in \Lambda (T_i)} \LL_\la(c) \Biggr)
\Val(\vartheta',c,\om) ,
\]
so that, if $\nu$ and $j$ are the momentum and the component label of the line $\ell_{T_p}'$
and $\om$ satisfies the Diophantine condition \eqref{diofantinoBIS}, we may bound, for some constant $C$,
\begin{subequations} \label{stimeorrende}
\begin{align}
& \left| \prod_{i=1}^p \matG_{\ell_{T_i}} \right| \le |x_{p+1}|^{-p} \le C^p 
\prod_{\substack{i\in \ZZZ \\ \nu_i \neq 0}} (1 + \jap{i}^2|\nu_i|^2)^{\tau p} ,
\label{stimeorrendea} \\
& | \Val(\vartheta';c,\om)| \le A_{\!N\!R}^p  e^{-s_1|\nu|_\al} e^{-s_2\jap{j}^\alpha} e^{-(s'-s_1)J(\vartheta')} ,
\label{stimeorrendeb}
\end{align}
\end{subequations}
where \eqref{stimeorrendeb} follows from \eqref{viodio} in Lemma \ref{convergerebbe1}.
On the other hand from \eqref{stimeorrendeb} we may extract a factor
\[
\prod_{\la \in \Lambda(\vartheta')} e^{-s_0 \jap{j_\la}^\al } \le
\prod_{\substack{i\in \ZZZ \\ \nu_i \neq 0}} e^{-s_0 \jap{i}^\al |\nu_i| },
\]
with $2s_0:=\min\{s_2,s'-s_1\}$, so that, by taking into account also \eqref{stimeorrendea}, we obtain
eventually a bound, for a different constant $C$,
\begin{equation} \label{stimeorrendec}
\left| \Val(\vartheta;c,\om) \right| \le %C^k \prod_{\substack{i\in \ZZZ \\ \nu_i \neq 0}} | i \nu_i|^{2\tau(p+1)} e^{-s_0 | i \nu_i|^{\al} } \le
C^k \prod_{\substack{i\in \ZZZ \\ \nu_i \neq 0}} p^{2\tau p /\al}  \le % e^{-s_0 | i \nu_i|^{\al} } \le
C^k \bigl( k^{2\tau k/\al} \bigr)^k 
\end{equation}
where we have used the order of $\vartheta$ is $k=2p$ and that the number of components $\nu_i \neq0$ is at most $p$.

%%%%%%%%%%%%%%%%%%%%%%%%%%%%%%%%%%%%%%%%%%%%%%%%%%%%%%%%%%%%%%%%%%%%%%%%%%
% FIGURA 
%%%%%%%%%%%%%%%%%%%%%%%%%%%%%%%%%%%%%%%%%%%%%%%%%%%%%%%%%%%%%%%%%%%%%%%%%%
\begin{figure}[ht]
\vspace{.2cm}
\centering
%\null
%\hspace{-.6cm}
\ins{005pt}{-033pt}{$\vartheta=$}
\ins{096pt}{-000pt}{$T_1$}
\ins{162pt}{-000pt}{$T_2$}
\ins{230pt}{-000pt}{$T_3$}
\ins{362pt}{-000pt}{$T_p$}
\ins{420pt}{-010pt}{$\vartheta'$}
\subfigure{\includegraphics*[width=5.8in]{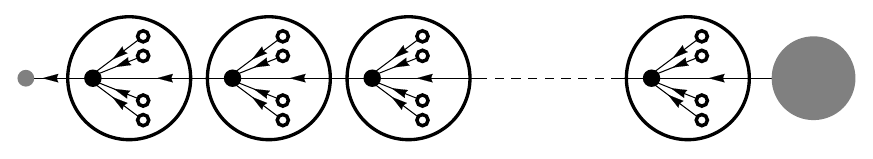}}
\caption{\small A tree $\vartheta\in\Theta$ of order $k=2p$ with a chain of $p$ resonant clusters.}
\label{pcatena}
\end{figure}
%%%%%%%%%%%%%%%%%%%%%%%%%%%%%%%%%%%%%%%%%%%%%%%%%%%%%%%%%%%%%%%%%%%%%%%%%%

The crucial point is that  -- as we shall see -- the value of $\h(c,\om)$ which makes the kernel equation \eqref{bif} to be satisfied
allows us also to exploit the symmetries due to the Hamiltonian nature of the equation and obtain a bound like 
\begin{equation}\label{cespero}
\left|\left[\matM_{j_1 j_{2}}^{(k_1)}(x_{2},c,\om,n_1) 
\prod_{i=2}^p   \Big( P_3 
\matM_{j_{i} j_{i+1}}^{(k_i)}(x_{i+1},c,\om,n_i)\Big)\right]_{\s_1\s_{p+1}} \right| 
\le {C^{k} }{|x_{p+1}|^{p-5}},
\end{equation}
for some constant $C>0$. This provides a gain factor that compensates the product of the propagators
of the resonant lines in \eqref{sommo}.
For instance, in the example of Figure \ref{pcatena}, once we have taken into account
all possible chains with $p$ RCs, we obtain a factor $|x_{p+1}|^{-6}$ instead of $|x_{p+1}|^{-p-1}$,
which in turn yields a final bound $C^k 6^{12\tau k/\al}$ instead of \eqref{stimeorrendec}.
However, to prove that a result  like \eqref{cespero} holds requires some work, because we cannot ignore that
the RTs may contain chains of RCs. Thus, we need a different iterative procedure, that we shall describe in
Subsections \ref{LR} and \ref{LFR}.

First of all, in analogy with Definition \ref{frodo}, we introduce the following sets of truncated trees.
The reason why we consider both expanded and unexpanded trees
is that on the one hand we need to work with the expanded trees in order to estimate the coefficients of fixed order,
on the other hand the symmetry properties responsible for the cancellations are better visualised in terms of the unexpanded trees --
as it will be discussed in Subsection \ref{simmetrie}.

%%%%%%%%%%%%%%%%%%%%%%%%%%%%%%%%%%%%%%%%%%%%%%%%%%%%%%%%%%%%%%%%%%%%%%%%%%
\begin{defi} [\textbf{Truncated unexpanded tree}] \label{Thetan}
An unexpanded tree $\vartheta\in \Theta$ is said to be \emph{truncated} at scale $n$ if
\begin{enumerate}[topsep=0ex]
\itemsep0em
\item $n_\ell< n$ for all $\ell\in L(\vartheta)$,
\item $J(\vartheta)<C_1r_{m_n-1}$.
\end{enumerate}
We call $\Theta^{(k)}_{j,\nu,\s}(n)$ the set of trees $\vartheta\in \Theta^{(k)}_{j,\nu,\s}$
which are truncated at scale $n$.
\end{defi}
%%%%%%%%%%%%%%%%%%%%%%%%%%%%%%%%%%%%%%%%%%%%%%%%%%%%%%%%%%%%%%%%%%%%%%%%%%

%%%%%%%%%%%%%%%%%%%%%%%%%%%%%%%%%%%%%%%%%%%%%%%%%%%%%%%%%%%%%%%%%%%%%%%%%%
\begin{defi} [\textbf{Truncated expanded tree}]\label{Tn}
An expanded tree $\vartheta \in \gotT$ is said to be \emph{truncated} at scale $n$ if
\begin{enumerate}[topsep=0ex]
\itemsep0em
\item $n_\ell<n$ for all $\ell\in L(\vartheta)$,
\item $J(\breve\vartheta)<C_1r_{m_n-1}$,
\item for all $v\in N_2(\vartheta)$ one has $J(\breve\vartheta_v)<C_1r_{m_n-1}$.
\end{enumerate}
We call $\gotT^{(k)}_{j,\nu,\s}(n)$ the set of trees $\vartheta\in\gotT^{(k)}_{j,\nu,\s}$
which are truncated at scale $n$.
\end{defi}
%%%%%%%%%%%%%%%%%%%%%%%%%%%%%%%%%%%%%%%%%%%%%%%%%%%%%%%%%%%%%%%%%%%%%%%%%%

%%%%%%%%%%%%%%%%%%%%%%%%%%%%%%%%%%%%%%%%%%%%%%%%%%%%%%%%%%%%%%%%%%%%%%%%%%
\begin{defi}[\textbf{Truncated unexpanded resonant tree}]\label{rifrodo}
An unexpanded resonant tree $\TT\in \SSSS$ is said to be \emph{truncated} at scale $n$ if
\begin{enumerate}[topsep=0ex]
\itemsep0em
\item $n_\ell< n$ for all $\ell\in L(\TT)$,
\item $J(\TT)<C_1r_{m_n-1}$.
\end{enumerate}
We call $\SSSS^{(k)}_{j,\nu,\s}(n)$ the set of unexpanded resonant trees $\TT\in \SSSS^{(k)}_{j,\nu,\s}$
which are truncated at scale $n$.
\end{defi}
%%%%%%%%%%%%%%%%%%%%%%%%%%%%%%%%%%%%%%%%%%%%%%%%%%%%%%%%%%%%%%%%%%%%%%%%%%

%%%%%%%%%%%%%%%%%%%%%%%%%%%%%%%%%%%%%%%%%%%%%%%%%%%%%%%%%%%%%%%%%%%%%%%%%%
\begin{rmk}\label{troncoj}
\emph{
Thanks to the constraints in Definitions \ref{Thetan} to \ref{rifrodo}
%that $n_\ell< n$ for all $\ell\in L(\vartheta)$ and $J(\vartheta)<C_1r_{m_n-1}$, 
the sets $\Theta^{(k)}_{j,\nu,\s}(n)$,  $\gotT^{(k)}_{j,\nu,\s}(n)$  and $\SSSS^{(k)}_{j,\nu,\s}(n)$ are finite as well as
$\gotS^{(k)}_{j,j',\s,\s'}(n)$, as pointed out in Remark \ref{fix}.
Note that, if we further require $\vartheta\in\gotT\cup\breve\gotT$ to satisfy the support property,
then the conditions 2 and 3 in Definition \ref{Tn} imply automatically the condition 1.
}
\end{rmk}
 %%%%%%%%%%%%%%%%%%%%%%%%%%%%%%%%%%%%%%%%%%%%%%%%%%%%%%%%%%%%%%%%%%%%%%%%%%

Set
\begin{subequations} \label{gin-ocean}
\begin{align}
G_{j,\s}^{(k)}(c,\om,\h,n) & := \sum_{\vartheta\in \Theta^{(k)}_{j,\gote_j,\s}(n)}\VV(\vartheta;c,\om,\h) ,
\label{gin} \\
(\h_j^{(k)}(c,\om,n))^{\s} & := - \frac{1}{c^\s_j} \sum_{\vartheta  \in \gotT^{(k)}_{ j ,\gote_j,\s}(n)} \Val (\vartheta;c,\om)  \,,
\label{ocean}
\end{align}
\end{subequations}
and  define
\begin{equation}\label{emmenne}
\MM_{j\s j'\s'}^{(k)}(x,c,\om,\h,n):= \sum_{\TT\in \SSSS^{(k)}_{j,\s,j',\s'}(n)}\VV(\TT;x,c,\om,\h) .
\end{equation}

%%%%%%%%%%%%%%%%%%%%%%%%%%%%%%%%%%%%%%%%%%%%%%%%%%%%%%%%%%%%%%%%%%%%%%%%%%
\begin{rmk}\label{troncoj2}
\emph{
By Remark \ref{troncoj}, the sums in \eqref{gin-ocean} and \eqref{emmenne} are finite.
}
\end{rmk}
 %%%%%%%%%%%%%%%%%%%%%%%%%%%%%%%%%%%%%%%%%%%%%%%%%%%%%%%%%%%%%%%%%%%%%%%%%%

%%%%%%%%%%%%%%%%%%%%%%%%%%%%%%%%%%%%%%%%%%%%%%%%%%%%%%%%%%%%%%%%%%%%%%%%%%
\begin{rmk}\label{sing}
\emph{
Lemma \ref{cijei} implies $(\h_j^{(k)}(c,\om,n))^{+}=(\h_j^{(k)}(c,\om,n))^{-}$, i.e.~$\h_j^{(k)}(c,\om,n)\in\RRR$.
Moreover,  by Lemma \ref{secondome} one has
\begin{equation}\label{connome}
G^{(k)}_{j,\s}(c,\om,\h(c,\om,n),n)\equiv0.
\end{equation}
}
\end{rmk}
%%%%%%%%%%%%%%%%%%%%%%%%%%%%%%%%%%%%%%%%%%%%%%%%%%%%%%%%%%%%%%%%%%%%%%%%%%

It is convenient to introduce the notation $\vartheta_{j,\s}^{(k)}$
for the tree in $\Theta^{(k)}_{j,\gote_j,\s}$ with only one node $v$ with $s_v=1$ and $k_v=k$ (see Figure \ref{thetajsk}), and define
\begin{equation}\label{tilden}
\tilde{G}_{j,\s}^{(k)} (c,\om,\h,n)  := \sum_{\vartheta\in \tilde{\Theta}^{(k)}_{j,\gote_j,\s}(n)}\VV(\vartheta;c,\om,\h),\qquad
\tilde{\Theta}^{(k)}_{j,\gote_j,\s}(n) :=  {\Theta}^{(k)}_{j,\gote_j,\s}(n)\setminus\{\vartheta_{j,\s}^{(k)}\} ,
\end{equation}
so as to write
\begin{equation}\label{posso}
G_{j,\s}^{(k)}(c,\om,\h,n) = \VV(\vartheta_{j,\s}^{(k)};c,\om,\h) + \tilde G_{j,\s}^{(k)}(c,\om,\h,n) , \qquad
\VV(\vartheta_{j,\s}^{(k)};c,\om,\h)=(\h^{(k)}_j)^{\s} c^{\s}_j.
\end{equation}
%

%%%%%%%%%%%%%%%%%%%%%%%%%%%%%%%%%%%%%%%%%%%%%%%%%%%%%%%%%%%%%%%%%%%%%%%%%%
% FIGURA 12
%%%%%%%%%%%%%%%%%%%%%%%%%%%%%%%%%%%%%%%%%%%%%%%%%%%%%%%%%%%%%%%%%%%%%%%%%%
\begin{figure}[ht]
\vspace{.4cm}
\centering
\ins{124pt}{-000pt}{$\vartheta_{j,\s}^{(k)} =$}
\ins{212pt}{0012pt}{$k,j,\s$}
\ins{284pt}{0010pt}{$j,\s$}
\ins{180pt}{-012pt}{$j,\gote_j,\s$}
\ins{248pt}{-012pt}{$j,\gote_j,\s$}
\subfigure{\includegraphics*[width=2in]{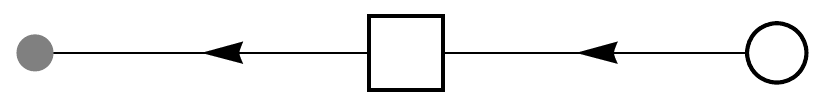}}
\caption{\small The tree $\vartheta_{j,\s}^{(k)}$.}
\label{thetajsk}
\end{figure}
%%%%%%%%%%%%%%%%%%%%%%%%%%%%%%%%%%%%%%%%%%%%%%%%%%%%%%%%%%%%%%%%%%%%%%%%%%

%%%%%%%%%%%%%%%%%%%%%%%%%%%%%%%%%%%%%%%%%%%%%%%%%%%%%%%%%%%%%%%%%%%%%%%%%%
\begin{rmk}\label{ririserve?}
\emph{
For any $\vartheta\in\tilde{\Theta}^{(k)}_{j,\gote_j,\s}(n)$,
in the light of Remark \ref{noetaleaf}, the lines exiting the leaves can only enter nodes $v\in N(\vartheta)$ with $s_v=5$.
}
\end{rmk}
%%%%%%%%%%%%%%%%%%%%%%%%%%%%%%%%%%%%%%%%%%%%%%%%%%%%%%%%%%%%%%%%%%%%%%%%%%

%%%%%%%%%%%%%%%%%%%%%%%%%%%%%%%%%%%%%%%%%%%%%%%%%%%%%%%%%%%%%%%%%%%%%%%%%%
\begin{rmk}\label{ordinieta}
\emph{
Given any $\vartheta\in\tilde\Theta^{(k)}_{j,\nu,\s}(n)$,
only the parameters $(\h_{j'}^{(h)})^{\s'}$ with $h<k$ contribute to $\VV(\vartheta;c,\om,\h)$.
Thus, in fact,
 $\tilde{G}_{j,\s}^{(k)}(c,\om,\h,n)$ depends only on $(\h_{j'}^{(h)})^{\s'}$ with $h<k$.
}
\end{rmk}
%%%%%%%%%%%%%%%%%%%%%%%%%%%%%%%%%%%%%%%%%%%%%%%%%%%%%%%%%%%%%%%%%%%%%%%%%%

As a consequence of Lemma \ref{cijei} one has
\begin{equation} \label{cijei2}
\tilde{G}_{j,\s}^{(k)}(c,\om, \h,n) = c_j^\s \widehat{G}_{j}^{(k)}(|c|^2,\om,\h,n),
\end{equation}
where $|c|^2:=\{|c_i|^2\}_{i\in\ZZZ}$ and
$\widehat{G}^{(k)}_{j}(|c|^2,\om,\h,n)$ is a  polynomial in $|c_i|^2$ of degree $\le 2k$ for all $i\in\ZZZ$.
Thus, if $\h_i^{(h)} \in \RRR$
for all $i\in\ZZZ$ and all $h < k$, then $\widehat{G}_{j}^{(k)}(|c|^2,\om,\h,n) \in \RRR$. 
In particular $\widehat{G}_{j,\s}^{(1)}(|c|^2,\om,\h,n)$ does not depend on $\h$ and hence it is real {even if $\eta$ is not}.
%}%remark
%\end{rmk}
%%%%%%%%%%%%%%%%%%%%%%%%%%%%%%%%%%%%%%%%%%%%%%%%%%%%%%%%%%%%%%%%%%%%%%%%%%
Therefore \eqref{connome} can be rewritten as
\begin{equation}\label{labase}
\h_j^{(k)}(c,\om,n) + \widehat{G}^{(k)}_{j}(|c|^2,\om,\h(c,\om,n),n)\equiv0.
\end{equation}

For all $j,j'\in\ZZZ$, let us introduce the matrix
\begin{equation}\label{rinotaz}
\MM_{j j'}^{(k)}(x,c,\om,\h,n) := 
\begin{pmatrix}
\MM_{j+ j'+}^{(k)}(x,c,\om,\h,n) & \MM_{j+ j'-}^{(k)}(-x,c,\om,\h,n) \cr
\MM_{j- j'+}^{(k)}(x,c,\om,\h,n) & \MM_{j- j'-}^{(k)}(-x,c,\om,\h,n)
\end{pmatrix}.
\end{equation}
%

%%%%%%%%%%%%%%%%%%%%%%%%%%%%%%%%%%%%%%%%%%%%%%%%%%%%%%%%%%%%%%%%%%%%%%%%%%
\begin{rmk}\label{giusto}
\emph{%remark
	As noted in Remark \ref{calcolo}, the sum \eqref{emmeta} equals the sum \eqref{emmenne}
	%as soon as one replaces the parameters $\h^{(k_v)}_{j_v}$
	%appearing in each $\Val(\vartheta;c,\om,\h)$ with $\h^{(k_v)}_{j_v}(c)$ in \eqref{etak}.
	as soon as, for any $\TT\in\SSSS^{(k)}_{j,\nu,\s}(n)$ and any node $v\in N(\TT)$ with $s_v=1$,
	one replaces the parameter $\h^{(k_v)}_{j_v}$
	with $\h^{(k_v)}_{j_v}(c,\om,n)$ in \eqref{ocean}. In other words one has
	\[
	{\matM}^{(k)}_{j,\s,j',\s'}(x,c,\om,n)=\MM^{(k)}_{j,\s,j',\s'}(x,c,\om,\h(c,\om,n),n) .
	\]
	}%remark
\end{rmk}
%%%%%%%%%%%%%%%%%%%%%%%%%%%%%%%%%%%%%%%%%%%%%%%%%%%%%%%%%%%%%%%%%%%%%%%%%%

%%%%%%%%%%%%%%%%%%%%%%%%%%%%%%%%%%%%%%%%%%%%%%%%%%%%%%%%%%%%%%%%%%%%%%%%%%
%%%%%%%%%%%%%%%%%%%%%%%%%%%%%%%%%%%%%%%%%%%%%%%%%%%%%%%%%%%%%%%%%%%%%%%%%%
\subsection{Symmetries and cancellations}
\label{simmetrie}
%%%%%%%%%%%%%%%%%%%%%%%%%%%%%%%%%%%%%%%%%%%%%%%%%%%%%%%%%%%%%%%%%%%%%%%%%%
%%%%%%%%%%%%%%%%%%%%%%%%%%%%%%%%%%%%%%%%%%%%%%%%%%%%%%%%%%%%%%%%%%%%%%%%%%

As anticipated in Subsection \ref{nuovo}, the resonant lines are no longer an obstruction to the convergence
if the counterterm is appropriately fixed. This is due to remarkable identities that we are going to discuss here.

%%%%%%%%%%%%%%%%%%%%%%%%%%%%%%%%%%%%%%%%%%%%%%%%%%%%%%%%%%%%%%%%%%%%%%%%%%
\begin{lemma}\label{cijeiSEC}
For all $j,j'\in\ZZZ$,  $k\ge1$ and $\s,\s'\in\{\pm\}$ and for any $\TT\in \SSSS^{(k)}_{j,\s,j',\s'}\cup \gotS^{(k)}_{j,\s,j',\s'}$ one has
\[
\begin{aligned}
& {|\Lambda^*_{j'',+}(\TT)| = |\Lambda^*_{j'',-}(\TT)|} , \qquad j''\ne j,j' , \phantom{\Big)} \\
&|\Lambda^*_{j,\s}(\TT)| = | \Lambda^*_{j,-\s}(\TT)|+1-\de_{jj'} \de_{\s\s'} \\
 &| \Lambda^*_{j',-\s'}(\TT)| = | \Lambda^*_{j',\s'}(\TT)|+1-\de_{jj'} \de_{\s\s'} \phantom{\Big)}
\end{aligned}
\]
\end{lemma}
%%%%%%%%%%%%%%%%%%%%%%%%%%%%%%%%%%%%%%%%%%%%%%%%%%%%%%%%%%%%%%%%%%%%%%%%%%

%%%%%%%%%%%%%%%%%%%%%%%%%%%%%%%%%%%%%%%%%%%%%%%%%%%%%%%%%%%%%%%%%%%%%%%%%%
\prova
One reasons as in Lemma \ref{cijei}, but uses \eqref{lefoglie*} {instead of \eqref{lefoglie}}.
\EP
%%%%%%%%%%%%%%%%%%%%%%%%%%%%%%%%%%%%%%%%%%%%%%%%%%%%%%%%%%%%%%%%%%%%%%%%%%

The matrix $\MM_{jj'}^{(k)}(x,c,\om,\h,n)$ in \eqref{rinotaz} satisfies the following symmetry properties.

%%%%%%%%%%%%%%%%%%%%%%%%%%%%%%%%%%%%%%%%%%%%%%%%%%%%%%%%%%%%%%%%%%%%%%%%%%
\begin{lemma}\label{pauli}
For all $k\ge 1$, $j,j' \in\ZZZ$, $n\ge -1$,
and for any $\h %=\{\h^{(k)}_j\}_{j\in\ZZZ,k\ge1}
\in\RRR^{\ZZZ\times\NNN}$,
there are real functions $\calA^{(k)}_{jj'\s}(x)=\calA^{(k)}_{jj'\s}(x,|c|^2,\om,\h,n)$, with $\s=\pm$,
%$\calA^{(k)}_{jj'-}(x)=\calA^{(k)}_{jj'}(x,|c|^2,\om,\h,n)$ and
and $\BB^{(k)}_{jj'}(x)=\BB^{(k)}_{jj'}(x,|c|^2,\om,\h,n)$, such that the matrix \eqref{notaz} can be written as
\begin{equation}\label{above}
\begin{aligned}
\null \hspace{-.2cm}
\MM^{(k)}_{jj'}(x,c,\om,\h,n) &= \CCCC_j \!
\begin{pmatrix}
\calA^{(k)}_{jj'+}(x) & \calA^{(k)}_{jj'-}(-x) \cr
\calA^{(k)}_{jj'-}(x) & \calA^{(k)}_{jj'+}(-x)
\end{pmatrix} \! \CCCC_{j'}^\dagger
% \\ & \qquad 
+ \de_{jj'} \!
\begin{pmatrix}
\BB^{(k)}_{jj'}(x) & 0 \cr
0 & \BB^{(k)}_{jj'}(-x)
\end{pmatrix},
\end{aligned}
\end{equation}
with $\CCCC_j$ defined in \eqref{matricette}.
\end{lemma}
%%%%%%%%%%%%%%%%%%%%%%%%%%%%%%%%%%%%%%%%%%%%%%%%%%%%%%%%%%%%%%%%%%%%%%%%%%

%%%%%%%%%%%%%%%%%%%%%%%%%%%%%%%%%%%%%%%%%%%%%%%%%%%%%%%%%%%%%%%%%%%%%%%%%%
\prova
Given  $\TT\in \SSSS^{(k)}_{j,\s,j',\s'}(n)$, 
using Lemma \ref{cijeiSEC} and the fact that both the sequence $\h$ and the propagators are real, 
one can write
\begin{equation}\label{aquesta}
\VV(\TT;x,c,\om,\h)=
\begin{cases}
c_{j}^{\s}c_{j'}^{-\s'}A(\TT;x,|c|^2,\om,\h) , &  \Lambda^*_{j',-\s'}(\TT)\ne\emptyset,\\
 \de_{jj'}\de_{\s\s'}B(\TT;x,|c|^2,\om,\h), &  \Lambda^*_{j',-\s'}(\TT)=\emptyset,
\end{cases}
\end{equation}
with suitable $A(\TT;x,|c|^2,\om,\h)$ and $B(\TT;x,|c|^2,\om,\h)$ in $\RRR$.

Moreover, if
 $\TT'\in\SSSS^{(k)}_{j,-\s,j',-\s'}(n)$ is the RT obtained from $\TT $ 
by changing the sign labels of all the leaves, including the special vertex $v_{in}$, 
by direct inspection one finds
\[
A(\TT;x,|c|^2,\om,\h)=A(\TT';x,|c|^2,\om,\h),\qquad B(\TT;x,|c|^2,\om,\h)=B(\TT';x,|c|^2,\om,\h),
\] 
so that
\[
\begin{aligned}
\calA^{(k)}_{jj'\s}(x,|c|^2,\om,\h,n)&:= \!\!\!\!\! 
\sum_{\TT\in \SSSS^{(k)}_{j,+ ,j', \s }(n)} \!\!\!\!\!  A(\TT;x,c,\om,\h) = \!\!\!\!\!  \sum_{\TT'\in \SSSS^{(k)}_{j,- ,j', -\s }(n)} \!\!\!\!\! A(\TT';x,c,\om,\h),\\
\BB^{(k)}_{jj'}(x,|c|^2,\om,\h,n)&:= \de_{jj'} \!\!\!\!\!  \sum_{\TT\in \SSSS^{(k)}_{j,+ ,j, + }(n)} \!\!\!\!\!  B(\TT;x,c,\om,\h) = \de_{jj'}
\!\!\!\!\!  \sum_{\TT'\in \SSSS^{(k)}_{j,- ,j, - }(n)} \!\!\!\!\!  B(\TT';x,c,\om,\h),\\
\end{aligned}
\]
with both $\calA_{jj'\s}(x,|c|^2,\om,\h,n)$ and $\BB_{jj'}(x,|c|^2,\om,\h,n)$ in $\RRR$. 
Therefore the assertion follows.
%depend on $|c|^2$ and $n$.
%%%%%Therefore, summing over all RTs of order $k$, one finds
%%%%%\[
%%%%%\begin{aligned}
%%%%%\MM^{(k)}_{j,\s,j',\s'}(x,c,\om,\h,n) &= c_{j}^{\s}c_{j'}^{-\s'} A^{(k)}_{jj'}(x,\s\s',|c|^2,\om,\h,n) + \de_{jj'}\de_{\s\s'} \matB^{(k)}_{jj'}(x,|c|^2,\om,\h,n), 
%%%%%%\\
%%%%%%\MM^{(k)}_{j,-\s,j',-\s'}(-x,c,\om,\h,n) &= c_{j}^{-\s}c_{j'}^{\s'} A^{(k)}_{jj'}(-x,\s\s',|c|^2,\om,\h,n) + \de_{jj'}\de_{\s\s'} \matB^{(k)}_{jj'}(-x,|c|^2,\om,\h,n),
%%%%%\end{aligned}
%%%%%\]
%%%%%and, if we set
%%%%%\[
%%%%%\begin{aligned}
%%%%%\calA^{(k)}_{jj'\s}(x,|c|^2,\om,\h,n) & = A^{(k)}_{jj'}(x,\s,|c|^2,\om,\h,n) , \\
%%%%%\BB^{(k)}_{jj'}(x,|c|^2,\om,\h,n) & = B^{(k)}_{jj'}(x,|c|^2,\om,\h,n) ,
%%%%%\end{aligned}
%%%%%\]
%%%%%then \eqref{above} follows.
\EP
%%%%%%%%%%%%%%%%%%%%%%%%%%%%%%%%%%%%%%%%%%%%%%%%%%%%%%%%%%%%%%%%%%%%%%%%%%

%%%%%%%%%%%%%%%%%%%%%%%%%%%%%%%%%%%%%%%%%%%%%%%%%%%%%%%%%%%%%%%%%%%%%%%%%%
\begin{rmk}\label{ham}
\emph{
The symmetry of Lemma \ref{pauli} is 
%a deep consequence of the ``real-on-real'' Hamiltonian nature of the problem, and it is
the analogous of \cite[Lemma 3.11]{GG} in complex instead of action-angle coordinates.
%It implies that, if $P_1$ is defined as in \eqref{matricette}, then
%
%\[
%\begin{aligned}
%P_1 \MM_{j j'}^{(k)}(x,c,\om,\h,n) P_1 & = \ol{\MM_{j j'}^{(k)}(-x,c,\om,\h,n)}, \\
%\CCCC_j\ol{\MM_{j j'}^{(k)}(x,c,\om,\h,n) } \, \CCCC_{j'}^\dagger & = \CCCC_j^\dagger \MM_{j j'}^{(k)}(x,c,\om,\h,n) \, \CCCC_{j'}.
%\end{aligned}
%\]
}
\end{rmk}
%%%%%%%%%%%%%%%%%%%%%%%%%%%%%%%%%%%%%%%%%%%%%%%%%%%%%%%%%%%%%%%%%%%%%%%%%%

%%%%%%%%%%%%%%%%%%%%%%%%%%%%%%%%%%%%%%%%%%%%%%%%%%%%%%%%%%%%%%%%%%%%%%%%%%
\begin{coro}\label{lineare}
For all $k\ge 1$, $j,j' \in\ZZZ$, $n\ge 0$, and for any $\h \in \RRR^{\ZZZ\times\NNN}$, there are
$\DD^{(k)}_{jj'\s}=\DD^{(k)}_{jj'\s}(|c|^2,\om,\h,n)\in\RRR$, 
with $\s=\pm$, and  $\calE^{(k)}_{jj}=\calE^{(k)}_{jj'}(|c|^2,\om,\h,n)\in\RRR$, such that
\begin{equation}\label{put}
\begin{aligned}
	\del_x \MM_{j j'}^{(k)}(0,c,\om,\h,n) = 
	\CCCC_j 
	&\begin{pmatrix}
	\DD^{(k)}_{jj'+} & -\DD^{(k)}_{jj'-} \cr
	\DD^{(k)}_{jj'-} & -\DD^{(k)}_{jj'+} \end{pmatrix}
	\CCCC_{j'}^\dagger %\\ &
	+ \de_{jj'} \calE^{(k)}_{jj} P_3.
\end{aligned}
\end{equation}
with $\CCCC_j$ and $P_3$ defined in \eqref{matricette}.
\end{coro}
%%%%%%%%%%%%%%%%%%%%%%%%%%%%%%%%%%%%%%%%%%%%%%%%%%%%%%%%%%%%%%%%%%%%%%%%%%

%%%%%%%%%%%%%%%%%%%%%%%%%%%%%%%%%%%%%%%%%%%%%%%%%%%%%%%%%%%%%%%%%%%%%%%%%%
\prova
It is a direct consequence of Lemma \ref{pauli}: differentiating the matrix in \eqref{above} w.r.t.~$x$ 
and computing at $x=0$, the assertion follows
with $\DD^{(k)}_{jj'\s} = \del_x \calA^{(k)}_{jj'\s}(0)$ and $\calE^{(k)}_{jj'} =\del_x \BB^{(k)}_{jj'}(0)$.
\EP
%%%%%%%%%%%%%%%%%%%%%%%%%%%%%%%%%%%%%%%%%%%%%%%%%%%%%%%%%%%%%%%%%%%%%%%%%%

%%%%%%%%%%%%%%%%%%%%%%%%%%%%%%%%%%%%%%%%%%%%%%%%%%%%%%%%%%%%%%%%%%%%%%%%%%
\begin{rmk}\label{rigiusto}
\emph{
Thanks to Remarks \ref{sing} and \ref{giusto}, the symmetry properties stated in Lemma \ref{pauli} -- and in Corollary \ref{lineare} -- hold
for ${\matM}^{(k)}_{j,\s,j',\s'}(x,c,\om,n)$ as well.
}
\end{rmk}
%%%%%%%%%%%%%%%%%%%%%%%%%%%%%%%%%%%%%%%%%%%%%%%%%%%%%%%%%%%%%%%%%%%%%%%%%%

Let $\TT_{j,\s}^{(k)}$ be the trivial RT $\TT\in\SSSS$ 
of order $k$ and with component label $j_{\ell_\TT^{out}}=j$ and sign label $\s_{\ell_\TT^{out}}=\s$ (see Definition \ref{RT} and Figure \ref{Tjsk}),
and define
\begin{equation}\label{splitto}
\tilde{\SSSS}^{(k)}_{j,\s,j',\s'}(n):=
\left\{
\begin{aligned}
&\SSSS^{(k)}_{j,\s,j,\s}(n)\setminus\{\TT_{j,\s}^{(k)}\} \quad \mbox{ if }j=j'\mbox{ and }\s=\s', \\
& \SSSS^{(k)}_{j,\s,j',\s'}(n) \quad\mbox{ otherwise},
\end{aligned}
\right.
\end{equation}
so that, recalling \eqref{emmenne},
\begin{equation}\label{separo}
\MM_{j\s j'\s'}^{(k)}(x,c,\om,\h,n) =  \de_{jj'}\de_{\s\s'}(\h^{(k)}_j)^\s + \tilde\MM_{j\s j'\s'}^{(k)}(x,c,\om,\h,n),
\end{equation}
where
\begin{equation}\label{emmetilde}
\tilde\MM_{j\s j'\s'}^{(k)}(x,c,\om,\h,n) := \sum_{\TT\in \tilde \SSSS^{(k)}_{j,\s,j',\s'}(n)}\VV(\TT;x,c,\om,\h).
\end{equation}

The next symmetry is the analogous of \cite[Lemma 4.12]{CG2} in the present context.

%%%%%%%%%%%%%%%%%%%%%%%%%%%%%%%%%%%%%%%%%%%%%%%%%%%%%%%%%%%%%%%%%%%%%%%%%%
\begin{prop}\label{emmeeladerivatadigi}
For all $k\ge 1$, $j,j' \in\ZZZ$, $n\ge-1$, and for any $\h \in \CCC^{\ZZZ\times\NNN}$, one has
\[
\MM_{j\s j'\s'}^{(k)}(0,c,\om,\h,n) =  \del_{c_{j'}^{\s'}} G_{j,\s}^{(k)}(c,\om,\h,n) ,
\]
with $G_{j,\s}^{(k)}(c,\om,\h,n)$ defined in \eqref{gin}.
\end{prop}
%%%%%%%%%%%%%%%%%%%%%%%%%%%%%%%%%%%%%%%%%%%%%%%%%%%%%%%%%%%%%%%%%%%%%%%%%%

%%%%%%%%%%%%%%%%%%%%%%%%%%%%%%%%%%%%%%%%%%%%%%%%%%%%%%%%%%%%%%%%%%%%%%%%%%
\prova
By \eqref{posso} and \eqref{cijei2}, one has
\[
\del_{c_{j'}^{\s'}} G_{j,\s}^{(k)}(c,\om,\h,n) = \de_{jj'}\de_{\s,\s'}( \h_j^{(k)})^\s + 
\del_{c_{j'}^{\s'}} (\tilde{G}_{j,\s}^{(k)}(c,\om,\h,n) )\,,%c_j^\s \hat{G}_{j}^{(k)}(|c|^2,\h))\,,
\]
so that by \eqref{separo} we only need to prove that
\[
\tilde\MM_{j\s j'\s'}^{(k)}(0,c,\om,\h,n) =  \del_{c_{j'}^{\s'}}\tilde{G}_{j,\s}^{(k)}(c,\om,\h,n) .
\]

%Define $\del_{c_{j'}^{\s'}} \tilde{\Theta}_{j,\gote_j,\s}^{(k)}(n)$ as the set of all trees obtained from the 
Given a tree $\vartheta\in\tilde{\Theta}_{j,\gote_j,\s}^{(k)}(n)$ such that $|\Lambda_{j',\s'}(\vartheta)|\ge1$,
choose a leaf $\la\in \Lambda_{j',\s'}(\vartheta) $ %(recall Definition \ref{indicetti})
and change its value
into $\LL_\la(c)=1$, so that the value of the tree is modified as well into a new value $\VV_\la(\vartheta;c,\om,\h)$;
we call such a leaf the \emph{fallen leaf} of the tree. By \eqref{tilden}  we have
\begin{equation}\label{low}
\null\hspace{-.3cm}
 \del_{c_{j'}^{\s'}}\tilde{G}_{j,\s}^{(k)}(c,\om,\h,n) %=\del_{c_{j'}^{\s'}} ( c_j^\s \hat{G}_{j}^{(k)}(|c|^2,\om,\h,n)) 
 = \sum_{\vartheta\in \tilde{\Theta}_{j,\gote_j,\s}^{(k)}(n)}\; \sum_{\la \in \Lambda_{j',\s'}(\vartheta)} \VV_\la (\vartheta;c,\om,\h).
\end{equation}
Clearly $ \tilde{\SSSS}^{(k)}_{j,\s,j',\s'}(n)$ is in one-to-one correspondence with the set of trees
$\vartheta\in\tilde{\Theta}_{j,\gote_j,\s}^{(k)}(n)$ with a fallen leaf: this is easily seen
%indeed each $\vartheta \in \del_{c_{j'}^{\s'}} \tilde{\Theta}_{j,\gote_j,\s}^{(k)}(n)$ corresponds to a RT $\TT\in \tilde \SSSS^{(k)}_{j,\s,j',\s'}(n)$
by identifying the root of $\vartheta$ with $v_{out}$ of the corresponding $\TT$, and the fallen leaf in $\vartheta$ with $v_{in}$ of $\TT$.
Finally for any 
$\vartheta,\TT$ as above one has  $\VV(\vartheta;c,\om,\h)=\VV(\TT;0,c,\om,\h)$, so that the assertion follows.
\EP
%%%%%%%%%%%%%%%%%%%%%%%%%%%%%%%%%%%%%%%%%%%%%%%%%%%%%%%%%%%%%%%%%%%%%%%%%%

A graphical representation of the identity in Proposition \ref{emmeeladerivatadigi}, for $k=1$, is provided in Figure \ref{figlow}.
If $\vartheta \in \tilde{\Theta}_{j,\gote_j,\s}^{(1)}(n)$ is such that
either $|\Lambda_{j',\s'}(\vartheta)|=1$ or $|\Lambda_{j',\s'}(\vartheta)|=2$, then the derivative
w.r.t.~$c_{j'}^{\s'}$, when acting on $\VV (\vartheta;c,\om,\h)$, produces one or two contributions, respectively,
to $\tilde\MM_{j\s j'\s'}^{(1)}(0,c,\om,\h,n)$. Note that, for $k=1$, $|\Lambda_{j',\s'}(\vartheta)|=3$ is possible only if $(j',\s')=(j,\s)$,
while $|\Lambda_{j',\s'}(\vartheta)|\ge 4$ is not possible at all.

%%%%%%%%%%%%%%%%%%%%%%%%%%%%%%%%%%%%%%%%%%%%%%%%%%%%%%%%%%%%%%%%%%%%%%%%%%
% FIGURA 15
%%%%%%%%%%%%%%%%%%%%%%%%%%%%%%%%%%%%%%%%%%%%%%%%%%%%%%%%%%%%%%%%%%%%%%%%%%
\begin{figure}[ht]
\vspace{.4cm}
\centering
\ins{220pt}{-032pt}{$= \; \partial_{c_{j'}^{\s'}}$}
\ins{278pt}{-115pt}{$= \;\;\; \partial_{c_{j'}^{\s'}}$}
\ins{328pt}{-032pt}{$j',\s'$}
\ins{396pt}{-115pt}{$j',\s'$}
\ins{396pt}{-103pt}{$j',\s'$}
\ins{142pt}{-115pt}{$+$}
%\ins{212pt}{0012pt}{$k,j,\s$}
%\ins{180pt}{-012pt}{$j,\gote_j,\s$}
%\ins{248pt}{-012pt}{$j,\gote_j,\s$}
\subfigure{\includegraphics*[width=1.6in]{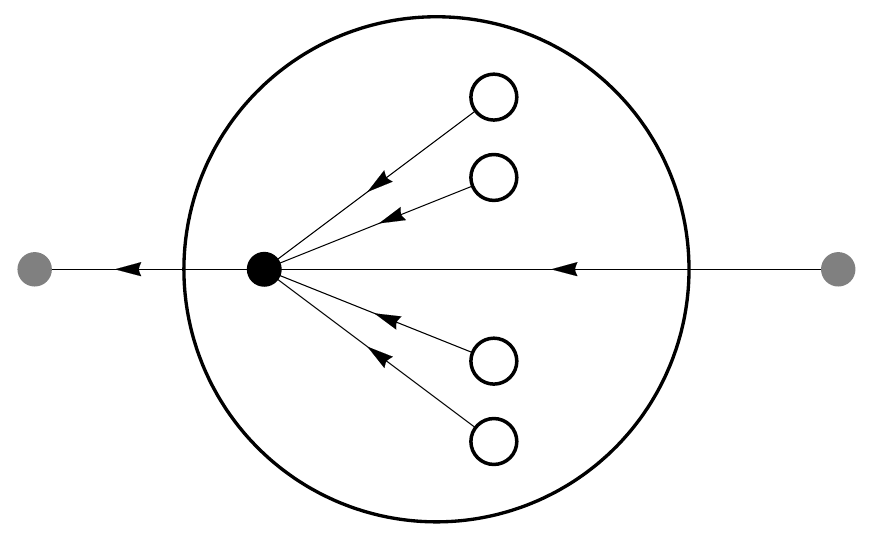}}
\hspace{1.5cm}
\subfigure{\includegraphics*[width=1.6in]{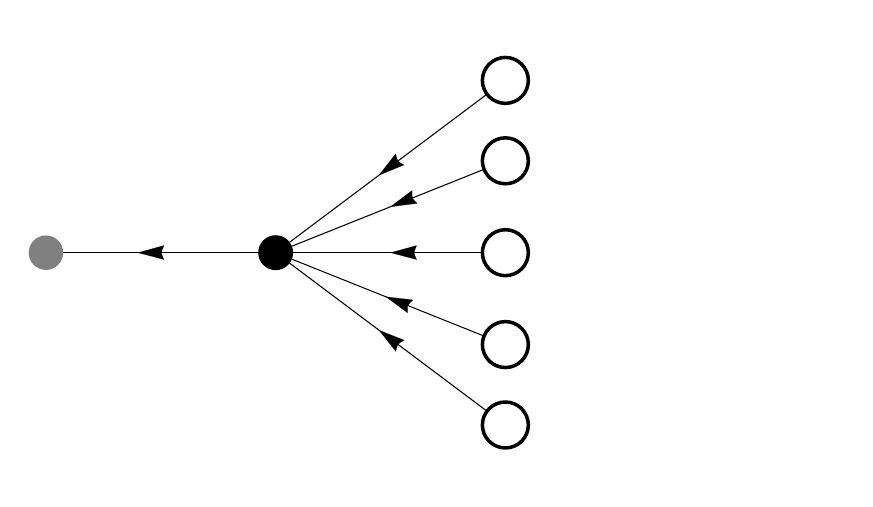}}
\\
\subfigure{\includegraphics*[width=1.6in]{figura15a}}
\hspace{.5cm}
\subfigure{\includegraphics*[width=1.6in]{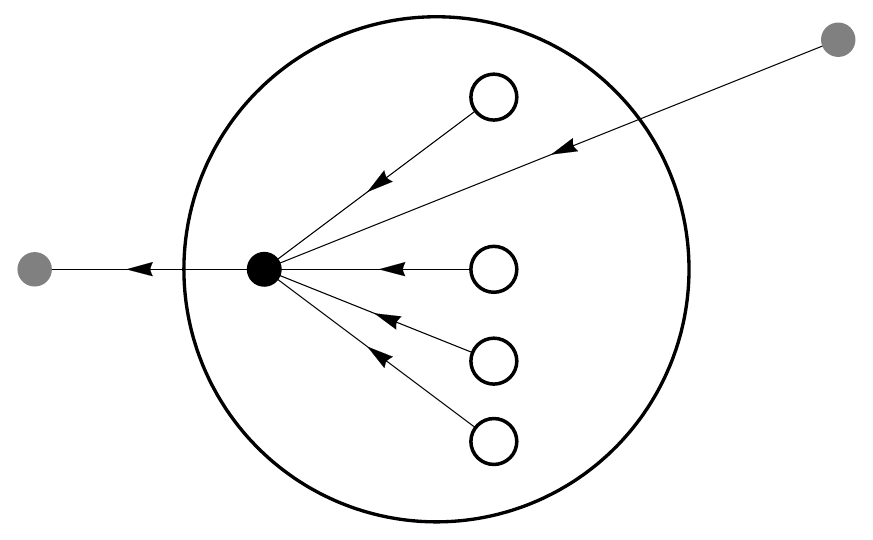}}
\hspace{1.5cm}
\subfigure{\includegraphics*[width=1.6in]{figura15b}}
\caption{\small Graphical representation of the identity in Proposition \ref{emmeeladerivatadigi} for $k=1$.}
\label{figlow}
\end{figure}
%%%%%%%%%%%%%%%%%%%%%%%%%%%%%%%%%%%%%%%%%%%%%%%%%%%%%%%%%%%%%%%%%%%%%%%%%%

%%%%%%%%%%%%%%%%%%%%%%%%%%%%%%%%%%%%%%%%%%%%%%%%%%%%%%%%%%%%%%%%%%%%%%%%%%
\begin{rmk}\label{primaderivo}
\emph{
From Proposition \ref{emmeeladerivatadigi} and Remark \ref{giusto}, it follows immediately  that
\[
\matM_{j\s j'\s'}^{(k)}(0,c,\om,n) = \left.\del_{c_{j'}^{\s'}}G^{(k)}_{j,\s}(c,\om,\h,n)\right|_{\h=\h(c,\om,n)}\, .
\]
}
\end{rmk}
%%%%%%%%%%%%%%%%%%%%%%%%%%%%%%%%%%%%%%%%%%%%%%%%%%%%%%%%%%%%%%%%%%%%%%%%%%

%%%%%%%%%%%%%%%%%%%%%%%%%%%%%%%%%%%%%%%%%%%%%%%%%%%%%%%%%%%%%%%%%%%%%%%%%%
\begin{coro}\label{emmecor}
For all $k\ge 1$, $j,j' \in\ZZZ$, $\s,\s'\in\{\pm\}$ and $n\ge -1$, and for any $\h \in \CCC^{\ZZZ\times\NNN}$, one has
\[
\MM_{j\s j'\s'}^{(k)}(0,c,\om,\h,n) = \de_{\s\s'}\de_{jj'} \bigl( (\h^{(k)}_j)^\s +\widehat{G}_{j}^{(k)}(|c|^2,\om,\h,n) \bigr) + 
c_j^\s c_{j'}^{-\s'}\del_{|c_{j'}|^2} \widehat{G}^{(k)}_{j}(|c|^2,\om,\h,n) .
% F^{(k)}_{j,j'}(|c|^2,\h)
\]
%for some function $F^{(k)}_{j,j'}(|c|^2,\om,\h,n)$ is a polynomial in $c_j^\s$, $ c_{j'}^{-\s'}$. 
\end{coro}
%%%%%%%%%%%%%%%%%%%%%%%%%%%%%%%%%%%%%%%%%%%%%%%%%%%%%%%%%%%%%%%%%%%%%%%%%%

%%%%%%%%%%%%%%%%%%%%%%%%%%%%%%%%%%%%%%%%%%%%%%%%%%%%%%%%%%%%%%%%%%%%%%%%%%
\prova
Proposition \ref{emmeeladerivatadigi}, combined with  \eqref{cijei2} and Remark \ref{posso}  %More specifically, the case $j\ne j'$ is obvious.
implies
\[
\MM_{j\s j'\s'}^{(k)}(0,c,\om,\h,n) = \de_{\s\s'}\de_{jj'} \bigl( (\h^{(k)}_j)^{\s} +\widehat{G}_{j}^{(k)}(|c|^2,\om,\h,n) \bigr) + 
c_j^\s \del_{c_{j'}^{\s'}} \widehat{G}^{(k)}_{j}(|c|^2,\om,\h,n),
\]
so that the assertion follows.
\EP
%%%%%%%%%%%%%%%%%%%%%%%%%%%%%%%%%%%%%%%%%%%%%%%%%%%%%%%%%%%%%%%%%%%%%%%%%%

From Remark \ref{primaderivo} and Corollary \ref{emmecor} we deduce the following result.

%%%%%%%%%%%%%%%%%%%%%%%%%%%%%%%%%%%%%%%%%%%%%%%%%%%%%%%%%%%%%%%%%%%%%%%%%%
\begin{coro}\label{simmetria}
For all $k\ge 1$, $j,j' \in\ZZZ$, $\s,\s'\in\{\pm\}$ and $n\ge -1$,  one has
\[
\matM_{j\s j'\s'}^{(k)}(0,c,\om,n) = c_{j}^{\s}c_{j'}^{-\s'} F_{j,j'}^{(k)}(|c|^2,\om,n),
\]
with
\[
F_{j,j'}^{(k)}(|c|^2,\om,n):= \left. \del_{|c_{j'}|^2} \widehat{G}^{(k)}_{j}(|c|^2,\om,\h,n)\right|_{\h=\h(c,\om,n)}.
\]
%with  $ F_{j,j'}(|c|^2,\h(c))$ having no singularity in $c_j^\s$, $ c_{j'}^{-\s'}$.
\end{coro}
%%%%%%%%%%%%%%%%%%%%%%%%%%%%%%%%%%%%%%%%%%%%%%%%%%%%%%%%%%%%%%%%%%%%%%%%%%

%%%%%%%%%%%%%%%%%%%%%%%%%%%%%%%%%%%%%%%%%%%%%%%%%%%%%%%%%%%%%%%%%%%%%%%%%%
\prova
By  Corollary \ref{emmecor} and Remark \ref{giusto} we have
\[
\begin{aligned}
\matM_{j\s j'\s'}^{(k)}(0,c,\om,n) =  \de_{\s\s'}\de_{jj'} &( \h^{(k)}_j(c,\om,n) +\widehat{G}_{j}^{(k)}(|c|^2,\h(c,\om,n),n) ) \\
&+
 \left.c_{j}^{\s}c_{j'}^{-\s'} \del_{|c_{j'}|^2} \widehat{G}^{(k)}_{j}(|c|^2,\om,\h,n)\right|_{\h=\h(c,\om,n)},
 \end{aligned}
\]
so that the assertion follows by \eqref{labase}.
\EP 
%%%%%%%%%%%%%%%%%%%%%%%%%%%%%%%%%%%%%%%%%%%%%%%%%%%%%%%%%%%%%%%%%%%%%%%%%%

%%%%%%%%%%%%%%%%%%%%%%%%%%%%%%%%%%%%%%%%%%%%%%%%%%%%%%%%%%%%%%%%%%%%%%%%%%
\begin{rmk} \label{pappapero}
\emph{
Both $\matM_{j\s j'\s'}^{(k)}(0,c,\om,n)$ and $\widehat{G}^{(k)}_{j}(|c|^2,\om,\h,n)$ can be expressed as sum of values of trees;
this follows from \eqref{emmeta} %and Remark \ref{giusto} 
for the first quantity and from \eqref{cijei2}
with \eqref{tilden} for the latter. Then Corollary \ref{simmetria} shows that, if we
\begin{itemize}[topsep=0.1ex]
\itemsep-0.2em
\item differentiate w.r.t.~$|c_{j'}|^2$ the value of a tree $\vartheta$ contributing to $\widehat{G}^{(k)}_{j}(|c|^2,\om,\h,n)$,
\item replace $\h^{(k_v)}$ with the value of a tree contributing to $\h^{(k_v)}(c,\om,n)$ for any $v\in N(\vartheta)$ with $s_v=1$,
%\item expand $\h(c,\om,n)$ according to \eqref{formale} and take a contribution of order $k_2=k-k_1$,
\end{itemize} 
we obtain one or more quantities which, multiplied times $c_{j}^{\s}c_{j'}^{-\s'}$, can be seen as values of trees contributing to
$\matM_{j\s j'\s'}^{(k)}(0,c,\om,n)$.
}
\end{rmk}
%%%%%%%%%%%%%%%%%%%%%%%%%%%%%%%%%%%%%%%%%%%%%%%%%%%%%%%%%%%%%%%%%%%%%%%%%%

As an example, consider the tree $\vartheta$ contributing to $c_j^\s\widehat{G}^{(4)}_{j}(|c|^2,\om,\h,n)$
represented in Figure \ref{fig-simmetria0},
where only the labels equal to $j,\pm\s$ and $j',\pm\s'$ are explicitly shown for the leaves, with $\s=\s'=+$ and $j\neq j'$,
and the node $v$ with $s_v=1$ has order label $k_v=1$.

%%%%%%%%%%%%%%%%%%%%%%%%%%%%%%%%%%%%%%%%%%%%%%%%%%%%%%%%%%%%%%%%%%%%%%%%%%
% FIGURA 17boh
%%%%%%%%%%%%%%%%%%%%%%%%%%%%%%%%%%%%%%%%%%%%%%%%%%%%%%%%%%%%%%%%%%%%%%%%%%
\begin{figure}[ht]
\centering
\ins{083pt}{-048pt}{$\vartheta=$}
\ins{202pt}{-075pt}{$j+$}
\ins{200pt}{-098pt}{$j'\!+$}
\ins{204pt}{000pt}{$j'\!-$}
\ins{314pt}{-098pt}{$j'\!-$}
\ins{356pt}{-066pt}{$j'\!+$}
\subfigure{\includegraphics*[width=3.6in]{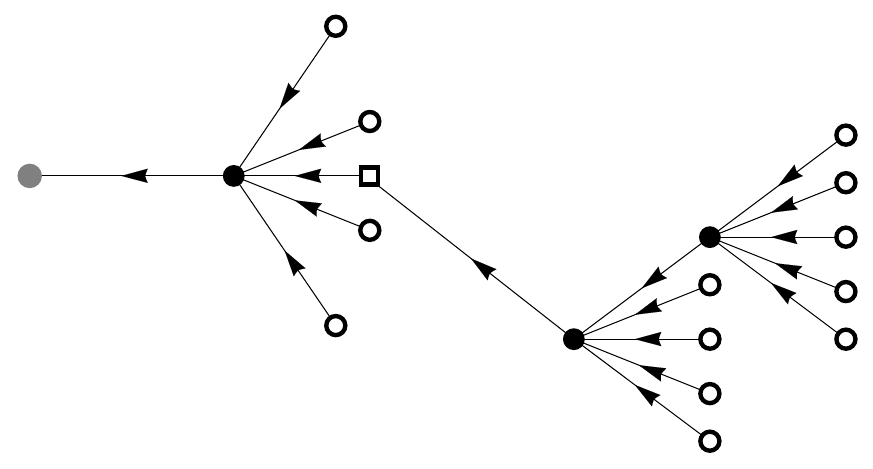}}
\caption{\small A tree $\vartheta$ whose value contributes to $c_j^\s\widehat{G}^{(4)}_{j}(|c|^2,\om,\h,n)$.}
\label{fig-simmetria0}
\end{figure}
%%%%%%%%%%%%%%%%%%%%%%%%%%%%%%%%%%%%%%%%%%%%%%%%%%%%%%%%%%%%%%%%%%%%%%%%%%

%%%%%%%%%%%%%%%%%%%%%%%%%%%%%%%%%%%%%%%%%%%%%%%%%%%%%%%%%%%%%%%%%%%%%%%%%%%
%% FIGURA 17a
%%%%%%%%%%%%%%%%%%%%%%%%%%%%%%%%%%%%%%%%%%%%%%%%%%%%%%%%%%%%%%%%%%%%%%%%%%%
%\begin{figure}[ht]
%\centering
%\ins{083pt}{-073pt}{$\vartheta=$}
%\ins{202pt}{-100pt}{$j\s$}
%\ins{200pt}{-123pt}{$j'\!\s'$}
%\ins{200pt}{-017pt}{$j'\!\!-\!\!\s'$}
%\ins{314pt}{-123pt}{$j'\!\!-\!\!\s'$}
%\ins{356pt}{-091pt}{$j'\!\s'$}
%\subfigure{\includegraphics*[width=3.6in]{figura17a}}
%\caption{\small A tree $\vartheta$ whose value contributes to $c_j^\s\widehat{G}^{(k)}_{j}(|c|^2,\om,\h,n)$ with $k=4$.}
%\label{fig-simmetria1}
%\end{figure}
%%%%%%%%%%%%%%%%%%%%%%%%%%%%%%%%%%%%%%%%%%%%%%%%%%%%%%%%%%%%%%%%%%%%%%%%%%%

Then, if we differentiate $\VV(\vartheta;c,\om,\h)$ with respect to $|c_{j'}|^2$, hence replace $\h^{(1)}$ with the value of a tree of order $1$
contributing to $\h^{(1)}(c,\om,n)$ and thence multiply the result times $c_j^\s c_{j'}^{-\s'}$,
we obtain the values, computed at $x=0$, of the two resonant trees $\TT_1$ and $\TT_2$
represented in Figure \ref{fig-simmetria2}: such values contribute to $\matM_{j\s j'\s'}^{(4)}(0,c,\om,n)$.
Note that the subtree whose value contributes to $\h^{(1)}(c,\om,n)$ -- the subgraph encircled by the dashed line, including the exiting line,
in Figure \ref{fig-simmetria2} -- may contain leaves with labels $j'$.
%
%%%%%%%%%%%%%%%%%%%%%%%%%%%%%%%%%%%%%%%%%%%%%%%%%%%%%%%%%%%%%%%%%%%%%%%%%%
% FIGURA 17b
%%%%%%%%%%%%%%%%%%%%%%%%%%%%%%%%%%%%%%%%%%%%%%%%%%%%%%%%%%%%%%%%%%%%%%%%%%
\begin{figure}[ht]
\centering
\ins{009pt}{-051pt}{$\TT_1\!=$}
\ins{235pt}{-053pt}{$\TT_2\!=$}
\ins{028pt}{-040pt}{$j+$}
\ins{256pt}{-042pt}{$j+$}
\ins{217pt}{-054pt}{$j'\!+$}
\ins{099pt}{-018pt}{$j'\!-$}
\ins{325pt}{-015pt}{$j'\!-$}
\ins{334pt}{-098pt}{$j'\!+$}
\ins{098pt}{-082pt}{$j'\!+$}
\ins{426pt}{-063pt}{$j'\!+$}
\ins{102pt}{-070pt}{$j+$}
\ins{323pt}{-072pt}{$j+$}
\ins{173pt}{-082pt}{$j'\!-$}
\ins{400pt}{-082pt}{$j'\!-$}
%\ins{212pt}{0012pt}{$k,j,\s$}
%\ins{180pt}{-012pt}{$j,\gote_j,\s$}
%\ins{248pt}{-012pt}{$j,\gote_j,\s$}
\subfigure{\includegraphics*[width=2.8in]{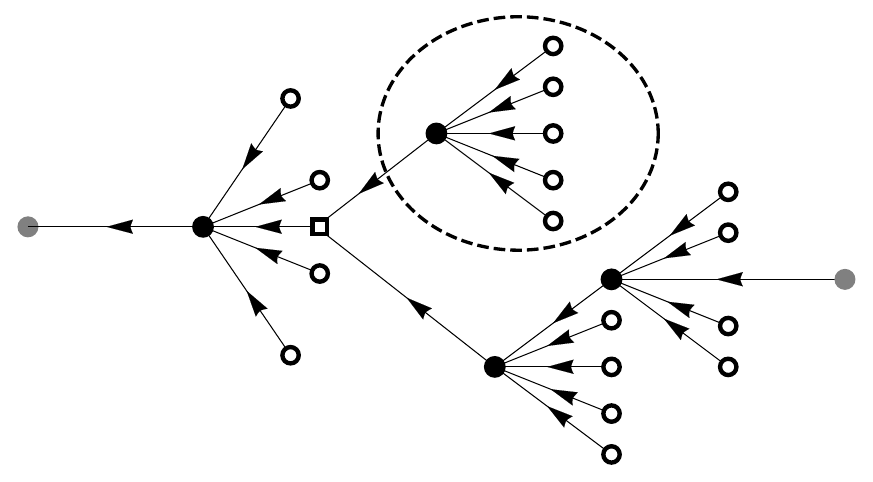}}
\hspace{0.6cm}
\subfigure{\includegraphics*[width=2.4in]{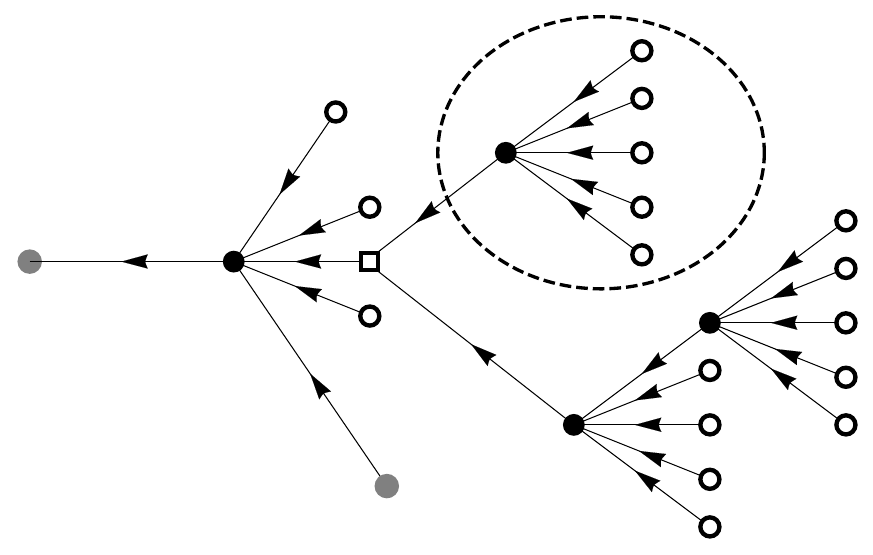}}
\caption{\small The resonant trees obtained from the tree $\vartheta$ in Figure \ref{fig-simmetria0}.}
\label{fig-simmetria2}
\end{figure}
%%%%%%%%%%%%%%%%%%%%%%%%%%%%%%%%%%%%%%%%%%%%%%%%%%%%%%%%%%%%%%%%%%%%%%%%%%

%%%%%%%%%%%%%%%%%%%%%%%%%%%%%%%%%%%%%%%%%%%%%%%%%%%%%%%%%%%%%%%%%%%%%%%%%%
\begin{coro}\label{pippero}
For all $k\ge 1$, $j,j' \in\ZZZ$ and $n\ge -1$, one has
\begin{equation}\label{lamad}
\matM_{j j'}^{(k)}(0,c,\om,n) =  F_{j,j'}^{(k)}(|c|^2,\om,n)
\begin{pmatrix}
c_j \ol{c}_{j'} &  c_jc_{j'} \cr
\ol{c}_j\ol{c}_{j'} & \ol{c}_j c_{j'}
\end{pmatrix}
=  F_{j,j'}^{(k)}(|c|^2,\om,n)\, \CCCC_j U \CCCC_{j'}^\dagger.
\end{equation}
\end{coro}
%%%%%%%%%%%%%%%%%%%%%%%%%%%%%%%%%%%%%%%%%%%%%%%%%%%%%%%%%%%%%%%%%%%%%%%%%%

%%%%%%%%%%%%%%%%%%%%%%%%%%%%%%%%%%%%%%%%%%%%%%%%%%%%%%%%%%%%%%%%%%%%%%%%%%
\proof
By explicit calculation, combining the definitions in \eqref{matricette} and Corollary \ref{simmetria}.
\EP
%%%%%%%%%%%%%%%%%%%%%%%%%%%%%%%%%%%%%%%%%%%%%%%%%%%%%%%%%%%%%%%%%%%%%%%%%%

%%%%%%%%%%%%%%%%%%%%%%%%%%%%%%%%%%%%%%%%%%%%%%%%%%%%%%%%%%%%%%%%%%%%%%%%%%
\begin{rmk} \label{notiamolo}
\emph{%remark
By comparing Corollary \ref{pippero} with Lemma \ref{pauli}, we deduce that
\[
\calA^{(k)}_{jj'+}(0) = \calA^{(k)}_{jj'-}(0) = F_{j,j'}^{(k)}(|c|^2,\om,n)
\]
when both $\calA^{(k)}_{jj'+}(0)$ and $\calA^{(k)}_{jj'-}(0)$ are computed at $\h=\h(c,\om,n)$.
}%remark
\end{rmk}
%%%%%%%%%%%%%%%%%%%%%%%%%%%%%%%%%%%%%%%%%%%%%%%%%%%%%%%%%%%%%%%%%%%%%%%%%%

%%%%%%%%%%%%%%%%%%%%%%%%%%%%%%%%%%%%%%%%%%%%%%%%%%%%%%%%%%%%%%%%%%%%%%%%%%
\begin{lemma}\label{cancellazione}
For all $p\ge2$ and all $k_1,\ldots,k_{p}\ge 1$, $n_1,\ldots,n_p \ge -1$, $j_1,\ldots,j_{p+1}\in\ZZZ$, one has
\begin{equation}\label{lei}
	\matM_{j_1 j_2}^{(k_1)}(0,c,\om,n_1)\,\Biggl(\, \prod_{i=2}^{p-1} P_3 \del_x\matM_{j_i j_{i+1}}^{(k_i)}(0,c,\om,n_i)\Biggr)\,P_3\, 
	\matM_{j_p j_{p+1}}^{(k_p)}(0,c,\om,n_p) = 0,
\end{equation}
where the product has to be interpreted as $1$ for $p=2$.
\end{lemma}
%%%%%%%%%%%%%%%%%%%%%%%%%%%%%%%%%%%%%%%%%%%%%%%%%%%%%%%%%%%%%%%%%%%%%%%%%%

%%%%%%%%%%%%%%%%%%%%%%%%%%%%%%%%%%%%%%%%%%%%%%%%%%%%%%%%%%%%%%%%%%%%%%%%%%
\proof
Combining Corollary \ref{pippero} and Lemma \ref{lineare}, since $\CCCC_j\CCCC_j^\dagger = |c_j|^2\uno$, we see that \eqref{lei} 
can be written as
\[
\g \, \CCCC_{j_1}
 U \Biggl(\,\prod_{i=2}^{p-1}P_3 D_i\Biggr) P_3 U \CCCC_{j_{p+1}}^\dagger =0,\qquad
D_i:=\begin{pmatrix}
a_i & -b_i \cr b_i & -a_i
\end{pmatrix},
\] 
where we have used the notation in \eqref{matricette} and defined
\[
\begin{aligned}
a_i&:=\DD^{(k_i)}_{j_i j_{i+1}+}(|c|^2,\om,\h(c,\om,n_i),n_i)|c_{j_i}|^2 + \de_{j_i j_{i+1}} \calE^{(k_i)}_{j_i j_{i+1}}(|c|^2,\om,\h(c,\om,n_i),n_i), \\
b_i&:= \DD^{(k_i)}_{j_i j_{i+1}-}(|c|^2,\om,\h(c,\om,n_i),n_i)|c_{j_i}|^2,\\
\g&:= F_{j_1j_2}^{(k_1)}(|c|^2,\om,n_1)  F_{j_pj_{p+1}}^{(k_p)}(|c|^2,\om,n_p) ,
\end{aligned}
\]
with $\DD^{(k)}_{jj'\s}$ and $\calE^{(k)}_{jj'\s}$ defined in Corollary \ref{lineare}, and $F^{(k)}_{jj'}$ as in Corollaruy \ref{simmetria}.

The case $p=2$ follows from the fact that $UP_3U=0$.
If $p\ge3$, by explicit calculation 
\[
UP_3D_i = (a_i-b_i) U,\qquad D_i P_3 U = (a_i+b_i)U,
\]
so that the result follows once more from the fact that $UP_3U=0$.
\EP
%%%%%%%%%%%%%%%%%%%%%%%%%%%%%%%%%%%%%%%%%%%%%%%%%%%%%%%%%%%%%%%%%%%%%%%%%%

%%%%%%%%%%%%%%%%%%%%%%%%%%%%%%%%%%%%%%%%%%%%%%%%%%%%%%%%%%%%%%%%%%%%%%%%%%
\begin{rmk} \label{oraeunremark}
\emph{%remark
For any $C^2$ function $f\!:\RRR\to\RRR$ we can write
\begin{equation}\label{taylor}
f(x)= f(0) + x \,\del f(0) + x^2\int_0^1 {\rm d}y\, (1-y) \, \del^2 f(yx),
\end{equation}
where $\del$ is the derivative w.r.t.~the argument.
%With this in mind
We introduce the \emph{localization operator} $\matL$ by setting
\begin{equation}\label{elle}
\matL f(x) := f(0) ,
\end{equation}
the \emph{derivative operator} $\matD$ by setting
\begin{equation}\label{der}
\matD f(x) := x \, \del f(0)  ,
\end{equation}
and the \emph{regularization operator} $\matR$ by setting 
\begin{equation}\label{erre}
\matR f(x) :=  x^2\int_0^1 {\rm d}y \,(1-y) \, \del^2 f(yx) , 
\end{equation}
so that $\matL+\matD+\matR=\uno$; such operators will be widely used in the rest of the section.
Then, we may reformulate \eqref{lei} as
\begin{equation}\label{leibis}
	\matL\matM_{j_1 j_2}^{(k_1)}(x,c,\om,n_1)\Bigl(\, \prod_{i=2}^{p-1}P_3\matD\matM_{j_i j_{i+1}}^{(k_i)}(x,c,\om,n_i)\Bigr)\, 
	 P_3\, 
	\matL\matM_{j_p j_{p+1}}^{(k_p)}(x,c,\om,n_p) = 0 ,
\end{equation}
where $\matO\matM_{j_i j_{i+1}}^{(k_i)}(x,c,\om,n_i)$
are obtained, for $\matO=\matL,\matD,\matR$, by applying the linear operators introduced above to the function 
$x \mapsto \matM_{j_i j_{i+1}}^{(k_i)}(x,c,\om,n_i)$.
The operators introducedabove will be widely used in the two forthcoming sections.
}%remark
\end{rmk}
\subsection{Renormalization procedure:
renormalized trees}
\label{LR}
%%%%%%%%%%%%%%%%%%%%%%%%%%%%%%%%%%%%%%%%%%%%%%%%%%%%%%%%%%%%%%%%%%%%%%%%%%
%%%%%%%%%%%%%%%%%%%%%%%%%%%%%%%%%%%%%%%%%%%%%%%%%%%%%%%%%%%%%%%%%%%%%%%%%%

Now we come back to the problem of proving a bound like \eqref{cespero}. The intuitive idea would be to consider a second order Taylor expansion
of the matrices appearing in \eqref{cespero}, as suggested by Remark \ref{oraeunremark},
and make use of Lemma \ref{cancellazione} to cancel out enough
summands and eventually obtain the bound. However, in order to make this argument rigorous, we need a different, iterative procedure.
%which allows us to deal with relevant resonant clusters containing other relevant resonant clusters.
We start by giving some preliminary definitions.
%Recall Definition \ref{rilevante} and, given any subgraph $S$ of a tree $\vartheta$, define $\fT_0(S)$ as the set of relevant RCs contained in $S$. 

%%%%%%%%%%%%%%%%%%%%%%%%%%%%%%%%%%%%%%%%%%%%%%%%%%%%%%%%%%%%%%%%%%%%%%%%%%
\begin{defi}[\textbf{Depth of a relevant resonant cluster}]
\label{profondita}
Given a tree $\vartheta\in\gotT\cup\breve\gotT$, we define the \emph{depth} $d(T)$ of a relevant RC $T$ in $\vartheta$ 
recursively as follows: we set $d(T)=0$ if there is no relevant RC containing $T$, and set
$d(T)=d(T')+1$ if $T$ is contained in a relevant RC $T'$ and no other relevant
RC contained in $T'$ (if any) contains $T$.
\end{defi}
%%%%%%%%%%%%%%%%%%%%%%%%%%%%%%%%%%%%%%%%%%%%%%%%%%%%%%%%%%%%%%%%%%%%%%%%%%

%%%%%%%%%%%%%%%%%%%%%%%%%%%%%%%%%%%%%%%%%%%%%%%%%%%%%%%%%%%%%%%%%%%%%%%%%%
\begin{rmk}\label{catenadepth}
\emph{
If $\gotC=\{T_1,\ldots,T_p\}$ is a chain of RCs, then $d(T_1)=\ldots=d(T_p)$.
Moreover, if $T$ is contained in $T'$, then either $\calP_{T}\subset \calP_{T'}$ or $L(T) \cap \calP_{T'}=\emptyset$.
}
\end{rmk}
%%%%%%%%%%%%%%%%%%%%%%%%%%%%%%%%%%%%%%%%%%%%%%%%%%%%%%%%%%%%%%%%%%%%%%%%%%

%%%%%%%%%%%%%%%%%%%%%%%%%%%%%%%%%%%%%%%%%%%%%%%%%%%%%%%%%%%%%%%%%%%%%%%%%%
\begin{defi}[\textbf{Depth of a line}]
\label{profondita2}
Given a tree $\vartheta\in\gotT\cup\breve\gotT$, 
we say that a line $\ell \in L(\vartheta)$ has \emph{depth} $d(\ell)=0$ if there is no relevant RC $T$ such that $\ell\in\calP_T$,
otherwise if $T_{d-1}\subset T_{d-2}\subset\ldots\subset T_0$ are all and only the relevant RCs such that 
$\ell\in\calP_{T_i} $ for $i=0,\ldots,d-1$ we set $d(\ell)=d\ge1$.
\end{defi}
%%%%%%%%%%%%%%%%%%%%%%%%%%%%%%%%%%%%%%%%%%%%%%%%%%%%%%%%%%%%%%%%%%%%%%%%%%

%%%%%%%%%%%%%%%%%%%%%%%%%%%%%%%%%%%%%%%%%%%%%%%%%%%%%%%%%%%%%%%%%%%%%%%%%%
\begin{defi}[\textbf{Cloud of a line}]
\label{nuvola}
Given a tree $\vartheta\in\gotT\cup\breve\gotT$ and a line $\ell \in L(\vartheta)$,
if $d(\ell)\ge1$ we define the \emph{cloud} of $\ell$ as the set 
$\calmC_{\ell}(\vartheta):= \{T_0,\ldots,T_{d(\ell)-1}\}$ of the relevant RCs such that $\ell\in\calP_{T_i} $ for $i=0,\ldots,d-1$,
and if $d(\ell)=0$ we set $\calmC_\ell(\vartheta):=\emptyset$.
Given a relevant RC $T$ in $\vartheta$ and a line $\ell\in \calP_T$, we call $\calmC_{\ell}(T):=\{T'\in\calmC_\ell(\vartheta): d(T')>d(T)\}$
the \emph{cloud of $\ell$ in $T$}.
%%
%%
%%
%%Given a relevant RC $T$ and a line $\ell\in \calP_T$ with depth $d(\ell)=d\ge1$,
%%let $T_{0},\ldots,T_{p-1}$, with $p\le d(\ell)$ be all and only the relevant RCs contained in $T$ such that
%% such that $T=T_{0} \supset T_{1} \supset \ldots \supset
%%T_{p-1}$ and $\ell\in \calP_{T_i}$ for $i=0,\ldots,p-1$ ;
%%we call $\calmC_{\ell}(T):= \{T_1,\ldots,T_{p-1}\}$ %\{T_{j}\}_{j=1}^{p-1}$
%%the \emph{cloud} of $\ell$ in $T$. If $p=1$ the  cloud of $\ell$ in $T$ is the empty set.
% and $\{T_{j}\}_{j=1}^{p-1}$ the \emph{internal cloud} of $\ell$.
%If $p=1$ then the internal cloud of $\ell$ is the empty set.
 \end{defi}
 %%%%%%%%%%%%%%%%%%%%%%%%%%%%%%%%%%%%%%%%%%%%%%%%%%%%%%%%%%%%%%%%%%%%%%%%%%

%%%%%%%%%%%%%%%%%%%%%%%%%%%%%%%%%%%%%%%%%%%%%%%%%%%%%%%%%%%%%%%%%%%%%%%%%%
\begin{rmk}\label{sbuffo}
\emph{
A line $\ell$ with depth $d\ge1$ is contained in $p\ge d$ relevant RCs, and it belongs to the resonant path of $d$ of them.
}
\end{rmk}
%%%%%%%%%%%%%%%%%%%%%%%%%%%%%%%%%%%%%%%%%%%%%%%%%%%%%%%%%%%%%%%%%%%%%%%%%%

In order to show that the cancellation mechanism described in Subsection \ref{simmetrie} still works
when dealing with relevant resonant clusters containing other relevant resonant clusters,
%make rigorous the argument outlined in Subsection \ref{nuovo}, 
it is useful to introduce a few new sets of trees.

%%%%%%%%%%%%%%%%%%%%%%%%%%%%%%%%%%%%%%%%%%%%%%%%%%%%%%%%%%%%%%%%%%%%%%%%%%
\begin{defi}[\textbf{Set $\boldsymbol{\fR}$ of the renormalized trees}] %\hspace{-.27cm}\boldsymbol{\fR}$}]
\label{nonepossopiu}
Let $\fR$ be the set of trees which can be obtained from any expanded tree  $\vartheta\in\gotT$  by assigning 
to each relevant RC $T$ in $\vartheta$ an \emph{operator label} $\matO_T\in\{\matL,\matD,\matR,\uno\}$
and an \emph{interpolation parameter}  $y_T\in[0,1]$,
and to each line $\ell\in L(\vartheta)$ 
a \emph{derivative label} $\del_\ell\in\{0,1,2\}$, 
a set of \emph{interpolation parameters}\footnote{Recall the definition of $\jap{\cdot}$ in \eqref{jap}.}
$\und{y}_\ell\in[0,1]^{\jap{d(\ell)}}$ and
an \emph{interpolated momentum} $\nu_\ell(\und{y}_\ell)$.
%For any relevant RC $T\in\fT_0(\vartheta)$, let $\fT(T)$ be the set of $T'\in\fT_0(T)$ with $\matO_{T'}\ne\uno$
%such that any relevant RC $T''\subset T$ containing $T'$ has $\matO_{T''}=\uno$.
The new labels are assigned according to the following iterative prescriptions.
\begin{enumerate}[topsep=0ex]
\itemsep0em
\item For any  tree $\vartheta\in\gotT$, let $\fT_0(\vartheta)$ be the set of the relevant RCs $T$ in $\vartheta$ with $d(T)=0$ (if any),
and let $\thetao$ denote the set of vertices and lines in $\vartheta$ which do
not belong to any $T\in\fT_0(\vartheta)$.
To each line $\ell\in L(\thetao)$ we assign $\del_\ell=0$.
\item\label{step3}  To any $T\in\fT_0(\vartheta)$ we assign an operator label $\matO_T\in\{\matL,\matD,\matR\}$,
with $\matO_T=\matL$ if $\calP_T=\emptyset$.
\begin{enumerate}[topsep=0ex]
\itemsep0em
\item[\ref{step3}.1.] \label{elleT}
If 
$\matO_T=\matL$,  let $\To$ denote the set of vertices and lines in $T$ which do not belong to any 
relevant RC  in $T$;
to each $\ell\in L(\To)$ we assign $\del_\ell=0$.
\item[\ref{step3}.2.] \label{dT}
If  $\matO_T=\matD$,  
we assign $\del_\ell=1$ to one and only one line $\ell\in\calP_T$.
Let $\To$ denote the set of vertices and lines in $T$ such that if there is a relevant RC $T'$ in $T$ in which they are contained,
then $T'\in\calmC_\ell(T)$;
%
%which do not belong to any relevant RC 
%in $\fT_0(T)\setminus \calmC_\ell(T)$; 
%
we set $\del_{\ell'}=0$ for each line $\ell'\in L(\To)\setminus\{\ell\}$,
and to each $T'\in \calmC_\ell(T)$ we assign $\matO_{T'}=\uno$.
\item[\ref{step3}.3.] \label{R12T}
If  $\matO_T=\matR$, either we assign 
 $\del_\ell=2$ to one and only one line $\ell\in\calP_T$ or we assign $\del_{\ell_1}=\del_{\ell_2}=1$ to two distinct lines $\ell_1,\ell_2\in\calP_T$.
\begin{enumerate}[topsep=0ex]
 \item[\ref{step3}.3.1.] \label{R1T}
 In the first case, let $\To$ denote the set of vertices and lines in $T$ such that if there is a relevant RC $T'$ in $T$ in which they are contained,
then $T'\in\calmC_\ell(T)$; 
% $\To$ denote the set of vertices and lines in $T$ which do not belong to any
% relevant RC 
%in $\fT_0(T)\setminus \calmC_\ell(T)$; 
we set $\del_{\ell'}=0$ for each line $\ell'\in L(\To)\setminus\{\ell\}$,
and  to each $T'\in \calmC_\ell(T)$
we assign $\matO_{T'}=\uno$. 
\item[\ref{step3}.3.2.] \label{R2T}
In the second case, let  $\To$ denote the set of vertices and lines in $T$ such that if there is a relevant RC $T'$ in $T$ in which they are contained,
then $T'\in\calmC_{\ell_1}(T)\cup\calmC_{\ell_2}(T)$;
%$\To$ denote the set of vertices and lines in $T$ which do not belong to any 
%relevant RC
%in $\fT_0(T)\setminus (\calmC_{\ell_1}(T)\cup\calmC_{\ell_2}(T))$;
 we set $\del_{\ell'}=0$ for each line $\ell'\in L(\To)\setminus\{\ell_1,\ell_2\}$, and 
 to each $T'\in \calmC_{\ell_1}(T)\cup\calmC_{\ell_2}(T)$ we assign $\matO_{T'}=\uno$. 
\end{enumerate}
\end{enumerate}
%Define $\fT_1(\vartheta)$ as the set of relevant RCs $T$ in $\vartheta$ with $d(T)=1$ that do not have an operator label assigned yet.

%
%Let $\fT(T)$ be the set of $T'\in\fT_0(T)$ with $\matO_{T'}\ne\uno$ such that any relevant RC $T''\subset T$ containing $T'$ has $\matO_{T''}=\uno$.
%Note that if $T'\in\fT(T)$,
%then it does not have an operator label assigned yet.
%
\item Next, for $d\ge1$, we define iteratively $\fT_d(\vartheta)$ as the set of relevant RCs $T$ in $\vartheta$ with $d(T)=d$ that do not have an operator 
label assigned yet, and, for each $T'\in\fT_d(\vartheta)$, we assign to $T'$ an operator label ${\matO_{T'}}\in\{\matL,\matD,\matR\}$
and then proceed as in Steps \ref{step3}.1 to \ref{step3}.3. 
%with $\fT_0(\vartheta)$ replaced with $\fT_d(\vartheta)$.

%consider the relevant RCs $T\in\fT_d(\vartheta)$ 
%
%
%The recursion is on the depth of $T'$: we assign to each $T'$ with $d(T')=1$  an operator label ${\matO_{T'}}\in\{\matL,\matD,\matR\}$
%and proceed as in Steps \ref{cinque}.1 to \ref{cinque}.3.
%
\item  For any relevant RC $T$ we set $y_T\in(0,1)$ if $\matO_T=\matR$, $y_T=0$ if 
$\matO_T\in\{\matL,\matD\}$ and $y_T=1$ if $\matO_T=\uno$.
\item For any $\ell\in L(\vartheta)$, if $d(\ell)=0$ we set $\und{y}_\ell=1$ and $\nu_\ell(1):=\nu_\ell$,
%%
%\item 
while, if $d=d(\ell)\ge 1$ and $T_{d-1}\subset T_{d-2}\subset\ldots\subset T_0$
are the relevant RCs such that $\ell\in\calP_{T_i}$ for $i=0,\ldots,d-1$, then, shortening
$y_i=y_{{T_i}}$, $\ell'_i=\ell'_{T_i}$, $j_i=j_{\ell'_i}$ and $\s_i=\s_{\ell'_i}$,
we set $\und{y}_\ell:=(y_{d-1},\ldots ,y_0)$ and  define $\nu_{\ell}(\underline{y}_\ell)$ so that
\begin{equation}\label{nuovonu}
\s_\ell\nu_{\ell}(\underline{y}_\ell):=
\nu_\ell^0(T_{d-1})+ \s_{d-1}\gote_{j_{d-1}} + \s_{d-1}y_{d-1}(\nu_{\ell'_{d-1}}(\und{y}_{\ell'_{d-1}})-\gote_{j_{d-1}}),
\end{equation}
with $\nu_\ell^0(T_{d-1})$ defined as in \eqref{ellezero}.
\item\label{dieci?} For any relevant RC $T$ in $\vartheta$ with $\matO_T \neq \uno$,
let $\fT(T)$ be the set of the relevant RCs $T'\subset T$ with $\matO_{T'}\ne\uno$
such that any relevant RC $T''\subset T$ containing $T'$ has $\matO_{T''}=\uno$.
\end{enumerate}
We call \emph{renormalized tree} any tree in $\fR$. 
Finally, we call  $\fR^{(k)}_{j,\nu,\s}$ the set of trees in $\fR$ with order
$k$ such that the root line has component $j$, momentum $\nu$ and sign $\s$.
\end{defi}
%%%%%%%%%%%%%%%%%%%%%%%%%%%%%%%%%%%%%%%%%%%%%%%%%%%%%%%%%%%%%%%%%%%%%%%%%%

%%%%%%%%%%%%%%%%%%%%%%%%%%%%%%%%%%%%%%%%%%%%%%%%%%%%%%%%%%%%%%%%%%%%%%%%%%
% FIGURA 16
%%%%%%%%%%%%%%%%%%%%%%%%%%%%%%%%%%%%%%%%%%%%%%%%%%%%%%%%%%%%%%%%%%%%%%%%%%
\begin{figure}[ht]
\vspace{.4cm}
\centering
\ins{224pt}{-032pt}{$T''$}
\ins{187pt}{-088pt}{$\ell_3$}
\ins{242pt}{-089pt}{$\ell_4$}
\ins{296pt}{-089pt}{$\ell_5$}
\ins{328pt}{-038pt}{$T'$}
\ins{358pt}{-020pt}{$T$}
\ins{136pt}{-076pt}{$\ell_2$}
\ins{108pt}{-079pt}{$\ell_1$}
\ins{357pt}{-089pt}{$\ell_6$}
%\ins{180pt}{-012pt}{$j,\gote_j,\s$}
%\ins{248pt}{-012pt}{$j,\gote_j,\s$}
\subfigure{\includegraphics*[width=6in]{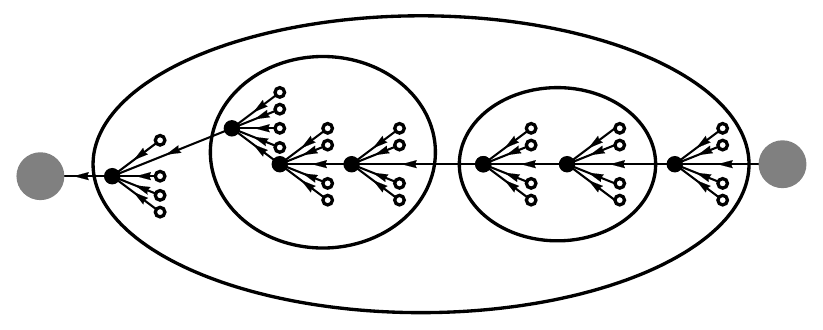}}
\caption{\small A RC $T$ containing two RCs $T'$ and $T''$.}
\label{RGT}
\end{figure}
%%%%%%%%%%%%%%%%%%%%%%%%%%%%%%%%%%%%%%%%%%%%%%%%%%%%%%%%%%%%%%%%%%%%%%%%%%

See Figure \ref{RGT} for an example illustrating the assignment of the operator labels in a resonant tree $\vartheta$ with a relevant
RC $T$ which contains two RCs $T'$ and $T''$ (as usual the grey circles represent the remaining parts of the tree).
Suppose that $\matO_T=\matR$ and that there are two distinct lines $\ell,\ell'\in \calP_T$
such that $\partial_\ell=\partial_{\ell'}=1$. If at least one of such lines is inside $T''$, then $\matO_{T''}=\uno$; if both of them
are inside $T''$ (i.e. $\ell=\ell_2$ and $\ell'=\ell_3$), then, again, $\matO_{T''}=\uno$, while $\matO_{T'}\in\{\matL,\matD,\matR\}$,
and moreover $\partial_{\ell_1}=\partial_{\ell_4}=\partial_{\ell_6}=0$, while $\partial_{\ell_5}=0$ if $\matO_{T'}=\matL$,
$\partial_{\ell_5}=1$ if $\matO_{T'}=\matD$ and $\partial_{\ell_5}=2$ if $\matO_{T'}=\matR$;
if $\ell\in L(T')$ and $\ell'\in L(T'')$, or vice versa,
then $\matO_{T'}=\matO_{T''}=\uno$ and all the other lines $\ell'' \in \calP_T \setminus\{\ell,\ell'\}$ have $\partial_{\ell''}=0$;
if $\ell$ and $\ell'$ belong both to the set $\{\ell_1,\ell_4,\ell_6\}$, then $\matO_{T'},\matO_{T''}\in\{\matL,\matD,\matR\}$.
Since we are supposing $\matO_T=\matR$, the interpolated momentum, say, of the line $\ell_3$ is such that,
using the notation in Remark \ref{nuzero},
\[
\s_{\ell_3} \nu_{\ell_3}(y_{T''},y_{T}) = 
\begin{cases}
\nu_{\ell_3}^0(T) + \s_{\ell_T'} \gote_{j_{\ell_T'}} + \s_{\ell_T'} y_T \, (\nu_{\ell_T'} - \gote_{j_{\ell_T'}}) , & \matO_{T''}=\uno , \\
\nu_{\ell_3}^0(T'') + \s_{\ell_4} \gote_{j_{\ell_4}}  , \phantom{(\nu_{\ell_T'} - \gote_{j_{\ell_T'}})} & \matO_{T''}=\matL, \, \matD , \\
%\nu_{\ell_3}^0(T) + \s_{\ell_4} \gote_{j_{\ell_4}}  , \phantom{(\nu_{\ell_T'} - \gote_{j_{\ell_T'}})} & \matO_{T''}=\matD, \\
\nu_{\ell_3}^0(T'') + \s_{\ell_4} \gote_{j_{\ell_4}} + \s_{\ell_4} y_{T''} \, (\nu_{\ell_4}(y_{T}) - \gote_{j_{\ell_4}}) , & \matO_{T''}=\matR, \\
\end{cases}
\]
with $\s_{\ell_4}\nu_{\ell_4}(y_{T}) \!=\! \nu_{\ell_4}^0(T) + \s_{\ell_{T}'} \gote_{j_{\ell_{T'}}} + \s_{\ell_{T}'} y_T \, (\nu_{\ell_T'} \!-\! \gote_{j_{\ell_T'}})$.
%while $\nu_{\ell_3}^0(T)$, $\nu_{\ell_3}^0(T'')$ and $\nu_{\ell_4}^0(T)$ are defined according to \eqref{ellezero}.
The momenta of the other lines are given by similar -- or simpler -- expressions.
Note that $y_{T''}=0$ if $\matO_{T''}=\matL,\matD$ and $y_{T''}=1$ if $\matO_{T''}=\uno$.

%%%%%%%%%%%%%%%%%%%%%%%%%%%%%%%%%%%%%%%%%%%%%%%%%%%%%%%%%%%%%%%%%%%%%%%%%%
\begin{rmk}\label{marmellatadinuovo}
\emph{%remark
Iterating \eqref{nuovonu} one finds, explicitly,
\[
\begin{aligned}
\s_\ell\nu_{\ell}(\underline{y}_\ell)&=
\nu_\ell^0(T_{d-1}) + \s_{d-1}\gote_{j_{d-1}}+\Big(\prod_{i=0}^{d-1}\s_iy_i\Big)(\nu_{\ell'_0} -\gote_{j_0}) \\
&%y_p\s_p(\nu_{\ell_p}^0+\s_{p-1}\gote_{j_{p-1}}-\gote_{j_p})\\
+\sum_{i=1}^{d-1}\Big(\prod_{h=i}^{d-1}\s_hy_h\Big)
 (\nu^0_{\ell'_i} (T_{i-1})+\s_{i-1} \gote_{j_{i-1}} - \gote_{j_i}) .
\end{aligned}
\]
}%remark
\end{rmk}
%%%%%%%%%%%%%%%%%%%%%%%%%%%%%%%%%%%%%%%%%%%%%%%%%%%%%%%%%%%%%%%%%%%%%%%%%%

%%%%%%%%%%%%%%%%%%%%%%%%%%%%%%%%%%%%%%%%%%%%%%%%%%%%%%%%%%%%%%%%%%%%%%%%%%
\begin{rmk}\label{marmellata}
\emph{%remark
By construction, for any relevant RC $T$ of a  tree $\vartheta\in\fR$ one has
\[
\sum_{\ell\in L(\To)} \del_\ell =
\begin{cases}
0, & \matO_T=\matL, \\
1, & \matO_T=\matD,\\
2, &\matO_T=\matR. 
\end{cases}
\]
In particular, for any chain $\gotC$ of RCs there are at most either two {links} $\ell,\ell'$ of 
$\gotC$ with $\del_\ell=\del_{\ell'}=1$ or one {link} $\ell$ of $\gotC$ with $\del_\ell=2$.
}%remark
\end{rmk}
%%%%%%%%%%%%%%%%%%%%%%%%%%%%%%%%%%%%%%%%%%%%%%%%%%%%%%%%%%%%%%%%%%%%%%%%%%

Define the {\emph{renormalized small divisor} }
\begin{equation}\label{simbolor}
x_{\ell}(\und{y}_\ell):= \om\cdot\nu_\ell(\und{y}_\ell)-\om_{j_\ell}\,,
\end{equation}
and the \emph{renormalized propagator} as
\begin{equation}\label{gr}
\matG^\ttR_\ell(\om) := 
\begin{cases}
\calG_{n_\ell}(x_{\ell}(\und{y}_\ell)),&\qquad \del_\ell=0, \\
\del\calG_{n_\ell} (x_{\ell}(\und{y}_\ell)), & \qquad\del_\ell=1,\\
\del^2\calG_{n_\ell} (x_{\ell}(\und{y}_\ell)), & \qquad\del_\ell=2,
\end{cases}
\end{equation}
where $\del$ denotes the derivative w.r.t. the argument
and $\calG_n(x)$ is defined in \eqref{procesi}.

%%%%%%%%%%%%%%%%%%%%%%%%%%%%%%%%%%%%%%%%%%%%%%%%%%%%%%%%%%%%%%%%%%%%%%%%%%
\begin{defi}[\textbf{Renormalized value of a resonant cluster}]\label{defvalue6-0}
Given a relevant RC $T\subset\vartheta\in\fR$ with $\matO_T\ne\uno$, setting
\begin{equation}\label{termouno}
\Val_0^\ttR(T;c,\om):= 
\Biggl(\prod_{\la\in \Lambda(\To)}\LL_\la(c)\Biggr) 
\Biggl(\prod_{v\in N(\To)}\calF_v^*(c)\Biggr) \Biggl(\prod_{\ell\in L(\To)\setminus \calP_T}\matG_\ell^\ttR(\om)\Biggr)
\end{equation}
and $\PPTo:=L(\To)\cap\calP_T$, the \emph{renormalized value} of $T$ is recursively defined, if $\matO_T=\matL$, as
\begin{equation} \label{termounoderiregola1}
\Val^\ttR(T;c,\om) :=
\Val_0^\ttR(T;c,\om) \Biggl(\prod_{T'\in \fT(T)} \Val^\ttR(T';c,\om) \Biggr) \!
\Biggl(\prod_{\ell\in \PPTo}\matG_\ell^\ttR(\om)\Biggr) ,
\end{equation}
if $\matO_T=\matD$, as
\begin{equation} \label{termounoderiregola2}
\Val^\ttR(T;c,\om) :=
x_{\ell'_T}(\und{y}_{\ell'_T}) \, \Val_0^\ttR(T;c,\om) \Biggl(\prod_{T'\in \fT(T)} \Val^\ttR(T';c,\om) \Biggr) \!
\Biggl(\prod_{\ell\in \PPTo}\matG_\ell^\ttR(\om)\Biggr) ,
\end{equation}
and, if $\matO_T=\matR$, as
\begin{equation} \label{termounoderiregola3}
\null\!\!\!
\Val^\ttR(T;c,\om) :=
(x_{\ell'_T}(\und{y}_{\ell'_T}))^2 \! \Val_0^\ttR(T;c,\om)\!\!\! \int_0^1 \!\! {\rm d}y_{T} (1 \!-\! y_T) 
\Biggl(\prod_{T'\in \fT(T)} \!\!\!\!\!\! \Val^\ttR(T';c,\om) \!\! \Biggr)  \!\!
\Biggl(\prod_{\ell\in  \PPTo} \!\! \matG_{\ell}^\ttR(\om) \!\!\Biggr) .
\end{equation}
\end{defi}
%%%%%%%%%%%%%%%%%%%%%%%%%%%%%%%%%%%%%%%%%%%%%%%%%%%%%%%%%%%%%%%%%%%%%%%%%%

%%%%%%%%%%%%%%%%%%%%%%%%%%%%%%%%%%%%%%%%%%%%%%%%%%%%%%%%%%%%%%%%%%%%%%%%%%
\begin{rmk} \label{wellposedrecursive}
\emph{
The recursive definitions \eqref{termounoderiregola1}--\eqref{termounoderiregola3}
are well-posed since any $T'\in\fT(T)$ has depth $d(T')>d(T)$. The renormalized value $\Val^\ttR(T;c,\om)$
depends on $y_{T''}$ for all $T''$ such that $\calP_{T''}\cap\calP_T\ne\emptyset$; see \eqref{nuovonu}.
 }
 \end{rmk}
%%%%%%%%%%%%%%%%%%%%%%%%%%%%%%%%%%%%%%%%%%%%%%%%%%%%%%%%%%%%%%%%%%%%%%%%%%
 
%%%%%%%%%%%%%%%%%%%%%%%%%%%%%%%%%%%%%%%%%%%%%%%%%%%%%%%%%%%%%%%%%%%%%%%%%%
\begin{defi}[\textbf{Gain factor}] \label{guadagno}
We call \emph{gain factor} both the factor $x_{\ell'_T}$ appearing in $\Val^\ttR(T;c,\om)$ when $\matO_T=\matD$
and the factor $(x_{\ell'_T}(\und{y}_{\ell'_T}))^2$ appearing in $\Val^\ttR(T;c,\om)$ when $\matO_T=\matR$.
\end{defi}
%%%%%%%%%%%%%%%%%%%%%%%%%%%%%%%%%%%%%%%%%%%%%%%%%%%%%%%%%%%%%%%%%%%%%%%%%%

%%%%%%%%%%%%%%%%%%%%%%%%%%%%%%%%%%%%%%%%%%%%%%%%%%%%%%%%%%%%%%%%%%%%%%%%%%
\begin{rmk}\label{nonevero}
\emph{%remark
In \eqref{termounoderiregola1}--\eqref{termounoderiregola3}
the product of the renormalized propagators depends on the label $\matO_T$:
indeed,  all propagators $\matG^\ttR_\ell(\om)$, with $\ell \in L(\To)$,
are given by the first line in \eqref{gr}, except one if $\matO_T=\matD$, which is given by the second line in \eqref{gr},
and one or two if $\matO_T=\matR$, which are given by the third line or both by the second line of \eqref{gr}, respectively.
}%remark
\end{rmk}
%%%%%%%%%%%%%%%%%%%%%%%%%%%%%%%%%%%%%%%%%%%%%%%%%%%%%%%%%%%%%%%%%%%%%%%%%%

%%%%%%%%%%%%%%%%%%%%%%%%%%%%%%%%%%%%%%%%%%%%%%%%%%%%%%%%%%%%%%%%%%%%%%%%%%
\begin{defi}[\textbf{Renormalized value of an renormalized tree}]\label{defvalue6-1}
Given a renormalized tree $\vartheta\in  \fR$, we define the \emph{renormalized value} of $\vartheta$ as
\begin{equation}\label{valR}
\null\hspace{-.2cm}
\Val^\ttR(\vartheta;c,\om) :=
\displaystyle{\Bigg(\!\prod_{\substack{T\subset \vartheta \\ d(T)=0}} \!\!\!  \Val^\ttR(T;c,\om) \!\!\Bigg) 
\!\! \Bigg(\prod_{\la\in \Lambda(\thetao)} \!\!\! \LL_\la (c)\Bigg) 
\!\! \Bigg(\prod_{v\in N(\thetao) } \!\!\! \calF_v^*(c) \Bigg)
\!\! \Bigg(\prod_{\ell\in L(\thetao)} \!\!\! \matG_\ell^\ttR(\om)\Bigg) } . 
\end{equation}
where $\LL_\la(c)$, $\calF_v^*(c)$ and $\matG_\ell^\ttR(\om)$ are defined 
in \eqref{foglie}, \eqref{nodiker} and \eqref{gr}, respectively, while $\Val^\ttR(T;c,\om)$
is recursively defined in \eqref{termounoderiregola1}--\eqref{termounoderiregola3},
\end{defi}
%%%%%%%%%%%%%%%%%%%%%%%%%%%%%%%%%%%%%%%%%%%%%%%%%%%%%%%%%%%%%%%%%%%%%%%%%%

%
%Similarly, for any tree $\vartheta \in \breve\fR^{(k)}_{j,\nu,\s}$ we define the \emph{renormalized value} of  $\vartheta$ as
%\begin{equation}\label{cio}
%\null\hspace{-.2cm}
%\Val^\ttR(\vartheta;c,\om):= \Bigg(\!\prod_{\substack{T\subset \vartheta \\ d(T)=0}} \Val^\ttR(T;c,\om) \!\!\Bigg) 
%\!\! \Bigg(\prod_{\la\in \Lambda(\thetao)}\LL_\la (c)\Bigg) 
%\!\! \Bigg(\prod_{\ell\in L(\thetao)} \matG_\ell^\ttR(\om)\Bigg).
%\end{equation}

%%%%%%%%%%%%%%%%%%%%%%%%%%%%%%%%%%%%%%%%%%%%%%%%%%%%%%%%%%%%%%%%%%%%%%%%%%
\begin{rmk}\label{lescaleR}
\emph{%remark
If a tree $\vartheta\in\fR$ is such that $\Val^\ttR(\vartheta;c,\om)\ne0$, then
for all $\ell\in L(\vartheta)$ with $n_\ell\ge 0$ one has, by definition of scale,
\begin{equation}\label{tagliaR}
%\frac{1}{16}\be_{\om}(r_{m_{n_\ell}}) \le |\om\cdot\nu_\ell - \om_{j_\ell}| < \frac{1}{8}\be_{\om}(r_{m_{n_\ell-1}}),
\frac{1}{16}\be({m_{n_\ell}}) < |x_\ell(\und{y}_\ell) | < \frac{1}{8}\be({m_{n_\ell-1}}) ,
\end{equation}
where, {as in Remark \ref{lescale}}, $\be(m_{-1})$ has to be interpreted as $+\io$.
}%remark
\end{rmk}
%%%%%%%%%%%%%%%%%%%%%%%%%%%%%%%%%%%%%%%%%%%%%%%%%%%%%%%%%%%%%%%%%%%%%%%%%%

Now, we want to extend the definitions above to the resonant trees.

%%%%%%%%%%%%%%%%%%%%%%%%%%%%%%%%%%%%%%%%%%%%%%%%%%%%%%%%%%%%%%%%%%%%%%%%%%
\begin{defi}[\textbf{Set $\boldsymbol{\fS}$ of the renormalized resonant trees}] %\hspace{-.27cm}\boldsymbol{\fS}$}]
\label{ultima?}
Let $\fS$ be the set of RTs
which can be obtained from any RT $\TT\in\gotS$ as follows.
\begin{enumerate}[topsep=0ex]
\itemsep0em
\item 
%Let $\fT_0(\TT)$ be the set of relevant RCs contained in $\TT$.
Define the depth of a relevant RC $T $ in $\TT$ as in Definition \ref{profondita},
and the depth and the cloud of a line $\ell\in L(\TT)$ as in Definitions \ref{profondita2} and \ref{nuvola}, respectively,
with $\vartheta$ replaced with $\TT$. 
\item\label{cinque} Assign to $\TT$ two labels, the \emph{operator label} $\matO_\TT\in\{\matL,\matD,\matR\}$,
with the constraint that $\matO_\TT=\matL$ if $\calP_\TT = \emptyset$,
and the \emph{interpolation label} $y_\TT\in[0,1]$, with the constraint that $y_\TT=0$ if $\matO_\TT\in\{\matL,\matD\}$
and $y_\TT\in(0,1)$ if $\matO_\TT=\matR$.
%\item Assign to each  $T\in\fT_0(\TT)$ an \emph{operator label} $\matO_T\in\{\matL,\matD,\matR,\uno\}$
%and an \emph{interpolation parameter}  $y_T\in[0,1]$,
%and to each line $\ell\in L(\TT)$  a \emph{derivative label} $\del_\ell\in\{0,1,2\}$, according to the following iterative procedure.
%
%\item \label{cinque} First we consider the RT $\TT$.
%
\begin{enumerate}[topsep=0ex]
\itemsep0em
\item[\ref{cinque}.1.] \label{elleTbis}
If $\matO_\TT=\matL$,  we define $\TTo$ as the set of vertices and lines in $\TT$ which do not belong to any 
relevant RC  in $\TT$, and to each $\ell\in L(\TTo)$ we assign a label $\del_\ell=0$.
\item[\ref{cinque}.2.] \label{dTbis}
If  $\matO_\TT=\matD$,  we assign a label $\del_\ell=1$ to one and only one line $\ell\in\calP_\TT$, we define
$\TTo$ as the set of vertices and lines in $\TT$ such that if there is a relevant RC $T$ in $\TT$ in which they are
contained, then $T\in\calmC_\ell(\TT)$; we set $\del_{\ell'}=0$ for each line $\ell'\in L(\TTo)\setminus\{\ell\}$,
and to each $T\in \calmC_\ell(\TT)$ we assign label $\matO_{T}=\uno$.
\item[\ref{cinque}.3.] \label{R12Tbis}
If  $\matO_\TT=\matR$, either we assign 
 $\del_\ell=2$ to one and only one line $\ell\in\calP_\TT$ or we assign $\del_{\ell_1}=\del_{\ell_2}=1$ to two distinct lines $\ell_1,\ell_2\in\calP_\TT$.
\begin{enumerate}[topsep=0ex]
 \item[\ref{cinque}.3.1.] \label{R1Tbis}
 In the first case, we define $\TTo$ as the set of vertices and lines in $\TT$ such that if there is a relevant RC $T$ in $\TT$ in which they are
contained, then $T\in\calmC_\ell(\TT)$; we set $\del_{\ell'}=0$ for each line $\ell'\in L(\TTo)\setminus\{\ell\}$,
and to each $T\in \calmC_\ell(\TT)$ we assign $\matO_{T}=\uno$. 
\item[\ref{cinque}.3.2.] \label{R2Tbis}
In the second case, we define $\TTo$ as the set of vertices and lines in $\TT$ such that if there is a relevant RC $T$ in $\TT$ in which they are
contained, then $T\in\calmC_{\ell_1}(\TT)\cup\calmC_{\ell_2}(\TT)$;
we set $\del_{\ell'}=0$ for each line $\ell'\in L(\TTo)\setminus\{\ell_1,\ell_2\}$, and 
to each $T\in \calmC_{\ell_1}(\TT)\cup\calmC_{\ell_2}(\TT)$ we assign a label $\matO_{T}=\uno$. 
\end{enumerate}
\end{enumerate}
\item\label{sei} Next for $d\ge0$ we define iteratively $\fT_d(\TT)$ as the set of the relevant RCs $T$ in $\TT$ with $d(T)=d$ (if any)
that do not have an operator label assigned yet.
% consider recursively the relevant RCs $T$ in $\TT$ with $\matO_{T'}$ not assigned yet,
%the recursion being on the depth $d(T')$ of $T'$.
\begin{itemize}[topsep=0ex]
\itemsep0em
\item[\ref{sei}.1]
We assign to any $T$ an operator label ${\matO_{T}}\in\{\matL,\matD,\matR\}$,
and proceed as in Steps \ref{step3}.1 to \ref{step3}.3 of Definition \ref{nonepossopiu}.
\item[\ref{sei}.2]
For any $T$ we assign to the RT $\TT_{T}$ associated with $T$
an operator label ${\matO_{\TT_{T}}}\in\{\matL,\matD,\matR\}$ by setting $\matO_{\TT_{T}}:=\matO_{T}$.
\end{itemize}
\item We define $\fT(\TT)$ as the set of relevant RCs $T$ with $\matO_T\ne \uno$
such that any relevant RC $T'$ containing $T$ has $\matO_{T'}=\uno$.
\item  For any relevant RC $T$ in $\TT$ we set $y_T\in(0,1)$ if $\matO_T=\matR$, $y_T=0$ if 
$\matO_T\in\{\matL,\matD\}$ and $y_T=1$ if $\matO_T=\uno$, and define $y_{\TT_T}=y_T$.
%
%\item If $d(\ell)=0$ we set 
%\[
%\und{y}_\ell:=\begin{cases}
%0 & \ell\notin \calP_\TT,\\
%y_\TT & \ell\in \calP_\TT,
%\end{cases}
%\]
%%
%and define the  \emph{renormalized symbol function} $\x^\ttR_\ell\!:\RRR^{2}\to\RRR$ given by
%\begin{equation}\label{lexir0}
%\x_{\ell}^\ttR(\und{y}_\ell,x) =\begin{cases}
%x_\ell , & \ell\notin \calP_\TT, \\
% \s_{\ell}(\om\cdot\nu^0_{\ell}(\TT) + \s_{\ell^{in}_{\TT}}\om_{j_{\ell^{in}_{\TT}}}+ \s_{\ell^{in}_{\TT}} y_\TT x ) -\om_{j_{\ell}}    , & \ell\in \calP_\TT ,
%\end{cases}
%\end{equation}
%%
%\item 
%If $d=d(\ell)\ge 1$ and $T_{d-1}\subset T_{d-2}\subset\ldots\subset T_0\subset\TT$
%are the relevant RCs such that $\ell\in\calP_{T_i}$ for $i=0,\ldots,d-1$, then, shortening
%$y_i=y_{{T_i}}$, $\ell'_i=\ell'_{T_i}$, $j_i=j_{\ell'_i}$ and $\s_i=\s_{\ell'_i}$,
%we set 
%\[
%\und{y}_\ell:=\begin{cases}
%(y_{d-1},\ldots ,y_0,0) & \ell\notin \calP_\TT,\\
%(y_{d-1},\ldots ,y_0,y_\TT) & \ell\in \calP_\TT.
%\end{cases}
%\]
%and define the  \emph{renormalized symbol function} $\x^\ttR_\ell\!:\RRR^{d+2}\to\RRR$ given by
%\begin{equation}\label{lexird}
%\x_{\ell}^\ttR(\und{y}_\ell,x) =
% \s_{\ell}(\om\cdot\nu^0_{\ell}(T_{d-1}) + \s_{\ell'_{d-1}}\om_{j_{\ell'_{d-1}}}+ \s_{\ell'_{d-1}} y_{d-1} \x_{\ell'_{d-1}}^\ttR(\und{y}_{\ell'_{d-1}},x) ) -\om_{j_{\ell}}     ,
%\end{equation}

%
\end{enumerate}
We call \emph{renormalized resonant trees} the trees $\TT\in\fS$, and, for $\matO=\matL,\matD,\matR$, we  indicate with
$\fS(\matO)$ the sets of RTs $\TT\in\fS$ with $\matO_\TT=\matO$.
\end{defi} 
%%%%%%%%%%%%%%%%%%%%%%%%%%%%%%%%%%%%%%%%%%%%%%%%%%%%%%%%%%%%%%%%%%%%%%%%%%

%%%%%%%%%%%%%%%%%%%%%%%%%%%%%%%%%%%%%%%%%%%%%%%%%%%%%%%%%%%%%%%%%%%%%%%%%%
\begin{rmk}
\emph{
For any $\TT\in\fS$, the set $\TTo$ is in fact the set
of vertices and lines in $\TT$ which do not belong to any relevant RC  in $\fT(\TT)$.
}
\end{rmk}
%%%%%%%%%%%%%%%%%%%%%%%%%%%%%%%%%%%%%%%%%%%%%%%%%%%%%%%%%%%%%%%%%%%%%%%%%%

%%%%%%%%%%%%%%%%%%%%%%%%%%%%%%%%%%%%%%%%%%%%%%%%%%%%%%%%%%%%%%%%%%%%%%%%%%
\begin{defi}[\textbf{Renormalized value of a renormalized resonant tree}]\label{defvalue7}
For$\,$any$\,$renor\-malized RT $\TT\in\fS$, set
\begin{equation} \label{argh}
\Val_0^\ttR(\TT;c,\om):= 
\Biggl(\prod_{\la\in \Lambda(\TTo)}\LL_\la(c)\Biggr) 
\Biggl(\prod_{v\in N(\TTo)}\calF_v^*(c)\Biggr) \Biggl(\prod_{\ell\in L(\TTo)\setminus \calP_\TT}
\matG_\ell^\ttR(\om)\Biggr) ,
\end{equation}
where  $\PPTTo:=L(\TTo)\cap\calP_\TT$.
Then %for any $k\ge 1$, $j\in\ZZZ$, $\nu\in\ZZZ^{\ZZZ}_f$ and $\s\in\{\pm\}$,
we define recursively the \emph{renormalized value} of $\TT$, %\in  \fS^{(k)}_{j,\nu,\s}$,
if $\matO_\TT=\matL$, as
\begin{equation} \nonumber \label{termounoderiregola4}
\Val^\ttR(\TT;x,c,\om)=\Val^\ttR(\TT;0,c,\om) :=
\Val_0^\ttR(\TT;c,\om) \Biggl( \prod_{T'\in \fT(\TT)} \Val^\ttR(\TT_{T'};\x_{\ell'_{T'}}(0), c, \om) \Biggr) \!
\Biggl(\prod_{\ell\in \PPTTo} \matG_\ell^\ttR(0,\om)\Biggr) ,
\end{equation}
if $\matO_\TT=\matD$, as
\begin{equation} \nonumber \label{termounoderiregola5}
\Val^\ttR(\TT;x,c,\om) :=
x \Val_0^\ttR(\TT;c,\om) \,
\Biggl( \prod_{T'\in \fT(\TT)} \Val^\ttR(\TT_{T'}; \x_{\ell'_{T'}}(0), c,\om) \Biggr)\!
\Biggl( \prod_{\ell\in \PPTTo}\matG_\ell^\ttR(0,\om) \Biggr) ,
\end{equation}
and if $\matO_\TT=\matR$, as
\begin{equation} \nonumber \label{termounoderiregola6}
\null\!\!\!\!\!\!\!\!\!
\Val^\ttR(\TT \! ; x, \! c, \! \om) \!:=\!
x^2 \!\Val_0^\ttR(\TT \! ; \! c , \!\om) \!
\int_0^1 \!\!
{\rm d}y_{\TT} (1 \!-\! y_\TT)
 \Biggl(\prod_{T'\in \fT(\TT)} \!\!\!\!\!\! \Val^\ttR(\TT_{T'}; \x_{\ell'_{T'}}(y_\TT x) ,  c,\om ) \!\! \Biggr)
\!\!\Biggl(\prod_{\ell\in \PPTTo} \!\!\! \matG_{\ell}^\ttR(y_\TT x,\om) \!\! \Biggr) \! ,
\end{equation}
%
%and, if $\matO_\TT=\uno_i$, with $i=1,2,3$, as
%%
%\begin{equation} \nonumber \label{termounoderiregola7}
%\Val^\ttR(\TT;x,c,\om) :=
%\Val_0^\ttR(\TT;c,\om) 
%\Biggl(\prod_{T'\in \fT(\TT)} \!\!\!\! \Val^\ttR(\TT_{T'};\x_{\ell'_{T'}}(x) , c,\om) \!\Biggr) \!
%\Biggl(\prod_{\ell\in \PPTTo} \!\! \matG_{\ell}^\ttR(x,\om) \! \Biggr) ,
%\end{equation}
%
where the {symbol function} $\x_\ell(x)$ is defined in \eqref{lexi}
%\begin{equation}\label{lexir}
%\x_{\ell}(x) =\begin{cases}
%x_\ell , & \ell'_{T'}\notin \calP_\TT, \\
% \s_{\ell}(\om\cdot\nu^0_{\ell} + \s_{\ell^{in}_{\TT}}\om_{j_{\ell^{in}_{\TT}}}+ \s_{\ell^{in}_{\TT}} x ) -\om_{j_{\ell}}    , & \ell\in \calP_\TT ,
%\end{cases}
%\end{equation}
and the propagator $\matG_{\ell}^\ttR(x,\om)$ is given by
\begin{equation}\label{vocette}
\matG^\ttR_\ell(x,\om) := 
\begin{cases}
\calG_{n_\ell}(\x_\ell(x) ),&\qquad \del_\ell=0, \\
 \del\calG_{n_\ell} (\x_\ell(x) ), & \qquad\del_\ell=1,\\
  \del^2\calG_{n_\ell} (\x_\ell(x) ), & \qquad\del_\ell=2 ,
\end{cases}
\end{equation}
with $\calG_n(x)$ defined in \eqref{procesi}.
\end{defi}
%%%%%%%%%%%%%%%%%%%%%%%%%%%%%%%%%%%%%%%%%%%%%%%%%%%%%%%%%%%%%%%%%%%%%%%%%%

 %%%%%%%%%%%%%%%%%%%%%%%%%%%%%%%%%%%%%%%%%%%%%%%%%%%%%%%%%%%%%%%%%%%%%%%%%%
\begin{defi}[\textbf{Trimmed tree associated with a renormalized tree}] 
\label{ritrimmoful}
{
%For all $k$, $\nu$, $j$, $\s$
For all $k\ge1$, $j\in\ZZZ$,  $\nu\in\ZZZ^{\ZZZ}_{f}$ and $\s\in\{\pm\}$ and for any $\vartheta\in\fR^{(k)}_{ j ,\nu,\s}$, set
\[
N_2(\vartheta) :=\{ v \in N(\vartheta) :s_{v}=2\}.% \hbox{ and } v \st{\preceq} r_{\vartheta}\}  . 
\]
For any $v\in N_2(\vartheta)$, if $\vartheta_v$ denotes the subtree of $\vartheta$ entering $v$
whose root line is an $\h$-line, set 
\[
F_2^*(\vartheta) := \{ \vartheta_v : v \in N_2(\vartheta) \hbox{ and } v \st{\preceq} r_{\vartheta}\}.
\]
Let $\breve\vartheta$ be the trimmed tree 
obtained from $\vartheta$ by cutting the subtrees in $F_2^*(\vartheta)$;
we call $\breve\vartheta$ the \emph{trimmed renormalized tree associated with the renormalized tree} $\vartheta$.}
\end{defi}
%%%%%%%%%%%%%%%%%%%%%%%%%%%%%%%%%%%%%%%%%%%%%%%%%%%%%%%%%%%%%%%%%%%%%%%%%%

%%%%%%%%%%%%%%%%%%%%%%%%%%%%%%%%%%%%%%%%%%%%%%%%%%%%%%%%%%%%%%%%%%%%%%%%%%
\begin{defi}[\textbf{Trimmed subgraph of a renormalized tree}] 
\label{trimmedsubgraph2}
{Given a subgraph $T$ of a renormalized tree $\vartheta$, we define $\breve T$ as the subgraph of $\breve\vartheta$
such that $V(\breve T) = V(T)\cap V(\breve\vartheta)$ and $L(\breve T)= L(T )\cap L(\breve\vartheta)$.
We call $\breve T$ the \emph{trimmed subgraph associated with the subgraph} $T$.}
\end{defi}
%%%%%%%%%%%%%%%%%%%%%%%%%%%%%%%%%%%%%%%%%%%%%%%%%%%%%%%%%%%%%%%%%%%%%%%%%%

%%%%%%%%%%%%%%%%%%%%%%%%%%%%%%%%%%%%%%%%%%%%%%%%%%%%%%%%%%%%%%%%%%%%%%%%%%
\begin{defi} [\textbf{Truncated renormalized tree}] \label{Rn}
A renormalized tree $\vartheta\in \fR$ is said to be \emph{truncated} at scale $n$ if
\begin{enumerate}[topsep=0ex]
\itemsep0em
\item $n_\ell< n$ for all $\ell\in L(\vartheta)$,
\item $J(\breve\vartheta)<C_1r_{m_n-1}$,
\item $J(\breve\vartheta_v)<C_1r_{m_n-1}$ for all $v\in N_2(\vartheta)$.
\end{enumerate}
We call $\fR^{(k)}_{j,\nu,\s}(n)$ the set of trees $\vartheta\in \fR^{(k)}_{j,\nu,\s}$
which are truncated at scale $n$.
\end{defi}
%%%%%%%%%%%%%%%%%%%%%%%%%%%%%%%%%%%%%%%%%%%%%%%%%%%%%%%%%%%%%%%%%%%%%%%%%%

%%%%%%%%%%%%%%%%%%%%%%%%%%%%%%%%%%%%%%%%%%%%%%%%%%%%%%%%%%%%%%%%%%%%%%%%%%
\begin{defi} [\textbf{Truncated renormalized resonant tree}] \label{Sn}
A renormalized resonant tree $\TT\in \fS$ is said to be \emph{truncated} at scale $n$ if
\begin{enumerate}[topsep=0ex]
\itemsep0em
\item $n_\ell< n$ for all $\ell\in L(\TT)$,
\item $J(\TT)<C_1r_{m_n-1}$,
\item $J(\breve\vartheta_v)<C_1r_{m_n-1}$ for all $v\in N_2(\TT)$.
\end{enumerate}
For any $k\ge 1$, $j,j'\in\ZZZ$, $\s,\s'\in\{\pm\}$, $n \ge -1$ and $\matO\in\{\matL,\matD,\matR\}$,
call $\fS^{(k)}_{j,\s,j',\s'}(n,\matO)$ the set of renormalized RTs  $\TT\in\fS(\matO)$ of order $k$
which are truncated at scale $n$ and are such that $j_{\ell^{out}_\TT}=j$,
$\s_{\ell^{out}_\TT}=\s$, $j_{\ell^{in}_\TT}=j'$, $\s_{\ell^{in}_\TT}=\s'$.
\end{defi}
%%%%%%%%%%%%%%%%%%%%%%%%%%%%%%%%%%%%%%%%%%%%%%%%%%%%%%%%%%%%%%%%%%%%%%%%%%

%%%%%%%%%%%%%%%%%%%%%%%%%%%%%%%%%%%%%%%%%%%%%%%%%%%%%%%%%%%%%%%%%%%%%%%%%%
\begin{lemma}\label{proc}
For all $k\ge 1$, $j\in\ZZZ$, $\nu\in\ZZZ^{\ZZZ}_f$, $\s\in\{\pm\}$ and $n\ge -1$, one has
 \[
 \sum_{\vartheta\in\gotT^{(k)}_{j,\nu,\s}(n)}\Val(\vartheta;c,\om) = \sum_{\vartheta\in\fR^{(k)}_{j,\nu,\s}(n)}\Val^\ttR(\vartheta;c,\om) .
%=  \sum_{\vartheta\in\gotF^{(k)}_{j,\nu,\s}(n)}\Val^\ttR(\vartheta;c)  .
 \]
\end{lemma}
%%%%%%%%%%%%%%%%%%%%%%%%%%%%%%%%%%%%%%%%%%%%%%%%%%%%%%%%%%%%%%%%%%%%%%%%%%

%%%%%%%%%%%%%%%%%%%%%%%%%%%%%%%%%%%%%%%%%%%%%%%%%%%%%%%%%%%%%%%%%%%%%%%%%%
\prova
By construction, the whole set $\fR^{(k)}_{j,\nu,\s}(n)$ is obtained by taking all trees $\vartheta\in\gotT^{(k)}_{j,\nu,\s}(n)$
and, for each of them, by assigning further labels so that the following operations are performed on the value
of each relevant RC $T$, iteratively of the depth.
\begin{enumerate}[topsep=0ex]
\itemsep0em
\item First, for any RC $T$ with depth $d(T)=0$,
one writes $\Val(T;c,\om) = \Val(\TT_T;x_{\ell'_T},c,\om)$ according to Remark \ref{pedante};
then, by considering $\Val(\TT_T;x_{\ell'_T},c,\om)$ as a function of $x_{\ell'_T}$, 
one applies, according to the value of the label $\matO_T$, 
an operator $\matO\in\{\matL,\matD , \matR, \uno\}$, as defined in Remark \ref{oraeunremark},
to the value $ \Val(\TT_T;x_{\ell'_T},c,\om)$.
\item Then, one considers any RC $T$ with depth $d(T)=1$, and reasons as before, with the only difference that, because of
the operations performed in the previous step, the value $\Val(T;c,\om)$ has been replaced by
$\Val^\ttR(T;c,\om)=\Val^\ttR(\TT_T;x_{\ell'_T}(\und{y}_\ell),c,\om)$.
\item One iterates the construction until no RC are left.
\end{enumerate}
Finally, one notes that $\matL+\matD+\matR=\uno$.
\EP
%%%%%%%%%%%%%%%%%%%%%%%%%%%%%%%%%%%%%%%%%%%%%%%%%%%%%%%%%%%%%%%%%%%%%%%%%%

Define, for $\matO\in\{\matL,\matD,\matR\}$,
 \begin{subequations}\label{emmetaO}
 \begin{align}
%\matM_{j\s j'\s'}^{(k)}(x,c,\om)  := \!\!\!\sum_{\TT\in \gotS^{(k)}_{j,j',\s,\s'}} \!\!\! \Val(\TT;c,x) , \qquad
\matM_{j\s j'\s'}^{(k)}(x,c,\om,n,\matO) &:= \!\!\! \sum_{\TT\in \fS^{(k)}_{j,\s,j',\s'}(n,\matO)} \!\!\!\!\!\! \Val^\ttR(\TT;x,c,\om),
\label{emmetaOa} \\
\matM_{j j'}^{(k)}(x,c,\om,n,\matO) &:=
\begin{pmatrix}
\matM_{j+ j'+}^{(k)}(x,c,\om,n,\matO)  & \matM_{j+ j'-}^{(k)}( -x,c,\om,n,\matO)  \cr
\matM_{j- j' +}^{(k)}(x,c,\om,n,\matO)  & \matM_{j- j'-}^{(k)}( -x,c,\om,n,\matO) 
\end{pmatrix} ,
\label{emmetaOb}
\end{align}
\end{subequations}
where the set $ \fS^{(k)}_{j,\s,j',\s'}(n,\matO)$ is introduced in Definition \ref{ultima?}.

%%%%%%%%%%%%%%%%%%%%%%%%%%%%%%%%%%%%%%%%%%%%%%%%%%%%%%%%%%%%%%%%%%%%%%%%%%
\begin{lemma}\label{dadim}
For all $k\ge 1$, $j,j'\in\ZZZ$, $\s\in\{\pm\}$ and $n\ge -1$, one has
\[
 \matM_{j\s j'\s'}^{(k)}(x,c,\om,n,\matO)=\matO \matM_{j\s j'\s'}^{(k)}(x,c,\om,n) , \qquad\qquad \matO\in\{\matL,\matD,\matR\}.
 \]
 \end{lemma}
 %%%%%%%%%%%%%%%%%%%%%%%%%%%%%%%%%%%%%%%%%%%%%%%%%%%%%%%%%%%%%%%%%%%%%%%%%%
 
 %%%%%%%%%%%%%%%%%%%%%%%%%%%%%%%%%%%%%%%%%%%%%%%%%%%%%%%%%%%%%%%%%%%%%%%%%%
 \prova
One reasons as in the proof of Lemma \ref{proc}.
 \EP
 %%%%%%%%%%%%%%%%%%%%%%%%%%%%%%%%%%%%%%%%%%%%%%%%%%%%%%%%%%%%%%%%%%%%%%%%%%
 
% %%%%%%%%%%%%%%%%%%%%%%%%%%%%%%%%%%%%%%%%%%%%%%%%%%%%%%%%%%%%%%%%%%%%%%%%%%
% \begin{rmk}\label{stomaco}
% \emph{%remark
% From Lemma \ref{dadim}, the Taylor formula \eqref{taylor} and the definition of the linear operators $\matL,\matD,\matR$,
% one has
%  \[
% \matM_{j\s j'\s'}^{(k)}(x,c,\om,n,\matL) +  \matM_{j\s j'\s'}^{(k)}(x,c,\om,n,\matD) +
%  \matM_{j\s j'\s'}^{(k)}(x,c,\om,n,\matR)= \matM_{j\s j'\s'}^{(k)}(x,c,\om,n).
% \]
%}%remark
%\end{rmk}
%%%%%%%%%%%%%%%%%%%%%%%%%%%%%%%%%%%%%%%%%%%%%%%%%%%%%%%%%%%%%%%%%%%%%%%%%%%

%%%%%%%%%%%%%%%%%%%%%%%%%%%%%%%%%%%%%%%%%%%%%%%%%%%%%%%%%%%%%%%%%%%%%%%%%%
\begin{lemma}\label{stanca}
Assume $\vartheta\in\fR$ to be such that $\Val^\ttR(\vartheta;c,\om)\ne0$. 
Then for any tree $\vartheta\in\fR$ and any line $\ell\in L(\vartheta)$, one has
\begin{equation} \label{iterative}
|x_\ell - x_\ell(\und{y}_\ell)| < \frac{1}{2}  |x_\ell (\und{y}_\ell) | ,
\end{equation}
where $x_\ell $ and $x_\ell(\und{y}_\ell)$ are defined in \eqref{simbolo} and \eqref{simbolor}, respectively.
For $n_{\ell}=0$ the bound \eqref{iterative} may be improved into
\begin{equation} \label{iterative2}
|x_\ell - x_\ell(\und{y}_\ell)| < \frac{3}{16}  \be(m_0) .
\end{equation}
\end{lemma}
%%%%%%%%%%%%%%%%%%%%%%%%%%%%%%%%%%%%%%%%%%%%%%%%%%%%%%%%%%%%%%%%%%%%%%%%%%

%%%%%%%%%%%%%%%%%%%%%%%%%%%%%%%%%%%%%%%%%%%%%%%%%%%%%%%%%%%%%%%%%%%%%%%%%%
\prova
We prove \eqref{iterative} by induction on the depth of the line.
If $d(\ell)=0$, the bound is trivially satisfied since
$\n_\ell(\und{y}_\ell)=\n_\ell(1)=\nu_\ell$ and hence $x_\ell = x_\ell(\und{y}_\ell)$.

Then, consider a line $\ell$ with depth $d(\ell)=p\ge 1$. If $T$ is
the relevant RC with highest depth of the cloud $\matC_\ell(\vartheta)$, 
then $\und{y}_\ell=(y_T,\und{y}_{\ell'_T})$, so that, setting 
\begin{equation}
\begin{aligned}
x^0_\ell & :=\om\cdot\nu_\ell^0+\s_{\ell'_{T}}\om_{j_{\ell'_{T}}} \!\! -\s_\ell\om_{j_\ell} , \\
x & := \om\cdot\nu_{\ell'_{T}} \!\! -\om_{j_{\ell'_{T}}} ,
\end{aligned}
\end{equation}
where $\nu_\ell^0=\nu_\ell^0(T)$ is defined accordingly to \eqref{ellezero}, one has
\[
\s_\ell x_\ell =x_\ell^0+ \s_{\ell_T'} x , \qquad \s_\ell x_\ell(\und{y}_\ell) = x^0_\ell+y_T \s_{\ell_T'} x(\und{y}_{\ell'_{T}}) .
\]

Since $\Val^\ttR(\vartheta;c,\om)\ne0$, Remark \ref{lescaleR} gives % and the inductive hypothesis give, respectively,
\begin{equation} \label{scales}
%\begin{aligned}
%&
\frac{1}{16}\be(m_n) < |x_\ell(\und{y}_{\ell}) |<  \frac{1}{8}\be(m_{n-1}) , \qquad
\frac{1}{16}\be(m_{n'}) < |x (\und{y}_{\ell'_{T}}) |<  \frac{1}{8}\be(m_{n'-1}) , %\\
%& |x - x(\und{y}_{\ell_T'})| < % \frac{1}{2} |x| , \qquad \frac{1}{32}\be(m_{n'}) <  | x| <  \frac{1}{4}\be(m_{n'-1})  ,
%\end{aligned}
\end{equation}
with $n':=n_{\ell'_T}>n_\ell=:n$, by definition of RC; moreover the inductive hypothesis implies
\[
|x-x (\und{y}_{\ell'_{T}})| \le \frac{1}{2}|x (\und{y}_{\ell'_{T}})| .
\]

The SEG $T'$ with $\ell$ as exiting line and $\ell'_{T}$ as entering line
satisfies the hypotheses of Lemma \ref{ovvio2} and
hence, using also Remark \ref{utile} and Definition \ref{coso}, one has
$$ C_1 r_{m_n-1} \le C_1  |\nu^0_\ell+\s_{\ell'_T}\gote_{j_{\ell'_T}}-\s_\ell\gote_{j_\ell} |_{\al/2}\le J(T') \le J(T) < C_1 r_{m_{n'}-1} , $$
so that $|x_\ell^0|\ge\be(m_{n'}-1)\ge \be(m_{n'-1})/2$, which yields
\begin{equation} \label{lower}
%\begin{aligned}
|x_\ell (\und{y}_\ell) | \ge |x_\ell^0 | - |y_T  x (\und{y}_{\ell'_{T}}) |  \ge \frac{1}{2} \be(m_{n'-1}) - \frac{1}{8} \be(m_{n'-1}) \ge
\frac{3}{8} \be(m_{n'-1}) .
%\end{aligned}
\end{equation}
On the other hand, by the inductive hypothesis, one has
\begin{equation} \label{uppero}
\begin{aligned}
|x_\ell - x_\ell (\und{y}_\ell) |
 & \le | x - y_T x(\und{y}_{\ell'_{T}}) | \le | x - x(\und{y}_{\ell'_{T}}) | +  (1-y_T) \, | x(\und{y}_{\ell'_{T}}) | \\
& \le \frac{1}{2} | x(\und{y}_{\ell'_{T}}) | + | x(\und{y}_{\ell'_{T}}) | 
\le \frac{3}{2} | x(\und{y}_{\ell'_{T}}) | \le \frac{3}{16} \be(m_{n'-1}) ,
\end{aligned}
\end{equation}
so that one obtains
\begin{equation} \nonumber
|x_\ell - x_\ell (\und{y}_\ell) | < \frac{1}{2} |x_\ell (\und{y}_\ell) | ,
\end{equation}
which shows that the line $\ell$, too, satisfies the bound \eqref{iterative}. 
This completes the inductive proof of \eqref{iterative}. 
To obtain \eqref{iterative2} simply note that if $n_{\ell} =0$ then $m_{n'-1} \ge m_0$, and use \eqref{uppero}.
\EP
%%%%%%%%%%%%%%%%%%%%%%%%%%%%%%%%%%%%%%%%%%%%%%%%%%%%%%%%%%%%%%%%%%%%%%%%%%

%%%%%%%%%%%%%%%%%%%%%%%%%%%%%%%%%%%%%%%%%%%%%%%%%%%%%%%%%%%%%%%%%%%%%%%%%%
\begin{coro}\label{cor:stanca}
Assume $\vartheta\in\fR$ to be such that $\Val^\ttR(\vartheta;c,\om)\ne0$. 
Then for all $\ell\in L(\vartheta)$ one has
\[
\frac{1}{32}\be(m_{n_\ell}) < |\om\cdot\nu_{\ell}-\om_{j_\ell}| < \frac{1}{4}\be(m_{n_\ell-1}).
\]
\end{coro}
%%%%%%%%%%%%%%%%%%%%%%%%%%%%%%%%%%%%%%%%%%%%%%%%%%%%%%%%%%%%%%%%%%%%%%%%%%

%%%%%%%%%%%%%%%%%%%%%%%%%%%%%%%%%%%%%%%%%%%%%%%%%%%%%%%%%%%%%%%%%%%%%%%%%%
\prova
It follows from Lemma \ref{stanca} and Remark \ref{lescaleR}.
\EP
%%%%%%%%%%%%%%%%%%%%%%%%%%%%%%%%%%%%%%%%%%%%%%%%%%%%%%%%%%%%%%%%%%%%%%%%%%

%%%%%%%%%%%%%%%%%%%%%%%%%%%%%%%%%%%%%%%%%%%%%%%%%%%%%%%%%%%%%%%%%%%%%%%%%%
\begin{rmk} \label{ripetiamolo}
\emph{%remark
One can rephrase Corollary \ref{cor:stanca} by saying that all $\vartheta\in\fR$ such that $\Val^\ttR(\vartheta;c,\om)\ne0$ 
 verify the support property of Definition \ref{supporto}.
}%remark
\end{rmk}
%%%%%%%%%%%%%%%%%%%%%%%%%%%%%%%%%%%%%%%%%%%%%%%%%%%%%%%%%%%%%%%%%%%%%%%%%% 
 
%%%%%%%%%%%%%%%%%%%%%%%%%%%%%%%%%%%%%%%%%%%%%%%%%%%%%%%%%%%%%%%%%%%%%%%%%% 
\begin{rmk}\label{latte}
\emph{%
{In \eqref{gr} one has, for $h=0,1,2$, 
\[
\del^h \calG_{n} (x)=
\sum_{p=0}^{h} \frac{\Psi^{h,p}_{n}(x)}{x^{1+p}} , \qquad
|\Psi^{h,p}_{n_\ell}(x) | \le \frac{C}{|x|^{h-p} }, \qquad p=0,\ldots,h, 
\]
for suitable functions $\Psi^{h,0}_{n},\ldots,\Psi^{h,h}_n$ and some positive constant $C$, so that}
\[
| \del^h \calG_{n} (x)| 
\le  \frac{a_0}{ (\be({{m_{n}}}))^{1+h}} , \qquad h=0,1,2,
\]
for some constant $a_0 \ge 1$. %where $\del$ denotes derivative with respect to the argument.
{An explicit calculation gives
\begin{equation} \label{miele}
\Psi^{h,p}_n(x) = (-1)^p \left( \begin{matrix}  h \\ p \end{matrix} \right) \del^{h-p}\Psi_{n}(x) .
\end{equation}
In particular $\Psi^{h,h}_n(x)=\Psi_{n}(x)$ and, for $h\ge 1$ and $ p=0,\ldots,h-1$, one may have $\Psi^{h,p}_0(x) \neq 0$
only if $|x|\in(\be(m_0)/16,\be(m_0)/8)$.}
}
\end{rmk}
%%%%%%%%%%%%%%%%%%%%%%%%%%%%%%%%%%%%%%%%%%%%%%%%%%%%%%%%%%%%%%%%%%%%%%%%%% 

%%%%%%%%%%%%%%%%%%%%%%%%%%%%%%%%%%%%%%%%%%%%%%%%%%%%%%%%%%%%%%%%%%%%%%%%%%
%%%%%%%%%%%%%%%%%%%%%%%%%%%%%%%%%%%%%%%%%%%%%%%%%%%%%%%%%%%%%%%%%%%%%%%%%%
\subsection{Renormalization procedure: fully renormalized trees}
\label{LFR}
%%%%%%%%%%%%%%%%%%%%%%%%%%%%%%%%%%%%%%%%%%%%%%%%%%%%%%%%%%%%%%%%%%%%%%%%%%
%%%%%%%%%%%%%%%%%%%%%%%%%%%%%%%%%%%%%%%%%%%%%%%%%%%%%%%%%%%%%%%%%%%%%%%%%%

Now we want to combine the cancellations due to the symmetries studied in Subsection \ref{simmetrie}
with the renormalization scheme described in Subsection \ref{LR}. In order to exploit both features we introduce a further set of trees,
with the aim of ruling out since the beginning all trees whose values cancel out when summed together.

%%%%%%%%%%%%%%%%%%%%%%%%%%%%%%%%%%%%%%%%%%%%%%%%%%%%%%%%%%%%%%%%%%%%%%%%%% 
 \begin{defi}[\textbf{Set $\boldsymbol{\gotF}$ of the fully renormalized trees}]\label{terribile}
Let $\gotF$ be the subset of $\fR$ consisting of the renormalized trees $\vartheta$ such that
if 
\begin{itemize}[topsep=0ex]
\itemsep0em
\item
$\gotC=\{T_1,\ldots,T_p\}$ is a chain of RCs in $\vartheta$, with $p\ge 2$,
\item there are $q_1,q_2\in\{1,\ldots,p\}$ such that
$1\le q_1<q_2\le p$ and
$\matO_{T_{q_1}}=\matO_{T_{q_2}}=\matL$, 
\end{itemize}
then $p\ge3$, $q_2\ge q_1 + 2$ and there is $h$ such that $q_1<h<q_2$ and
$\matO_{T_h}\in\{\matR,\uno\}$. 
{We call \emph{fully renormalized tree} any tree $\vartheta\in\gotF$,}
and  $\gotF^{(k)}_{j,\nu,\s}$ the set of trees in $\gotF$ with order
$k$ such that the root line has component $j$, momentum $\nu$ and sign $\s$.
\end{defi}
%%%%%%%%%%%%%%%%%%%%%%%%%%%%%%%%%%%%%%%%%%%%%%%%%%%%%%%%%%%%%%%%%%%%%%%%%% 

%%%%%%%%%%%%%%%%%%%%%%%%%%%%%%%%%%%%%%%%%%%%%%%%%%%%%%%%%%%%%%%%%%%%%%%%%%
\begin{defi} [\textbf{Truncated fully renormalized tree}] \label{Fn}
A fully renormalized tree $\vartheta\in \gotF$ is said to be \emph{truncated} at scale $n$ if
it is a renormalized tree truncated at scale $n$ which belongs to $\gotF$.
We call $\gotF^{(k)}_{j,\nu,\s}(n)$ the set of trees $\vartheta\in \gotF^{(k)}_{j,\nu,\s}$
which are truncated at scale $n$.
\end{defi}
%%%%%%%%%%%%%%%%%%%%%%%%%%%%%%%%%%%%%%%%%%%%%%%%%%%%%%%%%%%%%%%%%%%%%%%%%%

%%%%%%%%%%%%%%%%%%%%%%%%%%%%%%%%%%%%%%%%%%%%%%%%%%%%%%%%%%%%%%%%%%%%%%%%%% 
\begin{lemma}\label{effone}
One has
\[
\sum_{\vartheta\in\gotT^{(k)}_{j,\nu,\s}}\Val(\vartheta;c,\om) 
=  \sum_{\vartheta\in\gotF^{(k)}_{j,\nu,\s}}\Val^\ttR(\vartheta;c,\om)
\]
 \end{lemma}
 %%%%%%%%%%%%%%%%%%%%%%%%%%%%%%%%%%%%%%%%%%%%%%%%%%%%%%%%%%%%%%%%%%%%%%%%%% 
 
 %%%%%%%%%%%%%%%%%%%%%%%%%%%%%%%%%%%%%%%%%%%%%%%%%%%%%%%%%%%%%%%%%%%%%%%%%% 
 \prova
 It follows from Lemma \ref{proc}, since for all $n\ge0$ one has
 \[
  \sum_{\vartheta\in\fR^{(k)}_{j,\nu,\s}(n)}\Val^\ttR(\vartheta;c,\om) 
=  \sum_{\vartheta\in\gotF^{(k)}_{j,\nu,\s}(n)}\Val^\ttR(\vartheta;c,\om),
 \]
 by Lemma \ref{dadim} and Lemma \ref{cancellazione}. Then, thanks to Lemma \ref{welldefined},
 one can take the limit $n\to\io$ and the assertion follows.
 \EP
%%%%%%%%%%%%%%%%%%%%%%%%%%%%%%%%%%%%%%%%%%%%%%%%%%%%%%%%%%%%%%%%%%%%%%%%%% 
 
%%%%%%%%%%%%%%%%%%%%%%%%%%%%%%%%%%%%%%%%%%%%%%%%%%%%%%%%%%%%%%%%%%%%%%%%%%
\begin{defi}[\textbf{Renormalized value of a trimmed tree}]\label{defvalue8}
For all $k\!\ge\!1$, $j\in\ZZZ$,  $\nu\in\ZZZ^{\ZZZ}_{f}$ and $\s\in\{\pm\}$, and for any fully renormalized tree $\vartheta\in\gotF^{(k)}_{ j ,\nu,\s}$,
define the \emph{renormalized value} of the trimmed tree $\breve\vartheta$
associated with $\vartheta$ as
\begin{equation}\label{valRtrim}
\null\hspace{-.2cm}
\Val^\ttR(\breve\vartheta;c,\om) :=
\displaystyle{\Bigg(\!\prod_{\substack{\breve T\subset \breve\vartheta \\ d(T)=0}} \!\!\!  \Val^\ttR(\breve T;c,\om) \!\!\Bigg) 
\!\! \Bigg(\prod_{\la\in \Lambda(\thetabo)} \!\!\! \LL_\la (c)\Bigg) 
%\!\! \Bigg(\prod_{v\in N(\thetao) } \!\!\! \calF_v^*(c) \Bigg)
\!\! \Bigg(\prod_{\ell\in L(\thetabo)} \!\!\! \matG_\ell^\ttR(\om)\Bigg) } . 
\end{equation}
where $\LL_\la(c)$ and $\matG_\ell^\ttR(\om)$ are defined 
in \eqref{foglie} and \eqref{gr}, respectively, while $\Val^\ttR(\breve T;c,\om)$
is recursively defined in \eqref{termounoderiregola1}--\eqref{termounoderiregola3},
with the resonant clusters $T'$ replaced by the trimmed resonant clusters $\breve T'$.
\end{defi}
%%%%%%%%%%%%%%%%%%%%%%%%%%%%%%%%%%%%%%%%%%%%%%%%%%%%%%%%%%%%%%%%%%%%%%%%%%

%%%%%%%%%%%%%%%%%%%%%%%%%%%%%%%%%%%%%%%%%%%%%%%%%%%%%%%%%%%%%%%%%%%%%%%%%%
\begin{rmk}\label{madonne}
\emph{%remark
For any $\vartheta\in\gotF$, with the notation of Definition \ref{ritrimmoful}, one has
$$ 
k(\breve\vartheta) = k - \sum_{\substack{v\in N_2(\vartheta) \\ v \st{\preceq} r_{\vartheta}}} k (\vartheta_v )
$$
and
\begin{equation} \label{eqmadonne}
\Val^\ttR(\vartheta;c,\om) ={\Val}^\ttR(\breve{\vartheta};c,\om)\prod_{v \in N_2(\vartheta)}
\Big(- \frac{1}{c_{j_{v}}^{\sigma_{v}}}
\Val^\ttR(\breve\vartheta_v;c,\om)
\Big).
\end{equation}
}
\end{rmk}
%%%%%%%%%%%%%%%%%%%%%%%%%%%%%%%%%%%%%%%%%%%%%%%%%%%%%%%%%%%%%%%%%%%%%%%%%%

%%%%%%%%%%%%%%%%%%%%%%%%%%%%%%%%%%%%%%%%%%%%%%%%%%%%%%%%%%%%%%%%%%%%%%%%%%
\begin{rmk}\label{rmk:hidden1}
\emph{
When trimming a tree $\vartheta\in\gotF$, a subgraph $\breve T$ of $\breve \vartheta$ can be a RC even though $T$ is not a RC of $\vartheta$.
Indeed,  let $T$ be a subgraph of a tree $\vartheta$ which has only one entering and one exiting line and,
despite fulfilling the conditions 2 to 5 of Definition \ref{coso}, is not a RC because it contains at least
one node $v\in N_2(\vartheta)$ such that either there exists $\ell\in L(\vartheta_v)$ with $n_\ell \ge \und{n}_T$
or $J(\breve\vartheta_v) \ge C_1 r_{m_{\und{n}_T}-1}$ (see Figure \ref{fig-hidden}).
Since $J(\breve T)=J(T)$ by construction, if $n_\ell < \und{n}_T$ for all $\ell\in L(\breve T)$,
then $\breve T$ satisfies all the conditions 1 to 5 of Definition \ref{coso}, and hence it is a RC of $\breve\vartheta$.
}
\end{rmk}
%%%%%%%%%%%%%%%%%%%%%%%%%%%%%%%%%%%%%%%%%%%%%%%%%%%%%%%%%%%%%%%%%%%%%%%%%%

%%%%%%%%%%%%%%%%%%%%%%%%%%%%%%%%%%%%%%%%%%%%%%%%%%%%%%%%%%%%%%%%%%%%%%%%%%
% FIGURA 16
%%%%%%%%%%%%%%%%%%%%%%%%%%%%%%%%%%%%%%%%%%%%%%%%%%%%%%%%%%%%%%%%%%%%%%%%%%
\begin{figure}[ht]
\vspace{.4cm}
\centering
\ins{162pt}{-090pt}{$v$}
\ins{168pt}{-037pt}{$\vartheta_v$}
\ins{033pt}{-088pt}{$\vartheta =$}
%\ins{242pt}{-089pt}{$\ell_4$}
%\ins{296pt}{-089pt}{$\ell_5$}
\ins{328pt}{-038pt}{$T$}
%\ins{358pt}{-020pt}{$T$}
%\ins{136pt}{-076pt}{$\ell_2$}
%\ins{108pt}{-079pt}{$\ell_1$}
%\ins{357pt}{-089pt}{$\ell_6$}
%\ins{180pt}{-012pt}{$j,\gote_j,\s$}
%\ins{248pt}{-012pt}{$j,\gote_j,\s$}
\subfigure{\includegraphics*[width=5in]{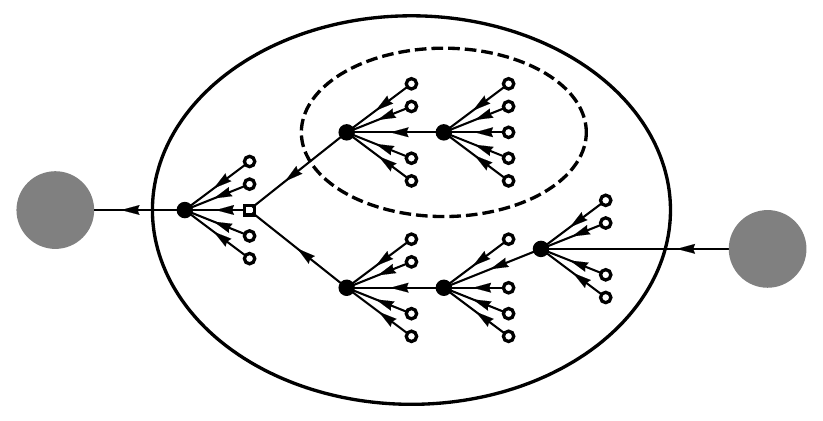}}
\caption{\small The hidden RC $\breve T$ of $\breve\vartheta$ obtained from $T$ by cutting $\vartheta_v$ from $\vartheta$.}
\label{fig-hidden}
\end{figure}
%%%%%%%%%%%%%%%%%%%%%%%%%%%%%%%%%%%%%%%%%%%%%%%%%%%%%%%%%%%%%%%%%%%%%%%%%%

In the light of Remark \ref{rmk:hidden1}, we introduce the following notions.

%%%%%%%%%%%%%%%%%%%%%%%%%%%%%%%%%%%%%%%%%%%%%%%%%%%%%%%%%%%%%%%%%%%%%%%%%%
\begin{defi}[\textbf{Hidden resonant cluster}] 
\label{defi:hidden1}
{We call \emph{hidden resonant cluster} a subgraph $\breve T$ of the trimmed tree $\breve\vartheta$
associated to a tree $\vartheta\in\gotF$, such that $\breve T$ is a RC of $\breve\vartheta$, whereas $T$ is not a RC of $\vartheta$.
A RC of $\breve\vartheta$ such that also $T$ is a RC of $\vartheta$ is said to be \emph{non-hidden} or
\emph{inherited from $\vartheta$}.}
\end{defi}
%%%%%%%%%%%%%%%%%%%%%%%%%%%%%%%%%%%%%%%%%%%%%%%%%%%%%%%%%%%%%%%%%%%%%%%%%%

%%%%%%%%%%%%%%%%%%%%%%%%%%%%%%%%%%%%%%%%%%%%%%%%%%%%%%%%%%%%%%%%%%%%%%%%%%
\begin{defi}[\textbf{Hidden resonant line}] 
\label{defi:hidden2}
{Let $\breve\vartheta$ be the trimmed tree associated to a tree $\vartheta\in\gotF$.
We say that a line $\ell$ is a hidden resonant line if it connects a hidden RC of $\vartheta$
to either a non-hidden RC or another hidden RC of $\vartheta$.}
\end{defi}
%%%%%%%%%%%%%%%%%%%%%%%%%%%%%%%%%%%%%%%%%%%%%%%%%%%%%%%%%%%%%%%%%%%%%%%%%%

%%%%%%%%%%%%%%%%%%%%%%%%%%%%%%%%%%%%%%%%%%%%%%%%%%%%%%%%%%%%%%%%%%%%%%%%%%
\begin{rmk}\label{rmk:tardi}
\emph{
As a consequence of Remarks \ref{troncoj} and \ref{rmk:hidden1}, if the tree $\vartheta\in \gotF$ satisfies the support property
and $\breve T$ is a hidden RC of $\breve\vartheta$, then there is $v\in N_2(T)$ such that
$J(\breve \vartheta_v) \ge C_1 r_{m_{\und{n}_T}-1}$.
 }
\end{rmk}
%%%%%%%%%%%%%%%%%%%%%%%%%%%%%%%%%%%%%%%%%%%%%%%%%%%%%%%%%%%%%%%%%%%%%%%%%%

%%%%%%%%%%%%%%%%%%%%%%%%%%%%%%%%%%%%%%%%%%%%%%%%%%%%%%%%%%%%%%%%%%%%%%%%%%
%\begin{rmk}\label{rmk:hidden2}
%\emph{
{The set of lines of a trimmed tree can be divided into three disjoint sets: the resonant lines, the hidden resonant lines and the non-resonant
lines.  
The definition of chain of RCs of a trimmed tree $\breve\vartheta$, with $\vartheta\in\gotF$, is adapted from Definition \ref{chain}
as follows.
%
%
%; in particular
%it may contain both hidden and non-hidden resonant clusters: if
%$\breve T$ is a hidden resonant cluster, then  the external lines of $\breve T$ are either hidden resonant lines or non-resonant lines.}
%}
%\end{rmk}
%%%%%%%%%%%%%%%%%%%%%%%%%%%%%%%%%%%%%%%%%%%%%%%%%%%%%%%%%%%%%%%%%%%%%%%%%%

%%%%%%%%%%%%%%%%%%%%%%%%%%%%%%%%%%%%%%%%%%%%%%%%%%%%%%%%%%%%%%%%%%%%%%%%%%
\begin{defi}[\textbf{Chain of trimmed resonant clusters}]
\label{chainbis}
A \emph{chain} of RCs for the trimmed tree $\breve \vartheta$ with $\vartheta \in \gotF$ is a set $\gotC=\{\breve T_1,\ldots,\breve T_p\}$, where 
$\breve T_1,\ldots,\breve T_p$, with $p \ge 2$, are RCs  of $\breve\vartheta$ such that $\ell_{\breve T_{i+1}}=\ell_{\breve T_{i}}'$ for 
$i=1,\ldots,p-1$ and both $\ell_{\breve T_1}$ and $\ell_{\breve T_p}'$ are non-resonant.
%We call $\ell_\gotC=\ell_{\breve T_p}$ and $\ell'_\gotC=\ell'_{\breve T_1}$ the 
%\emph{exiting line} and the \emph{entering line} of $\gotC$, respectively, and $\ell_{\breve T_2},\ldots,\ell_{\breve T_{p}}$ the \emph{links} of $\gotC$.
A RC $\breve T$ is called \emph{relevant} if there is a chain $\gotC$ such that $\breve T\in\gotC$:
we associate with $\breve T$ an operator label $\matO_{\breve T}=\uno$ if $\breve T$ is a hidden RC and
$\matO_{\breve T}=\matO_T$ if $\breve T$ is a inherited RC.
\end{defi}
%%%%%%%%%%%%%%%%%%%%%%%%%%%%%%%%%%%%%%%%%%%%%%%%%%%%%%%%%%%%%%%%%%%%%%%%%%

For convenience, we associate a depth label also to each hidden resonant cluster,
by adapting the definition given in the case of relevant resonant clusters (see Definition \ref{profondita}).

%%%%%%%%%%%%%%%%%%%%%%%%%%%%%%%%%%%%%%%%%%%%%%%%%%%%%%%%%%%%%%%%%%%%%%%%%%
\begin{defi}[\textbf{Depth of a relevant hidden resonant cluster}]
\label{profonditahidden}
{Given a tree $\breve\vartheta$, obtained by trimming a fully renormalized tree $\vartheta\in\gotF$,
we define the \emph{depth} $d(\breve T)$ of a relevant hidden RC $\breve T \subset \breve\vartheta$ 
recursively as follows: we set $d(\breve T)=0$ if there is no relevant hidden RC containing $\breve T$, and set
$d(\breve T)=d(\breve T')+1$ if $\breve T$ is contained in a relevant hidden RC $\breve T'$ and no other relevant hidden
RC contained in $\breve T'$ (if any) contains $\breve T$.}
\end{defi}
%%%%%%%%%%%%%%%%%%%%%%%%%%%%%%%%%%%%%%%%%%%%%%%%%%%%%%%%%%%%%%%%%%%%%%%%%%

The following result extends Proposition \ref{brjuno} to the set of the fully renormalized trees.

%%%%%%%%%%%%%%%%%%%%%%%%%%%%%%%%%%%%%%%%%%%%%%%%%%%%%%%%%%%%%%%%%%%%%%%%%%
\begin{lemma}\label{ribry}
Let $\vartheta\in\gotF$ and for all $n\ge1$ define $\gotN_n(\breve\vartheta)$ as in \eqref{Nn}.  Then one has
\[
\gotN_n(\breve\vartheta) \le \max\left\{
\frac{K}{r_{m_n-1}}J(\breve\vartheta)-2,0
\right\},
\]
with $K$ as in Proposition \ref{brjuno}.
\end{lemma}
%%%%%%%%%%%%%%%%%%%%%%%%%%%%%%%%%%%%%%%%%%%%%%%%%%%%%%%%%%%%%%%%%%%%%%%%%%

%%%%%%%%%%%%%%%%%%%%%%%%%%%%%%%%%%%%%%%%%%%%%%%%%%%%%%%%%%%%%%%%%%%%%%%%%%
\prova
Thanks to Lemma \ref{stanca} -- see also Remark \ref{ripetiamolo} -- the proof proceeds
essentially word by word as the proof of Proposition \ref{brjuno}.
\EP
%%%%%%%%%%%%%%%%%%%%%%%%%%%%%%%%%%%%%%%%%%%%%%%%%%%%%%%%%%%%%%%%%%%%%%%%%%

The following result is the analogue of \cite[Lemma 4.4]{CGP}.
 
 %%%%%%%%%%%%%%%%%%%%%%%%%%%%%%%%%%%%%%%%%%%%%%%%%%%%%%%%%%%%%%%%%%%%%%%%%% 
 \begin{lemma}\label{ilmare}
 Let $\ell$ be a resonant line with depth $d(\ell)=d\ge0$ and let $T$ be the relevant RC with highest depth containing $\ell$. Then there are in $T$
 at least two non-resonant lines with depth $d$ on scale $\ge n_{\ell}-1$.
 \end{lemma}
 %%%%%%%%%%%%%%%%%%%%%%%%%%%%%%%%%%%%%%%%%%%%%%%%%%%%%%%%%%%%%%%%%%%%%%%%%% 
 
 %%%%%%%%%%%%%%%%%%%%%%%%%%%%%%%%%%%%%%%%%%%%%%%%%%%%%%%%%%%%%%%%%%%%%%%%%% 
 \prova
 If $\ell$ is a resonant line, then it is a link of a chain of RCs. The exiting and entering lines of such a chain
 are non-resonant, are contained in $T$ and, having the same denominator as $\ell$, have the same minimum scale 
 $\nmin_\ell$ as $\ell$. Thus the assertion follows.
 \EP
 %%%%%%%%%%%%%%%%%%%%%%%%%%%%%%%%%%%%%%%%%%%%%%%%%%%%%%%%%%%%%%%%%%%%%%%%%% 
 
  %%%%%%%%%%%%%%%%%%%%%%%%%%%%%%%%%%%%%%%%%%%%%%%%%%%%%%%%%%%%%%%%%%%%%%%%%% 
\begin{rmk}\label{laspiaggia}
\emph{
If a tree $\vartheta\in\gotF$ is such that $\Val^\ttR(\vartheta;c,\om)\ne0$, then 
 Lemma \ref{ilmare}, combined together with Corollary \ref{cor:stanca}, implies that
the number of chains in $\breve\vartheta$ with a link $\ell$ on scale $n_\ell=n$ is bounded by $\frac{1}{2}\gotN_n(\breve\vartheta)$. 
}
\end{rmk}
  %%%%%%%%%%%%%%%%%%%%%%%%%%%%%%%%%%%%%%%%%%%%%%%%%%%%%%%%%%%%%%%%%%%%%%%%%% 

We are finally ready to prove the convergence of the series \eqref{formale+fout}.
 
%%%%%%%%%%%%%%%%%%%%%%%%%%%%%%%%%%%%%%%%%%%%%%%%%%%%%%%%%%%%%%%%%%%%%%%%%%
\begin{lemma}\label{converge1}
%Assume there is $C>0$ such  $|\h^{(k)}_j|\le C^k$ for all $k,j$.Then f
For all {$s_1 \ge 0$ and $s_2>0$ such that $s_1+s_2 = s$}, all $s'\in(s_1,s)$ and all $s_3\in(0,s')$, 
there is a positive constant {$A=A(s-s',s_3,\al,\om)$} such that,
for all $k\!\ge\!1$, $j\in\ZZZ$,  $\nu\in\ZZZ^{\ZZZ}_{f}$ and $\s\in\{\pm\}$,
and for any $\vartheta \in \gotF^{(k)}_{j,\nu,\s}$, one has, for $\nu\neq\gote_j$,
\begin{equation} \label{sudore}
|\Val^\ttR({\vartheta};c,\om)| \le 
A^{|N(\vartheta)|} 
e^{- s_1 |\nu|_\al} 
e^{- s_2 \jap{j}^\alpha}
e^{-(s'-s_1)J({\breve\vartheta})}\!\!\!
\prod_{v \in N_2(\vartheta)}  e^{- (s'-s_3) J({\breve\vartheta_v})} 
%\frac{1}{\be(m_{\und{n}_{\ell_v}})} \Bigr) }  ,
\end{equation}
and, for $\nu=\gote_j$,
\begin{equation} \label{sangue}
\left| \frac{1}{c_j} \Val^\ttR({\vartheta};c,\om) \right| \le 
A^{|N(\vartheta)|} 
e^{-(s'-s_1)J({\breve\vartheta})}\!\!\!
\prod_{v \in N_2(\vartheta)}  e^{- (s'-s_3) J({\breve\vartheta_v})} .
\end{equation}
\end{lemma}
%%%%%%%%%%%%%%%%%%%%%%%%%%%%%%%%%%%%%%%%%%%%%%%%%%%%%%%%%%%%%%%%%%%%%%%%%%
 
%%%%%%%%%%%%%%%%%%%%%%%%%%%%%%%%%%%%%%%%%%%%%%%%%%%%%%%%%%%%%%%%%%%%%%%%%%
\prova
Consider first the case $\nu\neq\gote_j$.
Given $\vartheta\in \gotF^{(k)}_{j,\nu,\s}$, let $\gotC=\{T_1,\ldots,T_p\}$ be a chain of resonant clusters in $\breve\vartheta$;
set $\matO_i=\matO_{T_i}$ and $\ell_i=\ell_{T_i}$, for $i=1,\ldots,p$, and $\ell_{p+1}=\ell'_{T_{p}}$, and write
\[
x_i:=x_{\ell_i}(\und{y}_{\ell_i})=\om\cdot\nu_{\ell_i}(\und{y}_{\ell_{i}}) - \om_{j_{\ell_{i}}} , \qquad \del_i=\del_{\ell_i} ,
\qquad i=1,\ldots,p+1 .
\]
In particular by definition of chain one has $|x_1|=\ldots=|x_{p+1}|$, and the lines $\ell_2,\ldots,\ell_p$
are {either resonant lines or hidden resonant lines}, whereas the lines $\ell_1$ and $\ell_{p+1}$ are non-resonant.
Moreover, by Remark \ref{latte}, one has
\begin{equation} \label{coli}
| \matG^\ttR_{\ell_i}(\om)| \le \frac{a_0}{|x_{1}|^{\del_i+1}} , \qquad i=1,\ldots,p+1.
\end{equation}
By \eqref{termounoderiregola1}--\eqref{termounoderiregola3}, for each $i=1,\ldots,p$, if $\matO_i=\matD$ the corresponding
value $\Val^\ttR(T_i;c,\om)$ contains a gain factor $x_{i+1}$, and
if $\matO_i=\matR$ it contains a gain factor $x_{i+1}^2$; {note that if $T_i$ is a hidden resonant cluster,
then one has $\matO_i=\uno$ (see Definition \ref{chainbis}).}

If we set, to simplify the notation, 
\[
\begin{aligned}
L & := \left| \{i \in \{1,\ldots,p\} : \matO_{i}= \matL \} \right|  , \qquad
D := \left| \{i \in \{1,\ldots,p\} : \matO_{i}= \matD \} \right|  , \\
R & := \left| \{i \in \{1,\ldots,p\} : \matO_{i}= \matR \} \right|  , \qquad
O := \left| \{i \in \{1,\ldots,p\} : \matO_{i}= \uno \} \right| , 
\end{aligned}
\]
{with $O$ splitted as $O=I + H$, where
\[
\begin{aligned}
I & = \left| \{i \in \{1,\ldots,p\} : \matO_{i}= \uno \hbox{ and } T_{i} \hbox{ is a inherited RC} \} \right| , \\
H & = \left| \{i \in \{1,\ldots,p\} : \matO_{i}= \uno \hbox{ and } T_{i} \hbox{ is a hidden RC} \} \right| , 
\end{aligned}
\]
and}
\[
d := \sum_{i=2}^{p} \del_i ,
\]
then Definitions \ref{nonepossopiu} and \ref{terribile} imply, respectively,  that
\begin{equation}\label{sempre2}
{I + d \le 2} , \qquad L-1 \le R + O .
%O + d \le 2 , \qquad L-1 \le R + O .
\end{equation} 
Since $L+D+R+O=p$, we can bound
\begin{equation} \label{stimacatena}
 \Biggl( \prod_{i=2}^{p} |\matG^\ttR_{\ell_i}(\om)| \Biggr)
%|x_{1}|^{|\{ i \in\{1,\ldots,p\} \,:\, \matO_{T_i}=\matD\}| + 2|\{ i \in\{1,\ldots,p\} \, : \, \matO_{T_i}=\matR\}|}
|x_{1}|^{D + 2R} \le {\frac{a_0^2}{|x_{1}|^{4+2H}} },
\end{equation}
where we have used that
\begin{equation} \label{nostimacatena}
\begin{aligned}
-d-(p-1) + D  + 2R 
& = -d-(p-1) + L + D + R + O - L - O + R  \\
& = - d - (L - 1 - R)  - O  \ge - 2 O - d \\
& \ge {-2(O+d) \ge -2 (I+d) - 2H \ge -4-2H } ,
%
%-d-(p-1) + D  + 2R 
%& = -d-(p-1) + L + D + R + O - L - O + R  \\
%& = - d - (L - 1 - R)  - O  \ge - 2 O - d \ge -2(d+O) \ge -4 ,
\end{aligned}
\end{equation}
which follows from \eqref{sempre2}.
Thus, in all cases
the product of the propagators of the links of the chain times
the gain factors is bounded by $a_0^2|x_{p+1}|^{-4}$, independently of the length of the 
chain, %\footnote{Recall the factor $|x_{p+1}|^{-4}$ in \eqref{cespero}.}
{times a further factor $|x_{1}|^{-2H}$ due to hidden resonant clusters.}

{Let $\fC(\breve\vartheta)$ denote the set of the chains of RCs  in $\breve\vartheta$.}
For any chain $\gotC \in \fC(\breve\vartheta)$, if $p(\gotC)$ denotes the number of RCs in $\gotC$,
call $\ell_2(\gotC),\ldots\ell_{p(\gotC)}(\gotC)$ the links of $\gotC$, and
$\ell_1(\gotC)$ and $\ell_{p(\gotC)+1}(\gotC) $ the exiting line and the entering line of $\gotC$, respectively.
Then set $x_{\ell_i(\gotC)} = x_{\ell_i(\gotC)} (\und{y}_{\ell_{i}(\gotC)})$,
for $i=1,\ldots,p(\gotC)+1$; by construction, one has $|x_{\ell_i(\gotC)}|=|x_{\ell_1(\gotC)}|$
for all $i=2,\ldots,p(\gotC)+1$. 
Then we proceed as in Lemma \ref{convergerebbe}, and -- using the bounds \eqref{carbonara} and
\eqref{stimacatena}, Corollary \ref{cor:stanca}, Remark \ref{laspiaggia} and Lemma \ref{ribry} -- {we find, for any $n_1\ge 1$},
\begin{equation}\label{telanumero}
\begin{aligned}
\null\hspace{-.3cm} 
& \Biggl(\prod_{\ell\in L_{\!N\!R}(\vartheta)}|\matG^{\ttR}_\ell(\om)|\Biggr)
\!\! \prod_{{\gotC \in \fC(\vartheta)}} \Biggl(
|x_{\ell_1(\gotC)}|^{{|\{ T\in \gotC \, : \, \matO_T=\matD\} + 2|\{ T \in\gotC \, : \, \matO_T=\matR\}|}}
 |\prod_{i=2}^{p(\gotC)}|\matG^{\ttR}_{\ell_i(\gotC)}(\om)|
\Biggr)\\
\null\hspace{-.3cm} 
& \quad \le
\prod_{n\ge0} \left(\frac{{2^{20}}a_0^2}{(\be(m_n))^5} \right)^{\gotN_n(\breve\vartheta)} \!\!\!\!
{\prod_{{\gotC \in \fC(\breve\vartheta)}} \!\!\! \Bigl( |x_{\ell_1(\gotC)}|^{-2H(\gotC)}\Bigr)} ,
%\\
%\null\hspace{-.3cm} 
%& \quad \le  
%(A_1(n_1))^{|N(\breve\vartheta)|} 
%\exp \!\! \left(5K J(\breve\vartheta){\sum_{n\ge n_1+1} \frac{1}{r_{m_{n}-1}} \log\left( \frac{2^4a_0^{2/5}}{\be({m_n})} \right)}\right) 
%\!\! {\prod_{{\gotC \in \fC(\breve\vartheta)}} \!\!\! \Bigl( |x_{\ell_1(\gotC)}|^{-2H(\gotC)}\Bigr) } ,
\end{aligned}
\end{equation}
where $H(\gotC)$ is the number of hidden resonant clusters contained in $\gotC$.
%and
%%
%\begin{equation}\label{guardata1}
%A_1(n):= \frac{2^{20}a_0^2}{(\be(m_n))^5} \,.
%\end{equation}
%%
In order to bound the last product in \eqref{telanumero} we reason as follows.
\begin{itemize}[topsep=0ex]
\itemsep0em
\item Let $\breve T$ be any relevant hidden RC of $\breve\vartheta$ with depth 0.
By construction there is at least one node $v \in N_2(\vartheta)$ such that $\breve \vartheta_v$ is contained in $T$ and
$J(\breve \vartheta_v) \ge C_1 r_{m_{\und{n}_T}-1}$ (see Remark \ref{rmk:tardi}). 
Take the node $v$ such that $\breve\vartheta_v$ contains the line with the highest scale in $T$;
if there are more nodes satisfying such a property, we chose any of them arbitrarily.
We say that the node $v$ is associated with $\breve T$. % and write $\breve T = \breve T(v)$.
Call $\calmT_{0}(\breve\vartheta)$ the set of the relevant hidden RCs of $\breve\vartheta$ with depth 0 and
$\calmN_{0}(\vartheta)$ the set of nodes $v\in N_2(\vartheta)$ associated to the relevant hidden RCs in $\calmT_0(\breve\vartheta)$.
By construction there is a one-to-one correspondence between the sets $\calmN_{0}(\vartheta)$ and $\calmT_0(\breve\vartheta)$: 
for any $v\in \calmN_{0}(\vartheta)$ we call $\breve T(v)$ the relevant hidden RC in $\calmT_0(\breve\vartheta)$ to which $v$ is associated.
\item Next, we consider any relevant hidden RC $\breve T'$ with depth $1$:
again there is at least  one subtree $\vartheta_{v} \in F_2^*(\vartheta)$ such that
$\vartheta_{v}\subset T'$ and  $J(\breve \vartheta_{v})\ge C_1 r_{m_{\und{n}_{T'}}-1}$. 
Consider all the nodes $v$ such that $\breve\vartheta_{v}$ contains a line with the highest scale in $T'$:
if any of them belongs to $\calmN_{0}(\vartheta)$, we say that
$\breve T'$ is not associated to any node of $N_2(\vartheta)$. 
If, instead, one has $v\notin \calmN_{0}(\vartheta)$ for all $v$ as above, 
we chose any of them arbitrarily, say $v'$, and say that $v'$ is associated with $\breve T'$. % and write $\breve T = \breve T(v)$.
Call $\calmT_{1}(\breve\vartheta)$ the set of the relevant hidden RCs of $\breve\vartheta$ with depth 1 to which a node has been associated, and
$\calmN_{1}(\vartheta)$ the set of nodes $v\in N_2(\vartheta)$ associated to the relevant hidden RCs in $\calmT_1(\breve\vartheta)$.
For any $v\in \calmN_{1}(\vartheta)$ we call $\breve T(v)$ the relevant hidden RC in $\calmT_1(\breve\vartheta)$ to which $v$ is associated.
Note that some relevant hidden RCs with depth 1 have not been associated to any node.
\item We iterate the construction above as many times as we need to reach the relevant hidden resonant clusters with maximum depth.
If $D$ denotes such a depth, set
$$ \calmT(\breve\vartheta) = \bigcup_{i=0}^D \calmT_i(\breve\vartheta) , \qquad \qquad
\calmN(\vartheta) = \bigcup_{i=0}^D \calmN_i(\vartheta) .$$
\item In this way, we may divide the set of all relevant hidden resonant clusters into two sets: the set containing the relevant hidden resonant clusters
which are associated to the nodes in $\calmN(\vartheta)$ and the set containing all the relevant hidden resonant clusters which are not associated
to any node in $\calmN(\vartheta)$.
\item We observe that, for each chain $\gotC$, there is at most one hidden resonant cluster which is not associated to any node $v\in \calmN(\vartheta)$.
This can be seen as follows. First of all, all relevant RCs with depth 0 are associated with some node in $\calmN(\vartheta)$,
thus we can confine ourselves to the chains $\gotC$ whose
RCs have depth $\ge1$.
For each relevant hidden resonant cluster $\breve T$ containing $\gotC$ the following cases are possible:
\begin{enumerate}[topsep=0ex]
\itemsep0em
\item either $\breve T$ is not associated with any node $v\in\calmN(\vartheta)$,
\item or $\breve T$ is associated with a node $v\in N(\breve T')$, for some $\breve T' \notin \gotC$,
\item or $\breve T$ is associated with a node $v \in N(\breve T')$, for some $\breve T' \in \gotC$.
\end{enumerate}
If the third case does not occur for any of the considered hidden resonant clusters, then
each hidden relevant cluster in $\gotC$ is associated with a node in $\calmN(\vartheta)$ by construction. %such that $\vartheta_v \subset T'$.
If the third case arises, then all the relevant hidden RCs $\breve T''$ contained in $\breve T$ and containing $\gotC$
are not associated with any node in $N_2(\vartheta)$,
because $v$ is such that $\vartheta_{v}\subset \breve T''$ and $\breve\vartheta_v$ is one of subtrees containing the line with the highest scale.
If $T'$ is the subgraph such that $\breve T'\in \gotC$ and $v\in N(T')$, all the other hidden resonant clusters of $\gotC$ are associated
to nodes in $\calmN(\vartheta)$ which are outside $T'$. %Call $\calmN(\theta;\gotC)$ the set of such nodes.
\end{itemize}
As a consequence of the discussion above, we may write
\begin{equation} \label{finiramai1}
\prod_{{\gotC \in \fC({\breve\vartheta})}} \Bigl( |x_{\ell_1(\gotC)}|^{-H(\gotC)}\Bigr) =
\prod_{{\gotC \in \fC({\breve\vartheta})}} \Bigl( |x_{\ell_1(\gotC)}|^{-1}\Bigr) 
\prod_{v \in \calmN(\vartheta)} \Bigl( |x_{\ell_{\breve T(v)}}|^{-1}\Bigr) ,
\end{equation}
where the last factor is bounded as
\begin{equation} \label{finiramai2}
\prod_{v \in \calmN(\vartheta)} \Bigl( |x_{\ell_{T(v)}}|^{-1}\Bigr) \le
%\prod_{v \in N_2(\vartheta)} \Bigl( |x_{\ell_v}|^{-1}\Bigr) \le
\prod_{v \in \calmN(\vartheta)} \frac{2^4}{\be
(m_{n(v)})} ,
%(m_{\und{n}_{\ell_{\breve T(v)}}})} ,
\end{equation}
where, to simplify the notation, we have shortened $n(v):= \und{n}_{\ell_{\breve T(v)}}$.
Therefore, bounding once more the number of chains with a link $\ell$ on scale $n_\ell=n$ by $\frac{1}{2}\gotN_n(\breve\vartheta)$,
combining \eqref{telanumero} with \eqref{finiramai1} and \eqref{finiramai2} gives, for any $n_1 \ge 1$,
\begin{equation}\label{telarinumero}
\begin{aligned}
\null\hspace{-.4cm} 
& \Biggl(\prod_{\ell\in L_{\!N\!R}(\vartheta)}|\matG^{\ttR}_\ell(\om)|\Biggr)
\!\! \prod_{{\gotC \in \matC(\vartheta)}} \Biggl(
|x_{\ell_1(\gotC)}|^{{|\{ T\in \gotC \, : \, \matO_T=\matD\} + 2|\{ T \in\gotC \, : \, \matO_T=\matR\}|}}
 |\prod_{i=2}^{p(\gotC)}|\matG^{\ttR}_{\ell_i(\gotC)}(\om)|
\Biggr)\\
\null\hspace{-.4cm} 
& \quad \le  
(A_2(n_1))^{|N(\breve\vartheta)|} 
\exp \! \Biggl(6K J(\breve\vartheta){\sum_{n\ge n_1+1} \frac{1}{r_{m_{n}-1}} \log\left( \frac{2^4 a_0^{1/3}}{\be({m_n})} \right)}\Biggr) 
\!\! {\prod_{v \in \calmN(\vartheta)} \Bigl( \frac{2^4}{\be
(m_{n(v)})} \Bigr)^2 } ,
%(m_{\und{n}_{\ell_{T(v)}}})} \Bigr)^2 } ,
\end{aligned}
\end{equation}
with
\begin{equation} \label{guardata2}
A_2(n):= \frac{2^{24}a_0^2}{\be(m_n)^6} = \Biggl( \frac{2^4 a_0^{1/3}}{\be({m_n})} \Biggr)^6 .
\end{equation}
Now given any $\de>0$ we define  $n_1(\de)$ the smallest integer such that
\begin{equation}\label{n1}
%\ka_*(n_1):= 
6 K \sum_{n\ge n_1(\de)+1} \frac{1}{r_{m_{n}-1}} \log\Biggl( \frac{2^4 a_0^{1/3}}{\be({m_n})} \Biggr) \le \de.
\end{equation}
The product of the node factors is $1$, while the product of the leaf factors is bounded as in \eqref{assai}.
Overall, we obtain
\begin{equation}\label{stastima}
|\Val^\ttR(\breve\vartheta;c,\om)| \le 
(A_2(n_1(s-s')))^{|N(\breve\vartheta)|} 
e^{-s'J(\breve\vartheta)}e^{-s\jap{j}^\alpha} 
{\prod_{v \in \calmN(\vartheta)} \Bigl( \frac{2^4}{\be
(m_{n(v)}) } \Bigr)^2 } .
%(m_{\und{n}_{\ell_v}})} \Bigr)^2 } .
\end{equation}

We can repeat the construction for any trimmed tree $\breve\vartheta_v$, with $v\in N_2(\vartheta)$,
by following the inclusion order -- this includes also the case of $\breve\vartheta$ with $\nu=\gote_j$ left out in the discussion above.
The only differences are the following.
\begin{itemize}[topsep=0ex]
\itemsep0em
\item The renormalized value of a tree $\vartheta\in\gotF^{(k)}_{j,\gote_j,\s}$ is multiplied by a further factor $-1/c_j^\s$,
but using Lemma \ref{cijei} -- which extends trivially to trees in $\gotF^{(k)}_{j,\gote_j,\s}$ -- one finds
\begin{equation} \label{servedopo}
\left| \frac{1}{c_j^\s}\prod_{\la\in \Lambda(\breve\vartheta)} \LL_\la(c)\right|\le e^{-s J(\breve\vartheta)} .
\end{equation}
%
%so that, following the proof of Proposition \ref{voltabuona} one sees that for all $0<s'<s$
%%
%\begin{equation}\label{stastima2}
%	\left| \frac{1}{c_j^\s}\Val^\ttR(\vartheta;c,\om)\right| \le 
%	A_1(n_1)^k 
%	e^{-s'J(\vartheta)},
%\end{equation}
%%
%(see \eqref{guardata} and \eqref{n1} for the definition of $A_1(n)$ and $n_1=n_1(s-s')$ respectively) 
%
\item When constructing the set $\calmN(\vartheta_v)$,
every time one considers a relevant hidden resonant cluster $\breve T$ of $\breve\vartheta_v$ and identifies
the subtree $\vartheta_{v'}$ with $v' \in N_2(\vartheta)$ such that $T$ contains $\breve \vartheta_{v'}$
and $\breve \vartheta_{v'}$ contains at least one line $\ell$ with scale $n_{\ell} > \und{n}_{\breve T}$,
if $v'$ belongs to some set $\calmN(\vartheta_{v''})$ constructed in a previous step, then
$T$ is not associated with any node in $N_2(\vartheta)$. In that way,
at most one factor $4/\be(m_{n(v)})$ arises for each node $v\in N_2(\vartheta)$
and, for any chain, at most two hidden resonant lines are not associated to any node $v\in N_2(\vartheta)$.
\end {itemize}

Define $\calmN_*(\vartheta)$ as the set of nodes in $N(\vartheta)$ which belong either
to $\calmN(\vartheta)$ or to $\calmN(\vartheta_v)$ for some $v\in N_2(\vartheta)$.  
Eventually we find
\begin{equation} \nonumber %label{eqmadonnestime}
| \Val^\ttR(\vartheta;c,\om) | \le
(A_2(n_1(s-s')))^{|N(\vartheta)|} 
e^{- s' J(\breve\vartheta)} 
e^{- s \jap{j}^\alpha} \!\!\!\!\!
\prod_{v \in N_2(\vartheta)} \!\!\!\! \left(
e^{-s'J(\breve\vartheta_v)} 
\right) \!\!
{\prod_{v \in \calmN_*(\vartheta)} \Bigl( \frac{4}{\be
(m_{n(v)}) } \Bigr)^2 } 
%(m_{\und{n}_{\ell_{T(v)}}})} \Bigr)^2 } 
\end{equation}
and hence, for any $s''\ge0$ such that $s''<s'$, we may write
\[
\begin{aligned}
|\Val^\ttR({\vartheta};c,\om)| &\le 
(A_2(n_1(s-s')))^{|N(\vartheta)|} 
e^{-s''|\nu|_\al} 
e^{-(s-s'')\jap{j}^\alpha}
e^{-(s'-s'')J({\breve\vartheta})}\\
&\qquad
\prod_{v \in N_2(\vartheta)} \!\!\!\! \left(
e^{-(s'-s_3)J(\breve\vartheta_v)} 
\right) \!\!
{\prod_{v \in \calmN_*(\vartheta)} \Bigl( \frac{4}{\be
(m_{n(v)})} \Bigr)^2
%(m_{\und{n}_{\ell_{T(v)}}})} \Bigr)^2
e^{-s_3J(\breve\vartheta_v)} } ,
\end{aligned}
\]
for some $s_3<s'$,
where we have used that $J(\breve\vartheta) + \jap{j}^\al\ge |\nu|_\al $ by \eqref{seserve}).
Finally, using that $J(\breve\vartheta_v)\ge C_1 r_{m_n-1}$ for some $n\ge n(v)$ %r_{m_{\und{n}_{\ell_{T(v)}}}-1}$ 
by construction, we have
\[
e^{-s_3 J(\breve\vartheta_v)}\Bigl( \frac{4}{\be
(m_{n(v)})} 
%(m_{\und{n}_{\ell_{T(v)}}})} 
\Bigr)^2\le
\exp \left( - C_1 r_{
m_{n(v)}-1
%m_{n_{\ell_{T(v)}}}-1
} \Biggl( s_3 - \frac{2}{r_{
m_{n(v)}-1
%m_{n_{\ell_{T(v)}}-1
}} \log \frac{4}{\be
(m_{n(v)})
%(m_{n_{\ell_{T(v)}}})
} \Biggr) \right) \le 1,
\]
if $n(v) > n_2(s_3)$, for a suitable $n_2(s_3)\le n_1(s_3)$, as a consequence
of the Bryuno condition  together with Lemma \ref{basta!}, and
\[
e^{-s_3 J(\breve\vartheta_v)}\Bigl( \frac{4}{\be
(m_{n(v)})
%(m_{\und{n}_{\ell_{T(v)}}})
} \Bigr)^2 \le \Bigl( \frac{4}{\be
(m_{n(v)})
%(m_{\und{n}_{\ell_{T(v)}}})
} \Bigr)^2  \le (A(n_2(s_3)))^2 \le A_2(n_1(s_3)) ,
\]
with $A(n)$ defined after \eqref{carbonara}, if $n(v) \le n_2(s_3)$.
Then the bound \eqref{sudore} follows, with 
\begin{equation}\label{pe}
A(s-s',s_3,\al,\om):=\max\{A_2(n_1(s-s')),A_2(n_1(s_3))\}=A_2(n_1(\min\{s-s',s_3\}))
\end{equation}
and $s_1=s''$.

The bound \eqref{sangue} follows immediately because, as already noted, if $\nu=\gote_j$ the trimmed tree $\breve\vartheta$ behaves as one of the
trimmed trees $\breve\vartheta_v$, with $v \in N_2(\vartheta)$, of the construction above.
\EP
%%%%%%%%%%%%%%%%%%%%%%%%%%%%%%%%%%%%%%%%%%%%%%%%%%%%%%%%%%%%%%%%%%%%%%%%%% 

With the notation of Definitions \ref{trimmed} and \ref{ritrimmoful}, for any tree $\vartheta\in \gotT \cup \gotF$ we define
\begin{equation} \label{allafine}
\calJ(\vartheta) := J(\breve\vartheta) + \sum_{v\in N_2(\vartheta)}J(\breve\vartheta_v) .
\end{equation}

%%%%%%%%%%%%%%%%%%%%%%%%%%%%%%%%%%%%%%%%%%%%%%%%%%%%%%%%%%%%%%%%%%%%%%%%%% 
\begin{coro}\label{trenta}
Let $s'$, $s_3$ and $A$ be as in Lemma \ref{converge1}.
For all $k\ge 1$, $j\in\ZZZ$, all $\hat s, \bar s >0$ such that $\hat s+\bar s=s'-s_3$, and
all $\vartheta\in\gotF^{(k)}_{j,\gote_j,+}$ such that either $|\Lambda_{j,+}(\vartheta)| \ge 2$ or
$|\Lambda_{j',+}(\vartheta)| \ge 1$ for some $j'> j$, one has
\begin{equation}\label{giunta}
\left| \frac{1}{c_j}\Val^\ttR (\vartheta;c,\om) \right| \le A^{|N(\vartheta)|} 
e^{-2 \hat s \jap{j}^{\al}}  e^{- \bar  s  \calJ(\vartheta)} .
 \end{equation}
\end{coro}
%%%%%%%%%%%%%%%%%%%%%%%%%%%%%%%%%%%%%%%%%%%%%%%%%%%%%%%%%%%%%%%%%%%%%%%%%% 

%%%%%%%%%%%%%%%%%%%%%%%%%%%%%%%%%%%%%%%%%%%%%%%%%%%%%%%%%%%%%%%%%%%%%%%%%% 
\prova
In \eqref{sangue}, taking $s_1 \le s_3$, use \eqref{allafine} and write $s'-s_3=\bar s + \hat s$.
By Lemma \ref{cijei}, if $|\Lambda_{j,+}(\vartheta)| \ge 2$, then $|\Lambda_{j,-}(\vartheta)| \ge 1$;
analogously, if $|\Lambda_{j',+}(\vartheta)| \ge 1$ for some $j'\neq j$, then $|\Lambda_{j',-}(\vartheta)| \ge 1$.
Thus, one uses that the sum $\hat s \calJ(\vartheta)$ contains al least two summands $\ge \hat s \jap{j}^\al$.
\EP
%%%%%%%%%%%%%%%%%%%%%%%%%%%%%%%%%%%%%%%%%%%%%%%%%%%%%%%%%%%%%%%%%%%%%%%%%% 

%%%%%%%%%%%%%%%%%%%%%%%%%%%%%%%%%%%%%%%%%%%%%%%%%%%%%%%%%%%%%%%%%%%%%%%%%% 
\begin{prop}\label{voltabuona}
For all $k\ge 1$, $j\in\ZZZ$, $\nu\in\ZZZ^\ZZZ_f$ and $\s\in\{\pm\}$,
and all $s_1 \ge 0$ and $s_2>0$ such that $s_1+s_2=s$, one has
\begin{subequations}\label{arghh}
\begin{align}
&\sum_{\vartheta \in \gotF^{(k)}_{j,\nu,\s}} 
%e^{s' J(\vartheta)}
 \! | \Val^\ttR (\vartheta;c,\om) | \le D_0^k %e^{-s' |\nu-\gote_j|_{\al/2}}
e^{-s_1 |\nu|_\al}e^{-s_2\jap{j}^\al } , &\nu\ne\gote_j,\label{ahia1} \\
&\sum_{\vartheta \in \gotF^{(k)}_{j,\gote_j,\s}} 
%e^{s' J(\vartheta)}
 \! \left|\frac{1}{c_j} \Val^\ttR (\vartheta;c,\om) \right| \le D_0^k , &\nu=\gote_j,\label{ahia2}
 \end{align}
\end{subequations}
for some positive constant $D_0=D_0(s_2,\al,\om)$.
\end{prop}
%%%%%%%%%%%%%%%%%%%%%%%%%%%%%%%%%%%%%%%%%%%%%%%%%%%%%%%%%%%%%%%%%%%%%%%%%% 

%%%%%%%%%%%%%%%%%%%%%%%%%%%%%%%%%%%%%%%%%%%%%%%%%%%%%%%%%%%%%%%%%%%%%%%%%% 
\prova
Fix $s'=s-s_2/5$ and $s_3=s_1+s_2/5$, so that $s'-s_1=4s_2/5$ and $s'-s_3=3s_2/5$, which inserted in \eqref{sudore} gives
\begin{equation} \label{cheppppp}
|\Val^\ttR({\vartheta};c,\om)| \le 
(A(n_1(s_2/5))^{|N(\vartheta)|} 
e^{- s_1 |\nu|_\al} 
e^{- s_2 \jap{j}^\alpha}
e^{-3s_2\calJ({\breve\vartheta})/5}.
\end{equation}
where $A_2(n)$ is defined in \eqref{guardata2} and we used that,
by \eqref{pe}, $A(s_2/5,s_1+s_2/5,\al,\om)= A_2(n_1(s_2/5))$.
Then we
apply Lemma \ref{convergerebbe2} with $s''$ replaced with $3s_2/5$ to perform the sum over the component labels, and,
%Finally, 
using that the number of unlabelled trees of order $k$ can be bounded  by {$2^{10k+2}$}, we arrive at
\[
\begin{aligned}
\sum_{\vartheta\in\gotF^{(k)}_{j,\nu,\s}} |\Val^\ttR({\vartheta};c,\om)|
& \le D_3^k (A_2(n_1(s_2/5)))^{4k}  e^{-s_1|\nu|_\al} e^{-s_2\jap{j}^\alpha} ,
\end{aligned}
\]
where we have bounded $|N(\vartheta)| \le 2k$, for some constant $D_3=D_3(s_2)$. One reasons in the same way for \eqref{ahia2}.
Then the assertion follows with $D_0=D_3(s_2) (A_2(n_1(s_2/5)))^{4}$.
\EP
%%%%%%%%%%%%%%%%%%%%%%%%%%%%%%%%%%%%%%%%%%%%%%%%%%%%%%%%%%%%%%%%%%%%%%%%%% 

%%%%%%%%%%%%%%%%%%%%%%%%%%%%%%%%%%%%%%%%%%%%%%%%%%%%%%%%%%%%%%%%%%%%%%%%%% 
\begin{rmk}\label{toner}
\emph{
The constant $D_0$ in Proposition \ref{voltabuona} is smooth in $s_2$  and decreasing. Moreover it  depends on $\om$ because the constant $A_2(n)$ in \eqref{guardata2} does,
through the sequence $\be(m_n)$.
More precisely, the larger the value of $\BB_\om(r)$, the smaller is the estimate
on the the radius of convergence of the series: indeed, if $\BB_\om(r)$ is very large one has to take a very large value $n_1$
in \eqref{n1}, which yields a very large value of the constant $A=A_2(n_1(\min\{s-s',s_3\}))$ in Lemma \ref{converge1}: this in turns means that
only values of $\e$ in a very small interval around the origin are allowed.
}
\end {rmk}
%%%%%%%%%%%%%%%%%%%%%%%%%%%%%%%%%%%%%%%%%%%%%%%%%%%%%%%%%%%%%%%%%%%%%%%%%% 

Let us now conclude the proof of Theorem \ref{moser}.

%%%%%%%%%%%%%%%%%%%%%%%%%%%%%%%%%%%%%%%%%%%%%%%%%%%%%%%%%%%%%%%%%%%%%%%%%% 
\begin{prop}\label{bellaprop}
For %all  $\al\in(0,1)$, all  $s>0$ and 
all $s_1\ge0$, $s_2>0$ such that $s_1+s_2=s$, the following holds:
\begin{enumerate}[topsep=0ex]
\itemsep0em
\item for all $k\ge 1$, $j\in\ZZZ$, $\nu\in\ZZZ^\ZZZ_f$ and $\s\in\{\pm\}$, one has 
\begin{subequations} \label{settotutto}
\begin{align}
\h_j^{(k)}(c,\om)&= -
\sum_{\vartheta \in {\gotF}^{(k)}_{j,\gote_j,+}}  \frac{1}{c_j}\Val^\ttR (\vartheta;c,\om), \label{settoeta}\\
u^{(k)}_{j,\nu}(c,\om)&=\sum_{\vartheta \in {\gotF}^{(k)}_{j,\nu,+}}  \Val^\ttR (\vartheta;c,\om),\qquad
\ol{u}^{(k)}_{j,\nu}(c,\om)=\sum_{\vartheta \in {\gotF}^{(k)}_{j,\nu,-}}  \Val^\ttR (\vartheta;c,\om),\qquad
 \nu\ne\gote_j; 
\label{settou}
\end{align}
\end{subequations}
\item
for all $k\ge 1$, $j\in\ZZZ$, $\nu\in\ZZZ^\ZZZ_f$ and $\s\in\{\pm\}$, one has 
\begin{equation} \label{poverestime}
|\h^{(k)}_j(c,\om)|\le D_0^k, \qquad |u^{(k)}_{j,\nu}(c,\om)| \le 
D_0^k e^{-s_1|\nu|_\al}e^{-s_2\jap{j}^\al},
\end{equation}
with $D_0=D_0(s_2,\al,\om)$ is as in Proposition \ref{voltabuona}.
\end{enumerate}
Then
\begin{equation}
	\label{converga}
 %\quad u_{j,\nu}(c,\e) = \begin{cases}
%	{ \sum_{k\ge1} \e^k u_{j,\nu}^{(k)}(c) \quad\mbox{if}\; \nu \ne \mathfrak e_j}\\
%	 c_j \quad \mbox{otherwise}\end{cases}
	U(x,\f; c,\om,\e):= \sum_{j\in\ZZZ}c_j e^{\ii (\varphi_j + jx)}+ \sum_{j\in\ZZZ}
	\sum_{\nu\ne \mathfrak e_j} u_{j,\nu}(c,\om,\e) e^{\ii(\varphi\cdot\nu+j x )}\,,
\end{equation}
where
\begin{equation}\label{purequesta}
 u_{j,\nu}(c,\om,\e) :=  \sum_{k\ge1} \e^k u_{j,\nu}^{(k)}(c,\om) ,
\end{equation}
and
\begin{equation}\label{ancheleidinuovo}
 \h_j(c,\om,\e) := \sum_{k\ge1} \e^k \h_j^{(k)}(c,\om) 
\end{equation}
are well-defined as absolutely convergent series for $|\e|<\e_0=\e_0(s,\al,\om):=(2D_0)^{-1}(s,\al,\om)$. % and solve \eqref{mtoro}.
Moreover  $U$ solves \eqref{mtoro} and $\eta$ satisfies \eqref{oh1}.
Finally both $U(x,\f; c,\om,\e)$ and $ \h(c,\om,\e)$ are separately analytic in $c,\ol{c}$ for $c\in\matU_1(\mathtt{g}(s,\al))$. 
\end{prop}
%%%%%%%%%%%%%%%%%%%%%%%%%%%%%%%%%%%%%%%%%%%%%%%%%%%%%%%%%%%%%%%%%%%%%%%%%% 

%%%%%%%%%%%%%%%%%%%%%%%%%%%%%%%%%%%%%%%%%%%%%%%%%%%%%%%%%%%%%%%%%%%%%%%%%% 
\prova
Proposition \ref{voltabuona} together with Lemma \ref{effone} ensures the absolute convergence of the series \eqref{converga}, 
\eqref{purequesta} and \eqref{ancheleidinuovo}, for $\e\in(-\e_0,\e_0)$, uniformly in $x,\f,c$. 
The  equalities \eqref{settoeta} and \eqref{settou}, follow from \eqref{veri} and from  Lemma \ref{effone}. 
Since the series are uniformly convergent toghether with their time and space derivatives, then $U$ defined above is a solution of \eqref{mtoro}. 
The bound \eqref{oh1} follows directly from \eqref{poverestime}., the symmetry property comes from Lemma \ref{cijei} and Remark \ref{singolare}.
It only remains  to show the analyiticity in $c\in\matU_1(\mathtt{g}(s,\al))$ of each coefficient $u_{j,\nu}^{(k)}(c,\om)$, $\h_j^{(k)}(c,\om)$. 
Noting that $u_{j,\nu}^{(k)}(c,\om)$, $\h_j^{(k)}(c,\om)$ are bounded polynomials in $c,\ol{c}$ of degree $\le 4k+1$,  the assertion follows.
\EP
%%%%%%%%%%%%%%%%%%%%%%%%%%%%%%%%%%%%%%%%%%%%%%%%%%%%%%%%%%%%%%%%%%%%%%%%%% 

%%%%%%%%%%%%%%%%%%%%%%%%%%%%%%%%%%%%%%%%%%%%%%%%%%%%%%%%%%%%%%%%%%%%%%%%%% 
\begin{rmk}\label{trimbo}
\emph{
From Proposition \ref{bellaprop} it follows that \eqref{purequesta} and \eqref{ancheleidinuovo} converge for 
$\e$ small enough, and, by comparing \eqref{settotutto} with \eqref{veri}, we conclude that they solve \eqref{mtoro}.
It is important to stress that, since $\inf\{\e_0(s,\al,\om):\om\in\gotB\}=0$, the radius of convergence is not uniform in $\om$.
}
\end{rmk}
%%%%%%%%%%%%%%%%%%%%%%%%%%%%%%%%%%%%%%%%%%%%%%%%%%%%%%%%%%%%%%%%%%%%%%%%%% 

%%%%%%%%%%%%%%%%%%%%%%%%%%%%%%%%%%%%%%%%%%%%%%%%%%%%%%%%%%%%%%%%%%%%%%%%%% 
\begin{coro}\label{contenti}
With the notation of Proposition \ref{bellaprop} let us define
\[
U_\perp(x,\f;c,\om,\e):=\sum_{j\in\ZZZ}\sum_{\nu\ne \mathfrak e_j} u_{j,\nu,+}(c,\om,\e) e^{\ii(\varphi\cdot\nu+j x )}
\]
For all $\e\in(-\e_0,\e_0)$ there exists $s'_0=s'_0(\e)\in(0,s)$ such that for all $s_2\in(s'_0,s)$ and $s_1=s-s_2$, there is $\Phi_2=\Phi_2(s_2,\al,\om)$
such that \eqref{oh2} holds.
\end{coro}
%%%%%%%%%%%%%%%%%%%%%%%%%%%%%%%%%%%%%%%%%%%%%%%%%%%%%%%%%%%%%%%%%%%%%%%%%% 

%%%%%%%%%%%%%%%%%%%%%%%%%%%%%%%%%%%%%%%%%%%%%%%%%%%%%%%%%%%%%%%%%%%%%%%%%% 
\prova
For any $\e\in(-\e_0,\e_0)$, set
\[
s'_0(\e):=\inf\{s_2\in(0,s) \, :\, 2D_0(s_2,\al,\om) < |\e|^{-1}\}.
\]
Then the assertion follows from Proposition \ref{bellaprop}, with $\Phi_2=2D_0(s_2,\al,\om)$.
\EP
%%%%%%%%%%%%%%%%%%%%%%%%%%%%%%%%%%%%%%%%%%%%%%%%%%%%%%%%%%%%%%%%%%%%%%%%%% 

This concludes the proof of Theorem \ref{moser}.

%%%%%%%%%%%%%%%%%%%%%%%%%%%%%%%%%%%%%%%%%%%%%%%%%%%%%%%%%%%%%%%%%%%%%%%%%% 
%%%%%%%%%%%%%%%%%%%%%%%%%%%%%%%%%%%%%%%%%%%%%%%%%%%%%%%%%%%%%%%%%%%%%%%%%% 
\section{Expansion in $\boldsymbol{j}$ of the counterterm $\boldsymbol{\h_j}$}
%\addcontentsline{toc}{subsection}{\thesection.\arabic{subsection}}{Expansion in ${j}$ of ${\h_j}$}
\label{labifsec}
\zerarcounters
%%%%%%%%%%%%%%%%%%%%%%%%%%%%%%%%%%%%%%%%%%%%%%%%%%%%%%%%%%%%%%%%%%%%%%%%%% 
%%%%%%%%%%%%%%%%%%%%%%%%%%%%%%%%%%%%%%%%%%%%%%%%%%%%%%%%%%%%%%%%%%%%%%%%%% 

%We proved the convergence of the series for $U(x,\f;c,\om,\e)$ and the counterterm $\h(c,\om,\e)$
%for all $\om\in\gotB$.
Let us now study more in detail the expression of the counterterm $\eta$. To this aim, we fix $N\ge1$ 
and consider, for any $\ze=(\ka,\x)\in\matK_N$, frequencies $\om(\ze)$ of the form \eqref{espome}:
our goal is to prove that, in such a case, $\eta=\h(c,\om(\ze),\e)$ admits 
an expansion of the form \eqref{nemmenounnome!}. 

%%%%%%%%%%%%%%%%%%%%%%%%%%%%%%%%%%%%%%%%%%%%%%%%%%%%%%%%%%%%%%%%%%%%%%%%%%

%%%%%%%%%%%%%%%%%%%%%%%%%%%%%%%%%%%%%%%%%%%%%%%%%%%%%%%%%%%%%%%%%%%%%%%%%%
%%%%%%%%%%%%%%%%%%%%%%%%%%%%%%%%%%%%%%%%%%%%%%%%%%%%%%%%%%%%%%%%%%%%%%%%%%
\subsection*{\ref{labifsec}.1\hspace{0.5cm}Expanded and fully renormalized $\boldsymbol{\h}$-trees}
\addcontentsline{toc}{subsection}{\ref{labifsec}.1\hspace{0.4cm}Expanded and fully renormalized ${\h}$-trees}
\label{etatrees}
%%%%%%%%%%%%%%%%%%%%%%%%%%%%%%%%%%%%%%%%%%%%%%%%%%%%%%%%%%%%%%%%%%%%%%%%%%
%%%%%%%%%%%%%%%%%%%%%%%%%%%%%%%%%%%%%%%%%%%%%%%%%%%%%%%%%%%%%%%%%%%%%%%%%%

To exploit of the form \eqref{espome} of the frequency vector and show that $\h_j(c,\om(\ze),\e)$ is of the form \eqref{nemmenounnome!}, we need to
Taylor-expand of the propagators in inverse powers of $j$. This leads to modify the rules for the graphic representation,
and hence to introduce two furter sets of trees.

%%%%%%%%%%%%%%%%%%%%%%%%%%%%%%%%%%%%%%%%%%%%%%%%%%%%%%%%%%%%%%%%%%%%%%%%%% 
\begin{defi}[\textbf{Special leaf and special path}] \label{special path}
Given any  $\vartheta\in \gotT^{(k)}_{ j ,\gote_j,\s}\cup\gotF^{(k)}_{ j ,\gote_j,\s}$, with $k\ge 1$, $j\in\ZZZ$ and $\s\in\{\pm\}$, set
\begin{equation}\label{diam}
\delta_{\vartheta} := \min\{|\calP(\ell_\la,\ell_\vartheta)|\ :\ \la\in\Lambda^*_{j,\s}(\vartheta)\},
\end{equation}
where the set $\Lambda^*_{j,\s}(\vartheta)$ is defined in \eqref{Lambdastar},
and let $\la_{\vartheta}\in \Lambda^*_{j,\s}(\vartheta)$ denote the uppermost leaf where the minimum in \eqref{diam} is attained.
Then set $\calP_{\vartheta}:=\calP(\ell_{\la_{\vartheta}},\ell_\vartheta)$.
We call $\la_{\vartheta}$
and $\calP_{\vartheta}$ the \emph{special leaf} and the \emph{special path}, respectively, of the tree $\vartheta$.
\end{defi}
%%%%%%%%%%%%%%%%%%%%%%%%%%%%%%%%%%%%%%%%%%%%%%%%%%%%%%%%%%%%%%%%%%%%%%%%%% 

%%%%%%%%%%%%%%%%%%%%%%%%%%%%%%%%%%%%%%%%%%%%%%%%%%%%%%%%%%%%%%%%%%%%%%%%%% 
\begin{rmk} \label{ovvissimo}
\emph{%remark
Recall that  for any $\vartheta\in \gotT^{(k)}_{ j ,\gote_j,\s}$ one has $\Lambda^*_{j,\s}(\vartheta)\ne\emptyset$ by Lemma \ref{cijei};
by construction the same happens for all $\vartheta\in  \gotF^{(k)}_{ j ,\gote_j,\s}$.
Hence the definition above is well-posed. 
}%remark
\end{rmk}
%%%%%%%%%%%%%%%%%%%%%%%%%%%%%%%%%%%%%%%%%%%%%%%%%%%%%%%%%%%%%%%%%%%%%%%%%% 

%%%%%%%%%%%%%%%%%%%%%%%%%%%%%%%%%%%%%%%%%%%%%%%%%%%%%%%%%%%%%%%%%%%%%%%%%% 
\begin{rmk}\label{itu}
\emph{
Since, if $\vartheta'$ is a subtree of $\vartheta\in \gotT^{(k)}_{ j ,\gote_j,\s}$ whose
root line is an $\h$-line, then $\vartheta' \in \gotT^{(k')}_{j',\gote_{j'},\s'}$
for some $k',j',\s'$, with $k'\le k-1$, we can define $\de_{\vartheta'}$, $\la_{\vartheta'}$ and $\calP_{\vartheta'}$ as done for $\vartheta$.
A similar argument applies if $\vartheta\in\gotF^{(k)}_{ j ,\gote_j,\s}$ and 
the root line of a subtree $\vartheta'$ of $\vartheta$  is an $\h$-line.
}
\end{rmk}
%%%%%%%%%%%%%%%%%%%%%%%%%%%%%%%%%%%%%%%%%%%%%%%%%%%%%%%%%%%%%%%%%%%%%%%%%% 

%%%%%%%%%%%%%%%%%%%%%%%%%%%%%%%%%%%%%%%%%%%%%%%%%%%%%%%%%%%%%%%%%%%%%%%%%%
\begin{rmk}\label{potato}
\emph{
By construction $\calP_\vartheta = \calP_{\breve\vartheta}$, where $\breve\vartheta$ is the trimmed tree associated with $\vartheta$.
}
\end{rmk}
%%%%%%%%%%%%%%%%%%%%%%%%%%%%%%%%%%%%%%%%%%%%%%%%%%%%%%%%%%%%%%%%%%%%%%%%%% 

%%%%%%%%%%%%%%%%%%%%%%%%%%%%%%%%%%%%%%%%%%%%%%%%%%%%%%%%%%%%%%%%%%%%%%%%%% 
\begin{defi}[\textbf{Ramification}]
\label{ramification}
Given an expanded tree
$$ \vartheta\in \Biggl(\,\bigcup_{k\ge1} \gotT^{(k)}_{j,\gote_{j},\s}\Biggr)\bigcup \Biggl(\,\bigcup_{k\ge1} \gotF^{(k)}_{j,\gote_{j},\s} \Biggr), $$
define recursively  the sets of paths $\calmP_k(\vartheta)$, with $k\ge 0$, as follows:
\begin{enumerate}[topsep=0ex]
\itemsep0em
\item set $\calmP_0(\vartheta)= \{\calP_{\vartheta}\}$,
\item for $k\ge 1$, if $\calmP_{k-1}(\vartheta) \neq \emptyset$, and there are $p$ subtrees
$\vartheta_1,\ldots,\vartheta_p \subset\vartheta$ whose root line is an $\h$-line entering a node in a path of $\calmP_{k-1}(\vartheta)$,
set $\calmP_k(\vartheta) :=\{ \calP_{\vartheta_1}, \ldots, \calP_{\vartheta_p}\}$, otherwise set $\calmP_k(\vartheta)= \emptyset$.
\end{enumerate}
If $d=d(\vartheta)=\max\{ k : \calmP_k(\vartheta) \neq \emptyset\}$,
we call $\calmP(\vartheta)=\calmP_0(\vartheta) \cup \ldots \cup \calmP_d(\vartheta)$ the \emph{ramification} of $\vartheta$,
and $d$ the \emph{length} of the ramifcation.
\end{defi}
%%%%%%%%%%%%%%%%%%%%%%%%%%%%%%%%%%%%%%%%%%%%%%%%%%%%%%%%%%%%%%%%%%%%%%%%%% 

%%%%%%%%%%%%%%%%%%%%%%%%%%%%%%%%%%%%%%%%%%%%%%%%%%%%%%%%%%%%%%%%%%%%%%%%%% 
\begin{rmk} \label{nonsisamai}
\emph{
By construction, both notions of special path and ramification make sense only for $\vartheta$ with $\nu_{\ell_\vartheta}=\gote_{j_\vartheta}$,
that is for $\vartheta$ whose value contribute to the counterterms.
}
\end{rmk}
%%%%%%%%%%%%%%%%%%%%%%%%%%%%%%%%%%%%%%%%%%%%%%%%%%%%%%%%%%%%%%%%%%%%%%%%%% 

%%%%%%%%%%%%%%%%%%%%%%%%%%%%%%%%%%%%%%%%%%%%%%%%%%%%%%%%%%%%%%%%%%%%%%%%%% 
\begin{defi}[\textbf{Set $\boldsymbol{\gotH}$ of the expanded $\boldsymbol{\h}$-trees}] 
%\hspace{-.27cm}\boldsymbol{\fR}$}]
\label{cihoprovato}
Let $\gotH$ be the set of all trees which are obtained from any expanded tree
%$\vartheta \in\gotF^b$ 
\[
\vartheta \in \bigcup_{\substack{k\ge1 \\ j\in\ZZZ \\ \s=\pm1}}\gotT^{(k)}_{j,\gote_j,\s} 
\]
by assigning to each line $\ell\in L(\vartheta)$ two further labels, the \emph{activity label} $h_{\ell}\in\{ {\bf a},{\bf n},\bf{s}\}$
and the \emph{power label} $\pow_{\ell}\in\{0,1,\ldots,N\}$, through the following procedure:
\begin{enumerate}[topsep=0ex]
\itemsep0em
\item\label{item1} for any line $\ell\in\calmP(\vartheta)$, if $\vartheta'\subseteq \vartheta$ denotes the subtree of $\vartheta$ such that $\ell\in\calP_{\vartheta'}$,
we set
\begin{itemize}[topsep=0ex]
\item[1.1.] $h_\ell=\bf{a}$ if $(j_\ell , \s_\ell) \neq (j_{\ell_{\vartheta'}},\s_{\ell_{\vartheta'}})$ and $n_\ell=0$,
\item[1.2.] $h_\ell=\bf{s}$ if $(j_\ell , \s_\ell) \neq (j_{\ell_{\vartheta'}},\s_{\ell_{\vartheta'}})$ and $n_\ell\ge1$,
\item[1.3.] $h_\ell=\bf{n}$ if $(j_\ell , \s_\ell) = (j_{\ell_{\vartheta'}},\s_{\ell_{\vartheta'}})${;}
\end{itemize}
\item for any line $\ell\in L(\vartheta)\setminus\calmP(\vartheta)$ we set $h_\ell=\bf{n}${;}
\item we say that a line $\ell\in L(\vartheta)$ is \emph{active} if $h_\ell=\bf{a}$,
\emph{neutral} if $h_\ell=\bf{n}$ and \emph{small} if $h_{\ell}=\bf{s}$, and call $L_{\bf{a}}(\vartheta)$ the set of active lines,
$L_{\bf{n}}(\vartheta)$ the set of neutral lines and $L_{\bf{s}}(\vartheta)$ the set of small lines;
%\item if $|L_{\bf a}(\vartheta)| > N$, we decompose $L_{\bf a}(\vartheta) = L^1_{\bf a}(\vartheta) \sqcup L^2_{\bf a}(\vartheta)$,
%where $L^1_{\bf a}(\vartheta)$ contains $N$ active lines, arbitrarily chosen, while,
%if $|L_{\bf a}(\vartheta)| \le N$, we define $L^1_{\bf a}(\vartheta) := L_{\bf a}(\vartheta)$ and $L^2_{\bf a}(\vartheta):=\emptyset$;
\item with any line $\ell\in L(\vartheta)$ we associate a label $\pow_\ell\in \NNN\cup\{0\}$,
the \emph{activity degree label},
with the constraint that $\pow_\ell=0$ if $h_\ell=\bf{n},\bf{s}$ and $1 \le \pow_\ell \le M(\vartheta)$ if $h_\ell=\bf{a}$, where
$$ M(\vartheta):= %N - L_{\bf a}^1(\vartheta)+1= 
\max\{ N - |L_{\bf a}(\vartheta)|,0 \} + 1.$$
\end{enumerate}
We call expanded \emph{$\eta$-tree} any tree $\vartheta\in\gotH$,
and  $\gotH^{(k)}_{j,\s}$ the set of the expanded $\eta$-trees in $\gotH$
with order $k$ such that the root line has component $j$ and sign $\s$.
Finally we call ${\gotH}^{(k)}_{j,\s}(\pow)$, for $\pow\in\{0,\ldots,N-1\}$, the set of the expanded $\eta$-trees $\vartheta\in{\gotH}^{(k)}_{j,\s}$ such that 
\begin{equation} \label{sommatau}
\sum_{\ell\in L(\vartheta)}\pow_\ell =\pow ,
\end{equation}
and  ${\gotH}^{(k)}_{j,\s}(N)$ the set of the expanded $\eta$-trees $\vartheta\in{\gotH}^{(k)}_{j,\s}$  such that 
\begin{equation} \label{sommaN}
\sum_{\ell\in L(\vartheta)}\pow_\ell \ge N.
\end{equation}
\end{defi}
%%%%%%%%%%%%%%%%%%%%%%%%%%%%%%%%%%%%%%%%%%%%%%%%%%%%%%%%%%%%%%%%%%%%%%%%%% 

%%%%%%%%%%%%%%%%%%%%%%%%%%%%%%%%%%%%%%%%%%%%%%%%%%%%%%%%%%%%%%%%%%%%%%%%%% 
\begin{defi}[\textbf{Set $\boldsymbol{{\gotF\gotH}}$ of the fully renormalized $\boldsymbol{\h}$-trees}] \label{purequesto}
Let ${\gotF\gotH}$ denote the set of all trees which are obtained as in Definition \ref{cihoprovato} but 
\begin{enumerate}[topsep=0ex]
\itemsep0em
\item starting with trees in the set $\gotF$ of the fully renormalized trees instead of the set $\gotT$ of the expanded trees,
%\item substituting $x_\ell(\und{y}_\ell)$ in place of $x_\ell$ in item \ref{item1}.
\item adding to any line $\ell\in L(\vartheta)$ a further label $p_{\ell}\in\{0,\ldots,\del_\ell\}$, 
with the constraint that $p_\ell=0$ if $h_\ell=\bf{n},\bf{s}$ and $0 \le p_\ell \le \del_{\ell}$ if $h_\ell=\bf{a}$.
\end{enumerate}
We call \emph{fully renormalized $\eta$-tree} any tree $\vartheta\in\gotF\gotH$,
and  $\gotF\gotH^{(k)}_{j,\s}$ the set of the fully renormalized $\eta$-trees in $\gotF\gotH$ with order
$k$ such that the root line has component $j$ and sign $\s$.
For $\pow\in\{0,\ldots,N-1\}$, we
call ${\gotF\gotH}^{(k)}_{j,\s}(\pow)$ the set of the fully renormalized $\eta$-trees $\vartheta\in{\gotF\gotH}^{(k)}_{j,\s}$ such that \eqref{sommatau} holds,
%
%\begin{equation} \label{sommatauFK}
%\sum_{\ell\in L(\vartheta)}\pow_\ell =\pow ,
%\end{equation}
%
while
${\gotF\gotH}^{(k)}_{j,\s}(N)$ denotes the set of the fully renormalized $\eta$-trees $\vartheta\in{\gotF\gotH}^{(k)}_{j,\s}$  such that \eqref{sommaN} holds.
%
%\begin{equation} \label{sommaNFK}
%\sum_{\ell\in L(\vartheta)}\pow_\ell \ge N.
%\end{equation}
\end{defi}
%%%%%%%%%%%%%%%%%%%%%%%%%%%%%%%%%%%%%%%%%%%%%%%%%%%%%%%%%%%%%%%%%%%%%%%%%% 

%%%%%%%%%%%%%%%%%%%%%%%%%%%%%%%%%%%%%%%%%%%%%%%%%%%%%%%%%%%%%%%%%%%%%%%%%%
\begin{defi}[\textbf{$\boldsymbol{\h}$-tree}] \label{hT}
We call $\h$-\emph{tree}  any element of $\gotH\cup\gotF\gotH$, that is any either expanded or fully renormalized
$\h$-tree.
\end{defi}
%%%%%%%%%%%%%%%%%%%%%%%%%%%%%%%%%%%%%%%%%%%%%%%%%%%%%%%%%%%%%%%%%%%%%%%%%%

%%%%%%%%%%%%%%%%%%%%%%%%%%%%%%%%%%%%%%%%%%%%%%%%%%%%%%%%%%%%%%%%%%%%%%%%%% 
\begin{rmk}\label{contiamole}
\emph{
If $\vartheta\in {\gotH}^{(k)}_{j,\s}(q)\cup{\gotF\gotH}^{(k)}_{j,\s}(q)$ with $q\le N-1$, then $q_{\ell} < M(\vartheta)$ for all $\ell\in L(\vartheta)$.
Indeed, if $\vartheta$ is such that there is $\bar{\ell}\in L_{\bf{a}}(\vartheta)$ with $q_{\bar{\ell}}=M(\vartheta)$ then
\[
\sum_{\ell\in L(\vartheta)}q_\ell = M(\vartheta) + \sum_{\substack{\ell\in L(\vartheta) \\ \ell\ne \bar{\ell}}}q_\ell \ge N,
\]
and hence $\vartheta\in {\gotH}^{(k)}_{j,\s}(N)\cup{\gotF\gotH}^{(k)}_{j,\s}(N)$.
}
\end{rmk}
%%%%%%%%%%%%%%%%%%%%%%%%%%%%%%%%%%%%%%%%%%%%%%%%%%%%%%%%%%%%%%%%%%%%%%%%%% 

 %%%%%%%%%%%%%%%%%%%%%%%%%%%%%%%%%%%%%%%%%%%%%%%%%%%%%%%%%%%%%%%%%%%%%%%%%%
\begin{defi}[\textbf{Trimmed tree associated with an $\boldsymbol{\h}$-tree}] 
\label{trih}
{
%For all $k$, $\nu$, $j$, $\s$
For  any $\vartheta\in\gotH \cup \gotF\gotH$, set
\[
N_2(\vartheta) :=\{ v \in N(\vartheta) :s_{v}=2\}.% \hbox{ and } v \st{\preceq} r_{\vartheta}\}  . 
\]
For any $v\in N_2(\vartheta)$, if $\vartheta_v$ denotes the subtree of $\vartheta$ entering $v$
whose root line is an $\h$-line, set 
\[
F_2^*(\vartheta) := \{ \vartheta_v : v \in N_2(\vartheta) \hbox{ and } v \st{\preceq} r_{\vartheta}\}.
\]
Let $\breve\vartheta$ be the  tree 
obtained from $\vartheta$ by cutting the subtrees in $F_2^*(\vartheta)$;
we call $\breve\vartheta$ the \emph{trimmed $\h$-tree associated with the $\h$-tree} $\vartheta$.}
\end{defi}
%%%%%%%%%%%%%%%%%%%%%%%%%%%%%%%%%%%%%%%%%%%%%%%%%%%%%%%%%%%%%%%%%%%%%%%%%%

Define
\begin{equation} \label{calJh}
\calJ(\vartheta):= J(\breve\vartheta) + \sum_{v\in N_2(\vartheta)}J(\breve\vartheta_v) , \qquad
\calJ_{\al}(\vartheta):= (\calJ(\vartheta))^{1/\al} .
\end{equation}

%%%%%%%%%%%%%%%%%%%%%%%%%%%%%%%%%%%%%%%%%%%%%%%%%%%%%%%%%%%%%%%%%%%%%%%%%% 
\begin{rmk} \label{unicoP}
\emph{
If $\ell\in L_{\bf{a}}(\vartheta)$ then either $\ell \in \calP_{\breve\vartheta}$ or, 
using the notation introduced in Definition \ref{trih}, there are 
$v_1,\ldots,v_p\in N_2(\vartheta)$,
with $p\ge1$, such that % $\vartheta'\subseteq\vartheta$ %(possibly $\vartheta'=\vartheta$) 
\begin{itemize}[topsep=0ex]
\itemsep0em
\item $\ell\in\calP_{\breve\vartheta_{1}}$,
\item $\ell_{{i}}\in\calP_{\breve\vartheta_{{i+1}}}$ for $i=1,\ldots,p-1$,
\item $\ell_{p} \in\calP_{\breve\vartheta}$,
\end{itemize}
where we have shorthened $\ell_i=\ell_{v_i}$ and $\vartheta_i=\vartheta_{v_i}$ for $i=1,\ldots,p$.
}
\end{rmk}
%%%%%%%%%%%%%%%%%%%%%%%%%%%%%%%%%%%%%%%%%%%%%%%%%%%%%%%%%%%%%%%%%%%%%%%%%% 

%%%%%%%%%%%%%%%%%%%%%%%%%%%%%%%%%%%%%%%%%%%%%%%%%%%%%%%%%%%%%%%%%%%%%%%%%% 
%\begin{rmk}
%\emph{%remark
%Of course for any $\vartheta\in\gotH$ we can define $\Val(\vartheta;c,\om)$ as in \eqref{valeta}, simply ignoring the extra labels $h_\ell,\pow_\ell$.
%Similarly for any $\vartheta\in\gotF\gotH$ we can define $\Val^\ttR(\vartheta;c,\om)$ as in \eqref{valR}.
%simply ignoring the extra labels $h_\ell,\pow_\ell$.
%In other words we have
%\begin{equation}\label{rietabis}
%\h_j^{(k)}(c,\om)= -
%\sum_{\vartheta \in {\gotH}^{(k)}_{j,\gote_j,+}}  \frac{1}{c_j}\Val (\vartheta;c,\om)=
%-\sum_{\vartheta \in {\gotF\gotH}^{(k)}_{j,\gote_j,+}}  \frac{1}{c_j}\Val^{\ttR} (\vartheta;c,\om),
%\end{equation}
%where the sums are absolutely convergent, reasoning as in the proof of Lemma \ref{welldefined} for the first sum,
%and Proposition \ref{voltabuona} for the second sum.
%}%remark
%\end{rmk}
%%%%%%%%%%%%%%%%%%%%%%%%%%%%%%%%%%%%%%%%%%%%%%%%%%%%%%%%%%%%%%%%%%%%%%%%%% 

%%%%%%%%%%%%%%%%%%%%%%%%%%%%%%%%%%%%%%%%%%%%%%%%%%%%%%%%%%%%%%%%%%%%%%%%%% 
\begin{rmk}\label{secattivo}
\emph{
For any $\vartheta\in\gotH \cup \gotF\gotH$ such that $\Val^\ttR(\vartheta;c,\om)\ne0$, if 
$\ell\in L_{\bf{a}}(\vartheta)$ is the external line of a RC $T$, then $\nmax_T=-1$ and hence $\matO_T=\matL$.
In particular, recalling that the trees in $\gotF\gotH$ are obtained from trees in $\gotF$ according to Definition \ref{terribile},
in a fully renormalized $\h$-tree no active line can be resonant 
since no chain $\gotC=\{T_1,\ldots,T_p\}$ can have $\matO_{T_i}=\matO_{T_{i+1}}=\LL$ for some $i=1,\ldots,p-1$.
}
\end{rmk}
%%%%%%%%%%%%%%%%%%%%%%%%%%%%%%%%%%%%%%%%%%%%%%%%%%%%%%%%%%%%%%%%%%%%%%%%%% 

%%%%%%%%%%%%%%%%%%%%%%%%%%%%%%%%%%%%%%%%%%%%%%%%%%%%%%%%%%%%%%%%%%%%%%%%%% 
%%%%%%%%%%%%%%%%%%%%%%%%%%%%%%%%%%%%%%%%%%%%%%%%%%%%%%%%%%%%%%%%%%%%%%%%%%
\subsection*{\ref{labifsec}.2\hspace{0.5cm}Values of the expanded and fully renormalized $\boldsymbol{\h}$-trees}}
\addcontentsline{toc}{subsection}{\ref{labifsec}.2\hspace{0.4cm}Values of the expanded and fully renormalized ${\h}$-trees}
\label{etatreevalues}
%%%%%%%%%%%%%%%%%%%%%%%%%%%%%%%%%%%%%%%%%%%%%%%%%%%%%%%%%%%%%%%%%%%%%%%%%%
%%%%%%%%%%%%%%%%%%%%%%%%%%%%%%%%%%%%%%%%%%%%%%%%%%%%%%%%%%%%%%%%%%%%%%%%%%

\noindent
Given a tree $\vartheta\in\gotF\gotH^{(k)}_{j,\s}$ with $j\ne0$ and a line $\ell \in \calmP(\vartheta)$, 
if $\vartheta'\subseteq\vartheta$ is the subtree such that $\ell\in\calP_{\vartheta'}$, we set
\begin{subequations} \label{nuj0}
\begin{align}
\nu_\ell^\flat(\und{y}_\ell) &  := 
\s_\ell\nu_\ell(\und{y}_\ell) - \s_{\ell_{\vartheta'}}\gote_{j_{\ell_{\vartheta'}}}  ,
\label{nuj0a} \\
j_\ell^\flat &: = \sum_{\substack{ \la\in \Lambda(\vartheta') \\ \la\st{\preceq}\ell}} \s_\la {j_\la} - \s j = \s_\ell j_{\ell} - \s j .
\label{nuj0b} \\
j_{\ell_{\vartheta'}}^\flat &: = \sum_{ \la\in \Lambda^*(\vartheta') } \s_\la {j_\la} - \s j = \s_{\ell_{\vartheta'}} j_{\ell_{\vartheta'}} - \s j .
 \label{nuj0c}
\end{align}
\end{subequations}
We write the renormalized small divisor of any line $\ell \in \calmP(\vartheta)$ in terms of \eqref{nuj0} as
\begin{equation}\label{ambrogio}
\begin{aligned}
x_{\ell}(\und{y}_\ell)  &= \s_{\ell} \, \om\cdot\nu_\ell^\flat(\und{y}_\ell) + \s_\ell \s_{\ell_{\vartheta'}} \om_{j_{\ell_{\vartheta'}}}
- \om_{\s_\ell j_\ell^\flat + \s_\ell \s j} \\
& = \s_{\ell} \, \om\cdot\nu_\ell^\flat(\und{y}_\ell) + \s_\ell \s_{\ell_{\vartheta'}}
\om_{\s_{\ell_{\vartheta'}} j^\flat_{\ell_{\vartheta'}} + \s_{\ell_{\vartheta'}} \s j}  - \om_{\s_\ell j_\ell^\flat + \s_\ell \s j} .
\end{aligned}
\end{equation}
%where the last line follows by writing also
%%
%\begin{equation} \label{nuj0c}
%j_{\ell_{\vartheta'}}^\flat : = \sum_{ \la\in \Lambda^*(\vartheta') } \s_\la {j_\la} - \s j = \s_{\ell_{\vartheta'}} j_{\ell_{\vartheta'}} - \s j ,
%\end{equation}
%%
%in analogy with \eqref{nuj0b}.
Then both $\om_{\s_{\ell_{\vartheta'}} j_{\ell^\flat_{\vartheta'}} + \s_{\ell_{\vartheta'}} \s j}$ %$\om_{j_{\ell_{\vartheta'}}}$ 
and $\om_{\s_\ell j_\ell^\flat + \s_\ell \s j}$ can be expanded according to \eqref{espome}.
If we set
\[
A_i(p):=\frac{1}{p!}\frac{d^p}{dz^p}\left.\frac{1}{(1+z)^i}\right|_{z=0},
\qquad i\ge1, \quad p \ge 0, 
\]
we obtain, for any $M=1,\ldots,N$
\begin{equation}\label{icselle}
x_{\ell}(\und{y}_\ell) = \bigl(\s_\ell \s_{\ell_{\vartheta'}} -1 \bigr) j^2 + 2 \bigl( \s_\ell \s_{\ell_{\vartheta'}} j_{\ell_{\vartheta'}}^\flat  - j_\ell^\flat \bigr) \s j 
+ \sum_{q=0}^{M-1} \frac{C_q(\ze,\ell,j)}{j^{q}} + {r_M(\ze,\ell,j)},
 \end{equation}
with\footnote{Here we are using the convention $0^0=1$.}
\begin{subequations}\label{gliC}
\begin{align}
%\hspace{-.5cm}
\null \!\!\!\!\!\!
&C_0{(\ze,\ell,j)}:=  \s_\ell \, \omega \cdot\nu_\ell^\flat(\und{y}_\ell) 
+ \bigl(\s_\ell \s_{\ell_{\vartheta'}}(j_{\ell_{\vartheta'}}^\flat)^2 - (j_\ell^\flat)^2 \bigr) + (\s_\ell \s_{\ell_{\vartheta'}} -1 )\kappa_0 ,
\phantom{\sum_{h=2}^q}\\
&C_1(\ze,\ell,j)=0 , \\
%\hspace{-.5cm}
\null \!\!\!\!\!\!
&C_q{(\ze,\ell,j)} := 
\sum_{h=2}^q  \s^q \kappa_h \, A_{h}(q-h) \Bigl(  \s_\ell\sigma_{\ell_{\vartheta'}} 
(j_{\ell_{\vartheta'}}^\flat)^{q-h}\s_{\ell_{\vartheta'}}^h 
%\\ &\qquad\qquad\qquad
- (j_\ell^\flat)^{q-h} \s_\ell^h \Bigr) , \quad 2 \le q \le N-1,
\label{ciqu}
\end{align}
\end{subequations}
and\footnote{Recall that we are using the convention that the sum over the empty set is $0$.}
\begin{equation}\label{restino}
\begin{aligned}
\null\!\!
r_M(\ze,\ell,j)  & := \sum_{q=2}^{N-1}\Big(\frac{\s_\ell \s_{\ell_{\vartheta'}}\ka_q}{(\s_{\ell_{\vartheta'}} j^\flat_{\ell_{\vartheta'}} + \s_{\ell_{\vartheta'}} \s j)^q} -
\frac{\ka_q}{(\s_\ell j_\ell^\flat + \s_\ell \s j)^q}\Big) \\
\null\!\! & +  
\s_\ell \s_{\ell_{\vartheta'}}\xi_{\s_{\ell_{\vartheta'}} j^\flat_{\ell_{\vartheta'}} + \s_{\ell_{\vartheta'}} \s j} - \xi_{\s_\ell j_\ell^\flat + \s_\ell \s j} 
- \sum_{q=2}^{M-1} \frac{C_q(\ze,\ell,j)}{j^{q}} .
\end{aligned}
\end{equation}
Of course \eqref{ciqu} has to be disregarded if $N=1,2$.

Similarly, if $\vartheta\in\gotH^{(k)}_{j,\s}$ and $\ell \in \calmP(\vartheta)$ we may repeat the construction above by introducing
\begin{equation} \label{nuj0bis}
\nu_\ell^\flat  := 
\s_\ell\nu_\ell  - \s_{\ell_{\vartheta'}}\gote_{j_{\ell_{\vartheta'}}}  ,
\end{equation}
%
%$\nu_\ell^\flat$ as in \eqref{nuj0a} but  with $\nu_\ell$ replacing $\nu_\ell(\und{y}_\ell)$
and writing  $x_\ell$ as in \eqref{icselle}
with the only difference that $C_0(\ze,\ell,j)$ is given as in \eqref{gliC} with $\nu_\ell^\flat$ replacing $\nu_\ell^\flat(\und{y}_\ell)$.

%%%%%%%%%%%%%%%%%%%%%%%%%%%%%%%%%%%%%%%%%%%%%%%%%%%%%%%%%%%%%%%%%%%%%%%%%% 
\begin{lemma}\label{in}
For any $k\ge1$, $j\in\ZZZ\setminus\{0\}$ and $\s\in\{\pm\}$, and for any
$\vartheta \in\gotH_{j,\s}^{(k)}\cup \gotF\gotH_{j,\s}^{(k)}$ and any $\ell\in \calmP(\vartheta)$, one has
\begin{equation}\label{stoclaim}
|\nu_\ell^\flat|_2 \le (\calJ_{\al}(\vartheta))^2 , \qquad |j_\ell^\flat| \le  \calJ_{\al} (\vartheta) , \qquad
|j^\flat_{\ell_{\vartheta'}}| <  \calJ_{\al} (\vartheta) , \qquad
|\om\cdot\nu_{\ell}^\flat(\und{y}_\ell)|\le 2 |\nu_\ell^\flat|_2 ,
\end{equation}
where $\vartheta'\subseteq\vartheta$ is the subtree such that $\ell\in\calP_{\vartheta'}$.
\end{lemma}
%%%%%%%%%%%%%%%%%%%%%%%%%%%%%%%%%%%%%%%%%%%%%%%%%%%%%%%%%%%%%%%%%%%%%%%%%% 

%%%%%%%%%%%%%%%%%%%%%%%%%%%%%%%%%%%%%%%%%%%%%%%%%%%%%%%%%%%%%%%%%%%%%%%%%% 
\prova
We start proving the first three bounds.
By Remark \ref{unicoP} we can distinguish between two cases.

\begin{itemize}[topsep=0ex]
\itemsep0em

\item
If $\ell\in\calP_{\breve\vartheta}$, so that $\vartheta'=\vartheta$, then one has
\[
\nu_\ell^\flat=\sum_{\substack{\la\st{\preceq}\ell \\ \la\ne\la_{\vartheta}}}\s_\la\gote_{j_\la},\qquad 
j_\ell^\flat = \sum_{\substack{\la\st{\preceq}\ell \\ \la\ne\la_{\vartheta}}}\s_\la j_\la,
\]
and hence
\[
|\nu^\flat_{\ell}|_2\le \sum_{\substack{\la\st{\preceq}\ell \\ \la\ne\la_{\vartheta}}} \jap{j_\la}^2 \le\Biggl(
\sum_{\substack{\la\st{\preceq}\ell \\ \la\ne\la_{\vartheta}}}\jap{j_\la}^\al
\Biggr)^{2/\al}\le (J(\breve\vartheta))^{2/\al},
\]
and
\[
|j_\ell^\flat| \le \sum_{\substack{\la\st{\preceq}\ell \\ \la\ne\la_{\vartheta}}}|j_\la| \le
\Biggl(
\sum_{\substack{\la\st{\preceq}\ell \\ \la\ne\la_{\vartheta}}}|j_\la|^\al
\Biggr)^{1/\al}\le (J(\breve\vartheta))^{1/\al} ,
\]
while the third bound holds trivially since $j^\flat_{\ell_\vartheta}=0$.

\item
Otherwise, with the notations of
%If there are $v_1,\ldots,v_p\in N_2(\vartheta)$, with $p\ge1$, such that 1--3 of 
Remark \ref{unicoP}, one has $\vartheta'=\vartheta_{1}$ and %, if $v_0$ denotes the node $\ell$ exits, so that $\ell=\ell_{v_0}$, 
\[
\nu_\ell^\flat=\sum_{\substack{\la\st{\preceq}\ell \\ \la\ne\la_{\vartheta_{1}}}}\s_\la\gote_{j_\la},\qquad
j_{\ell}^\flat =  
\sum_{\substack{\la\st{\preceq}\ell_{{p}} \\ \la\ne\la_{\vartheta}}}\s_\la j_\la +
\sum_{i=0}^{p-1} \sum_{\substack{\la\st{\preceq}\ell_{{i}} \\ \la\ne\la_{\vartheta_{{i+1}}}}}\s_\la j_\la ,\qquad
j^\flat_{\vartheta'}=j^\flat_{\ell_{1}},
\]
with $\ell_0=\ell$,
and hence
\[
|\nu^\flat_{\ell}|_2\le \sum_{\substack{\la\st{\preceq}\ell \\ \la\ne\la_{\vartheta_v}}} \jap{j_\la}^2 \le\Biggl(
\sum_{\substack{\la\st{\preceq}\ell \\ \la\ne\la_{\vartheta_v}}}\jap{j_\la}^\al
\Biggr)^{2/\al}\le (J(\breve\vartheta_v))^{2/\al},
\]
and
\[
\begin{aligned}
|j_{\ell}^\flat| &\le 
\sum_{\substack{\la\st{\preceq}\ell_{{p}} \\ \la\ne\la_{\vartheta}}}|j_\la| +
\sum_{i=0}^{p-1} \sum_{\substack{\la\st{\preceq}\ell_{{i}} \\ \la\ne\la_{\vartheta_{{i+1}}}}}| j_\la|
\\ & \le \Biggl(
\sum_{\substack{\la\st{\preceq}\ell_{{p}} \\ \la\ne\la_{\vartheta}}}|j_\la|^\al +
 \sum_{i=0}^{p-1} \sum_{\substack{\la\st{\preceq}\ell_{{i}} \\ \la\ne\la_{\vartheta_{{i+1}}}}}| j_\la|^\al \Biggr)^{1/\al} 
%\\ &
 \le
\Biggl( 
J(\breve\vartheta) + \sum_{v\in N_2(\vartheta)}J(\breve\vartheta_v) 
\Biggr)^{1/\al},
\end{aligned}
\]
and the same argument can be used for $j^\flat_{\vartheta'}$.

\end{itemize}
This concludes the proof of the first three bounds in \eqref{stoclaim}.

Finally, by Lemma \ref{stanca} and Remark \ref{serveunnome} one has
$$ |\om\cdot\nu_{\ell}^\flat(\und{y}_\ell)-\om\cdot\nu_{\ell}^\flat| =
|\om\cdot\nu_{\ell}(\und{y}_\ell)-\om\cdot\nu_{\ell}|\le \frac{3}{32}
$$
and hence the fourth of \eqref{stoclaim} follows as well.
\EP
%%%%%%%%%%%%%%%%%%%%%%%%%%%%%%%%%%%%%%%%%%%%%%%%%%%%%%%%%%%%%%%%%%%%%%%%%% 

%%%%%%%%%%%%%%%%%%%%%%%%%%%%%%%%%%%%%%%%%%%%%%%%%%%%%%%%%%%%%%%%%%%%%%%%%% 
\begin{lemma}\label{helsinki}
There exists a constant $K=K(\al,N)$ such that for any $k\ge1$, $j\in\ZZZ\setminus\{0\}$ and $\s\in\{\pm\}$, and
for any $\vartheta \in\gotH_{j,\s}^{(k)}\cup \gotF\gotH_{j,\s}^{(k)}$, any $\ell\in L_{\bf{a}}(\vartheta)$ and all $\ze\in\matW_N$ one has
\begin{subequations}
\begin{align}
|C_0(\ze,\ell,j)|	 & \le  K \, (\calJ_{\al}(\vartheta))^2 , 
%\Biggl( J(\breve\vartheta) + \sum_{v\in N_2(\vartheta)}J(\breve\vartheta_v) \Biggr)^{2/\alpha}\,, 
\label{helsinki0} \\
|C_q(\ze,\ell,j)|	 & \le  K  \, (\calJ_{\al}(\vartheta))^{q-2}, \qquad
q=2,\dots,N-1,
%\Biggl( J(\breve\vartheta) + \sum_{v\in N_2(\vartheta)}J(\breve\vartheta_v) \Biggr)^{(q-2)/\alpha}\,, 
\label{helsinki1} \\
|j^Mr_M(\ze,\ell,j)|& \le  K  \, (\calJ_{\al}(\vartheta))^{M}, \qquad 
M=1,\ldots, N,
% \Biggl( J(\breve\vartheta) + \sum_{v\in N_2(\vartheta)}J(\breve\vartheta_v) \Biggr)^{(M-1)/\alpha} ,
\label{helsinki2} 
\end{align}
\end{subequations}
where \eqref{helsinki1} has to be disregarded if $N=1,2$.
\end{lemma}
%%%%%%%%%%%%%%%%%%%%%%%%%%%%%%%%%%%%%%%%%%%%%%%%%%%%%%%%%%%%%%%%%%%%%%%%%% 

%%%%%%%%%%%%%%%%%%%%%%%%%%%%%%%%%%%%%%%%%%%%%%%%%%%%%%%%%%%%%%%%%%%%%%%%%% 
\prova
Let $\vartheta \in \gotF\gotH_{j,\s}^{(k)}$, $\ell\in L_{\bf{a}}(\vartheta)$ and $\vartheta'\subseteq\vartheta$ such that $\ell\in\calP_{\vartheta'}$.
Each coefficient $C_q(\ze,\ell,j)$  in \eqref{gliC} is a polynomial of degree  $\le q-2$ in
the variables $j_\ell^\flat$, $j_{\ell_{\vartheta'}}^\flat$ and $\omega\cdot\nu_\ell^\flat(\und{y}_\ell)$,
whose coefficients are uniformly bounded for $\ze\in\matW_N$.
Therefore by Lemma \ref{in} there is a constant $K_1=K_1(\al,N)$ such that the bound \eqref{helsinki0} and \eqref{helsinki1} hold
for any $K\ge K_1$.
In order to prove \eqref{helsinki2}, we distinguish between two cases.
\begin{itemize}[topsep=0ex]
\itemsep0em

\item If $|j| \le 2\calJ_{\al}(\vartheta)$,
%
%$$ |j| \le \left( J(\breve\vartheta) + \sum_{v\in N_2(\vartheta)}J(\breve\vartheta_v) \right)^{2/\al} , $$
%
then \eqref{restino}, together with \eqref{helsinki1}, gives
\[
| r_M(\ze,\ell,j) | \le N + \sum_{q=2}^{M-1} K_1 \frac{(\calJ_{\al}(\vartheta))^{q-2}}{j^q}  \le K_2 \frac{(\calJ_{\al}(\vartheta))^{M-1}}{|j|^{M-1}} ,
\]
for some constant $K_2 = K_2(\al,N)\ge K_1$.
\item If $|j| > 2\calJ_{\al}(\vartheta)$, from \eqref{restino} we obtain
\begin{equation}\label{restino2}
\begin{aligned}
& r_M(\ze,\ell,j) :=  \sum_{q=M}^{N-1}\Big(\frac{\ka_q}{(\s_{\ell_{\vartheta'}} j^\flat_{\ell_{\vartheta'}} + \s_{\ell_{\vartheta'}} \s j)^q} -
\frac{\ka_q}{(\s_\ell j_\ell^\flat + \s_\ell \s j)^q}\Big) \\
& \quad +  
\xi_{\s_{\ell_{\vartheta'}} j^\flat_{\ell_{\vartheta'}} + \s_{\ell_{\vartheta'}} \s j} - \xi_{\s_\ell j_\ell^\flat + \s_\ell \s j} \\
& \quad +\frac{\s^M}{j^M}\sum_{h=2}^{M-1}  \kappa_h \, A_h(M-h) \Bigl(  \s_\ell\sigma_{\ell_{\vartheta'}} 
\Big(\frac{j_{\ell_{\vartheta'}}^\flat}{1+z_1}\Big)^{M-h}\s_{\ell_{\vartheta'}}^h 
%\\ &\qquad\qquad\qquad
- \Big(\frac{j_\ell^\flat}{1+z_2}\Big)^{M-h} \s_\ell^h \Bigr) ,
%\Big[\frac{1}{(1+z_1)^{M-h}} - \frac{1}{(1+z_2)^{M-h}}\Big]
\end{aligned}
\end{equation}
for suitable $z_1$ and $z_2$ such that $|z_1| \le |j_{\ell_{\vartheta'}}^\flat/j|$ and $|z_2| \le |j_{\ell}^\flat/j|$. Therefore the last line of \eqref{restino2}
is bounded proportionally to $(\calJ_{\al}(\vartheta))^{M-2}/|j|^M$, whereas the first two lines are both bounded by a constant times $1/|j|^M$.
\end{itemize}
Collecting together the results, the bound \eqref{helsinki2} follows for suitable $K$.
\EP
%%%%%%%%%%%%%%%%%%%%%%%%%%%%%%%%%%%%%%%%%%%%%%%%%%%%%%%%%%%%%%%%%%%%%%%%%% 

If $\vartheta \in \gotF\gotH_{j,\s}^{(k)}$, for some $k\ge1$, $j\in\ZZZ\setminus\{0\}$ and $\s\in\{\pm\}$, $\ell\in L_{\bf a}(\vartheta)$
and $\vartheta'\subseteq\vartheta$ is such that $\ell\in\calP_{\vartheta'}$,
one has
%$|x_\ell(\und{y}_\ell)|\ge\be(m_0)$ and
either $(\s_\ell \s_{\ell_{\vartheta'}} -1 )\ne0$
or $( \s_\ell \s_{\ell_{\vartheta'}} j_{\ell_{\vartheta'}}^\flat  - j_\ell^\flat)\ne0$, and hence one can expand,
for any $M=1,\ldots,N$,
\begin{equation}\label{odiotuttidavvero}
\frac{1}{\bigl( x_\ell(\und{y}_\ell)\bigr)^{p_\ell+1}} = \sum_{\pow=1}^{M-1} %\frac{
\gotg_\ell(\ze,\pow,j) %}{j^\pow} 
+ {\gotG^{(M)}_\ell(\ze,j)}, 
\end{equation}
where, if $\s_\ell\ne\s_{\ell_{\vartheta'}}$ and $q=2(p_\ell+1),\ldots,N-1$,
one has\footnote{Here the sum over $q_1,\ldots,q_p$ are meant as 1 for $p=0$.}
\begin{equation} \label{mava}
\gotg_\ell(\ze,\pow,j):=
\frac{1}{j^q} \sum_{p=0}^{\io} \frac{A_{p_\ell+1}(p)}{(-2)^{p_\ell+1}}
\!\!\!\!\!\!\!\!\!
\sum_{ \substack{ q_1,\ldots,q_p \ge 1 \\ q_1+\ldots+q_p=\pow-2(p_\ell+1)}} \!\!\!\!\!\!\!\!\!
B_{q_1}(\ze,\ell,j) \ldots B_{q_p}(\ze,\ell,j) , 
\end{equation}
while, if $\s_\ell=\s_{\ell_{\vartheta'}}$ and $q=p_\ell+1,\ldots,N-1$, one has
\begin{equation} \label{mava1}
\gotg_\ell(\ze,\pow,j):=
\frac{1}{j^q} \sum_{p=0}^{\io} \frac{A_{p_\ell+1}(p)}{(2\s  (j_{\ell_{\vartheta'}}^\flat - j_\ell^\flat ))^{p_\ell+1}}\!\!\!\!\!\!
\!\!\!\!\!\!\!\!
\sum_{ \substack{ q_1,\ldots,q_p \ge 1 \\ q_1+\ldots+q_p=\pow-(p_\ell+1)}}\!\!\!\!\!\!\!\!\!\!\!\!\!\!
B'_{q_1}(\ze,\ell,j) \! \ldots \! B'_{q_p}(\ze,\ell,j), 
\end{equation}
with
\begin{equation} \label{B}
\begin{aligned}
\hspace{-.3cm}
& B_1(\ze,\ell,j) \!:=\! \s \bigl( j_\ell^\flat + j_{\ell_{\vartheta'}}^\flat \bigr) , \qquad
B_2(\ze,\ell,j) \!:=\!  - \frac{C_0(\ze,\ell,j)}{2}, \qquad 
B_3(\ze,\ell,j) \!:=\! 0, \\
& 
B_q(\ell,j ) \!:=\! -\frac{C_{q-2}(\ze,\ell,j)}{2}, \quad 4 \le q\le N+1 , \\
\hspace{-3.cm}
& B'_1(\ze,\ell,j) :=  \frac{C_0(\ze,\ell,j)}{2\s (j_{\ell_{\vartheta'}}^\flat - j_\ell^\flat)} , \qquad B'_2(\ze,\ell,j) :=0, \\
& B'_q(\ze,\ell,j) :=\frac{C_{q-1}(\ze,\ell,j)}{2\s  (j_{\ell_{\vartheta'}}^\flat - j_\ell^\flat)}, \quad 3 \le q\le N,
\end{aligned}
\end{equation}
while $\gotG^{(M)}_\ell(\ze,j)$ is a suitable, explicitly computable remainder.
%with $\{\gotG_\ell(\ze,j)\}_{j\in\ZZZ} \in \ell^{N,\io}$.
Similarly, if $\vartheta\in\gotH^{(k)}_{j,\s}$ and $\ell\in L_{\bf{a}}(\vartheta)$ we may 
%epeat the construction above;  by writing
write $1/x_\ell$ as in \eqref{odiotuttidavvero},
with $p_\ell=0$ and $\nu_\ell^\flat$  instead of $\nu_\ell^\flat(\und{y}_\ell)$. The next step is to find a bound for the coefficients 
$\gotg_\ell(\ze,\pow,j)$ and $\gotG^{(M)}_\ell(\ze,j)$ in \eqref{odiotuttidavvero}.

%%%%%%%%%%%%%%%%%%%%%%%%%%%%%%%%%%%%%%%%%%%%%%%%%%%%%%%%%%%%%%%%%%%%%%%%%% 
\begin{rmk}\label{cappaexi}
\emph{
The coefficients $C_q(\ze,\ell,j)$ in \eqref{gliC} depend on $\ka$ also through $\om=\om(\ze)$; thus
if $\ell$ is active, the same happens also to the functions $\gotg_\ell(\ze,\pow,j)$ and $\gotG^{(M)}_\ell(\ze,j)$.
}
\end{rmk}
%%%%%%%%%%%%%%%%%%%%%%%%%%%%%%%%%%%%%%%%%%%%%%%%%%%%%%%%%%%%%%%%%%%%%%%%%% 

%%%%%%%%%%%%%%%%%%%%%%%%%%%%%%%%%%%%%%%%%%%%%%%%%%%%%%%%%%%%%%%%%%%%%%%%%% 
\begin{lemma}\label{belfi}
There exists a constant $K_0=K_0(\al,N)$ such that for any $k\ge1$, $j\in\ZZZ\setminus\{0\}$ and $\s\in\{\pm\}$, and
for any $\vartheta \in\gotH_{j,\s}^{(k)}\cup \gotF\gotH_{j,\s}^{(k)}$, any $\ell\in L_{\bf{a}}(\vartheta)$, and all $\ze\in\matW_N$,
 one has
\begin{subequations}
\begin{align}
 \left| j^q\gotg_\ell(\ze,q,j) \right|	 & \le  K_0 \, \calJ_{\al}(\vartheta)^{2(q-1)} ,
 \qquad q=1,\dots,N-1,
\label{archig} \\
|j^M\gotG_\ell^{(M)}(\ze,j)|& \le  K_0  \, (\calJ_{\al}(\vartheta))^{2M} ,
\qquad
M=1,\dots,N,
\label{archir} 
\end{align}
\end{subequations}
where \eqref{archig} has to be disregarded for $N=1$.
%for all $(\ka,\x)\in\matW_N$ and for any active line $\ell$  of a tree $\vartheta \in\gotH_{j,\s}^{(k)}\cup \gotF\gotH_{j,\s}^{(k)}$.
\end{lemma}
%%%%%%%%%%%%%%%%%%%%%%%%%%%%%%%%%%%%%%%%%%%%%%%%%%%%%%%%%%%%%%%%%%%%%%%%%% 

%%%%%%%%%%%%%%%%%%%%%%%%%%%%%%%%%%%%%%%%%%%%%%%%%%%%%%%%%%%%%%%%%%%%%%%%%% 
\prova
By the definitions \eqref{B}, and the bounds in Lemmata \ref{in} and \ref{helsinki}, one has
\begin{subequations}
\begin{align}
|B_1(\ze,\ell,j)|	 & \le  2 K \, \calJ_{\al}(\vartheta) , \qquad |B_2(\ze,\ell,j)|,\, |B_1'(\ze,\ell,j)|\le K (\calJ_{\al}(\vartheta))^2,
%\Biggl( J(\breve\vartheta) + \sum_{v\in N_2(\vartheta)}J(\breve\vartheta_v) \Biggr)^{2/\alpha}\,, 
\label{archi12} \\
|B_q(\ze,\ell,j)|	 & \le  K  \, (\calJ_{\al}(\vartheta))^{q-4}, \qquad q \ge 4, \\
|B'_q(\ze,\ell,j)| & \le K (\calJ_{\al}(\vartheta))^{q-3}, \qquad q\ge 3 ,
%\Biggl( J(\breve\vartheta) + \sum_{v\in N_2(\vartheta)}J(\breve\vartheta_v) \Biggr)^{(q-2)/\alpha}\,, 
\label{archiq}
\end{align}
\end{subequations}
which, combined with \eqref{mava} immediately implies \eqref{archig}. In order to prove \eqref{archir} we distinguish among three cases.

\begin{itemize}[topsep=0ex]
\itemsep0em

\item If $|j| \le (M+5)(\calJ_{\al}(\vartheta))^2$, we have
\[
%\begin{aligned}
%|\gotG_\ell^{(M)}(\ka,\x,j)| &= \left|
%\frac{j^{M-1}}{P_2(j)+ \displaystyle{\sum_{q=0}^{M-1}}j^{M-1-q}C_q(\ze,\ell,j)+j^{M-1}r_M(\x,\ell,j)} - \sum_{q=1}^{M-1}\frac{\gotg_\ell(\ka,\x,q,j)}{j^q}\right| \\
%&\le K_0 (\calJ_{\al}(\vartheta))^{M-1},
|j^{M}\gotG_\ell^{(M)}(\ze,j)|  \le \left|
\frac{j^M}{(x_\ell(\und{y}_\ell))^{p_\ell+1}} \right| + \left|j^M\sum_{q=1}^{M-1} %\frac{
\gotg_\ell(\ze,q,j) %}{j^q}
\right| 
\le K_0 (\calJ_{\al}(\vartheta))^{2M},
%\end{aligned}
\]
where we have
%denoted
%\[
%P_2(j) = j^{M-1}\big(\bigl(\s_\ell \s_{\ell_{\vartheta'}} -1 \bigr) j^2 + 2 \bigl( \s_\ell \s_{\ell_{\vartheta'}} j_{\ell_{\vartheta'}}^\flat  - j_\ell^\flat \bigr) \s j\big)
%\]
%and 
used the bounds \eqref{archig} and the fact that $|x_\ell(\und{y}_\ell)|>\be(m_0)/16$ since $n_\ell=0$.

\item If $|j| > (M+5)(\calJ_{\al}(\vartheta))^2$ and
%, we distinguish between two cases.\begin{itemize}[topsep=0ex]\itemsep0em\item If 
$\s_\ell\s_{\ell_{\vartheta'}}=-1$, setting
\begin{equation}\label{djej1}
\de_j:= \frac{ \s(j^\flat_{\ell_{\vartheta'}} + j_{\ell}^\flat )}{j} - \frac{1}{2} \sum_{q=0}^{M-1}\frac{C_q(\ze,\ell,j)}{j^{q+2}},
\qquad\qquad
\eps_j:= - \frac{r_M(\ze,\ell,j)}{2j^2} ,
\end{equation}
we have
\[
\begin{aligned}
|\gotG^{(M)}_{\ell}(\ze,j) | & = \left| \frac{1}{(-2j^2)^{p_\ell+1}}\frac{1}{(1+\de_j+\eps_j)^{p_\ell+1}} - 
\sum_{q=1}^{M-1} %\frac{
\gotg_\ell(\ze,q,j) %}{j^q} 
\right| \\
& \qquad \le \frac{1}{(2j^2)^{p_\ell+1}} \left|\frac{1}{(1+\de_j+\eps_j)^{p_\ell+1}}  - \frac{1}{(1+\de_j)^{p_\ell+1}} \right| \\ 
 &  \qquad \qquad + \left|\frac{1}{(-2j^2)^{p_\ell+1}}\frac{1}{(1+\de_j)^{p_\ell+1}  } - 
 \sum_{q=1}^{M-1} %\frac{
 \gotg_\ell(\ze,q,j) %}{j^q} 
 \right| .
\end{aligned}
\]
Using that $|\de_j|+|\eps_j| \le 1/2$ in this case,
the second line is bounded by $|r_M(\ze,\ell,j)|/j^4$, whereas the third line equals,
for a suitable $z_*$ with $|z_*|<1/2$,
\begin{equation} \label{bohbohbohbohboh}
\left|
\frac{1}{j^M}\sum_{p=0}^{\io} \frac{A_{p_\ell+1}(p)}{(-2)^{p_\ell+1}}\left(\frac{1}{1+z_*}\right)^{p+p_\ell+1}
\!\!\!\!\!\!\!\!\!\!\!\!\!\!
\sum_{ \substack{ q_1,\ldots,q_p \ge 1 \\ q_1+\ldots+q_p=M-2(p_\ell+1)}} \!\!\!\!\!\!\!\!
B_{q_1}(\ze,\ell,j) \ldots B_{q_p}(\ze,\ell,j) 
\right| ,
\end{equation}
which is turn is bounded proportionally to $(\calJ_{\al}(\vartheta))^{2(M-1)}/|j|^M$.
\item If $|j| > (M+5)(\calJ_{\al}(\vartheta))^2$ and $\s_\ell\s_{\ell_{\vartheta'}}=1$, setting 
\begin{equation}\label{djej2}
\de_j:=   \frac{1}{2\s  (j_{\ell_{\vartheta'}}^\flat - j_\ell^\flat )} \sum_{q=0}^{M-1}\frac{C_q(\ze,\ell,j)}{j^{q+1}},
\qquad\qquad
\eps_j:= \frac{r_M(\ze,\ell,j)}{2\s j  (j_{\ell_{\vartheta'}}^\flat - j_\ell^\flat )} ,
\end{equation}
we have, reasoning as above,
\[
\begin{aligned}
|\gotG^{(M)}_{\ell}(\ze,j) | & \le \frac{1}{|2 (j_{\ell_{\vartheta'}}^\flat - j_\ell^\flat )j|^{p_\ell+1}} 
\left|\frac{1}{(1+\de_j+\eps_j)^{p_\ell+1}  }  - \frac{1}{(1+\de_j)^{p_\ell+1} } \right|  \\
& \qquad + \left|\frac{1}{(2 (j_{\ell_{\vartheta'}}^\flat - j_\ell^\flat )j)^{p_\ell+1}}\frac{1}{(1+\de_j)^{p_\ell+1}}  - 
 \sum_{q=1}^{M-1} %\frac{
 \gotg_\ell(\ze,q,j) %}{j^q} 
 \right| .
\end{aligned}
\]
with the contribution on the r.h.s.~in the first line still smaller than a constant times $|r_M(\ze,\ell,j)|/j^2$ and the second line equal to
\begin{equation} \label{tramonto}
\null\!\!\!\!
\left| 
\frac{1}{j^M}
 \sum_{p=0}^{\io} \frac{A_{p_\ell+1}(p)}{(2\s  (j_{\ell_{\vartheta'}}^\flat - j_\ell^\flat ))^{p_\ell+1}}
\left(\frac{1}{1+z_*}\right)^{p+p_\ell+1}
\!\!\!
\!\!\!\!\!\!\!\!\!\!\!\!\!\!
\sum_{ \substack{ q_1,\ldots,q_p \ge 1 \\ q_1+\ldots+q_p=M-(p_\ell+1)}}\!\!\!\!\!\!\!\!\!\!\!\!\!\!
B'_{q_1}(\ze,\ell,j) \ldots  B'_{q_p}(\ze,\ell,j) \right| , 
\end{equation}
for a suitable $z_*$ with $|z_*|<1/2$, so that the bound $(\calJ_{\al}(\vartheta))^{2(M-1)}/|j|^M$ follows once more.
\end{itemize}
Collecting together all the bounds above we obtain \eqref{archir}.
\EP
%%%%%%%%%%%%%%%%%%%%%%%%%%%%%%%%%%%%%%%%%%%%%%%%%%%%%%%%%%%%%%%%%%%%%%%%%% 

Given any expanded $\eta$-tree $\vartheta\in\gotH^{(k)}_{j,\s}$, with $j\ne0$,
%we want to expand in inverse powers of $j$ the propagators of the active lines, and associate with any of such line
%the contribution of order $\pow_\ell$ in $1/j$.
%To this aim, 
with each line $\ell\in L(\vartheta)$ we associate a propagator
$\overline\matG_\ell(\ze)$ as follows:
\begin{enumerate}[topsep=0ex]
\itemsep0em
\item if $\ell$ is either neutral or small, we set $\overline\matG_\ell(\ze):=
\matG_\ell(\om(\ze))$, with $\matG_\ell(\om)$ defined in \eqref{prop};
\item if $\ell$ is active and $|L_{\bf a}(\vartheta)|\le N-1$, we set %and $\vartheta'$ is the subtree such that $\ell\in\calP(\vartheta')$ we set
\begin{equation}\label{ffa}
\overline \matG_\ell(\ze):=
\left\{\begin{aligned} 
& \Psi_{n_\ell}(x_\ell) \, \gotg_{\ell}(\ze,\pow_\ell,j) ,\qquad \pow_\ell=1,\ldots,M(\vartheta)-1, \\
& \Psi_{n_\ell}(x_\ell)  \, \gotG_{\ell}^{(M(\vartheta))}(\ze,j), \qquad  \pow_{\ell}=M(\vartheta) ;
\end{aligned}
\right.
\end{equation}
\item if $\ell$ is active and $|L_{\bf a}(\vartheta)|\ge N$, so that $\pow_{\ell}= 1$,
we set % and $\vartheta'$ is the subtree such that $\ell\in\calP(\vartheta')$ we set
\begin{equation}\label{ffa1}
\overline \matG_\ell(\ze):=\matG_\ell(\om(\ze))=
%\left\{\begin{aligned} 
%& \Psi_{n_\ell}(x_\ell) \gotg_{\ell}(\ka,\x,\pow_\ell) ,\qquad \pow_\ell=1,\ldots,N-1, \\
%& 
\Psi_{n_\ell}(x_\ell)  \gotG^{(1)}_{\ell}(\ze,j) . %\qquad  \pow_{\ell}= 1.
%\end{aligned}
%\right.
\end{equation}
\end{enumerate}

Similarly, given any renormalized $\eta$-tree $\vartheta\in\gotF\gotH^{(k)}_{j,\s}$, with $j\ne0$,
%we want to expand in inverse powers of $j$ the propagators of the active lines, and associate with any of such line
%the contribution of order $\pow_\ell$ in $1/j$.
%To this aim, 
with each line $\ell\in L(\vartheta)$ we associate a renormalized propagator
$\overline\matG^\ttR_\ell(\ze)$ as follows:
\begin{enumerate}[topsep=0ex]
\itemsep0em
\item if $\ell$ is either neutral or small, we set $\overline\matG^\ttR_\ell(\ze):=
\matG^\ttR_\ell(\om(\ze))$, with $\matG^\ttR_\ell(\om)$ defined in \eqref{gr};
\item if $\ell$ is active and $|L_{\bf a}(\vartheta)|\le N-1$, we set 
\begin{equation}\label{ffa2}
\overline \matG_\ell^\ttR(\ze):=
\begin{cases} 
\Psi^{\del_\ell,p_\ell}_{n_\ell}(x_\ell(\und{y}_\ell))
\, \gotg_{\ell}(\ze,\pow_\ell,j) , & \quad \pow_\ell=1,\ldots,M(\vartheta)-1, \\
\Psi^{\del_\ell,p_\ell}_{n_\ell}(x_\ell(\und{y}_\ell))
  \, \gotG_{\ell}^{(M(\vartheta))}(\ze,j), & \quad  \pow_{\ell}=M(\vartheta) ,
\end{cases}
\end{equation}
with $\Psi^{\del_\ell,p_\ell}_{n_\ell}(x_\ell(\und{y}_\ell))$ as in \eqref{miele};
\item if $\ell$ is active and $|L_{\bf a}(\vartheta)|\ge N$, so that $\pow_{\ell}= 1$,
we set 
\begin{equation}\label{ffa2-no2}
\overline \matG_\ell^\ttR(\ze) := \matG^\ttR_\ell(\om(\ze))=
 \Psi^{\del_\ell,p_\ell}_{n_\ell}(x_\ell(\und{y}_\ell)) \, \gotG^{(1)}_{\ell}(\ze,j) . 
\end{equation}
\end{enumerate}

%%%%%%%%%%%%%%%%%%%%%%%%%%%%%%%%%%%%%%%%%%%%%%%%%%%%%%%%%%%%%%%%%%%%%%%%%% 
\begin{rmk}\label{ciaociao1}
\emph{%remark
In fact, for $\vartheta\in\gotF\gotH$, we use the expansion in \eqref{icselle} only if $\ell$ is active.
Indeed, given a tree $\vartheta\in {\gotF\gotH}^{(k)}_{j,+}$,
if $\ell \in L_{\bf{n}}(\vartheta)$ and $\vartheta'\subseteq \vartheta$ is such that $\ell\in\calP_{\vartheta'}$,
one has  $(j_\ell , \s_\ell)= (j_{\ell_{\vartheta'}},\s_{\ell_{\vartheta'}})$ 
which, inserted into \eqref{gliC}, yields that $C_q(\ze,\ell,j)=0$ for all $q \ge 1$ and
$x_\ell(\und{y}_\ell)=C_0(\ze,\ell,j)=\s_\ell \, \omega \cdot\nu_\ell^\flat(\und{y}_\ell) $.
On the other hand, as we shall see (see Corollary \ref{small} below), we will not need to expand
the propagators of the small lines in $L_{\bf{s}}(\vartheta)$.
Similar considerations hold for the trees in $\gotH$.
}%remark
\end{rmk}
%%%%%%%%%%%%%%%%%%%%%%%%%%%%%%%%%%%%%%%%%%%%%%%%%%%%%%%%%%%%%%%%%%%%%%%%%% 

%%%%%%%%%%%%%%%%%%%%%%%%%%%%%%%%%%%%%%%%%%%%%%%%%%%%%%%%%%%%%%%%%%%%%%%%%%
\begin{defi}[\textbf{Value of an $\boldsymbol{\eta}$-tree}]\label{defvalue9}
The value $\overline{\Val}(\vartheta;c,\ze)$ of an expanded $\h$-tree $\vartheta\in\gotH$ is defined as
\begin{equation}\label{vala}
\ol\Val(\vartheta;c,\ze):=\Big(\prod_{\la\in \Lambda(\vartheta)}\LL_\la (c)\Big) 
\Big(\prod_{v\in N(\vartheta)} \starF_v(c)\Big)
\Big(\prod_{\ell\in L(\vartheta)} \overline \matG_\ell(\ze) \Big) ,
\end{equation}
while the renormalized value $\overline{\Val}^\ttR(\vartheta;c,\ze)$ of a fully renormalized
$\h$-tree $\vartheta\in\gotF\gotH$ is defined as in \eqref{valR}, 
provided in \eqref{valR} as well as in \eqref{termounoderiregola1}--\eqref{termounoderiregola3} the renormalized propagators $\matG^\ttR_\ell(\om)$
are replaced with the new renormalized propagators $\overline \matG^\ttR_\ell(\ze)$.
\end{defi}
%%%%%%%%%%%%%%%%%%%%%%%%%%%%%%%%%%%%%%%%%%%%%%%%%%%%%%%%%%%%%%%%%%%%%%%%%%

Finally, for $j\ne0$ and $q=0,\ldots,N-1$, we define
\begin{subequations} \label{glia}
\begin{align}
{\gota}_\pow^{(k)}(c,\ze,j)&:=-\frac{j^q}{c_j}\sum_{\vartheta\in {\gotH}^{(k)}_{j,+}(\pow)}\!\!\!\!\! \overline{\Val}(\vartheta;c,\ze),
\label{gliaa} \\
\tilde{\gota}_\pow^{(k)}(c,\ze,j)&:=-\frac{j^q}{c_j}  \sum_{\vartheta\in {\gotF\gotH}^{(k)}_{j,+}(\pow)}\!\!\!\!\!  \overline{\Val}^\ttR(\vartheta;c,\ze),
\label{gliab}
\end{align}
\end{subequations}
and
\begin{subequations} \label{sisi}
\begin{align}
\gotA^{(k)}(c,\ze,j)&:=-\frac{1}{c_j}  \!\!\! \sum_{\vartheta\in {\gotH}^{(k)}_{j,+}(N)}  \overline{\Val}(\vartheta;c,\ze),
\label{sisia} \\
\tilde\gotA^{(k)}(c,\ze,j)&:=-\frac{1}{c_j} \!\!\! \sum_{\vartheta\in {\gotF\gotH}^{(k)}_{j,+}(N)}\overline{\Val}^\ttR(\vartheta;c,\ze) .
\label{sisib}
\end{align}
\end{subequations}
%

%%%%%%%%%%%%%%%%%%%%%%%%%%%%%%%%%%%%%%%%%%%%%%%%%%%%%%%%%%%%%%%%%%%%%%%%%% 
\begin{rmk}\label{uguali}
\emph{
By reasoning as in the proof of Lemma \ref{welldefined}, and using the bounds on the propagators provided by Lemma \ref{belfi},
 one finds, for $q=0,1,\ldots,N$, 
\begin{equation}\label{burp}
 \sum_{\vartheta\in {\gotH}^{(k)}_{j,+}(\pow)}\left| \frac{j^q}{c_j} \overline{\Val}(\vartheta;c,\ze)\right|,
  \sum_{\vartheta\in {\gotF\gotH}^{(k)}_{j,+}(\pow)}\left| \frac{j^q}{c_j} \overline{\Val}^\ttR (\vartheta;c,\ze)\right| \le C(k)^{k} ,
\end{equation}
with $C(k)$ diverging as $k\to\io$. Therefore, the series in \eqref{glia} as well as those in \eqref{sisi} are absolutely convergent. 
}
\end{rmk}
%%%%%%%%%%%%%%%%%%%%%%%%%%%%%%%%%%%%%%%%%%%%%%%%%%%%%%%%%%%%%%%%%%%%%%%%%% 

%%%%%%%%%%%%%%%%%%%%%%%%%%%%%%%%%%%%%%%%%%%%%%%%%%%%%%%%%%%%%%%%%%%%%%%%%% 
\begin{lemma} \label{atildea}
For all $j\in\ZZZ\setminus\{0\}$, one has
$${\gota}_\pow^{(k)}(c,\ze,j)=\tilde{\gota}_\pow^{(k)}(c,\ze,j) ,
\quad q=0,\ldots,N-1, \qquad \gotA^{(k)}(c,\ze,j)=\tilde\gotA^{(k)}(c,\ze,j) . $$
\end{lemma}
%%%%%%%%%%%%%%%%%%%%%%%%%%%%%%%%%%%%%%%%%%%%%%%%%%%%%%%%%%%%%%%%%%%%%%%%%% 

%%%%%%%%%%%%%%%%%%%%%%%%%%%%%%%%%%%%%%%%%%%%%%%%%%%%%%%%%%%%%%%%%%%%%%%%%% 
\prova
One reasons as in the proof of Lemma \ref{effone}, noting that the cancellation ensured by Lemma \ref{cancellazione}
is not affected by the presence of the active lines. Indeed, by Remark \ref{secattivo}, 
for any tree $\vartheta\in\gotH$, a line $\ell\in L_{\rm a}(\vartheta)$, if resonant, connects resonant clusters $T$
on scale $-1$, i.e.~trivial resonant clusters, whose value does not depend on the propagator of the entering line.
\EP
%%%%%%%%%%%%%%%%%%%%%%%%%%%%%%%%%%%%%%%%%%%%%%%%%%%%%%%%%%%%%%%%%%%%%%%%%% 

%%%%%%%%%%%%%%%%%%%%%%%%%%%%%%%%%%%%%%%%%%%%%%%%%%%%%%%%%%%%%%%%%%%%%%%%%% 
\begin{rmk}\label{motooscrivo}
In the following, when using ${\gota}_\pow^{(k)}(c,\ze,j)$ or $\tilde{\gota}_\pow^{(k)}(c,\ze,j)$ we highlight the fact that we are using the representation
\eqref{gliaa} or \eqref{gliab} respectively. In view of Lemma \ref{atildea}, one may wonder why we introduce both ${\gota}_\pow^{(k)}(c,\ze,j)$ and 
$\tilde{\gota}_\pow^{(k)}(c,\ze,j)$. The reason is that, similarly to Lemma \ref{converge1}, in Lemma \ref{converge2} below we shall provide the 
bound \eqref{sanguebis} which does not hold true for the values $ \overline{\Val}(\vartheta;c,\ze)$. Thus, we shall use the representation
$\tilde{\gota}_\pow^{(k)}(c,\ze,j)$ whenever we need to use such bound.
\end{rmk}
%%%%%%%%%%%%%%%%%%%%%%%%%%%%%%%%%%%%%%%%%%%%%%%%%%%%%%%%%%%%%%%%%%%%%%%%%% 

%%%%%%%%%%%%%%%%%%%%%%%%%%%%%%%%%%%%%%%%%%%%%%%%%%%%%%%%%%%%%%%%%%%%%%%%%% 
%%%%%%%%%%%%%%%%%%%%%%%%%%%%%%%%%%%%%%%%%%%%%%%%%%%%%%%%%%%%%%%%%%%%%%%%%%
\subsection*{\ref{labifsec}.3\hspace{0.5cm}Bounds on the values of the %expanded and fully renormalized 
$\boldsymbol{\h}$-trees}}
\addcontentsline{toc}{subsection}{\ref{labifsec}.3\hspace{0.4cm}Bounds on the values of the %expanded and fully renormalized
${\h}$-trees}
\label{boundsetatreevalues}
%%%%%%%%%%%%%%%%%%%%%%%%%%%%%%%%%%%%%%%%%%%%%%%%%%%%%%%%%%%%%%%%%%%%%%%%%%
%%%%%%%%%%%%%%%%%%%%%%%%%%%%%%%%%%%%%%%%%%%%%%%%%%%%%%%%%%%%%%%%%%%%%%%%%%

\noindent
Now, we discuss how to bound the values of the $\eta$-trees. We prove that
for any tree $\vartheta\in\gotH\gotF^{(k)}_{j,+}(q)$
the quantity $(j^q/c_j)\Val^\ttR({\vartheta};c,\ze) $ admits a bound
uniform in $j$ and for a large subset of trees the bound it can be improved into
a bound exponentially small in $\jap{j}^\al$.

%%%%%%%%%%%%%%%%%%%%%%%%%%%%%%%%%%%%%%%%%%%%%%%%%%%%%%%%%%%%%%%%%%%%%%%%%% 
\begin{lemma}\label{converge2}
For all $s'\in(0,s)$ and all $s_3\in(0,s')$, 
there is a positive constant $K_1=K_1(\al,\ze,N)$ such that,
for all $k\ge1$, all $j\in\ZZZ\setminus\{0\}$,
all $q=0,1,\ldots,N$ and any $\vartheta \in \gotH\gotF^{(k)}_{j,+}(q)$, one has
\begin{equation} \label{sanguebis}
\left| \frac{j^q}{c_j} \ol\Val^\ttR({\vartheta};c,\ze) \right| \le 
A^{|N(\vartheta)|} 
K_1^{|N(\vartheta)|}
(\calJ_{\al}(\vartheta))^{2N}e^{-(s'-s_3)\calJ(\vartheta)} ,
\end{equation}
with $A$ as in Lemma \ref{converge1}.
\end{lemma}
%%%%%%%%%%%%%%%%%%%%%%%%%%%%%%%%%%%%%%%%%%%%%%%%%%%%%%%%%%%%%%%%%%%%%%%%%% 

\prova
%%%%%%%%%%%%%%%%%%%%%%%%%%%%%%%%%%%%%%%%%%%%%%%%%%%%%%%%%%%%%%%%%%%%%%%%%% 
Consider any $\vartheta\in\gotF\gotH^{(k)}_{j,+}(q)$. 
Since the active lines are non-resonant, we can bound the product of the propagators of the neutral and small lines as in the proof of Lemma \ref{converge1},
by relying on Lemma \ref{ribry}, by taking $s_1=0$.

In order to bound the propagators of the active lines, still using the notations of Lemma \ref{converge1}
with $s_1 =0$, we distinguish among the following cases.
\begin{itemize}[topsep=0ex]
\itemsep0em
\item
If $q=0,\ldots,N-1$, for any $\ell\in L_{\bf{a}}(\vartheta)$,
thanks to Remark \ref{contiamole} we can use \eqref{archig} to bound
$|\overline \matG_\ell^\ttR(\ze)|\le K_0 |j|^{-q_\ell} (\calJ_\al(\vartheta))^{2(q_\ell-1)}$, so that
\begin{equation} \label{caso1}
\left| \frac{j^q}{c_j} \ol{\Val}^\ttR(\vartheta;c,\ze)\right| \le A^{|N(\vartheta)|} K_0^{N-1} (\calJ_{\al}(\vartheta))^{2(q-1)}e^{-(s'-s_3)\calJ(\vartheta)}.
\end{equation}
\item
If $q=N$ and $|L_{\bf a}(\vartheta)|<N$, then we bound the propagator of each active line with
$|\overline \matG_\ell^\ttR(\ze)|\le K_0 |j|^{-q_\ell} (\calJ_\al(\vartheta))^{2q_\ell}$ and use that $q_{\ell} \le M(\vartheta)<N$, so as to get
\begin{equation} \label{caso2}
\left| \frac{1}{c_j} \ol{\Val}^\ttR(\vartheta;c,\ze)\right| \le A^{|N(\vartheta)|} K_0^{N-1} |j|^{-N} (\calJ_{\al}(\vartheta))^{2N^2}e^{-(s'-s_3)\calJ(\vartheta)}.
\end{equation}
\item
If $q=N$ and $|L_{\bf a}(\vartheta)| \ge N$, then we bound the propagator of $N$, arbitrarily chosen active lines with 
$|\overline \matG_\ell^\ttR(\ze)|\le K_0 |j|^{-1} (\calJ_\al(\vartheta))^{2}$ , whereas we bound the propagators of the remaining active lines
by using that $|x_\ell(\und{y}_\ell)| \ge \be(m_0)/32$, which gives
\begin{equation} \label{caso3}
\left| \frac{1}{c_j} \ol{\Val}^\ttR(\vartheta;c,\ze)\right| \le A^{|N(\vartheta)|} K_1^{|N(\vartheta)|}
|j|^{-N} (\calJ_{\al}(\vartheta))^{2N}e^{-(s'-s_3)\calJ(\vartheta)},
\end{equation}
where $K_1=\max\{K_0,32/\be(m_0)\}$.
\end{itemize}
Collecting together the bounds above gives the bound \eqref{sanguebis}.
\EP
%%%%%%%%%%%%%%%%%%%%%%%%%%%%%%%%%%%%%%%%%%%%%%%%%%%%%%%%%%%%%%%%%%%%%%%%%% 

%%%%%%%%%%%%%%%%%%%%%%%%%%%%%%%%%%%%%%%%%%%%%%%%%%%%%%%%%%%%%%%%%%%%%%%%%% 
\begin{lemma}\label{trentuno}
Let $s'$, $s_3$, $A$ and $K_1$ be as in Lemma \ref{converge2}.
For all $k\ge1$, all $j\in\ZZZ\setminus\{0\}$, all $q=0,1,\ldots,N$, all $\hat s, \bar s>0$ such that $\hat s+\bar s=s'-s_3$, and
all $\vartheta\in\gotF\gotH^{(k)}_{j,+}(q)$ such that either $|\Lambda_{j,+}(\vartheta)| \ge 2$ or
$|\Lambda_{j',+}(\vartheta)| \ge 1$ for some $j'$ with $|j'|> |j|$, one has
\begin{equation}\label{giuntuno}
\left| \frac{j^q}{c_j}\ol\Val^\ttR (\vartheta;c,\ze) \right| \le 
A^{|N(\vartheta)|}K_1^{|N(\vartheta)|}  (\calJ_{\al}(\vartheta))^{2N^2}
e^{-2 \hat s \jap{j}^{\al}}  e^{- \bar  s  \calJ(\vartheta)}.
 \end{equation}
\end{lemma}
%%%%%%%%%%%%%%%%%%%%%%%%%%%%%%%%%%%%%%%%%%%%%%%%%%%%%%%%%%%%%%%%%%%%%%%%%% 

%%%%%%%%%%%%%%%%%%%%%%%%%%%%%%%%%%%%%%%%%%%%%%%%%%%%%%%%%%%%%%%%%%%%%%%%%% 
\prova
Reason as in the proof of Corollary \ref{trenta}, by relying on the bounds \eqref{caso1} to \eqref{caso3} in the proof of Lemma \ref{converge2}.
\EP
%%%%%%%%%%%%%%%%%%%%%%%%%%%%%%%%%%%%%%%%%%%%%%%%%%%%%%%%%%%%%%%%%%%%%%%%%% 

We have the following result.

%%%%%%%%%%%%%%%%%%%%%%%%%%%%%%%%%%%%%%%%%%%%%%%%%%%%%%%%%%%%%%%%%%%%%%%%%% 
\begin{coro}\label{stimo}
For all $k\ge1$ and all $j\in\ZZZ\setminus\{0\}$
there is a positive constant $C=C(s,\al,\ze,N)$ such that
\[
|\gota_\pow^{(k)}(c,\ze,j)| %, \; |\tilde\gota_\pow^{(k)}(c,\ze,j)|
\le C^k,\quad \pow=0,\ldots,N-1
,\qquad\qquad
|\gotA^{(k)}(c,\ze,j)| %, \; |\tilde\gotA^{(k)}(c,\ze,j)|
\le \frac{C^k}{|j|^N}.
\]
\end{coro}
%%%%%%%%%%%%%%%%%%%%%%%%%%%%%%%%%%%%%%%%%%%%%%%%%%%%%%%%%%%%%%%%%%%%%%%%%% 

%%%%%%%%%%%%%%%%%%%%%%%%%%%%%%%%%%%%%%%%%%%%%%%%%%%%%%%%%%%%%%%%%%%%%%%%%% 
\prova
Thanks to Lemma \ref{atildea}, we may consider $\tilde\gota_\pow^{(k)}(c,\ze,j)$ and  $\tilde\gotA^{(k)}(c,\ze,j)$.
Then the bounds
\[
|\tilde\gota_\pow^{(k)}(c,\ze,j)|\le C^k,\quad q=0,\ldots,N-1 , \qquad\qquad
|\tilde\gotA^{(k)}(c,\ze,j)|\le \frac{C^k}{|j|^N}.
\]
follow using Lemma \ref{converge2}.   We perform the sum over the trees as in Proposition \ref{voltabuona}, with $s_1=0$, $s_3=s/5$,   and $s'= 4 s/5$.  
The term $(\calJ_\al(\vartheta))^{2N}$ is compensated by the exponential term and gives an extra factor $(5/3 s)^{2\al N^2}  (2\al N^2)!$.
\EP
%%%%%%%%%%%%%%%%%%%%%%%%%%%%%%%%%%%%%%%%%%%%%%%%%%%%%%%%%%%%%%%%%%%%%%%%%% 

Thanks to Corollary \ref{stimo}, for $j\ne0$ we can define
\begin{equation}\label{momovedo}
\begin{aligned}
{\gota}_\pow(c,\ze,\e,j)&:=  \sum_{k\ge1}\e^k \gota_\pow^{(k)}(c,\ze, j) ,
\quad \pow=0,\ldots,N-1, \\
\gotA(c,\ze,\e,j)&:=  \sum_{k\ge1}\e^k \gotA^{(k)}(c,\ze,j) . 
\end{aligned}
\end{equation}
Collecting together all the results above leads to the following proposition.

%%%%%%%%%%%%%%%%%%%%%%%%%%%%%%%%%%%%%%%%%%%%%%%%%%%%%%%%%%%%%%%%%%%%%%%%%% 
\begin{prop}\label{stimeta}
		For all $N\ge1$, all $j\in\ZZZ\setminus\{0\}$ and all $\zeta\in\matK_N$ one has
		\[
		\h_j(c,\om(\ze),\e)=\sum_{\pow=0}^{N-1}\frac{\gota_\pow(c,\ze,\e,j)}{j^\pow} +\gotA(c,\ze,\e,j),
		\]
with $\gota_\pow(c,\ze,\e,j)$ and $\gotA(c,\ze,\e,j)$ satisfying the  bounds 
\begin{equation}
	\label{stimafinalea,A}
	|\gota_\pow(c,\ze,\e,j)| \le D|\e| , \qquad |\gotA(c,\ze,\e,j)|\le \frac{D|\e|}{{|j|}^N} ,
\end{equation}
for some positive constant $D=D(s,\al,\ze,N)$ and for all $\e\in (-C^{-1},C^{-1})$ with $C=C(s,\al,\ze,N)$ as in Corollary \ref{stimo}.
\end{prop}
%%%%%%%%%%%%%%%%%%%%%%%%%%%%%%%%%%%%%%%%%%%%%%%%%%%%%%%%%%%%%%%%%%%%%%%%%% 

%%%%%%%%%%%%%%%%%%%%%%%%%%%%%%%%%%%%%%%%%%%%%%%%%%%%%%%%%%%%%%%%%%%%%%%%%% 
%%%%%%%%%%%%%%%%%%%%%%%%%%%%%%%%%%%%%%%%%%%%%%%%%%%%%%%%%%%%%%%%%%%%%%%%%%
\setcounter{subsection}{3}
\subsection{Proof of Proposition \ref{decadeta}}
\label{proofdecadeta}
%%%%%%%%%%%%%%%%%%%%%%%%%%%%%%%%%%%%%%%%%%%%%%%%%%%%%%%%%%%%%%%%%%%%%%%%%%
%%%%%%%%%%%%%%%%%%%%%%%%%%%%%%%%%%%%%%%%%%%%%%%%%%%%%%%%%%%%%%%%%%%%%%%%%%

The bounds provided by Corollary \ref{stimo} consider together the contributions arising from all trees in $\gotF\gotH^{(k)}_{j,+}(q)$ 
-- or $\gotH^{(k)}_{j,\s}(q)$, by Lemma \ref{atildea}. To prove Proposition \ref{decadeta} we have to show
that $\gota_\pow(c,\ze,\e,j)$ is independent of $j$ up to corrections of order $j^{-N}$.
In fact, we prove that all contributions to $\tilde\gota^{(k)}_q(c,\ze,j)$  in \eqref{gliab} which depend explicitly on $j$
are stretched-exponentially small in $j$. 

As a preliminary step we prove the following result.

%%%%%%%%%%%%%%%%%%%%%%%%%%%%%%%%%%%%%%%%%%%%%%%%%%%%%%%%%%%%%%%%%%%%%%%%%% 
\begin{lemma}\label{typeb}
For any $k\ge1$, $j\in\ZZZ\setminus\{0\}$, $\pow =0,1,\ldots,N$ and any $\vartheta\in \gotF\gotH^{(k)}_{j,+}(\pow)$ such that 
${\ol{\Val}^\ttR} (\vartheta;c,\ze)\ne0$, if either $\ell \in L_{\bf s}(\vartheta)$ or $\ell \in L_{\bf a}(\vartheta)$ and $|x_\ell(\und{y}_\ell)|<1$, one has 
\begin{equation}\label{pa}
|j| \le 4 \max\{|j_\ell^\flat|^2 ,  |\nu^\flat_\ell |_2\} .
\end{equation}
\end{lemma}
%%%%%%%%%%%%%%%%%%%%%%%%%%%%%%%%%%%%%%%%%%%%%%%%%%%%%%%%%%%%%%%%%%%%%%%%%% 

%%%%%%%%%%%%%%%%%%%%%%%%%%%%%%%%%%%%%%%%%%%%%%%%%%%%%%%%%%%%%%%%%%%%%%%%%% 
\prova
Given $\vartheta\in  \gotF\gotH^{(k)}_{j,+}(q)$, if $\ell$ is a line in  $L_{\bf{s}}(\vartheta)$ one has $n_\ell \ge 1$ and
$(j_\ell , \s_\ell) \neq (j_{\ell_{\vartheta'}},\s_{\ell_{\vartheta'}})$, where $\vartheta'$ is the subtree of $\vartheta$ such that $\ell\in\calP_{\vartheta'}$.
Write $x_\ell(\und{y}_\ell)$ according to \eqref{ambrogio} and assume by contradiction that
\eqref{pa} does not hold.
Then, shortening $\om(\ze)=\om$, one finds, by using the fourth bound in \eqref{stoclaim},
\[
\begin{aligned}
| \s_\ell(\om\cdot\nu_{\ell}^\flat(\und{y}_\ell)
 + \s_{\ell_{\vartheta'}}\om_{j} ) %- \s_\ell (j_\ell^\flat+ j)} ) 
- \om_{j_\ell}| 
& \ge | \s_\ell(\om\cdot\nu_{\ell}^\flat(\und{y}_\ell) +\s_{\ell_{\vartheta'}} j^2) - \bigl( (j_\ell^\flat)^2 + j^2 + 2j_\ell^\flat j \bigr) | - 1 \\ %\phantom{\sum_{\ZZZ}} \\
& \ge |(\s_\ell \s_{\ell_{\vartheta'}} -1)j^2 - 2 j^\flat_\ell j | -  (j_\ell^\flat)^2  -  | \om\cdot\nu_{\ell}^\flat(\und{y}_\ell)| - 1 \\
%\sum_{i\in\ZZZ}i^2(\nu_\ell^\flat(\und{y}_\ell))_i \Bigr| \\
& \ge 2 |j| - 2 \bigl( (j_\ell^\flat)^2 + |\nu_\ell^\flat|_2 \bigr) \ge1,
\end{aligned}
\]
which is not possible if ${\ol{\Val}^\ttR} (\vartheta;c,\ze)\ne0$ and  $\ell\in L_{\bf{s}}(\vartheta)$.
If $\ell\in L_{\bf a}(\vartheta)$, then one has $(j_\ell , \s_\ell) \neq (j_{\ell_{\vartheta'}},\s_{\ell_{\vartheta'}})$, where $\vartheta'$ is the subtree of 
$\vartheta$ such that $\ell\in\calP_{\vartheta'}$,
so that, if additionally $|x_\ell(\und{y}_\ell)|<1$,  the same argument applies.
\EP
%%%%%%%%%%%%%%%%%%%%%%%%%%%%%%%%%%%%%%%%%%%%%%%%%%%%%%%%%%%%%%%%%%%%%%%%%% 
 
For $q=0,1,\ldots,N$, let ${\gotH}^{(k)}_{j,\s}(q,\bf{s})$ and ${\gotF\gotH}^{(k)}_{j,\s}(q,{\bf s})$
denote the set of trees $\vartheta\in{\gotH}^{(k)}_{j,\s}(q)$ and $\vartheta\in{\gotF\gotH}^{(k)}_{j,\s}(q)$, respectively,
such that $L_{\bf{s}}(\vartheta)\ne\emptyset$. By construction, one has
\begin{equation} \label{solos}
\sum_{\vartheta\in {\gotH}^{(k)}_{j,+}(q,\bf{s})} \overline{\Val}(\vartheta;c,\ze)= \sum_{\vartheta\in {\gotF\gotH}^{(k)}_{j,+}(q,\bf{s})}
{\overline{\Val}^\ttR} (\vartheta;c,\ze).
\end{equation}
%
%The bounds provided by Corollary \ref{stimo} consider together the contributions arising from all trees in $\gotH^{(k)}_{j,+}(q)$ or
%${\gotF\gotH}^{(k)}_{j,\s}(q)$. However, some contributions are stretched-exponentially small in $j$, as the following result shows.

%%%%%%%%%%%%%%%%%%%%%%%%%%%%%%%%%%%%%%%%%%%%%%%%%%%%%%%%%%%%%%%%%%%%%%%%%% 
\begin{coro}\label{small}
For any $s_4\in(0,s'-s_3)$, with $s'$ and $s_3$ as in Lemma \ref{converge2}, 
for all $k\ge1$, all $j\in\ZZZ\setminus\{0\}$, all $\pow =0,1,\ldots,N$ and all $\vartheta\in \gotF\gotH^{(k)}_{j,+}(\pow,{\bf{s}})$
one has $\calJ_\al(\vartheta) \ge 2^{-\al} \jap{j}^{\al/2}$ and
\begin{equation} \label{eq:small}
\left| \frac{j^q}{c_j} {\ol{\Val}^\ttR} (\vartheta;c,\ze) \right|  \le 
A^{|N(\vartheta)|}K_1^{N-1} (\calJ_{\al}(\vartheta))^{2N^2}
e^{-(s'-s_3-s_4)\calJ(\vartheta)} e^{-s_4 2^{-\al}\jap{j}^{\al/2}} .
\end{equation}
The same bounds hold for all $\vartheta \in \gotF\gotH^{(k)}_{j,+}(\pow)$ 
with a line $\ell\in L_{\bf a}(\vartheta)$ with $|x_\ell(\und{y}_\ell)| < 1$.
\end{coro}
%%%%%%%%%%%%%%%%%%%%%%%%%%%%%%%%%%%%%%%%%%%%%%%%%%%%%%%%%%%%%%%%%%%%%%%%%% 

\prova
Given $\vartheta\in  \gotF\gotH^{(k)}_{j,+}(q)$, if there is either a line $\ell\in L_{\bf{s}}(\vartheta)$
or a line $\ell\in L_{\bf a}(\vartheta)$ with $|x_\ell(\und{y}_\ell)|<1$,
then, by Lemma \ref{typeb} and the first two bounds of \eqref{stoclaim}, one has
\[
(\calJ_\al(\vartheta))^2 \ge \max\{|j_\ell^\flat|^2 ,  |\nu^\flat_\ell |_2 \} \ge \frac{1}{4}|j| .
\]
As a consequence, with respect to the bounds in the proof of Lemma \ref{trentuno}, the factor 
$e^{-(s'-s_3)\calJ(\vartheta)}$ can be replaced with $e^{-(s'-s_3-s_4)\calJ(\vartheta)} e^{-2^{-\al}s_4\jap{j}^{\al/2}}$.
\EP

%%%%%%%%%%%%%%%%%%%%%%%%%%%%%%%%%%%%%%%%%%%%%%%%%%%%%%%%%%%%%%%%%%%%%%%%%% 
\begin{rmk} \label{dovesseservire}
\emph{
Corollary \ref{small}, together with Remark \ref{ciaociao1}, implies that the only trees $\vartheta\in\gotF\gotH^{(k)}_{j,+}(\pow)$
whose values may not decay exponentially in $\jap{j}^{\al/2}$ are those which contain at least
one active line and all active lines $\ell$ are such that $|x_\ell(\und{y}_\ell)| \ge 1$.
On the other hand, if this happens, then $\Psi^{\del_\ell,0}_{n_\ell}(x_{\ell}(\und{y}_{\ell})) =1$  and
$\Psi^{\del_\ell,p_\ell}_{n_\ell}(x_{\ell}(\und{y}_{\ell})) =0$  for all $p>0$: thus, 
in such trees, for each active line $\ell$, one has 
$p_\ell=\del_\ell$ and $\overline \matG_\ell^\ttR(\ze)$ equals either
$\gotg_{\ell}(\ze,\pow_\ell,j)$, if $\pow_\ell=1,\ldots,M(\vartheta)-1$,
or $\gotG_{\ell}^{(M(\vartheta))}(\ze,j)$, if $\pow_{\ell}=M(\vartheta)$.
}
\end{rmk}
%%%%%%%%%%%%%%%%%%%%%%%%%%%%%%%%%%%%%%%%%%%%%%%%%%%%%%%%%%%%%%%%%%%%%%%%%%

%%%%%%%%%%%%%%%%%%%%%%%%%%%%%%%%%%%%%%%%%%%%%%%%%%%%%%%%%%%%%%%%%%%%%%%%%%
\begin{defi}[\textbf{$\boldsymbol{jj'}$-associated trees}] \label{conjtree}
Given a tree $\vartheta\in\gotF\gotH^{(k)}_{j,+}(q)$,  let $\vartheta' \in \gotF\gotH^{(k)}_{j',+}(q)$
be the tree -- if existing at all -- obtained from $\vartheta$ by 
\begin{enumerate}[topsep=0ex]
\itemsep0em
\item replacing $j_{\la_{\vartheta}}=j$ with $j'$ (and hence $\gote_{j_{\la_{\vartheta}}}=\gote_j$ with $\gote_{j'}$),
\item replacing $j_{\la_{\vartheta_v}}$ with $j_{\la_{\vartheta_v}}-j+j'$ for all $v\in N_2(\vartheta)$, 
\item leaving the component and momentum labels of  the other leaves unchanged,
\item modifying accordingly the labels of the lines of the path $\calmP({\vartheta})$ 
so as to satisfy the conservation law \eqref{lefoglie*}.
\end{enumerate}
If $T$ is a subgraph of $\vartheta$ we call $T'$ the subgraph
of $\vartheta'$ obtained from $T$ with the procedure described above. If there exists such a tree $\vartheta'$ we say that
$\vartheta'$ is \emph{$j'$-associated} with $\vartheta$, or that $\vartheta$ and $\vartheta'$ are \emph{$jj'$-associated}.
\end{defi}
%%%%%%%%%%%%%%%%%%%%%%%%%%%%%%%%%%%%%%%%%%%%%%%%%%%%%%%%%%%%%%%%%%%%%%%%%%

%%%%%%%%%%%%%%%%%%%%%%%%%%%%%%%%%%%%%%%%%%%%%%%%%%%%%%%%%%%%%%%%%%%%%%%%%%
\begin{rmk} \label{rmk:conjtree}
\emph{
Given a tree $\vartheta$, the tree $j'$-associated with $\vartheta$ may fail to exist:
indeed, if, because of the label changing, the labels of a line $\ell\in\calP_{\vartheta}$ with
$\nu_\ell \neq \gote_{j_\ell}$ assumed values $\nu_{\ell}'$ ad $j_{\ell}'$ such that $\nu_{\ell}'  = \gote_{j_\ell'}$,
then the constraint 7 in Definition \ref{gotT} would not be satisfied.
}
\end{rmk}
%%%%%%%%%%%%%%%%%%%%%%%%%%%%%%%%%%%%%%%%%%%%%%%%%%%%%%%%%%%%%%%%%%%%%%%%%%

%%%%%%%%%%%%%%%%%%%%%%%%%%%%%%%%%%%%%%%%%%%%%%%%%%%%%%%%%%%%%%%%%%%%%%%%%%
\begin{rmk} \label{rmk:neutralindj}
\emph{
As Remark \ref{ciaociao1} shows, the propagators of the neutral lines do not depend on $j$, in the sense
that if $\vartheta'$ is the tree $j'$-associated to a tree $\vartheta\in \gotF\gotH^{(k)}_{j,+}(q)$, then
for any $\ell\in L_{\bf n}(\vartheta)$ the propagator $\overline\matG^\ttR_\ell(\ze)$ does not change.
}
\end{rmk}
%%%%%%%%%%%%%%%%%%%%%%%%%%%%%%%%%%%%%%%%%%%%%%%%%%%%%%%%%%%%%%%%%%%%%%%%%% 

The next result is the key to conclude the proof of Proposition \ref{decadeta}.

%%%%%%%%%%%%%%%%%%%%%%%%%%%%%%%%%%%%%%%%%%%%%%%%%%%%%%%%%%%%%%%%%%%%%%%%%% 
\begin{lemma} \label{presqua}
Let $s'$ and $s_3$ be defined as in Lemma \ref{converge2}.
For all  $k\ge1$, all $j\in\ZZZ\setminus\{0\}$, all $\pow =0,1,\ldots,N$,
given a tree $\vartheta\in\gotF\gotH^{(k)}_{j,+}(q)$ and $j'\in\ZZZ\setminus\{0,j\}$, 
if there exists no tree $j'$-associated with $\vartheta$ then, for any $s_4 \in (0,s'-s_3)$, one has
\begin{equation} \label{notree}
\left| \frac{j^q}{c_j} {\overline{\Val}^\ttR} (\vartheta;c,\ze) \right| \le
A^{|N(\vartheta)|}K_1^{N-1} (\calJ_{\al}(\vartheta))^{2N^2}
e^{-(s'-s_3- s_4)\calJ(\vartheta)} e^{- s_4 \min\{\jap{j}^{\al/2},\jap{j'}^{\al/2}\}} ,
\end{equation}
while if there exists a tree $\vartheta'$ which is $j'$-associated with $\vartheta$ then either
\begin{equation} \label{diffe}
\frac{j^q}{c_j} {\overline{\Val}^\ttR} (\vartheta;c,\ze) =
\frac{(j')^q}{c_{j'}} {\overline{\Val}^\ttR} (\vartheta';c,\ze) 
\end{equation}
or both $|j^qc_j^{-1} {\overline{\Val}^\ttR} (\vartheta;c,\ze)|$ and $|(j')^qc_{j'}^{-1} {\overline{\Val}^\ttR} (\vartheta';c,\ze)|$
are bounded by
\begin{equation} \label{maximum}
A^{|N(\vartheta)|}K_1^{N-1} (\calJ_{\al}(\vartheta))^{2N^2}
e^{-(s'-s_3-s_4)\calJ(\vartheta)} e^{- s_4 2^{-\al}\min\{\jap{j}^{\al/2},\jap{j'}^{\al/2}\}} ,
\end{equation}
for any $s_4 \in (0,s'-s_3)$.
\end{lemma}
%%%%%%%%%%%%%%%%%%%%%%%%%%%%%%%%%%%%%%%%%%%%%%%%%%%%%%%%%%%%%%%%%%%%%%%%%% 

%%%%%%%%%%%%%%%%%%%%%%%%%%%%%%%%%%%%%%%%%%%%%%%%%%%%%%%%%%%%%%%%%%%%%%%%%% 
\prova
%Given any $\vartheta\in\gotF\gotH^{(k)}_{j,+}(q)$,  let $\vartheta' \in \gotF\gotH^{(k)}_{j',+}(q)$ be the tree -- if existing at all -- %\in\gotH^{(k)}_{j',+}$
%obtained from $\vartheta$ by replacing $j_{\la_{\vartheta}}=j$ with $j'$ (and hence $\gote_{j_{\la_{\vartheta}}}=\gote_j$ with $\gote_{j'}$)
%and by replacing $j_{\la_{\vartheta_v}}$ with $j_{\la_{\vartheta_v}}-j+j'$ for all $v\in N_2(\vartheta)$, 
%while leaving the component and momentum labels of 
%the other leaves unchanged, and modifying accordingly the labels of the lines of the path $\calmP({\vartheta})$ 
%so as to satisfy the conservation law \eqref{lefoglie*}. If $T$ is a subgraph of $\vartheta$ we call $T'$ the subgraph
%of $\vartheta'$ obtained from $T$ with the procedure described above.
%
%Note that such a tree may fail to exist: indeed, if, because of the label changing, the labels of a line $\ell\in\calP_{\vartheta}$ with
%$\nu_\ell \neq \gote_{j_\ell}$ assumed values $\nu_{\ell}'$ ad $j_{\ell}'$ such that $\nu_{\ell}'  = \gote_{j_\ell'}$,
%then the constraint 7 in Definition \ref{gotT} would not be satisfied.
%However, if this happens, we deduce from \eqref{lefoglie*} that there is at least one leaf
If there exists no tree $j'$-associated with $\vartheta$, we deduce from \eqref{lefoglie*} that there is at least one leaf
$\la_1\in\Lambda_{j',-}(\vartheta)$ with $\la_1\st{\preceq}\ell$, and, by Lemma \ref{cijei}, since $j\neq j'$,
there is necessarily also a leaf $\la_2\in\Lambda_{j',+}(\vartheta)$, so that
we can use Lemma \ref{trentuno} and obtain the bound \eqref{notree} with $\hat{s}=s_4$.

Otherwise, there exists a tree $\vartheta'$ which is $j'$-associated to $\vartheta$. %, we study the difference
%%
%\begin{equation} \label{diffe}
%\frac{j^q}{c_j} {\overline{\Val}^\ttR} (\vartheta;c,\ze)-
%\frac{(j')^q}{c_{j'}} {\overline{\Val}^\ttR} (\vartheta';c,\ze) . 
%\end{equation}
%%
First of all note that $\calJ(\vartheta')=\calJ(\vartheta)$ and, moreover,
for any $\ell$, if $(j_\ell,\s_\ell)$ and $(j'_\ell,\s'_\ell)$
denote the component and sign label of
$\ell$ as line in $\vartheta$ and $\vartheta'$, respectively, then
by construction $\s_\ell'j'_\ell=\s_\ell j_\ell-j+j'$ and $\s_\ell'=\s_\ell$. This implies that if $\ell$ is neutral as a line of $\vartheta$, 
then it is neutral as a line of $\vartheta'$ as well, and vice versa.
Moreover if $\ell$ belongs to the ramification of $\vartheta$ (and hence of $\vartheta'$), one has
\begin{equation}\label{bemolle}
(j'_\ell)^\flat = j_\ell^\flat.
\end{equation}
In particular $(j'_\ell)^\flat = j_\ell^\flat=0$ if $\ell$ is neutral.
We discuss separately various cases.

\begin{itemize}[topsep=0ex]
\itemsep0em
\item There exists a line $\ell$ which is small as a line of at least one of the two trees $\vartheta$ and $\vartheta'$.
Then, thanks to \eqref{bemolle} and Lemmata \ref{in} and \ref{typeb} one has
\begin{equation}\label{m5}
\min\{|j|,|j'|\}< 4\max\{ (j_\ell^\flat)^2 , |\nu^\flat_\ell|_2\} \le 4(\calJ_\al(\vartheta))^2 .
\end{equation}
Consequently, by reasoning as in the proof of Corollary \ref{small},
the values of both trees satisfy the bound \eqref{eq:small}, 
with the last factor replaced with 
$e^{- s_4 2^{-\al}\min\{\jap{j}^{\al/2},\jap{j'}^{\al/2}\}}$. 
\item
There are no small lines in both $\vartheta$ and $\vartheta'$.
Then, if $\nu_\ell(\und{y}_\ell)$ and $\nu_\ell'(\und{y}_\ell)$ denote the
renormalized momenta of $\ell$ as a line in $\vartheta$ and $\vartheta'$, respectively,
we start proving that, by induction on the depth of the line, either
\begin{equation}\label{nubemolle}
\nu_\ell^\flat(\und{y}_\ell)=\nu_\ell^{\prime\, \flat}(\und{y}_\ell),
\end{equation}
or the values of both trees are bounded again as in \eqref{eq:small}, 
with the last factor replaced with 
$e^{- s_4 2^{-\al}\min\{\jap{j}^{\al/2},\jap{j'}^{\al/2}\}}$. 
If $d(\ell)=0$ then $\s_\ell'\nu_\ell'(\und{y}_\ell) = \s_\ell'\nu_\ell' = \s_\ell\nu_\ell - \gote_j + \gote_{j'}$,
so that \eqref{nubemolle} follows. Assume \eqref{nubemolle} to hold up to depth $d-1$ and call $T$ the relevant RC
in $\vartheta$ with highest depth containing $\ell$.
We distinguish between two subcases.
\begin{itemize}[topsep=0ex]
\itemsep0em
\item %[(i)] %[2.1]
If $\lambda_\vartheta\notin \Lambda(T)$, then $T'$ is still a relevant RC, 
the entering lines of $T$ and $T'$ belongs to $\calmP(\vartheta)$ and $\calmP(\vartheta')$, respectively,
and are both neutral. Thus $(j_{\ell_T'},\s_{\ell_T'})=(j,+)$ and $(j_{\ell_T'},\s_{\ell_T'})=(j',+)$, so that, by the inductive hypothesis,
\[
\nu_\ell(\und{y}_\ell) =\nu_\ell^0(T)
 + \gote_j + y_{d-1} \nu_{\ell_T'}^\flat (\und{y}_{\ell_T'}) , \qquad
\nu_\ell'(\und{y}_\ell) =\nu_\ell^0(T') + \gote_{j'} + y_{d-1} \nu_{\ell_{T'}'}^\flat (\und{y}_{\ell_{T'}'}) ,
\]
which yields \eqref{nubemolle} once more because $\nu_\ell^0(T')=\nu_\ell^0(T)$.
\item %[(ii)] %[2.2]
If $\lambda_\vartheta\in \Lambda(T)$, then  the exiting line of $T$ belongs to $\calmP(\vartheta)$
and is neutral, so that $(j_{\ell_T},\s_{\ell_T})=(j,+)$. Moreover one has 
\[
\sum_{\la\in \Lambda^*(T')} \s_\la \gote_{j_\la} + \s_{\ell'_{T'}}\gote_{j_{\ell'_{T'}}} - \gote_{j'}=
\sum_{\la\in \Lambda^*(T) } \s_\la \gote_{j_\la} + \s_{\ell'_T}\gote_{j_{\ell'_T}} - \gote_j=0
\]
and $J(T')=J(T)$. If $\s_{\ell}(\nu_{\ell} -\gote_{j_{\ell}}) = \s_{\ell'_T}(\nu_{\ell'_T} - \gote_{j_{\ell'_T}})$,
then $|x_{\ell}(\und{y}_\ell)|<1$ because $\Psi_1(x_\ell)=\Psi_1(x_{\ell'_T}) \neq 0$: in that case, by Corollary \ref{typeb},
one has $\calJ(\vartheta)=\calJ(\vartheta') \ge 2^{-\al} \jap{j'}^{\al/2}$. If instead 
$\s_{\ell}(\nu_{\ell} -\gote_{j_{\ell}}) \neq \s_{\ell'_T}(\nu_{\ell'_T} - \gote_{j_{\ell'_T}})$, then $T'$ is a relevant
RC of $\vartheta'$ and hence
\[
\begin{aligned}
\nu_\ell(\und{y}_\ell) & =\nu_\ell^0(T) + \s_{\ell_T'} \gote_{j_{\ell_T'}} +
 y_{d-1} \s_{\ell_T'} ( \nu_{\ell_T'} (\und{y}_{\ell_T'}) - \gote_{j_{\ell_T'}}) , \\
\nu_\ell'(\und{y}_\ell) & =\nu_\ell^0(T') + \s_{\ell_T'} \gote_{j_{\ell_T'}} + 
y_{d-1} \s_{\ell_T'} ( \nu_{\ell_T'} (\und{y}_{\ell_T'}) - \gote_{j_{\ell_T'}}),
\end{aligned}
\]
where $\nu_\ell^0(T')=\nu_\ell^0(T)-\gote_j+\gote_{j'}$.
\end{itemize}
This concludes the proof of \eqref{nubemolle}.
Thus either both summands in \eqref{diffe} are bounded by \eqref{maximum}
%%
%\begin{equation} \nonumber
%A^{|N(\vartheta)|}K_1^{N-1} (\calJ_{\al}(\vartheta))^{2N^2}
%e^{-(s'-s_3-s_4)\calJ(\vartheta)} e^{- s_4 2^{-\al}\min\{\jap{j}^{\al/2},\jap{j'}^{\al/2}\}}
%\end{equation}
%
or, by collecting together the identities \eqref{bemolle} and \eqref{nubemolle}, one finds that
$C_q(\ze,\ell,j')=C_q(\ze,\ell,j)$ for all $q=0,1,\ldots,N-1$ 
and hence
\[
(j')^{q_\ell} \gotg_{\ell}(\ze,\pow_\ell,j')=j^{q_\ell} \gotg_{\ell}(\ze,\pow_\ell,j) 
\]
for any active line $\ell \in \calmP(\vartheta)$, while the propagator
of any neutral line  $\ell\in\calmP(\vartheta)$ does not change when considered as a line in $\vartheta'$.
\end{itemize}
Therefore, we deduce that whenever \eqref{diffe} does not hold, then
$|j^qc_j^{-1} {\overline{\Val}^\ttR} (\vartheta;c,\ze)|$ and $|(j')^qc_{j'}^{-1} {\overline{\Val}^\ttR} (\vartheta';c,\ze)|$
are exponentially small in $\min\{\jap{j}^{\al/2},\jap{j'}^{\al/2}\}$ according to \eqref{maximum}.
\EP
%%%%%%%%%%%%%%%%%%%%%%%%%%%%%%%%%%%%%%%%%%%%%%%%%%%%%%%%%%%%%%%%%%%%%%%%%% 

%%%%%%%%%%%%%%%%%%%%%%%%%%%%%%%%%%%%%%%%%%%%%%%%%%%%%%%%%%%%%%%%%%%%%%%%%% 
\begin{lemma}\label{squa}
For any $\pow=0,\ldots,N-1$, 
there are positive constants $C=C(s,\al,\ze,N)$, $D=D(s,\al,\ze,N)$, $\de_0$ and $\al_0$, such that for all $j,j'\in\ZZZ \setminus \{0\}$ and all $\e\in(-C^{-1},C^{-1})$
 one has
\begin{equation} \label{differenza}
\left|\gota_\pow(c,\ze,\e,j)-\gota_\pow(c,\ze,\e,j')\right| < D |\e| e^{-\de_0\min\{\jap{j}^{\al_0},\jap{j'}^{\al_0}\}} ,
\end{equation}
and hence there exists
\[
\gota_{\pow}(c,\ze,\e):=\lim_{|j|\to+\io}\gota_{\pow}(c,\ze,\e,j) .
\]
\end{lemma}
%%%%%%%%%%%%%%%%%%%%%%%%%%%%%%%%%%%%%%%%%%%%%%%%%%%%%%%%%%%%%%%%%%%%%%%%%% 

%%%%%%%%%%%%%%%%%%%%%%%%%%%%%%%%%%%%%%%%%%%%%%%%%%%%%%%%%%%%%%%%%%%%%%%%%% 
\prova
Fix $s'$ and $s_3$ as in the proof of Proposition \ref{voltabuona} with $s_1=0$, and $s_4=2s/5$.
Then $\bar{s}=s/5$, and the bound \eqref{differenza} follows from Lemma \ref{presqua}, with $\de_0=s/5$ and $\al_0=\al/2$.
This proves that $\{\gota_\pow(c,\ze,\e,j)\}_{j\in\ZZZ}$ is a Cauchy sequence.
%so that the second assertion follows immediately.
\EP
%%%%%%%%%%%%%%%%%%%%%%%%%%%%%%%%%%%%%%%%%%%%%%%%%%%%%%%%%%%%%%%%%%%%%%%%%% 

%%%%%%%%%%%%%%%%%%%%%%%%%%%%%%%%%%%%%%%%%%%%%%%%%%%%%%%%%%%%%%%%%%%%%%%%%% 
\begin{lemma}\label{auno}
Let $C$, $D$, $\al_0$ and $\de_0$ be the constants appearing in Lemma \ref{squa}.
If $N\ge2$ one has
\[
|\gota_1(c,\ze,\e,j) |\le D {\e^2 } e^{-\de_0 \jap{j}^{\al_0}},
\]
for all $\e\in(-C^{-1},C^{-1})$
\end{lemma}
%%%%%%%%%%%%%%%%%%%%%%%%%%%%%%%%%%%%%%%%%%%%%%%%%%%%%%%%%%%%%%%%%%%%%%%%%% 

%%%%%%%%%%%%%%%%%%%%%%%%%%%%%%%%%%%%%%%%%%%%%%%%%%%%%%%%%%%%%%%%%%%%%%%%%%
\prova
Consider a tree $\vartheta\in\gotF\gotH^{(k)}_{j,+}(1)$.
First of all we note that, since $q=1$, $\vartheta$ contains one and only one active line $\ell$.
Let $\vartheta_*$ be the tree obtained from $\vartheta$ by exchanging the root line $r_\vartheta$
with the special leaf $\la_\vartheta$ and, if $\ell$ is along a path $\calP_{\breve\vartheta'}$, with $\vartheta'\neq\vartheta$,
by exchanging also  $r_{\vartheta'}$ with $\la_{\vartheta'}$. This means that we reverse the orientation
of all the lines in $\calP_{\breve\vartheta}\cup\{\ell_\vartheta,\ell_{\la_{\vartheta}}\}$,
with $\ell_{\la_{\vartheta}}$ becoming the root line of $\vartheta_*$, and if $\ell$ is along the path $\calP_{\breve\vartheta'}$,so as to obtain a new tree $\vartheta_*'$ in such a way that $\ell_{\la_{\vartheta'}}$ becomes the root line of $\vartheta_*'$, 
and graft $\vartheta_*'$ again to the same node of $\vartheta$ which $\ell_{\vartheta'}$ entered.
The labels of $\vartheta_*$ are changed accordingly so as to satisfy \eqref{lefoglie*}.
Note that, as in the proof of Lemma \ref{squa}, $\vartheta_*$ may fail to exists,
but in such a case $\ol{\Val}^{\ttR}(\vartheta;c,\ze)$ has at least a factor $e^{-2s\jap{j}^\al}$.

The propagator of the line $\ell$, regarded as a line in $\vartheta$, according to \eqref{mava1}, is
we also cut the subtree $\vartheta'$ from $\vartheta$, reverse the orientation of 
all the lines in $\calP_{\breve\vartheta'}\cup\{\ell_{\vartheta'},\ell_{\la_{\vartheta'}}\}$ 
\begin{equation} \label{basta}
\gotg_\ell(\ze,1,j):=
\frac{1}{j} \frac{1}{(2 ( j_{\ell_{\vartheta'}}^\flat - j_\ell^\flat ))} .
\end{equation}
We distinguish between two cases:
\begin{itemize}[topsep=0ex]
\itemsep0em
\item If $\ell\in\calP(\vartheta)$, so that $\vartheta'=\vartheta$,
we write $j=j_{\ell_{\vartheta}}= \gotj + \s_\ell j_{\ell}$ and $\s_\ell j_{\ell}= \gotj' + j_{\la_{\vartheta}}$,
which implicitly define $\gotj$ and $\gotj'$, so that, using that $j_{\la_{\vartheta}}=j$ implies $\gotj'=-\gotj$. In that case
one has $j_{\ell_{\vartheta}}^\flat=j_{\ell_{\vartheta_*}}^\flat=0$, while
$j_\ell^\flat = \gotj'$ in $\vartheta$ and $j_\ell^\flat = \gotj=-\gotj'$ in $\vartheta_*$. Then the
propagator of $\ell$ in $\vartheta_*$ differs from \eqref{basta} only by a sign. 
\item If $\ell\in\calP(\vartheta')$, with $\vartheta' \neq \vartheta$, then 
$j=\gotj+\s_{\ell_{\vartheta'}} j_{\vartheta'}$ and
$\s_{\ell_{\vartheta'}} j_{\vartheta'}= \gotj'+j$, and, similarly,
$\s_{\ell_{\vartheta'}} j_{\ell_{\vartheta'}}=\gotj''+\s_\ell j_{\ell}$
and $\s_\ell j_\ell=\gotj''' + \s_{\ell_{\vartheta'}} j_{\ell_{\vartheta'}}$, 
which implicitly define $\gotj$, $\gotj'$, $\gotj''$ and $\gotj'''$. Therefore one deduces that $\gotj'=-\gotj$ and $\gotj'''=-\gotj''$.
Moreover, in $\vartheta$ one has $j_{\ell_{\vartheta'}}^\flat=\gotj'$ and $j_\ell^\flat=\gotj'''+ \s_{\ell_{\vartheta'}} j_{\ell_{\vartheta'}} - j= 
\gotj'''+\gotj'$, whereas in $\vartheta_*$ one has
$j_{\ell_{\vartheta'}}^\flat=\gotj=-\gotj'$ and $j_\ell^\flat=\gotj''+ \s_{\ell_{\vartheta'}} j_{\ell_{\vartheta'}} - j= \gotj''+\gotj=-\gotj'''-\gotj'$.
This means that once more the propagator just change sign.
\end{itemize}
The same argument shows that, for all the other lines $\ell$, the denominators
$x_\ell(\und{y}_{\ell})=\om\cdot (\nu_{\ell}(\und{y}_\ell)-\gote_{j_\ell})$
do not change when regarded as lines of $\vartheta$ or $\vartheta_*$.
Summarizing, if $\vartheta_*$ exists, one has
$$
\ol{\Val}^{\ttR}(\vartheta_*;c,\ze)=-\ol{\Val}^{\ttR}(\vartheta;c,\ze).
$$
Hence the assertion follows.
\EP
%%%%%%%%%%%%%%%%%%%%%%%%%%%%%%%%%%%%%%%%%%%%%%%%%%%%%%%%%%%%%%%%%%%%%%%%%% 

We are finally ready to conclude the proof of Proposition \ref{decadeta} and derive the asymptotic expansion of the counterterms.
%
%%%%%%%%%%%%%%%%%%%%%%%%%%%%%%%%%%%%%%%%%%%%%%%%%%%%%%%%%%%%%%%%%%%%%%%%%%% 
%\noindent {\it Proof of Proposition \ref{decadeta}}.
Proposition \ref{stimeta} and Lemmata \ref{squa} and \ref{auno} immediately give \eqref{nemmenounnome!},
%one has, for $j\ne0$,
%\[
%\begin{aligned}
%\h_j(c,\om(\ze),\e)&=\sum_{\pow=0}^{N-1}\frac{\gota_\pow(c,\ze,\e,j)}{j^\pow} + \gotA(c,\ze,\e,j) \\
%&= 
%\gota_0(c,\ze,\e)+\sum_{\pow=2}^{N-1}\frac{\gota_\pow(c,\ze,\e)}{j^\pow} + \gotr_j(c,\ze,\e)
%\end{aligned}
%\]
where, for $j\neq 0$,
\[
\gotr_j(c,\ze,\e) :=  \gotA (c,\ze,\e,j)+ 
%(\gota_0(c,\ze,\e,j)-\gota_0(c,\ze,\e)+ \frac{\gota_1(j)}{c,\ze,\e,j}+
\sum_{\pow=0}^{N-1}\frac{\gota_\pow(c,\ze,\e,j)-\gota_\pow(c,\ze,\e)}{j^\pow} , 
\qquad \gota_1(c,\ze,\e)=0,
\]
while $\gotr_0(c,\ze,\e):=  \h_0(c,\om(\ze),\e)-\gota_0(c,\ze,\e)$.
The bounds \eqref{tosse} follow from Proposition \ref{stimeta} and Lemma \ref{squa}. This concludes the proof 
of Proposition \ref{decadeta}, with $\e_1=\e_1(s,\al,\ze,N)=C^{-1}$ defined in Lemma \ref{squa}.

%%%%%%%%%%%%%%%%%%%%%%%%%%%%%%%%%%%%%%%%%%%%%%%%%%%%%%%%%%%%%%%%%%%%%%%%%% 
\subsection{Proof of Lemma \ref{carinello}}\label{provocarinello}
%%%%%%%%%%%%%%%%%%%%%%%%%%%%%%%%%%%%%%%%%%%%%%%%%%%%%%%%%%%%%%%%%%%%%%%%%%

Consistently with the notation \eqref{media} we set
\[
\av{|U|^4}_{\TTT^\ZZZ\times \TTT} = \sum_{k=0}^{\io} \e^k
\sum_{\substack{ \nu_1 -\nu_2 + \nu_3-\nu_4=0 \\ j_1-j_2+j_3 - j_4 = 0  \\ k_1+k_2+k_3+k_4=k}}
u_{j_1,\nu_1}^{(k_1)}\ol{u}_{j_2,\nu_2}^{(k_2)}u_{j_3,\nu_3}^{(k_3)}\ol{u}_{j_4,\nu_4}^{(k_4)} ,
\]
where $U$ is as in \eqref{converga} and we used the convention
\[
u^{(0)}_{j,\nu}=
\begin{cases}
c_j,\qquad \nu=\gote_j,& \\
0,\qquad \nu\ne\gote_j.&
\end{cases}
\]

%
%%%%%%%%%%%%%%%%%%%%%%%%%%%%%%%%%%%%%%%%%%%%%%%%%%%%%%%%%%%%%%%%%%%%%%%%%%% 
%\begin{lemma}\label{carinello}
%Let $\gota_0(c,\ze,\e)$ as in Proposition \ref{decadeta}, be defined as in Lemma \ref{squa}.
%One has $\gota_0(c,\ze,\e) = - 3\,\e\av{|U|^4}_{\TTT\times\TTT^\ZZZ}$.
%\end{lemma}
%%%%%%%%%%%%%%%%%%%%%%%%%%%%%%%%%%%%%%%%%%%%%%%%%%%%%%%%%%%%%%%%%%%%%%%%%%% 

%%%%%%%%%%%%%%%%%%%%%%%%%%%%%%%%%%%%%%%%%%%%%%%%%%%%%%%%%%%%%%%%%%%%%%%%%% 
We prove that there are constants $C_0,\de_0,\al_0>0$ such that
\begin{equation}\label{a}
\left| \gota_0(c,\ze,\e,j) + 3\e \av{|U|^4}_{\TTT\times\TTT^\ZZZ}\right| \le C_0{\e^2 }e^{-\de_0 \jap{j}^{\al_0}} ,
\end{equation}
where $\gota_0(c,\ze,\e,j)$ is defined in \eqref{momovedo}; recall that by Remark \ref{uguali}, 
the series in \eqref{gliaa} defining ${\gota}_0^{(k)}(c,\ze,j)$ is absolutely convergent.

%where, according to \eqref{gliaa}, $\gota_0(c,\ze,\e,j)$ is given by
%%
%\begin{equation} \nonumber
%\gota_0(c,\ze,\e,j)= \sum_{k\ge1} \e^k {\gota}_0^{(k)}(c,\ze,j),\qquad
%{\gota}_0^{(k)}(c,\ze,j) =-\frac{1}{c_j}\sum_{\vartheta\in {\gotH}^{(k)}_{j,+}(0)}\!\!\!\!\! \overline{\Val}(\vartheta;c,\ze) .
%\end{equation}
%%

By using \eqref{converga} and \eqref{settou} in Proposition \ref{bellaprop}, we can express the average
 $3\,\e\av{|U|^4}_{\TTT\times \TTT^\ZZZ}$
as a sum of values of renormalized trees, in the form
\begin{equation}\label{path0}
\begin{aligned}
3\,\e\av{|U|^4}_{\TTT\times \TTT^\ZZZ} &= 3\sum_{k\ge1}\e^k 
\sum_{\substack{k_1+k_2+k_3+k_4=k \\ j_1+j_2+j_3+j_4=0 \\ \nu_1+\nu_2+\nu_3+\nu_4=0 \\ \s_11+\s_21+\s_31+\s_41=0}}
\prod_{i=1}^4\Bigg(\sum_{\substack{\vartheta\in \gotF^{(k_i)}_{j_i,\nu_i,\s_i} \\ \calP_\vartheta=\emptyset }}\Val^{\ttR}(\vartheta;c,\om(\ze))\Bigg)\\
&=\frac{1}{c_j} \sum_{k\ge1}\e^k
\sum_{\substack{\vartheta\in \gotF^{(k)}_{j,\gote_j,+} \\ \calP_\vartheta=\emptyset }}\Val^{\ttR}(\vartheta;c,\om(\ze)) =
\frac{1}{c_j} \sum_{k\ge1}\e^k 
\sum_{\substack{\vartheta\in\gotH^{(k)}_{j,+}(0) \\ \calP_{\vartheta}=\emptyset }} \ol{\Val}(\vartheta;c,\ze) ,
\end{aligned}
\end{equation}
where the second equality follows Lemma \ref{effone} and by the observation that, with the notation of Definition \ref{ramification},
$\calmP(\vartheta)=\emptyset$ if $\calP_\vartheta=\emptyset$, so that any tree $\vartheta\in\gotH^{(k)}_{j,+}(0)$ with
$\calP_{\vartheta}=\emptyset$ is also an expanded tree in $\gotT^{(k)}_{j,\gote_j,+}$.

On the other hand, thanks to Remark \ref{uguali} and \eqref{solos}, and using that for any tree
$\vartheta\in\gotH^{(k)}_{j,+}(0)$ necessarily one has $L_{\bf{a}}(\vartheta)=\emptyset$,
we can write
\[
\sum_{\substack{\vartheta\in\gotH^{(k)}_{j,+}(0) \\ \calP_{\vartheta}\ne\emptyset }} \ol{\Val}(\vartheta;c,\ze)=
\sum_{\substack{\vartheta\in\gotH^{(k)}_{j,+}(0,{\bf{n}}) \\ \calP_{\vartheta}\ne\emptyset}} \ol{\Val}(\vartheta;c,\ze) +
\sum_{\substack{\vartheta\in \gotF\gotH^{(k)}_{j,+}(0,{\bf{s}}) \\ \calP_\vartheta\ne\emptyset}}\ol{\Val}^\ttR(\vartheta;c,\ze),
\]
where $\gotH^{(k)}_{j,+}(0,{\bf{n}})$ denotes the set of trees $\vartheta\in\gotH^{(k)}_{j,+}(0)$ with only neutral lines.

Since Corollary \ref{small} guarantees that
\[
\sum_{\vartheta\in \gotF\gotH^{(k)}_{j,+}(0,{\bf{s}})}|\ol{\Val}^\ttR(\vartheta;c,\ze)| \le C^k e^{-s_42^{-\al} \jap{j}^{\al/2}} ,
\]
for some positive constant $C$, the bound \eqref{a} follows if we show that
\begin{equation}\label{mahmah}
\sum_{\substack{\vartheta\in\gotH^{(k)}_{j,+}(0,{\bf{n}}) \\ \calP_{\vartheta}\ne\emptyset }} \ol{\Val}(\vartheta;c,\ze)=0.
\end{equation}

For any tree $\vartheta\in\gotH^{(k)}_{j,+}(0,{\bf{n}})$ such that $\calP_\vartheta \neq \emptyset$,
write $\calP_{\vartheta}=\{(v_1,v_2),\ldots,(v_{n-1},v_{n})\}$,
with $v_1$ such that $(\la_\vartheta,v_1)=\ell_{\la_\vartheta}$ and $v_n$ such that $(v_{n},r)=\ell_\vartheta$.
For $i=1,\ldots,n$ let $\Theta_i := \{\vartheta_{i,1},\ldots,\vartheta_{i,s_{v_i}-1}\}$ be the set of the
subtrees whose root lines enter $v_i$ and do not belong to $\calP_{\vartheta}\cup\{\ell_{\la_\vartheta}\}$.
Set also $\ell_i:=(v_{i},v_{i+1})$ for $i=1,\ldots,n-1$, while $\ell_0:=\ell_{\la_{\vartheta}}$ and $\ell_n:=\ell_{\vartheta}$.
Then for each $i=1,\ldots,n$ one has
\begin{equation} \label{ciuffi}
\s_{\ell_{i}} \nu_{\ell_i}=\nu_i + \s_{\ell_{i-1}}\nu_{\ell_{i-1}} , \qquad 
\s_{\ell_{i}} j_{\ell_i}=j_i + \s_{\ell_{i-1}} j_{\ell_{i-1}} ,
\end{equation}
where
\[
\nu_{i} := \sum_{k=1}^{4} \s_{\ell_{\vartheta_{i,k}}} \nu_{\ell_{\vartheta_{i,k}}} , \qquad
j_{i} := \sum_{k=1}^{4} \s_{\ell_{\vartheta_{i,k}}} j_{\ell_{\vartheta_{i,k}}} , \qquad \qquad i=1,\ldots,n ,
\]
if $s_{v_i}=5$ and $\nu_{i} =0$, $j_i=0$ if $s_v=2$.

Since all lines $\ell_1,\ldots,\ell_{n-1}$ are neutral one has $j_{\ell_i}=j$ and $\s_{\ell_i}=+$ for all $i=1,\ldots,n$, so that
\eqref{ciuffi} yields 
\[
\nu_{\ell_i}=\nu_i + \nu_{\ell_{i-1}} , \qquad j_{\ell_i} = j_{\ell_{i-1}} , \qquad j_i = 0 ,
\qquad \qquad i=1,\ldots,n .
\]
As a consequence, if we set $x_i:=\om\cdot\nu_i$ for $i=1,\ldots,n$, we obtain
\[
x_{\ell_i} = \om \cdot \nu_{\ell_i} -\om_{j_{\ell_i}} = \om \cdot (\nu_{\ell_i} - \gote_{j_{\ell_i}} )  =
\om \cdot (\nu_{\ell_{i-1}} - \gote_{j_{\ell_{i-1}}} ) + \om\cdot\nu_i 
= x_{\ell_{i-1}} + x_i ,
\]
for all $i=2,\ldots,n-1$, and hence, by iterating, we obtain $x_{\ell_i}=x_1+\ldots + x_i$, where we have used that $x_{\ell_1}=x_1$.

In \eqref{mahmah}, for any $\vartheta\in\gotH^{(k)}_{j,+}(0,{\bf{n}}) $ such that $\calP_{\vartheta}\ne\emptyset$ we can write
\begin{equation} \label{macheneso1}
\ol{\Val}(\vartheta;c,\ze)  
= c_j \prod_{i=1}^{n} \Biggl( \overline\matG_{\ell_i}(\ze) 
\starF_{v_i}(c) \prod_{k=1}^{s_{v_i}-1} \ol{\Val}(\vartheta_{i,k};c,\ze)  \Biggr) ,
\end{equation} 
so that if we sum over all possible assignments of the scale labels of the lines $\ell_1,\ldots,\ell_n$, we obtain
\begin{equation} \label{macheneso2}
G_n(x)
\Biggl( \prod_{i=1}^{n} \Biggl( \starF_{v_i}(c) \prod_{k=1}^{s_{v_i}-1} \ol{\Val}(\vartheta_{i,k};c,\ze)  \Biggr) \Biggr) c_j ,
\end{equation} 
with $x=(x_1,\ldots,x_n)$ and using the notation \eqref{GnD}. Let us distinguish the following cases.
\begin{itemize}[topsep=0ex]
\itemsep0em
\item If $x_i \neq 0$ for all $i=1,\ldots,n$, then if we consider together all the trees which differ from $\vartheta$
because of a permutation of the sets $\Theta_i$, while leaving unchanged all the other factors in \eqref{macheneso2},
and sum all the corresponding contributions, we obtain from \eqref{macheneso2}
\begin{equation} \nonumber %\label{macheneso3}
\ol{\Val}(\vartheta;c,\ze)  = \left( \sum_{x'\in \Pi_A(x)} G_n(x') \right)
\Biggl( \prod_{i=1}^{n} \Biggl( \starF_{v_i}(c) \prod_{k=1}^{s_{v_i}-1} \ol{\Val}(\vartheta_{i,k};c,\ze)  \Biggr) \Biggr) c_j ,
\end{equation} 
which vanishes because of Lemma \ref{stoc}.
\item If $x_i=0$ for some $i=1,\ldots,n$, then either $s_{v_i}=2$  or $s_{v_i}=5$ and $\nu_i=0$. In that case,
we consider together all such contributions. In particular, if $s_{v_i}=2$, then $\vartheta_{i,1}$ has a special
leaf $\la_i$ such that $j_{\la_i} = j_{\ell_{\vartheta_{v,1}}}=j_{\ell_i}=j$, where the latter identity follows
from the fact that the line $\ell_i$ is neutral. All the contributions such that $\calP_{\vartheta_{v_i,1}} =\emptyset$,
i.e.~such that the line $\ell_{\la_i}$ enter the node which the root line $\ell_{\ell_{\vartheta_{v_i,1}}}$ exits,
cancel out the contributions with $s_{v,i}=5$, again by \eqref{path0}.
\item Therefore we are left with trees such that there is at least one node $v_i$ with $s_{v_i}=2$ and
$\calP_{\vartheta_{v,1}} \neq \emptyset$. In that case, we consider $\ol{\Val}(\vartheta_{i,1};c,\ze)$:
by construction $\vartheta_{i,1}\in\gotH^{(k')}_{j,+}(0,{\bf{n}})$ for some $k'<k$ and $\calP_{\vartheta_{i,1}}\ne\emptyset$.
Therefore we can study $\vartheta_{i,1}$ as we have studied $\vartheta$, and we find that the only contributions
which do not vanish when summed together are those that correspond to trees containing at least one node $v'$
along the special path $\calP_{\vartheta_{i,1}}$ such that $s_{v'}=2$ and the subtree $\vartheta'$
grafted to $v'$ is such that $\calP_{\vartheta'} \neq \emptyset$.
\item We iterate the construction above: at each step, we have to deal with a tree
which belongs to $\gotH^{(k')}_{j,+}(0,{\bf{n}})$ for some $k'<k$ and has a non-empty special path.
However, at any step the order $k'$ decreases, so that at some point one reaches the value $k'=1$,
so that it is no longer possible to have a non-empty path.
\end{itemize}
In conclusion the bound \eqref{a} follows with $\al_0=\al/2$
and $\de_0=s_02^{-\al}$.
%%%%%%%%%%%%%%%%%%%%%%%%%%%%%%%%%%%%%%%%%%%%%%%%%%%%%%%%%%%%%%%%%%%%%%%%%% 

%%%%%%%%%%%%%%%%%%%%%%%%%%%%%%%%%%%%%%%%%%%%%%%%%%%%%%%%%%%%%%%%%%%%%%%%%% 
%%%%%%%%%%%%%%%%%%%%%%%%%%%%%%%%%%%%%%%%%%%%%%%%%%%%%%%%%%%%%%%%%%%%%%%%%% 
\section{Lipschitz extensions and measure estimates}
\label{provafinale}
\zerarcounters
%%%%%%%%%%%%%%%%%%%%%%%%%%%%%%%%%%%%%%%%%%%%%%%%%%%%%%%%%%%%%%%%%%%%%%%%%% 
%%%%%%%%%%%%%%%%%%%%%%%%%%%%%%%%%%%%%%%%%%%%%%%%%%%%%%%%%%%%%%%%%%%%%%%%%% 

%%%%%%%%%%%%%%%%%%%%%%%%%%%%%%%%%%%%%%%%%%%%%%%%%%%%%%%%%%%%%%%%%%%%%%%%%% 
%%%%%%%%%%%%%%%%%%%%%%%%%%%%%%%%%%%%%%%%%%%%%%%%%%%%%%%%%%%%%%%%%%%%%%%%%% 

We are finally ready to combine Theorem \ref{moser} and Proposition \ref{decadeta} with the compatibility equation \eqref{perv}, so as
to obtain the proof of Theorem \ref{main}.

%%%%%%%%%%%%%%%%%%%%%%%%%%%%%%%%%%%%%%%%%%%%%%%%%%%%%%%%%%%%%%%%%%%%%%%%%% 
\subsection{Lipschitz bounds and proof of Proposition \ref{fabrizio}}\label{lip}
%%%%%%%%%%%%%%%%%%%%%%%%%%%%%%%%%%%%%%%%%%%%%%%%%%%%%%%%%%%%%%%%%%%%%%%%%% 

Proposition \ref{bellaprop} ensures, for any fixed $\om\in\gotB^{(0)}$, the existence of both the solution $U$ and the counterterm $\h$ as 
convergent series. Moreover, as Proposition \ref{decadeta} shows, for $\ze=(\ka,\x)\in\matK_N$, with $N\ge1$, and $\om$ as in \eqref{espome},
each component $\h_j$ with $j\ne0$ of the counterterm can be expanded in inverse powers of $j$ as well as $\om_j$.
This allows us, in principle, to make \eqref{perv} to be satisfied with $V\in \ell^{N,\io}$ by solving the implicit function problem given by \eqref{system}.
However, as already pointed out (see Remark \ref{trimbo}), the radius of convergence depends on $\ze$, and $\ze$ varies in a set which does not 
contain any open set. 

Therefore, in order to conclude the proof of Theorem \ref{main}, the first task is to prove Proposition \ref{fabrizio}, with the aim
of obtaining a Lipschitz extension of the solution $U$, the counterterm $\h$ and the 
coefficients $\gota_\pow(c,\ze,\e)$, $\gotA(c,\ze,\e,j)$ for all $\ze$ in an open set. 
For $\g>0$ let $\matK_N(\g)$ be the set of good parameters as in Definition \ref{unibr}. Then
the proof of Proposition \ref{fabrizio} goes through two steps: for $\ze\in\matK_N(\g)$ we obtain, 
first, uniform bounds for the renormalized values of the trees,
next, we show that they verify Lipschitz estimates.

%The following results hold for  all $c\in\matU_1(\mathtt{g}(s,\al))$.
%%%
%%%Below we study explicitly the case $N\ge 2$. The case $N=1$ is easier and can be discussed exactly in the same way,
%%%up to notation simplifications: essentially, in what follows, $(\ze)$ has to be replaced with $\xi$ and there are no coefficients
%%%$\gota_q$ to be considered.

%%%%%%%%%%%%%%%%%%%%%%%%%%%%%%%%%%%%%%%%%%%%%%%%%%%%%%%%%%%%%%%%%%%%%%%%%% 
\begin{lemma}\label{uniforme}
For all $\g>0$ let $\matK_N(\g)$ be defined as in \eqref{bgamma}.
For all {$s_1 \ge 0$ and $s_2>0$ such that $s_1+s_2 = s$}, all $s'\in(s_1,s)$ and all $s_3\in [s_1,s')$, 
there is a positive constant {$A_*=A_*(s-s',s_3,\al,\g)$} such that
\begin{enumerate}[topsep=0ex]
\itemsep0em
\item
for all $k\ge1$, $j\in\ZZZ$, $\nu\in\ZZZ^\ZZZ_f$, $\s\in\{\pm\}$ and for any $\vartheta \in \gotF^{(k)}_{j,\nu,\s}$, one has %, if $\nu\neq\gote_j$,
\begin{subequations} \label{diobono+nonchanomestima}
\begin{align}
%\begin{equation}
\sup_{\ze\in \matK_N(\g)}
 \! | \Val^\ttR (\vartheta;c,\om(\ze)) |
& \le
A_*^{|N(\vartheta)|} 
e^{- s_1 |\nu|_\al} 
e^{- s_2 \jap{j}^\alpha}
e^{-(s'-s_3) \calJ(\vartheta)} , 
\quad \nu\neq\gote_j ,
\label{diobono} \\
%\end{equation}
%
%while, if $\nu=\gote_j$, one has 
%
%\begin{equation} \label{nonchanomestima}
\sup_{\ze\in \matK_N(\g)} \left| \frac{1}{c_j}\Val^\ttR (\vartheta;c,\om(\ze)) \right|
& \le
A_*^{|N(\vartheta)|} 
%e^{- s_2 \jap{j}^\alpha}
e^{-(s'-s_3) \calJ(\vartheta)} ,
\quad \qquad \quad\, \nu=\gote_j ,
\label{nonchanomestima}
\end{align}
\end{subequations}
%\end{equation}
%
\item
for all $k\ge1$, $j\in\ZZZ\setminus\{0\}$, $q=0,\ldots, N$ and for any $\vartheta \in \gotF\gotH^{(k)}_{j,+}(q)$, one has
\begin{subequations} \label{anchelei+ancheleianche}
\begin{align}
%\begin{equation} \label{anchelei}
\null\hspace{-1cm}
\sup_{\ze\in \matK_N(\g)} \left| \frac{j^q}{c_j}\ol{\Val}^\ttR (\vartheta;c,\ze) \right| 
& \le
A_*^{|N(\vartheta)|} K_0^N (\calJ_{\al}(\vartheta))^{2(q-1)}
e^{-(s'-s_3) \calJ(\vartheta)} ,
 \quad q=0,\ldots,N-1,
\label{anchelei} \\
%\end{equation}
%
%for $q=0,1,\ldots,N-1$, and
%
%\begin{equation} \label{ancheleianche}
\null\hspace{-1cm}
\sup_{\ze\in \matK_N(\g)} \left| \frac{j^N}{c_j}\ol{\Val}^\ttR (\vartheta;c,\ze) \right|
& \le
A_*^{|N(\vartheta)|} K_1^{|N(\vartheta)|}
 (\calJ_{\al}(\vartheta))^{2N^2}e^{-(s'-s_3)\calJ(\vartheta)} , \qquad\hspace{0.6cm} q=N,
\label{ancheleianche}
\end{align}
\end{subequations}
%\end{equation}
%
where $K_0$ and $K_1$ are as in the proof of Lemma \ref{converge2}.
\end{enumerate}
\end{lemma}
%%%%%%%%%%%%%%%%%%%%%%%%%%%%%%%%%%%%%%%%%%%%%%%%%%%%%%%%%%%%%%%%%%%%%%%%%% 

%%%%%%%%%%%%%%%%%%%%%%%%%%%%%%%%%%%%%%%%%%%%%%%%%%%%%%%%%%%%%%%%%%%%%%%%%% 
\prova
If $\ze\in \matK_N(\g)$ then $\om=\om(\ze)\in\gotB(\gamma,N+1)$, so we may 
take $r\in\gotR$, $r^*\in\gotR^*$ such that 
$\be_\om^{(0)}(r_m)\ge\be^*(r_m^*,\g)$  and $r_m\ge r_m^*$ for all $m\ge1$. 
Therefore for any 
$\ze\in \matK_N(\g)$ we may bound in \eqref{telanumero}
 $$
\frac{1}{r_{m_n-1}}\log\left(\frac{1}{ \be(m_n)}\right) \le \frac{1}{r^*_{m_n-1}}\log\left(\frac{1}{ \be^*(r_{m_n}^*,\g)}\right).
 $$
 This implies that the bounds in  Lemma \ref{converge1} hold
with
\[
A^*(s-s',s_3,\al,\g):= A_2^* (n_3(\min\{s-s',s_3\}))
\]
where
\begin{equation}\label{76}
A^*_2(n):= \frac{2^{24}a_0^2}{\be(r^*_{m_n})^6} = \Biggl( \frac{2^4 a_0^{1/3}}{\be^*(r^*_{m_n},\g)} \Biggr)^6 ,
\end{equation}
 and $n_3(\de)$ is the smallest integer so that 
\[
\sum_{n\ge n_3(\de)+1} \frac{1}{r^*_{m_n-1}}\log\left(\frac{1}{ \be^*(r_{m_n}^*,\g)}\right) \le \de.
\]
In this way the bounds are all uniform in $\om \in \gotB(\gamma,N+1)$
and hence uniform in $\ze\in \matK_N(\g)$. This proves \eqref{diobono} and \eqref{nonchanomestima}. Regarding
\eqref{anchelei} and \eqref{ancheleianche} one 
reasons in the same way, but relying on the proof of Lemma \ref{converge2}.
 \EP
 %%%%%%%%%%%%%%%%%%%%%%%%%%%%%%%%%%%%%%%%%%%%%%%%%%%%%%%%%%%%%%%%%%%%%%%%%% 
 
We now discuss  the Lipschitz regularity. First, we need some notation.
Introduce an ordering on $L_0(\vartheta):=\{\ell\in L(\vartheta) : n_{\ell} \ge 0 \}$ through the following iterative procedure:
\begin{enumerate}[topsep=0ex]
\itemsep0em
\item 
set $\ell_1=\ell_{\vartheta}$, 
\item given $\ell_1,\ldots,\ell_i$, with $1<i<|L_0(\vartheta)|$, set $\LLL_i(\vartheta):=\{\ell_1,\ldots,\ell_i\}$ and,
\begin{itemize}[topsep=0ex]
\itemsep0em
\item[2.1.] if $\ell_{i}$ is not contained in any relevant RC, take as $\ell_{i+1}$ any line in $L_0(\vartheta)$ connected to $\LLL_i(\vartheta)$,
\item[2.2.] if $T$ is the relevant RC with highest depth containing $\ell_i$ and $\LLL_i(\vartheta)\setminus L(T) \neq \emptyset$,
take as $\ell_{i+1}$ any line in $L_0(\vartheta)$ contained in $T$ connected to $\LLL_i(\vartheta)$;
\item[2.3.] if $T$ is the relevant RC with highest depth containing $\ell_i$ and $\LLL_i(\vartheta) \setminus L(T) =\emptyset$,
take $\ell_{i+1}=\ell_{T}'$.
\end{itemize}
\end{enumerate}
Set $\LLL_*(\vartheta):=\LLL_{|L_0(\vartheta)|}(\vartheta)$: by construction, $\LLL_*(\vartheta)$ contains the same lines as $L_0(\vartheta)$, the
only difference being that the lines of $\LLL_*(\vartheta)$ are ordered.
Any line $\ell_i\in\LLL_{*}(\vartheta)$ induces a natural splitting of $L_0(\vartheta)$ into three disjoint sets
$\LLL_{i-1}(\vartheta)$, $\{\ell_i\}$ and $\LLL_i^c(\vartheta):=L_*(\vartheta)\setminus\LLL_i(\vartheta)$.

\begin{figure}[H]
%\vspace{-.2cm}
\centering
%\null
%\hspace{-.6cm}
\ins{028pt}{-132pt}{$A$}
\ins{402pt}{-120pt}{$B$}
\ins{325pt}{-218pt}{$C$}
\ins{220pt}{-030pt}{$T$}
\ins{350pt}{-010pt}{$T'$}
\ins{367pt}{-050pt}{$T''$}
\ins{050pt}{-118pt}{$\ell_{i-4}$}
\ins{090pt}{-083pt}{$\ell_{i-3}$}
\ins{085pt}{-128pt}{$\ell_{i-2}$}
\ins{127pt}{-148pt}{$\ell_{i-1}$}
\ins{156pt}{-074pt}{$\ell_{i}$}
\subfigure{\includegraphics*[width=6.0in]{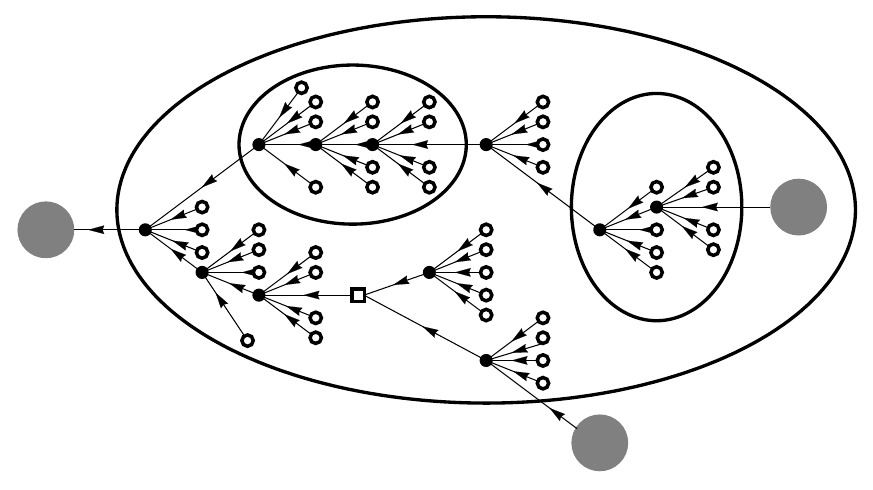}}
\caption{\small Examples of sets $\LLL_i(\vartheta)$ and lines $\ell_{i}$.}
\label{maiunaltrocosi}
\end{figure}

See Figure \ref{maiunaltrocosi} for an example. Assume that $T$, $T'$ and $T''$ are all RCs, and that
the set $\LLL_{i-5}(\vartheta)$ contains all the lines $\ell\in L_0(\vartheta)$ 
of the subgraph of the tree $\vartheta$ represented by the grey circle $A$ on the left. The sets $\LLL_{i-4}(\vartheta),\ldots,\LLL_{i}(\vartheta)$
are formed starting from $\LLL_{i-5}(\vartheta)$ and adding the lines $\ell_{i-4},\ldots,\ell_i$ one after the other.
Then, according to the constraints above, the line $\ell_{i+1}$ is necessarily the line exiting a node
and entering the node which $\ell_i$ exits, and the line $\ell_{i+2}$ is the entering line of $T$.
The following sets $\LLL_{i'}(\vartheta)$, with $i'\ge i+3$, are constructed in such a way that all the lines contained inside $T'$,
including the lines of the subgraph of $\vartheta$ represented by the grey circle $B$, must be added before
considering any line in the subgraph of $\vartheta$ represented by the grey circle $C$. In particular, once a set
$\LLL_{i''}(\vartheta)$ is formed, for some $i''\ge i+4$, by adding the ine exiting a node inside $T''$,
then the line $\ell_{i''+1}$ must be the line entering $T''$. Note that, while the set $\LLL_{i-1}(\vartheta)$ is connected,
the set $\LLL_i^c(\vartheta)$ in general is not.

%%%%%%%%%%%%%%%%%%%%%%%%%%%%%%%%%%%%%%%%%%%%%%%%%%%%%%%%%%%%%%%%%%%%%%%%%%
% FIGURA 19-2
%%%%%%%%%%%%%%%%%%%%%%%%%%%%%%%%%%%%%%%%%%%%%%%%%%%%%%%%%%%%%%%%%%%%%%%%%%
%%%%%%%%%%%%%%%%%%%%%%%%%%%%%%%%%%%%%%%%%%%%%%%%%%%%%%%%%%%%%%%%%%%%%%%%%%

%%%%%%%%%%%%%%%%%%%%%%%%%%%%%%%%%%%%%%%%%%%%%%%%%%%%%%%%%%%%%%%%%%%%%%%%%% 
\begin{lemma}\label{incartate}
For any  $j\in\ZZZ$, $\nu\in\ZZZ^\ZZZ_f$, $\s\in\{\pm\}$, $k\ge1$ and any $\vartheta\in \gotF^{(k)}_{j,\nu,\s}$, the function $\Val^\ttR(\vartheta;c,\om(\cdot))$
%$\tilde{\gota}^{(k)}_q(c,\cdot,j)$ and $\tilde{\gotA}^{(k)}(c,\cdot,j)$ are
is continuous w.r.t.~the product topology. Moreover 
for  all $s_1\ge0$ and $s_2>0$ such that $s_1+s_2=s$, all $s'\in(s_1,s)$ and $s_3\in [s_1,s')$,
there is $D_*=D_*(s-s',s_3,\al,\g)$
such that for all $\vartheta\in \gotF^{(k)}_{j,\nu,\s}$ one has
\begin{equation}\label{conla}
\sup_{\substack{\ze,\ze'\in\matK_N(\g)\\ \ze\ne\ze' }}\!\!\!\!
\frac{ |\Val^\ttR(\vartheta;c,\om(\ze))-\Val^\ttR(\vartheta;c,\om(\ze'))| }{\|\ze-\ze'\|_{\io}} \le 
D_*^{|N(\vartheta)|} e^{-s_1|\nu|_\al}e^{-s_2\jap{j}^\al} e^{-(s'-s_3) \calJ(\vartheta)} .
\end{equation}
%
%Same for the series \eqref{nonchanome}, and the second series in \eqref{glia} and \eqref{sisi} with  $\gotB(\gamma,1)$ replaced with
%$\gotB_N(\gamma)$. 
\end{lemma}
%%%%%%%%%%%%%%%%%%%%%%%%%%%%%%%%%%%%%%%%%%%%%%%%%%%%%%%%%%%%%%%%%%%%%%%%%% 

%%%%%%%%%%%%%%%%%%%%%%%%%%%%%%%%%%%%%%%%%%%%%%%%%%%%%%%%%%%%%%%%%%%%%%%%%% 
\prova
First of all note that $\om=\om(\ze)$ is continuous in $\ze$ w.r.t.~the product topology and it is Lipschitz-continuous
w.r.t. the $\ell^\io$-norm in $\matK_N(\g)$.  
For any $\vartheta\in\gotF^{(k)}_{j,\nu,\s}$ the renormalized value $\Val^\ttR(\vartheta;c,\om)$  depends only on finitely many
components of $\om$ and hence it is trivially continuous w.r.t.~the product topology in $\gotB(\g,N+1)$.
%Thanks to the bound  \eqref{diobono},
%continuity holds also for the uniformly convergent sum over all  $\vartheta\in\gotF^{(k)}_{j,\nu,\s}$.

Regarding the Lipschitz regularity of  $\Val^\ttR(\vartheta;c,\om(\cdot))$, again  we can confine ourselves to prove that
\begin{equation}\label{apa}
\sup_{\substack{\om,\om'\in\gotB(\gamma,N+1)\\ \om\ne\om' }}\frac{ |\Val^\ttR(\vartheta;c,\om)-\Val^\ttR(\vartheta;c,\om')| }{\|\om-\om'\|_{\io}} 
\le D^{|N(\vartheta)|} e^{-s_1|\nu|_\al}e^{-s_2\jap{j}^\al}
e^{-(s'-s_3) \calJ(\vartheta)},
\end{equation}
for some constant $D$, because \eqref{conla} immediately follows from \eqref{apa} by using the Lipschitz-continuity of the function $\om(\ze)$.

By looking at \eqref{valR}, one realizes that $\Val^\ttR(\vartheta;c,\om)$ and $\Val^\ttR(\vartheta;c,\om')$ differ only in
the product of the propagators, since the product of the leaf and node factors and the integrals do not depend on the frequency vector.
We can write
\begin{equation} \label{valRdif}
\begin{aligned}
& \Bigg(\prod_{\ell\in L(\vartheta)} \!\!\! \matG_\ell^\ttR(\om)\Bigg) -  \Bigg(\prod_{\ell\in L(\vartheta)} \!\!\! \matG_\ell^\ttR(\om')\Bigg) \\
& \qquad \qquad =
\sum_{i=1}^{|L_0(\vartheta)|} 
\Bigg(\prod_{\ell\in \LLL_{i-1}(\vartheta)} \!\!\! \matG_\ell^\ttR(\om)\Bigg) 
\Big( \matG_{\ell_i}^\ttR(\om) - \matG_{\ell_i}^\ttR(\om') \Big) 
\Bigg(\prod_{\ell\in \LLL_{i}^c(\vartheta)} \!\!\! \matG_\ell^\ttR(\om')\Bigg) ,
\end{aligned}
\end{equation}
where, according to \eqref{gr},
\begin{equation} \label{grdif}
\matG_{\ell_i}^\ttR(\om)-\matG_{\ell_i}^\ttR(\om') = 
\partial^{\del_{\ell_i}} \left( \calG_{n_{\ell_i}}(x_{\ell_i}(\und{y}_{\ell_i})) - \calG_{n_{\ell_i}}(x'_{\ell_i}(\und{y}_{\ell_i})) \right) ,
\end{equation}
with $x_{{\ell_i}}(\und{y}_{\ell_i}):= \om\cdot\nu_{\ell_i}(\und{y}_{\ell_i})-\om_{j_{\ell_i}}$ and $x_{{\ell_i}}'(\und{y}_{\ell_i}):= \om'\cdot\nu_{\ell_i}(\und{y}_{\ell_i})-\om'_{j_{\ell_i}}$.
In \eqref{grdif}, without loss of generality, we may assume that $\Psi_{n_{\ell_i}}(x_{\ell_i}(\und{y}_{\ell_i}))\neq0$; then,
by recalling the definition \eqref{procesi} of $\calG_{n_{\ell_i}}(x)$, we write
\begin{equation} \label{lostrozzo}
\begin{aligned}
& \calG_{n_{\ell_i}}(x_{\ell_i}(\und{y}_{\ell_i})) - \calG_{n_{\ell_i}}(x'_{\ell_i}(\und{y}_{\ell_i})) \\
& \qquad =
 \frac{ \Psi_{n_{\ell_i}}(x_{\ell_i}(\und{y}_{\ell_i})) - \Psi_{n_{\ell_i}}(x'_{\ell_i}(\und{y}_{\ell_i}))}{x_{\ell_i}(\und{y}_{\ell_i})} +
\Psi_{n_{\ell_i}}(x'_{\ell_i}(\und{y}_{\ell_i})) \Biggr( \frac{1}{x_{\ell_i}(\und{y}_{\ell_i})} - \frac{1}{x'_{\ell_i}(\und{y}_{\ell_i})} \Biggr) ,
\end{aligned}
\end{equation}
where 
\[
\begin{aligned}
& \left| \frac{ \Psi_{n_{\ell_i}}(x_{\ell_i}(\und{y}_{\ell_i})) - \Psi_{n_{\ell_i}}(x'_{\ell_i}(\und{y}_{\ell_i}))}{x_{\ell_i}(\und{y}_{\ell_i})} \right| \\
& \qquad \qquad = \left| \frac{x_{\ell_i}(\und{y}_{\ell_i}) -x'_{\ell_i}(\und{y}_{\ell_i})  }{x_{\ell_i}(\und{y}_{\ell_i})}
\int_0^1 \!\! \der t \, \partial \Psi_{n_{\ell_i}}(x'_{\ell_i}(\und{y}_{\ell_i})+t( x_{\ell_i}(\und{y}_{\ell_i}) - x'_{\ell_i}(\und{y}_{\ell_i})) \right| \\
& \qquad \qquad \qquad \qquad
\le \frac{C_1}{(\be^*(r_{m_{n_{\ell_i}}}^*,\g))^2} \left| x_{\ell_i}(\und{y}_{\ell_i}) -x'_{\ell_i}(\und{y}_{\ell_i})  \right| ,
\end{aligned} 
 \]
for some positive constant $C_1$, where we are using that $\om,\om'\in \gotB(\gamma,N+1)$.
On the other hand the second contribution either vanishes, if $\Psi_{n_{\ell_i}}(x'_{\ell_i}(\und{y}_{\ell_i}))=0$, or can be bounded as
\[
\begin{aligned}
& \left| \Psi_{n_{\ell_i}}(x'_{\ell_i}(\und{y}_{\ell_i})) \Biggr( \frac{1}{x_{\ell_i}(\und{y}_{\ell_i})} - \frac{1}{x'_{\ell_i}(\und{y}_{\ell_i})} \Biggr) \right| \\
& \qquad\qquad
\le \Psi_{n_{\ell_i}}(x'_{\ell_i}(\und{y}_{\ell_i})) \left| \frac{x_{\ell_i}(\und{y}_{\ell_i}) - x'_{\ell_i}(\und{y}_{\ell_i})}{x_{\ell_i}(\und{y}_{\ell_i}) \, x'_{\ell_i}(\und{y}_{\ell_i})} \right|
\le \frac{C_2}{(\be^*(r_{m_{n_{\ell_i}}}^*,\g))^2} \left| x_{\ell_i}(\und{y}_{\ell_i}) -x'_{\ell_i}(\und{y}_{\ell_i})  \right| ,
\end{aligned}
\]
for some positive constant $C_2$. More generally, by taking into account also the derivative $\del^{\del_{\ell_i}}$
and reasoning in a similar way, the difference of propagators \eqref{grdif} is proved to be bounded as
\begin{equation} \label{grdifbound}
\left| \matG_{\ell_i}^\ttR(\om)-\matG_{\ell_i}^\ttR(\om') \right| \le
\frac{C_0}{(\be^*(r_{m_{n_{\ell_i}}}^*,\g))^{2+\del_{\ell_i}}} \left| x_{\ell_i}(\und{y}_{\ell_i}) -x'_{\ell_i}(\und{y}_{\ell_i})  \right| ,
\end{equation}
for some positive constant $C_0$.

By construction, for any summand in \eqref{valRdif} with the exception of the cloud of the line $\ell_i$,
all the relevant RCs of $\vartheta$ contain only lines which either belong to the set $\LLL_i(\vartheta)$ or
belong to the set $\LLL_i^c(\vartheta)$. In particular, given any chain $\gotC$ of RCs, there is at most
one relevant RC $T$ such that $\ell_i \in L(T)$. %\cap \LLL_i(\vartheta) \neq \emptyset$ and $L(T) \cap \LLL_i^c(\vartheta) \neq \emptyset$.
If this happens, the line $\ell_i$ splits $\gotC$ into two disjoint (possibly empty) subsets $\gotC_1$ and $\gotC_2$ such that
all the lines of the RCs $T'\in \gotC_1$ belong to $\LLL_i(\vartheta)$, while
all the lines of the RCs $T''\in \gotC_2$ belong to $\LLL_i^c(\vartheta)$. Therefore,
to deal with the product of the propagators of the resonant lines connecting the RCs contained in $\gotC_1$ and in $\gotC_2$,
we can reason as in the proof of Lemma \ref{converge1}, because the propagators of all the RCs $T'\in \gotC_1$ have frequency vector $\om$
and the propagators of all the RCs $T''\in \gotC_2$ have frequency vector $\om'$.
With respect to the case considered when proving Lemma \ref{converge1} we have no gain corresponding to the two external lines
of  $T$, because one cannot exploit the cancellation mechanism of Lemma \ref{cancellazione}.

On the other hand, as already stressed, all the RCs for which there is no gain belong the same cloud, and hence,
for any fixed scale $n \ge 0$, there is at most one RC $T$ with $\overline{n}_T=n$ such that we lose an overall  gain factor
$x_{\ell_T'}(\und{y}_{\ell_T'})$ times $x'_{\ell_T'}(\und{y}_{\ell_T'})$, where the first factor corresponds to the exiting line of $T$,
whose propagator has frequency vector $\om$, while the second factor corresponds to the entering line of $T$,
whose propagator has frequency vector $\om'$.
We can take into account the propagators of the external lines of $T$ together with the entering line $\ell_\gotC'$ and the exiting line $\ell_\gotC$
of the chain containing $T$: since such lines are non-resonant and 
$\del_{\ell_\gotC}+\del_{\ell_\gotC'}+\del_{\ell_T}+\del_{\ell_T'}\le 2$ by construction, one can bound
\[
\left| \matG_{\ell_\gotC}^\ttR(\om) \, \matG_{\ell_T}^\ttR(\om) \, \matG_{\ell_T'}^\ttR(\om') \, \matG_{\ell_\gotC'}^\ttR(\om') 
\right| \le \frac{a^4_0}{(
\be^*(r_{m_{n_{\ell_\gotC}}}^*,\g))^6} .
\]

Next, observe that $|x_\ell(\und{y}_\ell) -x'_\ell(\und{y}_\ell)| \le \|\om-\om'\|_{\io}\|\nu_\ell-\gote_{j_\ell}\|_1$, with (see Remark \ref{stoper})
\[
\|\nu_\ell - \gote_{j_{\ell}} \|_1 \le  \| \nu_\ell \|_1  + 1 \le |\Lambda(\vartheta)| + 1 \le 4k + 2 .
\]
In conclusion, each summand in \eqref{valRdif} is found to admit a bound like \eqref{sudore} or \eqref{sangue},
with a power $7$ instead of $6$ in the definition \eqref{76} of $A_2^*(n)$, because of the bound \eqref{grdifbound} for the line $\ell_i$.
Thus, the bound \eqref{apa} follows.
\EP
%%%%%%%%%%%%%%%%%%%%%%%%%%%%%%%%%%%%%%%%%%%%%%%%%%%%%%%%%%%%%%%%%%%%%%%%%% 

%%%%%%%%%%%%%%%%%%%%%%%%%%%%%%%%%%%%%%%%%%%%%%%%%%%%%%%%%%%%%%%%%%%%%%%%%% 
\begin{lemma}\label{aristesso}
For any $N\ge1$, $k\ge1$, $j\in\ZZZ\setminus\{0\}$, the functions $\tilde{\gota}^{(k)}_q(c,\cdot,j)$, with $q=0,\ldots,N-1$, and $\tilde{\gotA}^{(k)}(c,\cdot,j)$ 
are continuous w.r.t.~the product topology.
Moreover for
 all $s'\in(0,s)$ and $s_3\in (0,s')$,
there is $E_*=E_*(s-s',s_3,\al,\g)$
such that for $q=0,1,\ldots,N$ and $\s=\pm$,
and for all $\vartheta\in \gotF\gotH^{(k)}_{j,\s}(q)$ one has
\begin{equation}\label{bip1}
%\begin{equation} \label{anchelei}
\null\hspace{-1cm}
\sup_{\substack{\ze,\ze'\in\matK_N(\g)\\ \ze\ne \ze' }}
\left| \frac{j^q}{c_j} \right|
\frac{ |\ol{\Val}^\ttR(\vartheta;c,\ze)-\ol{\Val}^\ttR(\vartheta;c,\ze'))| }{\|\ze-\ze'\|_{\io}}
 \le
E_*^{|N(\vartheta)|} K_0^N (\calJ_{\al}(\vartheta))^{2(q-1)}
e^{-(s'-s_3) \calJ(\vartheta)} ,
\end{equation}
for $q=0,\ldots,N-1$, and
 \\
\begin{equation} \label{bip2}
\null\hspace{-1cm}
\sup_{\substack{\ze,\ze'\in\matK_N(\g)\\ \ze\ne \ze' }}
\left| \frac{j^q}{c_j} \right|
\frac{ |\ol{\Val}^\ttR(\vartheta;c,\ze)-\ol{\Val}^\ttR(\vartheta;c,\ze'))| }{\|\ze-\ze'\|_{\io}}
 \le
E_*^{|N(\vartheta)|} K_1^{|N(\vartheta)|}
 (\calJ_{\al}(\vartheta))^{2N^2}e^{-(s'-s_3)\calJ(\vartheta)} ,
\end{equation}
%\end{equation}
%
for $q=N$,
where $K_0$ and $K_1$ are as in the proof of Lemma \ref{converge2}.

\end{lemma}
%%%%%%%%%%%%%%%%%%%%%%%%%%%%%%%%%%%%%%%%%%%%%%%%%%%%%%%%%%%%%%%%%%%%%%%%%% 

%%%%%%%%%%%%%%%%%%%%%%%%%%%%%%%%%%%%%%%%%%%%%%%%%%%%%%%%%%%%%%%%%%%%%%%%%% 
\prova
One reasons as in the proof of Lemma \ref{incartate}. 
The continuity w.r.t. the product topology is obtained as before.
Now, for any $\vartheta\in\gotF\gotH^{(k)}_{j,\s}(q)$,
one considers the difference of the products of the propagators 
\begin{equation} \label{valRdif2}
\begin{aligned}
& \Bigg(\prod_{\ell\in L(\vartheta)} \!\!\!\overline \matG_\ell^\ttR(\ze)\Bigg) -  \Bigg(\prod_{\ell\in L(\vartheta)} \!\!\! \overline \matG_\ell^\ttR(\ze')\Bigg) \\
& \qquad \qquad =
\sum_{i=1}^{|L_0(\vartheta)|} 
\Bigg(\prod_{\ell\in \LLL_{i-1}(\vartheta)} \!\!\! \overline \matG_\ell^\ttR(\ze)\Bigg) 
\Big( \overline \matG_{\ell_i}^\ttR(\ze)- \overline \matG_{\ell_i}^\ttR(\ze') \Big) 
\Bigg(\prod_{\ell\in \LLL_{i}^c(\vartheta)} \!\!\! \overline \matG_\ell^\ttR(\ze')\Bigg) ,
\end{aligned}
\end{equation}
and writes $ \overline \matG_{\ell_i}^\ttR(\ze)- \overline \matG_{\ell_i}^\ttR(\ze')
=\matG_{\ell_i}^\ttR(\om(\ze))-\matG_{\ell_i}^\ttR(\om(\ze'))$ as in \eqref{grdif},
if the line $\ell_i$ is either neutral or small, whereas, if $\ell_i$ is active, recalling that $n_{\ell_i}=0$ by definition,
one has,
by defining $x_\ell(\ze,\und{y}_\ell)$ and $x_\ell(\ze',\und{y}_\ell)$ as in \eqref{ambrogio} with $\om$ replaced with $\om(\ze)$
and $\om(\ze')$, respectively, 
\begin{equation} \label{grdifdifaa}
\ol\matG_{\ell_i}^\ttR(\ze)-\ol\matG_{\ell_i}^\ttR(\ze') = 
\Psi^{\del_{\ell_i},p_{\ell_i}}_{0}(x_{\ell_i}(\ze,\und{y}_{\ell_i}))
\gotg_{\ell_i}(\ze,\pow_{\ell_i},j) - 
\Psi^{\del_{\ell_i},p_{\ell_i}}_{0}(x_{\ell_i}(\ze',\und{y}_{\ell_i}))\gotg_{\ell_i}(\ze',\pow_{\ell_i},j)  ,
\end{equation}
if $q_{\ell_i} = 0,1,\ldots,M(\vartheta)-1$, and
\begin{equation} \label{grdifdifbb}
\ol\matG_{\ell_i}^\ttR(\ze)-\ol\matG_{\ell_i}^\ttR(\ze')
= \Psi^{\del_{\ell_i},p_{\ell_i}}_{0}(x_{\ell_i}(\ze,\und{y}_{\ell_i})) \gotG_{\ell_i}(\ze,j) - 
\Psi^{\del_{\ell_i},p_{\ell_i}}_{0}(x_{\ell_i}(\ze',\und{y}_{\ell_i})) \gotG_{\ell_i}(\ze',j) ,
\end{equation}
if $q_{\ell_i} = M(\vartheta)$, where $\gotg_\ell(\ze,\pow,j)$ is given by \eqref{mava}, while $\gotG_\ell(\ze,j)$ is implicitly defined by \eqref{odiotuttidavvero}.

In order to bound \eqref{grdifdifaa}, we proceed as in the proof of Lemma \ref{incartate}, by writing 
\begin{equation} \label{bibi}
\begin{aligned}
& \ol\matG_{\ell_i}^\ttR(\ze)-\ol\matG_{\ell_i}^\ttR(\ze') =
\big( \Psi^{\del_{\ell_i},p_{\ell_i}} _{0}(x_{\ell_i}(\ze,\und{y}_{\ell_i}))
 - \Psi^{\del_{\ell_i},p_{\ell_i}}_{0}(x_{\ell_i}(\ze',\und{y}_{\ell_i})) \big)
\gotg_{\ell_i}(\ze,\pow_{\ell_i},j) \\
& \qquad \qquad +
\Psi^{\del_{\ell_i},p_{\ell_i}}_{0}(x_{\ell_i}(\ze',\und{y}_{\ell_i}))
 \! \left( \gotg_{\ell_i}(\ze,\pow_{\ell_i},j) - \gotg_{\ell_i}(\ze',\pow_{\ell_i},j) \! \right) \! ,
\end{aligned}
\end{equation}
where, analogously to the proof of Lemma \ref{incartate}, we may assume that 
$\Psi^{\del_{\ell_i},p_{\ell_i}} _{0}(x_{\ell_i}(\ze,\und{y}_{\ell_i})) \neq 0$. Recalling \eqref{gliC},
%if $C_q(\ze,\ell,j)$ and $C_q'(\ell,j)$ denote the coefficients \eqref{gliC} corresponding to the frequencies $\om$ and $\om'$, respectively, then one has,
if $\ze=(\ka,\xi)$ and $\ze'=(\ka',\xi')$, one has
\begin{equation} \label{helsinki56}
\begin{aligned}
| C_0(\ze,\ell,j) - C_0(\ze',\ell,j) | & \le K_* \|\ze-\ze'\|_{\io} \left( \| \nu_\ell^\flat\|_{1} + 2 \right) , \\
| C_q(\ze,\ell,j) - C_q(\ze',\ell,j) | & \le K_* \|\ka-\ka'\|_{\io} (\calJ_{\al}(\vartheta))^{q-2}, \qquad q=2,\ldots,N-1 ,
\end{aligned}
\end{equation}
for some positive constant $K_*$ depending on $N$. Thus
the first contribution in \eqref{bibi} is dealt with
as in the proof of Lemma \ref{incartate}, so as to obtain the bound
\[
\begin{aligned}
&\left| \big( \Psi^{\del_{\ell_i},p_{\ell_i}} _{0}(x_{\ell_i}(\ze,\und{y}_{\ell_i}))
 - \Psi^{\del_{\ell_i},p_{\ell_i}}_{0}(x_{\ell_i}(\ze',\und{y}_{\ell_i})) \big) 
\gotg_\ell(\ze,\pow,j) \right| \\
&\qquad\qquad\qquad\le
C \frac{| x_{\ell_i}(\ze,\und{y}_{\ell_i}) -x_{\ell_i}(\ze',\und{y}_{\ell_i})  | }{(\be^*(r_{0}^*,\g))^{\del_{\ell_i} - p_{\ell_i} + 1}} 
\, |j|^{-q_{\ell_i}} \, (\calJ_{\al}(\vartheta))^{2(q_{\ell_i}-1)}\\
&\qquad\qquad\qquad\le C \, |j|^{-q_{\ell_i}}  \|\ze-\ze'\|_{\io} (\calJ_{\al}(\vartheta))^{2(q_{\ell_i}-1)},
\end{aligned}
\]
while the second one is found to be bounded as
\[
\left|  \Psi^{\del_{\ell_i},p_{\ell_i}}_{0}(x_{\ell_i}(\ze',\und{y}_{\ell_i}))\left( 
 \gotg_{\ell_i}(\ze,\pow_{\ell_i},j) - \gotg_{\ell_i}(\ze',\pow_{\ell_i},j) \right) \right| \le
C \, |j|^{-q_{\ell_i}}  \|\ze-\ze'\|_{\io} (\calJ_{\al}(\vartheta))^{2(q_{\ell_i}-1)} , 
\]
for some positive constant $C$.

By collecting together the two bounds we obtain that, if $q_{\ell_i}=0,1,\ldots,M(\vartheta)-1$,
then, for some positive constant $C$, one has
\begin{equation} \label{bound:g-g'}
\left|  \ol\matG_{\ell_i}^\ttR(\ze)-\ol\matG_{\ell_i}^\ttR(\ze')  \right|
\le C \, |j|^{-q_{\ell_i}}  \|\ze-\ze'\|_{\io} (\calJ_{\al}(\vartheta))^{2(q_{\ell_i}-1)} .
\end{equation}

Shorten $M=M(\vartheta)$ not to overwhelm the notation.
By reasoning as in the proof of Lemma \ref{helsinki} and using \eqref{helsinki56}, one finds
\[
|j|^M \left| r_M(\ze,\ell,j) - r_M(\ze',\ell,j) \right| \le  C \, \|\ze-\ze'\|_{\io} (\calJ_{\al}(\vartheta))^{M} ,
\]
for some positive constant $C$. Then, in order to study the difference in \eqref{grdifdifbb},
as in the proof of Lemma \ref{belfi}, we define $\vartheta'$ as the subtree such that $\ell_i\in\calP_{\vartheta'}$ and
distinguish among the following cases.

\begin{itemize}[topsep=0ex]
\itemsep0em
\item If $|j| \le (M+5)(\calJ_{\al}(\vartheta))^2$, we write
\[
\begin{aligned}
\left| j^{M} \, \left(\gotG_{\ell_i}^{(M)}(\ze,j)- \gotG_{\ell_i}^{(M)}(\ze',j) \right) \right|
& \le \left| \frac{j^M}{(x_{\ell_i}(\ze,\und{y}_{\ell_i}))^{p_{\ell_i}+1}} - \frac{j^M}{(x_{\ell_i}(\ze',\und{y}_{\ell_i}))^{p_{\ell_i}+1}} \right| \\
& + 
|j|^M\sum_{q=1}^{M-1} 
\left| \gotg_{\ell_i}(\ze,q,j) - \gotg_{\ell_i}(\ze',q,j) \right| .
\end{aligned}
\]
\item If $|j| > (M+5)(\calJ_{\al}(\vartheta))^2$ and
$\s_{\ell_i}\s_{\ell_{\vartheta'}}=-1$, defining $\de_j$ and $\e_j$ as in \eqref{djej1}, with $\ell={\ell_i}$, and similarly $\de_j'$ and $\e_j'$,
with $\ze'$ instead of $\ze$, we bound
\[
\begin{aligned}
\left| \gotG^{(M)}_{{\ell_i}}(\ze,j) - \gotG^{(M)}_{{\ell_i}}(\ze',j) \right| 
& \le \left| \frac{1}{(-2j^2)^{p_{\ell_i}+1}} \left( \frac{1}{(1+\de_j+\eps_j)^{p_{\ell_i}+1}}  - \frac{1}{(1+\de_j+\eps'_j)^{p_{\ell_i}+1}}\right) \right| \\
& + \left| \frac{1}{(-2j^2)^{p_{\ell_i}+1}} \left( \frac{1}{(1+\de_j+\eps'_j)^{p_{\ell_i}+1}}  - \frac{1}{(1+\de'_j+\eps'_j)^{p_{\ell_i}+1}} \right) \right| \\
& + \left|\frac{1}{(-2j^2)^{p_{\ell_i}+1}} \left( \frac{1}{(1+\de_j'+\eps_j')^{p_{\ell_i}+1}}  - \frac{1}{(1+\de_j')^{p_{\ell_i}+1}} \right)\right| \\
& + \Biggl| \Bigg( \frac{1}{(-2j^2)^{p_{\ell_i}+1}}\frac{1}{(1+\de_j)^{p_{\ell_i}+1}  } -  \sum_{q=1}^{M-1} \gotg_{\ell_i}(\ze,q,j) \Bigg) \\
& \qquad\qquad  - \Bigg( \frac{1}{(-2j^2)^{p_{\ell_i}+1}}\frac{1}{(1+\de'_j)^{p_{\ell_i}+1}  } +  \sum_{q=1}^{M-1} \gotg_{\ell_i}(\ze',q,j) \Bigg) \Biggr| ,
\end{aligned}
\]
and rewrite the last two lines as in \eqref{bohbohbohbohboh}, with 
\begin{equation} \label{B-B}
B_{q_1}(\ze,{\ell_i},j) \ldots B_{q_p}(\ze,{\ell_i},j) - B_{q_1}(\ze',{\ell_i},j) \ldots B_{q_p}(\ze',{\ell_i},j) 
\end{equation}
instead of $B_{q_1}(\ze,\ell,j) \ldots B_{q_p}(\ze,\ell,j)$.
\item If $|j| > (M+5)(\calJ_{\al}(\vartheta))^2$ and $\s_{\ell_i}\s_{\ell_{\vartheta'}}=1$, we obtain a bound
which differs from that of the previous case because
\begin{itemize}[topsep=0ex]
\itemsep0em
\item[1.]
$\de_j$ and $\e_j$, as well as $\de_j'$ and $\e_j'$, with $\ze'$ instead of $\ze$, are defined as in \eqref{djej2},
with $\ell={\ell_i}$;
\item[2.]
the factor $1/(-2j^2)^{p_{\ell_i}+1}$ is replaced with $1/(2 \s (j_{\ell_{\vartheta'}}^\flat - j_{\ell_i}^\flat ))^{p_{\ell_i}+1}$;
\item[3.]
the analogous of the last two lines of the previous bound are written as in
\eqref{tramonto}, with
\begin{equation} \label{B'-B'}
B'_{q_1}(\ze,{\ell_i},j) \ldots B'_{q_p}(\ze,{\ell_i},j) - B'_{q_1}(\ze',{\ell_i},j) \ldots B'_{q_p}(\ze',{\ell_i},j) 
\end{equation}
replacing $B'_{q_1}(\ze,\ell,j) \ldots B'_{q_p}(\ze,\ell,j)$.
\end{itemize}
\end{itemize}
By collecting together the results above we obtain that
\begin{equation} \label{G-G'}
|j|^M \left| \gotG_{\ell_i}^{(M)}(\ze,j) - \gotG_{\ell_i}^{(M)}(\ze',j) \right| \le C \, \|\ze-\ze'\|_\io (\calJ_{\al}(\vartheta))^{2M} ,
\end{equation}
for some constant $C$. Therefore, by writing \eqref{grdifdifbb} as
\begin{equation} \nonumber
\begin{aligned}
&  \ol\matG_{\ell_i}^\ttR(\ze)-\ol\matG_{\ell_i}^\ttR(\ze')  = 
\big( \Psi^{\del_{\ell_i},p_{\ell_i}} _{0}(x_{\ell_i}(\ze,\und{y}_{\ell_i}))
 - \Psi^{\del_{\ell_i},p_{\ell_i}}_{0}(x_{\ell_i}(\ze',\und{y}_{\ell_i})) \big)
\gotG_{\ell_i}^{(M)}(\ze,j) \\
&\qquad\qquad+
\Psi^{\del_{\ell_i},p_{\ell_i}}_{0}(x_{\ell_i}(\ze,\und{y}_{\ell_i}))
 \! \left( \gotG_{\ell_i}^{(M)}(\ze,j) - \gotG_{\ell_i}^{(M)}(\ze',j) \! \right) \! ,
\end{aligned}
\end{equation}
dealing with the first contribution as done for \eqref{grdifdifaa} and using \eqref{G-G'} to bound
the second one, finally we obtain, for $q_{\ell_i}=M$,
\begin{equation} \label{bound:G-G'}
\left|   \ol\matG_{\ell_i}^\ttR(\ze)-\ol\matG_{\ell_i}^\ttR(\ze')\right|
\le C \, |j|^{-M}  \|\ze-\ze'\|_{\io} (\calJ_{\al}(\vartheta))^{2M} ,
\end{equation}
for some positive constant $C$.

Now, we can proceed exactly as in the proofs of Lemma \ref{incartate},
but relying, in order to deal with the $\eta$-fully renormalized trees,
on the bounds of Lemma \ref{trentuno} rather than Lemma \ref{converge1}.
This concludes the proof.
\EP
\begin{rmk}\label{maieutica}
\emph{
Note that  $\xi\in\ell^{N,\io}\subset \ell^\io$, but in both Lemmata \ref{incartate} and \ref{aristesso} we consider Lipschitz 
semi-norms w.r.t.~the $\ell^\io$-norm, thus we require a stronger Lipschitz condition. 
} 
\end{rmk}
%%%%%%%%%%%%%%%%%%%%%%%%%%%%%%%%%%%%%%%%%%%%%%%%%%%%%%%%%%%%%%%%%%%%%%%%%% 

%It is now easy to prove Proposition \ref{fabrizio}:
%\medskip
%
%%%%%%%%%%%%%%%%%%%%%%%%%%%%%%%%%%%%%%%%%%%%%%%%%%%%%%%%%%%%%%%%%%%%%%%%%%% 
%
%\noindent {\it Proof of Proposition \ref{fabrizio}}.
The results above lead immediately to Proposition \ref{fabrizio}. Indeed the proof
 proceed as the proof of Proposition \ref{bellaprop} and Corollary \ref{contenti}, using Lemmata \ref{uniforme}, %\ref{stessoperglia}, 
\ref{incartate} and \ref{aristesso}. This fixes $\e_2=\e_2(s,\al,N,\g)$ and concludes the proof 
of Proposition \ref{fabrizio}.

%%%%%%%%%%%%%%%%%%%%%%%%%%%%%%%%%%%%%%%%%%%%%%%%%%%%%%%%%%%%%%%%%%%%%%%%%% 

%%%%%%%%%%%%%%%%%%%%%%%%%%%%%%%%%%%%%%%%%%%%%%%%%%%%%%%%%%%%%%%%%%%%%%%%%% 
\subsection{Lipschitz extensions}\label{lipext}
%%%%%%%%%%%%%%%%%%%%%%%%%%%%%%%%%%%%%%%%%%%%%%%%%%%%%%%%%%%%%%%%%%%%%%%%%% 

We are now ready to extend the functions appearing in Proposition \ref{fabrizio} to Lipschitz continuous functions
defined on the whole $\matW_N$.
To this aim we use the following version \cite{BMP3} of McShane's theorem \cite{mcshane}.

%%%%%%%%%%%%%%%%%%%%%%%%%%%%%%%%%%%%%%%%%%%%%%%%%%%%%%%%%%%%%%%%%%%%%%%%%% 
\begin{theo}[\textbf{McShane Theorem}]\label{mcsh}
	Let $X$ be a metric space, $A\subset X$ a non-empty subset and $f:A\to\ell^\io(\RRR)$ a Lipschitz function.
	Then there exists a Lipschitz extension $\bar{f}:X\to\ell^\io(\RRR)$ such that
	$\bar{f}|_A=f$ and $\bar{f}$ has the same Lipschitz constant as $f$. Moreover if $M:=\sup_A|f|<\io$, then $f^{\rm Ext}:=\max\{-M,\min\{\bar{f},M\}\}$
	is also a Lipschitz extension of $f$ with the same Lipschitz constant as $f$ and one has $\sup_X|f^{\rm Ext}|=M$. 
\end{theo}
%%%%%%%%%%%%%%%%%%%%%%%%%%%%%%%%%%%%%%%%%%%%%%%%%%%%%%%%%%%%%%%%%%%%%%%%%% 

Thanks to Theorem \ref{mcsh}, for any $\vartheta\in \gotF^{(k)}_{j,\nu,\s}$,
we can extend  $\Val^\ttR (\vartheta;c,\om(\ze))$ to a function $\Val^{\ttR, {\rm{Ext}} }(\vartheta;c,\ze)$
defined for all $\ze\in \matW_N $ and Lipschitz-continuous w.r.t.~the $\ell^\infty$ norm;
similarly, for $\vartheta\in {\gotF\gotH}^{(k)}_{j,+}(\pow)$, we can extend $\ol\Val^\ttR (\vartheta;c,\ze) $ to a function
$\ol{\Val}^{\ttR, {\rm{Ext}} }(\vartheta;c,\ze)$ defined for all $\ze\in \matW_N $ and Lipschitz-continuous w.r.t.~the $\ell^\infty$ norm.
Recalling that both $\Val^\ttR (\vartheta;c,\om(\ze)) $ and $\ol\Val^\ttR (\vartheta;c,\ze)$ depend only on a finite number of components 
of $\x$, the extended functions are also continuous w.r.t.~the product topology.
By construction, all possible extended functions coincide for all good parameters $\ze \in \matK_N(\g)$. 

For $(c,\ze,\e)\in  \matU_{1}(\mathtt g{(s,\alpha)})\times \matW_N\times (-{\e_2},{\e_2})$, with ${\e_2}$ as in Proposition \ref{fabrizio}, 
given any extension $\ol{\Val}^{\ttR, {\rm{Ext}} }(\vartheta;c,\ze)$ of  $\ol\Val^\ttR (\vartheta;c,\ze)$, resp. 
 ${\Val}^{\ttR, {\rm{Ext}} }(\vartheta;c,\ze)$ of  $\Val^\ttR (\vartheta;c,\ze)$,
define
\begin{subequations} \label{estendo}
\begin{align}
{\gota}_\pow^{(k,{\rm{Ext}})}(c,\ze,j) & := -\frac{1}{c_j}  \sum_{\vartheta\in {\gotF\gotH}^{(k)}_{j,+}(\pow)} \ol{\Val}^{\ttR, {\rm{Ext}} }(\vartheta;c,\ze),
\qquad j \neq 0 , \quad \pow =0,\ldots,N-1,
\label{estendoa} \\ 
%\notag
{\gotA}^{(k,{\rm{Ext}})}(c,\ze,j) & := -\frac{1}{c_j}  \sum_{\vartheta\in {\gotF\gotH}^{(k)}_{j,+}(N)} \ol{\Val}^{\ttR, {\rm{Ext}} }(\vartheta;c,\ze),
\qquad \!\! j \neq 0 ,
\label{estendob} \\
{\h}_0^{(k,{\rm{Ext}})}(c,\om(\ze)) & :=-\frac{1}{c_0}  \sum_{\vartheta\in {\gotF}^{(k)}_{0,\gote_0,+}} {\Val}^{\ttR, {\rm{Ext}} }(\vartheta;c,\ze) ,
\label{estendof}
\end{align}
\end{subequations}
and set
\begin{subequations} \label{estendosum}
\begin{align}
{\gota}^{{\rm{Ext}}}_\pow(c,\ze,\e,j) & :=  \sum_{k\ge1}\e^k {\gota}_\pow^{(k, {\rm{Ext}})}(c,\ze, j),
\qquad j \neq 0 , \quad q=0,\ldots,N-1,
\label{estendoc} \\ 
{\gotA}^{{\rm{Ext}}}(c,\ze,\e,j) &:=  \sum_{k\ge1}\e^k {\gotA}^{(k, {\rm{Ext}})}(c,\ze, j) ,
\qquad \! j \neq 0 ,
\label{estendod} \\
\h_0^{{\rm{Ext}}}(c,\om(\ze),\e) & := \sum_{k,\ge1} \e^k {\h}_0^{(k,{\rm{Ext}})}(c,\om(\ze)) .
\label{estendoe}
\end{align}
\end{subequations}
As a consequence of Theorem \ref{mcsh} the series \eqref{estendosum} are absolutely convergent, independently of
the extension. However, in order to prove that 
\begin{equation} \nonumber
\lim_{|j|\to \io} \gota^{{\rm{Ext}}}_\pow(c,\ze,\e,j)
\end{equation}
exists, and is continuous w.r.t.~the product topology and Lipschitz-continuous w.r.t.~the $\ell^\io$ norm,
we need to check that it is possible to define the extension in such a way that the sequence
\[
\{ {\gota}^{{\rm{Ext}}}_\pow(c,\ze,\e,j) \}_{j\in\ZZZ\setminus\{0\}}
\]
is a Cauchy sequence. Thus, we need that the cancellation mechanism exploited in the proof of Lemma \ref{squa} still works:
this can be achieved by choosing carefully the extended values $\ol{\Val}^{\ttR, {\rm{Ext}} }(\vartheta;c,\ze)$ resp. 
${\Val}^{\ttR, {\rm{Ext}} }(\vartheta;c,\ze)$.
 
%%%%%%%%%%%%%%%%%%%%%%%%%%%%%%%%%%%%%%%%%%%%%%%%%%%%%%%%%%%%%%%%%%%%%%%%%% 
\begin{lemma}\label{pieta}
Let $\e_2=\e_2(s,\al,N,\g)$ be as in Proposition \ref{fabrizio}.
There exist extensions $\gota_\pow^{\rm{Ext}}(c,\cdot,\e,j)$, $\gotA^{{\rm{Ext}}}(c,\cdot,\e,j)$ and $\h_0^{{\rm{Ext}}}(c,\om(\cdot),\e)$
to the whole $\matW_N$ of the functions $\gota_\pow(c,\cdot,\e,j)$, $\gotA(c,\cdot,\e,j)$ and $\h_0(c,\om(\cdot),\e)$, respectively,
such that, for all $(c,\e)\in\matU_1(\mathtt{g}(s,\al))\times(-{\e_2},{\e_2})$ and all $j\in\ZZZ\setminus\{0\}$,
\begin{enumerate}[topsep=0ex]
\itemsep0em
\item they are continuous w.r.t.~the product topology and Lipschitz-continuous w.r.t.~the $\ell^\io$ norm;
\item setting %$\gota_\pow^{{\rm{Ext}}}(c,\ze,\e,0):=0$ and 
$\gotA^{{\rm{Ext}}}(c,\ze,\e,0):=0$, 
and calling $\gotA^{{\rm{Ext}}}(c,\ze,\e):=\{\gotA^{{\rm{Ext}}}(c,\ze,\e,j)\}_{j\in\ZZZ}$,
there is $C>0$ such that
\begin{subequations} \label{apatiafin}
\begin{align}
\sup_{j\in\ZZZ\setminus\{0\}}|\gota^{{\rm{Ext}}}_\pow(c,\cdot,\e,j)|^{{\rm{Lip}}}_{\RRR} & \le C|\e|, 
\qquad q=0,\ldots,N-1,
\label{apatiafinb} \\
\phantom{\sup_{j\in\ZZZ\setminus\{0\}}}
|\gotA^{{\rm{Ext}}}(c,\cdot,\e)|^{{\rm{Lip}}}_{\ell^{N,\infty}} & \le C|\e| ,
\label{apatiafinc} \\
|\h_0^{{\rm{Ext}}}(c,\om(\cdot),\e)|^{{\rm{Lip}}}_\RRR & \le C|\e| ;
\label{apatiafina}
\end{align}
\end{subequations}
\item there are $\de_0,\al_0>0$ such that, for $q=0,\ldots,N-1$,
\begin{subequations} \label{apatialim}
\begin{align}
| \gota_\pow^{{\rm{Ext}}}(c,\cdot,\e,j)-\gota_\pow^{{\rm{Ext}}}(c,\cdot,\e,j')|^{{\rm{Lip}}}_\RRR	
&\le C|\e|e^{-\de_0\min\{\jap{j}^{\al_0},\jap{j'}^{\al_0}\}},
\label{diocane1} \\
|\gota^{{\rm{Ext}}}_1(c,\cdot,\e,j)|^{{\rm{Lip}}}_{\RRR} & \le C|\e|
e^{-\de_0\jap{j}^{\al_0}}.
\label{diocane2} 
\end{align}
\end{subequations}
\end{enumerate}
\end{lemma}
%%%%%%%%%%%%%%%%%%%%%%%%%%%%%%%%%%%%%%%%%%%%%%%%%%%%%%%%%%%%%%%%%%%%%%%%%% 

%%%%%%%%%%%%%%%%%%%%%%%%%%%%%%%%%%%%%%%%%%%%%%%%%%%%%%%%%%%%%%%%%%%%%%%%%% 
\prova
For $q=0,1,\ldots,N$, define
\[
\BBB^{(k)}_{j,+}(\pow):=\{\vartheta\in {\gotF\gotH}^{(k)}_{j,+}(\pow) :
|j_{\la}| < |j| \; \: \forall \la \in \Lambda(\vartheta) \setminus\{\la_{\vartheta} \} \}
\]
and set $\gotW^{(k)}_{j,+}(\pow):=\gotF\gotH^{(k)}_{j,+}(\pow) \setminus \BBB^{(k)}_{j,+}(\pow)$.
Recalling Remark \ref{motooscrivo}, we use the representation $\tilde{\gota}_\pow^{(k)}(c,\ze,j)$ in \eqref{gliab},
which we write as
\begin{equation} \label{gliabb}
\tilde{\gota}_\pow^{(k)}(c,\ze,j) = \tilde{\gotb}_\pow^{(k)}(c,\ze,j) + \tilde{\gotw}_\pow^{(k)}(c,\ze,j) ,
\end{equation}
where
\[
\tilde{\gotb}_\pow^{(k)}(c,\ze,j) := -\frac{j^q}{c_j}  
\sum_{\vartheta\in {\BBB}^{(k)}_{j,+}(\pow)}\!\!\!\!\!  \overline{\Val}^\ttR(\vartheta;c,\ze) , \qquad
\tilde{\gotw}_\pow^{(k)}(c,\ze,j) := -\frac{j^q}{c_j}  
\sum_{\vartheta\in {\gotW}^{(k)}_{j,+}(\pow)}\!\!\!\!\!  \overline{\Val}^\ttR(\vartheta;c,\ze) .
\]
Note that $\tilde{\gotw}_\pow^{(k)}(c,\ze,j) $ takes into account only contributions arising from trees
$\vartheta\in {\gotF\gotH}^{(k)}_{j,+}(\pow)$ such that $\overline{\Val}^\ttR(\vartheta;c,\ze)$ decays
exponentially in $\jap{j}^\al$, by Lemma \ref{trentuno}, because any $\vartheta\in{\gotW}^{(k)}_{j,+}(\pow)$ 
contains at least one leaf $\la\neq\la_{\vartheta}$ with $|j_\la| \ge |j_{\la_{\vartheta}}|$.

We extend the values of the fully renormalized $\eta$-trees recursively on $j$ as follows.
\begin{enumerate}[topsep=0ex]
\itemsep0em
\item For $j=1$ and for any tree $\vartheta\in\gotF\gotH^{(k)}_{1,+}$ we define the extension
${\overline{\Val}^{\ttR,{\rm Ext}}} (\vartheta;c,\ze)$ by using Theorem \ref{mcsh}.
\item Now, consider any $j>0$ and assume that the values of the trees $\vartheta\in\gotF\gotH^{(k)}_{j',+}$,
with $j'\neq0$ such that $-j+1 \le j' \le j$ have been extended. Recalling Definition \ref{conjtree},
by construction each tree $\vartheta '\in \BBB^{(k)}_{-j,+}(\pow)$
is $(-j)$-associated with a tree $\vartheta \in \BBB^{(k)}_{j,+}(\pow)$; see also the proof of Lemma \ref{presqua}. 
On the other hand
each tree $\vartheta' \in \BBB^{(k)}_{j+1,+}(\pow)$
either is $(j+1)$-associated with a tree $\vartheta \in \BBB^{(k)}_{j,+}(\pow)$ 
or contains at least one leaf  $\la'$ with $j_{\la'}=j$.
\begin{itemize}
\item[2.1.]
For $j'=-j,j+1$, in the case in which $\vartheta\in \BBB^{(k)}_{j',+}(\pow)$ is $j'$-associated with $\vartheta$, if
\begin{equation}\nonumber
\frac{(j')^q}{c_{j'}} {\overline{\Val}^\ttR} (\vartheta';c,\ze) =
\frac{j^q}{c_{j}} {\overline{\Val}^\ttR} (\vartheta;c,\ze) ,
\end{equation}
we set 
\begin{equation}\nonumber
{\overline{\Val}^{\ttR,{\rm Ext}}} (\vartheta';c,\ze) :=
\frac{c_{j'}} {(j')^q}
\frac{(j)^q}{c_{j}} {\overline{\Val}^{\ttR,{\rm Ext}}} (\vartheta;c,\ze) .
\end{equation}
Otherwise, from Lemma \ref{presqua}, we deduce that 
\begin{equation} \label{maximumval}
\begin{aligned}
& \left|
\frac{(j')^q}{c_{j'}} {\overline{\Val}^\ttR} (\vartheta';c,\ze) \right| , \,
\left| \frac{(j)^q}{c_{j}} {\overline{\Val}^\ttR} (\vartheta;c,\ze) \right| \\
& \qquad \qquad \le 
A^{|N(\vartheta)|}K_1^{N-1} (\calJ_{\al}(\vartheta))^{2N^2}
e^{-(s'-s_3-s_4)\calJ(\vartheta)} e^{- s_4 2^{-\al} \jap{j}^{\al/2}} .
\end{aligned}
\end{equation}
and we can extend ${\overline{\Val}^\ttR} (\vartheta';c,\ze)$ in such a way that still one has
\begin{equation} \nonumber %\label{maximumvalext}
%\begin{aligned}
%& 
\left|
\frac{(j')^q}{c_{j'}} {\overline{\Val}^{\ttR,{\rm Ext}}} (\vartheta';c,\ze) \right| 
% \\ & \qquad \qquad 
\le 
A^{|N(\vartheta)|}K_1^{N-1} (\calJ_{\al}(\vartheta))^{2N^2}
e^{-(s'-s_3-s_4)\calJ(\vartheta)} e^{- s_4 2^{-\al} \jap{j}^{\al/2}} .
%\end{aligned}
\end{equation}
\item[2.2.]
 If either $\vartheta' \in \BBB^{(k)}_{j+1,+}(\pow)$ and it contains at least one leaf  $\la'$ with $j_{\la'}=j$ or
$\vartheta' \in \gotW^{(k)}_{j',+}(\pow)$, with $j'=-j,j+1$, we construct
${\overline{\Val}^{\ttR,{\rm Ext}}} (\vartheta';c,\ze)$ by using Theorem \ref{mcsh}.
Again, the extension still satisfies the bound
\[
\left| \frac{(j')^q}{c_{j'}} {\overline{\Val}^{\ttR,{\rm Ext}}} (\vartheta';c,\ze) \right| \le
A^{|N(\vartheta)|}K_1^{|N(\vartheta)|}  (\calJ_{\al}(\vartheta))^{2N^2}
e^{-2 \hat s \jap{j}^{\al}}  e^{- \bar  s  \calJ(\vartheta)} . \]
\end{itemize}
\end{enumerate}
With the extensions defined as above,
fixing $s'$ and $s_3$ as in the proof of Corollary \ref{stimo},
 the bounds \eqref{apatiafin} follow immediately from Theorem \ref{mcsh}.
Moreover, Lemma \ref{presqua} still applies to the extended tree values, then the bounds \eqref{diocane1} are obtained 
by reasoning as in the proof of Lemma \ref{squa}, fixing $s_4$ and $\hat{s}$ in the same way.
Finally the bound \eqref{diocane2} is obtained by reasoning as in Lemma \ref{auno}.
\EP
%%%%%%%%%%%%%%%%%%%%%%%%%%%%%%%%%%%%%%%%%%%%%%%%%%%%%%%%%%%%%%%%%%%%%%%%%% 

%%%%%%%%%%%%%%%%%%%%%%%%%\mathtt g{(s,\alpha)}%%%%%%%%%%%%%%%%%%%%%%%%%%%%%%%%%%%%%%%%%%%%%%%%% 
\begin{coro} \label{noname}
For all $(c,\ze,\e)\in \matU_{1}(\mathtt g{(s,\alpha)})\times \matW_N\times (-{\e_2},{\e_2})$ and all $q=0,2,\dots,N-1$,
the sequence $\{\gota_\pow(c,\ze,\e,j)\}_{j\in\ZZZ}$ is a Cauchy sequence and
 \begin{equation}\label{Atau}
%\matW_N\ni(\ze)\mapsto 
A_\pow(c,\ze,\e):= \lim_{|j|\to \io} \gota^{{\rm{Ext}}}_\pow(c,\ze,\e,j)
 \end{equation}
 is continuous w.r.t.~the product topology and Lipschitz-continuous w.r.t.~the $\ell^\io$ norm, and one has
 \[
 |A_\pow(c,\cdot,\e)|^{{\rm{Lip}}}_{\RRR} \le C|\e|,
 \]
 with $C$ as in Lemma \ref{pieta}.
\end{coro}
%%%%%%%%%%%%%%%%%%%%%%%%%%%%%%%%%%%%%%%%%%%%%%%%%%%%%%%%%%%%%%%%%%%%%%%%%% 

We now define, for $(c,\ze,\e)\in \matU_{1}(\mathtt g{(s,\alpha)})\times \matW_N\times (-{\e_2},{\e_2})$ and $j\in\ZZZ\setminus\{0\}$,
\[
\eta_j^{{\rm{Ext}}}(c,\om(\ze),\e):= A_0(c,\ze,\e)+\sum_{\pow=2}^{N-1}\frac{A_\pow(c,\ze,\e)}{j^\pow} + R_j(c,\ze,\e),
\]
where  $A_0(c,\ze,\e),A_2(c,\ze,\e),\ldots, A_{N-1}(c,\ze,\e)$ are defined in \eqref{Atau}, while
\[
\begin{aligned}
R_j(c,\ze,\e):=  \gota_0^{{\rm{Ext}}}&(c,\ze,\e,j)-A_0(c,\ze,\e) + \frac{\gota_1^{{\rm{Ext}}}(c,\ze,\e,j)}{j} \\
&+\sum_{\pow=2}^{N-1}\frac{\gota_\pow^{{\rm{Ext}}}(c,\ze,\e,j)-A_\pow(c,\ze,\e)}{j^\pow} + \gotA^{{\rm{Ext}}}(c,\ze,\e,j)
\end{aligned}
\]
extends to $\matW_N$ the function $\gotr_j(c,\ze,\e)$ in Proposition \ref{decadeta}.
Finally we set
\[
R_0(c,\ze,\e):=  \eta_0^{{\rm{Ext}}}(c,\ze,\e)- A_0(c,\ze,\e) ,
\]
and define $R(c,\cdot,\e)=\{R_j(c,\cdot,\e)\}_{j\in\ZZZ}$.
Corollary \ref{noname} and a direct application of Lemma \ref{pieta} give the bounds
\begin{equation}
\label{stimAR}
|R(c,\cdot,\e)|^{{\rm{Lip}}}_{\ell^{N,\io}} \le C|\e|, \qquad
|A_q(c,\cdot,\e)|^{{\rm{Lip}}}_\RRR \le C|\e| , \quad q=0,2,\ldots,N-1.
\end{equation}

This concludes the extension argument.

%%%%%%%%%%%%%%%%%%%%%%%%%%%%%%%%%%%%%%%%%%%%%%%%%%%%%%%%%%%%%%%%%%%%%%%%%% 
%%%%%%%%%%%%%%%%%%%%%%%%%%%%%%%%%%%%%%%%%%%%%%%%%%%%%%%%%%%%%%%%%%%%%%%%%% 
\subsection{The implicit function equation}\label{implicit}

We are now ready to use a fixed point argument to solve the implicit function equation \eqref{system}.
For $(c,\e)\in \matU_1(\mathtt g(s,\al))\times (-{\e_2},{\e_2})$, define $\matF: [-1/2,1/2]^\ZZZ\times\matW_N\to \RRR^{N-1} \times \ell^\io $
as
\begin{equation} \label{kqinsieme}
\matF(V,\ze; c,\e)=(\{\matF_{\ka,q}(\ze; c,\e)\}_{q=0,2,\ldots,N-1}, \{\matF_{\x,j}(V,\ze; c,\e)\}_{j\in\ZZZ}),
\end{equation}
with
\begin{subequations} \label{mapF}
\begin{align}
\matF_{\ka,q}(V,\ze; c,\e) & := \kappa_q + A_q(c,\ze,\e)  , \qquad\qquad\quad\;\;\, q=0,2,\ldots, N-1 ,
\label{mapFa} \\
\matF_{\x,j}(V,\ze; c,\e) & := 
\x_j + R_j(c,\ze,\e) - V_j .
\label{mapFb}
\end{align}
\end{subequations}

Given any subset $\matU\subseteq \ell^\io$ and any map $f:\matU\to E$, 
with $(E,|\cdot|_E)$ some Banach space, we define the Lipschitz norm as
\begin{equation} \nonumber %\label{lipnorm2}
|f|_{\matU,E}^{{\rm{Lip}}}:= \sup_{V \in \matU} | f(V) |_E + 
\sup_{\substack{V,V'\in \matU \\ V \neq V'}} \frac{|f(V) -f(V') |_E}{\|V-V'\|_{\io}} .
\end{equation}

The map $\matF$ is well defined, continuous w.r.t.~the product topology and Lipschitz-continuous w.r.t.~the $\ell^\infty$ norm,
and it satisfies \eqref{stimAR}.
Moreover \eqref{system} can be written as $\matF(V,\kappa,\xi; c,\e)=0$.
\begin{lemma}\label{contraggo}
Fixed $(c,\e)\in \matU_1(\mathtt g(s,\al))\times (-{\e_2},{\e_2})$, let
%$\matF: [-1/2,1/2]^\ZZZ\times\matW_N\to \RRR^{N-1} \times \ell^\io$ be given by \eqref{mapF}.
$\matF\!:  \ol{\matU}_{1/4}(\ell^{N,\io})\times\matW_N\to \RRR^{N-1} \times \ell^\io$ be given by \eqref{mapF}.
%Then  there exists a map $\ze\!: [-1/4,1/4]^\ZZZ \to \matW_N$
Then  there exists a map $\ze\!: \ol{\matU}_{1/4}(\ell^{N,\io}) \to \matW_N$
%functions $ [-1/4,1/4]^\ZZZ \to \matU$ 
%	 \begin{equation}\label{mappa}
%	[-1/4,1/4]^\ZZZ\ni V\mapsto \ze(V) \in \matW_N ,
%	 \end{equation}
such that 
\begin{enumerate}[topsep=0ex]
\itemsep0em
\item it is continuous in the product topology and Lipschitz-continuous w.r.t.~the $\ell^{\infty}$ norm;
\item there is $C>0$ such that, setting $\Delta(V):= \x(V)-V$, one has
	\begin{equation}\label{item2}
	%\begin{aligned}
	%& 
	\left| \ka_q \right|^{\rm Lip}_{\ol{\matU}_{1/4}(\ell^{N,\io}),\RRR} 
	%= \sup_{V\in[-1/4,1/4]^\ZZZ} | \ka_q(V) | + \sup_{{\substack{V,V'\in[-1/4,1/4]^\ZZZ \\ V\ne V'}}} \frac{|\ka_q(V)-\ka_q(V') |}{\|V-V'\|_{\io}} 
	\le C|\e| ,
	\quad q=0,2,\ldots,N-1, \qquad\qquad
	%\\ & 
	\left| \Delta \right|^{\rm Lip}_{\ol{\matU}_{1/4}(\ell^{N,\io}),\ell^{N,\io}} 
	%= \sup_{V\in[-1/4,1/4]^\ZZZ} \|\Delta(V)\|_{N,\io} + \sup_{{\substack{V,V'\in[-1/4,1/4]^\ZZZ \\ V\ne V'}}} \frac{\|\Delta(V)-\Delta(V')\|_{N,\io}}{\|V-V'\|_{\io}} 
	\le C|\e| ;
	%\end{aligned}
	\end{equation}
	\item one has $\matF(V,\ze(V); c,\e)\equiv0$.
%	\item if $V\in \matU_{1/4}(\ell^{N,\io})$, then  $\xi(V) \in \matU_{1/2}(\ell^{N,\io}) $.
	\end{enumerate}
\end{lemma}

\prova
The proof is based on a fixed point argument for Lipschitz-continuous functions \cite{BMP3}.
Let us write $\matF(V,\ze; c,\e)=\ze - \matH(V,\ze; c,\e)$, with
\begin{subequations} \label{mapH}
\begin{align}
\matH_{\ka,q}(V,\ze; c,\e) & := - A_q(c,\ze,\e)  , \qquad q=0,2,\ldots, N-1 ,
\label{mapHa} \\
\matH_{\x,j}(V,\ze; c,\e) & := \begin{cases}
- R_0(c,\ze,\e) + V_ 0 , & \qquad j=0, \\
- R_j(c,\ze,\e) + V_j , & \qquad j \neq 0 , \end{cases}
\label{mapHb}
\end{align}
\end{subequations}
where we are using the same notation as in \eqref{kqinsieme} to describe the components of $\matH(V,\ze; c,\e)$.
The bounds \eqref{stimAR} imply that the function $\matH(V,\cdot; c,\e)$
is a contraction on $\matW_N$. Therefore, by \cite[Lemma B.1]{BMP3}, for any $V\in  \ol{\matU}_{1/4}(\ell^{N,\io})$
there exists a unique $\ze=\ze(V)\in\matW_N$ such that $\matH(V,\ze(V);c,\e)=V$ and $\ze(V)$ is continuous in the product topology
and Lipschitz-continuous w.r.t.~the $\ell^{\infty}$ norm, and satisfies the bounds in \eqref{item2}.
\EP
%%%%%%%%%%%%%%%%%%%%%%%%%%%%%%%%%%%%%%%%%%%%%%%%%%%%%%%%%%%%%%%%%%%%%%%%%% 

%%%%%%%%%%%%%%%%%%%%%%%%%%%%%%%%%%%%%%%%%%%%%%%%%%%%%%%%%%%%%%%%%%%%%%%%%% 
%%%%%%%%%%%%%%%%%%%%%%%%%%%%%%%%%%%%%%%%%%%%%%%%%%%%%%%%%%%%%%%%%%%%%%%%%% 
\subsection{Measure estimates and proof of Proposition \ref{bohboh}}
%\addcontentsline{toc}{subsection}{\thesection.\arabic{subsection}}{Expansion in ${j}$ of ${\h_j}$}
\label{misura}
%%%%%%%%%%%%%%%%%%%%%%%%%%%%%%%%%%%%%%%%%%%%%%%%%%%%%%%%%%%%%%%%%%%%%%%%%% 
%%%%%%%%%%%%%%%%%%%%%%%%%%%%%%%%%%%%%%%%%%%%%%%%%%%%%%%%%%%%%%%%%%%%%%%%%% 

Combining Lemma \ref{contraggo} with the definition \eqref{espome} we obtain
\begin{equation}
	\label{raccapriccio}
	\om_j(\ze(V))=\begin{cases}
	\kappa_0(V) + V_0 + \Delta_0(V),& \qquad j=0,\\
	\displaystyle{j^2 + \kappa_0(V)+\sum_{\pow=2}^{N-1}\frac{\kappa_\pow(V)}{j^\pow}+ V_j+ \Delta_j(V),}& \qquad j\ne0.
	\end{cases}
\end{equation}
Define also
\begin{equation}
	\label{raccapriccio2}
	\ol\om_j(V)=\begin{cases}
	V_0 + \Delta_0(V),& \qquad j=0,\\
	\displaystyle{j^2 +\sum_{\pow=2}^{N-1}\frac{\kappa_\pow(V)}{j^\pow}+ V_j+ \Delta_j(V),}& \qquad j\ne0 ,
	\end{cases}
\end{equation}
so that $\ol\om_j(V)=\om_j(\ze(V))-\ka_0(V)$ for all $j\in\ZZZ$.

%%%%%%%%%%%%%%%%%%%%%%%%%%%%%%%%%%%%%%%%%%%%%%%%%%%%%%%%%%%%%%%%%%%%%%%%%% 
\begin{lemma}\label{canali}
For any $\vartheta \in \gotF\cup\gotF\gotH$ and any line $\ell\in L(\vartheta)$ one has
\[
\om(\ze(V))\cdot\nu_{\ell} - \om_{j_\ell} (\ze(V)) = \ol\om(V)\cdot\nu_{\ell} - \ol\om_{j_\ell} (V) .
\]
\end{lemma}
%%%%%%%%%%%%%%%%%%%%%%%%%%%%%%%%%%%%%%%%%%%%%%%%%%%%%%%%%%%%%%%%%%%%%%%%%% 

%%%%%%%%%%%%%%%%%%%%%%%%%%%%%%%%%%%%%%%%%%%%%%%%%%%%%%%%%%%%%%%%%%%%%%%%%% 
\prova
For all $\ell\in L(\vartheta)$ one has
\[
\om(\ze(V))\cdot\nu_{\ell} - \om_{j_\ell} (\ze(V)) - \bigl( \ol\om(V)\cdot\nu_{\ell} - \ol\om_{j_\ell} (V) \bigr) =
\ka_0(V) \sum_{i\in\ZZZ} \Bigl( (\nu_\ell)_i - \de_{i,j_{\ell}} \Bigr) = 0 ,
\]
where the last equality follows from \eqref{tocchera}.
\EP
%%%%%%%%%%%%%%%%%%%%%%%%%%%%%%%%%%%%%%%%%%%%%%%%%%%%%%%%%%%%%%%%%%%%%%%%%% 

For all $\nu\in \ZZZ^{\ZZZ}_{f} \setminus\{0\}$  and $\de>0$, set 
 \[
\mathfrak R_\nu(\de):=\{V\in \ol{\matU}_{1/4}(\ell^{N,\io}) : |\ol\om(V)\cdot \nu|\le \de \}\,.
 \]
 
We note that  $\mathfrak R_\nu(\de)$ is measurable w.r.t.~the probability measure induced by the product measure on $\ell^\io$,
since it is the preimage of  a closed set via the  function
$V\mapsto \ol\om(V)\cdot\nu$, which by Lemma \ref{contraggo} is continuous in the product topology. 
Our aim is to give an upper bound on the probability measure $\meas(\mathfrak R_\nu(\de))$.

%%%%%%%%%%%%%%%%%%%%%%%%%%%%%%%%%%%%%%%%%%%%%%%%%%%%%%%%%%%%%%%%%%%%%%%%%% 
\begin{lemma}\label{misuretta}
For all $\de$ sufficiently small and all $\nu\in \ZZZ^{\ZZZ}_{f} \setminus\{0\}$ such that $\sum_{i\in\ZZZ}\nu_i=0$, one has
\begin{equation}
\label{upper}
\meas(\mathfrak R_\nu(\de)) \le 2 \jap{i_0(\nu)}^N \de ,
\end{equation}
where
\[
i_0(\nu):= \max\{i\in \ZZZ:   |\nu_i|=\|\nu\|_\infty\} .
\]
\end{lemma}
%%%%%%%%%%%%%%%%%%%%%%%%%%%%%%%%%%%%%%%%%%%%%%%%%%%%%%%%%%%%%%%%%%%%%%%%%% 

%%%%%%%%%%%%%%%%%%%%%%%%%%%%%%%%%%%%%%%%%%%%%%%%%%%%%%%%%%%%%%%%%%%%%%%%%% 
\prova
For all $\nu\in\ZZZ^{\ZZZ}_{f}$ such that $\sum_{i\in\ZZZ}\nu_i=0$ one has
\[
\ol\omega(V)\cdot \nu = \sum_{j\in \ZZZ} j^2 \nu_j +
V\cdot\nu + \Delta(V)\cdot\nu +  \sum_{q=2}^{N-1} \kappa_q(V) \sum_{j\in \ZZZ\setminus\{0\}} 
\frac{\nu_j}{j^q}.
\]

Fix $\nu\in\ZZZ^{\ZZZ}_{f}\setminus\{0\}$ with $\sum_{i\in\ZZZ}\nu_i=0$. Write $V= (V_{i_0},V')$, with $i_0=i_0(\nu)$, and $V'=\{V_j\}_{j\in\ZZZ\setminus\{i_0\}}$.
For $h\in\RRR$ small enough one has
\begin{align*}
& \frac{|\ol{\omega}(V_{i_0}+h,V')\cdot \nu - \overline\omega(V_{i_0},V')\cdot \nu |}{h} \\
& \qquad \ge | \nu_{i_0}| 
\Big( 1 - |\Delta|^{{\rm{Lip}}}_{\ol{\matU}_{1/4}(\ell^{N,\io}),\ell^{N,\io}}\sum_{j\in\ZZZ}\jap{j}^{-N} - 
\sum_{q=2}^{N-1} |\kappa_q|^{\text{Lip}}_{\ol{\matU}_{1/4}(\ell^{N,\io}),\RRR}  \sum_{j\in \ZZZ} \jap{j}^{-q}
\Big) \\
& \qquad \ge  |\nu_{i_0} |(1- C|\e| ) \ge \frac12,
\end{align*}
for some contant $C=C(N)>0$ such that $C(N)\to+\io$ as $N\to+\io$.
Reasoning as in \cite{BMP3}, by using the latter bound, we find that, for $\de$ small enough, the set
$\mathfrak R_\nu(\de)$ is contained in the normal domain
\[
\gotE:=\Biggl\{ V\in \ell^{N,\infty}:  a(V') \le V_{i_0} \le b(V') \,, \;  \sup_{\substack{j\in\ZZZ \\ j\ne i_0}}|V_j| \jap{j}^N \le \frac{1}{4} \Biggr\} ,
\]
where $a$ and $b$ are two  functions continuous w.r.t.~the product topology such that
%$a(V')\le b(V')$ for all $V'\in \matU_{1/4}(\ell^{N,\io})$ and
\[
0 \le b(V')- a (V') \le 2\de \qquad \forall V'\in \ol{\matU}_{1/4}(\ell^{N,\io}).
\]
Then, by Fubini's theorem,
\[
\meas(\mathfrak R_\nu(\de))\le \meas({\gotE}) \le 4 \jap{i_0}^N \de\,,
\]
thus concluding the proof.
\EP
%%%%%%%%%%%%%%%%%%%%%%%%%%%%%%%%%%%%%%%%%%%%%%%%%%%%%%%%%%%%%%%%%%%%%%%%%% 

%\noindent{\it Proof of Proposition \ref{bohboh}}.
Now, in order to prove Proposition \ref{bohboh}, i.e.~that the set $\widehat{\matC}_N(\e)$ defined in \eqref{rettangolino} has asymptotically full measure, 
thanks to Lemma \ref{diovsbrj}  it is enough to verify that the set
\[
\tilde{\matC}_N(\e):=\{ V\in \ol{\matU}_{1/4}(\ell^{N,\io}) : \om(V)\in \gD^{(0)}({\g,N+1})\} %\mathtt{D}_{\g,N+1}\}
%\{  \zeta\in B_{1,N}:\quad |\overline{\omega}(\zeta)\cdot \nu|\ge \gamma \mathtt d(\nu)^{-N-1} \,,\quad \forall \nu\in\ZZZ^{\ZZZ}_{f}: \sum_j\nu_j=0\}\,,\quad \mathtt d (\nu)= \prod_j (1+\nu_j^2 \jap{j}^2)
\]
has asymptotically full measure, where 
\begin{equation}\label{dio0}
\gD^{(0)}({\g,N+1} ) :=\Bigl\{ \omega\in {\gQ} : |\omega\cdot \nu|> \g
\prod_{i\in \ZZZ}\frac{1}{(1+ \jap{i}^{2}|\nu_i|^{2})^\tau} \;\ \forall \nu\in \ZZZ^{\ZZZ}_{f} \setminus\!\{0\},\ \sum_{i\in\ZZZ}\nu_i=0 \Bigr\}.
\end{equation}
 By construction $\tilde{\matC}_N(\e)$ is measurable, and $\omega(\ze(V))\in \gQ$ for all $V\in \ol{\matU}_{1/4}(\ell^{N,\io})$, hence  by Lemma \ref{misuretta} 
 \[
 \begin{aligned}
 \meas(\ol{\matU}_{1/4}(\ell^{N,\io})\setminus  \tilde{\matC}_N(\e))&= \meas\left( 
 \bigcup_{\substack{\nu\in\ZZZ^{\ZZZ}_{f}\setminus\{0\} \\ \sum_{i\in\ZZZ} \nu_i=0}}
  \mathfrak R_\nu\left(\g\prod_{i\in\ZZZ}\frac{1}{(1+\jap{i}^2|\nu_i|^2)^{N+1}}\right)\right)\\
 &\le 4\gamma \sum_{\nu\in\ZZZ^{\ZZZ}_{f}\setminus\{0\}} \jap{i_0(\nu)}^N \prod_{i\in\ZZZ}\frac{1}{(1+\jap{i}^2|\nu_i|^2)^{N+1}}\\
 &\le 4\gamma\sum_{\nu\in\ZZZ^{\ZZZ}_{f}} \prod_{i\in\ZZZ}\frac{1}{(1+\jap{i}^2|\nu_i|^2)} \le C\gamma,
 \end{aligned}
 \]
where the last inequality is proved in ref.~\cite[Lemma 4.1]{Bjfa}; see also refs.~\cite{BMP1,MP} and Lemma \ref{misurazza}.

%%%%%%%%%%%%%%%%%%%%%%%%%%%%%%%%%%%%%%%%%%%%%%%%%%%%%%%%%%%%%%%%%%%%%%%%%% 
%%%%%%%%%%%%%%%%%%%%%%%%%%%%%%%%%%%%%%%%%%%%%%%%%%%%%%%%%%%%%%%%%%%%%%%%%% 
\appendix
%%%%%%%%%%%%%%%%%%%%%%%%%%%%%%%%%%%%%%%%%%%%%%%%%%%%%%%%%%%%%%%%%%%%%%%%%% 
%%%%%%%%%%%%%%%%%%%%%%%%%%%%%%%%%%%%%%%%%%%%%%%%%%%%%%%%%%%%%%%%%%%%%%%%%% 

%%%%%%%%%%%%%%%%%%%%%%%%%%%%%%%%%%%%%%%%%%%%%%%%%%%%%%%%%%%%%%%%%%%%%%%%%%
%%%%%%%%%%%%%%%%%%%%%%%%%%%%%%%%%%%%%%%%%%%%%%%%%%%%%%%%%%%%%%%%%%%%%%%%%%
\zerarcounters 
\section{On the Gevrey regularity} 
\label{gevgev} 
%%%%%%%%%%%%%%%%%%%%%%%%%%%%%%%%%%%%%%%%%%%%%%%%%%%%%%%%%%%%%%%%%%%%%%%%%%
%%%%%%%%%%%%%%%%%%%%%%%%%%%%%%%%%%%%%%%%%%%%%%%%%%%%%%%%%%%%%%%%%%%%%%%%%%

%%%%%%%%%%%%%%%%%%%%%%%%%%%%%%%%%%%%%%%%%%%%%%%%%%%%%%%%%%%%%%%%%%%%%%%%%%
\subsection{Gevrey regularity of the solution} \label{gevgevgev} 
%%%%%%%%%%%%%%%%%%%%%%%%%%%%%%%%%%%%%%%%%%%%%%%%%%%%%%%%%%%%%%%%%%%%%%%%%%

In order to discuss the regularity of the functions in $\mathtt{G}(s_1,s_2,\al)$ it is convenient to introduce some notation.
Let us define two families of {\it thickened tori} by setting, for $a>0$,
\[
\TTT^\ZZZ_a:=\{\f\in \CCC^\ZZZ:\quad \mbox{Re}\f\in \TTT^\ZZZ\,,\quad \|\mbox{Im}\f\|_\infty < a \},
\]
and, for $s,\al>0$,
\[
\TTT^\ZZZ_{s,\al}:=\{\f\in \CCC^\ZZZ:\quad \mbox{Re}\f\in \TTT^\ZZZ\,,\quad \sup_{j\in\ZZZ}\jap{j}^\al|\mbox{Im}\f_j|< s \}\,.
\]
The thickened tori defined above, endowed with the metric
\[
{\rm dist}(\varphi,\psi):= \sup_{j\in\ZZZ} \bigl(|\mbox{Re}(\varphi_j-\psi_j)\,{\rm mod}{2\pi}|+|\mbox{Im}(\varphi_j-\psi_j)|\bigl)\,,
\]
 are Banach manifolds modelled on $\ell^{\io}$.
Thus one can define  the differential operators $\del_{\f_j}$ in the natural way, and characterize the 
 analytic functions as closure w.r.t.~the uniform topology of the set of trigonometric polynomials depending on a finite number of angles \cite{MP}.

%%%%%%%%%%%%%%%%%%%%%%%%%%%%%%%%%%%%%%%%%%%%%%%%%%%%%%%%%%%%%%%%%%%%%%%%%% 
\begin{lemma}\label{ginocchio}
	For any $s_1,s_2>0$, any $\al\in(0,1)$, any $\om\in \gQ$ and any $U\in \mathtt{G}(s_1,s_2,\al)$,
	writing $U(x,\f)$ as in \eqref{fury} the following holds:
	\begin{enumerate}[topsep=0ex]
\itemsep0em
		\item for all $x\in\RRR$ the map $t\mapsto U(x,\omega t)$ is Gevrey with index $2/\al$,
		\item for all $t\in\RRR$ the map $x\mapsto U(x,\omega t)$ is Gevrey with index $1/\al$,
		\item for all $s'< s_1$, $U$ defines an analytic  map $\mathfrak{i}_U:\TTT^\ZZZ_{s',\al}\to \mathtt{g}(s_2,\al)$ by setting
		\begin{equation}
		\label{varvara}
		\mathfrak{i}_U(\f):= \{	u_j(\f)\}_{j\in \ZZZ}\,,\qquad  u_j(\f):=\sum_{ \nu\in \ZZZ^\ZZZ_{f}}u_{j,\nu} e^{\ii \f\cdot\nu}\,.
		\end{equation}
	\end{enumerate}
\end{lemma}
%%%%%%%%%%%%%%%%%%%%%%%%%%%%%%%%%%%%%%%%%%%%%%%%%%%%%%%%%%%%%%%%%%%%%%%%%% 

%%%%%%%%%%%%%%%%%%%%%%%%%%%%%%%%%%%%%%%%%%%%%%%%%%%%%%%%%%%%%%%%%%%%%%%%%% 
 \prova
First note that for $t\in \RRR$ one has $\{\omega_j t\hbox{ mod }2\pi\}_{j\in\ZZZ}\in \TTT^\ZZZ_{s',\al}$ for all $s'\ge 0$. Thus, if we assume item 3, 
then for all $t\in \RRR$  one has $\mathfrak{i}_U(\omega t)\in \mathtt{g}(s_2,\al)$. Consequently 
if $\calF$ is the Fourier transform operator,  %we use the identity 
then  $\calF^{-1}\mathfrak{i}_U(\om t)$ is  Gevrey with index $1/\al$, so that
item 2 follows by the identity $U(\cdot,\f)=\calF^{-1}\mathfrak{i}_U(\f)$.
Therefore we are left to prove items 1 and 3 only.
 
 To prove item 1, note that 
 \[
\Biggl|\del^k_t \sum_{j\in\ZZZ}  \sum_{\nu\in\ZZZ^{\ZZZ}_{f}} u_{j,\nu} \, e^{\ii jx} e^{\ii\nu\cdot\om t }\Biggr|
  \le \sum_{j\in\ZZZ} \sum_{\nu\in\ZZZ^{\ZZZ}_{f}} |\om\cdot\nu|^k |u_{j,\nu} | 
 \le 2^k \sum_{j\in\ZZZ} \sum_{\nu\in\ZZZ^{\ZZZ}_{f}}  |u_{j,\nu}| |\nu|_\alpha^{2k/\alpha},
 \]
where we used that, for  $\om\in\gQ$, one has
 \[
 |\om\cdot\nu| \le 2 |\nu|_2 \le 2 |\nu|_\alpha^{2/\alpha} \,.
 \]
 On the other hand, using that
 $
 %2^k\sum_{j\in\ZZZ} \sum_{\nu\in\ZZZ^{\ZZZ}_{f}} 
 |u_{j,\nu}| %(|\nu|_\alpha^{2/\alpha})^k 
 \le %2^k 
 \|U\|_{s_1,s_2,\al} %\sum_{j\in\ZZZ} \sum_{\nu\in\ZZZ^{\ZZZ}_{f}} |\nu|_\alpha^{2 k /\alpha} 
 e^{-s_1|\nu|_\al}e^{-s_2\jap{j}^\al}$,
we find,  for a suitable constant $C=C(s_1,s_2,\al)$,
 \[
 |\del^k_t  U(x,\om t)|\le %C^k (2 k /\alpha)^{2k/\alpha } \le 
 C^k (k!)^{2/\alpha }\,,
 \]
which shows that the map $t\mapsto U(x,\omega t)$ is Gevrey with index $2/\al$.
 
 Regarding item 3, for any  
 %$\nu\in \ZZZ_{f}^\ZZZ$ and any $s'>0$ we have
 %\[
 %\sup_{\f\in \TTT^\ZZZ_{s',\al}}|e^{\ii \f\cdot\nu}|\le e^{s'|\nu|_\al}\,.
 %\]
% Then, for  
$s'\in(0,s_1)$ we have
 \[
\sup_{j\in\ZZZ} \sup_{\f\in \TTT^\ZZZ_{s',\al}}e^{s_2\jap{j}^\al}|u_j(\f)|
\le  \sup_{j\in\ZZZ} e^{s_2\jap{j}^\al}\sum_{ \nu\in \ZZZ^\ZZZ_{f}}|u_{j,\nu}|e^{s'|\nu|_\al}
\le \|U\|_{s_1,s_2,\al} \sum_{ \nu\in \ZZZ^\ZZZ_{f}}e^{-(s_1-s')|\nu|_\al}<\io\,,
 \]
implying analyticity of the map $\mathfrak{i}_U$.
\EP
%%%%%%%%%%%%%%%%%%%%%%%%%%%%%%%%%%%%%%%%%%%%%%%%%%%%%%%%%%%%%%%%%%%%%%%%%% 

%%%%%%%%%%%%%%%%%%%%%%%%%%%%%%%%%%%%%%%%%%%%%%%%%%%%%%%%%%%%%%%%%%%%%%%%%%   
\begin{lemma}\label{previous}
Fix $s>0$, $\al\in(0,1)$ and $\om\in\gotB^{(0)}$.
Let $\e_0=\e_0(s,\al,\om)$ be as in Theorem \ref{moser}, and
let $\mathfrak{i}_U$ be the map defined in \eqref{varvara}, with $U(x,\f)= U(x,\f;\e,\om,c)$ as in Theorem \ref{moser}. 
Then
for all $\e\in(-\e_0,\e_0)$ there exists $a=a(\e)$ such that $\mathfrak{i}_U$ is analytic as a map $\TTT^\ZZZ_a\to \mathtt{g}(s,\al)$.
\end{lemma}
%%%%%%%%%%%%%%%%%%%%%%%%%%%%%%%%%%%%%%%%%%%%%%%%%%%%%%%%%%%%%%%%%%%%%%%%%%   

%%%%%%%%%%%%%%%%%%%%%%%%%%%%%%%%%%%%%%%%%%%%%%%%%%%%%%%%%%%%%%%%%%%%%%%%%%   
\prova
Decompose $U=U_0+U_\perp$ as in Theorem \ref{moser}, with $U_0$ as in \eqref{sollinear}. Since $\mathfrak{i}_{U_0}$ is trivially analytic as a map 
$\TTT^\ZZZ_a\to \mathtt{g}(s,\al)$, we only need to prove that the series defining $\mathfrak{i}_{U_\perp}$ is totally convergent, i.e. that
%
%, we restrict ourself to study $\mathfrak{i}_{U_\perp}$ and, as in the proof of item 3.~in Lemma \ref{ginocchio}, we only need to prove total convergence, i.e. 
\[
\mathfrak{S} := \sup_{j\in\ZZZ} e^{s\jap{j}^\al}\sum_{ \nu\in \ZZZ^\ZZZ_{f}\setminus\{\gote_j\}}|u_{j,\nu}|e^{a\|\nu\|_{{1}}}<\io\,.
\]
Now \eqref{purequesta}, \eqref{settou} and the bound \eqref{sudore}, with $s_1=0$, $s'=4s/5$ and $s_3= {s'}/4$, imply 
\begin{align*}
\mathfrak{S} &\le  \sup_{j\in\ZZZ} e^{s\jap{j}^\al}\sum_{k\ge1}\e^k\sum_{ \nu\in \ZZZ^\ZZZ_{f}\setminus\{\gote_j\}} 
\sum_{ \vartheta \in \gotF^{(k)}_{j,\nu,+} }|\Val^\ttR({\vartheta};c,\om)|  e^{a|\Lambda(\breve{\vartheta})|}\\ 
& \le \sup_{j\in\ZZZ}  \sum_{k\ge1}\e^k\sum_{ \nu\in \ZZZ^\ZZZ_{f}\setminus\{\gote_j\}} \sum_{ \vartheta \in \gotF^{(k)}_{j,\nu,+} }
A^{|N(\vartheta)|}e^{a(4k+1)} e^{-s' J({\breve\vartheta})}\!\!\! \prod_{v \in N_2(\vartheta)}  e^{- \frac{3s'}{4} J({\breve\vartheta_v})} \\ 
& \le \sup_{j\in\ZZZ}  \sum_{k\ge1}\e^k A^{2k}  e^{a(4k+1)} 
\sum_{ \vartheta \in{{\bigcup_{\nu\ne \gote_j}} }\gotF^{(k)}_{j,\nu,+} } e^{- \frac{3s}5 \calJ(\vartheta)}.
\end{align*}

Now one procceds as in the proof of Lemma \ref{convergerebbe2}, with the following modifications.
For any $k\ge 1$, $j\in\ZZZ$, $\nu\in\ZZZ^\ZZZ_f\setminus\{\gote_j\}$ and $\s\in\{\pm\}$, %$k,j,\nu,\s$,
given a tree $\vartheta\in %\gotT^{(k)}_{j,\nu,\s}\cup
\gotF^{(k)}_{j,\nu,\s}$, 
let $\calmX^{\rm Tot}(\vartheta)$ denote the set of all trees $\vartheta'\in %\gotT^{(k)}_{j,\nu,\s}\cup
\cup_{\nu'\ne \gote_j}\gotF^{(k)}_{j,\nu',\s}$ which are obtained
from $\vartheta$ by  changing the component labels $j_\la$ but not the sign labels $\s_\la$ of the leaves $\la\in \Lambda(\vartheta)$, 
and changing all the other labels consistently. Following the proof of Lemma \ref{convergerebbe2} verbatim, with $s''=3 s/5$, one finds
\[
\sum_{\vartheta' \in \calmX^{\rm Tot}(\vartheta)} e^{- \frac{3s}5 \calJ(\vartheta')} \le D_2^k .
\]
Then one has to perform the sum over all the sign labels of the leaves, and over all the unlabelled trees of order $k$:
proceeding as in the proof of Lemma \ref{convergerebbe2}, with $s_2=s$, we obtain
\[
\mathfrak{S} \le \sum_{k\ge1}\e^k D_0^k e^{a(4k+1)}
\]
with $D_0=D_0(s,\al,\om)$, by using the notation in Proposition \ref{voltabuona}.
Since, by definition $\e_0= D_0^{-1}$, in conclusion the sum  converges provided that
$|\e| \e_0^{-1} e^{4a}  <1$. This fixes the value of $a=a(\e)$.
\EP
%%%%%%%%%%%%%%%%%%%%%%%%%%%%%%%%%%%%%%%%%%%%%%%%%%%%%%%%%%%%%%%%%%%%%%%%%% 

%%%%%%%%%%%%%%%%%%%%%%%%%%%%%%%%%%%%%%%%%%%%%%%%%%%%%%%%%%%%%%%%%%%%%%%%%%
\subsection{Normal separate analyticity of the solution} \label{gevgevgevgev} 
%%%%%%%%%%%%%%%%%%%%%%%%%%%%%%%%%%%%%%%%%%%%%%%%%%%%%%%%%%%%%%%%%%%%%%%%%%

Let us conclude by showing that, for all $\om\in \gotB^{(0)}$ and all $\e\in(-\e_0,\e_0)$,
 the map from $\matU_1(\mathtt{g}(s,\al)) $ to $ \mathtt{G}(s_1,s_2,\al)$ defined by $c\mapsto U(\cdot,\cdot;c,\om,\e)$
%$c\mapsto U(\cdot,\cdot;c,\om,\e)$ 
is normally separately analytic.
Once more, as in the proof of Lemma \ref{previous},
we can confine ourselves to consider $U_\perp$.
We start by writing the Taylor series expansion in the variables $c$ of the coefficients $u_{j,\nu}(\e,\om,c)$ in \eqref{purequesta}.

First of all, we set
 \[
 \ZZZ^{\io}_{+,f}:=\{ a \in \ZZZ^\io_f \,:\, a_j\ge0\ \forall\,j\in\ZZZ\},
 \]
and, for any $a,b\in\ZZZ^\io_{+,f}$, define
\[
j(a,b):= \sum_{i\in\ZZZ} i(a_i-b_i)\,,
%\quad 
%\gotI(a,b):=\{(j,+) \mbox{repeated $a_j$ times}\}\cup  \{(j,-) \mbox{repeated $b_j$ times}\}.
\]
Next, for any tree $\vartheta\in \gotF$ and $\s=\pm$, we set $\gotI_\s(\vartheta)=\{ |\Lambda_{j,\s} (\vartheta)|\}_{j\in\ZZZ}$.
Then we define %$U^{(a,b)}(\om,\e) =0$ if $a-b= \gote_{j}$ and
\[
U^{(a,b)}(\om,\e) := \begin{cases}
0  ,\phantom{\displaystyle{\sum_{N=M}}} & \qquad a-b= \gote_{j(a,b)} , \\
\displaystyle{ \sum_{k \ge 1}\e^k \sum_{ \substack{\vartheta\in \gotF^{(k)}_{j(a,b),a-b,+} \\ 
\gotI_+(\vartheta)=a , \, \gotI_-(\vartheta)=b }}
%\{j_\lambda,\s_\lambda\}_{\lambda\in \Lambda(\vartheta)}=\gotI(a,b)}}
\Val^\ttR(\vartheta;\uno,\om) } ,  & \qquad a-b \neq  \gote_{j(a,b)} ,
\end{cases}
\]
%otherwise,
where $\uno\in\ell^\io$ is such that $\uno_j=1$ for all $j\in\ZZZ$ (in other words, with respect to
$\Val^\ttR(\vartheta;c,\om)$, in $\Val^\ttR(\vartheta;\uno,\om)$ all the leaf factors are replaced with one).
Then, recalling \eqref{purequesta} and \eqref{settou}, one can write
\[
u_{j,\nu}(c,\om,\e)= \sum_{ \substack{a,b\in \ZZZ^\ZZZ_{+,f}\\ a-b = \nu , \, j(a,b)=j}} 
U^{(a,b)}(\om,\e) \, c^a \bar c^b \,,\qquad  c^a \bar c^b:=\prod_{i\in\ZZZ} c_i^{a_i}\bar c_i^{b_i}\,.
\] 

The function $U$ is normally separately analytic in $c,\ol{c}$ for $c\in\matU_1(\mathtt g(s,\al))$ 
if the Cauchy majorant 
\[
\und{U}(x,\f;c,\om,\e):= \sum_{j\in\ZZZ} \sum_{{\nu\in\ZZZ^{\ZZZ}_{f}\setminus{ \{\gote_j\}}}} 
\und{u}_{j,\nu}(c,\om,\e) \, e^{\ii (j x+\nu\cdot\f) },
\]
with
\[
\und{u}_{j,\nu}(c,\om,\e):= \sum_{ \substack{a,b\in \ZZZ^\ZZZ_{+,f}\\ a-b =\nu, \, j(a,b)=j}} |U^{(a,b)}(\om,\e)| c^a \bar c^b,
\]
belongs to $\mathtt G(s_1,s_2,\al)$ and is separately analytic in $c,\ol{c}$ for $c\in\matU_1(\mathtt g(s,\al)).$ 
To show that  the latter property holds, it is sufficient to prove that
\begin{equation} \label{normanal}
\und{\mathfrak{S}} :=\sup_{\substack{j\in\ZZZ \\ \nu\in\ZZZ^\io_f}} e^{s_1|\nu|_\al} e^{s_2\jap{j}^\al}
\sum_{ \substack{a,b\in \ZZZ^\ZZZ_{+,f}\\ a-b = \nu , \, j(a,b)=j}} 
|U^{(a,b)}(\e,\om)| \sup_{c,\ol{c}\in \matU_1(\mathtt{g}(s,\al))}|c^a \bar c^b |<\io.
\end{equation}
On the other hand one has
\begin{align*}
\und{\mathfrak{S}}
&\le \sup_{\substack{j\in\ZZZ \\\nu\in\ZZZ^\io_f}} e^{s_1|\nu|_\al} e^{s_2\jap{j}^\al} 
\sum_{ \substack{a,b\in \ZZZ^\ZZZ_{+,f}\\ a-b=\nu , \, j(a,b)=j}} \sum_{k \ge 1}|\e|^k 
\sum_{ \substack{\vartheta\in \gotF^{(k)}_{j(a,b),a-b,+} \\ 
\gotI_+(\vartheta)=a , \, \gotI_-(\vartheta)=b }}
%\{j_\lambda,\s_\lambda\}_{\lambda\in \Lambda(\vartheta)} =\gotI(a,b)}}
\sup_{c\in \matU_1(\mathtt{g}(s,\al))}|\Val^\ttR(\vartheta;c,\om)|
\\
&\le \sup_{\substack{j\in\ZZZ \\ \nu\in\ZZZ^\io_f}} e^{s_1|\nu|_\al} e^{s_2\jap{j}^\al}
\sum_{k \ge 1}|\e|^k \sum_{ \substack{\vartheta\in \gotF^{(k)}_{j,\nu,+} }}\sup_{c\in \matU_1(\mathtt{g}(s,\al))}|\Val^\ttR(\vartheta;c,\om)|.
\end{align*}
Using the bounds \eqref{sudore} and proceeding as in the proof of Lemma \ref{previous} we obtain \eqref{normanal}.

%%%%%%%%%%%%%%%%%%%%%%%%%%%%%%%%%%%%%%%%%%%%%%%%%%%%%%%%%%%%%%%%%%%%%%%%%%
%%%%%%%%%%%%%%%%%%%%%%%%%%%%%%%%%%%%%%%%%%%%%%%%%%%%%%%%%%%%%%%%%%%%%%%%%%
\zerarcounters 
\section{On the Bryuno condition} 
\label{confronto} 
%%%%%%%%%%%%%%%%%%%%%%%%%%%%%%%%%%%%%%%%%%%%%%%%%%%%%%%%%%%%%%%%%%%%%%%%%%
%%%%%%%%%%%%%%%%%%%%%%%%%%%%%%%%%%%%%%%%%%%%%%%%%%%%%%%%%%%%%%%%%%%%%%%%%%

%%%%%%%%%%%%%%%%%%%%%%%%%%%%%%%%%%%%%%%%%%%%%%%%%%%%%%%%%%%%%%%%%%%%%%%%%%
\subsection{The finite-dimensional case}
%%%%%%%%%%%%%%%%%%%%%%%%%%%%%%%%%%%%%%%%%%%%%%%%%%%%%%%%%%%%%%%%%%%%%%%%%%

The Bryuno condition was first introduced in ref.~\cite{Bry} for $\om\in\RRR^d$ as 
\begin{equation}\label{originale1}
\sum_{m\ge0}\frac{1}{2^m}\log \left( \frac{1}{\be_\om(2^m)} \right) <\io,
\end{equation}
with $\be_{\om}(x)$ defined as in \eqref{beta} with $\al=0$.\footnote{
Clearly in the finite-dimensional case the value of $\al$ is irrelevant in the light of the norm equivalence.}
In particular, for $\om\in\RRR^2$, up to a rotation and a dilatation, one can take $\om=(1,\la)$, with 
$\la\in\RRR\setminus\QQQ$, and, 
if $q_k$ denotes the denominator of the $k$-th convergent of the continuous fraction of $\la$, then
the condition \eqref{originale1} is equivalent to
\begin{equation}\label{originale}
\sum_{k\ge0}\frac{1}{q_k}\log(q_{k+1}) < \io.
\end{equation}

In refs.~\cite{Y1,Y2} Yoccoz shows that \eqref{originale} is the necessary and sufficient condition in order to 
have an analytically linearizable map $f:\CCC\to\CCC$, $f(z)=e^{2\pi i\theta}z + O(z^2)$.
For $d\ge3$ the sharpness of \eqref{originale1} in the analytic setting remains an open question; see also
ref.~\cite{Bounemoura}, where other related issues are discussed.

For $\om\in\RRR^d$,
in ref.~\cite{BF}, the so-called $\al$-Bryuno-R\"ussmann condition, i.e.~the condition that
\begin{equation}\label{buofe}
\calF_{\al}(\om):=\int_1^{+\io} \frac{1}{Q^{1+\frac{1}{\al}}}\log\left(\frac{1}{\be_\om(Q)}\right) \der Q <\io,
\qquad \al\ge 1,
\end{equation}
is used to obtain KAM tori in Gevrey class.  The set $BR_\al:=\{\om\in\RRR^d:\calF_\al(\om)<\io\}$ decreases
w.r.t.~$\al$, and for $\al=1$ one obtains KAM tori in the analytic case, as shown by R\"ussmann in ref.~\cite{Rus}.
In the finite-dimensional case, the condition \eqref{originale1} is equivalent to \eqref{buofe} with $\al=1$.
Indeed if in \eqref{buofe} one takes $\al=1$ and considers the substitution
$Q=2^{x}$, one finds
\[
\calF_1(\om)=\log 2 \int_0^{+\io} \frac{1}{2^x} \log\left(\frac{1}{\be_\om(2^x)}\right) \der x ,
\]
which converges if and only if \eqref{originale1} holds.

As already mentioned,
the condition \eqref{originale1} corresponds to \eqref{juno} with $\al=0$ and $r_m=2^{m}$,
thus, at least in the finite-dimensional case, the Bryuno condition in Definition \ref{br} is implied by the
standard Bryuno condition \eqref{originale1}.
The results in refs.~\cite{Y1,Y2} imply that for $d=2$ the condition \eqref{juno} is equivalent to \eqref{originale1}. 
In fact, the two conditions \eqref{juno} and \eqref{originale1} are equivalent
also in dimension $d\ge3$. Indeed, let $r\in\gotR$ be a sequence such that $B_{\om}(r)<\io$.
%By Remark \ref{brsubsequence} we may assume that $\be_\om(r_{m-1})< 2\be_\om(r_{m})$ for all $m\ge 1$.
Then we can write
\[ 
\begin{aligned}
\sum_{n\ge 0} \frac{1}{2^n} \log \frac{1}{\be_\om(2^n)} & =
\sum_{m\ge 1} \sum_{\substack{n\in \NNN \\ r_{m-1} \le 2^n < r_{m}}} \frac{1}{2^n} \log \frac{1}{\be_\om(2^n)}
\le \sum_{m\ge 1} \sum_{\substack{n\in \NNN \\ r_{m-1} \le 2^n < r_{m}}} \frac{1}{2^n} \log \frac{1}{\be_\om(r_{m})} \\
& \le \sum_{m\ge 1} \sum_{k \ge 0} \frac{1}{2^k r_{m-1}} \log \frac{1}{\be_\om(r_{m})} \le
2\BB_\om(r) .
\end{aligned}
\]

%%%%%%%%%%%%%%%%%%%%%%%%%%%%%%%%%%%%%%%%%%%%%%%%%%%%%%%%%%%%%%%%%%%%%%%%%%
\subsection{The infinite-dimensional case}
%%%%%%%%%%%%%%%%%%%%%%%%%%%%%%%%%%%%%%%%%%%%%%%%%%%%%%%%%%%%%%%%%%%%%%%%%%
Let us come back to the infinite-dimensional case. We
prove that the frequency vectors satisfying the Diophantine condition in Definition \ref{dio}
satisfy automatically the Bryuno condition too.
As a first step, the following result, adapted from~\cite[Lemma C.2]{MP},
provides a lower bound on $\be^*(r,\g)$, defined in \eqref{notte}, in terms of $r>0$.

%%%%%%%%%%%%%%%%%%%%%%%%%%%%%%%%%%%%%%%%%%%%%%%%%%%%%%%%%%%%%%%%%%%%%%%%%%
\begin{lemma} \label{eradiovsbrj} 
For any $\al\in(0,1)$ and any $r>0$, one has
\begin{equation} \label{stimobeta}
\be^*(r,\g) \ge \g \left(1+r\right)^{-C(\al)r^{\frac{1}{1+\al/2}}} ,
\end{equation}
for some constant $C(\al)>0$.
\end{lemma}
%%%%%%%%%%%%%%%%%%%%%%%%%%%%%%%%%%%%%%%%%%%%%%%%%%%%%%%%%%%%%%%%%%%%%%%%%%

%%%%%%%%%%%%%%%%%%%%%%%%%%%%%%%%%%%%%%%%%%%%%%%%%%%%%%%%%%%%%%%%%%%%%%%%%%
\prova
For $\nu\in\ZZZ^\ZZZ_f$ such that $0<|\nu|_\al \le r$, let $k(\nu)$ be the number of non-vanishing components of $\nu$;
in other words $k(\nu):=|\{ \nu_i \neq 0 : i  \in \ZZZ\}|$.
One has
\[ 
r \ge |\nu|_\al = \sum_{\substack{ i\in\ZZZ \\ \nu_i \neq 0}} \jap{i}^\al |\nu_i| \ge
 \sum_{\substack{ i\in\ZZZ \\ \nu_i \neq 0}} \jap{i}^\al \ge \sum_{i=1}^{k(\nu)} |i|^\al \ge C (k(\nu))^{1+\al} ,
 \]
for some $C>0$ independent of $\al$, and hence $k(\nu) \le C \, r^{\frac{1}{1+\al}}$, 
for some other $\al$-independent positive constant $C$. On the other hand, for $\al\in(0,1)$, one has
\[
\log( 1 + \jap{i}^2|\nu_i|^2)  \le \log( 1 + (\jap{i}^\al |\nu_i|)^{\frac2\al})  \le \log (1 + r^{\frac{2}{\al}}) \le K_0(\al)  \log (1+r) ,
\]
with $K_0(\al)=O(\al^{-1})$, so that, for any or $\nu\in\ZZZ^{\ZZZ}_{f}$ such that $|\nu|_\al \le r$ one finds
\begin{equation} \label{mp}
\log \prod_{i\in\ZZZ}(1+\jap{i}^{2} |\nu_i|^{2}) \le 
\sum_{\substack{ i\in\ZZZ \\ \nu_i \neq 0}} \log( 1 + \jap{i}^2|\nu_i|^2) \le K(\al) \, r^{\frac{1}{1+\al}} \log (1+r)
\end{equation}
and hence
\begin{equation} \label{taglio}
\sup_{0<|\nu|_{\al}<r} 
\prod_{i\in\ZZZ}(1+\jap{i}^{2} |\nu_i|^{2}) \le (1+r)^{K(\al)r^{\frac{1}{1+\al}}} ,
\end{equation}
for some constant $K(\al)>0$ proportional to $K_0(\al)$.
The bound \eqref{taglio} yields immediately \ref{stimobeta} with $C(\al)=K(\al/2)$.
\EP
%%%%%%%%%%%%%%%%%%%%%%%%%%%%%%%%%%%%%%%%%%%%%%%%%%%%%%%%%%%%%%%%%%%%%%%%%%

%%%%%%%%%%%%%%%%%%%%%%%%%%%%%%%%%%%%%%%%%%%%%%%%%%%%%%%%%%%%%%%%%%%%%%%%%%
\begin{rmk} \label{rmkaaa}
\emph{
In fact the bound \eqref{stimobeta} holds also for $\al \ge 1$, with slight changes in the proof~\cite{MP}.
By looking at the proof of Lemma \ref{eradiovsbrj} we see that $C(\al)\to\io$ proportionally to $\al^{-1}$ as $\al\to0$.
}
\end{rmk}
%%%%%%%%%%%%%%%%%%%%%%%%%%%%%%%%%%%%%%%%%%%%%%%%%%%%%%%%%%%%%%%%%%%%%%%%%%

Then, with the notation in Definitions \ref{br} and \ref{dio0}, one has the following result,
which is the analogous of Lemma 2 in ref.~\cite{Ge} in the finite-dimensional case.

%%%%%%%%%%%%%%%%%%%%%%%%%%%%%%%%%%%%%%%%%%%%%%%%%%%%%%%%%%%%%%%%%%%%%%%%%%
\begin{lemma}\label{diovsbrj} 
One has 
\[
\gD(\g,\tau)\subseteq \gD^{(0)}(\g,\tau)  \subseteq \gotB(\gamma,\tau),\qquad \gD \subseteq \gotB.
\]
\end{lemma}
%%%%%%%%%%%%%%%%%%%%%%%%%%%%%%%%%%%%%%%%%%%%%%%%%%%%%%%%%%%%%%%%%%%%%%%%%%

%%%%%%%%%%%%%%%%%%%%%%%%%%%%%%%%%%%%%%%%%%%%%%%%%%%%%%%%%%%%%%%%%%%%%%%%%%
\prova
Of course $\gD(\g,\tau)\subseteq \gD^{(0)}(\g,\tau)$. Now
if $\om\in\gD^{(0)}(\g,\tau)$ and $r^*\in\gotR^*$, one has
\begin{equation} \label{b*b}
\be_{\om}^{(0)}(r_m^*) = 
\!\!\!\! \inf_{\substack{\nu\in\ZZZ^{\ZZZ}_{f} \\ \sum_{i\in\ZZZ}\nu_i=0 \\ 0 <|\nu|_{\al/2}\le r_m^* }} \!\!\!\!
|\omega\cdot \nu| 
> \g 
\!\!\!\! \inf_{\substack{\nu\in\ZZZ^{\ZZZ}_{f} \\ \sum_{i\in\ZZZ}\nu_i=0 \\ 0 <|\nu|_{\al/2}\le r_m^*}}
\prod_{i\in \ZZZ}\frac{1}{(1+ \jap{i}^{2}|\nu_i|^{2})^\tau} = \be^*(r_m^*,\g,\tau) .
\end{equation}
By inserting \eqref{stimobeta} into \eqref{adessoce} and using \eqref{b*b} we obtain
\begin{equation} \label{adessoceanchequesto}
\BB_\om^{(0)}(r^*) < \BB(r^*,\g,\tau) \le \log\frac{1}{\g} \left( \sum_{m\ge1}\frac{1}{r^*_{m-1}} \right)
+ C(\al) \left( \tau \sum_{m\ge1}\frac{(r_m^*)^{\frac{1}{1+\al/2}}}{r^*_{m-1}}  \log( 1+r_m^*) \right) ,
\end{equation}
where both sums are convergent if $r^*$ is suitably chosen and $\al>0$. For instance one may take $r^*_m=2^m$.
Thus the second inclusion is proved.
To prove that $ \gD \subseteq \gotB$ one reasons in the same way, simply dropping the constraint $\sum_{i\in\ZZZ}\nu_i=0$.
%, and the second follows naturally.
%Therefore the first inclusion is proved. The same argument can be used to prove the second inclusion.
\EP
%%%%%%%%%%%%%%%%%%%%%%%%%%%%%%%%%%%%%%%%%%%%%%%%%%%%%%%%%%%%%%%%%%%%%%%%%%

%%%%%%%%%%%%%%%%%%%%%%%%%%%%%%%%%%%%%%%%%%%%%%%%%%%%%%%%%%%%%%%%%%%%%%%%%%
\begin{rmk} \label{primanoncera}
\emph{
The proof of Lemma \ref{diovsbrj} shows as a byproduct that the set $\gotR^*$ is non empty and that
$\{2^m\}\in\gotR^*$. Furthemore \eqref{adessoceanchequesto} shows that the bound on $\BB(r^*,\g,\tau)$ -- and hence
the bound on $\BB_\om^{(0)}(r^*)$ as well -- diverges in the limit $\g\to0$ or $\al\to0$.
}
\end{rmk}
%%%%%%%%%%%%%%%%%%%%%%%%%%%%%%%%%%%%%%%%%%%%%%%%%%%%%%%%%%%%%%%%%%%%%%%%%%
 
The following measure estimate is proved by reasoning as in ref.~\cite[Lemma 4.1]{Bjfa}; see also refs.~\cite{BMP1,MP}.

%%%%%%%%%%%%%%%%%%%%%%%%%%%%%%%%%%%%%%%%%%%%%%%%%%%%%%%%%%%%%%%%%%%%%%%%%% 
\begin{lemma}\label{misurazza}
For any $\g>0$ and any $\tau > 1/2$, there is an absolute constant $C>0$ such that
\[
\meas(\gQ\setminus\gotD(\gamma,\tau)) \le C\gamma.
\]
\end{lemma}
%%%%%%%%%%%%%%%%%%%%%%%%%%%%%%%%%%%%%%%%%%%%%%%%%%%%%%%%%%%%%%%%%%%%%%%%%% 

%%%%%%%%%%%%%%%%%%%%%%%%%%%%%%%%%%%%%%%%%%%%%%%%%%%%%%%%%%%%%%%%%%%%%%%%%% 
\prova
%Thanks to Lemma \ref{diovsbrj}, it is enough to prove that $\meas(\gQ\setminus\gD(\gamma,\tau)) \le C\gamma$. 
Given $\nu\in\ZZZ^\ZZZ_f\setminus\{0\}$, setting
\[
\Omega(\nu):=\left\{ \om\in\gQ\,:\, |\om\cdot\nu| \le  \g \prod_{i\in \ZZZ} \frac{1}{(1+ \jap{i}^{2}|\nu_i|^{2})^\tau} \right\} , 
\]
one has
\[
\meas(\Omega(\nu)) \le C \g  \prod_{i\in \ZZZ} \frac{1}{(1+ \jap{i}^{2}|\nu_i|^{2})^\tau}  ,
\qquad
\meas(\gQ\setminus\gD(\gamma,\tau) \le \sum_{\nu\in\ZZZ^{\ZZZ}_{f}} \meas(\Omega(\nu)) ,
\] 
for some constant $C$. For any $\nu\in\ZZZ^{\ZZZ}_{f}\setminus\{0\}$, set
\[
s(\nu)=\min\{ |i| : \nu_i \neq 0\} ,
\qquad S(\nu)=\max\{ |i| : \nu_i \neq 0\}
\]
and, shorten for notational simplicity, $a_{i}(x):=1/(1+ \jap{i}^{2}x^{2})$ and $b_i(x):=1/(\jap{i}^2x^2)$.
Then 
\[ 
\begin{aligned}
& \sum_{\nu\in\ZZZ^{\ZZZ}_{f}\setminus\{0\}} \meas(\Omega(\nu)) \le C \g \sum_{s\ge0} \sum_{S\ge s}
\sum_{\substack{ \nu\in\ZZZ^{\ZZZ}_{f}\setminus\{0\}\\ S(\nu) = S, \,s(\nu)=s}}  \prod_{i\in \ZZZ} (a_i(\nu_i))^\tau \\
& \quad \le C\g \sum_{s\ge0}
\Biggl( \; \sum_{\substack{h\in\ZZZ\setminus\{ 0\}} } (a_s(h))^{\tau}  \Biggr)
\Biggl( \; \sum_{h\in\ZZZ} (a_s(h))^ {\tau}  \Biggr) \Biggl[1+\\
&\qquad\qquad+
 \sum_{S\ge s+1} 
\Biggl( \; \sum_{\substack{h\in\ZZZ\setminus\{ 0\}} }(a_S(h))^{\tau}  \Biggr)
\Biggl( \; \sum_{h\in\ZZZ} (a_S(h))^ {\tau}  \Biggr)
\prod_{k=s+1}^{S-1} \Biggl( \; \sum_{h \in\ZZZ} (a_k(h))^\tau \Biggr)^{\!\!2}\Biggr]  \\
& \quad \le C\g  \sum_{s\ge0} \Biggl( \; \sum_{\substack{h\in\ZZZ\setminus\{ 0\}} } \frac{1}{(\jap{s}h)^{2\tau}} \Biggr)
\Biggl( 1 + \!\!\!  \sum_{\substack{h\in\ZZZ\setminus\{ 0\}} } \frac{1}{(\jap{s}h)^{2\tau}}\Biggr) 
\Biggl[ 1 +\\
&\qquad\qquad +  \sum_{S\ge s+1}
\Biggl( \;  \sum_{\substack{h\in\ZZZ\setminus\{ 0\}} } \frac{1}{({S}h)^{2\tau}} \Biggr)
\Biggl( 1 + \!\!\!  \sum_{\substack{h\in\ZZZ\setminus\{ 0\}} } \frac{1}{({S}h)^{2\tau}} \Biggr) 
\prod_{k=s+1}^{S-1} \Biggl( 1 + \!\!\!  \sum_{\substack{h\in\ZZZ\setminus\{ 0\}} } \frac{1}{({k}h)^{2\tau}} \Biggr)^{\!\!2} \; \Biggr]  \\
& \quad \le C \g \sum_{s\ge0} \frac{1}{\jap{s}^{2\tau}} \Biggl( 1 + \frac{C}{\jap{s}^{2\tau}} \Biggr)
\Biggl( 1 + \!\!\! \sum_{S\ge s+1} \frac{1}{S^{2\tau}} \Biggl( 1 + \frac{C}{S^{2\tau}} \Biggr)
\prod_{k=s+1}^{S-1} \Biggl( 1 + \frac{C}{k^{2\tau}} \Biggr)^2 \Biggr) \\
& \quad \le C \g \sum_{s\ge0} \frac{1}{\jap{s}^{2\tau}} 
\Biggl( 1 + \!\!\! \sum_{S\ge 1} \frac{1}{S^{2\tau}} \prod_{k\ge1} \Biggl( 1 + \frac{C}{k^{2\tau}} \Biggr)^2 \Biggr) 
 \le C \g ,
\end{aligned}
\]
where we have
\begin{itemize}[topsep=0ex]
\itemsep0em
\item taken into account that for any $k>0$ there are two label $\nu_i$ with $|i|=k$ and at least one
between $\nu_{-s}$ and $\nu_s$ does not vanish, and same for $S$;
\item distinguished between the cases $S(\nu)=s(\nu)$ and $S(\nu)>s(\nu)$, and used that in the latter case
$S=S(\nu) \neq 0$;
\item called $C$ any constant independently of its value;
\item noted that all sums and products converge for $\tau > 1/2$.
\end{itemize}
This concludes the proof.
\EP
%%%%%%%%%%%%%%%%%%%%%%%%%%%%%%%%%%%%%%%%%%%%%%%%%%%%%%%%%%%%%%%%%%%%%%%%%% 

%%%%%%%%%%%%%%%%%%%%%%%%%%%%%%%%%%%%%%%%%%%%%%%%%%%%%%%%%%%%%%%%%%%%%%%%%% 
%%%%%%%%%%%%%%%%%%%%%%%%%%%%%%%%%%%%%%%%%%%%%%%%%%%%%%%%%%%%%%%%%%%%%%%%%% 
\zerarcounters 
\section{Technical lemmata} 
\label{tecnici} 
%%%%%%%%%%%%%%%%%%%%%%%%%%%%%%%%%%%%%%%%%%%%%%%%%%%%%%%%%%%%%%%%%%%%%%%%%% 
%%%%%%%%%%%%%%%%%%%%%%%%%%%%%%%%%%%%%%%%%%%%%%%%%%%%%%%%%%%%%%%%%%%%%%%%%% 

%%%%%%%%%%%%%%%%%%%%%%%%%%%%%%%%%%%%%%%%%%%%%%%%%%%%%%%%%%%%%%%%%%%%%%%%%% 
\subsection{Bourgain's lemmata}\label{dibou}
%%%%%%%%%%%%%%%%%%%%%%%%%%%%%%%%%%%%%%%%%%%%%%%%%%%%%%%%%%%%%%%%%%%%%%%%%% 

The following results are due to Bourgain \cite{B96, CLSY}.
We report here the proofs for completeness.

%%%%%%%%%%%%%%%%%%%%%%%%%%%%%%%%%%%%%%%%%%%%%%%%%%%%%%%%%%%%%%%%%%%%%%%%%% 
\begin{lemma}\label{constance}
For $p\ge3$, let $\{n_i\}_{i=1}^{p}\subset \NNN$ be an ordered sequence $n_1\ge n_2\ge \ldots\ge n_{p} \ge1$
such that 
\begin{equation}\label{somma}
n_1+\sum_{i=2}^{p} \s_i n_i =0\,,
\end{equation}
for some sequence $\{\s_i\}_{i=2}^{p}\subset\{-1,0,1\}^p$.
Then for any $0<\al<1$ one has
\begin{equation}\label{co}
\sum_{i=2}^{p} n_i^\al - n_1^\al \ge ({2 - 2^\al}) \sum_{i=3}^{p} n_i^\al\,.
\end{equation}
\end{lemma}
%%%%%%%%%%%%%%%%%%%%%%%%%%%%%%%%%%%%%%%%%%%%%%%%%%%%%%%%%%%%%%%%%%%%%%%%%% 

%%%%%%%%%%%%%%%%%%%%%%%%%%%%%%%%%%%%%%%%%%%%%%%%%%%%%%%%%%%%%%%%%%%%%%%%%% 
\prova
First of all note that  from \eqref{somma} it follows that
\[
n_1^\al \le \Bigl(\sum_{i=2}^{p} n_i \Bigr)^\al,
\]
so we have
\begin{equation}\label{questa}
\sum_{i=2}^{p} n_i^\al - n_1^\al  \ge \sum_{i=2}^{p} n_i^\al - \Bigl( \; \sum_{i=2}^{p} n_i \Bigr)^\al .
%= n_2^\al + \sum_{i=3}^{p+1} n_i^\al - \left(n_2  + \sum_{i=3}^{p+1} n_i \right)^\al.
\end{equation}
%This implies that if $n_2=\ldots=n_{p+1}$ then the assertion holds trivially.

{
For $k=2,\ldots,p$, define $D_k:=\{(x_k,\ldots,x_{p}) \in \RRR^{p-k+1} : x_k \ge \ldots \ge x_{p}>0 \}$.
Let us now introduce the auxiliary function
\[
\begin{aligned}
f(x_2,\dots,x_{p}):= & 
\sum_{i=2}^{p} x_i^\al - \Bigl( \; \sum_{i=2}^{p} x_i \Bigr)^\al - ({2 - 2^\al}) \sum_{i=3}^{p} x_i^\al\\
=& \;
x_2^\al +(2^\al - 1) \sum_{i= 3}^{p} x_i^\al - \Bigl(x_2+\sum_{i= 3}^{p} x_i \Bigr)^\al \,,
\end{aligned}
\]
and note that $\partial_{x_2} f> 0$ on $D_2$.   Then 
$$
f(x_2,x_3,x_4,\dots,x_{p})\ge  f(x_3,x_3,x_4,\dots,x_{p})=: f_3(x_3,x_4,\dots,x_{p})\,.
$$ 
It $p=3$ the result follows by direct inspection, since $f(x_3,x_3)>0$ in $D_2$.

Otherwise, for all $k=3,\ldots,p-1$ define
\begin{equation} \nonumber
\begin{aligned}
f_k(x_k,\dots ,x_{p}) &:= f(\underbrace{x_k,\dots,x_k}_{(k-1)\text{ times}},x_{k+1},\dots,x_{p})\\ 
&= (1+ (2^\al -1)(k-2))x_k^\al  +(2^\al-1)\sum_{i\ge  k+1}x_i^\al- \Big((k-1)x_k+ \sum_{i\ge  k+1}x_i\Big)^\al.
\end{aligned}
\end{equation}
Assume inductively that, in $D_2$, for some $3\le k<p$, one has 
\[
f(x_2,\dots,x_{p})\ge f_3(x_3,\dots,x_{p}) \ge \dots \ge f_{k}(x_k,\dots ,x_{p})\,.
\] 
Since $ 1+(2^\al-1)(k-2) \ge (k-1)^\al$ for all $k\geq 2$, we see that
$\partial_{x_k} f_k >0$ on $D_k$ and hence the minimum is attained 
at $x_k= x_{k+1}$. In conclusion one has
\[
f(x_2,\dots,x_{p})\ge f_{p}(x_{p}) = f(x_{p},\dots,x_{p})> 0,
\]
where the last inequality is easily checked.
Therefore the assertion follows.}
%Now for $x\ge1$ consider the auxiliary function
%\[
%f(x):=(1+x)^\al -x^\al
%\]
%and note that, since $f'(x)<0$ one has
%\begin{equation}\label{au}
%(1+x)^\al -x^\al \le f(1) = 2^\al -1.
%\end{equation}
%This in turn implies that for any $a\ge b >0$ we have
%\[
%(2-2^\al)b^\al + (a+b)^\al -a^\al-b^\al = b^\al\Big( 
%(1+\frac{a}{b})^\al - (\frac{a}{b})^\al - (2^\al-1)
%\Big)\le0
%\]
%i.e.
%\begin{equation}\label{ab}
%a^\al + b^\al - (a+b)^\al \ge (2-2^\al)b^\al
%\end{equation}
%which iterated implies the assertion.
\EP
%%%%%%%%%%%%%%%%%%%%%%%%%%%%%%%%%%%%%%%%%%%%%%%%%%%%%%%%%%%%%%%%%%%%%%%%%% 

%%%%%%%%%%%%%%%%%%%%%%%%%%%%%%%%%%%%%%%%%%%%%%%%%%%%%%%%%%%%%%%%%%%%%%%%%% 
\begin{lemma}\label{costanza}
Let $\{n_i\}_{i=1}^{\widehat{p}}\subseteq \NNN$ be an ordered sequence $n_1\ge n_2\ge \ldots\ge n_{\widehat{p}} \ge1$,
with $\widehat{p}\ge3$, such that
%there are $\s_i=\pm1$ so that
\begin{equation}\label{sommazzabi}
\left| \sum_{i=1}^{\widehat{p}} \s_i n_i^2 \right| \le p\,,
\end{equation}
for some $p\ge\widehat{p}$ and some sequence $\{\s_i\}_{i=1}^{\widehat{p}}\subset\{\pm 1\}^{\widehat{p}}$
such that $\s_1=\s_2$ if $n_1=n_2$.
Then  one has
\begin{equation}\label{biasco}
n_2\le n_1 \le 3 \sum_{i=3}^{\widehat{p}} n_i^2 + p -\widehat{p} .
\end{equation}
%for some absolute constant $C>0$.\footnote{According to Biasco one has $C=31$. According to Procesi $C=7$ (tipo). } 
\end{lemma}
%%%%%%%%%%%%%%%%%%%%%%%%%%%%%%%%%%%%%%%%%%%%%%%%%%%%%%%%%%%%%%%%%%%%%%%%%% 

%%%%%%%%%%%%%%%%%%%%%%%%%%%%%%%%%%%%%%%%%%%%%%%%%%%%%%%%%%%%%%%%%%%%%%%%%% 
\prova
Assume without loss of generality that $\s_1=1$. Then of course by \eqref{sommazzabi}
\[
\left| 
n_1^2+\s_2 n_2^2 +\sum_{i=3}^{\widehat{p}}\s_i n_i^2
\right| \le p
\]
and hence
\[
n_1^2 +\s_2n_2^2 \le \left|
\sum_{i=3}^{\widehat{p}}\s_i n_i^2
\right| + p\le 2 \sum_{i=3}^{\widehat{p}} n_i^2 + 2 + p-\widehat{p}\,.
\]

Now if $\s_2=1$ then 
$$
n_1^2 + n_2^2 \ge n_1^2+1\ge n_1+1.
$$
On the other hand, if $\s_2=-1$ then $n_1 >n_2$,
and hence
$$
n_1^2-n_2^2 = (n_1+n_2)(n_1-n_2) \ge n_1 + n_2 \ge n_1+1.
$$
In both cases the assertion follows. % with $C=3$.
\EP
%%%%%%%%%%%%%%%%%%%%%%%%%%%%%%%%%%%%%%%%%%%%%%%%%%%%%%%%%%%%%%%%%%%%%%%%%% 

%%%%%%%%%%%%%%%%%%%%%%%%%%%%%%%%%%%%%%%%%%%%%%%%%%%%%%%%%%%%%%%%%%%%%%%%%% 
\subsection{Lemmata on permutations}
%%%%%%%%%%%%%%%%%%%%%%%%%%%%%%%%%%%%%%%%%%%%%%%%%%%%%%%%%%%%%%%%%%%%%%%%%% 

Let $S_n$ denote the set of all permutations $\s$ of $(1,\ldots,n)$, with $n\in \NNN$ such that $n\ge 2$.
Given $x=(x_1,\ldots,x_n) \in\RRR^n$ such that $x_i\neq0$ $\forall i=1,\ldots,n$,
let $\Pi(x)$ be the set of the vectors $\pi_\s(x)=(x_{\s(1)},\ldots,x_{\s(n)})$
obtained by permuting the components of $x$.  Set
\begin{equation} \label{GnD}
D_i(x) = x_{1}+\ldots + x_{i} , \quad i=1,\ldots,n  , \qquad \qquad
G_n(x) = \prod_{i=1}^{n-1} \frac{1}{D_i(x)} .
\end{equation}
%

%%%%%%%%%%%%%%%%%%%%%%%%%%%%%%%%%%%%%%%%%%%%%%%%%%%%%%%%%%%%%%%%%%%%%%%%%% 
%\begin{defi}[\textbf{Allowed vector}]\label{def1}
Given $x=(x_1,\ldots,x_n) \in\RRR^n$ such that $x_i\neq0$ $\forall i=1,\ldots,n$,
a vector $\pi_\s(x) \in \Pi(x)$ is said to be \emph{allowed} if $D_i(\pi_\s(x)) \neq 0$ $\forall i=1,\ldots,n -1$.
%and that a permutation $\s\in S$ is \emph{allowed} if $\pi_\s(x)$ is allowed.
%\end{defi}
%%%%%%%%%%%%%%%%%%%%%%%%%%%%%%%%%%%%%%%%%%%%%%%%%%%%%%%%%%%%%%%%%%%%%%%%%% 

Let $\Pi_A(x)$ be the set of all allowed vectors of $\Pi(x)$, and set $\Pi_N(x):=\Pi(x)\setminus\Pi_A(x)$.

%%%%%%%%%%%%%%%%%%%%%%%%%%%%%%%%%%%%%%%%%%%%%%%%%%%%%%%%%%%%%%%%%%%%%%%%%% 
\begin{lemma} \label{1}
If $x\in\RRR^n$ is such that $x_i\ne0$ for all $i=1,\ldots,n$, and $\Pi_A(x)=\Pi(x)$, then one has
$$ \sum_{x'\in \Pi(x)} G_n(x') = (x_1+\ldots+x_n) \prod_{i=1}^{n}\frac{1}{x_i} . $$
\end{lemma}
 %%%%%%%%%%%%%%%%%%%%%%%%%%%%%%%%%%%%%%%%%%%%%%%%%%%%%%%%%%%%%%%%%%%%%%%%%% 
 
 %%%%%%%%%%%%%%%%%%%%%%%%%%%%%%%%%%%%%%%%%%%%%%%%%%%%%%%%%%%%%%%%%%%%%%%%%% 
\proof By induction on $n$. \EP
%%%%%%%%%%%%%%%%%%%%%%%%%%%%%%%%%%%%%%%%%%%%%%%%%%%%%%%%%%%%%%%%%%%%%%%%%% 

%%%%%%%%%%%%%%%%%%%%%%%%%%%%%%%%%%%%%%%%%%%%%%%%%%%%%%%%%%%%%%%%%%%%%%%%%%  
%\begin{defi}[\textbf{Minimum sum}]\label{def2}
Given $x=(x_1,\ldots,x_n) \in\RRR^n$ such that $x_i\neq0$ $\forall i=1,\ldots,n$,
define the \emph{minimum sum} as
$$ \m := \min\{ D_i(\pi_\s(x))  : \s \in S_n , \, i\in \{1,\ldots , n-1 \} , \,D_i(\pi_\s(x)) \neq 0 \}  , $$
and set
$E:= (0,\mu/2n]$.  %$E:= [0,\mu/2n]$. 
For $\e=(\e_1,\ldots,\e_n)\in E^n$, write $x+\e=(x_1+\e_1,\ldots,x_n+\e_n)$.
%\end{defi}
%%%%%%%%%%%%%%%%%%%%%%%%%%%%%%%%%%%%%%%%%%%%%%%%%%%%%%%%%%%%%%%%%%%%%%%%%% 

%%%%%%%%%%%%%%%%%%%%%%%%%%%%%%%%%%%%%%%%%%%%%%%%%%%%%%%%%%%%%%%%%%%%%%%%%% 
\begin{rmk} \label{rmk0}
\emph{
By construction $\Pi_A(x+\e)=\Pi(x+\e)$ when
$\e_i\in E$
%E_0:=E\setminus\{0\}$ 
$\forall i=1,\ldots,n$.
}
\end{rmk}
%%%%%%%%%%%%%%%%%%%%%%%%%%%%%%%%%%%%%%%%%%%%%%%%%%%%%%%%%%%%%%%%%%%%%%%%%% 

%%%%%%%%%%%%%%%%%%%%%%%%%%%%%%%%%%%%%%%%%%%%%%%%%%%%%%%%%%%%%%%%%%%%%%%%%% 
%\begin{defi}[\textbf{Proper subset}]\label{def3}
Let $n,p\in\NNN$ such that $n\ge 3$ and $p\in\{2,\ldots,n-1\}$, and let $x=(x_1,\ldots,x_n)$.
A proper subset $\{x_{i_1},\ldots,x_{i_p}\}$ of the unordered set $\{x_1,\ldots,x_n\}$ is said to be %\emph{$p$-minimal} if
\emph{minimal} if
$x_{i_1} + \ldots + x_{i_p}=0$
and $x_{j_1}+\ldots+x_{j_{p'}} \neq 0$ for all $p'<p$ and all 
$\{x_{j_1},\ldots,x_{j_{p'}}\} \subset \{x_{i_1},\ldots,x_{i_p}\}$.
The set of all minimal subsets of $\{x_1,\ldots,x_n\}$ is denoted by $\Xi(x)$.
%%We say that a proper subset of $\{x_1,\ldots,x_n\}$ is \emph{minimal} if it is $p$-minimal for some $p=2,\ldots,n-2$.
%\end{defi}
%%%%%%%%%%%%%%%%%%%%%%%%%%%%%%%%%%%%%%%%%%%%%%%%%%%%%%%%%%%%%%%%%%%%%%%%%% 

Given any minimal %$p$-minimal 
subset $X=\{x_{i_1},\ldots,x_{i_p}\} \in \Xi(x)$, let $A(X)$ denote the subset of vectors $\pi_\s(x) \in \Pi(x)$
 such that $(\s(1),\ldots,\s(p))$ is a permutation of $(i_1,\ldots,i_p)$.

%%%%%%%%%%%%%%%%%%%%%%%%%%%%%%%%%%%%%%%%%%%%%%%%%%%%%%%%%%%%%%%%%%%%%%%%%% 
\begin{rmk} \label{rmk1}
\emph{
The minimal subsets $X\in\Xi(x)$ naturally induce a splitting of $\Pi_N(x)$ into disjoints sets $A(X)$,
in the sense that
$$ \Pi_N(x) = \bigsqcup_{X \in \Xi(x)} A(X).  $$
}
\end{rmk}
%%%%%%%%%%%%%%%%%%%%%%%%%%%%%%%%%%%%%%%%%%%%%%%%%%%%%%%%%%%%%%%%%%%%%%%%%% 
 
We are now ready to prove the following result.

%%%%%%%%%%%%%%%%%%%%%%%%%%%%%%%%%%%%%%%%%%%%%%%%%%%%%%%%%%%%%%%%%%%%%%%%%% 
\begin{lemma} \label{stoc}
If $x\in\RRR^n$ is such that $x_i\ne0$ for all $i=1,\ldots,n$, and $x_1+\ldots+x_n=0$, then one has
\begin{equation} \label{result1}
\sum_{x'\in \Pi_A(x)} G_n(x') = 0 .
\end{equation}
\end{lemma}
%%%%%%%%%%%%%%%%%%%%%%%%%%%%%%%%%%%%%%%%%%%%%%%%%%%%%%%%%%%%%%%%%%%%%%%%%% 

%%%%%%%%%%%%%%%%%%%%%%%%%%%%%%%%%%%%%%%%%%%%%%%%%%%%%%%%%%%%%%%%%%%%%%%%%% 
\prova
We prove the result by induction on $n$.
For $n=2,3$ one has $\Pi_A(x)=\Pi(x)$ hence the result follows by Lemma \ref{1}. 
Take any $n\ge4$ and
consider $x+\e$, with $\e\in E^n$.
%E_0^n$  as in Remark \ref{rmk0}. 
Then one has
$$ 
\sum_{x'\in \Pi_A(x+\e)} G_n(x') = \sum_{x'\in \Pi(x+\e)} G_n(x') = (\e_1+\ldots+\e_n)
\prod_{i=1}^{n}\frac{1}{x_i + \e_i}  
$$
by Lemma \ref{1}.

%Consider the set $\Pi_N(x)$ of the non-allowed vectors in $\Pi(x)$.
%By remark \ref{rmk1}, if $\e\in E_0^n$, witrh $E_0$ as in remark \ref{rmk0}, one has
%%
%$$ \Pi_N(x) = \bigsqcup_{X \in \Xi(x)} A(X).  $$
%%
Write
$$  \sum_{x'\in \Pi(x+\e)} G_n(x') = \sum_{x'\in \Pi_A(x)} G_n(x'+\e) + \sum_{x'\in \Pi_N(x)} G_n(x'+\e) . $$

By construction one has
\[
\begin{aligned}
 \sum_{x'\in \Pi_A(x)} &G_n(x') = \lim_{\e \to 0^+} \sum_{x'\in \Pi_A(x)} G_n(x'+\e) \\
&= \lim_{\e\to0+} \left(  \sum_{x'\in \Pi(x+\e)} G_n(x') - \sum_{x'\in \Pi_N(x)} G_n(x'+\e) \right)
= - \lim_{\e\to0+} \sum_{x'\in \Pi_N(x)} G_n(x'+\e), 
\end{aligned}
\]
so that, thanks to Remark \ref{rmk1}, the identity \eqref{result1} is equivalent to
\begin{equation} \label{result2}
\lim_{\e\to0^+} \sum_{x'\in \Pi_N(x)} G_n(x'+\e) =  \lim_{\e\to0^+}
\sum_{X \in \Xi(x)} \sum_{x' \in A(X) } G_n(x'+\e)=0.
\end{equation}
%

%We prove \eqref{result1} and \eqref{result2} together by induction on $n$,
%observing that they are easily checked to hold for the first non-trivial cases,
%which are $n=2$ and $n=3$, respectively.

Assume inductively \eqref{result2}
to be true up to $n-1\ge3$. Note that if $\Xi(x)=\emptyset$, then $\Pi_A(x)=\Pi(x)$ and again the
result follows from Lemma \ref{1}, so we may assume $\Xi(x)\ne\emptyset$.

Let $X=\{x_{i_1},\ldots,x_{i_p}\} \in \Xi(x)$, any $\pi_\s(x)\in A(X)$ is of the form 
$\pi_\s(x)=(x_{\s(1)},\ldots,x_{\s(n)})$,
with $(\s(1),\ldots,\s(p))$ a permutation of $({i_1},\ldots,{i_p})$ and 
$(\s({p+1}),\ldots,\s(n))$ a permutation of the others labels $(i_{p+1},\ldots,i_n)$.
Set $\xi'=(x_{i_1},\ldots,x_{i_p})$ and $\xi''=(x_{i_{p+1}},\ldots,x_{i_{n}})$; then, for all $\pi_\s(x)\in A(X)$, 
one can write
$x':=\pi_\s(x) =(\eta',\eta'')$ with $\eta'=(x_{\s(1)},\ldots,x_{\s(p)})$ and  $\eta''=(x_{\s({p+1})},\ldots,x_{\s(n)})$.

One has, by Lemma \ref{1},
$$ 
\begin{aligned}
& \sum_{x' \in A(X) } G_n(x'+\e) = \sum_{\eta' \in \Pi(\xi')} \sum_{\eta''\in \Pi(\xi'')} G_p(\eta'+\e') \, \frac{1}{D_p(\eta'+\e')} G_{n-p}(D_p(\eta'+\e')v + \eta''+\e'') \\
& \quad = \Biggl( \, \sum_{j=1}^p \left( x_{i_j} + \e_{i_j} \right) \Biggr) 
\Biggl( \, \prod_{k=1}^p \frac{1}{x_{i_k} +\e_{i_k} } \Biggr) \frac{1}{D_p(\xi'+\e')} \sum_{\eta'' \in \Pi(\xi'')} G_{n-p}(D_p(\xi' + \e ')v + \eta''+\e'') ,
\end{aligned}
$$
where 
$\e'=(\e_{\s(1)},\ldots,\e_{\s(p)})$, $\e''=(\e_{\s({p+1})},\ldots,\e_{\s(n)})$ and $v=(1,0,0,\ldots,0)\in\RRR^{n-p}$.

Since $ \left( x_{i_1} + \e_{i_1} + \ldots + x_{i_p} + \e_{i_p} \right) =D_p(\xi'+\e')$, we can set $\e'=0$ and obtain
$$
\begin{aligned}
\sum_{x' \in A(X) } G_n(x'+(0,\e'')) & = 
\Biggl( \prod_{k=1}^p \frac{1}{x_{i_k} } \Biggr) \sum_{\eta'' \in \Pi(\xi'')} G_{n-p}(D_p(\xi')v+ \eta''+\e'') \\
& = \Biggl( \prod_{k=1}^p \frac{1}{x_{i_k} } \Biggr) \sum_{\eta'' \in \Pi(\xi'')} G_{n-p}( \eta''+\e'') ,
\end{aligned}
$$
since one has $D_p(\xi')=0$.

By the inductive hypothesis one has
$$ \lim_{\e'' \to 0^+} \sum_{\eta'' \in \Pi(\xi'')} G_{n-p}( \eta''+\e'') =
\sum_{\eta'' \in \Pi_A(\xi'')} G_{n-p}( \eta'') = 0 , $$
where the first equality follows from \eqref{result2} and the second one from \eqref{result1} both with $n-p$ instead of $n$.
\EP
%%%%%%%%%%%%%%%%%%%%%%%%%%%%%%%%%%%%%%%%%%%%%%%%%%%%%%%%%%%%%%%%%%%%%%%%%% 

%%%%%%%%%%%%%%%%%%%%%%%%%%%%%%%%%%%%%%%%%%%%%%%%%%%%%%%%%%%%%%%%%%%%%%%%%%
%%%%%%%%%%%%%%%%%%%%%%%%%%%%%%%%%%%%%%%%%%%%%%%%%%%%%%%%%%%%%%%%%%%%%%%%%%
 
%%%%%%%%%%%%%%%%%%%%%%%%%%%%%%%%%%%%%%%%%%%%%%%%%%%%%%%%%%%%%%%%%%%%%%%%%%

\end{document}